\pdfoutput=1
\documentclass[leqno]{article}

\usepackage{amsmath,amsthm}

\usepackage[LGR,T1]{fontenc} % LGR allows for Greek text input
\usepackage[utf8]{inputenc} 
\usepackage[lcgreekalpha,upint]{stix2} % For older TeX distributions, use stix rather than stix2
\usepackage{inconsolata}
\renewcommand{\in}{\smallin}
\renewcommand{\emptyset}{\varnothing}
\usepackage{geometry}
\usepackage{afterpage}
\usepackage{setspace}
\usepackage{xspace}
\usepackage{CJKutf8} % Chinese, Japanese, and Korean

\newcommand{\Japanese}[1]{\begin{CJK}{UTF8}{min} \CJKfamily{goth} \mathchoice{\ensuremath{\text{\footnotesize #1}}}{\ensuremath{\text{\footnotesize #1}}}{\ensuremath{\text{\tiny #1}}}{\ensuremath{\text{\tiny #1}}} \begin{CJK}{UTF8}{min}} % Defines a macro to use Japanese (or any CJK) characters in math mode; font size adjusted to be closer to the standard size for math fonts, and mathchoice makes it so that subscripts etc. scale. It will probably be better to take a smaller size than "tiny" for the subsubscript style…though I can't forsee this being used much

\newcommand{\yo}{\begin{CJK}{UTF8}{min} \CJKfamily{goth} \mathchoice{\ensuremath{\text{\footnotesize よ}}}{\ensuremath{\text{\footnotesize よ}}}{\ensuremath{\text{\footnotesize ょ}}}{\ensuremath{\text{\tiny ょ}}} \end{CJK}}

\renewcommand{\setminus}{\smallsetminus}

\usepackage[shortlabels]{enumitem}
\setlist[enumerate]{
     itemsep  = 0.15cm,
       label  = {\upshape (\arabic*)},
         ref  = \arabic*,
    leftmargin  = *
}

\newcommand{\stlabel}[1]{{\upshape (\ref*{#1}.\arabic*)}}

\newcommand{\enumref}[2]{(\hyperref[#1.#2]{\ref*{#1}.\ref*{#1.#2}})}

\usepackage{imakeidx}
\usepackage{url} % Allows for printing of URLs with TeX symbols that normally wouldn't be printed such as "~"
\usepackage{hyperref}

\usepackage{xcolor}
\definecolor{love}{RGB}{128, 15, 37}
\definecolor{slate}{RGB}{42, 54, 59}

\hypersetup{
      colorlinks=true, %set true if you want colored links
      linktoc=all,     %set to all if you want both sections and sections linked
      linkcolor=love,  %choose some color if you want links to stand out
      citecolor = love,
      urlcolor=slate, % choose some color to set links to websites
      pdfencoding=unicode,
      bookmarksdepth=subsubsection,
 	  	pdfencoding=unicode,
      hypertexnames = true,
  }

\usepackage{bookmark}
\bookmarksetup{
  numbered, 
  open,
}

\usepackage[style=ieee-alphabetic,backend=biber,sorting=nyt,url=false,backref=true]{biblatex} % url=false removes '[Online] Available: _repetition of DOI_'
\AtEveryBibitem{\clearfield{issn}} % Removes ISSN (retains ISBN)
\AtEveryCitekey{\clearfield{issn}}
\AtBeginBibliography{\small} % Sets references font size
\renewcommand*{\bibnamedash}{\underline{\hspace{3em}}\kern 0.1em} % Uses an underline when an author has multiple publications (in line with the AMS style)
\DeclareFieldFormat{postnote}{#1} % Removes the "p." when just a number is referenced
\DeclareFieldFormat[article,inbook,incollection,inproceedings,patent,thesis,unpublished]
  {title}{\textit{#1\isdot}}     % Uniformly italicizes title of all works
\DeclareFieldFormat{journaltitle}{#1}  % Makes Journal title in normal text
\DeclareFieldFormat[book,inbook,incollection,inproceedings]{series}{#1} % Takes away the "ser." in front of Series Title that is part of the ieee styling

\newcommand{\HAsubseclink}[1]{\href{http://www.math.ias.edu/~lurie/papers/HA.pdf\#subsection.#1}{#1}}

\newcommand{\HTTthmlink}[1]{\href{http://www.math.ias.edu/~lurie/papers/HTT.pdf\#theorem.#1}{#1}}

\newcommand{\HTTsec}[1]{\href{http://www.math.ias.edu/~lurie/papers/HTT.pdf\#section.#1}{\S #1}}
\newcommand{\HAsec}[1]{\href{http://www.math.ias.edu/~lurie/papers/HA.pdf\#section.#1}{\S #1}}
\newcommand{\SAGsec}[1]{\href{http://www.math.ias.edu/~lurie/papers/SAG-rootfile.pdf\#section.#1}{\S #1}}

\newcommand{\SAGsubsec}[1]{\href{http://www.math.ias.edu/~lurie/papers/SAG-rootfile.pdf\#subsection.#1}{\S #1}}

\newcommand{\HTTthm}[2]{\href{http://www.math.ias.edu/~lurie/papers/HTT.pdf\#theorem.#2}{#1 #2}}
\newcommand{\HAthm}[2]{\href{http://www.math.ias.edu/~lurie/papers/HA.pdf\#theorem.#2}{#1 #2}}
\newcommand{\HAappthm}[2]{\href{http://www.math.ias.edu/~lurie/papers/HA.pdf\#atheorem.#2}{#1 #2}} % Appendix theorems are numbered differently
\newcommand{\SAGthm}[2]{\href{http://www.math.ias.edu/~lurie/papers/SAG-rootfile.pdf\#theorem.#2}{#1 #2}}

\newcommand{\HTTpage}[1]{\href{http://www.math.ias.edu/~lurie/papers/HTT.pdf\#page.#1}{p. #1}}

\newcommand{\HTT}[2]{\cite[\HTTthm{#1}{#2}]{HTT}}
\newcommand{\HA}[2]{\cite[\HAthm{#1}{#2}]{HA}}
\newcommand{\HAa}[2]{\cite[\HAappthm{#1}{#2}]{HA}} % For theorems in the appendix of HA because the numbering is different
\newcommand{\SAG}[2]{\cite[\SAGthm{#1}{#2}]{SAG}}

\usepackage[noabbrev,nameinlink]{cleveref}		% Uses cleveref package to cite things in the text with names
	\crefname{section}{Chapter}{Chapters}
	\crefname{subsection}{\S\!}{\S\S\!}
	\crefname{axioms}{Axiom}{Axioms}
	\crefname{exercise}{Exercise}{Exercises}
	\crefname{exercisenum}{Exercise}{Exercises}
	\crefname{construction}{Construction}{Constructions}
	\crefname{problem}{Problem}{Problems}
	\crefname{theorem}{Theorem}{Theorems}
	\crefname{definition}{Definition}{Definitions}
	\crefname{proposition}{Proposition}{Propositions}
	\crefname{lemma}{Lemma}{Lemmas}
	\crefname{lem}{Lemma}{Lemmas}
	\crefname{remark}{Remark}{Remarks}
	\crefname{example}{Example}{Examples}
	\crefname{examplealph}{Example}{Examples}
	\crefname{corollary}{Corollary}{Corollaries}
	\crefname{nonexample}{Nonexample}{Nonexamples}
	\crefname{subappendix}{\S\!}{\S\S\!}
	\crefname{part}{Part}{Parts}
	\crefname{notation}{Notation}{Notations}
	\crefname{thm}{Theorem}{Theorems}
	\crefname{defn}{Definition}{Definitions}
	\crefname{cor}{Corollary}{Corollaries}
	\crefname{rmk}{Remark}{Remarks}
	\crefname{observation}{Observation}{Observations}
	\Crefname{observation}{Observation}{Observations}
	\Crefname{rmk}{Remark}{Remarks}
	\Crefname{cor}{Corollary}{Corollaries}
	\crefname{conjecture}{Conjecture}{Conjectures}
	\Crefname{defn}{Definition}{Definitions}
	\Crefname{thm}{Theorem}{Theorems}
	\crefname{Notation}{Notation}{Notations}
	\Crefname{subappendix}{Section}{Sections}
	\Crefname{corollary}{Corollary}{Corollaries}
	\Crefname{axioms}{Axiom}{Axioms}
	\Crefname{exercise}{Exercise}{Exercises}
	\Crefname{exercisenum}{Exercise}{Exercises}
	\Crefname{construction}{Construction}{Constructions}
	\Crefname{problem}{Problem}{Problems}
	\Crefname{theorem}{Theorem}{Theorems}
	\Crefname{definition}{Definition}{Definitions}
	\Crefname{proposition}{Proposition}{Propositions}
	\Crefname{lemma}{Lemma}{Lemmas}
	\Crefname{lem}{Lemma}{Lemmas}
	\Crefname{remark}{Remark}{Remarks}
	\Crefname{example}{Example}{Examples}
	\Crefname{examplealph}{Example}{Examples}
	\Crefname{section}{Chapter}{Chapters}
	\Crefname{subsection}{Section}{Sections}
	\Crefname{summary}{Summary}{Summaries}
	\Crefname{attempt}{Attempt}{Attempts}
	\Crefname{part}{Part}{Parts}
	\Crefname{conjecture}{Conjecture}{Conjectures}

  \crefformat{nul}{(#2#1#3)} % Has cref parenthesize "nul" references
	\crefname{nul}{}{}
	\Crefformat{nul}{(#2#1#3)}
	\Crefname{nul}{}{}

\numberwithin{equation}{subsection}

\swapnumbers % Changes environment style from "Theorem #" to "# Theorem"
\usepackage{etoolbox}
\usepackage{mathtools}

\patchcmd{\swappedhead}{(#3)}{\normalfont #3}{}{} % patch 2 uses
\theoremstyle{definition}

\newtheorem{definition}[equation]{Definition}
\newtheorem{nul}[equation]{}

\newtheorem{construction}[equation]{Construction}

\newtheorem{notation}[equation]{Notation}
\newtheorem{recollection}[equation]{Recollection}
\newtheorem{convention}[equation]{Convention}
\newtheorem{goal}[equation]{Goal}

\newtheorem{warning}[equation]{Warning}
\newtheorem{question}[equation]{Question}
\newtheorem{questions}[equation]{Questions}
\newtheorem{idea}[equation]{Idea}
\newtheorem{remark}[equation]{Remark}

\newtheorem{observation}[equation]{Observation}

\newtheorem{observations}[equation]{Observations}

\newtheorem{example}[equation]{Example}

\newtheorem{exercise}[equation]{Exercise}

\newtheorem{properties}[equation]{Properties}
\newtheorem{attempt}[equation]{Attempt}

\newtheorem*{definition*}{Definition}
\newtheorem*{remark*}{Remark}

\theoremstyle{plain}

\newtheorem{theorem}[equation]{Theorem}

\newtheorem{claim}[equation]{Claim}
\newtheorem{proposition}[equation]{Proposition}
\newtheorem{lemma}[equation]{Lemma}

\newtheorem{conjecture}[equation]{Conjecture}

\newtheorem{corollary}[equation]{Corollary}

\newtheorem{cor}[equation]{Corollary}
\newtheorem{lem}[equation]{Lemma}

\newtheorem*{theorem*}{Theorem}

\theoremstyle{definition}

\newtheorem{quest}[equation]{Question}

\newtheorem{recall}[equation]{Recall}

\usepackage{tikz}
\usepackage{tikz-cd}
\usepackage{graphicx}
\usepackage[all,arc]{xy}
\SelectTips{cm}{10}

\usetikzlibrary{matrix,arrows,arrows.meta,bending,calc,decorations,decorations.pathmorphing,shapes}

  \tikzset{
    RRight/.tip={Glyph[glyph math command=arrowkitHeadright]},
    double distance={0.135em},
    commutative diagrams/.cd, 
    arrow style=tikz, 
  }

  \pgfdeclarearrow{
      name = minionto,
      parameters = { \the\pgfarrowlength\the\pgfarrowwidth},
      setup code = {
          \pgfarrowssettipend{\pgfarrowlength}
          \pgfarrowssetbackend{0.00512821\pgfarrowlength}
          \pgfarrowssetlineend{0.774369\pgfarrowlength}
          \pgfarrowssetvisualbackend{0.574359\pgfarrowlength}
          \pgfarrowssetvisualtipend{\pgfarrowlength}
          \pgfarrowsupperhullpoint{0cm}{0.914141\pgfarrowwidth}
          \pgfarrowsupperhullpoint{0cm}{0.469697\pgfarrowwidth}
          \pgfarrowsupperhullpoint{0.107692\pgfarrowlength}{\pgfarrowwidth}
          \pgfarrowshullpoint{0.574359\pgfarrowlength}{0cm}
          \pgfarrowshullpoint{\pgfarrowlength}{0cm}
          \pgfarrowssavethe\pgfarrowlength
          \pgfarrowssavethe\pgflinewidth
      },
      drawing code = {
          \pgfpathmoveto{\pgfpoint{0cm}{-0.914141\pgfarrowwidth}}
          \pgfpathlineto{\pgfpoint{0.574359\pgfarrowlength}{0cm}}
          \pgfpathlineto{\pgfpoint{0cm}{0.914141\pgfarrowwidth}}        
          \pgfpathlineto{\pgfpoint{0.107692\pgfarrowlength}{\pgfarrowwidth}}
          \pgfpathcurveto{\pgfpoint{0.528205\pgfarrowlength}{0.5\pgfarrowwidth}}
                         {\pgfpoint{0.548718\pgfarrowlength}{0.469697\pgfarrowwidth}}
                         {\pgfpoint{\pgfarrowlength}{0cm}}
          \pgfpathcurveto{\pgfpoint{0.528205\pgfarrowlength}{-0.5\pgfarrowwidth}}
                         {\pgfpoint{0.548718\pgfarrowlength}{-0.469697\pgfarrowwidth}}
                         {\pgfpoint{0.107692\pgfarrowlength}{-\pgfarrowwidth}}
          \pgfpathclose
          \pgfusepathqfill
      },
      defaults = {length=0.22em, width=0.19em}
  }

  \pgfarrowsdeclarereversed{miniontail}{miniontail}{minionto}{minionto}

  \pgfdeclarearrow{
      name = minionhook,
      parameters = { \the\pgflinewidth},
      setup code = {
          \pgfmathsetlength{\pgfarrowlength}{4.48\pgflinewidth}
          \pgfarrowssetlineend{\pgflinewidth}
          \pgfarrowssettipend{\pgfarrowlength}
          \pgfarrowsupperhullpoint{0cm}{0.914141\pgflinewidth}
          \pgfarrowsupperhullpoint{0cm}{0.469697\pgflinewidth}
          \pgfarrowsupperhullpoint{0.107692\pgfarrowlength}{\pgflinewidth}
          \pgfarrowshullpoint{0.574359\pgfarrowlength}{0cm}
          \pgfarrowshullpoint{\pgfarrowlength}{0cm}
          \pgfarrowssavethe\pgflinewidth
          \pgfarrowssavethe\pgfarrowlength
      },
      drawing code = {
          \pgfpathmoveto{\pgfpoint{0.25\pgfarrowlength}{-3.84\pgflinewidth}}
          \pgfpathlineto{\pgfpoint{0.25\pgfarrowlength}{-2.86\pgflinewidth}}
          \pgfpathcurveto{\pgfpoint{0.73660714\pgfarrowlength}{-2.86\pgflinewidth}}
                         {\pgfpoint{0.73660714\pgfarrowlength}{-1.9\pgflinewidth}}
                         {\pgfpoint{0.73660714\pgfarrowlength}{-1.6\pgflinewidth}}
          \pgfpathcurveto{\pgfpoint{0.73660714\pgfarrowlength}{-0.54\pgflinewidth}}
                         {\pgfpoint{0.35267857\pgfarrowlength}{-0.45\pgflinewidth}}
                         {\pgfpoint{0.0\pgfarrowlength}{-0.51\pgflinewidth}}
          \pgfpathlineto{\pgfpoint{0.0\pgfarrowlength}{0.5\pgflinewidth}}
          \pgfpathcurveto{\pgfpoint{0.97767857\pgfarrowlength}{0.7\pgflinewidth}}
                         {\pgfpoint{0.99107143\pgfarrowlength}{-1.2\pgflinewidth}}
                         {\pgfpoint{1.0\pgfarrowlength}{-1.6\pgflinewidth}}
          \pgfpathcurveto{\pgfpoint{0.95982143\pgfarrowlength}{-3.86\pgflinewidth}}
                         {\pgfpoint{0.39732143\pgfarrowlength}{-3.86\pgflinewidth}}
                         {\pgfpoint{0.25\pgfarrowlength}{-3.86\pgflinewidth}}
          \pgfpathclose
          \pgfusepathqfill
      },
      defaults = {length=0.22em, width=0.19em}
  }

  \pgfdeclarearrow{
      name = minionhook',
      parameters = { \the\pgflinewidth},
      setup code = {
          \pgfmathsetlength{\pgfarrowlength}{4.48\pgflinewidth}
          \pgfarrowssetlineend{\pgflinewidth}
          \pgfarrowssettipend{\pgfarrowlength}
          \pgfarrowsupperhullpoint{0cm}{0.914141\pgflinewidth}
          \pgfarrowsupperhullpoint{0cm}{0.469697\pgflinewidth}
          \pgfarrowsupperhullpoint{0.107692\pgfarrowlength}{\pgflinewidth}
          \pgfarrowshullpoint{0.574359\pgfarrowlength}{0cm}
          \pgfarrowshullpoint{\pgfarrowlength}{0cm}
          \pgfarrowssavethe\pgflinewidth
          \pgfarrowssavethe\pgfarrowlength
      },
      drawing code = {
          \pgfpathmoveto{\pgfpoint{0.25\pgfarrowlength}{3.84\pgflinewidth}}
          \pgfpathlineto{\pgfpoint{0.25\pgfarrowlength}{2.86\pgflinewidth}}
          \pgfpathcurveto{\pgfpoint{0.73660714\pgfarrowlength}{2.86\pgflinewidth}}
                         {\pgfpoint{0.73660714\pgfarrowlength}{1.9\pgflinewidth}}
                         {\pgfpoint{0.73660714\pgfarrowlength}{1.6\pgflinewidth}}
          \pgfpathcurveto{\pgfpoint{0.73660714\pgfarrowlength}{0.54\pgflinewidth}}
                         {\pgfpoint{0.35267857\pgfarrowlength}{0.45\pgflinewidth}}
                         {\pgfpoint{0.0\pgfarrowlength}{0.51\pgflinewidth}}
          \pgfpathlineto{\pgfpoint{0.0\pgfarrowlength}{-0.5\pgflinewidth}}
          \pgfpathcurveto{\pgfpoint{0.97767857\pgfarrowlength}{-0.7\pgflinewidth}}
                         {\pgfpoint{0.99107143\pgfarrowlength}{1.2\pgflinewidth}}
                         {\pgfpoint{1.0\pgfarrowlength}{1.6\pgflinewidth}}
          \pgfpathcurveto{\pgfpoint{0.95982143\pgfarrowlength}{3.86\pgflinewidth}}
                         {\pgfpoint{0.39732143\pgfarrowlength}{3.86\pgflinewidth}}
                         {\pgfpoint{0.25\pgfarrowlength}{3.86\pgflinewidth}}
          \pgfpathclose
          \pgfusepathqfill
      },
      defaults = {length=0.22em, width=0.19em}
  }

  \tikzset{
    >/.tip={minionto},
    >->/.style = {miniontail-minionto},
    <-</.style = {minionto-miniontail},
    <->/.style = {minionto-miniontail},
    hooked/.style = {minionhook-minionto},
    hooked'/.style = {minionhook'-minionto},
  }

  \tikzset{
  RRight/.tip={Glyph[glyph math command=arrowkitHeadright]},
}

\pgfdeclaredecoration{single line}{initial}{\state{initial}[width=\pgfdecoratedpathlength-1sp]{\pgfpathmoveto{\pgfpointorigin}}\state{final}{\pgfpathlineto{\pgfpointorigin}}}

\tikzcdset{  crossing over/.style={
    /tikz/preaction={
      /tikz/draw,
      /tikz/color=\pgfkeysvalueof{/tikz/commutative diagrams/background color},
      /tikz/arrows=-,
      /tikz/line width=\pgfkeysvalueof{/tikz/commutative
      diagrams/crossing over clearance}}}}

\tikzset{curve/.style={settings={#1},to path={(\tikztostart)
    .. controls ($(\tikztostart)!\pv{pos}!(\tikztotarget)!\pv{height}!270:(\tikztotarget)$)
    and ($(\tikztostart)!1-\pv{pos}!(\tikztotarget)!\pv{height}!270:(\tikztotarget)$)
    .. (\tikztotarget)\tikztonodes}},
    settings/.code={\tikzset{quiver/.cd,#1}
        \def\pv##1{\pgfkeysvalueof{/tikz/quiver/##1}}},
    quiver/.cd,pos/.initial=0.35,height/.initial=0}

\newcommand{\andeq}{\text{\qquad and\qquad}}
\newcommand{\period}{\rlap{\ .}}
\newcommand{\comma}{\rlap{\ ,}}
\newcommand{\semicolon}{\rlap{\ ;}}
\renewcommand{\ast}{\mathbin{*}}

\renewcommand{\square}{\mdwhtsquare} 

\newcommand{\D}{\categ{D}}

\newcommand{\CC}{\mathbb{C}}
\newcommand{\NN}{\mathbb{N}}
\newcommand{\RR}{\mathbb{R}}
\newcommand{\QQ}{\mathbb{Q}}
\newcommand{\ZZ}{\mathbb{Z}}

\newcommand{\FF}{\mathbb{F}}

\newcommand{\colonequals}{\coloneq}
\newcommand{\cross}{\times}
\newcommand{\tensor}{\otimes}

\newcommand{\tensorhat}{\otimeshat}
\newcommand{\isomorphic}{\cong}
\newcommand{\union}{\cup}
\newcommand{\of}{\circ}
\newcommand{\idnosub}{\mathrm{id}}
\newcommand{\id}[1]{\operatorname{id}_{#1}}

\newcommand{\crosslimits}{\operatornamewithlimits{\cross}}

\newcommand{\Directsumhat}{\operatornamewithlimits{\widehat{\bigoplus}}}

\newcommand{\categ}{\msf}

\newcommand{\TT}{\bb{T}}

\newcommand{\Ucal}{\cal{U}}

\DeclareMathOperator{\basic}{basic}

\DeclareMathOperator{\Conf}{Conf}
\DeclareMathOperator{\Cent}{Cent}
\DeclareMathOperator{\Sym}{Sym}
\newcommand{\CS}{\mathrm{CS}}
\newcommand{\cs}{\mathrm{cs}}

\DeclareMathOperator{\pt}{pt}

\DeclareMathOperator{\ku}{ku}
\DeclareMathOperator{\ko}{ko}
\newcommand{\kuhat}{\widehat{\ku}}

\DeclareMathOperator{\dlog}{dlog}

\DeclareMathOperator{\odd}{odd}

\DeclareMathOperator{\rank}{rank}

\DeclareMathOperator{\Hol}{Hol}

\newcommand{\HZZ}{\mathrm{H}\ZZ}
\newcommand{\HZZhat}{\widehat{\HZZ}}

\newcommand{\HRR}{\mathrm{H}\RR}
\newcommand{\HR}{\mathrm{H}R}

\newcommand{\HCC}{\mathrm{H}\CC}

\newcommand{\Omegabullet}{\Omega^{\bullet}}
\newcommand{\HOmega}{\mathrm{H}\Omega}
\newcommand{\HOmegabullet}{\mathrm{H}\Omegabullet}

\DeclareMathOperator{\cl}{cl}
\newcommand{\Omegacl}{\Omega_{\cl}}
\newcommand{\OmegadR}{\Omega_{\dR}}
\newcommand{\HOmegacl}{\mathrm{H}\Omegacl}

\DeclareMathOperator*{\colim}{colim}
\newcommand{\fib}{\ensuremath{\textup{fib}}}
\newcommand{\cofib}{\ensuremath{\textup{cofib}}}

\newcommand{\forget}{\mathrm{forget}}

\newcommand{\cupprod}{\smallsmile}

\newcommand{\uppi}{\mathrm{\pi}}

\newcommand{\counit}{\varepsilon}
\newcommand{\unit}{\eta}

\newcommand{\conn}{A}
\newcommand{\curvature}[1]{F_{#1}}

\DeclareMathOperator{\GL}{GL}
\DeclareMathOperator{\PGL}{PGL}
\DeclareMathOperator{\PSL}{PSL}
\DeclareMathOperator{\SL}{SL}
\DeclareMathOperator{\SO}{SO}
\DeclareMathOperator{\Or}{O}
\DeclareMathOperator{\PU}{PU}
\DeclareMathOperator{\Spin}{Spin}
\DeclareMathOperator{\SU}{SU}
\DeclareMathOperator{\LSU}{LSU}
\DeclareMathOperator{\LU}{LU}
\newcommand{\SLtilde}{\widetilde{\SL}}
\DeclareMathOperator{\Uup}{U}

\newcommand{\LM}{\mathrm{L}M}
\newcommand{\LG}{\mathrm{L}G}
\newcommand{\LGtilde}{\widetilde{\mathrm{L}}G}
\newcommand{\LGtildeplus}{\widetilde{\mathrm{L}}G^{+}}

\newcommand{\g}{\mathfrak{g}}
\newcommand{\gtilde}{\tilde{\g}}

\newcommand{\ufrak}{\mathfrak{u}}
\newcommand{\pfrak}{\mathfrak{p}}
\newcommand{\afrak}{\mathfrak{a}}
\newcommand{\ofrak}{\mathfrak{o}}
\newcommand{\kfrak}{\mathfrak{k}}
\newcommand{\kdual}{\kfrak^{\vee}}

\newcommand{\psl}{\mathfrak{psl}}
\newcommand{\slfrak}{\mathfrak{sl}}
\newcommand{\gl}{\mathfrak{gl}}
\newcommand{\so}{\mathfrak{so}}
\newcommand{\su}{\mathfrak{su}}

\DeclareMathOperator{\heis}{Heis}
\DeclareMathOperator{\witt}{Witt}
\DeclareMathOperator{\vir}{Vir}
\newcommand{\wittCC}{\witt_{\CC}}
\newcommand{\virCC}{\vir_{\CC}}
\newcommand{\wittRR}{\witt_{\RR}}
\newcommand{\virRR}{\vir_{\RR}}
\DeclareMathOperator{\Fock}{Fock}

\DeclareMathOperator{\Lg}{L\g}
\newcommand{\Lgtilde}[1]{\widetilde{\Lup}_{#1}\g}
\newcommand{\Lgtildeno}{\widetilde{\Lup}\g}

\DeclareMathOperator{\Cas}{Cas}

\newcommand{\Univ}{\mathcal{U}}

\newcommand{\MC}{\mathrm{MC}}

\DeclareMathOperator{\BGL}{BGL}
\DeclareMathOperator{\BO}{BO}
\newcommand{\BP}{\mathrm{B}P}
\DeclareMathOperator{\BSO}{BSO}
\DeclareMathOperator{\BSL}{BSL}
\DeclareMathOperator{\BU}{BU}
\DeclareMathOperator{\BSU}{BSU}

\newcommand{\BbulletGL}[2]{\Bun_{\GL_{#1}(#2)}}
\newcommand{\BbulletGLp}[2]{\Bun_{\GL_{#1}^+(#2)}}

\newcommand{\EGamma}{\mathrm{E}\Gamma}
\newcommand{\BGamma}{\mathrm{B}\Gamma}
\newcommand{\BbulletGamma}{\Bun_{\Gamma}}
\newcommand{\EbulletGamma}{\mathrm{E}_{\bullet}\Gamma}
\newcommand{\EG}{\mathrm{E}G}
\newcommand{\BG}{\mathrm{B}G}

\newcommand{\BK}{\mathrm{B}K}

\newcommand{\BbulletG}{\BunG}

\newcommand{\EnablaG}{\mathrm{E}_{\nabla}G}
\newcommand{\BunGnablatriv}{\BunG^{\kern-0.15em\nabla, \mathrm{triv}}}

\newcommand{\disc}{\mathrm{disc}}
\newcommand{\Gdisc}{G^{\disc}}
\newcommand{\BGdisc}{\mathrm{B}G^{\disc}}

\DeclareMathOperator{\triv}{triv}
\newcommand{\trivG}{\triv_{G}}
\newcommand{\trivGnabla}{\triv_{G}^{\kern-0.15em\nabla}}

\newcommand{\VectRR}{\operatorname{Vect}_{\RR}}

\newcommand{\Bun}{\mathrm{Bun}}
\newcommand{\BunG}{\Bun_{G}}

\newcommand{\BunGnabla}{\BunG^{\kern-0.15em\nabla}}
\newcommand{\Bunnabla}{\Bun^{\kern-0.15em\nabla}}

\newcommand{\BunOnabla}[1]{\Bun_{\Or_{#1}}^{\kern-0.15em\nabla}}
\newcommand{\BunUnabla}[1]{\Bun_{\Uup_{#1}}^{\kern-0.15em\nabla}}

\DeclareMathOperator{\MTH}{MTH}
\DeclareMathOperator{\KO}{KO}
\newcommand{\KOhat}{\widehat{\KO}{}}
\newcommand{\IZ}{\mathrm{I}_{\ZZ}}
\newcommand{\IZhat}{\hat{\mathrm{I}}_{\ZZ}}

\newcommand{\Ahat}{\widehat{A}}

\newcommand{\Gr}{\mathrm{Gr}}
\newcommand{\St}{\mathrm{St}}

\newcommand{\ev}{\mathrm{ev}}

\newcommand{\sLine}{\categ{sLine}}
\newcommand{\Line}{\categ{Line}}

\newcommand{\arrow}{\ar}
\newcommand{\msf}{\mathsf}
\newcommand{\bb}{\mathbb}
\newcommand{\mrm}{\mathrm}
\newcommand{\del}{\partial}
\renewcommand{\cal}{\mathcal}
\newsavebox{\pullbacksq}
\sbox\pullbacksq{%
\begin{tikzpicture}%
\draw (0,0) -- (1ex,0ex);%
\draw (1ex,0ex) -- (1ex,1ex);%
\end{tikzpicture}}

\DeclareMathOperator{\PT}{PT}

\let\originalleft\left
\let\originalright\right
\renewcommand{\left}{\mathopen{}\mathclose\bgroup\originalleft}
\renewcommand{\right}{\aftergroup\egroup\originalright}

\DeclareMathOperator\Tot{Tot}
\DeclareMathOperator\Kos{Kos}
\newcommand\Vect{\mathsf{Vect}}
\DeclareMathOperator{\pr}{pr}
\newcommand{\modmod}{\mathbin{\sslash}}
\newcommand{\exterior}{\Lambda}
\newcommand{\exteriorhat}{\widehat{\exterior}}

\newcommand{\Fcal}{\cal{F}}

\newcommand{\Vdual}{V^{\vee}}
\newcommand{\Wdual}{W^{\vee}}

\DeclareMathOperator{\Jet}{Jet}

\DeclareMathOperator{\rot}{rot}
\newcommand{\TTrot}{\TT_{\rot}}

\DeclareMathOperator{\ce}{cent}
\newcommand{\TTce}{\TT_{\ce}}

\DeclareMathOperator{\tmf}{tmf}

\renewcommand{\H}{\mathrm{H}}
\newcommand{\Hc}{\H_{\mathrm{c}}}

\newcommand{\Har}{\mathrm{Har}}

\DeclareMathOperator{\dR}{dR}
\newcommand{\HdR}{\H_{\dR}}

\newcommand{\Hhat}{\hat{\H}{}}
\newcommand{\Ehat}{\hat{E}}
\newcommand{\Phat}{\hat{P}}

\newcommand{\HG}{\H_{G}}
\newcommand{\HGdR}{\H_{G,\dR}}

\DeclareMathOperator{\sm}{sm}
\DeclareMathOperator{\sing}{sing}
\newcommand{\Zsm}{\operatorname{Z}^{\sm}}
\newcommand{\Csm}{\operatorname{C}^{\sm}}
\newcommand{\Csing}{\operatorname{C}_{\sing}}

\DeclareMathOperator{\cont}{cont}
\newcommand{\Hcont}{\H_{\cont}}

\renewcommand{\Cup}{\operatorname{C}}

\DeclareMathOperator{\Lie}{Lie}
\newcommand{\HLie}{\H_{\Lie}}

\newcommand{\Omegac}{\Omega_{\operatorname{c}}}

\newcommand{\dup}{\mathrm{d}}
\renewcommand{\d}{\mathrm{d}}

\newcommand{\curv}{\mathrm{curv}}
\newcommand{\im}{\mathrm{im}}
\newcommand{\image}{\mathrm{im}}

\DeclareMathOperator{\ch}{ch}
\DeclareMathOperator{\cc}{cc}
\DeclareMathOperator{\Arf}{Arf}

\DeclareMathOperator{\K}{K}
\DeclareMathOperator{\KU}{K}
\newcommand{\lowerstar}{_{\ast}}
\newcommand{\upperstar}{^{\ast}}
\newcommand{\fupperstar}{f\upperstar}
\newcommand{\gupperstar}{g\upperstar}

\newcommand{\fuppersharp}{f^{\sharp}}
\newcommand{\xupperstar}{x\upperstar}
\newcommand{\prupperstar}{\pr\upperstar}

\newcommand{\phiupperstar}{\phi\upperstar}

\newcommand{\uupperstar}{u\upperstar}
\newcommand{\jupperstar}{j\upperstar}
\newcommand{\jlowerstar}{j\lowerstar}
\newcommand{\jlowershriek}{j_!}

\newcommand{\hupperstar}{h\upperstar}

\newcommand{\iupperstar}{i\upperstar}
\newcommand{\ilowerstar}{i\lowerstar}

\newcommand{\iuppershriek}{i^!}

\newcommand{\jhat}{\hat{\jmath}}

\newcommand{\Fhat}{\hat{F}}

\newcommand{\ulowershriek}{u_{!}}

\DeclareMathOperator{\vol}{vol}
\DeclareMathOperator{\Gal}{Gal}

\newcommand{\ctilde}{\tilde{c}}
\newcommand{\ptilde}{\tilde{p}}
\newcommand{\chitilde}{\tilde{\chi}}
\newcommand{\iotatilde}{\tilde{\iota}}

\DeclareMathOperator{\Ex}{Ex}

\DeclareMathOperator{\Ext}{Ext}
\newcommand{\paren}[1]{\left( #1 \right)}

\DeclareMathOperator{\Lhi}{L_{hi}}
\DeclareMathOperator{\Rhi}{R_{hi}}
\newcommand{\LRR}{\Lup_{\RR}}
\renewcommand{\AA}{\mathbb{A}}
\newcommand{\PP}{\mathbb{P}}
\DeclareMathOperator{\Nis}{Nis}
\DeclareMathOperator{\mot}{mot}
\newcommand{\LAA}{\Lup_{\AA^1}}
\newcommand{\Lmot}{\Lup_{\mot}}
\newcommand{\Sm}{\categ{Sm}}

\newcommand{\Gammainverse}{\Gamma^{-1}}
\newcommand{\Gammaupperstar}{\Gamma^{\ast}}
\newcommand{\Gammalowerstar}{\Gamma_{\ast}}
\newcommand{\Gammauppersharp}{\Gamma^{\sharp}}
\newcommand{\Gammalowersharp}{\Gamma_{\sharp}}
\newcommand{\Piinf}{\Pi_{\infty}}
\newcommand{\Bup}{\mathrm{B}}

\newcommand{\Fup}{\mathrm{F}}
\newcommand{\Hup}{\mathrm{H}}
\newcommand{\Iup}{\mathrm{I}}
\newcommand{\Rup}{\mathrm{R}}
\newcommand{\Sup}{\mathrm{S}}

\newcommand{\EM}{\mathrm{H}}

\newcommand{\SMan}{\Sup_{\Mfld}}

\DeclareMathOperator{\swap}{swap}
\DeclareMathOperator{\mult}{mult}

\newcommand{\real}[1]{|#1|}
\newcommand{\setbar}[2]{\{\, #1 \, | \, #2 \,\}}

\DeclareMathOperator{\alg}{alg}
\newcommand{\Deltaalg}{\Delta_{\alg}}
\newcommand{\Deltaalgdot}{\Deltaalg^{\bullet}}

\newcommand{\flowerstar}{f_{\ast}}
\newcommand{\finverse}{f^{-1}}
\newcommand{\fhat}{\hat{f}}
\newcommand{\fbar}{\bar{f}}

\newcommand{\pupperstar}{p^{\ast}}
\newcommand{\plowerstar}{p_{\ast}}
\newcommand{\pinverse}{p^{-1}}

\newcommand{\sigmainverse}{\sigma^{-1}}

\newcommand{\el}{\ell}
\newcommand{\rbar}{\bar{r}}
\newcommand{\elbar}{\bar{\el}}

\newcommand{\coproduct}{\sqcup}
\newcommand{\directsum}{\oplus}

\newcommand{\restrict}[2]{#1|_{#2}}

\newcommand{\all}{{\ensuremath{\textup{all}}}}
\newcommand{\fcov}{{\ensuremath{\textup{fcov}}}}
\newcommand{\loc}{{\ensuremath{\textup{loc}}}}
\newcommand{\inj}{{\ensuremath{\textup{inj}}}}
\newcommand{\foldtext}{{\ensuremath{\textup{fold}}}}

\newcommand{\PShhi}{\PSh_{\RR}}
\newcommand{\Shhi}{\Sh_{\RR}}
\DeclareMathOperator{\pu}{pu}
\newcommand{\pure}{\pu}
\newcommand{\Shpu}{\Sh_{\pu}}
\newcommand{\Shpure}{\Shpu}

\newcommand{\equivalent}{\simeq}

\newcommand{\Uu}{\mathcal{U}}
\newcommand{\cA}{\mathcal A}

\newcommand{\intersect}{\cap}

\newcommand{\Union}{\bigcup}

\newcommand{\Euc}{\categ{Euc}}
\newcommand{\Eucop}{\Euc^{\op}}

\newcommand{\Xdot}{X_{\bullet}}
\newcommand{\Ydot}{Y_{\bullet}}

\newcommand{\X}{\categ{X}}
\newcommand{\Y}{\categ{Y}}

\newcommand{\Ucat}{\categ{U}}
\newcommand{\Xcat}{\categ{X}}

\newcommand{\Zcat}{\categ{Z}}

\DeclareMathOperator{\can}{can}
\DeclareMathOperator{\post}{post}
\DeclareMathOperator{\hyp}{hyp}

\newcommand{\trun}{\mathrm{\tau}}

\DeclareMathOperator{\Bock}{Bock}

\DeclareMathOperator{\Spine}{Spine}

\newcommand{\Eff}{\ensuremath{\textup{Eff}}}

\newcommand{\PrL}{\categ{Pr}^{\mathrm{L}}}

\newcommand{\Nerve}{\mathrm{N}}
\renewcommand{\Bar}{\operatorname{Bar}}
\newcommand{\cn}{\mathrm{cn}}

\newcommand{\Top}{\categ{Top}}

\newcommand{\Fre}{\categ{Fr\acute{e}}}

\newcommand{\Mfld}{\categ{Mfld}}

\newcommand{\GrAb}{\categ{GrAb}}
\newcommand{\Ab}{\categ{Ab}}
\newcommand{\Mod}{\categ{Mod}}
\newcommand{\LMod}{\categ{LMod}}
\newcommand{\Sp}{\categ{Spt}}
\newcommand{\Spt}{\Sp}
\newcommand{\Spc}{\categ{Spc}}
\newcommand{\Cat}{\categ{Cat}}
\newcommand{\Catalg}{\categ{Cat}^{\alg}}
\newcommand{\Gpdalg}{\categ{Gpd}^{\alg}}
\newcommand{\Space}{\Spc}
\newcommand{\PSh}{\categ{PSh}}
\newcommand{\Sh}{\categ{Sh}}
\newcommand{\LC}{\categ{LC}}
\newcommand{\Shv}{\categ{Sh}}
\newcommand{\sSet}{\categ{sSet}}
\newcommand{\Set}{\categ{Set}}
\newcommand{\Ch}{\categ{Ch}}
\newcommand{\LieAlg}{\categ{Lie}}
\newcommand{\Gpd}{\categ{Gpd}}
\newcommand{\Grp}{\categ{Grp}}
\newcommand{\Mon}{\categ{Mon}}
\newcommand{\DDelta}{\mathbf{\Delta}}
\newcommand{\Deltaplus}{\DDelta_{+}}
\newcommand{\Deltaplusop}{\DDelta_{+}^{\op}}

\newcommand{\Cor}{\categ{Cor}}
\newcommand{\Span}{\categ{Span}}
\DeclareMathOperator{\Mor}{Mor}
\DeclareMathOperator{\target}{target}
\newcommand{\Corfcov}{\Cor_{\fcov}(\Mfld)}
\newcommand{\Corfold}{\Cor_{\foldtext}(\Mfld)}

\newcommand{\CMon}{\categ{CMon}}
\newcommand{\Bord}{\categ{Bord}}

\DeclareMathOperator{\Open}{Open}

\newcommand{\fC}{\categ{C}}

\newcommand{\Fin}{\categ{Fin}}
\newcommand{\Finstar}{\Fin_{\ast}}
\newcommand{\Finop}{\Fin^{\op}}

\DeclareMathOperator{\Fun}{Fun}
\DeclareMathOperator{\Funlim}{Fun^{lim}}

\newcommand{\Funcross}{\Fun^{\cross}}
\DeclareMathOperator{\op}{op}
\newcommand{\Mfldop}{\Mfld^{\op}}
\newcommand{\Deltaop}{\DDelta^{\op}}
\newcommand{\Deltabf}{\DDelta}

\DeclareMathOperator{\RHom}{RHom}
\DeclareMathOperator{\Hom}{Hom}
\DeclareMathOperator{\Map}{Map}
\newcommand{\Mapsm}{\Map_{\sm}}
\newcommand{\Mapcont}{\Map_{\cont}}
\DeclareMathOperator{\Nm}{Nm}
\newcommand{\QCoh}{\categ{QCoh}}
\newcommand{\Coh}{\categ{Coh}}

\newcommand{\RRP}{\RR\mathrm{P}}
\newcommand{\CCP}{\CC\mathrm{P}}
\newcommand{\HHP}{\mathbb{H}\mathrm{P}}

\newcommand{\Eone}{\mathbb{E}_{1}}

\newcommand{\Einf}{\mathbb{E}_{\infty}}

\newcommand{\orthogonal}{\ensuremath{\rotatebox[origin=c]{90}{$\scriptstyle\vdash$}}}
\newcommand{\smallorthogonal}{\ensuremath{\rotatebox[origin=c]{90}{$\scriptscriptstyle\vdash$}}}
\newcommand{\Zlorth}{^{\orthogonal}\Zcat}
\newcommand{\Zrorth}{\Zcat^{\orthogonal}}

\newcommand{\Cinf}{\mathrm{C}^{\infty}}

\newcommand{\gdual}{\g^{\vee}}

\DeclareMathOperator{\Kil}{Kil}

\newcommand{\Tup}{\mathrm{T}}
\newcommand{\Tan}{\Tup}
\newcommand{\TanM}{\Tan M}

\newcommand{\Norm}{\mathrm{N}}
\newcommand{\NormM}{\Norm M}

\newcommand{\Lup}{\mathrm{L}}
\newcommand{\Kup}{\mathrm{K}}
\newcommand{\Tstar}{\Tup^{\ast}}

\DeclareMathOperator{\Aut}{Aut}
\DeclareMathOperator{\End}{End}
\DeclareMathOperator{\ad}{ad}
\DeclareMathOperator{\Ad}{Ad}
\DeclareMathOperator{\tr}{tr}
\DeclareMathOperator{\Tr}{tr}
\DeclareMathOperator{\Det}{det}
\DeclareMathOperator{\Pf}{pf}

\DeclareMathOperator{\Tors}{Tors}

\newcommand{\Def}{\mathrm{Def}}
\newcommand{\Cyc}{\mathrm{Cyc}}

\newcommand{\Diff}{\ensuremath{\textup{Diff}}}
\newcommand{\Diffplus}{\ensuremath{\textup{Diff}}{}^{\,+}}
\newcommand{\Sph}[1]{\mathrm{S}^{#1}}
\newcommand{\Circ}{\Sph{1}}
\newcommand{\Vir}{{\widetilde{\Diff}{}^{\,+}}(\Circ)}
\DeclareMathOperator{\chg}{chg}

\newcommand{\ang}[1]{\langle #1\rangle}

\providecommand{\vcentcolon}{\mathrel{\mathop{:}}}
\newcommand{\normord}[1]{\vcentcolon\mathrel{#1}\vcentcolon}

\newcommand{\leftadjoint}{\shortlefttack}

\newcommand{\fromto}[2]{{#1}\to{#2}}
\newcommand{\isomto}[2]{{#1}\similarrightarrow{#2}}
\newcommand{\equivto}[2]{{#1}\similarrightarrow{#2}}
\newcommand{\surjto}[2]{{#1}\twoheadrightarrow{#2}}
\newcommand{\incto}[2]{{#1}\hookrightarrow{#2}}
\newcommand{\into}[2]{{#1}\rightarrowtail{#2}}
\newcommand{\goesto}[2]{{#1}\mapsto{#2}}
\def\reflectedop#1#2{\mathop{\reflectbox{$#1{#2}$}}}
\def\noloc{{\mathpalette\reflectedop\colon}}
\newcommand{\adjto}[4]{{#1}\colon{#2} \rightleftarrows {#3}\noloc{#4}}

\newcommand{\isomorphism}{\similarrightarrow}
\newcommand{\equivalence}{\similarrightarrow}
\newcommand{\inclusion}{\hookrightarrow}

\newcommand{\RRinvariant}{\xspace{$\RR$-in\-vari\-ant}\xspace}
\newcommand{\RRinvariance}{\xspace{$\RR$-in\-vari\-ance}\xspace}
\newcommand{\RRlocalization}{\xspace{$\RR$-lo\-cal\-iza\-tion}\xspace}

\newcommand{\Ktheory}{\xspace{$\Kup$-the\-o\-ry}\xspace}

\newcommand{\category}{\xspace{$\infty$-cat\-e\-go\-ry}\xspace}
\newcommand{\categories}{\xspace{$\infty$-cat\-e\-gories}\xspace}
\newcommand{\categorical}{\xspace{$\infty$-cat\-e\-gor\-i\-cal}\xspace}

\newcommand{\acategory}{{an $\infty$-cat\-e\-go\-ry}\xspace}
\newcommand{\Acategory}{{An $\infty$-cat\-e\-go\-ry}\xspace}

\newcommand{\topos}{\xspace{$\infty$-topos}\xspace}
\newcommand{\topoi}{\xspace{$\infty$-topoi}\xspace}

\newcommand{\atopos}{{an $\infty$-topos}\xspace}
\newcommand{\Atopos}{{An $\infty$-topos}\xspace}

\newcommand{\asite}{{an $\infty$-site}\xspace}

\newcommand{\Cech}{\check{\operatorname{C}}{}}
\newcommand{\Hcech}{\check{\operatorname{H}}{}}

\newcommand{\chat}{\hat{c}}
\newcommand{\phat}{\hat{p}}
\newcommand{\ehat}{\hat{e}}

\usepackage{todonotes}

\indexsetup{toclevel=part} % Adds indices to table of contents as 'parts'
\makeindex[name=terminology,title={Index of Terminology},columns=2,intoc]
\makeindex[name=notation,title={Index of Notation},columns=2,intoc]

\AtBeginEnvironment{appendx}{\crefalias{section}{appsec}}

\makeatletter
	\newcounter{savedsection}% for remembering the last section number
	\preto\appendix{\setcounter{savedsection}{\arabic{section}}}% remembering!
	\newcommand\resumesections{% the \appendix command with some tweakB∇s
	  \setcounter{section}{\arabic{savedsection}}% restore section number
	  \setcounter{section}{\thesavedsection}% reset section counter
	  \gdef\@secpapp{\sectionname}% reset section name
	  \gdef\thesection{\@arabic\c@section}% make section numbers arabic
	  \crefalias{appendix}{section}%
	  \crefname{appendix}{Chapter}{Chapters}%
	  \Crefname{appendix}{Chapter}{Chapters}%
	}
\makeatother

\defbibheading{references}[\refname]{%
	\section*{#1}%
	\addcontentsline{toc}{part}{References} %
	\markboth{#1}{#1}
	}

\title{Differential Cohomology \\ \large Categories, Characteristic Classes, and Connections}

\date{\today}

\author{Edited by Araminta Amabel, Arun Debray, and Peter Haine}

\begin{document}

\bookmark[page=1,level=1]{Cover} % Adds a bookmark in the pdf for the table of contents
%\usetikzlibrary{fadings}
%\tikzfading[name=fade out, inner color=transparent!0,
%  outer color=transparent!50]

\DeclareFixedFont{\titlefont}{T1}{ppl}{b}{}{0.6in}
\DeclareFixedFont{\subtitlefont}{T1}{ppl}{b}{}{0.35in}
\DeclareFixedFont{\subsubtitlefont}{T1}{ppl}{b}{}{0.2in}
\newgeometry{left=0.6in, right=0.6in,top=0in, bottom=0in}
\afterpage{\restoregeometry}
%\definecolor{mytan}{HTML}{FFEEDD}
\pagecolor{black}\afterpage{\nopagecolor}

\thispagestyle{empty}
\begin{center}
\begin{tikzpicture}
\node[%scope fading=south,
inner sep=0pt, outer sep=0pt]{
 \makebox[\textwidth]{\includegraphics[width=\paperwidth]{hexagons2}}
};
\end{tikzpicture}
%\vfill
\vspace{1cm}

\textcolor{white}{
%\begin{flushleft}
\titlefont Differential Cohomology \\[0.3cm]
\subtitlefont Categories, Characteristic\\Classes, and Connections\\[0.75cm]
\subsubtitlefont Edited by Araminta Amabel, Arun Debray, and Peter Haine
%\end{flushleft}
}
\end{center}

\newgeometry{margin=1.5in}

\maketitle

\bookmark[page=2,level=1]{Title Page} % Adds a bookmark in the pdf for the table of contents

\begin{abstract}
	We give an overview of differential cohomology from a modern, homotopy-theoretic perspective in terms of sheaves on manifolds. 
	Although modern techniques are used, we base our discussion in the classical precursors to this modern approach, such as Chern--Weil theory and differential characters, and include the necessary background to increase accessibility. 
	Special treatment is given to differential characteristic classes, including a differential lift of the first Pontryagin class. 
	Multiple applications, including to configuration spaces, invertible field theories, and conformal immersions, are also discussed. 
	This book is based on talks given at MIT's Juvitop seminar run jointly with UT Austin in the Fall of 2019.  
\end{abstract}

\newpage

\bookmark[page=3,level=1]{Contents} % Adds a bookmark in the pdf for the table of contents
\setcounter{tocdepth}{2}
\tableofcontents

%!TEX root = ../diffcoh.tex

\section{Preface}
% stuff Arun added, very much a draft for now
Differential cohomology begins with the observation that many naturally occurring differential forms have
integrality properties. One example is the curvature $\Omega$ of a connection on a complex vector bundle over a
closed manifold $M$; if $N\subset M$ is a closed, oriented, two-dimensional submanifold, then $\int_N \Omega$ is
an integer multiple of $2\pi$. Analogous statements are true, though with different normalization constants, for
other Chern--Weil forms of a vector bundle with connection. The first explanation given is typically that the
cohomology classes represented by these forms are in the image of the map $\H^*(-;\ZZ)\to \H^*(-;
\RR)$, but in a way this fails to capture the entire picture: that the de Rham class of the Chern--Weil form has a
canonical lift to $\H^*(-;\ZZ)$. For example, $(1/2\pi) \Omega$ lifts to the first Chern class of a complex
vector bundle. Differential cohomology is built to house this kind of data: a closed differential form, an
integer-valued cohomology class, and an identification of their images in de Rham cohomology.

A similar situation can happen in quantum physics: abelian
gauge fields give rise to differential forms such as field strengths and currents, and quantization imposes strong
integrality properties on these objects. For example, in the classical theory of electromagnetism, the electric
field $E$ is a $1$-form, and the magnetic field $B$ is a $2$-form. Maxwell's equations on a closed $4$-manifold $M$
imply that the field strength\index[terminology]{field strength} $F = B - \d t\wedge E$ is a closed $2$-form. But
in the quantum theory, the possible values of electric and magnetic fluxes and charges are discretized; there is a
minimum magnetic charge $q_B$, and the integral of $F$ on a closed, oriented surface must be an integer multiple of
$2\pi q_B$. Again we have closed forms with integrality conditions, and so the field strength $B$ refines to a
cocycle representative of a differential cohomology class $\hat B\in\Hhat^2(M; q_B\ZZ)$.

Another perspective on differential cohomology is that it does for geometric objects what ordinary cohomology does
for their topological analogues. Vector bundles and principal bundles have characteristic classes in cohomology;
vector bundles with connection and principal bundles with connection have characteristic classes in differential
cohomology.
%An orientation induces a pushforward map in cohomology, and an orientation and a Riemannian metric
%induce a pushforward map in differential cohomology.
Analogously, topological \Ktheory is built out of vector bundles, and differential
\Ktheory is built out of vector bundles with connection. 

The goal of this book is to provide an introduction to differential cohomology, including both foundational aspects
of generalized differential cohomology theories and applications. We follow Bunke--Nikolaus--Völkl, defining
differential (generalized) cohomology theories as sheaves of spectra on the site of smooth manifolds. We go over the
basics of the theory, including defining the cup product and integration maps. We spend time with characteristic
classes: as hinted above, Chern--Weil forms refine to characteristic classes in differential cohomology, but there
are additional classes which have no topological counterparts. We also go over several applications of differential
cohomology. Often, these are geometric analogues of a well-known application of cohomology to topological
questions. For example, characteristic classes obstruct smooth embeddings of manifolds into $\RR^n$, and
differential characteristic classes can obstruct conformal embeddings into $\RR^n$. Some of these applications are
angled towards physics; for example, we revisit the idea above that differential cohomology has something to say
about quantization.

This book began as lectures given in a graduate student seminar joint between MIT and UT Austin in fall 2019,
initiated by Dan Freed and Mike Hopkins. Most chapters are notes from talks given by various speakers at the
seminar and a few chapters were written afterwards.

%-------------------------------------------------------------------%
%-------------------------------------------------------------------%
%  Assumed background                                               %
%-------------------------------------------------------------------%
%-------------------------------------------------------------------%

\subsection{Assumed background}

We hope that these notes are accessible to readers with a wide range of background knowledge.
The talks included here were part of a topology seminar,
and are therefore biased toward the homotopy theoretic perspective. 
This is evidenced by the fact that we review the definition of a connection
and not that of an \category.
However, knowledge of \categories is not a prerequisite for making use of these notes.
Comfort with sheaves, spectra, and simplicial sets will make reading easier.
The reader will also benefit, both in motivation and understanding,
from a familiarity with basic differential geometry;
this includes connections, curvature, and de Rham cohomology.
\Cref{part:applications} of these notes includes talks on several different applications of differential cohomology.
Enjoyment of these sections should not require any background other than interest in the section title.

%-------------------------------------------------------------------%
%-------------------------------------------------------------------%
%  Linear overview                                                  %
%-------------------------------------------------------------------%
%-------------------------------------------------------------------%

\subsection{Linear overview}

We give a brief overview of the three parts of these notes. A more detailed introduction is given at the beginning of each part.

%-------------------------------------------------------------------%
%  Part I: Basics of the theory                                     %
%-------------------------------------------------------------------%

\subsubsection{\texorpdfstring{\Cref{part:basics}}{Part \ref*{part:basics}}: Basics of the theory}

The purpose of this part is to introduce the basics of and develop the general theory behind differential cohomology.
In \Cref{Introduction}, we start with some motivation to the approach we take to differential cohomology coming from work of Cheeger--Simons \cite{MR827262} and Simons--Sullivan \cite{MR2365651} on differential characters and ordinary differential cohomology.
The perspective we take on differential cohomology theories is as sheaves of spectra on the category $ \Mfld $ of manifolds; since we also want to consider sheaves that come from chain complexes, we'll work in the framework of sheaves with values in a general \category.
While this might sound somewhat daunting, there are many familiar examples: 
\begin{enumerate}
	\item The functor sending a manifold $ M $ to the complex $ \Omegabullet(M) $ of de Rham cochains on $ M $.

	\item The functor sending a manifold $ M $ to the complex $ \Csing^\bullet(M) $ of singular cochains on $ M $.

	\item Given a Lie group $ G $, the functor sending a manifold $ M $ to the groupoid $ \BunG(M) $ (or $ \BunGnabla(M) $) of principal $ G $-bundles on $ M $ (with connection). 
\end{enumerate}
The new example of differential cohomology is essentially built from these ones in a nontrivial way.

In \Cref{sec:basicsetup}, we introduce the basics of sheaves on the category of manifolds, how to manipulate sheaves on $ \Mfld $, and any the category of sheaves (of sets) on $ \Mfld $ contains the standard category of infinite-dimensional manifolds (\textit{Fréchet manifolds}) as a full subcategory.
One important class of sheaves on $ \Mfld $ are those that invert all homotopy equivalences of manifolds.
\Cref{sec:hisheaves} is dedicated to explaining why all sheaves with this property have a very simple and concrete description.
In \Cref{sec:localization}, we explain how to resolve a sheaf by one that inverts all homotopy equivalences of manifolds.
This provides a way of decomposing a sheaf of spectra on $ \Mfld $ into one that inverts all homotopy equivalences
and another that ``comes from geometry''.
\Cref{sec:stable} explains this decomposition as well as how this gives rise to the Simons--Sullivan ``differential
cohomology hexagon'' \cite[\S1]{MR2365651}) relating ordinary cohomology, differential forms, and differential cohomology.

The remainder of this part is dedicated to important examples of differential cohomology theories and refining important constructions with ordinary cohomology.
\Cref{sec:examples} explains Cheeger--Simons differential characters, differential $ \K $-theory, and examples coming from $ G $-bundles in the framework of sheaves on $ \Mfld $.
\Cref{sec:Delignecup} refines the cup product to differential cohomology and explains how to calculate it in many examples.
\Cref{FiberIntegration} refines fiber integration to differential cohomology.
\Cref{sec:digressiononTransferConjecture} finishes the main text of this part with a digression proving Quillen's Transfer Conjecture.
Though not directly related to differential cohomology, this result states that connective spectra can be realized as homotopy-invariant sheaves on the category of correspondences of manifolds where the backwards maps are finite covering maps (i.e., connective spectra have natural \textit{transfers} along finite covering maps).
Our exposition follows work of Bachmann--Hoyois \cite[Appendix C]{MotivicNorms:BachmannHoyois}.

\Cref{part:basics} also has an appendix (\Cref{app:technicaldeatails}).
In this appendix, we prove a few technical category theory results that we need to get the foundations of sheaves on $ \Mfld $ on a solid framework in \Cref{sec:basicsetup,sec:hisheaves}.

%-------------------------------------------------------------------%
%  Part II: Characteristic classes                                  %
%-------------------------------------------------------------------%

\subsubsection{\texorpdfstring{\Cref{part:charclasses}}{Part \ref*{part:charclasses}}: Characteristic classes}

Just as one ordinary cohomology is a natural home for characteristic classes, differential cohomology offers its own invariants of bundles. These invariants, known as ``differential characteristic classes,'' are refinements of the classical characteristic classes in cohomology. More explicitly, we will investigate lifts of well-known characteristic classes, such as Chern classes, under the map from differential cohomology to ordinary cohomology.

This part begins be reviewing a few classical techniques and results that will be useful in studying differential characteristic classes, see \Cref{ChernWeilTheory} and \Cref{EquivariantdeRhamCohomology}.

Differential characteristic classes where first studied by Cheeger--Simons \cite{Cheeger-Simons}. %sentence about what these are
We discuss differential characters in \Cref{DifferentialCharacteristicClasses}. 
Building on work of Bott \cite{BottPaper}, Freed and Hopkins \cite{FreedHopkins} classified all differential characteristic classes for bundles equipped with a flat connection. 
This refines the classical Chern--Weil story, which we review in \Cref{ChernWeilTheory}. 
The contents of \cite{FreedHopkins} are covered in \Cref{WorkofFreedHopkins}.  
A closer look at the methods used in \cite{BottPaper} reveal that one can remove the connection data with some alterations. 
In \Cref{BottsMethod}, we delve into Bott's paper and the theorems it relies upon. In particular, we discuss van Est's theorem relating continuous cohomology to Lie algebra cohomology. 
Using the results of \cite{BottPaper}, Hopkins, in \Cref{LiftsofChernClasses}, discusses how to lift ordinary Chern classes to a form of differential cohomology, without the presence of a connection. 
The existence of a differential version of the Cartan formula is also considered. 

This part of the notes concludes with an interesting application of differential lifts of Chern classes to a possible construction of the Virasoro group. The Virasoro group is a certain central extension of orientation preserving diffeomorphisms $\Diffplus(\Circ)$ of $\Circ$. As Hopkins outlines in \Cref{LiftsofChernClasses}, one can obtain central extensions of $\Diffplus(\Circ)$ from a certain differential cohomology group. The details of this construction, as well as a review of the Virasoro algebra and group, appear in \Cref{VirasoroAlgebra}.

%-------------------------------------------------------------------%
%  Part III: Applications                                           %
%-------------------------------------------------------------------%

\subsubsection{\texorpdfstring{\Cref{part:applications}}{Part \ref*{part:applications}}: Applications}
In this part we discuss some uses of differential cohomology in topology, geometry, and physics. Some, but not all,
of these applications are part of the idea that what ordinary cohomology can do for topological questions,
differential cohomology can do for geometric ones, and many of these applications are related to various aspects of
quantum field theory.

One of the key links between differential cohomology and geometry is through Chern--Simons invariants, invariants
of connections which can be defined either in terms of integration of differential characteristic classes or
directly using geometric information. Because of this, several applications of differential cohomology to geometry
or physics pass through Chern--Simons theory. We introduce and apply Chern--Simons invariants in
\Cref{config_spaces} and also use them in \Cref{conformal_immersions}.

Our first two applications of differential cohomology are in geometry and topology. In \Cref{config_spaces}, we
discuss work of Evans-Lee--Saveliev \cite{deletedsquare}, who use Chern--Simons invariants to study the homotopy
types of two-point configuration spaces of lens spaces.
%Longoni--Salvatore \cite{LS05} had shown that the homotopy
%type of the two-point configuration space is in general a stronger invariant than the homotopy type of the
%underlying manifold, and using Chern--Simons invariants and a few other techniques, Evans-Lee--Saveliev are able to
%provide many more examples of this phenomenon.
Then in \Cref{conformal_immersions}, we use differential Pontryagin classes and Chern--Simons forms to obstruct
conformal immersions of conformal manifolds into Euclidean space, following Chern--Simons \cite{cs}; along the way
we spend some time getting to know the geometry of Chern--Weil and Chern--Simons forms.

The next two applications are to physics. \Cref{field_theory} applies differential cohomology to the quantization
of abelian gauge fields, using electromagnetism as an example. In classical physics, the field strength of an
abelian gauge field is a closed differential form; quantization lifts from closed forms to
cocycles for a differential cohomology group. The other physics
application we discuss, in \Cref{invertible_field_theories}, is quite different: a conjecture of
Freed--Hopkins \cite{FH21} using differential generalized cohomology to classify invertible, non-topological field
theories. This is a geometric conjecture modeled on a topological theorem of Freed--Hopkins (\textit{ibid.})
classifying invertible topological field theories using Madsen--Tillmann spectra.
We discuss this conjecture and several examples, including
classical Chern--Simons theory.

Our final two chapters are about the representation theory of loop groups. Loop groups are infinite-dimensional Lie
groups whose representation theory is strikingly similar to that of compact Lie groups, so long as one works with
what are called positive energy representations. In \Cref{loop_groups}, we survey this theory, defining and
motivating positive energy representations and sketching a proof of a theorem of Pressley--Segal \cite{loop}, which
says that positive energy representations admit projective intertwining actions of $\Diffplus(\Circ)$.
%We also
%discuss connections between the representation theory of loop groups and differential cohomology.
In
\Cref{segal_sugawara}, we study the Pressley--Segal theorem at the Lie algebra level, where this intertwining
projective action can be made more explicit. Since projective representations are equivalent to representations of
a central extension, the Virasoro algebra makes an appearance here.
%The math in
%\Cref{loop_groups,segal_sugawara} is related to aspects of two-dimensional conformal field theory, and we discuss
%some of the connections.

%-------------------------------------------------------------------%
%-------------------------------------------------------------------%
%  What's not included                                              %
%-------------------------------------------------------------------%
%-------------------------------------------------------------------%

\subsection{What's not included}

One %the?
original approach to differential cohomology is presented by Hopkins and Singer in \cite{HopkinsSinger}. 
While we look to this reference for motivation and intuition, 
we do not take this as our definition of a differential cohomology theory. 
Instead, we work with the more modern approach using sheaves on manifolds. 
We also make use of \cite{HopkinsSinger} for constructions of the cup product and fiber integration in differential cohomology, see \Cref{sec:Delignecup,FiberIntegration}. 

Several examples of differential cohomology theories, such as differential K-theory, are discussed in \Cref{sec:examples}; 
however, there are many more examples that we do not mention. %add reference to other examples in literature
Moreover, for most of these notes, we focus our attention on the specific example of the differential version of ordinary cohomology. 
This leaves several interesting areas of study, such as differential K-theory characteristic classes, untouched.

We do not present Schreiber's elegant and very general theory of \textit{differential cohomology in a cohesive \topos} \cite{Urs}.
Schreiber's work requires background that we do not assume; we decided to stick with the setting of sheaves on the category of manifolds to make the material accessible to the graduate students attending the seminar.

There are also many applications of differential cohomology to physics which we do not discuss in detail here. See \Cref{applications_part} for a discussion of related work. 

%-------------------------------------------------------------------%
%-------------------------------------------------------------------%
%  Cover image                                                      %
%-------------------------------------------------------------------%
%-------------------------------------------------------------------%

\subsection{Cover image}

One of the theses of this book is that differential cohomology has applications to physics. It therefore seems apt
to choose a cover image of another example of hexagons in the real world. Our cover image is a picture of Giant's
Causeway, a part of the coastline in Northern Ireland consisting of tens of thousands of tessellating hexagonal
basalt columns. This image is by Giuseppe Milo and can be found at
\href{https://www.flickr.com/photos/giuseppemilo/46587488041/in/photostream/}{\nolinkurl{flickr.com/photos/giuseppemilo/46587488041/in/photostream/}}; we cropped it slightly. It is
licensed under the \href{https://creativecommons.org/licenses/by/2.0/}{CC BY 2.0} license.

%-------------------------------------------------------------------%
%  Acknowledgements                                                 %
%-------------------------------------------------------------------%

\subsection{Acknowledgements}

These notes are a compilation of talks given in the Juvitop seminar at MIT in the fall of 2019. 
The seminar was run jointly with UT Austin and we thank the participants of both cities for their comments and discussion. 
%In particular, thanks to specific participants, 
%I don't remember who Arun said was particularly helpful, maybe it was Charlie?
%maybe something about people who talked in the slack
%maybe thanks to Dylan for discussions about organization: We thank Dylan Wilson for helping us figuring out how to structure the seminar as well as for answering a number of mathematical questions.
We would also like to thank the speakers, 
Dexter Chua, Sanath Devalapurkar, Dan Freed, Mike Hopkins, Greg Parker, Charlie Reid, and Adela Zhang,
both for volunteering to speak, as well as writing up notes to appear here. 

The seminar, and these notes, could not have existed without the help of Dan Freed and Mike Hopkins. 
We thank both Freed and Hopkins for their mathematical help, their organizational help, and for giving talks in the seminar. 
We also would like to express our appreciation for Freed and Hopkins %and teleman?
generosity in sharing new mathematical ideas and older insights along the way.

Extra thanks are due to Hopkins for buying us fancy video equipment to help with our half-virtual seminar.

We thank Kiran Luecke and Zhouhang Mao for pointing out a number of typos and errors in the first draft of this text.

PH gratefully acknowledges support from the MIT Dean of Science Fellowship, UC President's Postdoctoral Fellowship, and NSF Mathematical Sciences Postdoctoral Research Fellowship under Grant \#DMS-2102957. 
PH and AA acknowledge support from the National Science Foundation Graduate Research Fellowship under Grant \#112237.
This text is partially based upon work supported by the National Science Foundation under Grant \#DMS-1440140, while AD and PH were in residence at the Mathematical Sciences Research Institute in Berkeley, California, during the Spring 2020 semester.

\newpage

%-------------------------------------------------------------------%
%  Basics of the theory                                             %
%-------------------------------------------------------------------%

\part{Basics of the theory}

\numberwithin{equation}{part} % Changes numbering to I.# for the part intro
%!TEX root = ../diffcoh.tex

%-------------------------------------------------------------------%
%-------------------------------------------------------------------%
%  Part I Introduction                                              %
%-------------------------------------------------------------------%
%-------------------------------------------------------------------%

\label{part:basics}

The goal of this first part of the text is to introduce and study \textit{differential cohomology theories}. 
The term ``differential cohomology'' was first coined by Hopkins and Singer in \cite{HopkinsSinger}. 
In \cref{Introduction}, we introduce the ideas of differential cohomology theories following Cheeger--Simons \cite{MR827262} and Simons--Sullivan \cite{MR2365651}.
The basic point is that given a manifold $ M $, we can consider both the ``homotopy-theoretic'' complex of singular
cochains on $ M $, and the ``geometric'' complex of differential forms on $ M $.  These are related by the de Rham
isomorphism, and we would like to combine them together into a ``cohomology theory'' that captures both the
features of $ M $ as a homotopy type as well as the geometry of $ M $.  The thing to notice is that both the complex
of singular cochains and differential forms are \textit{sheaves} (in the homotopy-theoretic sense) on the category
of all manifolds.  So this category of sheaves on manifolds is the setting in which both these homotopy-theoretic
and geometric objects live.

Thus the perspective that we take in this text is that differential cohomology theories are \textit{sheaves of spectra on the category $ \Mfld $ of manifolds}.
It will also be useful to consider sheaves of spaces on $ \Mfld $ or sheaves with values in the derived \category of a ring; \Cref{sec:basicsetup} starts with introducing sheaves on the category of manifolds with values in any \category.
While the phrase ``sheaf on $ \Mfld $'' may sound somewhat daunting, it is surprisingly concrete: a sheaf $ F $ on $ \Mfld $ consists of a functor $ \fromto{\Mfld}{C} $ such that for each manifold $ M $, the restriction of $ F $ to open subsets of $ M $ defines a sheaf on $ M $. 

Let $ C $ be a presentable \category (e.g, spaces, spectra, or the derived \category of a ring).
One of the basic features of the category $ \Sh(\Mfld;C) $ of $ C $-valued sheaves on $ \Mfld $ is that the full subcategory $ \Shhi(\Mfld;C) $ spanned by those sheaves that invert homotopy equivalences is already familiar: 

\begin{theorem}[(\Cref{prop:Dugger})]\label{introthm:Dugger}
	Evaluation on the point defines an equivalence
	\begin{align*}
		\Gammalowerstar \colon \Shhi(\Mfld;C) &\equivalence C \\
		F &\mapsto F(*) \period
	\end{align*}
	Moreover, the inverse equivalence is given by the \emph{constant sheaf functor} $ \Gammaupperstar \colon \fromto{C}{\Sh(\Mfld;C)} $.
	That is, $ \Shhi(\Mfld;C) $ coincides with the full subcategory of $ \Sh(\Mfld;C) $ spanned by the constant sheaves. 
\end{theorem}

We call objects of $ \Shhi(\Mfld;C) $ \textit{\RRinvariant} sheaves.
\Cref{sec:hisheaves} is dedicated to proving \Cref{introthm:Dugger}. 
In \Cref{sec:hisheaves} we also give an explicit formula for the constant sheaf functor $ \fromto{C}{\Sh(\Mfld;C)} $:

\begin{proposition}[(\Cref{lem:constantishi})]\label{introprop:constantishi}
	The constant sheaf functor 
	\begin{equation*}
		\Gammaupperstar \colon C \equivalence \Shhi(\Mfld;C) \subset \Sh(\Mfld;C)
	\end{equation*}
	is given by the assignment
	\begin{equation*}
		X \mapsto \big[M \mapsto X^{\Piinf(M)}\big] \period 
	\end{equation*}
	Here, $ X^{\Piinf(M)} $ denotes the \emph{cotensor} of the object $ X \in C $ by the underlying homotopy type $ \Piinf(M) $ of the manifold $ M $ (see \Cref{rec:cotensoring}).
\end{proposition}

The cotensor in \Cref{introprop:constantishi} might look a bit mystifying, but it is actually a familiar object in
the specific values of $ C $ that we're most interested in:
\begin{enumerate}
	\item Let $ C = \Spc $ be the \category of spaces.
	In this case, the constant sheaf functor is given by
	\begin{equation*}
		X \mapsto \big[M \mapsto \Map_{\Spc}(\Piinf(M),X) \big] \period 
	\end{equation*}

	\item Let $ C = \Sp $ be the \category of spectra.
	In this case, the constant sheaf functor is given by
	\begin{equation*}
		E \mapsto \big[M \mapsto \Hom_{\Sp}(\Sigma_{+}^{\infty} \Piinf(M),E) \big] \comma 
	\end{equation*}
	where $ \Hom_{\Sp} $ is the mapping \textit{spectrum}.

	\item Let $ R $ be a ring and let $ C = \D(R) $ be the \textit{derived \category} of $ R $ obtained from the
	category of chain complexes of $ R $-modules by formally inverting the quasi-isomorphisms
	\cite[\HAthm{Definition}{1.3.5.8}, \HAthm{Proposition}{1.3.5.15}, \& \HAthm{Remark}{7.1.1.16}]{HA}.
	In this case, the constant sheaf functor is given by
	\begin{equation*}
		A_{*} \mapsto \big[M \mapsto \RHom_{R}(\Cup_*(M;R),A_{*}) \big] \period 
	\end{equation*}
	Here $ \Cup_*(M;R) $ is the complex of singular chains on $ M $, and $ \RHom_{R} $ is the derived $ \Hom $
	functor of chain complexes of $ R $-modules.
\end{enumerate}

As a consequence of \Cref{introprop:constantishi} (and some simple observations), we show that there is a chain of four adjoints
\begin{equation}\label{diag:intromanyadjoints}
	\begin{tikzcd}[sep=4em]
		\Sh(\Mfld;C) \arrow[r, shift left=2ex, "\Gammalowersharp"] \arrow[r, shift right=0.75ex, "\Gammalowerstar"{description, xshift=0.5em}] & C \period \arrow[l, shift left=2ex, "\Gammauppersharp", hooked'] \arrow[l, shift right=0.75ex, "\Gammaupperstar"{description, xshift=-0.5em}, hooked']
	\end{tikzcd}
\end{equation}
Here functors lie above their right adjoints.
The extreme right adjoint $ \Gammauppersharp $ has an explicit formula (see \Cref{lem:Gammauppershriekpsh}), but is not particularly useful.
On the other hand, under the identification
\begin{equation*}
	\Gammaupperstar \colon \equivto{C}{\Shhi(\Mfld;C)}
\end{equation*}
the extreme left adjoint $ \Gammalowersharp $ corresponds to the left adjoint to the inclusion
\begin{equation*}
	\incto{\Shhi(\Mfld;C)}{\Sh(\Mfld;C)}
\end{equation*}
We initially construct the left adjoint $ \Gammalowersharp $ abstractly via the Adjoint Functor Theorem, but since it plays a very important role throughout this text, it is useful to have an explicit formula for $ \Gammalowersharp $.
\Cref{sec:localization} is dedicated to showing that $ \Gammalowersharp(F) $ is computed by a simple geometric realization.
Write $ \Deltaalg^n $ for the hyperplane
\begin{equation*}
	\Deltaalg^n \colonequals \setbar{(t_0,\ldots,t_{n}) \in \RR^{n+1}}{t_0 + \cdots + t_{n} = 1} \subset \RR^{n+1} \period
\end{equation*}

\begin{theorem}[(\Cref{cor:formula_for_Gammalowersharp})]\label{introthm:Gammalowersharp}
	The left adjoint $ \Gammalowersharp \colon \fromto{\Sh(\Mfld;C)}{C} $ is given by the formula
	\begin{equation*}
		\Gammalowersharp(F) \equivalent \real{F(\Deltaalgdot)} \period
	\end{equation*}
\end{theorem}

\noindent \Cref{sec:localization} also explores some important consequences of \Cref{introthm:Gammalowersharp}.
For example, we give differential refinements of classifying spaces for $ G $-bundles (see \cref{subsubsec:consequencesofMSV}).

Some of the proofs in \Cref{sec:basicsetup,sec:hisheaves,sec:localization} rely on technical results about \topoi or presentable \categories.
To avoid distracting the reader from the main point of the text, we have relegated many of these details to \Cref{app:technicaldeatails}.

\Cref{sec:stable} specializes to sheaves with values in a presentable stable \category like spectra or the derived \category of a ring.
Using the many adjoint functors \eqref{diag:intromanyadjoints} constructed in \Cref{sec:hisheaves}, we prove the existence of a \textit{fracture square} that shows that every sheaf on $ \Mfld $ can be glued together from an \RRinvariant sheaf and a sheaf with vanishing global sections (\cref{subsec:fracture}).
Using this fracture square, we provide a version of the Simons--Sullivan differential cohomology diagram (\Cref{thm:SimonsSullivanunique}) for any differential cohomology theory (\cref{subsec:diffcohdiagram}).
We also begin the study of differential refinements of spectra (\cref{subsec:refinements}).

With the basic foundations out of the way, \Cref{sec:examples} is dedicated to examples of differential cohomology theories.
These include \textit{ordinary differential cohomology} after Cheeger--Simons and Delgine (\cref{subsec:CheegerSimonsDiffchar}), and \textit{differential \Ktheory} after Hopkins--Singer (\cref{subsec:diffKtheory}).

In \Cref{sec:Delignecup} we further analyze ordinary differential cohomology by giving it a product structure called the \textit{Deligne cup product}.

\begin{definition}
	Let $ k \geq 0 $ be an integer.
	The \textit{Deligne complex} $ \ZZ(k) $ is the pullback
	\begin{equation*}
		\begin{tikzcd}
			\ZZ(k)\arrow[r]\arrow[d] \arrow[dr, phantom, "\lrcorner"{description, very near start}] & \ZZ\arrow[d] \\
			\Sigma^k\Omegacl^k\arrow[r] & \RR 
		\end{tikzcd}
	\end{equation*}
	in the \category $ \Sh(\Mfld;\D(\ZZ)) $ of sheaves on $ \Mfld $ with values in the derived \category of $ \ZZ $.
\end{definition}

The Deligne complex $ \ZZ(k) $ used to define ordinary differential cohomology.
The Deligne cup product
\begin{equation*}
	\fromto{\ZZ(m) \tensor_{\ZZ} \ZZ(n)}{\ZZ(m+n)}
\end{equation*}
is defined by combining the cup product on integral cohomology with the wedge product on differential forms.
We conclude \Cref{sec:Delignecup} with an analysis of the Deligne cup product in detail in the lowest dimensions (\cref{subsec:Delignecupexamples}).

\Cref{FiberIntegration} refines fiber integration to differential cohomology.
After reviewing fiber integration for ordinary cohomology, we introduce differential versions of Thom classes and orientations (\cref{subsec:differentialintegration}).
We then use these notions to define \textit{differential fiber integration} and explain how this works for $ \Circ $-bundles.

\Cref{sec:digressiononTransferConjecture} is a digression explaining Bachmann and Hoyois' proof of Quillen's \textit{Transfer Conjecture} \cite[Appendix C]{MotivicNorms:BachmannHoyois}.
This identifies the \category of $ \Einf $-spaces with \RRinvariant sheaves on a $ 2 $-category of manifolds with morphisms \textit{correspondences}
\begin{equation*}
	\begin{tikzcd}[row sep=1.5em, column sep=0.75em]
		& N \arrow[dr] \arrow[dl] & \\
		M_0 & & M_1 \comma
	\end{tikzcd}
\end{equation*}
where the ``backwards'' maps are finite covering maps.
Restricting to grouplike objects on both sides gives a description of the \category of connective spectra in terms of sheaves on this $ 2 $-category of manifolds and correspondences.
\Cref{sec:digressiononTransferConjecture} is not used later in the text; we have included it because of its connection to \Cref{sec:basicsetup,sec:hisheaves,sec:localization}, but the uninterested reader can safely skip it.

\numberwithin{equation}{subsection} % Changes numbering back to normal
\newpage
%!TEX root = ../diffcoh.tex

%-------------------------------------------------------------------%
%-------------------------------------------------------------------%
%  Introduction                                                     %
%-------------------------------------------------------------------%
%-------------------------------------------------------------------%

\section{Introduction}\label{Introduction}
\textit{by Peter Haine}

The purpose of this chapter is to give some motivation for the perspective we take on differential cohomology.
We do this by giving an overview of the work of Cheeger--Simons \cite{MR827262}, Deligne \cites[\S2.2]{MR498551}[\S12.3]{MR2451566}, and Simons--Sullivan \cite{MR2365651} on differential cohomology.

%-------------------------------------------------------------------%
%-------------------------------------------------------------------%
%  Motivation for differential cohomology                           %
%-------------------------------------------------------------------%
%-------------------------------------------------------------------%

\subsection{Motivation for differential cohomology}

\begin{observation}[{(Simons--Sullivan \cite[\S1]{MR2365651})}]
	Let $ M $ be a manifold.
	Then we have exact sequences
	\begin{equation}\label{diag:diffcohpre}
		\begin{tikzcd}[column sep={10ex,between origins}, row sep={8ex,between origins}]
			& \H^{k-1}(M;\RR/\ZZ) \arrow[rr, "-\Bock"] & & \H^k(M;\ZZ) \arrow[dr] & \\
			\HdR^{k-1}(M) \arrow[ur] \arrow[dr] & & & & \HdR^k(M) \\
			& \Omega^{k-1}(M)/\im(\d) \arrow[rr, "\d"'] & & \Omegacl^k(M) \arrow[ur] & \phantom{\HdR^k(M)} \comma
		\end{tikzcd}
	\end{equation}
	where the top sequence is the Bockstein sequence associated to the short exact sequence
	\begin{equation*}
		\begin{tikzcd}[sep=1.5em]
			0 \arrow[r] & \ZZ \arrow[r, hook] & \RR \arrow[r, ->>] & \RR/\ZZ \arrow[r] & 0 \comma 
		\end{tikzcd}
	\end{equation*}
	and we are identifying singular and de Rham cohomology via the de Rham isomorphism
	\begin{equation*}
		 \HdR\upperstar(M) \isomorphic \H\upperstar(M;\RR) \period
	\end{equation*}

	The top sequence is ``purely homotopy-theoretic'' in nature, while the bottom sequence is ``purely geometric'' in nature (i.e., the functor $ \Omegacl^k $ is not homotopy-invariant).
\end{observation}

\begin{question}\label{qst:origin}
	Can we fill \eqref{diag:diffcoh} in with an invariant \textcolor{love}{$ \Hhat^k(M;\ZZ) $ in maroon} 
	\begin{equation}\label{diag:diffcoh}
		\begin{tikzcd}[column sep={10ex,between origins}, row sep={8ex,between origins}]
			& \H^{k-1}(M;\RR/\ZZ) \arrow[rr, "-\Bock"] \arrow[dr, love] & & \H^k(M;\ZZ) \arrow[dr] & \\
			\HdR^{k-1}(M) \arrow[ur] \arrow[dr] & & \textcolor{love}{\Hhat^k(M;\ZZ)} \arrow[ur, love] \arrow[dr, love] & & \HdR^k(M) \\
			& \Omega^{k-1}(M)/\im(\d) \arrow[rr, "\d"'] \arrow[ur, love] & & \Omegacl^k(M) \arrow[ur] & \phantom{\HdR^k(M)} \comma
		\end{tikzcd}
	\end{equation}
	that better blends homotopy theory and geometry, and makes the diagonals exact?
\end{question}

Now let us attempt to provide a satisfactory answer to \Cref{qst:origin} when $ k = 1 $.

\begin{attempt}[(for $ k = 1$)]
	Let $ M $ be a manifold.
	Consider the abelian group $ \Cinf(M,\RR/\ZZ) $ of smooth functions to the circle (with the group structure defined pointwise).
	We should really think of $ \Cinf(M,\RR/\ZZ) $ as an infinite-dimensional abelian Lie group.
	Recall that the inclusion
	\begin{equation*}
		 \Cinf(M,\RR/\ZZ) \subset \Map(M,\RR/\ZZ)
	\end{equation*}
	from the space of smooth maps to the space of continuous maps is a homotopy equivalence.
	Since the circle is $ 1 $-truncated,\footnote{I.e., only has nontrivial homotopy groups in degrees $ \leq 1$.} 
	this implies that 
	$ \Cinf(M,\RR/\ZZ) $ is also $ 1 $-truncated.

	Since $ \RR/\ZZ $ is a $ \K(\ZZ,1) $, we see that
	\begin{equation*}
		\uppi_0 \Cinf(M,\RR/\ZZ) \isomorphic \H^1(M;\ZZ) \period
	\end{equation*}
	In particular, we have a surjection $ \uppi_0 \colon \surjto{\Cinf(M,\RR/\ZZ)}{\H^1(M;\ZZ)} $.
	Also notice that
	\begin{align*}
		\uppi_1\Cinf(M,\RR/\ZZ) &\isomorphic \uppi_0\Map_{*}(\Circ,\Cinf(M,\RR/\ZZ)) \\
		&\isomorphic \uppi_0\Map_{*}(\Circ,\Map(M,\RR/\ZZ)) \\
		&\isomorphic \uppi_0\Map(M,\Map_{*}(\Circ,\RR/\ZZ)) \\
		&\isomorphic \uppi_0\Map(M,\Omega(\RR/\ZZ)) \\
		&\isomorphic \H^0(M;\ZZ) \period
	\end{align*}
\end{attempt}

\begin{construction}
	Let $ \vol $ denote the standard volume form on the circle $ \Circ \isomorphic \RR/\ZZ $.
	Define a \emph{curvature} map $ \curv \colon \fromto{\Cinf(M,\RR/\ZZ)}{\Omegacl^1(M)} $ by
	\begin{equation*}
		\curv(f) \colonequals \fupperstar(\vol) \period
	\end{equation*}
\end{construction}

\begin{nul}
	The kernel of $ \curv $ consists of the locally constant maps $ \fromto{M}{\RR/\ZZ} $, i.e.,
	\begin{equation*}
		\ker(\curv) \isomorphic \H^0(M;\RR/\ZZ) \period
	\end{equation*}
	Note that the curvature map is not surjective:
	\begin{equation*}
		\image(\curv) = \{ \alpha \in \Omegacl^1(M) \, | \, \textstyle\int_{\Circ} \alpha \in \ZZ \text{ for every embedding } \incto{\Circ}{M} \} \period
	\end{equation*}
	That is, the image of $ \curv $ is the group of \textit{closed $ 1 $-forms with integral periods}.
\end{nul}

\begin{definition}[(integral periods)]
	Let $ M $ be a manifold and $ k \geq 0 $ an integer. 
	A closed $ k $-form $ \omega $ on $ M $ \emph{has integral periods} if for every smooth $ k $-cycle $ c $ in $ M $ the integral $ \int_c \omega $ is an integer.
	We write
	\begin{equation*}
		\Omegacl^k(M)_{\ZZ} \subset \Omegacl^k(M)
	\end{equation*}
	for the subgroup of $ k $-forms with integral periods.
\end{definition}

\begin{remark}
	A closed $ k $-form $ \omega $ has integral periods if and only if the class of $ \omega $ lies in the image of the change-of-coefficients map
	\begin{equation*}
		\fromto{\H^k(M;\ZZ)}{\H^k(M;\RR) \isomorphic \HdR^k(M)} \period 
	\end{equation*}
\end{remark}

\begin{nul}
	We also have a map
	\begin{equation*}
		\iota \colon \fromto{\Omega^0(M) = \Cinf(M,\RR)}{\Cinf(M,\RR/\ZZ)}
	\end{equation*}
	given by post-composition with the quotient map $ \surjto{\RR}{\RR/\ZZ} $.
	The map $ \iota $ has kernel the integer-valued smooth functions $ \fromto{M}{\RR} $, i.e., the locally constant functions with integer values. 
	That is, $ \image(\iota) = \Omegacl^0(M)_{\ZZ} $.
\end{nul}

\begin{nul}
	These maps give rise to a commutative diagram with exact diagonals
	\begin{equation*}
		\begin{tikzcd}[column sep={10ex,between origins}, row sep={8ex,between origins}]
			& \H^{0}(M;\RR/\ZZ) \arrow[rr, "-\Bock"] \arrow[dr, hook] & & \H^1(M;\ZZ) \arrow[dr] & \\
			\HdR^{0}(M) \arrow[ur] \arrow[dr, hook] & & \Cinf(M,\RR/\ZZ) \arrow[ur, ->>, "\uppi_0" description] \arrow[dr, "\curv" description] & & \HdR^1(M) \\
			& \Omega^{0}(M) \arrow[rr, "\d"'] \arrow[ur, "\iota" description] & & \Omegacl^1(M) \arrow[ur] & \phantom{\HdR^1(M)} \period
		\end{tikzcd}
	\end{equation*}
	The diagonals become short exact sequences if we replace $ \Omega^0(M) $ by $ \Omega^0(M)/\Omegacl^0(M)_{\ZZ} $ and $ \Omegacl^1(M) $ by $ \Omegacl^1(M)_{\ZZ} $:
	\begin{equation*}
		\begin{tikzcd}[column sep={10ex,between origins}, row sep={8ex,between origins}]
			0 \arrow[dr] & & & & 0 \\
			& \H^{0}(M;\RR/\ZZ) \arrow[rr, "-\Bock"] \arrow[dr, hook] & & \H^1(M;\ZZ) \arrow[dr] \arrow[ur]  & \\
			\HdR^{0}(M) \arrow[ur] \arrow[dr, hook] & & \Cinf(M,\RR/\ZZ) \arrow[ur, ->>, "\uppi_0" description] \arrow[dr, "\curv" description, ->>] & & \HdR^1(M) \\
			& \Omega^0(M)/\Omegacl^0(M)_{\ZZ} \arrow[rr, "\d"'] \arrow[ur, >->, "\iota" description] & & \Omegacl^1(M)_{\ZZ} \arrow[ur] \arrow[dr] & \\
			0 \arrow[ur] & & & & 0 \period
		\end{tikzcd}
	\end{equation*}
\end{nul}

\begin{nul}
	The takeaway is that in \Cref{qst:origin}, we should really replace $ \Omega^{k-1}(M)/\im(d) $ by $ \Omega^{k-1}(M)/\Omegacl^{k-1}(M)_{\ZZ} $ and $ \Omegacl^k(M) $ by $ \Omegacl^0(M)_{\ZZ} $ and ask for the diagonal sequences to be short exact.
\end{nul}

% \begin{attempt}[for $ k = 2 $]
% 	Let $ M $ be a manifold and consider the group $ \Buncon(M) $ of isomorphism classes of principal $ \mathup{U}(1) $-bundles with connection under tensor product.
% \end{attempt}

%-------------------------------------------------------------------%
%-------------------------------------------------------------------%
%  Differential characters                                          %
%-------------------------------------------------------------------%
%-------------------------------------------------------------------%

\subsection{Differential characters}

We now present a unified approach to defining the ``differential cohomology'' groups $ \Hhat^{*}(M;\ZZ) $ due to Cheeger--Simons \cite{MR827262}.
We follow Bär and Becker's exposition on \textit{differential characters} \cite[Part I, \S5]{MR3237728}.

\begin{notation}
	Let $ M $ be a manifold and $ i \geq 0 $ an integer.
	We write $ \Csm_i(M;\ZZ) $ for the abelian group of smooth (integer-valued) chains on $ M $.
	We write $ \Zsm_i(M;\ZZ) \subset \Csm_i(M;\ZZ) $ for the subgroup of smooth cycles.
\end{notation}

\begin{definition}[{(Cheeger--Simons \cite[\S1]{MR827262})}]
	Let $ k \geq 1 $ be an integer and $ M $ a manifold.
	A \emph{degree $ k $ differential character} on $ M $ is a homomorphism $ \chi \colon \fromto{\Zsm_{k-1}(M;\ZZ)}{\RR/\ZZ} $ such that there exists a $ k $-form $ \omega(\chi) \in \Omega^k(M) $ with the property that for every $ c \in \Csm_k(M;\ZZ) $, 
	\begin{equation*}
		\chi(\partial c) = \int_c \omega(\chi) \quad \text{mod } \ZZ \period
	\end{equation*} 
	We write
	\begin{equation*}
		\Hhat^k(M;\ZZ) \subset \Hom_{\ZZ}(\Zsm_{k-1}(M;\ZZ),\RR/\ZZ)
	\end{equation*}
	for the abelian group of degree $ k $ differential characters on $ M $.
	
	It follows that $ \omega(\chi) $ is unique and closed.
	Moreover, $ \omega(\chi) $ has integral periods.
	The form $ \omega(\chi) $ is called the \emph{curvature} of $ \chi $, and we have a curvature map
	\begin{align*}
		\curv \colon \Hhat^k(M;\ZZ) &\to \Omega^k(M) \\
		\chi &\mapsto \omega(\chi)
	\end{align*} 
	with image $ \Omegacl^k(M)_{\ZZ} $ those closed $ k $-forms with integral periods.
\end{definition}

\begin{warning}
	The indexing convention used here is off by $ 1 $ from the indexing convention in \cite[\S1]{MR827262}.
	However, this indexing convention is what was later adopted by Simons--Sullivan \cite[\S1]{MR2365651}.
	See also \Cref{rem:Delignecohomology} for why $ k $ is the `right' index rather than $ k - 1 $.
\end{warning}

\begin{remark}
	When $ k = 0 $, the diagram \eqref{diag:diffcoh} is quite degenerate, and it will be convenient to define $ \Hhat^0(M;\ZZ) \colonequals \H^0(M;\ZZ) $.
\end{remark}

Now let us construct maps to fill in the ``differential cohomology'' diagram \eqref{diag:diffcoh}.

\begin{construction}[(characteristic class)]\label{cons:charclassmap}
	There is a \emph{characteristic class} map 
	\begin{equation*}
		\cc \colon \fromto{\Hhat^k(M;\ZZ)}{\H^k(M;\ZZ)}
	\end{equation*}
	defined as follows.
	Since $ \Zsm_{k-1}(M;\ZZ) $ is a free $ \ZZ $-module and the quotient map $ \surjto{\RR}{\RR/\ZZ} $ is an epimorphism, any homomorphism $ \chi \colon \fromto{\Zsm_{k-1}(M;\ZZ)}{\RR/\ZZ} $ lifts to a homomorphism
	\begin{equation*}
		\chitilde \colon \fromto{\Zsm_{k-1}(M;\ZZ)}{\RR} \period
	\end{equation*}
	Now define a homomorphism $ I(\chitilde) \colon \fromto{\Csm_k(M;\ZZ)}{\ZZ} $ by the assignment
	\begin{equation*}
		c \mapsto - \chitilde(\partial c) + \int_{c} \curv(\chi)  \period
	\end{equation*}

	Since $ \curv(\chi) $ is closed, $ I(\chitilde) $ defines a cocycle.
	Moreover, $ I(\chitilde) $ takes integral values, and the cohomology class $ [I(\chitilde)] \in \H^k(M;\ZZ) $ does not depend on the choice of lift $ \chitilde $.
	We define the characteristic class map $ \cc $ by the assignment
	\begin{align*}
		\cc \colon \Hhat^k(M;\ZZ) &\to \H^k(M;\ZZ) \\
		\chi &\mapsto [I(\chitilde)] \phantom{;\ZZ)} \period
	\end{align*} 
\end{construction}

\begin{warning}
	Simons and Sullivan \cite{MR2365651} denote the characteristic class map $ \cc $ by `$ \ch $'.
\end{warning}

\begin{construction}
	Consider the universal coefficient sequence
	\begin{equation*}
		\begin{tikzcd}[sep=2em]
			0 \arrow[r] & \Ext_{\ZZ}^1(\H_{i-1}(M;\ZZ),\RR/\ZZ) \arrow[r] & \H^i(M;\RR/\ZZ) \arrow[r, "\ang{-,-}"] & \Hom_{\ZZ}(\H_i(M;\ZZ),\RR/\ZZ) \arrow[r] & 0 \comma 
		\end{tikzcd}
	\end{equation*}
	where the morphism $ \ang{-,-} $ is given by sending the class of a cocycle $ u $ to the homomorphism
	\begin{align*}
		\ang{u,-} \colon \H_i(M;\ZZ) &\to \RR/\ZZ \\
		[z] &\mapsto u(z) \period
	\end{align*}
	Since the circle $ \RR/\ZZ $ is an injective $ \ZZ $-module, for any $ \ZZ $-module $ A $ and integer $ j > 0 $, we have \smash{$ \Ext_{\ZZ}^j(A,\RR/\ZZ) = 0 $}.
	In particular, $ \ang{-,-} $ is an isomorphism.

	Setting $ i = k - 1 $, precomposition with the quotient map $ \surjto{\Zsm_{k-1}(M;\ZZ)}{\H_{k-1}(M;\ZZ)} $ defines an injection
	\begin{equation*}
		\begin{tikzcd}[sep=1.5em]
			\H^i(M;\RR/\ZZ) \arrow[r, "\sim"{yshift=-0.1em}] & \Hom_{\ZZ}(\H_i(M;\ZZ),\RR/\ZZ) \arrow[r, hook] & \Hom_{\ZZ}(\Zsm_{k-1}(M;\ZZ),\RR/\ZZ) \period 
		\end{tikzcd}
	\end{equation*}
	It follows from the definitions that this factors through $ \Hhat^k(M;\ZZ) $.
	We simply denote this composite by $ \ang{-,-} \colon \incto{\H^{k-1}(M;\RR/\ZZ)}{\Hhat^k(M;\ZZ)} $.
\end{construction}

\begin{construction}
	Define a map $ \iota \colon \fromto{\Omega^{k-1}(M)}{\Hhat^k(M;\ZZ)} $ by setting
	\begin{equation*}
		\iota(\omega)(z) \colonequals \exp\paren{\textstyle 2\pi i \int_z \omega} 
	\end{equation*}
	for every smooth $ (k-1) $-cycle $ z $.
	By Stokes' Theorem, we see that $ \curv(\iota(\omega)) = \d\omega $.

	We have an $ \RR $-valued lift of $ \iota(\omega) $ given by setting
	\begin{equation*}
		\iotatilde(\omega)(z) \colonequals \int_z \omega
	\end{equation*}
	for every smooth $ (k-1) $-cycle $ z $.
	So by Stokes' Theorem we have
	\begin{align*}
		I(\iotatilde(\omega))(c) &= -\,\iotatilde(\omega)(\partial c) + \int_c \curv(\iota(\omega)) \\ 
		&= -\int_{\partial c} \omega + \int_{c} \d\omega = 0
	\end{align*}
	for every smooth $ k $-chain $ c $.
	Hence $ \cc \of \iota = 0 $.

	We see that $ \iota \colon \fromto{\Omega^{k-1}(M)}{\Hhat^k(M;\ZZ)} $ has kernel those closed forms $ \omega $ such that $ \int_{z} \omega $ is an integer for all $ z \in \Zsm_{k-1}(M;\ZZ) $. 
	That is,
	\begin{equation*}
		\ker(\iota) = \Omegacl^{k-1}(M)_{\ZZ}
	\end{equation*}
	is the group of closed $ (k-1) $-forms with integral periods.
	Hence $ \iota $ descends to an injection
	\begin{equation*}
		\iota \colon \into{\Omega^{k-1}(M)/\Omegacl^{k-1}(M)_{\ZZ}}{\Hhat^{k}(-;\ZZ)} \period
	\end{equation*}
\end{construction}

%-------------------------------------------------------------------%
%-------------------------------------------------------------------%
%  The differential cohomology hexagon                              %
%-------------------------------------------------------------------%
%-------------------------------------------------------------------%

\subsection{The differential cohomology hexagon}

\begin{notation}
	Write $ \Mfld $ for the category of smooth manifolds and $ \GrAb $ for the category of graded abelian groups. 
\end{notation}

\begin{theorem}[{(Simons--Sullivan \cite[Theorem 1.1]{MR2365651})}]\label{thm:SimonsSullivanunique}
	There is an essentially unique functor
	\begin{equation*}
		\Hhat^{*}(-;\ZZ) \colon \fromto{\Mfldop}{\GrAb}
	\end{equation*}
	equipped with natural transformations
	\begin{enumerate}[label=\stlabel{thm:SimonsSullivanunique}, ref=\arabic*]
		\item $ \ang{-,-} \colon \fromto{\H^{*-1}(-;\RR/\ZZ)}{\Hhat^{*}(-;\ZZ)} $, 

		\item $ \iota \colon \fromto{\Omega^{*-1}(M)/\Omegacl^{*-1}(M)_{\ZZ}}{\Hhat^{*}(-;\ZZ)} $,

	 	\item $ \cc \colon \fromto{\Hhat^{*}(-;\ZZ)}{\H^{*}(-;\ZZ)} $, 

	 	\item and $ \curv \colon \fromto{\Hhat^{*}(-;\ZZ)}{\Omegacl^{*}(-)_{\ZZ}} $
	\end{enumerate}
	filling in the ``differential cohomology hexagon''
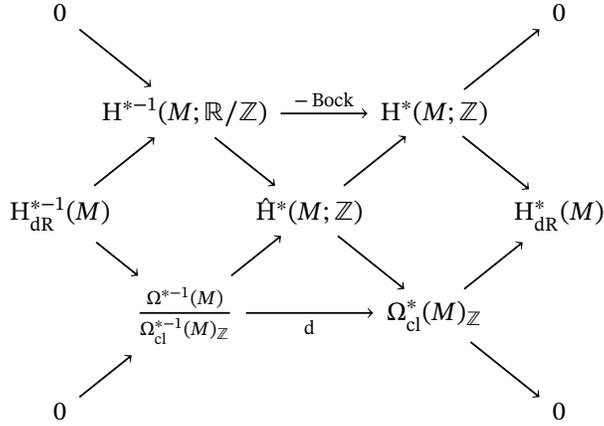
	\begin{equation*}
		\begin{tikzcd}[column sep={10ex,between origins}, row sep={8ex,between origins}]
			0 \arrow[dr] & & & & 0 \\
			& \H^{*-1}(M;\RR/\ZZ) \arrow[rr, "-\Bock"] \arrow[dr] & & \H^*(M;\ZZ) \arrow[dr] \arrow[ur] & \\
			\HdR^{*-1}(M) \arrow[ur] \arrow[dr] & & \Hhat^*(M;\ZZ) \arrow[ur] \arrow[dr] & & \HdR^*(M) \\
			& \frac{\Omega^{*-1}(M)}{\Omegacl^{*-1}(M)_{\ZZ}} \arrow[rr, "\d"'] \arrow[ur] & & \Omegacl^*(M)_{\ZZ} \arrow[ur] \arrow[dr] & \\
			0 \arrow[ur] & & & & 0 
		\end{tikzcd}
	\end{equation*}
	so that the diagonal sequences are exact.
\end{theorem}

\noindent Any functor $ \Hhat^{*}(-;\ZZ) \colon \fromto{\Mfldop}{\GrAb} $ satisfying these properties is called \emph{ordinary differential cohomology}. 

\begin{remark}[(Deligne's model)]\label{rem:Delignecohomology}
	Motivated by Deligne cohomology in Hodge theory \cites[\S2.2]{MR498551}[\S12.3]{MR2451566}, we can consider the smooth version of the Deligne complex on a manifold $ M $.
	Write $ \ZZ(k) $ for the complex of sheaves on $ M $
	\begin{equation*}
		\begin{tikzcd}[sep=1.5em]
			0 \arrow[r] & \ZZ \arrow[r, hook] & \Omega^0 \arrow[r, "\d"] & \Omega^1 \arrow[r, "\d"] & \cdots \arrow[r] & \Omega^{k-1} \arrow[r, ""] & 0 \comma
		\end{tikzcd}
	\end{equation*}
	where $ \Omega^i $ is in degree $ i + 1 $.
	The \emph{$ k $-th smooth Deligne cohomology group} of $ M $ is the sheaf cohomology (i.e., hypercohomology) group $ \H^k(M;\ZZ(k)) $.
	We will see later that smooth Deligne cohomology agrees with ordinary differential cohomology (see \Cref{lem:Delignemodel}).
\end{remark}

\begin{questions}\label{qst:extensions}
	There are a number of questions that naturally arise
	\begin{enumerate}[label=\stlabel{qst:extensions}, ref=\arabic*]
		\item Is there differential $ \K $-theory?

		Yes! Hopkins--Singer \cite{HopkinsSinger} define differential $\K$-theory. 
		Simons--Sullivan \cites{MR2732065}{MR3220448} tell a similar story, and define differential $ \K $-theory in terms of vector bundles with connection.
		We study this in \cref{subsec:diffKtheory}. 

		\item What about differential [favorite cohomology theory]? 

		Also yes, but the theory is more complicated. 
		The fundamental observation is that everything we've considered comes from a sheaf of abelian groups or chain complexes (which we regard as spectra) on the category of \textit{all} smooth manifolds.
		We begin to set up this theory in \Cref{sec:basicsetup}.

		Moreover, the \category $ \Sh(\Mfld;\Sp) $ of sheaves of spectra on the category of manifolds has rich structure that gives rise to a ``differential cohomology hexagon''
		associated to every object.
		We study this in \Cref{sec:stable}.
	\end{enumerate}
\end{questions}

\begin{remark}
	The category $ \Sh(\Mfld;\Set) $ is really the right place for moduli spaces of manifolds to live, and Fréchet manifolds embed as a full subcategory of $ \Sh(\Mfld;\Set) $.
	See \cref{subsec:infdimMan}.
\end{remark}

There are many applications of this perspective on differential cohomology that we study throughout this book.
See, in particular, \Cref{part:applications}.

\newpage
%!TEX root = ../diffcoh.tex

%-------------------------------------------------------------------%
%-------------------------------------------------------------------%
%  Basics of sheaves on manifolds                                 %
%-------------------------------------------------------------------%
%-------------------------------------------------------------------%

\section{Basics of sheaves on manifolds}\label{sec:sheavesonMan}\label{sec:basicsetup}
\textit{by Peter Haine}

The purpose of this chapter is to begin to set up the basics of differential cohomology theories as sheaves on the category of all manifolds.
\Cref{subsec:ShManDef} starts with the basic definitions.
\Cref{subsec:DRreminder} gives a reminder on derived \categories and their relation to spectra so that we can give examples of sheaves on the category of manifolds in \cref{subsec:firstexamples}.
In \cref{subsec:checkonstalks}, we explain why in all situations of interest, we can check equivalences of
differential cohomology theories ``on stalks''.
\Cref{subsec:Cartsheaves} gives an alternative description of the \category of sheaves on manifolds in terms of sheaves on the smaller category of Euclidean spaces.
\Cref{subsec:excision} is a digression giving a reformulation of the sheaf condition in terms of an
\textit{excision} condition (or \textit{Mayer--Vietoris} property) and a ``finiteness'' condition.
We finish the chapter with a digression explaining Losik and Hain's results embedding infinite dimensional manifolds into sheaves of sets on the category of (finite dimensional) manifolds (\cref{subsec:infdimMan}).

%-------------------------------------------------------------------%
%  Definitions                                                      %
%-------------------------------------------------------------------%

\subsection{Definitions}\label{subsec:ShManDef}

\begin{notation}
	We write $ \Mfld $ for the (ordinary) category of smooth manifolds, including the empty manifold.
	The category $ \Mfld $ has a Grothendieck topology where the covering families are families of open embeddings
	\begin{equation*}
		\{ j_{\alpha} \colon \incto{U_{\alpha}}{M} \}_{\alpha \in A}
	\end{equation*}
	such that the family of open sets $ \{j_{\alpha}(U_{\alpha})\}_{\alpha \in A} $ is an open cover of $ M $. 
	Whenever we regard $ \Mfld $ as a site, we use this topology.
\end{notation}

\begin{remark}
	Since the category $ \Mfld $ is equivalent to the category of manifolds with a fixed embedding into $ \RR^\infty $, the category $ \Mfld $ is essentially small.
\end{remark}

\begin{definition}\label{def:sheafonMan}
	Let $ C $ be a presentable \category.
	We write 
	\begin{equation*}
		\PSh(\Mfld;C) \colonequals \Fun(\Mfldop,C)
	\end{equation*}
	and write
	\begin{equation*}
		\Sh(\Mfld;C) \subset \PSh(\Mfld;C)
	\end{equation*}
	for the full subcategory spanned by the $ C $-valued sheaves on the site $ \Mfld $ with respect to the Grothendieck topology given by open covers.

	Explicitly, a $ C $-valued presheaf $ E \colon \fromto{\Mfldop}{C} $ is a sheaf if and only if for each manifold $ M $, the restriction $ \restrict{E}{\Open(M)} $ of $ E $ to the site $ \Open(M) $ of open submanifolds of $ M $ is a sheaf on the topological space $ M $.
\end{definition}

\begin{remark}[(on presentability)]
	Perhaps somewhat surprisingly, if $ C $ is \acategory with all limits and colimits, then the inclusion of $ C $-valued sheaves into $ C $-valued presheaves need not admit a left adjoint: presentability is the standard assumption which guarantees that this left adjoint exists.
	For this reason, we essentially always work with sheaved valued in a presentable \category.
\end{remark}

\begin{notation}
	We write $ \SMan \colon \fromto{\PSh(\Mfld;C)}{\Sh(\Mfld;C)} $ for the left adjoint to the inclusion, that is, the \textit{sheafification} functor.
\end{notation}

\begin{notation}
	We write $ \Set $ for the category of sets, $ \Spc $ for the \category of spaces, $ \Sp $ for the \category of spectra, and $ \Cat_{\infty} $ for the \category of \categories.
\end{notation}

\begin{example}\label{ex:representable_sheaf}
	Let $ \yo \colon \incto{\Mfld}{\PSh(\Mfld;\Set)} $ denote the Yoneda embedding.
	For each manifold $ M $, the representable presheaf $ \yo(M) $ is a sheaf.
	Unless noted otherwise, we simply write $ M $ for the sheaf $ \yo(M) $.
\end{example}

The following is the fundamental definition of this text:

\begin{definition}\label{def:diffcoh}
	The \category of \textit{differential cohomology theories} is the \category $ \Sh(\Mfld;\Sp) $ of sheaves of spectra on $ \Mfld $. 
\end{definition}

\noindent For most of this text we work in the generality of sheaves with values in a general presentable \category, or stable presentable \category. 
The main reason for doing this is because we have reason to consider sheaves of spaces, sheaves of chain complexes, and sheaves of spectra, and want to treat them on the same footing.

\begin{remark}
	We take the approach of Freed--Hopkins \cite{MR3049871} and consider sheaves on the category of smooth manifolds.
	The general setup here is very robust, and one can take the basic objects to be manifolds with corners without
	essential change to how theory works; this is the approach taken by Hopkins--Singer \cite{MR2192936} and
	Bunke--Nikolaus--Völkl \cite{MR3462099}.
\end{remark}

The first basic property we prove about sheaves on $ \Mfld $ is that morphism is an equivalence if and only if it is when evaluated on each Euclidean space.
For this, we use the fact that manifolds admit \textit{good covers}.

\begin{recollection}[(good covers)]
	Let $ M $ be an $ n $-manifold.
	An open cover $ \Uu $ of $ M $ is \textit{good} if for every finite set $ U_1,\ldots,U_m \in \Uu $ of opens in $ \Uu $, the intersection $ U_1 \intersect \cdots \intersect U_m $ is either empty or diffeomorphic to $ \RR^n $.
\end{recollection}

\begin{notation}
	Let $ T $ be a topological space and $ U \subset T $ be open.
	For every open cover $ \Uu $ of $ U $, write $ \Iup(\Uu) \subset \Open(T) $ for the full subposet  consisting of all nonempty finite intersections of elements in $ \Uu $.
\end{notation}

\begin{lemma}\label{lem:checkequivonRn}
	Let $ C $ be a presentable \category.
	A morphism $ f \colon \fromto{E}{E'} $ in $ \Sh(\Mfld;C) $ is an equivalence if and only if for each integer $ n \geq 0 $, the morphism $ f(\RR^n) \colon \fromto{E(\RR^n)}{E'(\RR^n)} $ is an equivalence in $ C $.
\end{lemma}

\begin{proof}
	Let $ M $ be a manifold and $ \Uu $ a good cover of $ M $.
	The morphism $ f $ induces a commutative square
	\begin{equation*}\label{eq:Cechequivs}
		\begin{tikzcd}
			E(M) \arrow[d, "f(M)"'] \arrow[r, "\sim"{yshift=-0.2em}] & \lim_{U \in \Iup(\Uu)^{\op}} E(U) \arrow[d] \\
			E'(M) \arrow[r, "\sim"{yshift=-0.2em}] & \lim_{U \in \Iup(\Uu)^{\op}} E'(U) \comma
		\end{tikzcd}
	\end{equation*}
	where the horizontal morphisms are equivalences because $ E $ and $ E' $ are sheaves.
	Since the cover $ \Uu $ is good and $ f $ is an equivalence on Euclidean spaces, we see that the induced morphism
	\begin{equation*}
		f \colon \fromto{\restrict{E}{\Iup(\Uu)^{\op}}}{\restrict{E'}{\Iup(\Uu)^{\op}}}
	\end{equation*}
	of $ \Iup(\Uu)^{\op} $-indexed diagrams in $ C $ is an equivalence, which proves the claim.
\end{proof}

%-------------------------------------------------------------------%
%  Reminder on derived ∞-categories and Eilenberg–MacLane spectra   %
%-------------------------------------------------------------------%

\subsection{Reminder on derived \texorpdfstring{$\infty$}{∞}-categories and Eilenberg--MacLane spectra}\label{subsec:DRreminder}

In order to give some important examples of sheaves on $ \Mfld $, we review the basics of derived \categories of rings and their relation to spectra.

\begin{notation}[(derived \categories)]
	Let $ R $ be a ring.
	We write $ \Ch(R) $ for the category of chain complexes of $ R $-modules.
	We write $ \D(R) $ for the \emph{derived \category} of $ R $ obtained from the category $ \Ch(R) $ by formally inverting the quasi-isomorphisms \cite[\HAthm{Definition}{1.3.5.8}, \HAthm{Propositon}{1.3.5.15}, \& \HAthm{Remark}{7.1.1.16}]{HA}.
	That is, $ \D(R) $ is the universal \category equipped with a functor $ \fromto{\Ch(R)}{\D(R)} $ carrying quasi-isomorphisms in $ \Ch(R) $ to equivalences in the \category $ \D(R) $.
	Note that for every map of rings $ \fromto{S}{R} $, the forgetful functor $ \fromto{\Ch(R)}{\Ch(S)} $ preserves quasi-isomorphisms, hence induces a forgetful functor $ \fromto{\D(R)}{\D(S)} $.
\end{notation}

\begin{recollection}[(Eilenberg--MacLane spectra)]\label{rec:EMspectra}
	The inclusion $ \Ab \subset \Sp $ of the category of abelian groups into the category of spectra as those spectra with homotopy groups in degree $ 0 $ (i.e., ordinary cohomology theories) extends to a right adjoint functor
	\begin{equation*}
		\EM \colon \fromto{\D(\ZZ)}{\Sp} \period
	\end{equation*}
	The functor $ \EM $ is called the \emph{Eilenberg--MacLane} functor \HA{Example}{1.3.3.5}.
	For a ring $ R $, we also simply write $ \EM $ for the composite
	\begin{equation*}
		\begin{tikzcd}
			\D(R) \arrow[r] & \D(\ZZ) \arrow[r, "\EM"] & \Sp
		\end{tikzcd}
	\end{equation*}
	for the composite of the forgetful functor $ \fromto{\D(R)}{\D(\ZZ)} $ with the Eilenberg--MacLane functor.
	The spectrum $ \HR $ represents ordinary cohomology with coefficients in $ R $.
\end{recollection}

\begin{recollection}[{($ \HR $-modules)}]
	Every spectrum in the image of $ \EM \colon \fromto{\D(R)}{\Sp} $ is a module over the Eilenberg--MacLane spectrum $ \HR $ representing ordinary cohomology with coefficients in $ R $.
	Moreover, the Eilenberg--MacLane functor induces an equivalence
	\begin{equation*}
		\equivto{\D(R)}{\Mod(\HR)}
	\end{equation*}
	between the derived \category $ \D(R) $ and the \category $ \Mod(\HR) $ of $ \HR $-module spectra \HA{Proposition}{7.1.4.6}.
	Under this equivalence $ \equivto{\D(R)}{\Mod(\HR)} $, the functor $ \EM \colon \fromto{\D(R)}{\Sp} $ corresponds to the forgetful functor $ \fromto{\Mod(\HR)}{\Sp} $
\end{recollection}

%-------------------------------------------------------------------%
%  First examples                                                   %
%-------------------------------------------------------------------%

\subsection{First examples}\label{subsec:firstexamples}

Now we give some examples of sheaves on manifolds coming from topological spaces, complexes of differential forms, and bundles.

%-------------------------------------------------------------------%
%  Topological spaces                                               %
%-------------------------------------------------------------------%

\subsubsection{Topological spaces}

\begin{notation}\label{ntn:Top}
	Write $ \Top $ for the category of topological spaces.
\end{notation}

\begin{construction}\label{ex:Topyoneda}
	Define a restricted Yoneda functor $ y_{\Top} $ by
	\begin{align*}
		y_{\Top} \colon \Top &\to \PSh(\Mfld;\Set) \\
		T &\mapsto [M \mapsto \Map_{\Top}(M,T)] \period
	\end{align*}
	Since continuous functions glue over open covers, the assignment $ M \mapsto \Map_{\Fre}(M,T) $ is a sheaf on $ \Mfld $.
	That is, $ y_{\Top} $ factors through $ \Sh(\Mfld;\Set) $.
	Hence every topological space defines a sheaf on $ \Mfld $.
\end{construction}

%-------------------------------------------------------------------%
%  Differential forms                                               %
%-------------------------------------------------------------------%

\subsubsection{Differential forms}

\begin{example}[(differential forms)]\label{ex:Omegai_is_a_sheaf}
	Let $ i \geq 0 $ be an integer.
	The functor
	\begin{equation*}
		\Omega^i \colon \fromto{\Mfldop}{\Vect(\RR)}
	\end{equation*} 
	sending manifold $ M $ to vector space $ \Omega^i(M) $ of differential $ i $-forms on $ M $ with functoriality given by pullback of bundles is a sheaf.
	Note that by the Yoneda Lemma, there is a natural isomorphism
	\begin{equation*}
		\Map_{\Sh(\Mfld;\Set)}(M,\Omega^i) \isomorphic \Omega^i(M) \period
	\end{equation*}
\end{example}

\begin{example}[(de Rham complex)]
	Putting togther all $ i $ at once, the functor
	\begin{equation*}
		\Omegabullet \colon \fromto{\Mfldop}{\Ch(\RR)}
	\end{equation*} 
	sending manifold $ M $ to its de Rham complex $ \Omegabullet(M) $ is a sheaf of chain complexes on $ \Mfld $.

	Even better, $ \Omegabullet $ is a sheaf in the derived sense: the composite
	\begin{equation*}
		\begin{tikzcd}
			\Mfldop \arrow[r, "\Omegabullet"] & \Ch(\RR) \arrow[r] & \D(\RR)
		\end{tikzcd}
	\end{equation*}
	with the localization functor $ \fromto{\Ch(\RR)}{\D(\RR)} $ is a sheaf valued in the \textit{\category} $ \D(\RR) $.
\end{example}

%-------------------------------------------------------------------%
%  Bundles & sheaves                                                %
%-------------------------------------------------------------------%

\subsubsection{Bundles \& sheaves}\label{subsec:bundles_and_sheaves}

\begin{example}[(vector bundles)]
	Write
	\begin{equation*}
		\VectRR \colon \fromto{\Mfldop}{\Gpd}
	\end{equation*} 
	for the functor sending a manifold $ M $ to the groupoid of (finite dimensional) real vector bundles on $ M $ and bundle isomorphisms, with functoriality given by pullback of bundles.
	Again, the local nature of the definition of a vector bundle ensures that $ \VectRR  $ is a sheaf of groupoids on $ \Mfld $.
\end{example}

\begin{example}[(principal bundles)]\label{ex:BunG}
	Let $ G $ be a Lie group.
	Write
	\begin{equation*}
		\BunG \colon \fromto{\Mfldop}{\Gpd}
	\end{equation*} 
	for the functor sending a manifold $ M $ to the groupoid of (smooth) principal $ G $ bundles on $ M $ and bundle isomorphisms, with functoriality given by pullback of bundles.
	The locally triviality of principal bundles implies that $ \BunG $ is a sheaf of groupoids on $ \Mfld $.
\end{example}

\begin{example}[(principal bundles with connection)]\label{ex:BunGnabla}
	Let $ G $ be a Lie group.
	Write
	\begin{equation*}
		\BunGnabla \colon \fromto{\Mfldop}{\Gpd}
	\end{equation*} 
	for the functor sending a manifold $ M $ to the groupoid of (smooth) principal $ G $ bundles on $ M $ \textit{with connection} and bundle isomorphisms respecting connections, with functoriality given by pullback of bundles.
	Explicitly, an object of $ \BunGnabla(M) $ consists of a pair $ (P,\theta) $ of a principal $ G $-bundle $ \fromto{P}{M} $ and a connection $ 1 $-form $ \theta \in \Omega^1(P;\g) $.
	See \Cref{ChernWeilTheory} for background on connections.
	A morphism
	\begin{equation*}
		\fromto{(P,\theta)}{(P',\theta')}
	\end{equation*}
	in $ \BunGnabla(M) $ consists of an isomorphism of principal $ G $-bundles $ f \colon \isomto{P}{P'} $ such that $ \fupperstar(\theta') = \theta $.

	Again, the local nature of the definition of a bundle with connection ensures that $ \BunGnabla  $ is a sheaf of groupoids on $ \Mfld $.
\end{example}

\begin{warning}
	In \Cref{ex:BunG,ex:BunGnabla}, it is very important that we have not passed to \textit{isomorphism classes} of principal $ G $-bundles (with connection).
	The reason is that isomorphism classes do not glue, i.e., do not form a sheaf.
\end{warning}

The sheaves $ \BunG $ and $ \BunGnabla $ are of great importance.
In \Cref{char_class_part}, we study them in detail.

\begin{example}[(sheaves)]\label{ex:ShLC}
	For each manifold $ M $, write $ \Sh(M) $ for the \category of sheaves of spaces on $ M $, and $ \LC(M) \subset \Sh(M) $ for the full subcategory spanned by the locally constant sheaves of spaces.
	The assignment $ \goesto{M}{\Sh(M)} $ extends to a functor
	\begin{equation*}
		\Sh \colon \fromto{\Mfldop}{\Cat_{\infty}}
	\end{equation*} 
	with functoriality given by pullback of sheaves.
	The functor $ \Sh $ is a sheaf of (large) \categories on $ \Mfld $ \HTT{Theorem}{6.1.3.9}.
	Since locally constant sheaves are preserved by sheaf pullback and local constancy is a local condition, the subfunctor $ \LC \subset \Sh $ is also a sheaf of (large) \categories on $ \Mfld $.
\end{example}

%-------------------------------------------------------------------%
%  Checking equivalences on stalks                                  %
%-------------------------------------------------------------------%

\subsection{Checking equivalences on stalks}\label{subsec:checkonstalks}

We now explain that equivalences of sheaves on $ \Mfld $ with values in a \textit{compactly generated} \category
(e.g., $ \Spc $, $ \Sp $, $ \D(R) $) can be checked on ``stalks'' at the origins in $ \RR^n $ for $ n \geq 0 $.
The proof of this requires a few technical detours which we defer to \Cref{subsec:points}.

\begin{notation}
	Let $ M $ be a manifold and $ x \in M $.
	We write $ \Open_x(M) \subset \Open(X) $ for the full subposet spanned by the open neighborhoods of $ x \in M $.
\end{notation}

\begin{definition}\label{def:stalkCG}
	Let $ C $ be a compactly generated \category, $ E \in \Sh(\Mfld;C) $ a $ C $-valued sheaf on $ \Mfld $, $ M $ a manifold, and $ x \in M $.
	The \textit{stalk} of $ E $ at $ x \in M $ is the filtered colimit
	\begin{equation}\label{eq:stalkfilteredcolim}
		\xupperstar (E) \colonequals \colim_{U \in \Open_x(M)^{\op}} E(U)
	\end{equation}
	in $ C $.
\end{definition}

\begin{warning}
	It is important that we have phrased \Cref{def:stalkCG} only for compactly generated coefficients.
	It is true that for any presentable \category $ C $, manifold $ M $, and point $ x \in M $, there is a stalk functor $ \xupperstar \colon \fromto{\Sh(\Mfld;C)}{C} $ (see \Cref{cons:stalk}).
	However, if $ C $ is not compactly generated then $ \xupperstar $ need not be computed by the filtered colimit \eqref{eq:stalkfilteredcolim}.
\end{warning}

\begin{notation}\label{ntn:origin}
	For each integer $ n \geq 0 $ and number $ r \in \RR_{>0} $, write $ 0_n \in \RR^n $ for the origin, and write
	\begin{equation*}
		\Bup_{\RR^n}(r) \subset \RR^n
	\end{equation*}
	for the open ball in $ \RR^n $ of radius $ r $ centered at the $ 0_n $.
\end{notation}

\begin{nul}
	Let $ E \colon \fromto{\Mfldop}{C} $ be a sheaf on $ \Mfld $.
	Note that the stalk $ 0_n\upperstar(E) $ can be computed as the colimit
	\begin{equation*}
		0_n\upperstar(E) \equivalent \colim_{k \in \NN} E(\Bup_{\RR^n}(1/k)) \period
	\end{equation*}
\end{nul}

The following result comes from the functoriality of a sheaf on $ \Mfld $ in \textit{all} manifolds, the fact that for ever $ n $-manifold $ M $ and point $ x \in M $, there exists an open embedding $ j \colon \incto{\RR^n}{M} $ such that $ j(0_n) = x $, and that equivalences in sheaves on $ M $ can be checked on stalks.
In \Cref{subsec:points} we provide a detailed proof.

\begin{proposition}[(\Cref{app.prop:hypercompleteness})]\label{prop:hypercompleteness}
	Let $ C $ be a compactly generated \category.
	A morphism $ f $ in $ \Sh(\Mfld;C) $ is an equivalence if and only if for each integer $ n \geq 0 $, the morphism $ 0_n\upperstar(f) $ is an equivalence in $ C $.
\end{proposition}

\begin{remark}\label{rem:FreedHopkinsmodel}
	\Cref{prop:hypercompleteness} is important from our perspective.
	Freed and Hopkins work with differential cohomology theories using the language of simplicial sheaves and model categories \cite{MR3049871}.  
	Combining \Cref{prop:hypercompleteness} with \cite[\HTTthm{Remark}{6.5.2.2} \& \HTTthm{Proposition}{6.5.2.14}]{HTT} shows that the model structure on simplicial presheaves on $ \Mfld $ considered in \cite[§5]{MR3049871} pre\-sents the \category $ \Sh(\Mfld;\Spc) $.
\end{remark}

\begin{warning}
	\Cref{prop:hypercompleteness} does \textit{not} hold when $ C $ is replaced by an arbitrary presentable \category.
\end{warning}

%-------------------------------------------------------------------%
%  Sheaves on the Euclidean site                                    %
%-------------------------------------------------------------------%

\subsection{Sheaves on the Euclidean site}\label{subsec:Cartsheaves}

In this section, we refine \Cref{lem:checkequivonRn} in the following manner.
Since every manifold admits an open cover by Euclidean spaces, the category of sheaves of \textit{sets} on $ \Mfld $ is equivalent to sheaves of sets on the full subcategory spanned by the Euclidean spaces.
We prove an analogous result for sheaves of \textit{spaces}; this is not immediate in the higher-categorical setting \SAG{Counterexample}{20.4.0.1}.
The reason for this subtlety is exactly the failure of Whitehead's Theorem to hold in an arbitrary \category of sheaves of spaces.
However, \Cref{prop:hypercompleteness} implies that Whitehead's Theorem holds in $ \Sh(\Mfld;\Spc) $; a general result \cites[Appendix A]{arXiv:2001.00319}[Corollary 3.12.13]{exodromy} implies that sheaves on the site of Euclidean spaces and sheaves on $ \Mfld $ coincide.

\begin{definition}\label{def:Euclideansite}
	The \textit{Euclidean site} is the full subcategory $ \Euc \subset \Mfld $ spanned by the Euclidean spaces $ \RR^n $ for $ n \geq 0 $, with the induced Grothendieck topology.
\end{definition}

The proof of the following is quite short. 
However, it involves some technical tools we have not yet introduced, so we defer it to \cref{subsec:points}.

\begin{notation}
	Let $ C $ be a presentable \category. 
	Write $ j \colon \incto{\Eucop}{\Mfldop} $ for the inclusion, and
	\begin{equation*}
		\jlowerstar \colon \incto{\PSh(\Euc;C)}{\PSh(\Mfld;C)}
	\end{equation*}
	for the fully faithful functor given by right Kan extension along $ j $.
	Restriction of presheaves is left adjoint to $ \jlowerstar $; we denote this left adjoint by $ \jupperstar $ or $ \restrict{(-)}{\Eucop} $.
\end{notation}

\begin{lemma}[(\Cref{app.lem:hypersheavesonCart})]\label{lem:hypersheavesonCart}
	Let $ C $ be a presentable \category and $ F \in \PSh(\Mfld;C) $.
	Then:
	\begin{enumerate}[label=\stlabel{lem:hypersheavesonCart}, ref=\arabic*]
		\item\label{lem:hypersheavesonCart.1} The functors $ \jlowerstar \colon \incto{\PSh(\Euc;C)}{\PSh(\Mfld;C)} $ and $ \jupperstar \colon \incto{\PSh(\Mfld;C)}{\PSh(\Euc;C)} $ preserve sheaves.

		\item\label{lem:hypersheavesonCart.2} The adjoint functors $ \adjto{\jupperstar}{\Sh(\Mfld;C)}{\Sh(\Euc;C)}{\jlowerstar} $ are inverse equivalences of \categories.

		\item\label{lem:hypersheavesonCart.3} If $ \jupperstar(F) $ is a sheaf on $ \Euc $, then the unit $ \fromto{F}{\jlowerstar\jupperstar(F)} $ exhibits $ \jlowerstar\jupperstar(F) $ as the sheafification $ \SMan(F) $.

		\item\label{lem:hypersheavesonCart.4} If $ \jupperstar(F) $ is a sheaf on $ \Euc $, then for all $ n \geq 0 $, the unit $ \fromto{F}{\SMan(F)} $ is an equivalence when evaluated on $ \RR^n $.
		In particular, the unit $ \fromto{F}{\SMan(F)} $ induces an equivalence on global sections.
	\end{enumerate}
\end{lemma}

%-------------------------------------------------------------------%
%  Digression: excision & the sheaf condition                       %
%-------------------------------------------------------------------%

\subsection{Digression: excision \& the sheaf condition}\label{subsec:excision}

The goal of this section is to prove a convenient reformulation of the sheaf condition in terms of an \textit{excision} property.
We do not make use of the reformation in this text, but present it here because it is the manifold analogue of \textit{Nisnevich excision} from algebraic geometry \cites[\SAGthm{Proposition}{B.5.1.1}]{SAG}[\S3.2]{MR3679884}[\S3.1, Proposition 1.4]{MR1813224}.
Another way of explaining the following result is that it says that a presheaf on $ \Mfld $ is a sheaf if and only if it satisfies the \textit{Mayer--Vietoris} property and transforms countable increasing chains of open submanifolds to limits.

\begin{theorem}[{\cite[Theorem 5.1]{BEBdBP:Classifying}}]\label{thm:excision}
	Let $ C $ be a presentable \category.
	A $ C $-valued presheaf $ F \colon \fromto{\Mfldop}{C} $ on $ \Mfld $ is a sheaf if and only if $ F $ satisfies the following conditions:
	\begin{enumerate}[label=\stlabel{thm:excision}, ref=\arabic*]
		\item\label{thm:excision.1} The object $ F(\emptyset) $ is terminal in $ C $.

		\item\label{thm:excision.2} For every manifold $ M $ and pair of open subsets $ U, V \subset M $ such that $ U \union V = M $, the induced square 
		\begin{equation*}
			\begin{tikzcd}
				F(M) \arrow[r] \arrow[d] & F(V) \arrow[d] & \\
				F(U) \arrow[r] & F(U \intersect V) 
			\end{tikzcd}
		\end{equation*}
		is a pullback square in $ C $.

		\item\label{thm:excision.3} For every manifold $ M $ and $ \NN $-indexed sequence of open sets
		\begin{equation*}
			U_0 \subset U_1 \subset \cdots \subset M
		\end{equation*}
		such that $ \Union_{n \geq 0} U_n = M $, the induced morphism
		\begin{equation*}
			\fromto{F(M)}{\lim_{n \geq 0} F(U_n)}
		\end{equation*}
		is an equivalence in $ C $.
	\end{enumerate}
\end{theorem}

\begin{remark}
	Pavlov~\cite[Theorem 1.7]{Pav22} proves that for $\sSet$-valued presheaves on $\Mfld$, these conditions can be considerably simplified to two conditions on two-element open covers and pairwise disjoint open covers.
\end{remark}

We do not have occasion to use \Cref{thm:excision} in this text, but include it for completeness and because it is useful.
For example, \Cref{thm:excision} is crucial to work of Berwick-Evans--Boavida de Brito--Pavlov \cite{BEBdBP:Classifying} extending results of Madsen--Weiss \cite[Appendix A]{MR2335797}.
See \Cref{remark:BE-BdB-P} for more details.

The idea of \Cref{thm:excision} is as follows.
Conditions \enumref{thm:excision}{1} and \enumref{thm:excision}{2} guarantee that $ F $ satisfies the sheaf condition with respect to finite open covers.
Given descent with respect to finite open covers, by writing a countable cover as a union of a sequence of finite covers of smaller subspaces, \enumref{thm:excision}{3} implies descent with respect to \textit{countable} open covers.
Note that implicit in \Cref{thm:excision} is the claim that descent with respect to countable open covers implies descent with respect to arbitrary open covers.

Since the sheaf condition on $ \Mfld $ is defined after restriction to each manifold, \Cref{thm:excision} follows from an analogous rephrasing of the sheaf condition for a presheaf on an individual manifold (\Cref{prop:Lindelofexcision}).
The manifold structure isn't really used here; all that is necessary is that an open cover of an open subset of a manifold admits a countable subcover.
Hence we work at this level of generality.

\begin{observation}\label{obs:descfinitecovers}
	Let $ T $ be a topological space and $ C $ a presentable \category.
	Since limits of finite cubes can be written as iterated pullbacks, the following are equivalent for a presheaf $ F \in \PSh(T;C) $ on $ T $:
	\begin{enumerate}[label=\stlabel{obs:descfinitecovers}, ref=\arabic*]
		\item The presheaf $ F $ satisfies descent with respect to nonempty finite covers.

		\item For all opens $ U, V \subset T $, the induced square 
		\begin{equation*}
			\begin{tikzcd}
				F(U \union V) \arrow[r] \arrow[d] & F(V) \arrow[d] & \\
				F(U) \arrow[r] & F(U \intersect V) 
			\end{tikzcd}
		\end{equation*}
		is a pullback square in $ C $.
	\end{enumerate}
\end{observation}

\begin{recollection}\label{rec:Lindelof}
	A topological space $ T $ is \textit{Lindelöf} if every open cover of $ T $ has a countable subcover.

	The following conditions are equivalent for a topological space $ T $:
	\begin{enumerate}[label=\stlabel{rec:Lindelof}, ref=\arabic*]
		\item\label{rec:Lindelof.1} Every open subspace of $ T $ is Lindelöf.

		\item\label{rec:Lindelof.2} Every subspace of $ T $ is Lindelöf.
	\end{enumerate}
	We say that $ T $ is \textit{hereditarily Lindelöf} if $ T $ satisfies the equivalent conditions \enumref{rec:Lindelof}{1}--\enumref{rec:Lindelof}{2}.

	Note that every second-countable topological space (e.g., manifold) is hereditarily Lindelöf.
\end{recollection}

\begin{lemma}\label{lem:sheafLindelof}
	Let $ T $ be a hereditarily Lindelöf topological space and $ C $ a presentable \category.
	The following are equivalent for a presheaf $ F \in \PSh(T;C) $ on $ T $:
	\begin{enumerate}[label=\stlabel{lem:sheafLindelof}, ref=\arabic*]
		\item\label{lem:sheafLindelof.1} The presheaf $ F $ is a sheaf on $ T $.

		\item\label{lem:sheafLindelof.2} The presheaf $ F $ satisfies descent with respect to countable open covers.
	\end{enumerate}
\end{lemma}

\begin{proof}
	Clearly \enumref{lem:sheafLindelof}{1} $ \Rightarrow $ \enumref{lem:sheafLindelof}{2}.
	To see that \enumref{lem:sheafLindelof}{2} $ \Rightarrow $ \enumref{lem:sheafLindelof}{1}, let $ U \subset T $ be open and let $ \Uu $ be an open cover of $ U $
	Since $ T $ is hereditarily Lindelöf, there exists a countable subset $ \Uu_0 \subset \Uu $ that also covers $ U $. 
	To conclude, note that have a commutative triangle
	\begin{equation*}
		\begin{tikzcd}[sep=1.5em]
			& F(U) \arrow[dl] \arrow[dr, "\sim"{sloped, yshift=-0.2em}] & \\
			\displaystyle\lim_{V \in \Iup(\Uu)^{\op}} F(V) \arrow[rr, "\sim"{yshift=-0.2em}] & & \displaystyle\lim_{V \in \Iup(\Uu_0)^{\op}} F(V) \comma
		\end{tikzcd}
	\end{equation*}
	where the right-hand diagonal morphism is an equivalence by \enumref{lem:sheafLindelof}{2} and the horizontal morphism is an equivalence because the inclusion $ \Iup(\Uu_0)^{\op} \subset \Iup(\Uu)^{\op} $ is limit-cofinal.
\end{proof}

Now we provide a characterization of sheaves on a hereditarily Lindelöf topological space in terms of an excision property.
This characterization immediately implies \Cref{thm:excision}.

\begin{proposition}\label{prop:Lindelofexcision}
	Let $ T $ be a hereditarily Lindelöf topological space and $ C $ a presentable \category.
	A $ C $-valued presheaf $ F \in \PSh(T;C) $ on $ T $ is a sheaf if and only if $ F $ satisfies the following conditions:
	\begin{enumerate}[label=\stlabel{prop:Lindelofexcision}, ref=\arabic*]
		\item\label{prop:Lindelofexcision.1} The object $ F(\emptyset) $ is terminal in $ C $.

		\item\label{prop:Lindelofexcision.2} For all opens $ U, V \subset T $, the induced square 
		\begin{equation*}
			\begin{tikzcd}
				F(U \union V) \arrow[r] \arrow[d] & F(V) \arrow[d] & \\
				F(U) \arrow[r] & F(U \intersect V) 
			\end{tikzcd}
		\end{equation*}
		is a pullback square in $ C $.

		\item\label{prop:Lindelofexcision.3} For every $ \NN $-indexed sequence of open sets $ U_0 \subset U_1 \subset \cdots \subset T $, the induced morphism
		\begin{equation*}
			\fromto{F\paren{\textstyle\Union_{n \geq 0} U_n}}{\lim_{n \geq 0} F(U_n)}
		\end{equation*}
		is an equivalence in $ C $.
	\end{enumerate}
\end{proposition}

\begin{proof}
	First note that \enumref{prop:Lindelofexcision}{1} and \enumref{prop:Lindelofexcision}{2} are equivalent to saying that $ F $ satisfies descent with respect to finite covers.
	By \Cref{lem:sheafLindelof}, it suffices to show that $ F $ satisfies descent with respect to countable covers.

	Let $ V \subset T $ be open and $ \Uu = \{V_i\}_{i \in \NN} $ a countable open cover of $ V $.
	For each $ n \in \NN $, define
	\begin{equation*}
		U_n \colonequals \Union_{i=0}^{n} V_i \andeq \Uu_n \colonequals \{V_0,\ldots,V_n\} \period
	\end{equation*}
	Then $ \Uu_n $ is a finite open cover of $ U_n $ and we have inclusions $ U_n \subset U_{n+1} $ and $ \Uu_n \subset \Uu_{n+1} $.
	Note that the poset $ \Iup(\Uu) $ is the filtered union
	\begin{equation*}
		\Iup(\Uu) = \colim_{n \geq 0} \Iup(\Uu_n) \period
	\end{equation*}
	Since $ F $ satisfies descent with respect to finite covers, by \enumref{prop:Lindelofexcision}{3} we see that we have natural equivalences
	\begin{align*}
		F(V) &\equivalence \lim_{n \geq 0} F(U_n) \\ 
		&\equivalence \lim_{n \geq 0} \lim_{U \in \Iup(U_n)^{\op}} F(U) \\ 
		&\equivalent \lim_{U \in \Iup(\Uu)^{\op}} F(U) \period
	\end{align*}
	Hence $ F $ satisfies descent with respect to the countable cover $ \Uu $, as desired.
\end{proof}

\begin{proof}[Proof of \Cref{thm:excision}]
	Since manifolds are second-countable and open subsets of manifolds are manifolds, the claim is immediate from \Cref{prop:Lindelofexcision} and the definition of what it means to be a sheaf on $ \Mfld $ (\Cref{def:sheafonMan}).
\end{proof}

%-------------------------------------------------------------------%
%  Digression: relation to infinite dimensional manifolds           %
%-------------------------------------------------------------------%

\subsection{Digression: relation to infinite dimensional manifolds}\label{subsec:infdimMan}

We finish this chapter by describing a ``Yoneda embedding'' of infinite dimensional manifolds into sheaves of \textit{sets} on $ \Mfld $.

\begin{recollection}[(infinite dimensional manifolds)]
	There are two classes of possibly infinite dimensional manifolds that are commonly considered: \emph{Banach manifolds} and \emph{Fréchet manifolds} \cites[Chapter III, \S1]{MR0341518}[\S I.4]{MR656198}.  
	Banach spaces are examples of Fréchet spaces, and the category of Banach manifolds is a full subcategory of the
	category of Fréchet manifolds. 
\end{recollection}

One reason to consider Fréchet manifolds is that the (smooth) free loop space of a manifold naturally has the structure of a Fréchet manifold:

\begin{example}\label{ex:mappingFrechet}
	If $ M $ and $ N $ are manifolds, and $ M $ is compact, then the topological space $ \Cinf(M,N) $ of smooth maps $ \fromto{M}{N} $ has a natural Fréchet manifold structure.  
	See \cite[Chapter III, §1]{MR0341518}, in particular \cite[Chapter III, Theorem 1.11]{MR0341518}, for details.
\end{example}

\begin{notation}
	We write $ \Fre $ for the category of Fréchet manifolds.
	Note that $ \Mfld $ is the full subcategory of $ \Fre $ spanned by the finite dimensional Fréchet manifolds.
\end{notation}

\begin{construction}\label{construction:embedBanFre}
	Define a restricted Yoneda functor $ y_{\Fre} $ by
	\begin{align*}
		y_{\Fre} \colon \Fre &\to \PSh(\Mfld;\Set) \\
		F &\mapsto [M \mapsto \Map_{\Fre}(M,F)] \period
	\end{align*}
	Notice that since morphisms of Fréchet manifolds are defined locally, for each Fréchet manifold $ F $, the presheaf $ y_{\Fre}(F) $ is a sheaf.
	That is, $ y_{\Fre} $ factors through $ \Sh(\Mfld;\Set) $.
\end{construction}

\begin{theorem}[{(Hain \cite{MR539632}, Losik \cites{MR1213569}[Theorem 3.1.1]{MR1307261}[Theorem A.1.5]{MR2905777})}]
	The functor $ y_{\Fre} \colon \fromto{\Fre}{\Sh(\Mfld;\Set)} $ is fully faithful. 
\end{theorem}

The next result about infinite dimensional manifolds is that the embedding $ y_{\Fre} $ sends Fréchet manifold of smooth maps from a compact manifold to an arbitrary manifold (\Cref{ex:mappingFrechet}) to the internal-$ \Hom $ in $ \Sh(\Mfld;\Set) $.
In particular, free loop spaces are correctly represented in $ \Sh(\Mfld;\Set) $.
To state this result, let us first recall the internal-$ \Hom $ in sheaves on $ \Mfld $.

\begin{recollection}[(cartesian closedness)]
	Like any topos, the category $ \Sh(\Mfld;\Set) $ of sheaves of sets on $ \Mfld $ is \textit{cartesian closed}.
	In particular, $ \Sh(\Mfld;\Set) $ has an internal-$ \Hom $ defined by
	\begin{align*}
		\Hom_{\Sh(\Mfld;\Set)}(E,E') \colon \Mfldop &\to \Set \\
		M &\mapsto \Map_{\Sh(\Mfld;\Set)}(E \cross M, E') \period
	\end{align*}
\end{recollection}

\begin{theorem}[{(Waldorf \cite[Lemma A.1.7]{MR2905777})}]
	Let $ M $ and $ N $ be manifolds.
	If $ M $ is compact, then there is a natural isomorphism
	\begin{equation*}
		y_{\Fre}(\Cinf(M,N)) \isomorphic \Hom_{\Sh(\Mfld;\Set)}(M,N) \period
	\end{equation*}
\end{theorem}

We finish this section by explaining how a commonly used enlargement of the category of Fréchet manifolds fits into the category $ \Sh(\Mfld;\Set) $.

\begin{remark}[(diffeological spaces)]
	Souriau introduced \cite{MR607688} \emph{diffeological spaces} as generalization of manifolds to include infinite dimensional manifolds as well as manifold-like spaces appearing in mathematical physics.
	Diffeological spaces have been extensively studied by Iglesias-Zemmour and collaborators \cites{MR799609}{MR844156}{MR906897}{MR906897}{MR2346968}{MR2424313}{MR3025051}{MR2846342}{MR2592936}.

	To explain how diffeological spaces fit into sheaves on manifolds, write 
	\begin{equation*}
		\yo \colon \incto{\Euc}{\Sh(\Euc;\Set) \equivalent \Sh(\Mfld;\Set)}
	\end{equation*}
	for the Yoneda embedding.
	A \emph{diffeological space} is a sheaf $ E $ on $ \Euc $ such that for each $ n \geq 0 $, the natural map
	\begin{equation*}
		\begin{tikzcd}[sep=2.5em]
			E(\RR^n) \isomorphic \Map_{\Sh(\Euc;\Set)}(\yo(\RR^n),E) \arrow[r] & \Map_{\Set}(\yo(\RR^n)(*),E(*)) = \Map_{\Set}(\RR^n,E(*))
		\end{tikzcd}
	\end{equation*}
	is injective.
	This injectivity condition allows a diffeological space to be described as a set $ X $ equipped with a
	collection of ``plots''
	\begin{equation*}
	 	\Cinf(\RR^n,X) \subset \Map_{\Set}(\RR^n,X)
	\end{equation*}
	for each $ n \geq 0 $, subject to a collection of explicit conditions that are equivalent to saying that the assignment 
	\begin{equation*}
		\goesto{\RR^n}{\Cinf(\RR^n,X)}
	\end{equation*}
	is a sheaf on $ \Euc $.
	(To match up notation, $ X = E(*) $ and $ \Cinf(\RR^n,X) = E(\RR^n) $.)
\end{remark}

\newpage
%!TEX root = ../diffcoh.tex

%-------------------------------------------------------------------%
%-------------------------------------------------------------------%
%  ℝ-invariant sheaves                                              %
%-------------------------------------------------------------------%
%-------------------------------------------------------------------%

\section{\texorpdfstring{$ \RR $}{ℝ}-invariant sheaves}\label{sec:hisheaves}
\textit{by Peter Haine}

In this chapter, we investigate \textit{\RRinvariant} (or \textit{homotopy invariant}) sheaves on $ \Mfld $.
These are the sheaves that invert homotopy equivalences of manifolds.
The main result of this chapter is Dugger's observation that the global sections functor induces an equivalence from the subcategory $ \Sh(\Mfld;C) $ of \RRinvariant sheaves to $ C $ (\Cref{prop:Dugger}).
In the case where $ C = \Spc $, we show that the constant sheaf functor $ \Gammaupperstar \colon \fromto{\Spc}{\Sh(\Mfld;\Spc)} $ is given by the the assignment
\begin{equation*}
	X \mapsto [M \mapsto \Map_{\Spc}(\Piinf(M),X)] \comma
\end{equation*} 
where $ \Piinf(M) $ denotes the underlying homotopy type of the manifold $ M $.
More generally, the constant sheaf functor $ \Gammaupperstar \colon \fromto{C}{\Sh(\Mfld;C)} $ is given by the assignment
\begin{equation*}
	X \mapsto [M \mapsto X^{\Piinf(M)}] \comma
\end{equation*} 
where $ X^{\Piinf(M)} $ denotes the \textit{cotensor} of $ X \in C $ by $ \Piinf(M) \in \Spc $; see \Cref{rec:cotensoring}.
These results imply that there exists a chain of four (explicit) adjoints
\begin{equation*}
	\begin{tikzcd}[sep=4em]
		\Sh(\Mfld;C) \arrow[r, shift left=2ex, "\Gammalowersharp"] \arrow[r, shift right=0.75ex, "\Gammalowerstar"{description, xshift=0.5em}] & C \arrow[l, shift left=2ex, "\Gammauppersharp", hooked'] \arrow[l, shift right=0.75ex, "\Gammaupperstar"{description, xshift=-0.5em}, hooked']
	\end{tikzcd}
\end{equation*}
relating $ \Sh(\Mfld;C) $ and $ C $ \Cref{nul:Gammaadjunctions}.

Looking forward, in \Cref{sec:localization}, we give an explicit formula for $ \Gammalowersharp $ as a geometric realization.
In \Cref{sec:stable}, we use these adjoints and relations between them to construct a ``differential cohomology
diagram'' for sheaves on $ \Mfld $ with values in any presentable stable \category.

\Cref{sec:prelimsglobalsec} starts with some preliminary observations about the global sections and constant sheaf functors.
In \Cref{sec:RRinvarbasics}, we define \RRinvariant sheaves and explore some of their basic properties.
\Cref{sec:constantsheafdescript} is dedicated to proving that the global sections functor restricts to an equivalence on \RRinvariant sheaves.
In \cref{sec:conequences_of_Dugger}, we explore some immediate consequences of this result.

%-------------------------------------------------------------------%
%-------------------------------------------------------------------%
%  Preliminaries on global sections and constant sheaves            %
%-------------------------------------------------------------------%
%-------------------------------------------------------------------%

\subsection{Preliminaries on global sections and constant sheaves}\label{sec:prelimsglobalsec}

We begin by fixing some notation that we use throughout the rest of this text.

\begin{notation}
	Write $ \Gammalowerstar \colon \fromto{\PSh(\Mfld;C)}{C} $ for the \textit{global sections} functor, defined by
	\begin{equation*}
		\Gammalowerstar(E) \colonequals E(*) \period
	\end{equation*}
	Write $ \Gammainverse \colon \fromto{C}{\PSh(\Mfld;C)} $ for the \textit{constant presheaf} functor.
	That is, $ \Gammainverse $ is left adjoint to global sections functor $ \Gammalowerstar \colon \fromto{\PSh(\Mfld;C)}{C} $.

	Write $ \Gammaupperstar \colon \fromto{C}{\Sh(\Mfld;C)} $ for the \textit{constant sheaf} functor.
	Then $ \Gammaupperstar $ is the sheafification of $ \Gammainverse $.
	Moreover, $ \Gammaupperstar $ is left adjoint to the restriction $ \Gammalowerstar \colon \fromto{\Sh(\Mfld;C)}{C} $ of the global sections functor to sheaves.

	We use the same notations for the constant (pre)sheaf and global sections functors for the Euclidean site $ \Euc \subset \Mfld $ introduced in \cref{subsec:Cartsheaves}.
\end{notation}

\begin{observation}
	Note that we have a natural identification $ \Gammalowerstar \Gammainverse \equivalent \id{} $.
	Since $ \Gammalowerstar $ is right adjoint to $ \Gammainverse $, we conclude that $ \Gammainverse $ is fully faithful \cite[Lemma 3.3.1]{arXiv:2007.13089}.
\end{observation}

\noindent The global sections functor also has a right adjoint.

\begin{lemma}\label{lem:Gammauppershriekpsh}
	Let $ C $ be a presentable \category.
	Then the functor $ \Gammauppersharp \colon \fromto{C}{\PSh(\Mfld;C)} $ defined by the formula
	\begin{equation*}
		\Gammauppersharp(X)(M) \colonequals \prod_{m \in M} X 
	\end{equation*}
	is fully faithful and right adjoint to the global sections functor $ \Gammalowerstar \colon \fromto{\PSh(\Mfld;C)}{C} $.
	(Here the product is over the underlying set of the manifold $ M $.)
\end{lemma}

\begin{proof}
	We define the unit and counit of the adjunction.
	The unit $ \unit_{F} \colon \fromto{F}{\Gammauppersharp \Gammalowerstar(F)} $ is defined by the natural map
	\begin{equation*}
		F(M) \to \prod_{m \in M} F(\{m\}) \equivalence \Gammauppersharp \Gammalowerstar(F)(M)
	\end{equation*}
	induced by the inclusions $ \incto{\{m\}}{M} $ for all $ m \in M $.
	The counit $ \counit_X \colon \fromto{\Gammalowerstar \Gammauppersharp(X)}{X} $ is given by the natural identification $ \prod_{*} X \equivalent X $.
	The triangle identities are immediate from the definitions.

	To conclude, note that since the counit $ \counit $ is an equivalence, the functor $ \Gammauppersharp $ is fully faithful.
\end{proof}

Before recording the consequences of \Cref{lem:Gammauppershriekpsh} on the level of sheaves, we recall a basic fact from category theory.
For a proof see, for example, \cite[Chapter VII, \S4, Lemma 1]{MR1300636}.

\begin{lemma}\label{lem:ladjradj}
	Let $ \flowerstar \colon \fromto{A}{B} $ be a functor between \categories.
	Assume that $ \flowerstar $ admits a left adjoint $ \fupperstar $ and right adjoint $ \fuppersharp $.
	Then $ \fupperstar $ is fully faithful if and only if $ \fuppersharp $ is fully faithful.
\end{lemma}

\begin{corollary}\label{cor:GammauppershriekSh}
	Let $ C $ be a presentable \category.
	\begin{enumerate}[label=\stlabel{cor:GammauppershriekSh}, ref=\arabic*]
		\item\label{cor:GammauppershriekSh.1} The functor $ \Gammauppersharp $ factors through $ \Sh(\Mfld;C) $.

		\item\label{cor:GammauppershriekSh.2} The global sections functor $ \Gammalowerstar \colon \fromto{\Sh(\Mfld;C)}{C} $ is left adjoint to $ \Gammauppersharp \colon \fromto{C}{\Sh(\Mfld;C)} $.
		
		\item\label{cor:GammauppershriekSh.3} The constant sheaf functor $ \Gammaupperstar $ is fully faithful.
	\end{enumerate}
\end{corollary}

\begin{proof}
	For \enumref{cor:GammauppershriekSh}{1}, note that is immediate from \Cref{def:sheafonMan} that for each $ X \in C $, the presheaf $ \Gammauppersharp(X) $ is a sheaf on $ \Mfld $.
	\Cref{lem:Gammauppershriekpsh} and \enumref{cor:GammauppershriekSh}{1} immediately imply \enumref{cor:GammauppershriekSh}{2}.
	Finally, \enumref{cor:GammauppershriekSh}{3} is a consequence of \Cref{lem:Gammauppershriekpsh}, \enumref{cor:GammauppershriekSh}{2}, and the full faithfulness of $ \Gammauppersharp $.
\end{proof}

%-------------------------------------------------------------------%
%-------------------------------------------------------------------%
%  Basics of ℝ-invariant sheaves                                    %
%-------------------------------------------------------------------%
%-------------------------------------------------------------------%

\subsection{Basics of \texorpdfstring{$ \RR $}{ℝ}-invariant sheaves}\label{sec:RRinvarbasics}

We start by introducing an important subcategory of $ \Sh(\Mfld;C) $.

\begin{definition}[(\RRinvariant sheaves on $ \Mfld $)]\label{def:RR-invariance}
	Let $ C $ be a presentable \category.
	We say that a $ C $-valued presheaf
	\begin{equation*}
		F \colon \fromto{\Mfldop}{C}
	\end{equation*}
	is \textit{\RRinvariant}, \textit{homotopy-invariant}, or \textit{concordance-invariant} if for each manifold $ M $, the first projection $ \pr_M \colon \fromto{M \cross \RR}{M} $ induces an equivalence
	\begin{equation*}
		\prupperstar_M \colon \isomto{F(M)}{F(M \cross \RR)} \period
	\end{equation*}
	Write
	\begin{equation*}
		\Shhi(\Mfld;C) \subset \Sh(\Mfld;C) \andeq \PShhi(\Mfld;C) \subset \PSh(\Mfld;C)
	\end{equation*}
	for the full subcategories spanned by the \RRinvariant $ C $-valued sheaves and presheaves, respectively.
\end{definition}

We now give some reformulations of \RRinvariance.
The main one is that a presheaf is \RRinvariant if and only if it inverts all homotopy equivalences between manifolds.

\begin{notation}
	Let $ M $ be a manifold and $ t \in \RR $.
	We write $ i_{M,t} \colon \incto{M}{M \cross \RR} $ for the closed embedding defined by $ \goesto{x}{(x,t)} $.
\end{notation}

\begin{observation}\label{obs:iMcrossR}
	For each manifold $ M $ and $ t \in \RR $, the map $ i_{M \cross \RR,t} $ is given by the composite 
	\begin{equation*}
		\begin{tikzcd}[column sep=4.5em]
			M \cross \RR \arrow[r, hooked, "i_{M,t} \cross \id{\RR}"] & M \cross \RR \cross \RR \arrow[r, "\id{M} \cross \swap", "\sim"' ] & M \cross \RR \cross \RR \comma
		\end{tikzcd}
	\end{equation*}
	where $ \swap \colon \isomto{\RR \cross \RR}{\RR \cross \RR} $ is the map that swaps the two factors.
\end{observation}

\begin{proposition}\label{prop:differentRinvariance}
	Let $ C $ be \acategory and $ F \colon \fromto{\Mfldop}{C} $ a presheaf.
	The following are equivalent:
	\begin{enumerate}[label=\stlabel{prop:differentRinvariance}, ref=\arabic*]
		\item\label{prop:differentRinvariance.1} The presheaf $ F $ is \RRinvariant.

		\item\label{prop:differentRinvariance.2} For each homotopy equivalence of manifolds $ f \colon \fromto{M}{N} $, the map $ \fupperstar \colon \fromto{F(N)}{F(M)} $ is an equivalence in $ C $.

		\item\label{prop:differentRinvariance.3} For each manifold $ M $ and $ t \in \RR $, the induced map $ \iupperstar_{M,t} \colon \fromto{F(M \cross \RR)}{F(M)} $ is an equivalence.

		\item\label{prop:differentRinvariance.4} For each manifold $ M $, the induced maps $ \iupperstar_{M,0}, \iupperstar_{M,1} \colon \fromto{F(M \cross \RR)}{F(M)}  $ are equivalent.
	\end{enumerate}
\end{proposition}

\begin{proof}
	Since the embeddings
	\begin{equation*}
		i_{M,t} \colon \incto{M}{M \cross \RR}
	\end{equation*}
	are sections of the projection $ \pr_M \colon \fromto{M \cross \RR}{M} $, it is clear that \enumref{prop:differentRinvariance}{1} $ \Leftrightarrow $ \enumref{prop:differentRinvariance}{3} and \enumref{prop:differentRinvariance}{1} $ \Rightarrow $ \enumref{prop:differentRinvariance}{4}.
	It is also clear that \enumref{prop:differentRinvariance}{2} $ \Rightarrow $ \enumref{prop:differentRinvariance}{1}

	To see that \enumref{prop:differentRinvariance}{1} $ \Rightarrow $ \enumref{prop:differentRinvariance}{2}, let $ g \colon \fromto{N}{M} $ be a homotopy inverse to $ f $, so that there are homotopies of $ gf $ and $ fg $ with the respective identities fitting into commutative diagrams
 	\begin{equation*}
	 		\begin{tikzcd}[column sep=3em]
		 			M \arrow[d, "i_{M,0}"', hooked] \arrow[dr, equals] & \\ 
		 			M \cross [0,1] \arrow[r, "h_0" description] & M \\ 
		 			M \arrow[u, "i_{M,1}", hooked'] \arrow[ur, "gf"'] &
		 		\end{tikzcd} 
	 		\andeq
	 		\begin{tikzcd}[column sep=3em]
		 			N \arrow[d, "i_{N,0}"', hooked] \arrow[dr, equals] & \\ 
		 			N \cross [0,1] \arrow[r, "h_1" description] & N \\ 
		 			N \arrow[u, "i_{N,1}", hooked'] \arrow[ur, "fg"'] & \phantom{N} \period
		 		\end{tikzcd}
	 	\end{equation*}
 	Since the morphisms $ \iupperstar_{N,0} $ and $ \iupperstar_{M,0} $ are equivalences and the upper triangles commute, by the $ 2 $-of-$ 3 $ property both $ \hupperstar_0 $ and $ \hupperstar_1 $ are equivalences.
 	Since the morphisms $ \iupperstar_{N,1} $, $ \iupperstar_{M,1} $, $ \hupperstar_0 $, and $ \hupperstar_1 $ are equivalences and the lower triangles commute, we see that
 	\begin{equation*}
 		(fg)\upperstar \equivalent \gupperstar \fupperstar \andeq (gf)\upperstar \equivalent \fupperstar\gupperstar
 	\end{equation*}
 	are equivalences.
 	By the $ 2 $-of-$ 6 $ property we deduce that $ \fupperstar $ and $ \gupperstar  $ are equivalences.

	To complete the proof, we show that \enumref{prop:differentRinvariance}{4} $ \Rightarrow $ \enumref{prop:differentRinvariance}{1}.
	Assuming \enumref{prop:differentRinvariance}{4}, since $ i_{M,0} $ is a section or the projection $ \pr_M \colon \fromto{M \cross \RR}{M} $, it suffices to show that we have an equivalence
	\begin{equation*}
		\prupperstar_M \iupperstar_{M,0} \equivalent \id{F(M \cross \RR)} \period
	\end{equation*}
	To see this, write $ \mult \colon \fromto{\RR \cross \RR}{\RR} $ for the multiplication map, and notice that we have a commutative diagram in $ \Mfld $
	\begin{equation}\label{diag:multcom}
		\begin{tikzcd}[row sep=3em, column sep=4em]
			M \cross \RR \arrow[r, hooked, "i_{M,1} \cross \id{\RR}"] \arrow[dr, equals] & M \cross \RR \cross \RR \arrow[d, "\id{M} \cross \mult" description] & M \cross \RR \arrow[l, hooked', "i_{M,0} \cross \id{\RR}"'] \arrow[d, "\pr_M"] \\
			 & M \cross \RR & M \arrow[l, hooked', "i_{M,0}"]
		\end{tikzcd}
	\end{equation}
	Combining the assumption that $ \iupperstar_{M \cross \RR,0} \equivalent \iupperstar_{M \cross \RR,1} $ with \Cref{obs:iMcrossR} shows that
	\begin{equation}\label{eq:prodequiv}
		(i_{M,0} \cross \id{\RR})\upperstar \equivalent (i_{M,1} \cross \id{\RR})\upperstar \period
	\end{equation}
	\Cref{eq:prodequiv} and the commutativity of the diagram \eqref{diag:multcom} now show that
	\begin{align*}
		\prupperstar_M \iupperstar_{M,0} &\equivalent (i_{M,0} \cross \id{\RR})\upperstar \of (\id{M} \cross \mult)\upperstar \\
		&\equivalent (i_{M,1} \cross \id{\RR})\upperstar \of (\id{M} \cross \mult)\upperstar \\
		&\equivalent \id{F(M \cross \RR)} \comma
	\end{align*}
	as desired.
\end{proof}

\begin{remark}
	The reformulation of \RRinvariance given in \enumref{prop:differentRinvariance}{4} is due to Voevodsky \cite[Lemma 2.16]{MR2242284}.
\end{remark}

Equivalences of \RRinvariant sheaves can be checked on global sections:

\begin{lemma}\label{lem:checkequivforhionpt}
	Let $ C $ be a presentable \category.
	A morphism $ f \colon \fromto{E}{E'} $ in $ \Shhi(\Mfld;C) $ is an equivalence if and only if $ \Gammalowerstar(f) $ is an equivalence in $ C $.
\end{lemma}

\begin{proof}
	This follows from \Cref{lem:checkequivonRn} and the assumption that $ E $ and $ E' $ are \RRinvariant.
\end{proof}

%-------------------------------------------------------------------%
%  ℝ-invariant sheaves on the Euclidean site                        %
%-------------------------------------------------------------------%

\subsubsection{\texorpdfstring{$ \RR $}{ℝ}-invariant sheaves on the Euclidean site}

It is convenient to describe \RRinvariant (pre)sheaves in terms of the Euclidean site $ \Euc \subset \Mfld $.
We start with the following reformulation of \RRinvariance:

\begin{lemma}\label{cor:differentRinvariance}
	Let $ C $ be a presentable \category and $ F \in \Sh(\Mfld;C) $ a sheaf.
	The following are equivalent:
	\begin{enumerate}[label=\stlabel{cor:differentRinvariance}, ref=\arabic*]
		\item\label{cor:differentRinvariance.1} The sheaf $ F $ is \RRinvariant.

		\item\label{cor:differentRinvariance.2} The functor $ \restrict{F}{\Eucop} \colon \fromto{\Eucop}{C} $ carries every morphism to an equivalence.
	\end{enumerate}
\end{lemma}

\begin{proof}
	Since Euclidean spaces are contractible, every map between Euclidean spaces is a homotopy equivalence.
	Thus applying \Cref{prop:differentRinvariance} shows that \enumref{cor:differentRinvariance}{1} $ \Rightarrow $ \enumref{cor:differentRinvariance}{2}.

	To see that \enumref{cor:differentRinvariance}{2} $ \Rightarrow $ \enumref{cor:differentRinvariance}{1}, assume that $ \restrict{F}{\Eucop} \colon \fromto{\Eucop}{C} $ carries every morphism to an equivalence.
	Let $ M $ be a manifold.
	To show that $ \prupperstar_M $ is an equivalence, fix a \textit{good} cover $ \Uu $ if $ M $.
	Then the collection $ \{U \cross \RR\}_{U \in \Uu} $ is a good cover of $ M \cross \RR $.
	The projection maps $ \pr_U \colon \fromto{U \cross \RR}{\RR} $ induce a commutative square
	\begin{equation}\label{eq:Cechequivs_again}
		\begin{tikzcd}[sep=3em]
			F(M) \arrow[d, "\prupperstar_M"'] \arrow[r, "\sim"{yshift=-0.2em}] & \displaystyle\lim_{U \in \Iup(\Uu)^{\op}} F(U) \arrow[d, "\lim\limits_{U \in \Iup(\Uu)^{\op}} \prupperstar_{U}"] \\
			F(M \cross \RR) \arrow[r, "\sim"{yshift=-0.2em}] & \displaystyle\lim_{U \in \Iup(\Uu)^{\op}} F(U \cross \RR) \comma
		\end{tikzcd}
	\end{equation}
	where the horizontal morphisms are equivalences because $ F $ is a sheaf.
	Since $ \restrict{F}{\Eucop} $ carries every morphism to an equivalence and the cover $ \Uu $ consists of Euclidean spaces, we see that the right-hand vertical morphism in \eqref{eq:Cechequivs_again} is an equivalence.
	By the $ 2 $-of-$ 3 $ property, the left-hand vertical morphism in \eqref{eq:Cechequivs_again} is an equivalence, as desired.
\end{proof}

\Cref{cor:differentRinvariance} motivates the following variant of \Cref{def:RR-invariance}:

\begin{definition}[(\RRinvariant sheaves on $ \Euc $)]
	Let $ C $ be a presentable \category.
	We say that a $ C $-valued presheaf
	\begin{equation*}
		F \colon \fromto{\Eucop}{C}
	\end{equation*}
	is \textit{\RRinvariant} if $ F $ carries every morphism in $ \Eucop $ to an equivalence.
	Write
	\begin{equation*}
		\Shhi(\Euc;C) \subset \Sh(\Euc;C) \andeq \PShhi(\Euc;C) \subset \PSh(\Euc;C)
	\end{equation*}
	for the full subcategories spanned by the \RRinvariant $ C $-valued sheaves and presheaves, respectively.
\end{definition}

\begin{observation}
	Let $ F $ be an \RRinvariant presheaf on $ \Mfld $.
	Then the restriction $ \restrict{F}{\Eucop} $ is an \RRinvariant presheaf on $ \Euc $.
	By \Cref{cor:differentRinvariance}, if $ G $ is an \RRinvariant \textit{sheaf} on $ \Euc $, then the right Kan extension $ \jlowerstar(G) \colon \fromto{\Mfldop}{C} $ is an \RRinvariant sheaf on $ \Mfld $.
	Hence the inverse equivalences of \categories
	\begin{equation*}
		\adjto{\restrict{(-)}{\Eucop}}{\Sh(\Mfld;C)}{\Sh(\Euc;C)}{\jlowerstar}
	\end{equation*}
	of \Cref{lem:hypersheavesonCart} restrict to inverse equivalences 
	\begin{equation*}
		\adjto{\restrict{(-)}{\Eucop}}{\Shhi(\Mfld;C)}{\Shhi(\Euc;C)}{\jlowerstar} \period
	\end{equation*}
\end{observation}

Similarly to \Cref{lem:checkequivforhionpt}, a morphism of $ \RR $-invariant presheaves on $ \Euc $ is an equivalence if and only if it induces an equivalence on global sections.

\begin{lemma}\label{lem:checkequivforhionpt_Euc}
	Let $ C $ be a presentable \category.
	A morphism $ f \colon \fromto{E}{E'} $ in $ \PShhi(\Euc;C) $ is an equivalence if and only if $ \Gammalowerstar(f) $ is an equivalence in $ C $.
\end{lemma}

\begin{proof}
	For the nontrivial direction, assume that $ \Gammalowerstar(f) $ is an equivalence, and fix an integer $ n \geq 0 $; we need to show that $ f(\RR^n) $ is an equivalence.
	Consider the commutative square
	\begin{equation}\label{eq:quivs_on_RRn}
		\begin{tikzcd}
			E(\ast) \arrow[d, "f(\ast)"'] \arrow[r] & E(\RR^n) \arrow[d, "f(\RR^n)"] \\
			E'(\ast) \arrow[r] & E'(\RR^n) \comma
		\end{tikzcd}
	\end{equation}
	where the horizontal morphisms are induced by the unique morphism $ \fromto{\RR^n}{\ast} $.
	Since $ E $ and $ E' $ carry all morphisms in $ \Eucop $ to equivalences, the horizontal morphisms in \eqref{eq:quivs_on_RRn} are equivalences.
	By assumption $ f(*) = \Gammalowerstar(f) $ is an equivalence.
	Hence the $ 2 $-of-$ 3 $ property shows that $ f(\RR^n) $ is an equivalence, as desired. 
\end{proof}

\begin{remark}
	Note that in \Cref{lem:checkequivforhionpt_Euc} we only require that $ E $ and $ E' $ be \textit{pre}sheaves.
\end{remark}

%-------------------------------------------------------------------%
%-------------------------------------------------------------------%
%  The constant sheaf functor                                       %
%-------------------------------------------------------------------%
%-------------------------------------------------------------------%

\subsection{The constant sheaf functor}\label{sec:constantsheafdescript}

The goal of this section is to prove the following result, originally sketched for sheaves of spaces by Dugger \cites[Theorem 3.4.3]{Duggersheaves}[Proposition 8.3]{MR1870515} and Morel--Voevodsky \cite[Proposition 3.3.3]{MR1813224}.
See the work of Bunk \cite{arXiv:2007.06039} for another model category-theoretic argument.

\begin{proposition}\label{prop:Dugger}
	Let $ C $ be a presentable \category.
	Then:
	\begin{enumerate}[label=\stlabel{prop:Dugger}, ref=\arabic*]
		\item The constant sheaf functor $ \Gammaupperstar \colon \fromto{C}{\Sh(\Mfld;C)} $ factors through $ \Shhi(\Mfld;C) $.

		\item The global sections functor
		\begin{equation*}
			\Gammalowerstar \colon \fromto{\Shhi(\Mfld;C)}{C}
		\end{equation*}
		is an equivalence with inverse given by $ \Gammaupperstar $.

		\item A sheaf $ F $ on $ \Mfld $ is \RRinvariant if and only if $ F $ is constant.

		\item The constant sheaf functor $ \Gammaupperstar \colon \fromto{C}{\Sh(\Mfld;C)} $ admits a left adjoint.
	\end{enumerate}
\end{proposition}

\begin{remark}
	An analogue of \Cref{prop:Dugger} holds where the category of manifolds is replaced by the category of smooth complex analytic spaces, and $ \RR $ is replaced by the open unit disk in $ \CC $; see \cite[Remarque 1.9]{MR2602027}.
	Similarly, there are many variants of this result where $ \Mfld $ is replaced by any reasonable category of locally contractible spaces.
\end{remark}

\begin{remark}
	See \cite[\S2]{arXiv:2010.06473} for related results about constant (hyper)sheaves on locally weakly contractible topological spaces.
\end{remark}

%-------------------------------------------------------------------%
%  Background on cotensors                                          %
%-------------------------------------------------------------------%

\subsubsection{Background on cotensors}

In order to prove \Cref{prop:Dugger}, we give a concrete description of the constant sheaf functor. 
To do this, we first recall the natural cotensoring of a presentable \category over $ \Spc $.

\begin{recollection}[(cotensoring over $ \Spc $)]\label{rec:cotensoring}
	Every presentable \category $ C $ is naturally \textit{cotensored over} the \category $ \Spc $ of spaces \HTT{Remark}{5.5.2.6}.
	That is, there is a functor 
	\begin{align*}
		(-)^{(-)} \colon \Spc^{\op} \cross C &\to C \\
		(K,X) &\mapsto X^K \comma
	\end{align*}
	along with natural equivalences
	\begin{equation*}
		\Map_{C}(X',X^K) \equivalent \Map_{\Spc}(K,\Map_C(X',X)) \period
	\end{equation*}
\end{recollection}

\begin{example}
	If $ C = \Spc $ is the \category of spaces, then the cotensoring is given by
	\begin{equation*}
		X^K \colonequals \Map_{\Spc}(K,X) \period
	\end{equation*}
\end{example}

\begin{example}\label{ex:Sptcotensor}
	If $ C = \Sp $ is the \category of spectra, then the cotensoring is given by
	\begin{equation*}
		X^K \colonequals \Hom_{\Sp}(\Sigma_+^\infty K,X) \comma
	\end{equation*}
	where $ \Hom_{\Sp} $ denotes the mapping \textit{spectrum} in $ \Sp $.
\end{example}

\begin{example}
	If $ R $ is a ring and let $ C = \D(R) $ be the derived \category of $ R $, then the cotensoring is given by
	\begin{equation*}
		A_{*}^K \colonequals \RHom_{R}(\Cup_*(K;R),A_{*}) \period
	\end{equation*}
	Here $ \Cup_*(K;R) $ is the complex of singular chains on $ K $, and $ \RHom_{R} $ is the derived $ \Hom $
	functor of chain complexes of $ R $-modules.
	If $ M $ is an ordinary $ R $-module regarded as an object of $ \D(R) $ concentrated in degree $ 0 $, then the complex
	\begin{equation*}
		\RHom_{R}(\Cup_*(K;R),M)
	\end{equation*}
	is the complex $ \Cup^{-*}(K;M) $ of \textit{singular cochains} on $ K $ with coefficients in $ M $.
\end{example}

%-------------------------------------------------------------------%
%  Description of the constant sheaf functor                        %
%-------------------------------------------------------------------%

\subsubsection{Description of the constant sheaf functor}

We now give an explicit formula for the constant sheaf functor.

\begin{notation}
	Recall that we write $ \Top $ for the category of topological spaces (\Cref{ntn:Top}).
	Write $ \Piinf \colon \fromto{\Top}{\Spc} $ for the functor sending a topological space $ T $ to the underlying homotopy type  $ \Piinf(T) $ of $ T $.
\end{notation}

\begin{construction}
	Write $ \Piinf \colon \fromto{\PSh(\Mfld;\Spc)}{\Spc} $ for the left Kan extension of the functor 
	\begin{equation*}
		\begin{tikzcd}[sep=2.75em]
			\Mfld \arrow[r, "\forget"] & \Top \arrow[r, "\Piinf"] & \Spc 
		\end{tikzcd}
	\end{equation*}
	along the Yoneda embedding $ \incto{\Mfld}{\PSh(\Mfld;\Spc)} $.
	By the universal property of the \category of presheaves, the functor $ \Piinf \colon \fromto{\PSh(\Mfld;\Spc)}{\Spc} $ is a left adjoint with right adjoint $ \sm \colon \fromto{\Spc}{\PSh(\Mfld;\Spc)} $ given by the formula
	\begin{equation*}
		X \mapsto \big[M \mapsto \Map_{\Spc}(\Piinf(M),X) \big] \period 
	\end{equation*}
	By the van Kampen Theorem \HAa{Proposition}{A.3.2}, the functor $ \sm \colon \fromto{\Spc}{\PSh(\Mfld;\Spc)} $ factors through $ \Sh(\Mfld;\Spc) $.
	We use the same notation for the resulting adjunction
	\begin{equation*}
		\adjto{\Piinf}{\Sh(\Mfld;\Spc)}{\Spc}{\sm} \period
	\end{equation*}

	Given a presentable \category $ C $, we write $ \Piinf \colon \fromto{\Sh(\Mfld;C)}{C} $ for the tensor product
	\begin{equation*} 
		\begin{tikzcd}[sep=5em]
			\Sh(\Mfld;C) \equivalent \Sh(\Mfld;\Spc) \tensor C \arrow[r, "\Piinf \tensor \id{C}"] & \Spc \tensor C \equivalent C \period
		\end{tikzcd}
	\end{equation*}
	We write
	\begin{equation*}
		\sm \colon \fromto{C}{\Sh(\Mfld;C)}
	\end{equation*}
	for the right adjoint of $ \Piinf $.
	Concretely, $ \sm $ is defined by sending $ X \in C $ to the sheaf
	\begin{equation*} 
		M \mapsto X^{\Piinf(M)} \period 
	\end{equation*}
\end{construction}

\begin{observation}\label{obs:sm_RR-invariant}
	For each $ X \in C $, the sheaf $ \sm(X) $ is \RRinvariant. 
\end{observation}

\begin{construction}[(comparison natural transformation)]
	Recall that we write 
	\begin{equation*}
		\Gammainverse \colon \fromto{C}{\PSh(\Mfld;C)}
	\end{equation*}
	for the constant presheaf functor.
	Define a natural transformation $ \alpha \colon \fromto{\Gammainverse}{\sm} $ as follows. 
	For each $ X \in C $, the component of $ \alpha(X) $ at a manifold $ M \in \Mfld $ is the map 
	\begin{equation*}
		\begin{tikzcd}
			\Gammainverse(X)(M) = X \arrow[r] & X^{\Piinf(M)} = \sm(X)(M) 
		\end{tikzcd}
	\end{equation*}
	induced by the unique map $ \fromto{\Piinf(M)}{\ast} $.
\end{construction}

\begin{proposition}[(formula for the constant sheaf functor)]\label{lem:constantishi}
	Let $ C $ be a presentable \category.
	Then: 
	\begin{enumerate}[label=\stlabel{lem:constantishi}, ref=\arabic*]
		\item\label{lem:constantishi.1} The natural transformation $ \alpha \colon \fromto{\Gammainverse}{\sm} $ is an equivalence when restricted to the subcategory $ \Eucop \subset \Mfldop $.
		That is, the composite
		\begin{equation*}
			\begin{tikzcd}[sep=3.5em]
				C \arrow[r, "\sm"] & \Sh(\Mfld;C) \arrow[r, "\restrict{(-)}{\Eucop}", "\sim"'] & \Sh(\Euc;C) 
			\end{tikzcd}
		\end{equation*}
		is naturally identified with the constant \emph{presheaf} functor.

		\item\label{lem:constantishi.2} Every constant $ C $-valued presheaf on $ \Euc $ is a sheaf.

		\item\label{lem:constantishi.3} The natural transformation $ \alpha \colon \fromto{\Gammainverse}{\sm} $ exhibits $ \sm $ as the sheafification of $ \Gammainverse $.
		That is, there is a natural equivalence $ \sm \equivalent \Gammaupperstar $.

		\item\label{lem:constantishi.4} The functor $ \Piinf \colon \fromto{\Sh(\Mfld;C)}{C} $ is left adjoint to $ \Gammaupperstar $.

		\item\label{lem:constantishi.5} The functor $ \Gammaupperstar \colon \incto{C}{\Sh(\Mfld;C)} $ factors through the subcategory
		\begin{equation*}
			\Shhi(\Mfld;C) \subset \Sh(\Mfld;C) \period
		\end{equation*}
		That is, every constant sheaf on $ \Mfld $ is \RRinvariant.
	\end{enumerate}
\end{proposition}

\begin{proof}
	To see \enumref{lem:constantishi}{1}, note that since each Euclidean space $ \RR^n $ is contractible, the unique map $ \fromto{\Piinf(\RR^n)}{\ast} $ is an equivalence. Hence the claim follows from the definition of $ \alpha $.

	For \enumref{lem:constantishi}{2}, combine \enumref{lem:constantishi}{1} with the fact that for each $ X \in C $, the presheaf
	\begin{equation*}
		\sm(X) \colon \fromto{\Mfldop}{C}
	\end{equation*}
	is a sheaf.
	Note that \enumref{lem:constantishi}{3} follows from \enumref{lem:constantishi}{1} and \Cref{lem:hypersheavesonCart}.
	Statement \enumref{lem:constantishi}{4} is immediate from \enumref{lem:constantishi}{3} and the definition of $ \sm $.
	Finally, \enumref{lem:constantishi}{5} is immediate from \enumref{lem:constantishi}{3} and \Cref{obs:sm_RR-invariant}.
\end{proof}

\Cref{prop:Dugger} now follows with the facts that $ \Gammaupperstar $ is fully faithful and $ \Gammalowerstar $ is conservative when restricted to the \RRinvariant sheaves (\Cref{lem:checkequivforhionpt}), combined with the following basic lemma from category theory.

\begin{lemma}\label{lem:conservativeradjffladj}
	Let $ \adjto{\fupperstar}{A}{B}{\flowerstar} $ be an adjunction between \categories, and assume that the left adjoint $ \fupperstar $ is fully faithful.
	Then $ \fupperstar $ is an equivalence if and only if $ \flowerstar $ is conservative.
\end{lemma}

\begin{proof}
	If $ \fupperstar $ is an equivalence, then $ \flowerstar $ is also an equivalence, hence conservative.

	On the other hand, assume that $ \flowerstar $ is conservative. 
	Since the left adjoint $ \fupperstar $ is fully faithful, the unit $ \fromto{\id{A}}{\flowerstar\fupperstar} $ is an equivalence.
	Hence $ \fupperstar $ is an equivalence if and only if for each object $ X \in B $, the counit $ \counit_{X} \colon \fromto{\fupperstar\flowerstar(X)}{X} $ is an equivalence.
	Since $ \flowerstar $ is conservative, the counit $ \counit_{X} $ is an equivalence if and only if
	\begin{equation*}
		\flowerstar(\counit_{X}) \colon \fromto{\flowerstar\fupperstar\flowerstar(X)}{\flowerstar(X)}
	\end{equation*}
	is an equivalence.
	The claim now follows from the fact that the unit $ \fromto{\id{A}}{\flowerstar\fupperstar} $ is an equivalence and the triangle identity.
\end{proof}

\begin{proof}[Proof of \Cref{prop:Dugger}]
	Since $ \Gammaupperstar \colon \incto{C}{\Shhi(\Mfld;C)} $ is fully faithful and
	\begin{equation*}
		\Gammalowerstar \colon \fromto{\Shhi(\Mfld;C)}{C}
	\end{equation*}
	is conservative (\Cref{lem:checkequivforhionpt}), we conclude by \Cref{lem:conservativeradjffladj}.
\end{proof}

\begin{corollary}\label{cor:RR-invariant_on_Euc}
	Let $ C $ be a presentable \category.
	\begin{enumerate}[label=\stlabel{cor:RR-invariant_on_Euc}, ref=\arabic*]
		% \item The constant presheaf functor $ \Gammainverse \colon \fromto{C}{\PSh(\Euc;C)} $ factors through $ \Shhi(\Euc;C) $.

		\item\label{cor:RR-invariant_on_Euc.1} The global sections functor $ \fromto{\PShhi(\Euc;C)}{C} $ is an equivalence with inverse given by the constant presheaf functor $ \Gammainverse \colon \fromto{C}{\PShhi(\Euc;C)} $.

		\item\label{cor:RR-invariant_on_Euc.2} The inclusion $ \Shhi(\Euc;C) \subset \PShhi(\Euc;C) $ is an equality.
		That is, every \RRinvariant presheaf on $ \Euc $ is automatically a sheaf.
	\end{enumerate}
\end{corollary}

\begin{proof}
	For \enumref{cor:RR-invariant_on_Euc}{1}, note that \Cref{lem:constantishi} shows that the the constant presheaf functor
	\begin{equation*}
		\Gammainverse \colon \fromto{C}{\PSh(\Euc;C)}
	\end{equation*}
	is fully faithful and factors through $ \PShhi(\Euc;C) $.
	Since the global sections functor
	\begin{equation*}
		\Gammalowerstar \colon \fromto{\PShhi(\Euc;C)}{C}
	\end{equation*}
	is conservative (\Cref{lem:checkequivforhionpt_Euc}), the claim follows from \Cref{lem:conservativeradjffladj}.

	Note that \enumref{cor:RR-invariant_on_Euc}{2} follows from \enumref{cor:RR-invariant_on_Euc}{1} and the fact that, again by \Cref{lem:constantishi}, the functor $ \Gammainverse \colon \fromto{C}{\PSh(\Euc;C)} $ factors through the subcategory $ \Shhi(\Euc;C) $.
\end{proof}

%-------------------------------------------------------------------%
%-------------------------------------------------------------------%
%  Consequences of \Cref{prop:Dugger}                               %
%-------------------------------------------------------------------%
%-------------------------------------------------------------------%

\subsection{Consequences of \texorpdfstring{\Cref{prop:Dugger}}{Proposition \ref*{prop:Dugger}}}\label{sec:conequences_of_Dugger}

%-------------------------------------------------------------------%
%  The four adjoints                                                %
%-------------------------------------------------------------------%

\subsubsection{The four adjoints}

\Cref{prop:Dugger} gives us a chain of four adjoints relating $ \Sh(\Mfld;C) $ and $ C $.
Let us start by giving a new name to the left adjoint $ \Piinf \colon \fromto{\Sh(\Mfld;C)}{C} $ that is consistent with all of the other functors relating $ \Sh(\Mfld;C) $ and $ C $.

\begin{definition}
	Let $ C $ be a presentable \category.
	We write
	\begin{equation*}
		\Gammalowersharp \colon \fromto{\Sh(\Mfld;C)}{C}
	\end{equation*}
	for the left adjoint to $ \Gammaupperstar \colon \fromto{C}{\Sh(\Mfld;C)} $.
\end{definition}

\begin{nul}\label{nul:Gammaadjunctions}
	Combining \Cref{cor:GammauppershriekSh,prop:Dugger}, we have a chain of four adjoints
	\begin{equation}\label{diag:4adjoints}
		\begin{tikzcd}[sep=4em]
			\Sh(\Mfld;C) \arrow[r, shift left=2ex, "\Gammalowersharp"] \arrow[r, shift right=0.75ex, "\Gammalowerstar"{description, xshift=0.5em}] & C \comma \arrow[l, shift left=2ex, "\Gammauppersharp", hooked'] \arrow[l, shift right=0.75ex, "\Gammaupperstar"{description, xshift=-0.5em}, hooked']
		\end{tikzcd}
	\end{equation}
	where functors lie above their right adjoints.
	Moreover, $ \Gammauppersharp $ and $ \Gammaupperstar $ are fully faithful.
\end{nul}

\begin{observation}
	The composite
	\begin{equation*}
		\Gammaupperstar \Gammalowersharp \colon \fromto{\Sh(\Mfld;C)}{\Shhi(\Mfld;C)}
	\end{equation*}
	is left adjoint to the inclusion $ \Shhi(\Mfld;C) \subset \Sh(\Mfld;C) $. 
	Similarly, the composite
	\begin{equation*}
		\Gammaupperstar \Gammalowerstar \colon \fromto{\Sh(\Mfld;C)}{\Shhi(\Mfld;C)}
	\end{equation*}
	is right adjoint to the inclusion $ \Shhi(\Mfld;C) \subset \Sh(\Mfld;C) $. 
\end{observation}

\begin{definition}[(homotopification)]\label{homotopification_definition}
	Let $ C $ be a presentable \category.
	We write
	\begin{equation*}
		\Lhi \colonequals \Gammaupperstar \Gammalowersharp \andeq \Rhi \colonequals \Gammaupperstar \Gammalowerstar 
	\end{equation*}
	for the left and right adjoint to the inclusion $ \Shhi(\Mfld;C) \subset \Sh(\Mfld;C) $, respectively. 
	We call $ \Lhi $ the \emph{homotopification} functor.
\end{definition}

\begin{observation}[(formulas for $ \Lhi $ and $ \Rhi $)]\label{obs:formulas_for_Lhi_and_Rhi}
	\Cref{lem:constantishi,cor:formula_for_Gammalowersharp} show that $ \Lhi $ and $ \Rhi $ is given by the formulas
	\begin{equation*}
		\Lhi(E)(M) \equivalent \Gammalowersharp(E)^{\Piinf(M)} \andeq \Rhi(E)(M) \equivalent E(\ast)^{\Piinf(M)} \period
	\end{equation*}
	Also note that we have identifications
	\begin{equation*}
		\Gammalowerstar\Lhi \equivalent \Gammalowersharp \andeq \Gammalowerstar\Rhi \equivalent \Gammalowerstar \period
	\end{equation*}
	In particular, when $ C = \Spc $, the functors $ \Lhi $ and $ \Rhi $ are given by the formulas
	\begin{equation*}
		\Lhi(E)(M) \equivalent \Map_{\Spc}(\Piinf(M),\Gammalowersharp(E)) \andeq \Rhi(E)(M) \equivalent \Map_{\Spc}(\Piinf(M),E(\ast)) \period
	\end{equation*}
\end{observation}

\begin{remark}[(cohesion)]
	Much of the structure of sheaves on $ \Mfld $ that we are interested in for studying differential cohomology (particularly \Cref{sec:stable}) only depends on the existence of the chain of four adjoints \eqref{diag:4adjoints}.
	In the case where $ C = \Spc $, the existence of these extra adjoints for the global sections geometric morphism (along with the condition that the extreme left adjoint $ \Gammalowersharp $ preserve finite products; see \Cref{cor:formula_for_Gammalowersharp}) is what Schreiber calls a \textit{cohesive \topos} \cite[Definition 3.4.1]{Schreiber:cohesive}.
	The primary examples of cohesive \topoi are \textit{global spaces} \cite{Rezk:global}, \textit{orbispaces} \cite[Chapter 3]{Ell-III}, and variants of sheaves on $ \Mfld $.
	Cohesive \topoi are a very general setting in which one can talk about a generalized form of ``differential
	cohomology''.

	Many of the ideas about cohesive \topoi go back to work of Lawvere \cites{MR1257681}{MR2125786}{MR2369017}[\S C.1]{MR1965482} as well as Simpson--Teleman \cite{SimpsonTeleman:deRham}.
\end{remark}

%-------------------------------------------------------------------%
%  ℝ-localization                                                   %
%-------------------------------------------------------------------%

\subsubsection{\texorpdfstring{$ \RR $}{ℝ}-localization}

We now observing that the inclusion of \RRinvariant (pre)sheaves into all (pre)sheaves admits a left adjoint.

\begin{observation}\label{obs:LRRandLhi}
	Notice that the full subcategory $ \PShhi(\Mfld;C) \subset \PSh(\Mfld;C) $ is closed under both limits and colimits.
	Hence $ \PShhi(\Mfld;C) $ is presentable and by the Adjoint Functor Theorem, the inclusion
	\begin{equation*}
		\PShhi(\Mfld;C) \subset \PSh(\Mfld;C)
	\end{equation*}
	admits both a left and a right adjoint.
	We write $ \LRR \colon \fromto{\PSh(\Mfld;C)}{\PShhi(\Mfld;C)} $ for the left adjoint to the inclusion.
\end{observation}
	
\begin{definition}[(\RRlocalization)]\label{def:RR-localization}
	Let $ C $ be a presentable \category.
	We refer to the left adjoint
	\begin{equation*}
		\LRR \colon \fromto{\PSh(\Mfld;C)}{\PShhi(\Mfld;C)} 
	\end{equation*}
	as the \emph{\RRlocalization} functor.
\end{definition}

\begin{nul}
	\Cref{sec:localization} is dedicated to providing an explicit formula for the functor $ \LRR $.
\end{nul}

It is not hard to see that the homotopification functor $ \Lhi \colon \fromto{\Sh(\Mfld;C)}{\Shhi(\Mfld;C)} $ is obtained by sheafifying $ \LRR $:

\begin{corollary}\label{cor:sheafification_preserves_RR-invariance}
	Let $ C $ be a presentable \category.
	If $ F \in \PSh(\Mfld;C) $ is \RRinvariant, then the sheafification $ \SMan(F) $ is \RRinvariant and the unit $ \fromto{F}{\SMan(F)} $ is an equivalence when restricted to $ \Eucop \subset \Mfldop $.
\end{corollary}

\begin{proof}
	Since the presheaf $ \restrict{F}{\Eucop} $ is \RRinvariant and a sheaf \enumref{cor:RR-invariant_on_Euc}{2}, the claim follows from \enumref{lem:hypersheavesonCart}{4}.
\end{proof}

\begin{corollary}\label{cor:hilowershriekasacomposite}
	Let $ C $ be a presentable \category.
	Then:
	\begin{enumerate}[label=\stlabel{cor:hilowershriekasacomposite}, ref=\arabic*]
		\item\label{cor:hilowershriekasacomposite.1} The composite
		\begin{equation*}
			\SMan \LRR \colon \fromto{\Sh(\Mfld;C)}{\Sh(\Mfld;C)}
		\end{equation*}
		factors through $ \Shhi(\Mfld;C) $ and is left adjoint to the inclusion $ \incto{\Shhi(\Mfld;C)}{\Sh(\Mfld;C)} $.
		That is, there is a natural equivalence
		\begin{equation*}
			\Lhi \equivalent \SMan\LRR \period
		\end{equation*}

		\item\label{cor:hilowershriekasacomposite.2} For each $ F \in \Sh(\Mfld;C) $, the natural transformation
		\begin{equation*}
			\fromto{\restrict{\LRR(F)}{\Eucop}}{\restrict{\Lhi(F)}{\Eucop}}
		\end{equation*}
		is an equivalence.
	\end{enumerate}
\end{corollary}

\begin{proof}
	For \enumref{cor:hilowershriekasacomposite}{1}, first note that \Cref{cor:sheafification_preserves_RR-invariance} immediately implies that $ \SMan \LRR $ factors through $ \Shhi(\Mfld;C) $.
	To see that $ \SMan \LRR $ is left adjoint to the inclusion, let $ F,G \in \Sh(\Mfld;C) $, and assume that $ G $ is \RRinvariant.
	Using the fact that $ \LRR \colon \fromto{\PSh(\Mfld;C)}{\PShhi(\Mfld;C)} $ is left adjoint to the inclusion and that $ F $ is a sheaf, we compute
	\begin{align*}
		\Map_{\Shhi(\Mfld;C)}(\SMan \LRR(F),G) &\equivalent \Map_{\PShhi(\Mfld;C)}(\LRR(F),G) \\
		&\equivalent \Map_{\Sh(\Mfld;C)}(F,G)
	\end{align*}

	Item \enumref{cor:hilowershriekasacomposite}{2} follows from \Cref{cor:sheafification_preserves_RR-invariance} and \enumref{cor:hilowershriekasacomposite}{1}.
\end{proof}

We finish this section with some remarks on the difference between $ \LRR $ and $ \Lhi $ and the notations we have chosen.

\begin{remark}[{($ \LRR $ vs. $ \Lhi $)}]\label{remark:BE-BdB-P}
	For a general presentable \category $ C $ and $ C $-valued sheaf $ E $ on $ \Mfld $, the presheaf $ \LRR(E) $ need not be a sheaf.
	Hence $ \Lhi $ is \textit{not} given by simply restricting $ \LRR $ to sheaves.
	However, the main result of the work of Berwick-Evans--Boavida de Brito--Pavlov \cite{BEBdBP:Classifying} shows that when $ C = \Spc $, the functor $ \LRR $ \textit{does} preserve sheaves.
	That is, for each sheaf $ E \in \Sh(\Mfld;\Spc) $, the natural morphism $ \fromto{\LRR(E)}{\Lhi(E)} $ is an equivalence.
	The keys to their proof are the reformulation of the sheaf condition given in \Cref{thm:excision} and technical results about when geometric realizations commute with infinite products and pullbacks akin to the results in \cite[\SAGsubsec{A.5.4}]{SAG}.
	We do not have occasion to use Berwick-Evans, Boavida de Brito, and Pavlov's result in this text.
\end{remark}

\begin{remark}[{(notations)}]
	Our notations $ \LRR $ and $ \Lhi $ are chosen in analogy with the standard notations in \textit{unstable motivic homotopy theory} \cites[\S2.2]{MotivicNorms:BachmannHoyois}{MR2242284}{MR2275634}{MR2934577}{MR1648048}{MR1813224}.
	To explain this, let us give an overview of how motivic spaces are defined.

	Let $ S $ be a scheme.
	We say that a presheaf $ F $ on the category $ \Sm_S $ of smooth schemes of finite type over $ S $ is \emph{$ \AA^1 $-invariant} if for every $ X \in \Sm_S $, the projection $ \pr_X \colon X \times_S \AA_S^1 \to X $ induces an equivalence
	\begin{equation*}
	 \prupperstar_X \colon \isomto{F(X)}{F(X \cross_S \AA_S^1)} \period
	\end{equation*}
	Write $ \PSh_{\AA^1}(\Sm_S) \subset \PSh(\Sm_S) $ for the full subcategory spanned by the $ \AA^1 $-invariant presheaves of spaces on $ \Sm_S $.
	The inclusion $ \PSh_{\AA^1}(\Sm_S) \subset \PSh(\Sm_S) $ admits a left adjoint, typically denoted by $ \LAA $ and called \emph{$ \AA^1 $-localization}.
	The \category of motivic spaces over $ S $ is defined as the \category
	\begin{equation*}
		\Sh_{\Nis,\AA^1}(\Sm_S) \colonequals \Sh_{\Nis}(\Sm_S) \intersect \PSh_{\AA^1}(\Sm_S)
	\end{equation*}
	of presheaves of spaces on $ \Sm_S $ that are $ \AA^1 $-invariant as well as sheaves for the \emph{Nisnevich topology} on $ \Sm_S $.
	The inclusion
	\begin{equation*}
		\Sh_{\Nis,\AA^1}(\Sm_S) \subset \Sh_{\Nis}(\Sm_S) 
	\end{equation*}
	of motivic spaces into Nisnevich sheaves on $ \Sm_S $ also admits a left adjoint, typically denoted by $ \Lmot $ and called \emph{motivic localization}.
	An important point is that the functor
	\begin{equation*}
		\LAA \colon \fromto{\PSh(\Sm_S)}{\PSh_{\AA^1}(\Sm_S)}
	\end{equation*}
	does not carry Nisnevich sheaves to Nisnevich sheaves, so $ \Lmot $ is not given by simply restricting $ \LAA $ to Nisnevich sheaves.

	Here, we should think $ \Mfld $ as analogous to $ \Sm_S $ and $ \Sh(\Mfld;\Spc) $ as analogous to $ \Sh_{\Nis}(\Sm_S) $.
	In analogy with $ \LAA $, we have chosen to use the notation $ \LRR $ for the functor inverting $ \RR $ at the level of presheaves.
	Similarly, we have used letters for the sheaf version of inverting $ \RR $.
	The ``hi'' in $ \Lhi $ stands for ``homotopy invariant''.
\end{remark}

%-------------------------------------------------------------------%
%  Changing coefficients                                            %
%-------------------------------------------------------------------%

\subsubsection{Changing coefficients}

We conclude this chapter by explaining how changing the coefficient presentable \category $ C $ interacts with the four adjoints $ \Gammalowersharp \leftadjoint \Gammaupperstar \leftadjoint \Gammalowerstar \leftadjoint \Gammauppersharp $.

\begin{observation}
	Let $ \adjto{\fupperstar}{D}{C}{\flowerstar} $ be an adjunction between presentable \categories.
	Since $ \flowerstar $ preserves limits, pointwise application of $ \flowerstar $ defines a functor
	\begin{equation*}
		\fromto{\Sh(\Mfld;C)}{\Sh(\Mfld;D)}
	\end{equation*}
	that we also denote by $ \flowerstar $.
	The functor $ \flowerstar \colon \fromto{\Sh(\Mfld;C)}{\Sh(\Mfld;D)} $ preserves limits and accessible, hence admits a left adjoint that we also denote by 
	\begin{equation*}
		\fupperstar \colon \fromto{\Sh(\Mfld;D)}{\Sh(\Mfld;C)} \period
	\end{equation*}
\end{observation}

\begin{warning}
	Pointwise application of $ \fupperstar \colon \fromto{D}{C} $ need \textit{not} preserve sheaves, hence the functor
	\begin{equation*}
		\fupperstar \colon \fromto{\Sh(\Mfld;D)}{\Sh(\Mfld;C)}
	\end{equation*}
	is \textit{not} generally given by pointwise application of $ \fupperstar \colon \fromto{D}{C} $.
	The left adjoint to the functor $ \flowerstar \colon \fromto{\Sh(\Mfld;C)}{\Sh(\Mfld;D)} $ is given by the composite 
	\begin{equation*}
		\begin{tikzcd}[sep=3em]
			\Sh(\Mfld;D) \arrow[r, hooked] & \PSh(\Mfld;D) \arrow[r, "\fupperstar \of -"] & \PSh(\Mfld;C) \arrow[r, "\SMan"] & \Sh(\Mfld;C)
		\end{tikzcd}
	\end{equation*}
	of pointwise application of $ \fupperstar \colon \fromto{D}{C} $ with sheafification. 
\end{warning}

The tensor product of presentable \categories (see \cref{sec:Cvaluedsheaves}) gives us an alternative description of these functors:

\begin{observation}\label{obs:changing_coeffs_tensoring}
	Let $ \adjto{\fupperstar}{D}{C}{\flowerstar} $ be an adjunction between presentable \categories.
	The equivalences 
	\begin{equation*}
		\Sh(\Mfld;\Spc) \tensor C \equivalence \Sh(\Mfld;C) \andeq \Sh(\Mfld;\Spc) \tensor D \equivalence \Sh(\Mfld;D)
	\end{equation*}
	of \Cref{ex:ShMan_tensor} fit into canonically commutative squares
	\begin{equation*}
		\begin{tikzcd}
			\Sh(\Mfld;\Spc) \tensor C \arrow[r, "\sim"{yshift = -0.2em}] \arrow[d, "{\id{} \tensor \flowerstar}"'] & \Sh(\Mfld;C) \arrow[d, "\flowerstar"] \\ 
			\Sh(\Mfld;\Spc) \tensor D \arrow[r, "\sim"{yshift = -0.2em}] & \Sh(\Mfld;D)
		\end{tikzcd}
		\andeq
		\begin{tikzcd}
			\Sh(\Mfld;\Spc) \tensor D \arrow[r, "\sim"{yshift = -0.2em}] \arrow[d, "{\id{} \tensor \fupperstar}"'] & \Sh(\Mfld;D) \arrow[d, "\fupperstar"] \\ 
			\Sh(\Mfld;\Spc) \tensor C \arrow[r, "\sim"{yshift = -0.2em}] & \Sh(\Mfld;C) \period
		\end{tikzcd}
	\end{equation*}
\end{observation}

\begin{observation}\label{obs:Gamma_tensoring}
	Let $ C $ be a presentable \category.
	Under the identifications
	\begin{equation*}
		\Sh(\Mfld;\Spc) \tensor C \equivalent \Sh(\Mfld;C) \andeq \Spc \tensor C \equivalent C \comma
	\end{equation*}
	the chain of four adjoints
	\begin{equation*}
		\begin{tikzcd}[sep=4em]
			\Sh(\Mfld;C) \arrow[r, shift left=2ex, "\Gammalowersharp"] \arrow[r, shift right=0.75ex, "\Gammalowerstar"{description, xshift=0.5em}] & C \arrow[l, shift left=2ex, "\Gammauppersharp", hooked'] \arrow[l, shift right=0.75ex, "\Gammaupperstar"{description, xshift=-0.5em}, hooked']
		\end{tikzcd}
	\end{equation*}
	is given by applying the tensor product of presentable \categories $ (-) \tensor C $ to the chain of adjoints
	\begin{equation*}
		\begin{tikzcd}[sep=4em]
			\Sh(\Mfld;\Spc) \arrow[r, shift left=2ex, "\Gammalowersharp"] \arrow[r, shift right=0.75ex, "\Gammalowerstar"{description, xshift=0.5em}] & \Spc \period \arrow[l, shift left=2ex, "\Gammauppersharp", hooked'] \arrow[l, shift right=0.75ex, "\Gammaupperstar"{description, xshift=-0.5em}, hooked']
		\end{tikzcd}
	\end{equation*}
\end{observation}

The general yoga of tensor products of presentable \categories implies that changing coefficients is compatible with the four operations $ \Gammalowersharp \leftadjoint \Gammaupperstar \leftadjoint \Gammalowerstar \leftadjoint \Gammauppersharp $.

\begin{proposition}\label{prop:changing_coeffs}
	Let $ \adjto{\fupperstar}{D}{C}{\flowerstar} $ be an adjunction between presentable \categories.
	Then the squares
	\begin{equation*}
		\begin{tikzcd}
			C \arrow[r, "\Gammauppersharp", hooked] \arrow[d, "\flowerstar"'] & \Sh(\Mfld;C) \arrow[d, "\flowerstar"] \\ 
			D \arrow[r, "\Gammauppersharp"', hooked] & \Sh(\Mfld;D) \comma
		\end{tikzcd} 
		\qquad
		\begin{tikzcd}
			\Sh(\Mfld;C) \arrow[r, "\Gammalowerstar"] \arrow[d, "\flowerstar"'] & C \arrow[d, "\flowerstar"] \\ 
			\Sh(\Mfld;D) \arrow[r, "\Gammalowerstar"'] & D\comma
		\end{tikzcd} 
		\andeq
		\begin{tikzcd}
			C \arrow[r, "\Gammaupperstar", hooked] \arrow[d, "\flowerstar"'] & \Sh(\Mfld;C) \arrow[d, "\flowerstar"] \\ 
			D \arrow[r, "\Gammaupperstar"', hooked] & \Sh(\Mfld;D)
		\end{tikzcd}
	\end{equation*}
	canonically commute.
	Equivalently, the squares
	\begin{equation*}
		\begin{tikzcd}
			\Sh(\Mfld;D) \arrow[r, "\Gammalowerstar"] \arrow[d, "\fupperstar"'] & D \arrow[d, "\fupperstar"] \\ 
			\Sh(\Mfld;C) \arrow[r, "\Gammalowerstar"'] & C \comma
		\end{tikzcd} 
		\qquad
		\begin{tikzcd}
			D \arrow[r, "\Gammaupperstar", hooked] \arrow[d, "\fupperstar"'] & \Sh(\Mfld;D) \arrow[d, "\fupperstar"] \\ 
			C \arrow[r, "\Gammaupperstar"', hooked] & \Sh(\Mfld;C) \comma
		\end{tikzcd}
		\andeq
		\begin{tikzcd}
			\Sh(\Mfld;D) \arrow[r, "\Gammalowersharp"] \arrow[d, "\fupperstar"'] & D \arrow[d, "\fupperstar"] \\ 
			\Sh(\Mfld;C) \arrow[r, "\Gammalowersharp"'] & C
		\end{tikzcd}
	\end{equation*}
	canonically commute.
\end{proposition}

\begin{proof}
	The claim for the first collection of three squares follows from the tensor product description of \Cref{obs:changing_coeffs_tensoring,obs:Gamma_tensoring} combined with \cite[Observation 1.15]{arXiv:2108.03545}.
\end{proof}

There are a few specific instances of \Cref{prop:changing_coeffs} that we are particularly interested in. 
Most of them concern the relationship between stable and unstable coefficients.

\begin{notation}\label{ntn:Sigmainf-Omegainf}
	We write
	\begin{equation*}
		\Sigma_+^{\infty} \colon \fromto{\Spc}{\Spt}
	\end{equation*}
	for the functor sending a space $ X $ to the suspension spectrum of the space $ X_+ \colonequals X \coproduct \ast $ obtained by adding a disjoint basepoint to $ X $.
	We write
	\begin{equation*}
		\Omega^{\infty} \colon \fromto{\Spt}{\Spc}
	\end{equation*}
	for the right adjoint to $ \Sigma_+^\infty $.
	Given an $ \Eone $-ring spectrum $ R $, we write $ \Omega_R^\infty \colon \fromto{\Mod(R)}{\Spc} $ for the composite
	\begin{equation*}
		\begin{tikzcd}
			\Mod(R) \arrow[r] & \Spt \arrow[r, "\Omega^{\infty}"] & \Spc
		\end{tikzcd}
	\end{equation*}
	of the forgetful functor with $ \Omega^{\infty}$
\end{notation}

\begin{example}\label{cor:suspension_spectra_ShMfld}
	The square
	\begin{equation*}
		\begin{tikzcd}
			\Sh(\Mfld;\Spc) \arrow[r, "\Gammalowersharp"] \arrow[d, "\Sigma_+^{\infty}"'] & \Spc \arrow[d, "\Sigma_+^{\infty}"] \\ 
			\Sh(\Mfld;\Spt) \arrow[r, "\Gammalowersharp"'] & \Spt
		\end{tikzcd}
	\end{equation*}
	canonically commutes.
\end{example}

\begin{example}
	Let $ R $ be an $ \Eone $-ring spectrum.
	Then the squares
	\begin{equation*}
		\begin{tikzcd}
			\Sh(\Mfld;\Mod(R)) \arrow[r, "\Gammalowerstar"] \arrow[d, "\Omega_R^\infty"'] & \Mod(R) \arrow[d, "\Omega_R^\infty"] \\ 
			\Sh(\Mfld;\Spc) \arrow[r, "\Gammalowerstar"'] & \Spc \comma
		\end{tikzcd} 
		\andeq 
		\begin{tikzcd}
			\Mod(R) \arrow[r, "\Gammaupperstar", hooked] \arrow[d, "\Omega_R^\infty"'] & \Sh(\Mfld;\Mod(R)) \arrow[d, "\Omega_R^\infty"] \\ 
			\Spc \arrow[r, "\Gammaupperstar"', hooked] & \Sh(\Mfld;\Spc) \comma
		\end{tikzcd}
	\end{equation*}
	canonically commute.
\end{example}

\newpage
%!TEX root = ../diffcoh.tex

%-------------------------------------------------------------------%
%-------------------------------------------------------------------%
%  ℝ-localization                                                   %
%-------------------------------------------------------------------%
%-------------------------------------------------------------------% 

\section{\texorpdfstring{$ \RR $}{ℝ}-localization}\label{sec:localization}
\textit{by Peter Haine}

The purpose of this chapter is to provide formulas for the \RRlocalization functor
\begin{equation*}
	\LRR \colon \fromto{\PSh(\Mfld;C)}{\PShhi(\Mfld;C)}
\end{equation*}
(\Cref{homotopification_definition}) and left adjoint $ \Gammalowersharp \colon \fromto{\Sh(\Mfld;C)}{C} $ to the constant sheaf functor \cref{nul:Gammaadjunctions}.
Specifically, write 
\begin{equation*}
	\Deltaalg^n \colonequals \setbar{(t_0,\ldots,t_{n}) \in \RR^{n+1}}{t_0 + \cdots + t_{n} = 1} \subset \RR^{n+1} 
\end{equation*}
for the algebraic $ n $-simplex; the assignment $ \goesto{[n]}{\Deltaalg^n} $ defines a cosimplicial manifold.
We show that $ \LRR $ and $ \Gammalowersharp $ are computed by the geometric realizations
\begin{equation*}
	\LRR(F)(M) \equivalent \real{F(M \cross \Deltaalgdot)} \andeq \Gammalowersharp(E) \equivalent \real{E(\Deltaalgdot)}
\end{equation*}
(\Cref{prop:MorelVoevodsky,cor:formula_for_Gammalowersharp}).

In \cref{sec:MVSconstruction}, we give a precise statement of the main result of this section (\Cref{prop:MorelVoevodsky}), but do not prove it.
We then explain some consequences of these formulas (\cref{subsubsec:consequencesofMSV}). 
In \cref{sec:simplicialhomotopy}, we recall some background on simplicial homotopies in \categories that we need to prove the formula for $ \LRR $.
\Cref{subsec:proofofMVS} is dedicated to proving this formula.

%-------------------------------------------------------------------%
%-------------------------------------------------------------------%
%  The Morel–Suslin–Voevodsky construction                          %
%-------------------------------------------------------------------%
%-------------------------------------------------------------------%

\subsection{The Morel--Suslin--Voevodsky construction}\label{sec:MVSconstruction}

%-------------------------------------------------------------------%
%  The construction                                                 %
%-------------------------------------------------------------------%

\subsubsection{The construction}

\begin{notation}[(algebraic simplices)]
	Let $ n \geq 0 $ be an integer.
	Write $ \Deltaalg^n $ for the hyperplane in $ \RR^{n+1} $ defined by
	\begin{equation*}
		\Deltaalg^n \colonequals \setbar{(t_0,\ldots,t_{n}) \in \RR^{n+1}}{t_0 + \cdots + t_{n} = 1} \subset \RR^{n+1} \comma
	\end{equation*}
	so that as a smooth manifold $ \Deltaalg^n $ is diffeomorphic to $ \RR^{n} $.
	We call $ \Deltaalg^n $ the \textit{algebraic $ n $-simplex}.

	In the usual way, the algebraic $ n $-simplices for $ n \geq 0 $ assemble into a cosimplicial manifold
	\begin{equation*}
		\Deltaalgdot \colon \fromto{\Deltabf}{\Mfld} \period
	\end{equation*}
\end{notation}

\begin{proposition}[(Morel--Suslin--Voevodsky construction)]\label{prop:MorelVoevodsky}
	Let $ C $ be a presentable \category. 
	The left adjoint
	\begin{equation*}
		\LRR \colon \fromto{\PSh(\Mfld;C)}{\PShhi(\Mfld;C)}
	\end{equation*}
	is given by the geometric realization
	\begin{equation*}
		\LRR(F)(M) \equivalent \real{F(M \cross \Deltaalgdot)} \period
	\end{equation*}
\end{proposition}

\begin{remark}
	We call the construction $ \goesto{F}{\real{F(- \cross \Deltaalgdot)}} $ the \textit{Morel--Suslin--Voevodsky} construction.
	Morel and Voevodsky provide a very general version of the Morel--Suslin--Voevodsky construction for ``sites
	with an interval object'' \cite[\S 2.3]{MR1813224}, which covers the site $ \Mfld $ with $ \RR $ as the interval object (see also \cites[\S 4.3]{MR3727503}[\S 4]{MR3679884}).
	They attribute this argument to Suslin.

	Their arguments are model category-theoretic and apply to a more specific coefficient \categories $ C $ than we're interested in.
	Hence we provide separate argument.
	See the work of Bunk \cite[Proposition 2.5]{arXiv:2008.12263} for another model category-theoretic argument.

	So as to not take us too far afield, we settle for working with the site of manifolds rather than a general site with an interval object.
	Our proof of \Cref{prop:MorelVoevodsky} takes the approach used in Brazelton's notes on motivic homotopy theory \cite[\S 3]{Brazelton:A1}.
\end{remark}

%-------------------------------------------------------------------%
%  Consequences of the Morel–Suslin–Voevodsky construction          %
%-------------------------------------------------------------------%

\subsubsection{Consequences of the Morel--Suslin--Voevodsky construction}\label{subsubsec:consequencesofMSV}

We defer the proof of \Cref{prop:MorelVoevodsky} to \cref{subsec:proofofMVS,sec:simplicialhomotopy} and first explain the consequences of \Cref{prop:MorelVoevodsky} for the functors $ \Gammalowersharp $ and $ \Lhi $.

\begin{corollary}\label{cor:formula_for_Gammalowersharp}
	Let $ C $ be a presentable \category and $ F \in \Sh(\Mfld;C) $.
	Then:
	\begin{enumerate}[label=\stlabel{cor:formula_for_Gammalowersharp}, ref=\arabic*]
		\item\label{cor:formula_for_Gammalowersharp.1} For each $ n \geq 0 $ we have
		\begin{equation*}
			\Lhi(F)(\RR^n) \equivalent \real{E(\Deltaalgdot \cross \RR^n)} \period
		\end{equation*}

		\item\label{cor:formula_for_Gammalowersharp.2} The left adjoint $ \Gammalowersharp \colon \fromto{\Sh(\Mfld;C)}{C} $ to the constant sheaf functor is given by
		\begin{equation*}
			\Gammalowersharp(E) \equivalent \real{E(\Deltaalgdot)} \period
		\end{equation*}
	\end{enumerate}
\end{corollary}

\begin{proof}
	For \enumref{cor:formula_for_Gammalowersharp}{1}, note that by \Cref{cor:hilowershriekasacomposite}, the natural map
	\begin{equation*}
		\fromto{\restrict{\LRR(E)}{\Eucop}}{\restrict{\Lhi(E)}{\Eucop}}
	\end{equation*}
	is an equivalence.
	Hence the claim is an immediate consequence of \Cref{prop:MorelVoevodsky}.
	Item \enumref{cor:formula_for_Gammalowersharp}{2} follows from \enumref{cor:formula_for_Gammalowersharp}{1}
	and the identification $ \Gammalowersharp \equivalent \Gammalowerstar \Lhi $ (\Cref{obs:formulas_for_Lhi_and_Rhi}).
\end{proof}

\begin{nul}
	Since $ \Lhi \equivalent \Gammaupperstar \Gammalowersharp $, \Cref{lem:constantishi,cor:formula_for_Gammalowersharp} show that $ \Lhi $ is given by the formula
	\begin{equation*}
		\Lhi(E)(M) \equivalent \real{E(\Deltaalgdot)}^{\Piinf(M)} \period
	\end{equation*}
	In particular, when $ C = \Spc $, the functor $ \Lhi $ is given by the formula
	\begin{equation*}
		\Lhi(E)(M) \equivalent \Map_{\Spc}(\Piinf(M),\real{E(\Deltaalgdot)}) \period
	\end{equation*}
\end{nul}

\Cref{cor:formula_for_Gammalowersharp} also reproves the formula for the underlying homotopy type of a manifold using smooth simplices:

\begin{corollary}\label{cor:formula_for_Gammalowersharpofamanifold}
	Let $ M $ be a manifold.
	There is a natural equivalence
	\begin{equation*}
		\Piinf(M) \equivalent \real{\Map_{\Mfld}(\Deltaalgdot,M)}
	\end{equation*}
	in the \category $ \Spc $.
\end{corollary}

\begin{proof}
	Recall that we write $ \yo \colon \incto{\Mfld}{\Sh(\Mfld;\Spc)} $ for the Yoneda embedding (\Cref{ex:representable_sheaf}).
	By definition, $ \Gammalowersharp \colon \fromto{\Sh(\Mfld;\Spc)}{\Spc} $ is the left Kan extension of $ \Piinf \colon \fromto{\Mfld}{\Spc} $.
	Hence we compute 
	\begin{align*}
		\Piinf(M) &\equivalent \Gammalowersharp(\yo(M)) \\
		&\equivalent \real{\yo(M)(\Deltaalgdot)} && (\text{\Cref{cor:formula_for_Gammalowersharp}}) \\ 
		&= \real{\Map_{\Mfld}(\Deltaalgdot,M)} \period && \qedhere
	\end{align*}
\end{proof}

The formula for $ \Gammalowersharp $ also proves:

\begin{corollary}\label{cor:formula_for_Gammalowersharp_preserves_finite_products}
	Let $ C $ be a presentable \category.
	If geometric realizations commute with finite products in $ C $ (e.g., $ C $ is \atopos or stable), then the functor
	\begin{equation*}
		\Gammalowersharp \colon \fromto{\Sh(\Mfld;C)}{C}
	\end{equation*}
	preserves finite products.
\end{corollary}

\begin{remark}
	The functors $ \LRR $, $ \Lhi $, and $ \Gammalowersharp $ do not generally commute with finite limits.
	However, general category theory \cite[Proposition 3.4]{MR3570135} shows that the functor
	\begin{equation*}
		\LRR \colon \fromto{\PSh(\Mfld;\Spc)}{\PShhi(\Mfld;\Spc)}
	\end{equation*}
	is \emph{locally cartesian}: for any cospan $ E \to G \leftarrow F $ with $ E, G \in \PShhi(\Mfld;\Spc) $, the natural morphism
	\begin{equation*}
		\fromto{\LRR(E \cross_G F)}{E \cross_G \LRR(F)}
	\end{equation*}
	is an equivalence.
	Since the sheafification functor $ \SMan \colon \fromto{\PSh(\Mfld;\Spc)}{\Sh(\Mfld;\Spc)} $ is left exact, \Cref{cor:hilowershriekasacomposite} shows that $ \Lhi $ and $ \Gammalowersharp $ are locally cartesian as well.
\end{remark}

%-------------------------------------------------------------------%
%-------------------------------------------------------------------%
%  Background on simplicial homotopies in ∞-categories              %
%-------------------------------------------------------------------%
%-------------------------------------------------------------------%

\subsection{Background on simplicial homotopies in \texorpdfstring{$ \infty $}{∞}-categories}\label{sec:simplicialhomotopy}

In order to prove the Morel--Suslin--Voevodsky formula (\Cref{prop:MorelVoevodsky}), we need to use homotopies of simplicial objects in an arbitrary \category.
Since we're working natively to \categories and not in simplicial sets or simplicial presheaves, doing so requires a reformulation of the usual definition of a simplicial homotopy.

%-------------------------------------------------------------------%
%  Motivation from simplicial sets                                  %
%-------------------------------------------------------------------%

\subsubsection{Motivation from simplicial sets}

Recall that a \textit{simplicial homotopy} between morphisms of simplicial sets $ f_0,f_1 \colon \fromto{\Xdot}{\Ydot} $ consists of a morphism $ h \colon \fromto{\Xdot \cross \Delta^1}{\Ydot} $ along with identifications of the restriction of $ h $ to $ \Xdot \cross \{0\} $ with $ f_0 $ and the restriction of $ h $ to $ \Xdot \cross \{1\} $ with $ f_1 $.
First we reformulate this notion in terms of morphisms in the overcategory $ \sSet_{/\Delta^1} $.

\begin{notation}
	Write $ \uupperstar \colon \fromto{\sSet}{\sSet_{/\Delta^1}} $ for the functor $ \goesto{\Xdot}{\Xdot \cross \Delta^1} $. 
	Note that $ \uupperstar $ is right adjoint to the forgetful functor $ \ulowershriek \colon \fromto{\sSet_{/\Delta^1}}{\sSet} $.
\end{notation}

\begin{lemma}\label{lem:reformshtpy}
	Let $ \Xdot $ and $ \Ydot $ be simplicial sets.
	There is a natural bijection 
	\begin{equation*}
		\Map_{\sSet}(\Xdot \cross \Delta^1,\Ydot) \isomorphic \Map_{\sSet_{/\Delta^1}}(\uupperstar(\Xdot),\uupperstar(\Ydot)) \period 
	\end{equation*}
\end{lemma}

\begin{proof}
	Since $ \ulowershriek $ is left adjoint to $ \uupperstar $, we have natural bijections
	\begin{align*}
		\Map_{\sSet_{/\Delta^1}}(\uupperstar(\Xdot),\uupperstar(\Ydot)) &\isomorphic \Map_{\sSet}(\ulowershriek\uupperstar(\Xdot),\Ydot) \\
		&= \Map_{\sSet}(\Xdot \cross \Delta^1,\Ydot) \period \qedhere
	\end{align*}
\end{proof}

In order to use \Cref{lem:reformshtpy} to generalize simplicial homotopies to arbitrary \categories, notice that the functor $ \uupperstar $ admits an alternative interpretation that makes sense for simplicial objects in any \category.

\begin{observation}[(presheaf categories and slice categories)]
	Let $ S $ be a small category and $ s \in S $.
	Write $ \yo \colon \incto{S}{\Fun(S^{\op},\Set)} $ for the Yoneda embedding.
	The colimit-preserving extension of the ``sliced Yoneda embedding''
	\begin{align*}
		S_{/s} &\inclusion \Fun(S^{\op},\Set)_{/\yo(s)} \\ 
		[s' \to s] &\mapsto [\yo(s') \to \yo(s)] \\
		\intertext{defines an equivalence of categories}
		\Fun((S_{/s})^{\op},\Set) &\equivalence \Fun(S^{\op},\Set)_{/\yo(s)} \period
	\end{align*}
	Under this identification, the functor $ \fromto{\Fun(S^{\op},\Set)}{\Fun((S_{/s})^{\op},\Set)} $ given by precomposition with the forgetful functor $ \fromto{(S_{/s})^{\op}}{S^{\op}} $ is identified with the functor
	\begin{equation*}
		\yo(s) \cross (-) \colon \fromto{\Fun(S^{\op},\Set)}{\Fun(S^{\op},\Set)_{/\yo(s)}} \period
	\end{equation*}
	Moreover, the functor $ \yo(s) \cross (-) $ is right adjoint to the forgetful functor
	\begin{equation*}
		\fromto{\Fun(S^{\op},\Set)_{/\yo(s)}}{\Fun(S^{\op},\Set)} \period
	\end{equation*}
\end{observation}

\begin{nul}
	Specializing to the case $ S = \Deltabf $ and $ s = [1] $ shows that the functor $ \uupperstar \colon \fromto{\sSet}{\sSet_{/\Delta^1}} $ is identified with the functor
	\begin{equation*}
		\Fun(\Deltaop,\Set) \to \Fun((\Deltabf_{/[1]})^{\op},\Set)
	\end{equation*}
	given by precomposition with the forgetful functor $ \fromto{(\Deltabf_{/[1]})^{\op}}{\Deltaop} $.
	We also write
	\begin{equation*}
		\uupperstar \colon \Fun(\Deltaop,\Set) \to \Fun((\Deltabf_{/[1]})^{\op},\Set)
	\end{equation*}
	for this functor.
\end{nul}

Thus, we have a further reformulation of what a simplicial homotopy is:

\begin{corollary}\label{cor:reformshtpy}
	Let $ \Xdot $ and $ \Ydot $ be simplicial sets.
	There is a natural bijection 
	\begin{equation*}
		\Map_{\sSet}(\Xdot \cross \Delta^1,\Ydot) \isomorphic \Map_{\Fun((\Deltabf_{/[1]})^{\op},\Set)}(\uupperstar(\Xdot),\uupperstar(\Ydot)) \period 
	\end{equation*}
\end{corollary}

\noindent The benefit of \Cref{cor:reformshtpy} is that the right-hand side makes sense in \textit{any} \category.

\begin{notation}
	Write $ u \colon \fromto{(\Deltabf_{/[1]})^{\op}}{\Deltaop} $ for the forgetful functor.
	For $ i \in [1] $, write
	\begin{equation*}
		j_i \colon \incto{\Deltaop}{(\Deltabf_{/[1]})^{\op}}
	\end{equation*}
	for the fully faithful functor given on objects by the assignment
	\begin{equation*}
		\goesto{[n]}{\big[[n] \to \{i\} \inclusion [1]\big]} \comma
	\end{equation*}
	with the obvious assignment on morphisms.
	Given \acategory $ D $, write
	\begin{equation*}
		\uupperstar \colon \fromto{\Fun(\Deltaop,D)}{\Fun((\Deltabf_{/[1]})^{\op},D)} \andeq \jupperstar_i \colon \fromto{\Fun((\Deltabf_{/[1]})^{\op},D)}{\Fun(\Deltaop,D)}
	\end{equation*}
	for the functors given by precomposition with $ u $ and $ j_i $, respectively.
\end{notation}

\begin{observation}\label{obs:jadjoint}
	For each $ i \in [1] $, the fully faithful functor $ j_i \colon \incto{\Deltaop}{(\Deltabf_{/[1]})^{\op}} $ is a section of the forgetful functor $ u \colon \fromto{(\Deltabf_{/[1]})^{\op}}{\Deltaop} $.
	That is $ uj_0 = uj_1 = \id{} $.
\end{observation}

\begin{definition}[{(simplicial homotopy \HA{Definition}{7.2.1.6})}]\label{def:simplicialhomotopy}
	Let $ D $ be \acategory and let
	\begin{equation*}
		f_0,f_1 \colon \fromto{\Xdot}{\Ydot}
	\end{equation*}
	be morphisms in the \category $ \Fun(\Deltaop,D) $ of simplicial objects in $ D $.
	A \textit{simplicial homotopy} from $ f_0 $ to $ f_1 $ consists of the following data:
	\begin{enumerate}[label=\stlabel{def:simplicialhomotopy}, ref=\arabic*]
		\item A morphism $ h \colon \fromto{\uupperstar(\Xdot)}{\uupperstar(\Ydot)} $ in $ \Fun((\Deltabf_{/[1]})^{\op},D) $.

		\item Equivalences $ \jupperstar_0(h) \equivalent f_0 $ and $ \jupperstar_1(h) \equivalent f_1 $ of morphisms $ \fromto{\Xdot}{\Ydot} $ in $ \Fun(\Deltaop,D) $.
	\end{enumerate}

	We often write $ h \colon \fromto{\uupperstar(\Xdot)}{\uupperstar(\Ydot)} $ for the entire data of a simplicial homotopy from $ f_0 $ to $ f_1 $.
\end{definition}

%-------------------------------------------------------------------%
%  Cofinality & sifted ∞-categories                                 %
%-------------------------------------------------------------------%

\subsubsection{Cofinality \& sifted \texorpdfstring{$\infty$}{∞}-categories}

The fact that we need about simplicial homotopies is that if $ h \colon \fromto{\uupperstar(\Xdot)}{\uupperstar(\Ydot)} $ is a simplicial homotopy from $ f_0 $ to $ f_1 $, then $ f_0 $ and $ f_1 $ induce the same map $ \fromto{\real{\Xdot}}{\real{\Ydot}} $ on geometric realizations.
To prove this, we need to use some special properties of category $ \Deltaop $, namely that it is \textit{sifted}.
The purpose of this subsection is to review the theory of sifted \categories.

\begin{definition}\label{def:cofinality}
	Let $ p \colon \fromto{I'}{I} $ be a functor between \categories.
	We say that $ p $ is \emph{colimit-cofinal} if for every \category $ D $ and diagram $ F \colon \fromto{I}{D} $, the following conditions hold:
	\begin{enumerate}[label=\stlabel{def:cofinality}, ref=\arabic*]
		\item The colimit $ \colim_I F $ exists if and only if the colimit $ \colim_{I'} Fp $ exists.

		\item If the colimit $ \colim_I F $ exists, then the natural map
		\begin{equation*}
			\textstyle \colim_I F \to \colim_{I'} Fp
		\end{equation*}
		is an equivalence.
	\end{enumerate}
\end{definition}

\begin{definition}
	\Acategory $ I $ is \emph{sifted} if $ I $ is nonempty and the diagonal functor $ \fromto{I}{I \cross I} $ is colimit-cofinal.
\end{definition}

That is, siftedness means that we can compute colimits indexed over $ I \cross I $ as colimits indexed over the diagonal copy of $ I $.
There are two key examples of sifted \categories: 

\begin{example}
	Every filtered \category is sifted \HTT{Example}{5.5.8.3}.
	The category $ \Deltaop $ is sifted \HTT{Lemma}{5.5.8.4}
\end{example}

It is useful to have an explicit condition to verify in order to see that \acategory is sifted.
Quillen's Theorem A \HTT{Theorem}{4.1.3.1} implies the following reformulations of siftedness.
For these, recall that \acategory $ C $ is \textit{weakly contractible} if the classifying space of $ C $ is contractible.

\begin{lemma}\label{lem:sifted_reformulations}
	Let $ I $ be a nonempty \category.
	The following conditions are equivalent:
	\begin{enumerate}[label=\stlabel{lem:sifted_reformulations}, ref=\arabic*]
		\item The \category $ I $ is sifted.

		\item For all $ i,i' \in I $, the \category $ I_{i/} \cross_{I} I_{i'/} $ is weakly contractible.

		\item For each $ i \in I $, the forgetful functor $ \fromto{I_{i/}}{I} $ is colimit-cofinal.
	\end{enumerate}
\end{lemma}

\begin{example}\label{ex:u_colimit-cofinal}
	Since the category $ \Deltaop $ is sifted, the forgetful functor 
	\begin{equation*}
		u \colon (\Deltabf_{/[1]})^{\op} = (\Deltaop)_{[1]/} \to \Deltaop
	\end{equation*}
	is colimit-cofinal.
\end{example}	

%-------------------------------------------------------------------%
%  Realizations of simplicial homotopies                            %
%-------------------------------------------------------------------%

\subsubsection{Realizations of simplicial homotopies}

Finally, we prove that simplicially homotopic maps induce equivalent maps on geometric realizations.
We begin with an observation.

\begin{observation}\label{obs:ji_inverse_equivalence}
	Let $ D $ be \acategory that admits geometric realizations of simplicial objects and let $ \Xdot $ be a simplicial object in $ D $.
	For $ i \in [1] $, the induced maps
	\begin{equation*}
		\begin{tikzcd}
			\alpha_i \colon \displaystyle \colim_{\Deltaop} \Xdot = \colim_{\Deltaop} (uj_i)\upperstar(\Xdot) \arrow[r] & \displaystyle \colim_{(\Deltabf_{/[1]})^{\op}} \uupperstar(\Xdot)
		\end{tikzcd} 
	\end{equation*}
	are equivalences and are homotopic.
	Specifically, they are both inverses of the equivalence
	\begin{equation*}
		\begin{tikzcd}
			\displaystyle \colim_{(\Deltabf_{/[1]})^{\op}} \uupperstar(\Xdot) \arrow[r, "\sim"{yshift=-0.25em}] & \displaystyle \colim_{\Deltaop} \Xdot 
		\end{tikzcd} 
	\end{equation*}
	provided by the fact that $ u $ is colimit-cofinal (\Cref{ex:u_colimit-cofinal}).
\end{observation}

\begin{lemma}\label{lem:simplicialhomotopyequivalentmaps}
	Let $ D $ be \acategory that admits geometric realizations of simplicial objects.
	Let $ f_0, f_1 \colon \fromto{\Xdot}{\Ydot} $ be morphisms of simplicial objects in $ D $ and let $ h $ be a simplicial homotopy from $ f_0 $ to $ f_1 $.
	Then the simplicial homotopy $ h $ induces an equivalence $ \real{f_0} \equivalent \real{f_1} $ between the induced morphisms
	\begin{equation*}
		\real{f_0}, \real{f_1} \colon \fromto{\real{\Xdot}}{\real{\Ydot}}
	\end{equation*}
	on geometric realizations.
\end{lemma}

\begin{proof}
	In light of \Cref{obs:ji_inverse_equivalence}, for $ i \in [1] $ we have a commutative square
	\begin{equation*}
		\begin{tikzcd}[sep=3em]
			{\real{\Xdot}} \arrow[d, "\real{\jupperstar_i(h)}"'] \arrow[r, "\alpha_i", "\sim"'{yshift=0.15em}] & \displaystyle \colim_{(\Deltabf_{/[1]})^{\op}} \uupperstar(\Xdot) \arrow[d, "\colim\limits_{(\Deltabf_{/[1]})^{\op}} h"] \\
			{\real{\Ydot}} \arrow[r, "\sim"{yshift=-0.2em}, "\alpha_i"'] & \displaystyle\colim_{(\Deltabf_{/[1]})^{\op}} \uupperstar(\Ydot) \period
		\end{tikzcd}
	\end{equation*}
	Thus the equivalences 
	\begin{equation*}
		\jupperstar_0(h) \equivalent f_0 \andeq \jupperstar_1(h) \equivalent f_1
	\end{equation*}
	provided by the simplicial homotopy $ h $ induce equivalences
	\begin{equation*}
		\real{f_0} \equivalent \real{\jupperstar_0(h)} \equivalent \colim_{(\Deltabf_{/[1]})^{\op}} h \equivalent \real{\jupperstar_1(h)} \equivalent \real{f_1}
	\end{equation*}
	in the arrow \category of $ D $.
\end{proof}

%-------------------------------------------------------------------%
%-------------------------------------------------------------------%
%  Proof of the Morel–Suslin–Voevodsky formula                      %
%-------------------------------------------------------------------%
%-------------------------------------------------------------------%

\subsection{Proof of the Morel--Suslin--Voevodsky formula}\label{subsec:proofofMVS}

We prove \Cref{prop:MorelVoevodsky} by applying the following recognition principle for localization functors.

\begin{proposition}[{\HTT{Proposition}{5.2.7.4}}]\label{prop:HTT.5.2.7.4}
	Let $ C $ be \acategory and $ L \colon \fromto{D}{D} $ a functor with essential image $ LD \subset D $.
	Then the following are equivalent:
	\begin{enumerate}[label=\stlabel{prop:HTT.5.2.7.4}, ref=\arabic*]
		\item\label{prop:HTT.5.2.7.4.1} There exists a functor $ F \colon \fromto{D}{D'} $ with fully faithful right adjoint $ G \colon \incto{D'}{D} $ such that $ GF \equivalent L $.

		\item\label{prop:HTT.5.2.7.4.2} The functor $ L \colon \fromto{D}{LD} $ is left adjoint to the inclusion $ \incto{LD}{D} $.

		\item\label{prop:HTT.5.2.7.4.3} There is a natural transformation $ \unit \colon \fromto{\id{D}}{L} $ such that for all $ d \in D $, the morphisms
		\begin{equation*}
			\unit_{L(d)}, L(\unit_d) \colon \fromto{L(d)}{L(L(d))}
		\end{equation*}
		are equivalences.
	\end{enumerate}
\end{proposition}

\begin{notation}
	Let us temporarily write $ \Hup \colon \fromto{\PSh(\Mfld;C)}{\PSh(\Mfld;C)} $ for the Morel--Suslin--Voevodsky construction
	\begin{equation*}
		\Hup(F)(M) \colonequals \real{F(M \cross \Deltaalgdot)} \period
	\end{equation*}
\end{notation}

\begin{construction}
	Let $ C $ be a presentable \category.
	Define a natural transformation
	\begin{equation*}
		\unit \colon \fromto{\id{\PSh(\Mfld;C)}}{\Hup} 
	\end{equation*}
	as follows.
	Let $ M $ be a manifold, and also simply write $ M $ for the constant cosimplicial manifold at $ M $.
	Projection onto the first factor defines a morphism of cosimplicial manifolds
	\begin{equation*}
		\pr_M \colon \fromto{M \cross \Deltaalgdot}{M}
	\end{equation*}
	from the product cosimplicial manifold $ M \cross \Deltaalgdot $ to the constant cosimplicial manifold at $ M $.
	For each $ C $-valued presheaf $ F \in \PSh(\Mfld;C) $, the morphism $ \unit_F \colon \fromto{F}{\Hup(F)} $ is defined as the geometric realization
	\begin{equation*}
		\unit_F(M) \colonequals \real{\prupperstar_M} \colon F(M) \equivalence \real{F(M)} \to \real{F(M \cross \Deltaalgdot)} = \Hup(F)(M) \period
	\end{equation*}

	Equivalently, the morphism $ \unit_F(M) $ is the composite
	\begin{equation*}
		F(M) \equivalent F(M \cross \Deltaalg^0) \to \real{F(M \cross \Deltaalgdot)}
	\end{equation*}
	of the equivalence $ \equivto{F(M)}{F(M \cross \Deltaalg^0)} $ induced by the projection $ \isomto{M \cross \Deltaalg^0}{M} $ with the induced map
	\begin{equation*}
		\fromto{F(M \cross \Deltaalg^0)}{\real{F(M \cross \Deltaalgdot)}}
	\end{equation*}
	from the $ 0 $-simplices of the simplicial object $ F(M \cross \Deltaalgdot) $ to its geometric realization. 
\end{construction}

%-------------------------------------------------------------------%
%  Proof of ℝ-invariance                                            %
%-------------------------------------------------------------------%

\subsubsection{Proof of \texorpdfstring{$ \RR $}{ℝ}-invariance}

In order to apply \Cref{prop:HTT.5.2.7.4}, the we first check:

\begin{lemma}\label{lem:Hishi}
	Let $ C $ be a presentable \category.
	For any presheaf $ F \colon \fromto{\Mfldop}{C} $, the presheaf $ \Hup(F) $ is \RRinvariant.
\end{lemma}

\noindent To prove \Cref{lem:Hishi}, we apply the technology of simplicial homotopies.

\begin{lemma}\label{lem:simplicialhomotopyindpr}
	Let $ M $ be a manifold.
	There is a natural simplicial homotopy in $ \Mfldop $ from the map
	\begin{equation*}
		i_{M \cross \Deltaalgdot,0}\, \of \pr_{M \cross \Deltaalgdot} \colon M \cross \Deltaalgdot \cross \RR \to M \cross \Deltaalgdot \cross \RR
	\end{equation*}
	to the identity.
\end{lemma}

\begin{proof}
	Define a simplicial homotopy
	\begin{equation*}
		h \colon \fromto{\uupperstar(M \cross \Deltaalgdot \cross \RR)}{\uupperstar(M \cross \Deltaalgdot \cross \RR)}
	\end{equation*}
	as follows.
	For each map $ \sigma \colon \fromto{[n]}{[1]} $ in $ \Deltabf $, write $ h'_{\sigma} \colon \fromto{\Deltaalg^n \cross \RR}{\Deltaalg^n \cross \RR} $ for the smooth map defined by the formula
	\begin{equation*}
		h'_{\sigma}(t_0,\ldots,t_n,x) \colonequals \paren{t_0,\ldots,t_n, \textstyle x\sum_{k \in \sigmainverse(1)} t_k} \period
	\end{equation*}
	Define $ h_{\sigma} \colon \fromto{M \cross \Deltaalg^n \cross \RR}{M \cross \Deltaalg^n \cross \RR} $ by setting $ h_{\sigma} \colonequals \id{M} \cross h'_{\sigma} $.
	It is immediate from the definitions that $ h $ defines a simplicial homotopy
	\begin{equation*}
		\fromto{\uupperstar(M \cross \Deltaalgdot \cross \RR)}{\uupperstar(M \cross \Deltaalgdot \cross \RR)} \comma
	\end{equation*}
	and, moreover,
	\begin{equation*}
		\jupperstar_0(h) = i_{M \cross \Deltaalgdot,0}\, \of \pr_{M \cross \Deltaalgdot} \andeq \jupperstar_1(h) = \id{M \cross \Deltaalgdot \cross \RR} \period \qedhere
	\end{equation*}
\end{proof}

\begin{proof}[Proof of \Cref{lem:Hishi}]
	Let $ M $ be a manifold.
	Since $ \pr_M i_{M,0} = \id{M} $, to see that
	\begin{equation*}
		\prupperstar_M \colon \fromto{\Hup(F)(M)}{\Hup(F)(M \cross \RR)}
	\end{equation*}
	is an equivalence, it suffices to show that $ \prupperstar_M \iupperstar_{M,0} \equivalent \id{\Hup(F)(M \cross \RR)} $.
	This follows from combining \Cref{lem:simplicialhomotopyindpr,lem:simplicialhomotopyequivalentmaps}.
\end{proof}

%-------------------------------------------------------------------%
%  Proof that the unit is an equivalence                            %
%-------------------------------------------------------------------%

\subsubsection{Proof that the unit is an equivalence}

The second thing to check is that for every presheaf $ G $, the morphism $ \unit_{\Hup(G)} $ is an equivalence.
Combined with \Cref{lem:Hishi} this guarantees that the essential image of the functor
\begin{equation*}
	\Hup \colon \fromto{\PSh(\Mfld;C)}{\PSh(\Mfld;C)}
\end{equation*}
is $ \PShhi(\Mfld;C) $.

\begin{lemma}\label{lem:Hofhiishi}
	Let $ C $ be a presentable \category.
	If $ F \colon \fromto{\Mfldop}{C} $ is \RRinvariant, then the map $ \unit_F \colon \fromto{F}{\Hup(F)} $ is an equivalence.
\end{lemma}

\begin{proof}
	Let $ M $ be a manifold.
	Since $ F $ is \RRinvariant and $ \Deltaalg^n \isomorphic \RR^n $ for each $ n \geq 0 $, the projection  $ \pr_M \colon \fromto{M \cross \Deltaalgdot}{M} $ from the cosimplicial manifold $ M \cross \Deltaalgdot $ to the constant cosimplicial manifold at $ M $ induces an equivalence
	\begin{equation*}
		\prupperstar_M \colon \equivto{F(M)}{F(M \cross \Deltaalgdot)} 
	\end{equation*}
	of simplicial objects in $ C $.
	The claim now follows by passing to geometric realizations.
\end{proof}

\begin{corollary}\label{cor:essentialimageofH}
	Let $ C $ be a presentable \category.
	The essential image of the functor
	\begin{equation*}
		 \Hup \colon \fromto{\PSh(\Mfld;C)}{\PSh(\Mfld;C)}
	\end{equation*}
	is $ \PShhi(\Mfld;C) $. 
\end{corollary}

Now we complete the proof of \Cref{prop:MorelVoevodsky} by showing that $ \Hup(\unit_F) $ is an equivalence.

\begin{lemma}\label{lem:alphaequivs}
	Let $ C $ be a presentable \category.
	For all $ F \in \PSh(\Mfld;C) $, the maps
	\begin{equation*}
		\unit_{\Hup(F)}, \Hup(\unit_F)  \colon \fromto{\Hup(F)}{\Hup(\Hup(F))}
	\end{equation*}
	are equivalences.
\end{lemma}

\begin{proof}
	By \Cref{lem:Hofhiishi,cor:essentialimageofH}, the morphism $ \unit_{\Hup(F)} $ is an equivalence.
	To see that $ \Hup(\unit_F) \colon \fromto{\Hup(F)}{\Hup(\Hup(F))} $ is an equivalence, note that for each manifold $ M $ we have
	\begin{align*}
		\Hup(F)(M) &= \colim_{[m] \in \Deltaop} F(M \cross \Deltaalg^m) \\
		\shortintertext{and} 
		\Hup(\Hup(F))(M) &= \colim_{[m] \in \Deltaop} \colim_{[n] \in \Deltaop} F(M \cross \Deltaalg^m \cross \Deltaalg^n) \\ 
		&\equivalent \colim_{([m],[n]) \in \Deltaop \cross \Deltaop} F(M \cross \Deltaalg^m \cross \Deltaalg^n) \period
	\end{align*} 
	Moreover, the map $ \Hup(\unit_F) \colon \fromto{\Hup(F)}{\Hup(\Hup(F))} $ is induced by restriction of diagrams along the fully faithful functor
	\begin{align*}
		\Deltaop &\inclusion \Deltaop \cross \Deltaop \\ 
		[m] &\mapsto ([m],[0]) \period
	\end{align*}
	First taking the colimit over the variable $ [m] \in \Deltaop $, we see that the map $ \Hup(\unit_F)(M) $ is induced by the map from the $ 0 $-simplices $ \Hup(F)(M) $ of the simplicial object $ \Hup(F)(M \cross \Deltaalgdot) $ to its geometric realization.
	Since $ \Hup(F) $ is \RRinvariant (\Cref{lem:Hishi}), the simplicial object $ \Hup(F)(M \cross \Deltaalgdot) $ is equivalent to the constant simplicial object at $ \Hup(F)(M) $, hence the induced map 
	\begin{equation*}
		\Hup(F)(M) \to \colim_{[n] \in \Deltaop} \Hup(F)(M \cross \Deltaalg^n) 
	\end{equation*}
	from the $ 0 $-simplices is an equivalence.
\end{proof}

\begin{proof}[Proof of \Cref{prop:MorelVoevodsky}]
	Combine \Cref{cor:essentialimageofH,lem:alphaequivs,prop:HTT.5.2.7.4}.
\end{proof}

\newpage
%!TEX root = ../diffcoh.tex

%-------------------------------------------------------------------%
%-------------------------------------------------------------------%
%  Structures in the stable case                                    %
%-------------------------------------------------------------------%
%-------------------------------------------------------------------%

\section{Structures in the stable case}\label{sec:stable}
\textit{by Peter Haine}

In ordinary differential cohomology, we had the Simons--Sullivan ``differential cohomology hexagon''
\begin{equation*}\label{diag:SimonsSullivan}
	\begin{tikzcd}[column sep={10ex,between origins}, row sep={8ex,between origins}]
		0 \arrow[dr] & & & & 0 \\
		& \H^{*-1}(M;\RR/\ZZ) \arrow[rr, "-\Bock"] \arrow[dr] & & \H^*(M;\ZZ) \arrow[dr] \arrow[ur] & \\
		\HdR^{*-1}(M) \arrow[ur] \arrow[dr] & & \Hhat^*(M;\ZZ) \arrow[ur] \arrow[dr] & & \HdR^*(M) \\
		& \frac{\Omega^{*-1}(M)}{\Omegacl^{*-1}(M)_{\ZZ}} \arrow[rr, "\d"'] \arrow[ur] & & \Omegacl^*(M)_{\ZZ} \arrow[ur] \arrow[dr] & \\
		0 \arrow[ur] & & & & 0 \comma
	\end{tikzcd}
\end{equation*}
which actually characterized ordinary differential cohomology (\Cref{thm:SimonsSullivanunique}).
We want to be able to reproduce an analogue of the differential cohomology hexagon for \textit{any} sheaf of spectra on $ \Mfld $.
To do this, we need to identify how cohomology with coefficients in $ \RR/\ZZ $, $ \ZZ $, and $ \RR $ as well as $ \Omega^{*-1}(M)/\Omegacl^{*-1}(M)_{\ZZ} $ and $ \Omegacl^*(M)_{\ZZ} $ fit into the story.

One general machine for producing diagrams aesthetically similar to the differential cohomology hexagon is the
theory of \textit{recollements}, or ways of ``gluing'' a category together out of two pieces.
It turns out that the differential cohomology hexagon falls exactly into this framework: one of the subcategories that we build $ \Sh(\Mfld;\Sp) $ from is the subcategory $ \Shhi(\Mfld;C) $ of \RRinvariant sheaves, and the other piece is the subcategory of sheaves with vanishing global sections.
Since this whole story is a special case of the theory of recollements, the first half of the section (\cref{subsec:recollement}) gives a quick introduction to the theory of recollements and the key results.
In \cref{subsec:fracture}, we apply this general machinery to sheaves on manifolds to obtain the a version of differential cohomology hexagon for any sheaf of spectra on $ \Mfld $ (see \cref{nul:spectraldiffcohdiag}).
We finish the section by making precise what it means for a sheaf of spectra on $ \Mfld $ to ``refine'' a cohomology theory.

%-------------------------------------------------------------------%
%-------------------------------------------------------------------%
%  Background on recollements                                       %
%-------------------------------------------------------------------%
%-------------------------------------------------------------------%

\subsection{Background on recollements}\label{subsec:recollement} 

\textit{Recollements}%
\footnote{Roughly, the French verb \textit{recoller} means ``to glue back together''.} %
were introduced by Grothendieck and Verdier in the context of topoi \cites[Exposé IV, \S9]{MR50:7130} and by
Beĭlinson--Bernstein--Deligne in the context of triangulated categories \cite[\S1.4]{MR751966} to ``glue'' together sheaves over open-closed decompositions of a space.
However, there are many other situations in which \acategory can be ``glued together'' from two subcategories that are in some sense complementary.
For example, if $ R $ is a ring and $ I \subset R $ is a finitely generated ideal, then the derived \category of $ R $ can be clued together from its subcategories of $ I $-nilpotent and $ I $-local objects.

The goal of this section is to explain this general theory and how it can be applied to the context of sheaves of spectra on the category of manifolds.
The key insight is that given a stable \category $ \Xcat $ and a full subcategory $ \ilowerstar \colon \incto{\Zcat}{\Xcat} $ that is both localizing and colocalizing
\begin{equation*}
	\begin{tikzcd}[sep=3.5em]
		\Zcat \arrow[r, "\ilowerstar" description, hooked] & \Xcat \comma \arrow[l, shift left=1.25ex, "\iuppershriek"] \arrow[l, shift right=1.25ex, "\iupperstar"'] 
	\end{tikzcd}
\end{equation*}
the \category $ \Xcat $ can be glued together from the subcategory $ \Zcat $ and the subcategory $ \Zrorth \subset \Xcat $ \textit{right orthogonal} to $ \Xcat $ (\Cref{prop:orthogonaladjoints,cor:stablerecollement}).
That is, $ \Zrorth $ is the subcategory of objects of $ \Xcat $ that admit no nontrivial maps \textit{from} objects of $ \Zcat $.
This applies to the situation of interest because we have both a left and right adjoint
\begin{equation*}
	\begin{tikzcd}[sep=3.5em]
		\Shhi(\Mfld;\Sp) \arrow[r, hooked] & \Sh(\Mfld;\Sp) \arrow[l, shift left=1.25ex, "\Rhi"] \arrow[l, shift right=1.25ex, "\Lhi"']
	\end{tikzcd}
\end{equation*}
to the inclusion of \RRinvariant sheaves on $ \Mfld $ into all sheaves \Cref{nul:Gammaadjunctions}.
We'll apply the general theory studied in this section to the context of sheaves on $ \Mfld $ in \cref{subsec:fracture}.
 
%-------------------------------------------------------------------%
%  Motivation                                                       %
%-------------------------------------------------------------------%

\subsubsection{Motivation}\label{subsec:recollementmotivation}

To explain the motivation for recollements, let $ X $ be a topological space and $ Z \subset X $ a closed subspace.
Write $ U \colonequals X \smallsetminus Z $ for the open complement of $ Z $ in $ X $, and write 
\begin{equation*}
	i \colon \incto{Z}{X} \andeq j\colon \incto{U}{X}
\end{equation*}
for the inclusions.
Any sheaf $ F $ of sets on $ X $ pulls back to sheaves
\begin{equation*}
	F_Z \colonequals \iupperstar(F) \andeq F_U \colonequals \jupperstar(F)
\end{equation*}
on $ Z $ and $ U $, respectively.
Moreover, the sheaf $ F $ is completely determined by the sheaves $ F_Z $ and $ F_U $ in the following sense.
Applying $ \iupperstar $ to the unit $ \unit \colon \fromto{F}{\jlowerstar\jupperstar(F)} $, we obtain a natural morphism
\begin{equation*}
	u \colon F_Z = \iupperstar(F) \to \iupperstar\jlowerstar\jupperstar(F) = \iupperstar\jlowerstar(F_U) \period
\end{equation*}
The triangle identities imply that there is a commutative square
\begin{equation}\label{sq:fracturemotivation}
	\begin{tikzcd}
		F \arrow[r] \arrow[d] & \jlowerstar(F_U) \arrow[d] \\
		\ilowerstar(F_Z) \arrow[r, "\ilowerstar(u)"'] & \ilowerstar\iupperstar\jlowerstar(F_U) \comma
	\end{tikzcd}
\end{equation}
where the three morphisms
\begin{equation*}
	F \to \ilowerstar\iupperstar(F) = \ilowerstar(F_Z) \, , \quad F \to \jlowerstar\jupperstar(F) = \jlowerstar(F_U) \, , \andeq \jlowerstar(F_U) \to \ilowerstar\iupperstar\jlowerstar(F_U)
\end{equation*}
are all unit morphisms.
One can show that the square \eqref{sq:fracturemotivation} is in a \textit{pullback square}.
This provides an explicit way to reconstruct $ F $ from the data of the sheaves $ F_Z $ and $ F_U $ along with the morphism $ u \colon \fromto{F_Z}{\iupperstar \jlowerstar(F_U)} $.

In fact, even more is true.
The whole \textit{category} $ \Sh(X;\Set) $ can be reconstructed from the categories $ \Sh(Z;\Set) $ and $ \Sh(U;\Set) $ together with the functor $ \iupperstar\jlowerstar \colon \fromto{\Sh(U;\Set)}{\Sh(Z;\Set)} $ in the following sense.
Write $ [1] $ for the ``walking arrow'' poset $ \{0 < 1\} $.
There is a pullback square of categories
\begin{equation}\label{sq:fractureCatmotivation}
	\begin{tikzcd}
		\Sh(X;\Set) \arrow[r] \arrow[d, "\jupperstar"'] & \Fun([1],\Sh(Z;\Set)) \arrow[d, "\target"] \\
		\Sh(U;\Set) \arrow[r, "\iupperstar\jlowerstar"'] & \Sh(Z;\Set) \period
	\end{tikzcd}
\end{equation}
Here the unlabeled top horizontal arrow sends a sheaf $ F \in \Sh(X;\Set) $ to the morphism given by applying $ \iupperstar $ to the unit $ \fromto{F}{\jlowerstar\jupperstar(F)} $.
More explicitly, an object of $ \Sh(X;\Set) $ is equivalent to the data of a sheaf $ F_{Z} $ on $ Z $, a sheaf $ F_{U} $ on $ U $, and a \textit{gluing morphism} $ \fromto{F_Z}{\iupperstar\jlowerstar(F_U)} $.
Morphisms are morphisms of sheaves on $ Z $ and $ U $ commuting with the specified gluing morphisms.

In the rest of this section, we explain the general categorical framework for decompositions of this form.
We do not explain the proofs of the results presented in this section; for those, the reader should consult \cites[\HAsec{A.8}]{HA}[\SAGsec{7.2}]{SAG}{arXiv:1607.02064}.

%-------------------------------------------------------------------%
%  Definitions and general results                                  %
%-------------------------------------------------------------------%

\subsubsection{Definitions and general results} 

Now we generalize the situation for sheaves explained in \cref{subsec:recollementmotivation}.
The following are the key features of the situation.

\begin{definition}\label{def:recollement}
	Let $ \Xcat $ be \acategory with finite limits.
	Fully faithful functors
	\begin{equation*}
		\ilowerstar \colon \incto{\Zcat}{\Xcat} \andeq \jlowerstar \colon \incto{\Ucat}{\Xcat}
	\end{equation*}
	exhibit $ \Xcat $ as the \textit{recollement} of $ \Zcat $ and $ \Ucat $ if:
	\begin{enumerate}[label=\stlabel{def:recollement}, ref=\arabic*]
		\item\label{def:recollement.1} The functors $ \ilowerstar $ and $ \jlowerstar $ admit left exact left adjoints $ \iupperstar $ and $ \jupperstar $, respectively.

		\item\label{def:recollement.2} The functor $ \jupperstar \ilowerstar \colon \fromto{\Zcat}{\Ucat} $ is constant at the terminal object of $ \Ucat $.

		\item\label{def:recollement.3} The functors $ \iupperstar \colon \fromto{\Xcat}{\Zcat} $ and $ \jupperstar \colon \fromto{\Xcat}{\Ucat} $ are jointly conservative.
		That is, a morphism $ f $ in $ \Xcat $ is an equivalence if and only if both $ \iupperstar(f) $ and $ \jupperstar(f) $ are equivalences.
	\end{enumerate}

	We refer to the subcategory $ \Zcat \subset \Xcat $ as the \textit{closed} subcategory, and $ \Ucat \subset \Xcat $ as the \textit{open} subcategory.
\end{definition}

\begin{remark}
	Note that \enumref{def:recollement}{2} in particular implies that the are \textit{no} nontrivial maps from objects in $ \Zcat \subset \Xcat $ to objects in $ \Ucat \subset \Xcat $.
\end{remark}

\begin{warning}
	Note that the condition that $ \Xcat $ be the recollement of $ \Zcat $ and $ \Ucat $ is \textit{not} symmetric: if $ \Xcat $ is the recollement of $ \Zcat $ and $ \Ucat $, then $ \Xcat $ need note be the recollement of $ \Ucat $ and $ \Zcat $.
	For example, the composite $ \iupperstar \jlowerstar $ is not usually constant at the terminal object of $ \Zcat $.  
\end{warning}

The two most important examples of recollements from topology and algebraic geometry are the following:

\begin{example}\label{ex:Shrecollement}
	Let $ X $ be a topological space, $ i \colon \incto{Z}{X} $ a closed subspace, and $ j \colon \incto{U}{X} $ the open complement of $ Z $ in $ X $.
	Let $ C $ be a presentable \category that is compactly generated or stable.
	Then the pushforward functors
	\begin{equation*}
		\ilowerstar \colon \incto{\Sh(Z;C)}{\Sh(X;C)} \andeq \jlowerstar \colon \incto{\Sh(U;C)}{\Sh(X;C)}
	\end{equation*}
	exhibit $ \Sh(X;C) $ as the recollement of $ \Sh(Z;C) $ and $ \Sh(U;C) $.
	See \cites[\HAappthm{Remark}{A.8.16}]{HA}[Corollaries 2.12 \& 2.23]{arXiv:2108.03545}
\end{example}

\begin{example}\label{ex:QCohrecollement}
	Let $ X $ be a scheme, $ \incto{Z}{X} $ a closed subscheme, and $ \incto{U}{X} $ the complementary open subscheme in $ X $.
	Assume that $ U $ is quasicompact.
	We write $ \QCoh(X) $ and $ \QCoh(U) $ for the stable \categories of quasicoherent sheaves on $ X $ and $ U $, respectively.
	We write
	\begin{equation*}
		\QCoh_{Z}(X) \subset \QCoh(X)
	\end{equation*}
	for the full subcategory spanned by those quasicoherent sheaves that are set-theoretically supported on $ Z $.
	In this setting, the pushforward functor $ \incto{\QCoh(U)}{\QCoh(X)} $ and the inclusion $ \QCoh_{Z}(X) \subset \QCoh(X) $ exhibit $ \QCoh(X) $ as the recollement of $ \QCoh(U) $ and $ \QCoh_{Z}(X) $.
	See, for example, \SAG{Proposition}{7.2.3.1}.
\end{example}

\begin{warning}
	In \Cref{ex:QCohrecollement}, note that the subcategory $ \QCoh(U) $ is the \textit{closed} subcategory, and the subcategory $ \QCoh_{Z}(X) $ is the \textit{open} subcategory.
	There are thus two competing naming conventions for the ``closed'' and ``open'' subcategories: one coming from the theory of sheaves on topological spaces (\Cref{ex:Shrecollement}), and one coming from quasicoherent sheaves on schemes (\Cref{ex:QCohrecollement}).
	Both are used in the literature, depending on whether one is working in a ``topological'' or
	``algebro-geometric'' context. 
	In this text we use the ``topological'' convention.
\end{warning}

The following result explains how to reconstruct a recollement from the closed and open subcategories together with \textit{gluing functor} $ \iupperstar\jlowerstar \colon \fromto{\Ucat}{\Zcat} $.

\begin{theorem}[{\cites[\HAappthm{Corollary}{A.8.13}, \HAappthm{Remark}{A.8.5}, \& \HAappthm{Proposition}{A.8.17}]{HA}[1.17]{arXiv:1909.03920}}]\label{thm:recollfracture}
	Let $ \ilowerstar \colon \incto{\Zcat}{\Xcat} $ and $ \jlowerstar \colon \incto{\Ucat}{\Xcat} $ be functors that exhibit $ \Xcat $ as the recollement of $ \Zcat $ and $ \Ucat $.
	There is a pullback square of \categories
	\begin{equation*}
		\begin{tikzcd}
			\Xcat \arrow[r] \arrow[d, "\jupperstar"'] & \Fun([1],\Zcat) \arrow[d, "\target"] \\
			\Ucat \arrow[r, "\iupperstar\jlowerstar"'] & \Zcat \period
		\end{tikzcd}
	\end{equation*}
	Here the unlabeled top horizontal arrow sends an object $ F \in \Xcat $ to the morphism given by applying $ \iupperstar $ to the unit $ \fromto{F}{\jlowerstar\jupperstar(F)} $.

	As a consequence, there is a pullback square of endofunctors of $ \Xcat $
	\begin{equation}\label{sq:fracturegeneral}
		\begin{tikzcd}
			\id{\Xcat} \arrow[r] \arrow[d] & \jlowerstar\jupperstar \arrow[d] \\
			\ilowerstar\iupperstar \arrow[r] & \ilowerstar\iupperstar\jlowerstar\jupperstar \period
		\end{tikzcd}
	\end{equation}
	Here the top horizontal and left vertical morphisms are the unit morphisms, the bottom horizontal morphism is obtained by applying $ \ilowerstar\iupperstar $ to the unit morphism $ \fromto{\id{\X}}{\jlowerstar\jupperstar} $, and the right vertical morphism is obtained by precomposing the unit morphism $ \fromto{\id{\X}}{\ilowerstar\iupperstar} $ with $ \jlowerstar\jupperstar $.
\end{theorem}

\begin{definition}
	Let $ \ilowerstar \colon \incto{\Zcat}{\Xcat} $ and $ \jlowerstar \colon \incto{\Ucat}{\Xcat} $ be functors that exhibit $ \Xcat $ as the recollement of $ \Zcat $ and $ \Ucat $.
	The pullback square \eqref{sq:fracturegeneral} is referred to as the \textit{fracture square} of the recollement.
\end{definition}

Often the functors $ \ilowerstar $ and $ \jupperstar $ admit further adjoints.

\begin{theorem}[{\cites[\HAappthm{Corollary}{A.8.7}, \HAappthm{Remark}{A.8.8}, \& \HAappthm{Proposition}{A.8.11}]{HA}[Corollary 1.10]{arXiv:1909.03920}}]\label{thm:stabrecollement}
	Let $ \ilowerstar \colon \incto{\Zcat}{\Xcat} $ and $ \jlowerstar \colon \incto{\Ucat}{\Xcat} $ be functors that exhibit $ \Xcat $ as the recollement of $ \Zcat $ and $ \Ucat $.
	\begin{enumerate}[label=\stlabel{thm:stabrecollement}, ref=\arabic*]
		\item\label{thm:stabrecollement.1} If the \category $ \Zcat $ has an initial object, then $ \jupperstar $ admits a fully faithful left adjoint $ \jlowershriek \colon \incto{\Ucat}{\Xcat} $.

		\item\label{thm:stabrecollement.2} If, moreover, $ \Xcat $ has a zero object, then $ \ilowerstar $ admits a right adjoint $ \iuppershriek \colon \fromto{\Xcat}{\Zcat} $ characterized by the property that 
		\begin{equation*}
			\ilowerstar\iuppershriek \equivalent \fib(\eta \colon \id{\Xcat} \to \jlowerstar\jupperstar) \period
		\end{equation*}
		In particular, applying $ \iupperstar $, there is a fiber sequence
		\begin{equation*}
			\begin{tikzcd}
				\iuppershriek \arrow[r] & \iupperstar \arrow[r, "\iupperstar\unit"] & \iupperstar\jlowerstar\jupperstar \period
			\end{tikzcd}
		\end{equation*}

		\item\label{thm:stabrecollement.3} If $ \Xcat $ is stable, then $ \Zcat $ and $ \Ucat $ are also stable.
		Moreover, there is a canonical fiber sequence
		\begin{equation*}
			\begin{tikzcd}
				\jlowershriek \jupperstar \arrow[r] & \id{\Xcat} \arrow[r] & \ilowerstar \iupperstar \comma
			\end{tikzcd}
		\end{equation*}
		where the first morphism is the counit and the second is the unit.

		\item\label{thm:stabrecollement.4} If $ \Xcat $ is presentable and the gluing functor $ \iupperstar \jlowerstar $ is accessible, then $ \Zcat $ and $ \Ucat $ are presentable.
	\end{enumerate}
\end{theorem}

\begin{nul}\label{nul:stabelrecadjunction}
	Thus, if $ \Xcat $ is stable, there is a chain of adjunctions
	\begin{equation*}
		\begin{tikzcd}[sep=3.5em]
			\Zcat \arrow[r, "\ilowerstar" description, hooked] & \Xcat \arrow[l, shift left=1.25ex, "\iuppershriek"] \arrow[l, shift right=1.25ex, "\iupperstar"'] \arrow[r, "\jupperstar" description]  & \Ucat \period \arrow[l, shift left=1.25ex, hooked', "\jlowerstar"] \arrow[l, shift right=1.25ex, hooked', "\jlowershriek"']
		\end{tikzcd}
	\end{equation*}
\end{nul}

We're interested in applying this to the situation where $ \ilowerstar $ is the inclusion of $ \Shhi(\Mfld;\Sp) $ into $ \Sh(\Mfld;\Sp) $, $ \iupperstar $ is $ \Lhi $, and $ \iuppershriek $ is $ \Rhi $.
To get an analogue of the ``differential cohomology hexagon'', we need to enlarge the fracture square \eqref{sq:fracturegeneral} using the fiber sequences from \enumref{thm:stabrecollement}{2} and \enumref{thm:stabrecollement}{3} along with one more.

\begin{construction}[(norm map)]
	Let $ \Xcat $ and $ \Ucat $ be \categories, and suppose we are given adjunctions
	\begin{equation*}
		\begin{tikzcd}[sep=3.5em]
			\Xcat \arrow[r, "\jupperstar" description] & \Ucat  \arrow[l, shift left=1.25ex, hooked', "\jlowerstar"] \arrow[l, shift right=1.25ex, hooked', "\jlowershriek"']
		\end{tikzcd}
	\end{equation*}
	where left adjoint $ \jlowershriek $ and right adjoint $ \jlowerstar $ are fully faithful.
	Write $ \counit \colon \fromto{\jupperstar\jlowerstar}{\id{\Ucat}} $ for the counit.
	Since $ \jlowershriek $ is left adjoint to $ \jupperstar $ and the counit $ \counit $ is an equivalence, we have equivalences
	\begin{equation}\label{eq:normequiv}
		\begin{tikzcd}[sep=3.5em]
			\Map(\jlowershriek,\jlowerstar) \equivalent \Map(\id{\Ucat},\jupperstar\jlowerstar) \arrow[r, "\counit \of -"', "\sim"] & \Map(\id{\Ucat},\id{\Ucat}) \period
		\end{tikzcd} 
	\end{equation}
	The \textit{norm} natural transformation
	\begin{equation*}
		\Nm \colon \fromto{\jlowershriek}{\jlowerstar}
	\end{equation*}
	is the natural transformation corresponding to the identity $ \fromto{\id{\Ucat}}{\id{\Ucat}} $ under the equivalence \eqref{eq:normequiv}.
\end{construction}

\begin{theorem}\label{thm:fracturehexagon}
	Let $ \Xcat $ be a stable \category and let $ \ilowerstar \colon \incto{\Zcat}{\Xcat} $ and $ \jlowerstar \colon \incto{\Ucat}{\Xcat} $ be functors that exhibit $ \Xcat $ as the recollement of $ \Zcat $ and $ \Ucat $.
	Then the sequence
	\begin{equation*}
		\begin{tikzcd}
			\jlowershriek \jupperstar \arrow[r, "\Nm\jupperstar"] & \jlowerstar\jupperstar \arrow[r] & \ilowerstar \iupperstar \jlowerstar\jupperstar 
		\end{tikzcd}
	\end{equation*}
	is a fiber sequence.
	As a consequence, the fracture square fits into a commutative diagram
	\begin{equation}\label{sq:fracturehexagon}
		\begin{tikzcd}[sep=2.5em]
			 &  \jlowershriek \jupperstar \arrow[r, equals] \arrow[d] & \jlowershriek \jupperstar \arrow[d, "\Nm\jupperstar"] \\
			\ilowerstar\iuppershriek \arrow[r] \arrow[d, equals] & \id{\Xcat} \arrow[r] \arrow[d] & \jlowerstar\jupperstar \arrow[d] \\
			\ilowerstar\iuppershriek \arrow[r] & \ilowerstar\iupperstar \arrow[r] & \ilowerstar\iupperstar\jlowerstar\jupperstar 
		\end{tikzcd}
	\end{equation}
	where all rows and columns are fiber sequences.
\end{theorem}

Aside from the explicit identification of the first map in the lower horizontal fiber sequence of \eqref{sq:fracturegeneral} with the norm map, \Cref{thm:fracturehexagon} can be deduced by applying the following characterization of pullback squares of stable \categories horizontally and vertically to the fracture square \eqref{sq:fracturegeneral}.

\begin{recollection}\label{rec:stablepullbackviafibers}
	Let $ C $ be a pointed \category and 
	\begin{equation}\label{sq:stablepullback}
		\begin{tikzcd}
			W \arrow[r, "\fbar"] \arrow[d] & Y \arrow[d] \\
			X \arrow[r, "f"'] & Z  
		\end{tikzcd}
	\end{equation}
	a commutative square in $ C $.
	Then there is a natural equivalence
	\begin{equation*}
		\fib(W \to X \cross_Z Y) \equivalent \fib(\fib(\fbar) \to \fib(f)) \period
	\end{equation*}
	In particular, if $ C $ is stable, then $ \equivto{\fib(\fbar)}{\fib(f)} $ if and only if the square \eqref{sq:stablepullback} is a pullback square.
	See \cites[\S2]{Chromfracture:BarthelAntolin}{MO:333239} for more details.
\end{recollection}

%-------------------------------------------------------------------%
%  Orthogonal complements & the stable situation                    %
%-------------------------------------------------------------------%

\subsubsection{Orthogonal complements \& the stable situation} 

In the stable case, it turns out that the data of a recollement of $ \Xcat $ is equivalent to the data of the closed subcategory $ \Zcat \subset \Xcat $.
The open subcategory $ \Ucat \subset \Xcat $ can be recovered as an \textit{orthogonal complement} to $ \Zcat $ in the following sense.

\begin{definition}\label{def:orthogonal}
	Let $ \Xcat $ be \acategory and $ \Zcat \subset \Xcat $ a full subcategory.
	\begin{enumerate}[label=\stlabel{def:orthogonal}, ref=\arabic*]
		\item We say that an object $ X \in \Xcat $ is \textit{right orthogonal} to the subcategory $ \Zcat $ if for each $ Z \in \Zcat $, the mapping space 
		$ \Map_{\Xcat}(Z,X) $ is contractible.

		\item We say that an object $ X \in \Xcat $ is \textit{left orthogonal} to the subcategory $ \Zcat $ if for each $ Z \in \Zcat $, the mapping space 
		$ \Map_{\Xcat}(X,Z) $ is contractible.
	\end{enumerate} 
	The \textit{right orthogonal complement} of $ \Zcat $ is the full subcategory $ \Zrorth \subset \Xcat $ spanned by those objects right orthogonal to $ \Zcat $.
	The \textit{left orthogonal complement} of $ \Zcat $ is the full subcategory $ \Zlorth \subset \Xcat $ spanned by those objects right orthogonal to $ \Zcat $.
\end{definition}

\begin{proposition}[{\cites[\SAGthm{Proposition}{7.2.1.10}]{SAG}[Lemmas 2 \& 5 and Proposition 7]{arXiv:1607.02064}}]\label{prop:orthogonaladjoints}
	Let $ \Xcat $ be a stable \category, and $ \ilowerstar \colon \incto{\Zcat}{\Xcat} $ a full subcategory.
	Assume that the inclusion $ \ilowerstar $ admits a left adjoint $ \iupperstar $ and a right adjoint $ \iuppershriek $.
	Then:
	\begin{enumerate}[label=\stlabel{prop:orthogonaladjoints}, ref=\arabic*]
		\item\label{prop:orthogonaladjoints.1} The inclusion $ \Zrorth \subset \Xcat $ admits a left adjoint $ j^{\orthogonal} \colon \fromto{\Xcat}{\Zrorth} $ defined as the cofiber
		\begin{equation*}
			j^{\orthogonal} \colonequals \cofib(\counit \colon \fromto{\ilowerstar\iuppershriek}{\id{\Xcat}}) \period
		\end{equation*}

		\item\label{prop:orthogonaladjoints.2} The inclusion $ \Zlorth \subset \Xcat $ admits a right adjoint $ ^{\orthogonal}j \colon \fromto{\Xcat}{\Zlorth} $ defined as the fiber
		\begin{equation*}
			^{\orthogonal}j \colonequals \fib(\unit \colon \fromto{\id{\Xcat}}{\iupperstar\ilowerstar}) \period
		\end{equation*}

		\item\label{prop:orthogonaladjoints.3} The composite functors
		\begin{equation*}
			\begin{tikzcd}
				\Zrorth \arrow[r, hooked] & \Xcat \arrow[r, "^{\smallorthogonal}j"] & \Zlorth
			\end{tikzcd}
			\andeq
			\begin{tikzcd}
				\Zlorth \arrow[r, hooked] & \Xcat \arrow[r, "j^{\smallorthogonal}"] & \Zrorth
			\end{tikzcd}
		\end{equation*}
		are inverse equivalences of \categories.

		\item\label{prop:orthogonaladjoints.4} The stable \category $ \Xcat $ is the recollement of the stable subcategories $ \Zcat $ and $ \Zrorth $.
	\end{enumerate}
\end{proposition}

\begin{corollary}\label{cor:stablerecollement}
	Let $ \Xcat $ be a stable \category, and let $ \ilowerstar \colon \incto{\Zcat}{\Xcat} $ and $ \jlowerstar \colon \incto{\Ucat}{\Xcat} $ be functors that exhibit $ \Xcat $ as the recollement of $ \Zcat $ and $ \Ucat $.
	Then the essential image of the fully faithful functor $ \jlowerstar $ is the right orthogonal complement $ \Zrorth $ of $ \Zcat $.
\end{corollary}

\noindent Said differently, \textit{every} stable recollement arises via \Cref{prop:orthogonaladjoints}.

\begin{remark}[(semiorthognal decompositions)]
	\Cref{prop:orthogonaladjoints,cor:stablerecollement} say that recollements are special types of \textit{semiorthogonal decompositions} of \categories.
	Semi\-orth\-ogonal decompositions were originally introduced (in the context of triangulated categories) by Bondal and Kapranov \cite{MR1039961} to break apart stable \categories arising in algebraic geometry into more simple pieces.
	There are many beautiful examples (namely, Beĭlinson's celebrated semiorthogonal decomposition of $ \Coh(\PP^n) $ \cites{MR509388}{MR863137}) and connections to other important algebraic structures such as $ \mathrm{t} $-structures.
	The interested reader is encouraged to consult \cite[\SAGsec{7.2}]{SAG} as well as Antieau and Elmanto's recent work \cite{MR4205113}.
\end{remark}

%-------------------------------------------------------------------%
%-------------------------------------------------------------------%
%  Decomposing sheaves on manifolds                                 %
%-------------------------------------------------------------------%
%-------------------------------------------------------------------%

\subsection{Decomposing sheaves on manifolds}\label{subsec:fracture}

We now apply the framework of recollements introduced in \cref{subsec:recollement} to the case where
\begin{equation*}
	\Xcat = \Sh(\Mfld;\Sp) \andeq \Zcat = \Shhi(\Mfld;\Sp) \period
\end{equation*}
Since we can do so at no extra cost, we'll work in the more general setting of sheaves valued in a presentable stable \category. 
First, let's align our notation with \Cref{prop:orthogonaladjoints}.

\begin{nul}
	Let $ C $ be a presentable stable \category.
	Writing $ \Xcat = \Sh(\Mfld;C) $ and $ \Zcat = \Shhi(\Mfld;C) $, in the notation of \Cref{prop:orthogonaladjoints} we have $ \iupperstar = \Lhi $ and $ \iuppershriek = \Rhi $.
\end{nul}

\begin{definition}\label{def:pursheaf}
	Let $ C $ be a stable presentable \category.
	A sheaf $ \Ehat \colon \fromto{\Mfldop}{C} $ is \textit{pure} if $ \Ehat $ is right orthogonal to $ \Shhi(\Mfld;C) $. 
	We write
	\begin{equation*}
		\Shpure(\Mfld;C) \colonequals \Shhi(\Mfld;C)^{\orthogonal} \subset \Sh(\Mfld;C)
	\end{equation*}
	for the full subcategory spanned by the pure sheaves.
\end{definition}

\begin{observation}
	Recall that the subcategory $ \Shhi(\Mfld;C) $ is the essential image of the constant sheaf functor $ \Gammaupperstar \colon \incto{C}{\Sh(\Mfld;C)} $ (\Cref{prop:Dugger}).
	Let $ X \in C $ and $ \Ehat \in \Sh(\Mfld;C) $.
	Then 
	\begin{equation*}
		\Map_{\Sh(\Mfld;C)}(\Gammaupperstar(X),\Ehat) \equivalent \Map_C(X,\Gammalowerstar(\Ehat)) \period
	\end{equation*}
	Thus $ \Ehat $ is right orthogonal to $ \Shhi(\Mfld;C) $ if and only if 
	\begin{equation*}
		\Gammalowerstar(\Ehat) = \Ehat(*) = 0 \period
	\end{equation*}
	Said differently, $ \Shpure(\Mfld;C) $ is the kernel of the constant sheaf functor $ \Gammalowerstar \colon \fromto{\Sh(\Mfld;C)}{C} $.

	Also note that since the global sections functor $ \Gammalowerstar $ preserves all limits and colimits, the subcategory of pure sheaves is stable under limits and colimits.
\end{observation}

Now we introduce the left adjoint to the inclusion $ \Shpure(\Mfld;C) \subset \Sh(\Mfld;C) $ following the prescription of \enumref{prop:orthogonaladjoints}{1}.
In the following, we think of $ \Rhi(\Ehat) $ as playing the role of cohomology with coefficients in $ \RR/\ZZ $ in the differential cohomology hexagon (\Cref{thm:SimonsSullivanunique}).

\begin{definition}\label{differential_cycles}
	Let $ C $ be a stable presentable \category.
	Define a functor
	\begin{equation*}
		\Cyc \colon \fromto{\Sh(\Mfld;C)}{\Sh(\Mfld;C)}
	\end{equation*}
	and a \textit{curvature} natural transformation $ \curv \colon \fromto{\id{}}{\Cyc} $ by the cofiber sequence
	\begin{equation*}
		\begin{tikzcd}[sep=2em]
			\Rhi \arrow[r, "\varepsilon"] & \id{} \arrow[r, "\curv"] & \Cyc \comma
		\end{tikzcd}
	\end{equation*} 
	where $ \varepsilon \colon \fromto{\Rhi}{\id{}} $ is the counit.
	For a $ C $-valued sheaf $ \Ehat $ on $ \Mfld $, we call $ \Cyc(\Ehat) $ the sheaf of \textit{differential cycles} associated to $ \Ehat $.
\end{definition}

\begin{nul}
	As a consequence of \Cref{prop:orthogonaladjoints}, $ \Cyc $ factors through $ \Shpure(\Mfld;C) $ and is left adjoint to the inclusion $ \Shpure(\Mfld;C) \subset \Sh(\Mfld;C) $.
\end{nul}

\begin{observation}
	Since the global sections functor $ \Gammalowerstar $ preserves all limits and colimits, the subcategory of pure sheaves is stable under \textit{both} limits and colimits.
	Since $ \Shpure(\Mfld;C) $ is presentable, the inclusion $ \incto{\Shpure(\Mfld;C)}{\Sh(\Mfld;C)} $ also admits a right adjoint.
\end{observation}

To do this, we identify the left adjoint to the functor $ \Cyc \colon \fromto{\Sh(\Mfld;C)}{\Shpure(\Mfld;C)} $.

\begin{definition}
	Let $ C $ be a stable presentable \category.
	Define a functor
	\begin{equation*}
		\Def \colon \fromto{\Sh(\Mfld;C)}{\Sh(\Mfld;C)}
	\end{equation*}
	by the fiber sequence
	\begin{equation*}
		\begin{tikzcd}[sep=1.5em]
			\Def \arrow[r] & \id{} \arrow[r, "\eta"] & \Lhi \comma
		\end{tikzcd}
	\end{equation*} 
	where $ \eta \colon \fromto{\id{}}{\Lhi} $ is the unit.
	For a $ C $-valued sheaf $ \Ehat $ on $ \Mfld $, we call $ \Def(\Ehat) $ the sheaf of \textit{differential deformations} associated to $ \Ehat $.
\end{definition}

\begin{observations}\label{obs:Acomposite}
	In light of \Cref{thm:stabrecollement}, the functor 
	\begin{equation*}
		\Def \colon \fromto{\Shpure(\Mfld;C)}{\Sh(\Mfld;C)}
	\end{equation*}
	is left adjoint to the functor $ \Cyc $.
	In particular, $ \Def \colon \fromto{\Shpure(\Mfld;C)}{\Sh(\Mfld;C)} $ is fully faithful (\Cref{lem:ladjradj}).
	% \hfill
	% \begin{enumerate}[label=\stlabel{obs:Acomposite}, ref=\arabic*]
	% 	\item Since $ \Lhi $ is idempotent and exact, we see that $ \Lhi \of \Def \equivalent 0 $.

	% 	\item Since $ \Lhi \Rhi \equivalent \Rhi $, for any $ C $-valued sheaf $ \Ehat $ on $ \Mfld $, the fiber sequence defining $ \Def $ gives a fiber sequence
	% 	\begin{equation*}
	% 		\begin{tikzcd}[sep=1.5em]
	% 			\Def\Rhi(\Ehat) \arrow[r] & \Rhi(\Ehat) \arrow[r, "\sim"{yshift=-0.25em}] & \Lhi\Rhi(\Ehat) \comma
	% 		\end{tikzcd}
	% 	\end{equation*} 
	% 	hence $ \Def \of \Rhi \equivalent 0 $.  
	% \end{enumerate}
\end{observations}

\begin{nul}\label{nul:diffcohadjunctions}
	We have chains of adjunctions
	\begin{equation*}
		\begin{tikzcd}[sep=3.5em]
			\Shhi(\Mfld;C) \arrow[r, description, hooked] & \Sh(\Mfld;C) \arrow[l, shift left=1.25ex, "\Rhi"] \arrow[l, shift right=1.25ex, "\Lhi"'] \arrow[r, "\Cyc" description]  & \Shpure(\Mfld;C) \period \arrow[l, shift left=1.25ex, hooked'] \arrow[l, shift right=1.25ex, hooked', "\Def"']
		\end{tikzcd}
	\end{equation*}
	To align notation with \cref{nul:stabelrecadjunction}, we have $ \Xcat = \Sh(\Mfld;C) $, $ \Zcat = \Shhi(\Mfld;C) $, and $ \Ucat = \Shpure(\Mfld;C) $.
	The functors $ \ilowerstar \colon \incto{\Zcat}{\Xcat} $ and $ \jlowerstar \colon \incto{\Ucat}{\Xcat} $ are the two unlabeled inclusions.
	We also have $ \iuppershriek = \Rhi $, $ \iupperstar = \Lhi $, $ \jupperstar = \Cyc $, and $ \jlowershriek = \Def $.
\end{nul}

%-------------------------------------------------------------------%
%  The differential cohomology hexagon                              %
%-------------------------------------------------------------------%

\subsubsection{The differential cohomology hexagon}\label{subsec:diffcohdiagram}

Now we explain how the extended fracture diagram of a stable recollement (\Cref{thm:fracturehexagon}) gives rise to
a ``differential cohomology hexagon''.

\begin{notation}
	We write $ \dup \colon \fromto{\Def}{\Cyc} $ for the composite 
	\begin{equation*}
		\begin{tikzcd}[sep=2em]
			\dup \colon \Def \arrow[r] & \id{} \arrow[r, "\curv"] & \Cyc \period
		\end{tikzcd}
	\end{equation*} 
\end{notation}

\begin{corollary}[(fracture square)]\label{thm:fracturesquare}
	Let $ C $ be a stable presentable \category.
	The \category $ \Sh(\Mfld;C) $ is the recollement of the subcategories $ \Shhi(\Mfld;C) $ and $ \Shpure(\Mfld;C) $.
	In particular, there is a commutative diagram
	\begin{equation}\label{diag:fiberseqs}
		\begin{tikzcd}[sep=2.5em]
			 &  \Def \arrow[r, equals] \arrow[d] & \Def \arrow[d, "\d"] \\
			\Rhi \arrow[r] \arrow[d, equals] & \id{\Sh(\Mfld;C)} \arrow[r] \arrow[d] \arrow[dr, phantom, "\square" description] & \Cyc \arrow[d] \\
			\Rhi \arrow[r] & \Lhi \arrow[r] & \Lhi\Cyc 
		\end{tikzcd}
	\end{equation}
	of functors $ \fromto{\Sh(\Mfld;C)}{\Sh(\Mfld;C)} $, where the lower right-hand square is a pullback and all rows and columns are fiber sequences.
\end{corollary}

\begin{nul}
	Informally, $ \Sh(\Mfld;C) $ is the \category of triples 
	\begin{equation*}
		(\Ehat_{\RR}, \Ehat_{\pure}, \phi \colon \fromto{\Ehat_{\RR}}{\Lhi\Ehat_{\pure}}) \comma
	\end{equation*}
	where $ \Ehat_{\RR} $ is a \RRinvariant sheaf, $ \Ehat_{\pure} $ is a pure sheaf, and $ \phi $ is any morphism.
\end{nul}

\begin{nul}[(differential cohomology hexagon)]\label{nul:spectraldiffcohdiag}
	With some rearrangement, \Cref{thm:fracturesquare} and the fact that pullback squares compose, we see that there is a diagram of pullback squares
	\begin{equation}\label{diag:firstdiffcoh}
		\begin{tikzcd}[sep=2.5em]
			\Sigma^{-1} \Lhi\Cyc \arrow[r] \arrow[d] \arrow[dr, phantom, "\square" description] & \Rhi \arrow[r] \arrow[d] \arrow[dr, phantom, "\square" description] & 0 \arrow[d] \\
			\Def \arrow[r] \arrow[d] \arrow[dr, phantom, "\square" description] & \id{\Sh(\Mfld;C)} \arrow[r] \arrow[d] \arrow[dr, phantom, "\square" description] & \Cyc \arrow[d] \\
			0 \arrow[r] & \Lhi \arrow[r] & \Lhi\Cyc \period
		\end{tikzcd}
	\end{equation}
	Rearranging the diagram \eqref{diag:firstdiffcoh}, for each $ \Ehat \in \Sh(\Mfld;C) $ we get the following
	\textit{differential cohomology hexagon}
	\begin{equation}\label{diag:generaldiffcohhex}
		\begin{tikzcd}[column sep={10ex,between origins}, row sep={8ex,between origins}]
			& \Rhi(\Ehat) \arrow[rr] \arrow[dr] & & \Lhi(\Ehat) \arrow[dr] & \\
			\Sigma^{-1} \Lhi \Cyc(\Ehat) \arrow[ur] \arrow[dr] & & \Ehat \arrow[ur] \arrow[dr, "\curv" description] & & \Lhi \Cyc(\Ehat) \\
			& \Def(\Ehat) \arrow[rr, "\dup"'] \arrow[ur] & & \Cyc(\Ehat) \arrow[ur]  & \phantom{\Lhi \Cyc(\Ehat)} \period
		\end{tikzcd}
	\end{equation}
	Here the diagonals are fiber sequences, the top and bottom rows are extensions of fiber sequences by one term, and both squares are pullback squares.
	The ``top row'' consists of \RRinvariant sheaves, whereas the ``bottom row'' consists of sheaves that are, in some sense, more geometric.

	Since $ \Lhi \equivalent \Gammaupperstar \Gammalowersharp $ and $ \Rhi \equivalent \Gammaupperstar \Gammalowerstar $ \Cref{nul:Gammaadjunctions}, the differential cohomology hexagon \eqref{diag:generaldiffcohhex} can be rewritten as
	\begin{equation*}\label{diag:generaldiffcohhexGamma}
		\begin{tikzcd}[column sep={10ex,between origins}, row sep={8ex,between origins}]
			& \Gammaupperstar\Gammalowerstar\Ehat\arrow[rr]\arrow[dr] & & \Gammaupperstar \Gammalowersharp \Ehat \arrow[dr] & \\
			\Sigma^{-1}\Gammaupperstar\Gammalowersharp\Cyc(\Ehat)\arrow[dr]\arrow[ur] & & \Ehat\arrow[dr]\arrow[ur] & & \Gammaupperstar\Gammalowersharp \Cyc(\Ehat)\\
			& \Def(\Ehat)\arrow[ur]\arrow[rr, "\d"'] & & \Cyc(\Ehat)\arrow[ur] & \phantom{\Gammaupperstar\Gammalowersharp \Cyc(\Ehat)} \period
		\end{tikzcd}
	\end{equation*}
\end{nul}

%-------------------------------------------------------------------%
%  Differential refinements                                         %
%-------------------------------------------------------------------%

\subsubsection{Differential refinements}\label{subsec:refinements}

We finish this section by making precise what it means for a differential cohomology theory $ \Ehat \in \Sh(\Mfld;\Sp) $ to refine a cohomology theory $ E \in \Sp $.

\begin{definition}\label{def:refinement}
	Let $ C $ be a presentable stable \category.
	A \textit{differential refinement} of a an object $ E \in C $ is pair $ (\Ehat,\phi) $ of a sheaf $ \Ehat \in \Sh(\Mfld;C) $ together with an equivalence $ \phi \colon \equivto{\Gammalowersharp(\Ehat)}{E} $ in $ C $.
\end{definition}

\begin{nul}\label{nul:altdiffrefinement}
	From the fracture square (\Cref{thm:fracturesquare}), a differential refinement of $ E \in C $ is equivalently the data of a pure sheaf $ \Phat \in \Shpure(\Mfld;C) $ along with a morphism $ \fromto{E}{\Gammalowersharp(\Phat)} $ in $ C $.
	Given this data, we can construct a differential refinement $ \Ehat $ in the sense of \Cref{def:refinement} as the pullback
	\begin{equation*}
		\begin{tikzcd}[column sep={11ex,between origins}, row sep={11ex,between origins}]
			\Ehat \arrow[r] \arrow[d] \arrow[dr, phantom, "\square" description] & \Phat \arrow[d] & \\
			\Gammaupperstar(E) \arrow[r] & \Gammaupperstar\Gammalowersharp(\Phat) \period
		\end{tikzcd}
	\end{equation*}

	In this case, we have:
	\begin{enumerate}[label=\stlabel{nul:altdiffrefinement}, ref=\arabic*]
		\item $ \equivto{\Def(\Ehat)}{\Def(\Phat)} $.

		\item $ \equivto{\Cyc(\Ehat)}{\Phat} $.

		\item $ \Gammalowerstar(\Ehat) $ fits into a fiber sequence
		\begin{equation*}
			\begin{tikzcd}[sep=1.5ex]
				\Gammalowerstar(\Ehat) \arrow[r] & E \arrow[r] & \Gammalowersharp(\Phat) \period
			\end{tikzcd}
		\end{equation*}
	\end{enumerate}
\end{nul}

\begin{construction}[(pullback of a differential refinement)]\label{constr:pullbackrefinement}
	Let $ C $ be a presentable stable \category, $ f \colon \fromto{E}{E'} $ a morphism in $ C $, and $ (\Ehat',\phi') $ a differential refinement of $ E' $. 
 	Form the pullback
	\begin{equation}\label{sq:differential}
		\begin{tikzcd}[column sep={11ex,between origins}, row sep={11ex,between origins}]
			\Ehat \arrow[r, "\fhat"] \arrow[d] \arrow[dr, phantom, "\square" description] & \Ehat' \arrow[d] & \\
			\Gammaupperstar(E) \arrow[r, "\Gammaupperstar(f)"'] & \Gammaupperstar(E') \comma
		\end{tikzcd}
	\end{equation}
	where the morphism $ \fromto{\Ehat'}{\Gammaupperstar(E')} $ is adjoint to the given equivalence $ \phi' \colon \equivto{\Gammalowersharp(\Ehat')}{E'} $.
	Since $ \Gammalowersharp $ is exact, applying $ \Gammalowersharp $ to the square \eqref{sq:differential} gives a pullback square
	\begin{equation*}
		\begin{tikzcd}[column sep={11ex,between origins}, row sep={11ex,between origins}]
			\Gammalowersharp(\Ehat) \arrow[r] \arrow[d, "\phi"'] \arrow[dr, phantom, "\square" description] & \Gammalowersharp(\Ehat') \arrow[d, "\phi'"', "\sim"{sloped, yshift=-0.2em}] & \\
			E \arrow[r, "f"'] & E' \comma
		\end{tikzcd}
	\end{equation*}
	which provides an equivalence $ \phi \colon \equivto{\Gammalowersharp(\Ehat)}{E} $.
	The \textit{pullback differential refinement} of $ (\Ehat',\phi') $ along $ f $ is the differential refinement $ (\Ehat,\phi) $ of $ E $.
\end{construction}

\begin{lemma}\label{lem:pullbackrefinements}
	In the notation of \Cref{constr:pullbackrefinement}, the following
	\begin{enumerate}[label=\stlabel{lem:pullbackrefinements}, ref=\arabic*]
		\item The morphism $ \Def(\fhat) \colon \fromto{\Def(\Ehat)}{\Def(\Ehat')} $ is an equivalence.

		\item The morphism $ \Cyc(\fhat) \colon \fromto{\Cyc(\Ehat)}{\Cyc(\Ehat')} $ is an equivalence.

		\item The global sections of $ \Ehat $ is given by the pullback
		\begin{equation*}
			\begin{tikzcd}[column sep={11ex,between origins}, row sep={11ex,between origins}]
				\Gammalowerstar(\Ehat) \arrow[r, "\Gammalowerstar(\fhat)"] \arrow[d] \arrow[dr, phantom, "\square" description] & \Gammalowerstar(\Ehat') \arrow[d] & \\
				E \arrow[r, "f"'] & E' \period
			\end{tikzcd}
		\end{equation*}
	\end{enumerate} 
\end{lemma}

\newpage
%!TEX root = ../diffcoh.tex

\section{Examples}\label{sec:examples}
\textit{by Araminta Amabel}

The purpose of this section is to construct examples of differential cohomology theories, i.e., sheaves of spectra on the category $ \Mfld $.
We'll construct these examples by using the method of \textit{differential refinements} introduced in \cref{subsec:refinements}.
Note that given a spectrum $E$, there are possibly many differential refinements of $E$. 
We will construct differential cohomology theories refining the cohomology theory $ E $ by the following process:
\begin{enumerate}
	\item Choose a pure sheaf $\Phat $ (\Cref{def:pursheaf}).

	\item Compute $\Gammalowersharp\Phat $ using the formula $ \Gammalowersharp\Phat  = \colim_{\Deltaop} \Phat (\Deltaalgdot) $ of \Cref{cor:formula_for_Gammalowersharp}.

	\item Find a map of spectra $ f \colon \fromto{E}{\Gammalowersharp\Phat} $.

	\item Define $\Ehat$ as in the pullback
	\begin{equation*}
		\begin{tikzcd}[column sep={15ex,between origins}, row sep={11ex,between origins}]
			\Ehat \arrow[r] \arrow[d] \arrow[dr, phantom, "\square" description] & \Phat \arrow[d] & \\
			\Gammaupperstar(E) \arrow[r, "\Gammaupperstar(f)"'] & \Gammaupperstar\Gammalowersharp(\Phat) \period
		\end{tikzcd}
	\end{equation*} 
\end{enumerate}

We start in \cref{subsec:moststimple} with differential refinements of $ 0 $ and what the differential cohomology hexagon looks like in this case.
In \cref{subsec:simplefilt}, we refine this simplest example by adding a filtration.
\Cref{subsec:CheegerSimonsDiffchar} explains how the Cheeger--Simons theory of differential characters fits into this story, and \cref{subsec:diffKtheory} studies differential refinements of \Ktheory.

%-------------------------------------------------------------------%
%-------------------------------------------------------------------%
%  The simplest example                                             %
%-------------------------------------------------------------------%
%-------------------------------------------------------------------%

\subsection{The simplest example}\label{subsec:moststimple}

To start off, let's try to construct a differential refinement where the pure sheaf $ \Phat $ is zero.
That is, $ \Phat = 0 = \Gammaupperstar  0 $. 
In this case, since the functor $ \Gammalowersharp $ is exact, $ \Gammalowersharp(\Phat) = 0 $.
Any spectrum $E$ maps uniquely to $ 0 $.
Thus for any spectrum $ E $ we have a differential refinement $ \Ehat $ defined by the pullback
\begin{equation*}
	\begin{tikzcd}
		\Ehat\arrow[r]\arrow[d] & \Gammaupperstar E\arrow[d]\\
		\Gammaupperstar 0\arrow[r, equals] & \Gammaupperstar 0
	\end{tikzcd}
\end{equation*}
Since the bottom horizontal arrow is an equivalence, the top horizontal arrow is as well: $ \Ehat = \Gammaupperstar E $. 
The rest of the differential cohomology diagram looks as follows
\begin{equation*}
	\begin{tikzcd}[column sep={10ex,between origins}, row sep={8ex,between origins}]
		& \Gammaupperstar E\arrow[rr]\arrow[dr, equals] & & \Gammaupperstar E\arrow[dr] & \\
		0 \arrow[dr]\arrow[ur] & & \Gammaupperstar E\arrow[ur, equals]\arrow[dr] & & 0\\
		& \Def(\Gammaupperstar E)\arrow[rr]\arrow[ur] & & 0 \arrow[ur] & \phantom{\Gammaupperstar E} \period
	\end{tikzcd}
\end{equation*}
Since the upwards diagonal sequence is a fiber sequence, we also have $\Def(\Gammaupperstar E) = 0 $.

This example is just saying that that $ E $-cohmology is a special case of differential cohomology.
We're really just reformulating the fact that the constant sheaf functor $ \Gammaupperstar \colon \fromto{\Sp}{\Sh(\Mfld;\Sp)} $ is fully faithful with essential image the \RRinvariant sheaves (\Cref{prop:Dugger}). 

%-------------------------------------------------------------------%
%-------------------------------------------------------------------%
%  The simplest example, but with a filtration                      %
%-------------------------------------------------------------------%
%-------------------------------------------------------------------%

\subsection{The simplest example, but with a filtration}\label{subsec:simplefilt} 

We give an alternative differential refinement of the zero spectrum which comes with a natural filtration.

\begin{nul}
	Let $\Omegabullet\in\Shv(\Mfld;\D(\RR))$ the sheaf of de Rham forms with cohomological grading; so $ \Omega^k $ is in degree $ -k $. 
	Consider the resulting functor of spectra, $\HOmegabullet$. 
	By the Poincaré Lemma, $\Omegabullet$ is quasi-isomorphic to the constant sheaf at $\RR[0]$. 
	Thus $ \equivto{\HOmegabullet}{\Gammaupperstar\HRR} $. 
	In particular, $\HOmegabullet$ is not pure. 
\end{nul}

However, since $ \HOmegabullet \equivalent \Gammaupperstar\HRR $ is \RRinvariant, the purification $
\Cyc(\HOmegabullet) $ is equivalent to zero.
Now $\Omegabullet$ has a filtration by degree. 
For $k\in\NN$, let $\Omega^{\geq k}$ denote the stunted piece of the chain complex $\Omegabullet$ where we have replaced everything in degrees $<k$ by $0$. 
We get induced filtrations of $\HOmegabullet$ and of $\Cyc(\HOmegabullet)\simeq\Gammaupperstar 0$.

For $k\geq 1$, there is an equivalence $\Omega^{\geq k}(*)\simeq 0$ of chain complexes. 
Thus the global sections of $\HOmega^{\geq k}$ is 0,
\begin{equation*}
	\Gammalowerstar \HOmega^{\geq k}= \HOmega^{\geq k}(*) = 0 \period
\end{equation*}
By definition, this means that $\HOmega^{\geq k}$ is a pure sheaf if (and only if) $k\geq 1$. 
The purification functor $\Cyc$ is the identity on pure sheaves, so we obtain a filtration of the pure sheaf $\Gammaupperstar 0$ by pure sheaves
\begin{equation*}
	\Gammaupperstar 0 \to \HOmega^{\geq 1}\to\cdots\to \HOmega^{\geq k}\to\cdots \period
\end{equation*}
Now for each $k\geq 1$, we can choose the pure sheaf $\HOmega^{\geq k}$ and follow our procedure.

We need to compute the homotopification of our chosen pure sheaf.

\begin{lemma}\label{2.1}
	For any $k\in\NN$, there is an equivalence $\Gammalowersharp \HOmega^{\geq k}\simeq \HRR$.
\end{lemma}

\begin{proof}
	For $k=0$ ,we have seen that $\HOmega^{\geq 0}\simeq \Gammaupperstar \HRR$, which is already homotopy invariant. Thus $\Gammalowersharp\Gammaupperstar \HRR\simeq \HRR$. 
	For $k\geq 1$, see \cite[Lemma 7.15]{MR3462099}.
\end{proof}

The following family of differential refinements was introduced by Hopkins and Singer, \cite{HopkinsSinger}.

\begin{definition}
	Let $ E $ be a spectrum and $ f \colon \fromto{E}{\HRR} $ a map of spectra.
	For each $ k \geq 1 $, write $ \Ehat(k) $ for the pullback 
	\begin{equation*}
		\begin{tikzcd}
			\Ehat(k)\arrow[r]\arrow[d] & \Gammaupperstar E\arrow[d, "f"]\\
			\HOmega^{\geq k}\arrow[r] & \HRR \period
		\end{tikzcd}
	\end{equation*}
\end{definition}

The differential cohomology diagram \eqref{diag:generaldiffcohhexGamma} for $\Ehat(k)$ looks like
\begin{equation*}
	\begin{tikzcd}[column sep={10ex,between origins}, row sep={8ex,between origins}]
		& \Gammaupperstar \Gammalowerstar\Ehat\arrow[rr]\arrow[dr] & & \Gammaupperstar E\arrow[dr] & \\
		\Sigma^{-1}\Gammaupperstar \HRR\arrow[dr]\arrow[ur] & & \Ehat(k)\arrow[dr]\arrow[ur] & & \Gammaupperstar \HRR\\
		& \HOmega^{\leq k-1}[-1]\arrow[rr] \arrow[ur] & & \HOmega^{\geq k}\arrow[ur] & \phantom{\Gammaupperstar \HRR} \period
	\end{tikzcd}
\end{equation*}

%-------------------------------------------------------------------%
%-------------------------------------------------------------------%
%  Cheeger–Simons Differential Characters                           %
%-------------------------------------------------------------------%
%-------------------------------------------------------------------%

\subsection{Ordinary differential cohomology}\label{subsec:CheegerSimonsDiffchar}

Take $E=\HZZ$ and the map $\HZZ\to \HRR$ induced from the inclusion $\ZZ\subset\RR$.

\begin{definition}
	The $k$-th \emph{ordinary differential cohomology group} of a manifold $M$, denoted $\Hhat^k(M)$ is the $(-k)$-th homotopy group
	\begin{equation*}
		\Hhat^k(M)=\uppi_{-k}\HZZhat(k)(M)
	\end{equation*}
	where $\HZZhat(k)$ is defined by the homotopy pullback square
	\begin{equation*}
		\begin{tikzcd}
			\HZZhat(k)\arrow[r]\arrow[d] & \Gammaupperstar  \HZZ\arrow[d]\\
			\Cyc(\HOmega^{\geq k})\arrow[r] & \Gammaupperstar \HRR \period
		\end{tikzcd}
	\end{equation*}
\end{definition}

\begin{nul}
	Note that $\Cyc(\HOmega^{\geq k})\simeq \HOmega^{\geq k}$ if $k\geq 1$ and is $\HRR$ if $k=0$.
\end{nul}

\begin{remark}
	The group $\Hhat^k(M)$ is also known as the \emph{Cheeger--Simons differential characters}, or the \emph{smooth Deligne cohomology}. 
\end{remark}

The following gives an explicit complex computing ordinary differential cohomology.
This complex first appeared in the setting of complex manifolds in Deligne's work on Hodge theory (see \cites[\S2.2]{MR498551}[\S12.3]{MR2451566}), and is why differential cohomology is also called smooth Deligne cohomology.

\begin{lemma}\label{lem:Delignemodel}
	Let $k\geq 1$. 
	The sheaf of spectra $\HZZhat(k)$ is given by applying the Eilenberg--MacLane functor $ \EM \colon \fromto{\D(\ZZ)}{\Sp} $ (\Cref{rec:EMspectra}) pointwise to the sheaf of chain complexes
	\begin{equation*}
		(\Gammaupperstar \ZZ\to\Omega^0\to\Omega^1\to\cdots\to\Omega^{k-1}\to 0\to \cdots) \period
	\end{equation*}
	Here $\Omega^i$ is in degree $-i-1$. 
	Moreover, the group $\Hhat^k(M)$, for a manifold $M$, can be computed as the $k$-th sheaf cohomology group of this sheaf of chain complexes.
\end{lemma}

\begin{proof}
	By construction, $\HZZhat(k)$ comes from applying $\EM$ of the sheaf of chain complexes $F$ given by the homotopy pullback
	\begin{equation*}
		\begin{tikzcd}
			F\arrow[r]\arrow[d] & \Gammaupperstar \ZZ[0]\arrow[d]\\
			\Omega^{\geq k}\arrow[r] & \Omegabullet \period
		\end{tikzcd}
	\end{equation*}
	Since the bottom horizontal arrow is an inclusion, its cofiber is given by the cokernel. 
	We have a cofiber sequence in $ \D(\ZZ) $
	\begin{equation*}
		\Omega^{\geq k}\to\Omegabullet\to\Omega^{\leq k-1}
	\end{equation*}
	where $\Omega^{\leq k-1}$ has $\Omega^i$ in degree $-i$, and $0$ above $k-1$. 
	The cofiber of the top horizontal map is equivalent to the cofiber of the bottom horizontal map. 
	Since we are in a stable setting, these cofiber sequences are also fiber sequences. 
	Thus, we have a fiber sequence
	\begin{equation*}
		F\to\Gammaupperstar \ZZ[0]\to\Omega^{\leq k-1}
	\end{equation*}
	where $\ZZ[0]\to\Omega^{\leq k-1}$ includes $\ZZ$ and $\Omega^0$. 
	The fiber of this inclusion is a shift of the mapping cone, which is
	\begin{equation*}
		(\Gammaupperstar \ZZ\to\Omega^0\to\Omega^1\to\cdots\to\Omega^{k-1}\to 0\to\cdots)
	\end{equation*}
	Finally, note that $ \uppi_{-k}(\EM F) = \H^k(F) $.
\end{proof}

\begin{example}
	Take $k=0$. 
	Then $\HZZhat(k)\simeq\Gammaupperstar \HZZ$ and 
	\begin{equation*}
		\Gammaupperstar \HZZ(M) = \Hom_{\Sp}(\Sigma_{+}^{\infty} \Piinf(M),\HZZ)
	\end{equation*}
	(\Cref{ex:Sptcotensor,lem:constantishi}).
	Hence $ \Gammaupperstar \HZZ(M) $ has $ 0 $-th homotopy group $\H^0(M;\ZZ)$.	
\end{example}

The following two computations from Kumar's notes \cite{Nilay}.

\begin{example}
	Take $k=1$. We compute $\Hhat^1(M)$. 
	By \Cref{lem:Delignemodel}, we can compute $\Hhat^1(M)$ as the $ 1 $-st sheaf cohomology group of the sheaf of chain complexes $(\Gammaupperstar \ZZ\to\Omega^0)$. 
	After choosing a good cover of $M$, we can compute this sheaf cohomology as Čech cohomology. 
	The Čech cohomology will be the cohomology of the total complex of the following bicomplex,
	\begin{equation*}
		\begin{tikzcd}
			\Cech^0(\Gammaupperstar \ZZ)\arrow[r]\arrow[d] & \Cech^0(\Omega^0)\arrow[d]\\
			\Cech^1(\Gammaupperstar \ZZ)\arrow[r]\arrow[d] & \Cech^1(\Omega^0)\arrow[d]\\
			\Cech^2(\Gammaupperstar \ZZ)\arrow[r]\arrow[d] & \Cech^1(\Omega^0)\arrow[d]\\
			\vdots & \vdots
		\end{tikzcd}
	\end{equation*}
	with $\Cech^i(\Gammaupperstar \ZZ)$ in bidegree $(0,-i)$ and $\Cech^i(\Omega^0)$ in bidgree $(-1,-i)$. 
	The differential on this bicomplex is
	\begin{equation*}
		D = d^{\mrm{hor}} + (-1)^p d^{\mrm{ver}} \comma
	\end{equation*}
	where $ p $ is the horizontal degree. 
	The piece of the total complex that we are interested looks like
	\begin{equation*}
		\begin{tikzcd}
			\Cech^0(\Gammaupperstar \ZZ) \arrow[r, "D_0"] & \Cech^0(\Omega^0)\oplus\Cech^1(\Gammaupperstar \ZZ) \arrow[r, "D_1"] & \Cech^1(\Omega^0)\oplus\Cech^2(\Gammaupperstar \ZZ) \period 
		\end{tikzcd}
	\end{equation*}
	If our good cover of $M$ is $\{U_\alpha\}$ with intersections $U_{\alpha\beta}$, then an element of $\Cech^0(\Omega^0)\oplus\Cech^1(\Gammaupperstar \ZZ)$ looks like a collection of smooth maps $f_\alpha\colon U_\alpha\to\RR$ and integers $n_{\alpha\beta}\in\ZZ$. The map $D_1$ sends
	\begin{equation*}D_1(f_\alpha, n_{\alpha\beta})=(f_\alpha-f_\beta+n_{\alpha\beta},n_{\beta\gamma}-n_{\alpha\gamma}+n_{\alpha\beta})\end{equation*}
	In particular, an element of $\ker D_1$ consists of maps $f_\alpha$ that agree on intersections up to an integer. These glue together to give a (smooth) map $f\colon M\to \Circ=\Uup_{1}$.

	The map $D_0$ sends a collection $(n_\alpha)$ to
	\begin{equation*}D_0(n_\alpha)= (c_{n_\alpha},n_\alpha-n_\beta)\end{equation*}
	where $c_{n_\alpha}$ is the constant function $U_\alpha\to\RR$ at the integer $n_\alpha$. As a map $M\to \Circ$, these glue together to the constant map at the base point.

	Thus we have an isomorphism
	\begin{equation*}
		\Hhat^1(M)\isomorphic\Mapsm(M,\Uup_{1}) \period
	\end{equation*}
	In ordinary cohomology, we have
	\begin{equation*}
		\H^1(M;\ZZ)= \uppi_{0} \Map_{\Spc}(M,\Kup(\ZZ,1)) = \uppi_{0} \Map_{\Spc}(M,\Uup_{1}) \period
	\end{equation*}
	In this sense, differential cohomology replaced homotopy maps with smooth maps.
\end{example}

\begin{example}
	Take $k=2$. 
	Then we have an isomorphism
	\begin{equation*}
		\Hhat^2(M) \isomorphic \{\text{line bundles on $M$ with connection}\}/\!\sim \period
	\end{equation*}
	In ordinary cohomology, we have
	\begin{align*}
		\H^2(M;\ZZ) &= \uppi_{0} \Map_{\Spc}(M,\Kup(\ZZ,2)) \\
		&= \uppi_{0} \Map_{\Spc}(M,\BU_{1}) \\
		&= \{\text{line bundles on $M$}\}/\!\sim \period 
	\end{align*}
	In this sense, the new geometric information encoded in differential cohomology is the connection.
\end{example}

%-------------------------------------------------------------------%
%-------------------------------------------------------------------%
%  Differential K-theory                                            %
%-------------------------------------------------------------------%
%-------------------------------------------------------------------%

\subsection{Differential \texorpdfstring{$\Kup$}{K}-theory}\label{subsec:diffKtheory}

\begin{nul}
	Consider de Rham forms with $\CC[u^{\pm 1}]$ coefficients, with $u$ in degree $2$. 
	We obtain a family of pure sheaves $\HOmega^{\geq k}(-;\CC[u^{\pm 1}])$. 
	As in \Cref{2.1}, we have an equivalence,
	\begin{equation*}
		\Gammalowersharp\HOmega^{\geq k}(-;\CC[u^{\pm 1}])\simeq \HCC[u^{\pm 1}]
	\end{equation*}
\end{nul}

\begin{nul}
	Take $E=\ku$ to be the spectrum defining connective complex \Ktheory. 
	The \textit{Chern character} defines a map of spectra
	\begin{equation*}
		\ch \colon \ku \to \HCC[u^{\pm 1}] \period
	\end{equation*}
	The resulting family of differential cohomology theories defined by pullback squares,
	\begin{equation*}
		\begin{tikzcd}
			\kuhat(k) \arrow[r] \arrow[d] & \Gammaupperstar  \HOmega^{\geq k}(-;\CC[u^{\pm 1}]) \arrow[d]\\
			\Gammaupperstar(\ku) \arrow[r, "\Gammaupperstar(\ch)"'] & \Gammaupperstar \HCC[u^{\pm 1}] \period
		\end{tikzcd}
	\end{equation*}
	first studied by Hopkins and Singer in \cite{HopkinsSinger} is called \textit{differential \Ktheory}.
\end{nul}

\begin{nul}
	There are other interesting differential refinements of $\ku$ that do not arise from the pure sheaves $\HOmega^{\geq k}(-;\CC[u^{\pm 1}])$.
\end{nul}

\newpage
%!TEX root = ../diffcoh.tex

%-------------------------------------------------------------------%
%-------------------------------------------------------------------%
%  The Deligne cup product                                          %
%-------------------------------------------------------------------%
%-------------------------------------------------------------------%

\section{The Deligne cup product}\label{sec:Delignecup}
\textit{by Araminta Amabel}

Let $M$ be a manifold. 
Recall that the Deligne complex $\ZZ(k)$ is the homotopy pullback
\begin{equation*}
	\begin{tikzcd}
		\ZZ(k)\arrow[r]\arrow[d] & \ZZ\arrow[d] \\
		\Sigma^k\Omegacl^k\arrow[r] & \RR \period
	\end{tikzcd}
\end{equation*}
The goal of this section is to combine the cup product on $\HZZ$ and the wedge product on differential forms to put a ring structure on differential cohomology.

%-------------------------------------------------------------------%
%-------------------------------------------------------------------%
%  Combining the cup and wedge products                             %
%-------------------------------------------------------------------%
%-------------------------------------------------------------------%

\subsection{Combining the cup and wedge products}

\begin{nul}
	Notice that the cup product on $ \HZZ $ and $ \HRR $ and the wedge product on differential forms fit into a commutative digram.
	\begin{equation*}
		\begin{tikzcd}
			\ZZ(n)\otimes\ZZ(m)\arrow[r]\arrow[d] & \HZZ[n]\otimes \HZZ[m]\arrow[r, "\cupprod"] \arrow[d] & \HZZ[m+n]\arrow[dd]\\
			\Omegacl^n\otimes\Omegacl^m \arrow[d, "\wedge"'] \arrow[r] & \HRR[n]\otimes \HRR[m]\arrow[dr, "\cupprod" description] & \\
			\Omegacl^{n+m}\arrow[rr] & & \HRR[m+n] \period 
		\end{tikzcd}
	\end{equation*}
\end{nul}

\begin{nul}
	By the definition of $ \ZZ(k) $ as a pullback, we can represent $\ZZ(k)(M)$ as a triple $(c,h,\omega)$ where $c$ is an integral degree $k$ cocycle on $M$, $\omega$ is a closed $k$ form on $M$, and $h$ is a degree $k-1$ real cochain on $M$ so that $\d x=\omega-c$.
\end{nul}

\begin{nul}
	In particular, if we represent an element of $\Cup^n(M;\ZZ(n)) $ by a triple $(c_1,h_1,\omega_1)$ and an element of $\Cup^m(M;\ZZ(m))$ by a triple $(c_2,h_2,\omega_2)$ we would like the product to be a triple
	\[
		(c_1,h_1,\omega_1)\cupprod(c_2,h_2,\omega_2)=(c_3,h_3,\omega_3)\in \Cup^{m+n}(M;\ZZ(m+n)) \period
	\]
	Saying that this product comes from combining the cup product and the wedge product, means that $c_3=c_1\cupprod c_2$ and $\omega_3=\omega_1\wedge\omega_2$. We are only left with figuring out what $h_3$ should be. 
	Heuristically, $h_3$ should be a homotopy between $c_3$ and $\omega_3$; i.e., a homotopy between the cup product and the wedge product. 
\end{nul}

\begin{nul}
	Given forms $\omega\in\Omega^n(M)$ and $\eta\in\Omega^m(M)$, we can form the wedge product $\omega\wedge \eta\in\Omega^{n+m}(M)$ and view that as a real cochain under the map
	\begin{equation*}
		\Omega^{n+m}(M)\to \Cup^{n+m}(M;\RR) \period
	\end{equation*}
	We could also map the forms $\omega,\eta$ to real cochains on $M$ and then take their cup product. 
	Let $B(\omega,\eta)\in \Cup^{n+m-1}(M;\RR)$ be a choice of natural homotopy between these two cochains so that 
	\[
		\d B(\omega,\eta)+B(d\omega,\eta)+(-1)^{|\omega|} B(\omega,d\eta)=\omega\wedge\eta-\omega\cupprod \eta \period
	\]
	Note that we can take $B(\omega,0)=0$. 
\end{nul}

\begin{nul}
	Then the product of $(c_1,h_1,\omega_1)\in \Cup^n(M;\ZZ(n))$ and $(c_2,h_2,\omega_2)\in \Cup^m(M;\ZZ(m))$ is given by
	\[
		(c_3,h_3,\omega_3)=(c_1\cupprod c_2,(-1)^{|c_1|}c_1\cupprod h_2+h_1\cupprod\omega_2+B(\omega_1,\omega_2),\omega_1\wedge\omega_2) \period 
	\]
	For this to be a differential cocycle, we need to have 
	\[
		\d((-1)^{|c_1|}c_1\cupprod h_2+h_1\cupprod \omega_2+B(\omega_1,\omega_2)=\omega_1\wedge\omega_2-c_1\cupprod c_2 \period 
		\]
	This will only work if $(c_1,h_1,\omega_1)$ and $(c_2,h_2,\omega_2)$ are themselves cocycles; i.e., $\d c_i=0=\d\omega_i$. 
	In this case, we have
	\[
		\omega_1\wedge\omega_2-\omega_1\cupprod\omega_2=\d B(\omega_1,\omega_2)=B(0,\omega_2)+(-1)^{|\omega_1|}B(\omega_1,0)=\d B(\omega_1,\omega_2) \period
	\]
	Thus
	\begin{align*}
		\d\paren{(-1)^{|c_1|}c_1\cupprod h_2+h_1\cupprod \omega_2+B(\omega_1,\omega_2)} &= (-1)^{|c_1|}\d(c_1\cupprod h_2)+\d(h_1\cupprod\omega_2)+\d B(\omega_1,\omega_2) \\
		&=(-1)^{|c_1|}\left(\d c_1\cupprod h_2+(-1)^{|c_1|}c_1\cupprod \d h_2\right)+\d h_1\cupprod\omega_2 \\ 
		&\phantom{=} \qquad +(-1)^{|h_1|}h_1\cupprod \d\omega_2+\d B(\omega_1,\omega_2) \\
		&= c_1\cupprod \d h_2+\d h_1\cupprod\omega_2+dB(\omega_1,\omega_2) \\
		&= c_1\cupprod(\omega_2-c_2)+(\omega_1-c_1)\cupprod\omega_2+\omega_1\wedge\omega_2-\omega_1\cupprod\omega_2 \\
		&= c_1\cupprod\omega_2-c_1\cupprod c_2+\omega_1\cupprod\omega_2-c_1\cupprod\omega_2 \\ 
		&\phantom{=} \qquad +\omega_1\wedge\omega_2-\omega_1\cupprod\omega_2 \\ 
		&=\omega_1\wedge\omega_2-c_1\cupprod c_2 \period
	\end{align*}
\end{nul}

\begin{remark}
	In fact we can get $\Einf$-structure from the homotopy pullback diagram.
	% \todo{It might be worth expanding on this remark more.} 
	View $\HZZ$ as a (trivially) filtered $\Einf$-algebra. 
	View the de Rham complex $ \Omegabullet$ as a filtered $\Einf$-algebra with filtration $ \{\Omega^{\geq k}\}_{k \geq 0} $. 
	Then the homotopy pullback of two $\Einf$-algebras is again an $\Einf$-algebra.
\end{remark}

%-------------------------------------------------------------------%
%-------------------------------------------------------------------%
%  The Deligne cup product                                          %
%-------------------------------------------------------------------%
%-------------------------------------------------------------------%

\subsection{The Deligne cup product}

Recall that we have an identification of the homotopy pullback $\HZZhat(k)$ with the complex of sheaves
\begin{equation*}
	\ZZ(k)=\bigg(
	\begin{tikzcd}[sep=1.5em]
		\Gamma^*\ZZ \arrow[r, "\iota"] & \Omega^0 \arrow[r, "\d"] & \Omega^1 \arrow[r, "\d"] & \cdots \arrow[r, "\d"] & \Omega^{k-1} 
	\end{tikzcd}\bigg) \period
\end{equation*}
Under this identification, we can describe the product in differential cohomology more explicitly. This is
sometimes called the \emph{Deligne cup product}.

Let $M$ be a manifold and $U\subset M$ an open set. 
Then $\ZZ(k)(U)$ is a chain complex that is $\Cup^0(U;\ZZ)$ in degree $ 0 $ and $\Omega^p(U)$ in degree $p+1$. 

\begin{proposition}\label{formula}
	The Deligne cup product
	\begin{equation*}
		\cupprod \colon \ZZ(k)(U)\otimes\ZZ(\ell)(U)\to\ZZ(k+\ell)(U)
	\end{equation*}
	is given by 
	\begin{equation*}
		x \cupprod y =
		\begin{cases}
			x\cdot y \comma & \deg(x)=0\\
			x\wedge \iota y \comma  &\deg(x)>0, \deg(y)=0\\
			x\wedge \d y \comma & \deg(x)>0,\deg(y)=\ell>0\\
			0 \comma & \textup{otherwise } \phantom{0,\deg(y)=\ell>0} \period
		\end{cases}
	\end{equation*}
\end{proposition}

\begin{remark}
	This is only commutative up to homotopy.
\end{remark}

%-------------------------------------------------------------------%
%-------------------------------------------------------------------%
%  Examples                                                         %
%-------------------------------------------------------------------%
%-------------------------------------------------------------------%

\subsection{Examples}\label{subsec:Delignecupexamples}

We analyze the Deligne cup product in detail in the lowest dimensions. Let $M$ be a manifold. Recall the following computations.
\begin{itemize}
	\item $\ZZ(0)=\Gamma^*\ZZ[0]$ is the complex with $\Gamma^*\ZZ$ in degree zero. 
	Thus $\Hcech^0(M)=\H^0(M;\ZZ)$.

	\item $\Hcech^1(M)=\Mapsm(M,\Uup_{1})$.

	\item $\Hcech^2(M)=\{\text{line bundles on $M$ with connection}\}/\!\sim$.
\end{itemize}
Let $\ZZ(k)^\ell$ denote the degree $\ell$ term of the complex $\ZZ(k)$. 
For example, $\ZZ(3)^2=\Omega^1$. 
Let $\Ucal$ be a good cover for $M$. 
Using Čech cohomology for this good cover, the Deligne cup product gives a map
\begin{equation*}
	\paren{\bigoplus_{i+j=k}\Cech^i(\Ucal;\ZZ(k)^j)} \tensor \paren{\bigoplus_{i+j=l}\Cech^i(\Ucal;\ZZ(\ell)^j)} \longrightarrow \paren{\bigoplus_{i+j=k+\ell}\Cech^i(\Ucal;\ZZ(k+\ell)^j)} \period
\end{equation*}

\begin{example}
	The Deligne cup product
	\[
		\ZZ(0)\otimes\ZZ(0)\to\ZZ(0)
	\]
	should give us a way of taking two locally constant functions of $M\to\ZZ$ and producing a third.
	By \Cref{formula}, the Deligne cup product of two elements in degree 0 agrees with the ordinary cup product in $\H^0(M;\ZZ)$; i.e., the product of the two locally constant functions.
\end{example}

\begin{example}
	The Deligne cup product
	\[
		\ZZ(0)\otimes\ZZ(1)\to\ZZ(1)
	\]
	should give us a way of taking a locally constant function $M\to\ZZ$ and a smooth map $g\colon M\to \Uup_{1}$ and producing a new smooth map $M\to \Uup_{1}$.  
	In the Čech complex, we are looking at a map
	\[
		\Cech^0(\Ucal;\ZZ(0)^0)\otimes \left(\Cech^0(\Ucal;\ZZ(1)^1)\oplus \Cech^1(\Ucal;\ZZ(1)^0)\right)\to \left(\Cech^0(\Ucal;\ZZ(1)^1)\oplus \Cech^1(\Ucal;\ZZ(1)^0)\right)
	\]
	Identifying these terms, we have
	\[
		\Cech^0(\Ucal;\ZZ)\otimes \left(\Cech^0(\Ucal;\Omega^0)\oplus \Cech^1(\Ucal;\ZZ)\right)\to \left(\Cech^0(\Ucal;\Omega^0)\oplus \Cech^1(\Ucal;\ZZ)\right)
	\]
	This sends $n\otimes(f,m)$ to $(n\cdot f,n\cdot m)$.
\end{example}

\begin{example}
	The Deligne cup product
	\[
		\ZZ(1)\otimes\ZZ(0)\to\ZZ(1)
	\]
	should give us a way of taking a locally constant function $M\to \ZZ$ and a smooth map $g\colon M\to \Uup_{1}$ and producing a new smooth map $M\to \Uup_{1}$. In the Čech complex, we are looking at a map
	\[ 
		\left(\Cech^0(\Ucal;\Omega^0)\oplus \Cech^1(\Ucal;\ZZ)\right)\otimes\Cech^0(\Ucal;\ZZ)\to \left(\Cech^0(\Ucal;\Omega^0)\oplus \Cech^1(\Ucal;\ZZ)\right)
	\]
	This map sends $(f,m)\otimes n)$ to $(f\cdot \iota n,m\cdot n)$.
\end{example}

More geometrically, we can describe the Deligne cup product as follows. 
Given a pair $(n,f)$ where $n\colon M\to\ZZ$ is a locally constant function and $f\colon M\to \Circ$ is a smooth map, the Deligne cup product of $n$ with $f$ is the smooth function $g\colon M\to \Circ$ given by $g(x)=e^{2\pi i n(x)}f(x)$. 

\begin{remark}
	We can note that the Deligne cup product commutes up to homotopy,
	\begin{equation*}
		\begin{tikzcd}
			\ZZ(1)\otimes\ZZ(0)\arrow[r]\arrow[d] & \ZZ(1)\\
			\ZZ(0)\otimes\ZZ(1)\arrow[ur] & 
		\end{tikzcd}
	\end{equation*}
	since $(f\cdot \iota n=n\cdot f)$ as functions to $\RR$.
\end{remark}

\begin{example}
	The Deligne cup product
	\[
		\ZZ(1)\otimes\ZZ(1)\to\ZZ(2)
	\]
	should give us a way of taking two smooth maps $M\to \Uup_{1}$ and producing a line bundle on $M$ with connection. In the Čech complex, we are looking at a map
	\[
		\left(\Cech^0(\Ucal;\ZZ(1)^1)\oplus \Cech^1(\Ucal;\ZZ(1)^0)\right)^{\otimes 2}\to\left(\Cech^0(\Ucal;\ZZ(2)^2)\oplus\Cech^1(\Ucal;\ZZ(2)^1)\oplus\Cech^2(\Ucal;\ZZ(2)^0)\right) \period
	\]
	Then the Deligne cup product sends 
	\[
		(f,n)\otimes(g,m)\mapsto(n_{\alpha\beta}\cdot m_{\beta\gamma},n_\alpha\beta\cdot g_\beta+0,f_\alpha \d g_\alpha) \period
	\]
	If we think of $(f,n)$ and $(g,m)$ as smooth maps $M\to \Uup_{1}$, then $(n_{\alpha\beta}\cdot m_{\beta\gamma}, n_{\alpha\beta}\cdot g_\beta,f_\alpha \d g_\alpha)$ corresponds to the line bundle with transition function $n_{\alpha\beta}\cdot g_\beta$ and connection given by one form $(2\pi i)f_\alpha \d g_\alpha$.

	By \cite[Lemma 1.3.1]{MR760999}, the curvature of $f\cupprod g$ is $\dlog (f)\wedge\dlog(g)$.
\end{example}

\newpage
%!TEX root = ../diffcoh.tex

%-------------------------------------------------------------------%
%-------------------------------------------------------------------%
%  Fiber integration                                                %
%-------------------------------------------------------------------%
%-------------------------------------------------------------------%

\section{Fiber integration}\label{FiberIntegration}
\textit{by Araminta Amabel}

The goal of this section is to define a refinement of \textit{fiber integration} (along with its usual properties) in the setting of differential cohomology.
In ordinary cohomology, we get a fiber integration map from combining the Thom isomorphism and the suspension isomorphism. 
Let $E\to B$ be an oriented fiber bundle with fiber a compact manifold of dimension $k$. 
Let $ E\hookrightarrow\RR^{N}$ be an embedding with normal bundle $\nu$, and let $E^\nu$ denote the Thom space of $ \nu $. 
Then fiber integration is given by the composite
\begin{equation*}
	\begin{tikzcd}
		\H^{q+k}(E) \arrow[r, "\sim"{yshift=-0.25em}] & \H^{q+N}(E^\nu) \arrow[r, "\PT"] & \H^{q+N}(B_+\wedge \Sph{N}) \simeq \H^{q}(B_+) \comma
	\end{tikzcd}
\end{equation*}
where the first map is the Thom isomorphism, the second map is the Pontryagin--Thom collapse map, and the third map is the suspension isomorphism. Recall that the Thom isomorphism is given by taking the cup product with the Thom class. 

To do fiber integration in differential cohomology, we need to provide differential refinements of the following:
\begin{enumerate}
	\item Thom classes/orientations.

	\item The suspension isomorphism.
\end{enumerate}
To do this, we combine fiber integration in ordinary cohomology with integration of forms. 

%-------------------------------------------------------------------%
%-------------------------------------------------------------------%
%  Differential integration                                         %
%-------------------------------------------------------------------%
%-------------------------------------------------------------------%

\subsection{Differential integration}\label{subsec:differentialintegration}

The input will be a fiber bundle of manifolds
\begin{equation*}
	M\to E\to X \comma
\end{equation*}
where $ M $ is a closed, smooth manifold of dimension $ d $.
The output will be a map of %ring? 
spectra
\begin{equation*}
	\ZZ(k)(E)\to\Sigma^{d}\ZZ(k-d)(X)
\end{equation*}
where $\ZZ(k)$ is the pullback 
\begin{equation*}
	\begin{tikzcd}
		\ZZ(k)\arrow[r]\arrow[d] & \Gamma^*\HZZ\arrow[d]\\
		\Sigma^{-k}\HOmegacl^k\arrow[r] & \Gamma^*\HRR
	\end{tikzcd}
\end{equation*}
in $\Sh(\Mfld;\Sp)$ and, similarly, $\ZZ(k-d)$ is the pullback
\begin{equation*}
	\begin{tikzcd}
		\ZZ(k-d)\arrow[r]\arrow[d] & \Gamma^*\HZZ\arrow[d]\\
	\Sigma^{d-k}\HOmegacl^{k-d}\arrow[r] & \Gamma^*\HRR \period
	\end{tikzcd}
\end{equation*}
To produce a map $\ZZ(k)\to\Sigma^d\ZZ(k-d)$, it therefore suffices to produce maps $\HZZ\to\Sigma^d \HZZ $ and $\Omegacl^k\to\Omegacl^{k-d}$
 together with a path between their images in $\Sigma^{d}\Gamma^*\HRR$.

%-------------------------------------------------------------------%
%  Differential Thom classes and orientations                       %
%-------------------------------------------------------------------%

\subsection{Differential Thom classes and orientations}

\begin{definition}
	Let $M$ be a smooth compact manifold and $V\to M$ a real vector bundle of dimension $k$. 
	A \emph{differential Thom cocycle} on $V$ is a cocycle
	\begin{equation*}
		U=(c,h,\omega)\in \check{Z}(k)^k_{\mathrm{c}}(V)
	\end{equation*}
	such that, for each $m\in M$
	\begin{equation*}
		\int_{V_m}\omega =\pm 1
	\end{equation*}
\end{definition}

\begin{remark}
	A differential Thom class determines a ordinary Thom class in integral cohomology $\Hc^k(V;\ZZ)$.
\end{remark}

\begin{definition}[{\cite[Definition 2.9]{HopkinsSinger}}]
	An \emph{$\Hhat$-orientation} of $p\colon E\to B$ consists of the following three pieces of data:
	\begin{enumerate}
		\item A smooth embedding $E\subset B\times\RR^N$ for some $N$.

		\item A tubular neighborhood $W\subset B\times\RR^N$.

		\item A differential Thom cocycle $U$ on $W$.
	\end{enumerate}
\end{definition}

%-------------------------------------------------------------------%
%  Differential fiber integration                                   %
%-------------------------------------------------------------------%

\subsection{Differential fiber integration}

Our hope is to get an analogue of the suspension isomorphism
\begin{equation*}
	\Hc^{q+N}(B\times\RR^N)\simeq \H^q(B) \period
\end{equation*}
To understand the correct analogue of the suspension isomorphism in the differential setting, let us consider the most simple case.

\begin{example}
	Consider the case when $B$ is a point and $N=1$.
	Then the ordinary suspension isomorphism says that 
	\begin{equation*}
		\H^1(\Circ;\ZZ)\cong  \H^0(\pt;\ZZ)\simeq \ZZ
	\end{equation*}
	The calculation $\H^1(\Circ;\ZZ) \isomorphic \ZZ$ is by degree:
	\begin{equation*}
		\begin{tikzcd}[sep=3em]
			\H^1(\Circ;\ZZ) = \uppi_0 \Map_{\Spc}(\Circ,\K(\ZZ,1)) = \uppi_0 \Map_{\Spc}(\Circ,\Circ) \arrow[r, "\sim"{yshift=-0.25em}, "\deg"'] & \ZZ \period 
		\end{tikzcd}
	\end{equation*}
	In differential cohomology, we have an isomorphism
	\begin{equation*}
		\Hhat^1(\Circ) \isomorphic \Mapsm(\Circ,\Circ) \period
	\end{equation*}
	We still have a degree map
	\begin{equation*}
		\deg \colon \Mapsm(\Circ,\Circ) \to \ZZ \comma
	\end{equation*}
	but it is no longer an isomorphism.
\end{example}

The upshot is that we are looking for a suspension \textit{map} not an isomorphism.

\begin{nul}
	We start by working with the trivial bundle $B\times\RR^N\to B$ and defining integration for compactly-supported forms. 
	This is \cite[\S 3.4]{HopkinsSinger}. 
	Define the map
	\begin{equation*}
		\int_{B\times\RR^N/B} \colon \Cech(p+N)_{\mathrm{c}}^{q+N}(B\times\RR^N)\to\Cech(p)^q(B)
	\end{equation*}
	by the slant product with a fundamental cycle $Z_N\in \Cup_N(\RR^N;\ZZ)$,
	\begin{equation*}
		(c,h,\omega)\mapsto\left(c/Z_N,h/Z_N,\int_{B\times\RR^N/B}\omega\right)
	\end{equation*}
	Note that this is simply a map, \textit{not} an isomorphism.
\end{nul}

\begin{remark}
	Checking that the slant product goes through to differential cohomology seems to require some work. 
	See \cite[\S3.4]{HopkinsSinger}.
\end{remark}

\begin{definition}[{\cite[Definition 3.11]{HopkinsSinger}}]
	Let $p\colon E\to B$ be an $\Hhat$-oriented map of smooth manifolds with boundary of relative dimension $k$. 
	The \emph{integration map} is the map
	\begin{equation*}
		\int_{E/B}\colon\Cech(p+k)^{q+k}(E)\to \Cech(p)^q(B)
	\end{equation*}
	given by the composite
	\begin{equation*}
		\begin{tikzcd}[sep=3em]
			\Cech(p+k)^{q+k}(E) \arrow[r, "\cupprod U"] & \Cech(p+N)_{\mathrm{c}}^{q+N}(B\times\RR^N) \arrow[r, "\int_{\RR^N}(-)"] & \Cech(p)_{\mathrm{c}}^q(B) \period
		\end{tikzcd}
	\end{equation*}
\end{definition}

\begin{example}
	In dimension $ 1 $, the only closed manifold is $\Circ$. 
	If $E\to B$ is an oriented $\Circ$-bundle, then integration along the fibers defines a map
	\begin{equation*}
		\int_{E/B}\colon \Hhat^2(E)\to\Hhat^1(E)
	\end{equation*}
	If $ x \in \Hhat^2(E) $ corresponds to a line bundle with connection, then 
	\begin{equation*}
		\int_{E/B}x
	\end{equation*}
	represents the function $ B \to \Circ $ sending $ b \in B $ to the monodromy of $ x $ computed around the fiber $ E_b $.
\end{example}
\newpage
%!TEX root = ../diffcoh.tex

%-------------------------------------------------------------------%
%-------------------------------------------------------------------%
%  Digression: the Transfer Conjecture                              %
%-------------------------------------------------------------------%
%-------------------------------------------------------------------%

\section{Digression: the Transfer Conjecture}\label{sec:digressiononTransferConjecture}
\textit{by Peter Haine}

%-------------------------------------------------------------------%
%-------------------------------------------------------------------%
%  Introduction                                                     %
%-------------------------------------------------------------------%
%-------------------------------------------------------------------%

\subsection{Introduction}

Let $ X $ be a space.
We have seen that the constant sheaf of spaces $ \Gammaupperstar(X) $ on $ \Mfld $ is given by the formula
\begin{equation*}
	 \Gammaupperstar(X) = \Map_{\Space}(\Piinf(M),X)
\end{equation*}
(\Cref{lem:constantishi}).
If $ X = \Omega^{\infty} E $ is the infinite loop space of a spectrum $ E $, then the sheaf $ \Gammaupperstar(X) $ acquires additional functoriality: for any finite covering map between manifolds $ f \colon \fromto{N}{M} $, the \textit{Becker--Gottlieb transfer}
\begin{equation*}
	\fromto{\Sigma_{+}^{\infty} \Piinf(M)}{\Sigma_{+}^{\infty} \Piinf(N)}
\end{equation*}
\cite[Definition 3.11]{Haugseng:Becker-Gottlieb} induces a \textit{transfer} map
\begin{equation*}
	\flowerstar \colon \fromto{\Gammaupperstar(X)(N)}{\Gammaupperstar(X)(M)} \period
\end{equation*} 
This enhanced functoriality can be used to make $ \Gammaupperstar(X) $ into a copresheaf on a $ 2 $-category $ \Corfcov $ with objects smooth manifolds and morphisms \textit{correspondences}
\begin{equation*}
	\begin{tikzcd}[row sep=1.5em, column sep=0.75em]
		& N \arrow[dr, "f"] \arrow[dl] & \\
		M_0 & & M_1 \comma
	\end{tikzcd}
\end{equation*}
where $ f $ is a finite covering map.
Composition in $ \Corfcov $ is given by pullback. 

For a sheaf $ F $ on $ \Mfld $, Quillen conjectured that an extension of $ F $ to $ \Corfcov $ is just another way of encoding an $ \Einf $-structure on $ F $.
However, when Quillen originally formulated this \textit{Transfer Conjecture}, the language to express the higher coherences necessary for the validity of the result was not available.
Moreover, Quillen's original formulation was \textit{dis}proven by Kraines and Lada \cites{MR557187}{MR3243401}.

The goal of this section is to explain how to deduce the following corrected version of the Transfer Conjecture from very general results of Bachmann--Hoyois on commutative algebras and \categories of spans \cite[Appendix C]{MotivicNorms:BachmannHoyois}.

\begin{theorem}[{(Transfer Conjecture; \Cref{cor:transferconjecture,cor:transferconjecturehi})}]\label{thm:transfer}
	Let $ C $ be a presentable \category.
	There is an equivalence of \categories
	\begin{equation*}
		\Fun_{\loc}(\Corfcov,C) \equivalence \Sh(\Mfld;\CMon(C))
	\end{equation*}
	between functors $ \fromto{\Corfcov}{C} $ whose restriction to $ \Mfldop $ is a sheaf and sheaves of commutative monoids in $ C $.
	This further restricts to an equivalence
	\begin{equation*}
		\Fun_{\loc,\RR}(\Corfcov,C) \equivalence \CMon(C)
	\end{equation*}
	between functors $ \fromto{\Corfcov}{C} $ whose restriction to $ \Mfldop $ is an \RRinvariant sheaf and commutative monoids in $ C $.
\end{theorem}

\begin{example}
	Setting $ C = \Spc $ in \Cref{thm:transfer} gives an equivalence between functors
	\begin{equation*}
		\fromto{\Corfcov}{\Spc}
	\end{equation*}
	whose restriction to $ \Mfldop $ is an \RRinvariant sheaf and $ \Einf $-spaces.
	Restricting to grouplike objects on both sides and applying the Segal's Recognition Principle for connective spectra \HA{Remark}{5.2.6.26} provides an equivalence between grouplike objects of
	\begin{equation*}
		\Fun_{\loc,\RR}(\Corfcov,\Space)
	\end{equation*}
	and the \category $ \Sp_{\geq 0} $ of connective spectra. 
\end{example}

\begin{remark}
	The Becker--Gottlieb transfer is defined in more generality than finite covering maps; for example, for proper submersions.
	It is possible to modify \Cref{thm:transfer} to encode this additional generality.
	However, since pullbacks along proper submersions do not exist in the category of manifolds, in order for composition of correspondences where one leg is proper to be defined, one needs to work with \textit{derived} manifolds \cites{arXiv:1905.06195}{MR2641940}.
	For the sake of simplicity, we will satisfy ourselves with just working with manifolds and finite covering maps.
\end{remark}

In order to give a more precise formulation of \Cref{thm:transfer}, we'll first review constructing $ 2 $-categories of \textit{correspondences} or \textit{spans} from $ 1 $-categories (\cref{subsec:spans}). 
We then briefly recall the role that \categories of spans play in encoding $ \Einf $-structures (\cref{subsec:spansCMon}).
Finally, we walk through \cite[Appendix C]{MotivicNorms:BachmannHoyois} in the case of interest and explain how to deduce the Transfer Conjecture from their results (\cref{subsec:transferproof}).  

%-------------------------------------------------------------------%
%  Categories of spans                                              %
%-------------------------------------------------------------------%

\subsection{Categories of spans}\label{subsec:spans}

In this section we explain how to construct the $ 2 $-category $ \Corfcov $ of correspondences of manifolds appearing in the Transfer Conjecture.
This is a special case of a general construction for \categories due to Barwick \cite[§§3--5]{MR3558219}.
If $ D $ is an $ n $-category, then Barwick's \category of spans in $ D $ is an $ (n+1) $-category. 
In order to avoid explaining how to deal with the homotopy coherence problems that arise, we only present the $ 1 $-categorical case as we can give a simple definition as a $ 2 $-category.

\begin{construction}[($ 2 $-category of spans)]\label{def:spancat}
	Let $ D $ be a $ 1 $-category, and let $ L, R \subset \Mor(D) $ be two classes of morphisms in $ D $ satisfying the following properties:
	\begin{enumerate}[label=\stlabel{def:spancat}, ref=\arabic*]
		\item The classes $ L $ and $ R $ contain all isomorphisms.

		\item The classes $ L $ and $ R $ are each stable under composition.

		\item Given a morphism $ \el \colon \fromto{X}{Z} $ in $ L $ and morphism $ r \colon \fromto{Y}{Z} $ in $ R $, there exists a pullback diagram
		\begin{equation*}
			\begin{tikzcd}
				W \arrow[r, "\rbar"] \arrow[d, "\elbar"'] \arrow[dr, phantom, very near start, "\lrcorner"] & Y \arrow[d, "r"] \\
             	X \arrow[r, "\el"'] & Z 
			\end{tikzcd}
		\end{equation*}
		in $ D $ where $ \elbar \in L $ and $ \rbar \in R $.
	\end{enumerate}

	Define a $ 2 $-category $ \Span(D;L,R) $ as follows.
	The objects of $ \Span(D;L,R) $ are the objects of $ D $.
	Given objects $ X_0, X_1 \in D $, the groupoid $ \Map_{\Span(D;L,R)}(X_0,X_1) $ has objects diagrams 
	\begin{equation*}
		\begin{tikzcd}[row sep=1.5em, column sep=0.75em]
			& Y \arrow[dr, "r"] \arrow[dl, "\el"'] & \\
			X_0 & & X_1 \comma
		\end{tikzcd}
	\end{equation*}
	in $ D $ where $ \el \in L $ and $ r \in R $, and morphisms isomorphisms of diagrams.
	Composition is given by pullback of spans: given morphisms $ X_0 \to X_1 $ and $ X_1 \to X_2 $ corresponding to spans
	\begin{equation*}
		\begin{tikzcd}[row sep=1.5em, column sep=0.75em]
			& Y \arrow[dr] \arrow[dl] & \\
			X_0 & & X_1 
		\end{tikzcd} \andeq
		\begin{tikzcd}[row sep=1.5em, column sep=0.75em]
			& Z \arrow[dr] \arrow[dl] & \\
			X_1 & & X_2 \comma
		\end{tikzcd}
	\end{equation*}
	the composite morphism $ X_0 \to X_2 $ in $ \Span(D;L,R) $ is defined as the large pullback span
	\begin{equation*}
		\begin{tikzcd}[column sep={6ex,between origins}, row sep={6ex,between origins}]
			&  & \displaystyle Y \cross_{X_1} Z \arrow[dr] \arrow[dl] & & \\
			& Y \arrow[dr] \arrow[dl] &  & Z \arrow[dr] \arrow[dl] & \\
			X_0 & & X_1 & & X_2 \period
		\end{tikzcd}
	\end{equation*}
\end{construction}

\begin{notation}
	Let $ D $ be a $ 1 $-category.
	We write $ \all \colonequals \Mor(D) $ for the class of all morphisms in $ D $.
	If $ D $ has pullbacks, we write
	\begin{equation*}
		\Span(D) \colonequals \Span(D;\all,\all)
	\end{equation*}
	for the $ 2 $-category of spans of arbitrary morphisms in $ D $.
\end{notation}

\begin{observation}
	Let $ D $ be a category and $ R $ a class of morphisms in $ D $ such that the pullback of a morphism in $ R $ along an arbitrary morphism of $ D $ exists, and the class $ R $ is stable under pullback.
	Then there is a natural faithful functor
	\begin{equation*}
		\fromto{D^{\op}}{\Span(D;\all,R)}
	\end{equation*}
	given by the identity on objects, and on morphisms by sending a morphism $ f \colon \fromto{X}{Y} $ to the span
	\begin{equation*}
		\begin{tikzcd}[row sep=1.5em, column sep=0.75em]
			& X \arrow[dl, "f"'] \arrow[dr, equals] & \\
			Y & & X \period
		\end{tikzcd}
	\end{equation*}
\end{observation}

\begin{example}
	Write $ \fcov \subset \Mor(\Mfld) $ for the class of finite covering maps of manifolds.
	Note that the pullback of a finite covering map of manifolds along any morphism exists, and the class of finite covering maps is stable under pullback.
	We write
	\begin{equation*}
		\Corfcov \colonequals \Span(\Mfld;\all,\fcov)
	\end{equation*}
	for the $ 2 $-category with objects manifolds and morphisms \textit{correspondences}\footnote{The term
	``correspondence'' is just another name for a span; ``correspondence'' seems to be the more common term in geometry.} of manifolds 
	\begin{equation*}
		\begin{tikzcd}[row sep=1.5em, column sep=0.75em]
			& N \arrow[dr, "f"] \arrow[dl] & \\
			M_0 & & M_1 \comma
		\end{tikzcd}
	\end{equation*}
	where $ f $ is a finite covering map.
\end{example}

\begin{example}
	Write $ \foldtext \subset \Mor(\Mfld) $ for the class of maps that are finite coproducts of fold maps of manifolds, i.e., finite coproducts of fold maps $ \nabla \colon \fromto{M^{\coproduct i}}{M} $ from a finite disjoint union of copies of $ M $ to $ M $.
	Note that coproduct decompositions are stable under all pullbacks that exist in the category of manifolds, hence the class $ \foldtext $ is stable under pullback.
	We write
	\begin{equation*}
		\Corfold \colonequals \Span(\Mfld;\all,\foldtext)
	\end{equation*}
	for the $ 2 $-category with objects manifolds and morphisms correspondences of manifolds 
	\begin{equation*}
		\begin{tikzcd}[row sep=1.5em, column sep=0.75em]
			& N \arrow[dr, "f"] \arrow[dl] & \\
			M_0 & & M_1 \comma
		\end{tikzcd}
	\end{equation*}
	where $ f $ is a finite coproduct of fold maps.

	Note that $ \foldtext \subset \fcov $, so that $ \Corfold $ defines a subcategory of $ \Corfcov $ that contains all objects.
\end{example}

%-------------------------------------------------------------------%
%  Spans and commutative monoids                                    %
%-------------------------------------------------------------------%

\subsection{Spans and commutative monoids}\label{subsec:spansCMon}

In this section we briefly recall the role that \categories of spans play in encoding $ \Einf $-structures.
We begin by introducing the relevant $ 2 $-category of spans.

\begin{notation}
	Write $ \Fin $ for the category of finite sets.
	Given \acategory $ C $ with a terminal object, we write $ * $ for the terminal object.
\end{notation}

\begin{recollection}
	Let $ C $ be \acategory with finite products.
	A \textit{commutative monoid} or \textit{$ \Einf $-monoid} in $ C $ is a functor $ M \colon \fromto{\Finstar}{C} $ such that $ \equivto{M(*)}{*} $ and for each integer $ n \geq  1 $, the collapse maps $ \fromto{\{1,\ldots,n\}_+}{\{i\}_+} $ induce an equivalence
	\begin{equation*}
		\equivto{M(\{1,\ldots,n\}_+)}{\prod_{i=1}^n M(\{i\}_+)} \period
	\end{equation*}
	We write $ \CMon(C) \subset \Fun(\Finstar,C) $ for the full subcategory spanned by the commutative monoids.
\end{recollection}

\begin{remark}
	By induction, a functor $ M \colon \fromto{\Finstar}{C} $ is a commutative monoid if and only if $ \equivto{M(*)}{*} $ and for every pair $ S,T \in \Finstar $, the functor $ M $ carries the pushout square
	\begin{equation*}
		\begin{tikzcd}[sep=3em]
			S \vee T \arrow[r, "\ast \vee \id{T}"] \arrow[d, "\id{S} \vee \ast"'] \arrow[dr, phantom, "\ulcorner"{description, very near end}] & T \arrow[d] \\ 
			S \arrow[r] & \ast 
		\end{tikzcd}
	\end{equation*}
	to a \textit{pullback} square in $ C $.
\end{remark}

\begin{observation}
	The $ 2 $-category $ \Span(\Fin) $ is \textit{semiadditive}: the direct sum in $ \Span(\Fin) $ is given by disjoint union of finite sets.
	See \cites[Lemma C.3]{MotivicNorms:BachmannHoyois}[Proposition 4.3]{MR3558219} for more general results on the semiadditivity of \categories of spans.
\end{observation}

\begin{observation}[(finite pointed sets via spans)]
	Write $ \inj $ for the class of injective maps in $ \Fin $.
	The functor
	\begin{equation*}
		\fromto{\Finstar}{\Span(\Fin;\inj,\all)}
	\end{equation*}
	given by sending $ \goesto{X_{+}}{X} $ and a morphism $ f \colon \fromto{X_{+}}{Y_{+}} $ to the span
	\begin{equation*}
		\begin{tikzcd}[row sep=1.5em, column sep=0.75em]
			& \finverse(Y) \arrow[dr, "f"] \arrow[dl, hooked'] & \\
			X & & Y 
		\end{tikzcd}
	\end{equation*}
	is an equivalence of categories.

	The category $ \Span(\Fin;\inj,\all) $ is often referred to as the category of finite sets and \textit{partially defined maps}.
\end{observation}

The importance of transfers in $ \Einf $-structures is explained by the following universal property of the $ 2 $-category $ \Span(\Fin) $ of spans of finite sets.

\begin{proposition}[{(Cranch \cites[Proposition C.1]{MotivicNorms:BachmannHoyois}[§5]{Cranch:Thesis})}]\label{prop:BH.C.1}
	Let $ C $ be \acategory with finite products.
	Then the restriction 
	\begin{equation*}
		\fromto{\Fun(\Span(\Fin),C)}{\Fun(\Finstar,C)}
	\end{equation*}
	along the inclusion $ \fromto{\Finstar}{\Span(\Fin)} $ restricts to an equivalence between:
	\begin{enumerate}[label=\stlabel{prop:BH.C.1}, ref=\arabic*]
		\item Functors $ M \colon \fromto{\Span(\Fin)}{C} $ that preserve finite products (equivalently, $ \restrict{M}{\Finop} $ preserves finite products).

		\item Commutative monoids in $ C $.
	\end{enumerate}
	The inverse is given by right Kan extension.
\end{proposition}

The $ 2 $-category $ \Span(\Fin) $ has a second (related) universal property: $ \Span(\Fin) $ is the free semiadditive \category generated by a single object.

\begin{proposition}[{(Harpaz \cite[Theorem 1.1]{Harpaz:Spans})}]
	Let $ C $ be a semiadditive \category.
	Then evaluation at $ \ast \in \Span(\Fin) $ defines an equivalence
	\begin{equation*}
		\equivto{\Fun^{\directsum}(\Span(\Fin),C)}{C} \period
	\end{equation*}
\end{proposition}

%-------------------------------------------------------------------%
%  The Transfer Conjecture after Bachmann–Hoyois                    %
%-------------------------------------------------------------------%

\subsection{The Transfer Conjecture after Bachmann--Hoyois}\label{subsec:transferproof}

In this section we outline work of Bachmann--Hoyois that implies the Transfer Conjecture \cite[Appendix C]{MotivicNorms:BachmannHoyois}. 
The perspective on commutative monoids in $ C $ as finite product-preserving functors $ \fromto{\Span(\Fin)}{C} $ (\Cref{prop:BH.C.1}) is fundamental to proving the Transfer Conjecture.

The first step is to relate finite product-preserving functors $ \fromto{\Corfold}{C} $ to presheaves of commutative monoids on $ \Mfld $.
Then we impose the sheaf condition to pass from $ \Corfold $ to $ \Corfcov $.

\begin{notation}
	Write $ \Theta \colon \fromto{\Mfldop \cross \Span(\Fin)}{\Corfold} $ for the functor given on objects by the assignment
	\begin{equation*}
		(M,I) \mapsto M^{\coproduct I}
	\end{equation*}
	and on morphisms by the assignment
	\begin{equation*}
		(M \to N, I_0 \leftarrow J \to I_1) \quad \mapsto
		\quad
		\begin{tikzcd}[row sep=1.5em, column sep=0.75em]
			& M^{\coproduct J} \arrow[dl] \arrow[dr] & \\
			N^{\coproduct I_0} & & M^{\coproduct I_1} \period
		\end{tikzcd} 
	\end{equation*}
\end{notation}

The functor $ \Theta $ is the universal functor that preserves finite products in each variable:

\begin{proposition}[{\cite[Proposition C.5]{MotivicNorms:BachmannHoyois}}]\label{prop:BH.C.5}
	Let $ C $ be \acategory with finite products.
	Then the restriction functor
	\begin{equation*}
		\Theta\upperstar \colon \Fun(\Corfold,C) \to \Fun(\Mfldop \cross \Span(\Fin), C)
	\end{equation*}
	restricts to an equivalence 
	\begin{equation*}
		\Funcross(\Corfold,C) \equivalence \Funcross(\Mfldop, \CMon(C)) \period
	\end{equation*}
	The inverse is given by right Kan extension along $ \Theta $.
\end{proposition}

Since every finite covering map is locally a fold map, we see:

\begin{proposition}[{\cite[Proposition C.11]{MotivicNorms:BachmannHoyois}}]\label{prop:BH.C.11}
	Let $ C $ be \acategory with finite products.
	Then the restriction functor
	\begin{equation*}
		\Fun(\Corfcov,C) \to \Fun(\Corfold,C) 
	\end{equation*}
	restricts to an equivalence between the full subcategories of those functors whose restrictions to $ \Mfldop $ are sheaves.
	The inverse is given by right Kan extension.
\end{proposition}

\begin{notation}
	Write
	\begin{equation*}
		\Fun_{\loc}(\Corfcov,C) \subset \Fun(\Corfcov,C)
	\end{equation*}
	for the full subcategory spanned by those functors $ F $ whose restrictions to $ \Mfldop $ are sheaves.
\end{notation}

We now arrive at Quillen's Transfer Conjecture:

\begin{corollary}[(Transfer Conjecture)]\label{cor:transferconjecture}
	Let $ C $ be \acategory with all limits.
	Restriction along the inclusion $ \incto{\Mfldop}{\Corfcov} $ defines an equivalence of \categories 
	\begin{equation*}
		\Fun_{\loc}(\Corfcov,C) \equivalence \Sh(\Mfld;\CMon(C)) \period
	\end{equation*}
\end{corollary}

\noindent Combining \Cref{prop:Dugger} and \Cref{cor:transferconjecture} shows:

\begin{corollary}\label{cor:transferconjecturehi}
	Let $ C $ be a presentable \category.
	Restriction along the inclusion
	\begin{equation*}
		\incto{\Mfldop}{\Corfcov}
	\end{equation*}
	defines an equivalence of \categories 
	\begin{equation*}
		\Fun_{\loc,\RR}(\Corfcov,C) \equivalence \Shhi(\Mfld;\CMon(C)) \period
	\end{equation*}
	Post-composing with the global sections functor $ \Gammalowerstar $ defines an equivalence
	\begin{equation*}
		\Fun_{\loc,\RR}(\Corfcov,C) \equivalence \CMon(C) \period
	\end{equation*} 
\end{corollary}

\begin{nul}
	Unwinding the definitions we see that restriction along the inclusion
	\begin{equation*}
		\Span(\Fin) \subset \Corfcov
	\end{equation*}
	defines an equivalence
	\begin{equation*}
		\Fun_{\loc,\RR}(\Corfcov,C) \equivalence \Funcross(\Span(\Fin),C) \equivalent \CMon(C) \period
	\end{equation*}
\end{nul}

\newpage
\appendix
%!TEX root = ../diffcoh.tex

%-------------------------------------------------------------------%
%-------------------------------------------------------------------%
%  Appendix: Technical details from topos theory                    %
%-------------------------------------------------------------------%
%-------------------------------------------------------------------%

\section{Technical details from topos theory}\label{app:technicaldeatails}
\textit{by Peter Haine}

The purpose of this appendix is to prove a number of technical results used throughout the text.
We have relegated these proofs to this appendix because of one of the following reason:
\begin{enumerate}
	\item They are lengthy and, while the result is important, the proof is not important to know.

	\item They require some knowledge from the theory of \topoi.
\end{enumerate}
In \cref{sec:Cvaluedsheaves}, we explain a formal procedure to get from sheaves of spaces to sheaves valued in another presentable \category $ C $.
This lets us deduce many results about sheaves on $ \Mfld $ valued in a general presentable \category $ C $ from the case $ C = \Spc $.
\Cref{sec:restriction} explains the important properties of the functor given by restricting a sheaf defined on $ \Mfld $ to a sheaf defined on only a single manifold.
\Cref{sec:sheafification} explains why this restriction procedure commutes with sheafification.
In \cref{sec:completeness}, we give some background on notions of ``completeness'' for \topoi.
\Cref{subsec:points} shows that equivalences in $ \Sh(\Mfld;\Spc) $ can be checked on stalks and uses this to show that $ \Sh(\Mfld;\Spc) $ satisfies the strongest of these completeness notions (\Cref{app.prop:postnikovcompleteness}).
This also implies that $ \Sh(\Mfld;C) $ is equivalent to the category of $ C $-valued sheaves on the subcategory $ \Euc \subset \Mfld $ spanned by the Euclidean spaces (\Cref{app.lem:hypersheavesonCart}).

Since we are mostly interested in sheaves of spaces in this appendix, we adopt the following notational convention.

\begin{notation}
	We write $ \Sh(\Mfld) \colonequals \Sh(\Mfld;\Spc) $ for the \topos of sheaves of spaces on $ \Mfld $.
\end{notation}

\begin{remark}\label{rem:geometric_morphism}
	For this appendix, it is sufficient to know that the \category of sheaves of spaces on a site is \atopos, and that a \textit{geometric morphism} of \topoi is a right adjoint functor $ \flowerstar \colon \fromto{\X}{\Y} $ whose left adjoint $ \fupperstar $ is left exact.
\end{remark}

%-------------------------------------------------------------------%
%  From sheaves of spaces to C-valued sheaves                       %
%-------------------------------------------------------------------%

\subsection{From sheaves of spaces to \texorpdfstring{$ C $}{C}-valued sheaves}\label{sec:Cvaluedsheaves}

Let $ C $ be a presentable \category.
In this section we explain a formal procedure that allows us to pass from the \category $ \Sh(\Mfld;\Spc) $ of sheaves of spaces on $ \Mfld $ to the \category $ \Sh(\Mfld;C) $ of $ C $-valued sheaves on $ \Mfld $.
We'll also recall the basics of tensor products of presentable \categories and explain how to describe $ \Sh(\Mfld;C) $ as the tensor product
\begin{equation*}
	\Sh(\Mfld;C) \equivalent \Sh(\Mfld;\Spc) \tensor C \period
\end{equation*} 

The first thing to observe is that if $ G \colon \fromto{\Sh(\Mfld;\Spc)^{\op}}{C} $ is a functor that preserves limits, then the restriction $ G \colon \fromto{\Mfldop}{C} $ is a sheaf. 
It turns out that all $ C $-valued sheaves arise in this way.

\begin{proposition}[{\SAG{Proposition}{1.3.1.7}}]
	Let $ (S,\tau) $ be \asite and $ C $ an \category with all limits.
	Write $ \yo_{\tau} \colon \fromto{S}{\Sh_{\tau}(S;\Spc)} $ for the $ \tau $-sheafification of the Yoneda embedding.
	Then pre-composition with $ \yo_{\tau} $ defines an equivalence
	\begin{equation*}
		\Funlim(\Sh_{\tau}(S;\Spc)^{\op},C) \equivalence \Sh_{\tau}(S;C) \period
	\end{equation*} 
\end{proposition}

Now we give the \category $ \Funlim(\Sh(\Mfld)^{\op},C) $ a description in terms of a universal property of presentable \categories.

\begin{recollection}[{\HA{Proposition}{4.8.1.17}}]
	Let $ C $ and $ D $ be presentable \categories.
	The \textit{tensor product} of presentable \categories $ C \tensor D $ along with the functor $ \tensor \colon \fromto{C \cross D}{C \tensor D} $ are characterized by the following universal property: for any presentable \category $ E $, restriction along $ \tensor $ defines an equivalence 
	\begin{equation*}
		\Fun^{\colim}(C \tensor D,E) \equivalence \Fun^{\colim,\colim}(C \cross D,E) \period
	\end{equation*}
	Here the right-hand side is the full subcategory of $ \Fun(C \cross D,E) $ spanned by those functors $ \fromto{C \cross D}{E} $ that preserve colimits separately in each variable.
	The tensor product of presentable \categories defines a functor
	\begin{equation*}
		\tensor \colon \fromto{\PrL \cross \PrL}{\PrL}
	\end{equation*}
	and can be used to equip $ \PrL $ with the structure of a symmetric monoidal \category.

	The tensor product $ C \tensor D $ admits the following useful (seemingly asymmetric) description:
	\begin{equation*}
		C \tensor D \equivalent \Fun^{\lim}(C^{\op},D) \period
	\end{equation*}
	If $ F \colon \fromto{D}{D'} $ is a right adjoint functor of presentable \categories, then the induced right adjoint
	\begin{equation*}
		\id{C} \tensor F \colon C \tensor D \equivalent \Fun^{\lim}(C^{\op},D) \to \Fun^{\lim}(C^{\op},D') \equivalent C \tensor D'
	\end{equation*}
	is given by post-composition with $ F $.
	Unfortunately, the left adjoint to $ \id{C} \tensor F $ does not generally admit a simple description.
	However, if $ C $ is compactly assembled and the left adjoint to $ F $ is left exact, then the left adjoint to $ \id{C} \tensor F $ admits a simple description; see \cite[\S2.2]{arXiv:2108.03545}.
\end{recollection}

\begin{example}\label{ex:ShMan_tensor}
	For any presentable \category $ C $, we have a natural equivalence
	\begin{equation*}
		\equivto{\Sh(\Mfld) \tensor C}{\Sh(\Mfld;C)} \period
	\end{equation*}
\end{example}
	
%-------------------------------------------------------------------%
%  Restriction to a manifold                                        %
%-------------------------------------------------------------------%

\subsection{Restriction to a manifold}\label{sec:restriction}

We now give an alternative description of the functor $ \fromto{\Sh(\Mfld;C)}{C} $ that sends a sheaf to its value on a manifold $ M $.

\begin{notation}
	Let $ T $ be a topological space and $ C $ a presentable \category.
	Write
	\begin{equation*}
		\PSh(T;C) \colonequals \Fun(\Open(T)^{\op},C)
	\end{equation*}
	and write $ \Sh(T;C) \subset \PSh(T;C) $ for the \category of $ C $-valued sheaves on $ T $.
	Write
	\begin{equation*}
		\Gamma_{T,*} \colon \fromto{\Sh(T;C)}{C}
	\end{equation*}
	for the \textit{global sections} functor, defined by $ \Gamma_{T,*}(F) \colonequals F(T) $, and write $ \Gammaupperstar_T \colon \fromto{C}{\Sh(T;C)} $ for the left adjoint to $ \Gamma_{T,*} $, i.e., the \textit{constant sheaf} functor.
\end{notation}

\begin{observation}\label{obs:restrictiontoamanifold}
	Let $ C $ be a presentable \category and $ M $ a manifold.
	The forgetful functor $ \fromto{\Open(M)}{\Mfld} $ preserves finite limits and is a morphism of sites.
	Moreover, the forgetful functor satisfies the \emph{covering lifting property} \cite[Definition A.12]{arXiv:1803.01804}.
	In particular: 
	\begin{enumerate}[label=\stlabel{obs:restrictiontoamanifold}, ref=\arabic*]
		\item The presheaf retriction functor $ \restrict{(-)}{M} \colon \fromto{\PSh(\Mfld;C)}{\PSh(M;C)} $ carries sheaves to sheaves.
		
		\item The functor $ \restrict{(-)}{M} \colon \fromto{\Sh(M;C)}{\Sh(\Mfld;C)} $ is both a left and right adjoint \cite[Proposition A.18]{arXiv:1803.01804}. 
	\end{enumerate}
\end{observation}

\begin{nul}
	Note that the functor given by sending a sheaf $ E $ on $ \Mfld $ to its value on $ M $ is given by the composite
	\begin{equation*}
		\begin{tikzcd}
			\Sh(\Mfld;C) \arrow[r, "\restrict{(-)}{M}"] & \Sh(M;C) \arrow[r, "{\Gamma_{M,*}}"] & C \period
		\end{tikzcd}
	\end{equation*} 
\end{nul}

\begin{nul}
	Moreover, if $ p \colon \fromto{N}{M} $ is a morphism in $ \Mfld $, then there is a canonical natural transformation fitting into the triangle 
	\begin{equation*}
		\begin{tikzcd}[column sep={12ex,between origins}, row sep={8ex,between origins}]
			\Sh(\Mfld;C) \arrow[rr, "\restrict{(-)}{M}"] \arrow[dr, "\restrict{(-)}{N}"'] & \phantom{\Sh(N;C)} & \Sh(M;C) \\
			& \Sh(N;C) \arrow[ur, "\plowerstar"'] \arrow[u, phantom, "\Longrightarrow"{sloped, near end}, "{\, \scriptstyle\can_p}"{xshift=1em, yshift=.25em}]
		\end{tikzcd}
	\end{equation*} 
	defined as follows: given a sheaf $ E $ on $ \Mfld $ and an open subset $ U \subset M $, the morphism
	\begin{equation*}
		\fromto{E(U)}{E(\pinverse(U))}
	\end{equation*}
	is induced by the projection $ \fromto{\pinverse(U)}{U} $ by the functoriality of $ E $.
	In particular, upon taking global sections, the morphism 
	\begin{equation*}
		\can_p \colon E(M) = \Gamma_{M,*}(\restrict{E}{M}) \to  \Gamma_{M,*}(\plowerstar(\restrict{E}{N})) = E(N)
	\end{equation*}
	is the morphism $ \fromto{E(M)}{E(N)} $ induced by $ p $ by the functoriality of $ E $.
\end{nul}

%-------------------------------------------------------------------%
%  Sheafification                                                   %
%-------------------------------------------------------------------%

\subsection{Sheafification}\label{sec:sheafification}

Next we show that restriction from $ \Sh(\Mfld;C) $ to $ \Sh(M;C) $ commutes with sheafification.

\begin{nul}
	Consider the commutative square
	\begin{equation*}
		\begin{tikzcd}
			\Sh(\Mfld;C) \arrow[d, "\restrict{(-)}{M}"'] \arrow[r, hooked] & \PSh(\Mfld;C) \arrow[d, "\restrict{(-)}{M}"] \\ 
			\Sh(M;C) \arrow[r, hooked] & \PSh(M;C) \period
		\end{tikzcd}
	\end{equation*} 
	Using the unit and counit of the sheafification-inclusion adjunctions for $ \Mfld $ and $ M $, one can define an \textit{exchange transformation}
	\begin{equation*}
		\Ex \colon \fromto{\Sup_{M} \of \restrict{(-)}{M}}{\restrict{(-)}{M} \of \SMan} \period
	\end{equation*}
	See \cites[\HAthm{Definition}{4.7.4.13}]{HA}[Definition 1.1]{arXiv:2108.03545}.
	The exchange morphism $ \Ex $ fits into a diagram
	\begin{equation*}
		\begin{tikzcd}
			\PSh(\Mfld;C) \arrow[d, "\restrict{(-)}{M}"'] \arrow[r, "\SMan"] & \Sh(\Mfld;C) \arrow[d, "\restrict{(-)}{M}"] \arrow[dl, phantom, "\scriptstyle \Ex" above left, "\Longrightarrow" sloped] \\ 
			\PSh(M;C) \arrow[r, "\Sup_{M}"'] & \Sh(M;C) \period
		\end{tikzcd}
	\end{equation*}
\end{nul}

\begin{lemma}\label{lem:rescommuteswithsheaf}
	Let $ C $ be a presentable \category and $ M $ a manifold.
	Then the exchange transformation 
	\begin{equation*}
		\Ex \colon \fromto{\Sup_{M} \of \restrict{(-)}{M}}{\restrict{(-)}{M} \of \SMan}
	\end{equation*}
	is an equivalence.
	That is, there is a commuative square of \categories
	\begin{equation*}
		\begin{tikzcd}
			\PSh(\Mfld;C) \arrow[d, "\restrict{(-)}{M}"'] \arrow[r, "\SMan"] & \Sh(\Mfld;C) \arrow[d, "\restrict{(-)}{M}"] \\ 
			\PSh(M;C) \arrow[r, "\Sup_{M}"'] & \Sh(M;C) \period
		\end{tikzcd}
	\end{equation*} 
\end{lemma}

\begin{proof}
	In the case $ C = \Spc $, the claim follows from the fact that the forgetful functor
	\begin{equation*}
		\fromto{\Open(M)}{\Mfld}
	\end{equation*}
	satisfies the covering lifting property; see \cites[Proposition 7.1]{ClausenMathew:hyperdescent}[Proposition A.12]{arXiv:1803.01804}.
	The claim for general $ C $ follows from the claim for sheaves of spaces by applying the tensor product of presentable \categories and \cite[Lemma 1.18]{arXiv:2108.03545}.
\end{proof}

\begin{corollary}\label{cor:restrictionofconstantsheaf}
	Let $ C $ be a presentable \category, $ X \in C $, and $ M $ a manifold.
	Then we have a natural identification $ \restrict{\Gammaupperstar(X)}{M} = \Gammaupperstar_M(X) $ of the restriction of $ \Gammaupperstar(X) $ to $ M $ with the constant sheaf on $ M $ at $ X $.
\end{corollary}

\begin{proof}
	Note that by tensoring with the presentable \category $ C $, it suffices to prove the claim for $ C = \Spc $.
	In this case, note that by \Cref{lem:rescommuteswithsheaf} the functors
	\begin{equation*}
		\restrict{(-)}{M} \of \Gammaupperstar, \Gammaupperstar_M \colon \fromto{\Spc}{\Sh(M)}
	\end{equation*}
	are both left exact left adjoints.
	The claim follows from the fact that for \atopos $ \X $, the constant sheaf functor is the unique left exact left adjoint $ \fromto{\Spc}{\X} $ \HTT{Proposition}{6.3.4.1}.
\end{proof}

%-------------------------------------------------------------------%
%  Background on notions of completeness for higher topoi           %
%-------------------------------------------------------------------%

\subsection{Background on notions of completeness for higher topoi}\label{sec:completeness}

There are three notions of ``completeness'' for \atopos $ \X $:
\begin{enumerate}
	\item \textit{Hypercompleteness:} Whitehead's Theorem holds in $ \X $.

	\item \textit{Convergence of Postnikov towers:} Every object of $ \X $ is the limit of its Postnikov tower.

	\item \textit{Postnikov completeness:} $ \X $ can be recovered as the limit $ \lim_{n} \X_{\leq n} $ of its subcategories $ \X_{\leq n} \subset \X $ of $ n $-truncated objects along the truncation functors $ \trun_{\leq n} \colon \fromto{\X_{\leq n+1}}{\X_{\leq n}} $.
\end{enumerate}
While all of these properties hold for the \topos $ \Spc $ of spaces, they need not hold for a general \topos.
We have implications (3) $ \Rightarrow $ (2) $ \Rightarrow $(1), and none of the implications are reversible in general.
In this section we give a brief overview of hypercompletness as it plays a role in relating the Freed--Hopkins approach to differential cohomology from \cite{MR3049871} to the \categorical approach we have taken here.
Detailed accounts of hypercompleteness and Postnikov completeness can be found in \cite[\HTTsec{6.5}]{HTT} and \cite[\SAGsec{A.7}]{SAG}, respectively.

\begin{definition}
	Let $ \X $ be \atopos.
	An object $ U \in \X $ is \textit{$ \infty $-connected} if for every integer $ n \geq -2 $ the $ n $-truncation $ \trun_{\leq n}(U) $ of $ U $ is the terminal object of $ \X $.
	A morphism $ f \colon \fromto{U}{V} $ is \textit{$ \infty $-connected} if $ f \colon \fromto{U}{V} $ is an $ \infty $-connected object of the \topos $ \X_{/V} $.
\end{definition}

\begin{definition}\label{def:hypercompleteness}
	Let $ \X $ be \atopos.
	An object $ U \in \X $ is \textit{hypercomplete} if $ U $ is local with respect to the class of $ \infty $-connected morphisms in $ \X $.
	We write $ \X^{\hyp} \subset \X $ for the full subcategory spanned by the hypercomplete objects of $ \X $.
	\Atopos is \textit{hypercomplete} if $ \X^{\hyp} = \X $.
\end{definition}

\begin{nul}
	The \category $ \X^{\hyp} \subset \X $ is a left exact localization of $ \X $, hence \atopos \cite[\HTTpage{699}]{HTT}.
	Moreover, the \topos $ \X^{\hyp} $ is hypercomplete \HTT{Lemma}{6.5.2.12}.
\end{nul}

\begin{nul}
	The \topos $ \X^{\hyp} $ is the universal hypercomplete \topos equipped with a geometric morphism to $ \X $ \HTT{Proposition}{6.5.2.13}.
	For this reason we call $ \X^{\hyp} $ the \textit{hypercompletion} of $ \X $.
\end{nul}

\begin{observation}\label{obs:hypercompletecons}
	Let $ \X $ be \atopos.
	Then $ \X $ is hypercomplete if and only if the pullback functor $ \pupperstar \colon \fromto{\X}{\X^{\post}} $ is conservative.
\end{observation}

The standard way of working with sheaves of spaces on a site $ (S,\tau) $ in the language of model-categories is to use the Brown--Joyal--Jardine model structure on simplicial presheaves \cites{MR341469}{MR906403}.
However, this model structure only presents the hypercompletion of the \topos of sheaves of spaces on $ (S,\tau) $.

\begin{proposition}[{\HTT{Proposition}{6.5.2.14}}]\label{app.prop:HTT.6.5.2.14}
	Let $ (S,\tau) $ be a site.
	Then the underlying \category of the category of simplicial presheaves on $ S $ in the Brown--Joyal--Jardine model structure is equivalent to the \topos $ \Sh_{\tau}(S;\Spc)^{\hyp} $ of hypercomplete sheaves of spaces on $ S $.
\end{proposition}

\begin{definition}
	Let $ \X $ be \atopos.
	A \textit{point} of $ \X $ is a left exact left adjoint $ \xupperstar \colon \fromto{\X}{\Spc} $.
	Given an object $ U \in \X $ and point $ \xupperstar $ of $ \X $, we call $ \xupperstar(U) $ the \textit{stalk} of $ U $ at $ \xupperstar $.
\end{definition}

\begin{example}
	Let $ T $ be a topological space and $ t \in T $. 
	Then the stalk functor
	\begin{equation*}
		(-)_t \colon \fromto{\Sh(T)}{\Spc}
	\end{equation*}
	defines a point of $ \Sh(T) $.
\end{example}

\begin{definition}
	\Atopos $ \X $ \textit{has enough points} if a morphism $ f $ in $ \X $ is an equivalence if and only if for every point $ \xupperstar $ of $ \X $, the stalk $ \xupperstar(f) $ is an equivalence in $ \Spc $.
\end{definition}

\begin{example}\label{exm:enoughpointsishypercomplete}
	\Atopos with enough points is hypercomplete.
\end{example}

\begin{remark}
	The existence of enough points is incomparable with the convergence of Postnikov towers and is also incomparable with Postnikov completeness (both of which imply hypercompleteness).
\end{remark}

\begin{example}\label{ex:ShMPostnikovcomplete}
	Let $ M $ be a manifold. 
	Then the \topos $ \Sh(M) $ is Postnikov complete \cite[\HTTthm{Proposition}{7.2.1.10} \& \HTTthm{Theorem}{7.2.3.6}]{HTT}.
\end{example}

%-------------------------------------------------------------------%
%  A conservative family of points                                  %
%-------------------------------------------------------------------%

\subsection{A conservative family of points}\label{subsec:points}

In this section we show the stalks at the origins in $ \RR^n $ for $ n \geq 0 $ form a conservative family of points for the \topos $ \Sh(\Mfld) $ (\Cref{app.prop:hypercompleteness}).
This implies that the model structure on simplicial presheaves on $ \Mfld $ considered by Freed--Hopkins in \cite[\S 5]{MR3049871} presents the \topos $ \Sh(\Mfld) $.
We also present an observation of Hoyois that shows that the \topos $ \Sh(\Mfld) $ is Postnikov complete (\Cref{app.prop:postnikovcompleteness}).

We begin by discussing the stalk of a sheaf on $ \Mfld $ at a point of a manifold.

\begin{construction}\label{cons:stalk}
	Let $ M $ be a manifold and $ x \in M $.
	In light of \Cref{lem:rescommuteswithsheaf}, the composition of the restriction to $ M $ with the stalk at $ x $ defines a left exact left adjoint
	\begin{equation*}
		\begin{tikzcd} 
			\Sh(\Mfld;C) \arrow[r, "\restrict{(-)}{M}"] & \Sh(M;C) \arrow[r, "{(-)_x}"] & C \comma
		\end{tikzcd}
	\end{equation*}
	which we denote by $ \xupperstar $.
	Given a sheaf $ E $ on $ \Mfld $, we call $ \xupperstar(E) $ the \textit{stalk} of $ E $ at $ x \in M $.
\end{construction}

\begin{observation}\label{obs:restricttoopen}
	Let $ M $ be a manifold and $ j \colon \incto{U}{M} $ an open embedding.
	Then, by definition, the triangle
	\begin{equation*}
		\begin{tikzcd}[sep=3em]
			\Sh(\Mfld;C) \arrow[r, "\restrict{(-)}{M}"] \arrow[dr, "\restrict{(-)}{U}"'] & \Sh(M;C) \arrow[d, "\jupperstar"] \\
			& \Sh(U;C) 
		\end{tikzcd}
	\end{equation*} 
	commutes.
	Thus for any $ x \in U $, then there is a canonical identification of the stalk functor $ \fromto{\Sh(\Mfld;C)}{C} $ at $ x \in U $ with the stalk functor at $ j(x) \in M $. 
\end{observation}

Recall that for each integer $ n \geq 0 $, write $ 0_n \in \RR^n $ for the origin (\Cref{ntn:origin}).

\begin{proposition}\label{app.prop:hypercompleteness}
	Let $ C $ be a compactly assembled \category.
	Then the set of stalk functors
	\begin{equation*}
		\{0_n\upperstar \colon \Sh(\Mfld;C) \to C \}_{n \geq 0}
	\end{equation*}
	is jointly conservative. 
	In particular, the \topos $ \Sh(\Mfld) $ is hypercomplete.
\end{proposition}

\begin{proof}
	In light if \cite[Lemma 2.12]{arXiv:2108.03545}, it suffices to treat the case $ C = \Spc $.
	In this case, first note that the family of restriction functors
	\begin{equation*}
		\restrict{(-)}{M} \colon \fromto{\Sh(\Mfld)}{\Spc} 
	\end{equation*}
	for $ M \in \Mfld $ is jointly conservative (\Cref{obs:restrictiontoamanifold}).
	For each manifold $ M $, the \topos $ \Sh(M) $ is a hypercomplete \topos and the points of $ M $ provide conservative family of points for $ \Sh(M) $ \HTT{Corollary}{7.2.1.17}.
	Thus the stalk functors
	\begin{equation*}
		\xupperstar \colon \Sh(\Mfld) \to \Spc
	\end{equation*}
	for all $ M \in \Mfld $ and $ x \in M $ form a conservative family of points for $ \Sh(\Mfld) $.
	To conclude, note that for every manifold $ M $ and point $ x \in M $, there exists an open embedding $ j \colon \incto{\RR^n}{M} $ such that $ j(0_n) = x $ and apply \Cref{obs:restricttoopen}.
\end{proof}

We now give a quick argument showing that the \topos $ \Sh(\Mfld) $ is Postnikov complete.
We learned the following argument from Hoyois; it is a slight refinement of the argument for the convergence of Postnikov towers that Hoyois gave in \cite{MO:130999}.

\begin{proposition}\label{app.prop:postnikovcompleteness}
	The \topos $ \Sh(\Mfld) $ is Postnikov complete.
\end{proposition}

\begin{proof}
	Since $ \Sh(\Mfld) $ is hypercomplete, by \Cref{obs:hypercompletecons} it suffices to show that the right adjoint $ \plowerstar \colon \fromto{\Sh(\Mfld)^{\post}}{\Sh(\Mfld)} $ is fully faithful.
	That is, we need to show that for every collection of objects $ \{F_n\}_{n \geq -2} $ of $ \Sh(\Mfld)$ equipped with compatible equivalences $ \equivto{\trun_{\leq n}(F_{n+1})}{F_n} $, and integer $ k \geq -2 $, the natural morphism
	\begin{equation}\label{eq:truncationmor}
		\fromto{\trun_{\leq k}\paren{\lim_{n \geq -2} F_n}}{F_k}
	\end{equation}
	is an equivalence.
	To see this, note that since the restriction functors 
	\begin{equation*}
		\{\restrict{(-)}{M} \colon \fromto{\Sh(\Mfld)}{\Sh(M)}\}_{M \in \Mfld}
	\end{equation*}
	are jointly conservative and commute with limits and truncations, it suffices to show that the morphism \eqref{eq:truncationmor} becomes an equivalence after restriction to each manifold $ M $.
	This last claim follows from the fact that the \topos $ \Sh(M) $ is Postnikov complete (\Cref{ex:ShMPostnikovcomplete}).
\end{proof}

We finish this section by proving that $ \Sh(\Mfld) $ is equivalent to the \topos of sheaves on the smaller site $ \Euc \subset \Mfld $ spanned by the Euclidean spaces (\Cref{def:Euclideansite}).
Note that since every manifold admits a cover by Euclidean spaces, the Euclidean site is a \textit{basis} for the Grothendieck topology on $ \Mfld $ (see \cite[\S B.6]{Ultracategories} for more about bases for Grothendieck topologies).

\begin{corollary}\label{app.lem:hypersheavesonCart}
	Let $ C $ be a presentable \category.
	Then restriction of presheaves
	\begin{equation*}
		\restrict{(-)}{\Eucop} \colon \fromto{\Sh(\Mfld;C)}{\Sh(\Euc;C)}
	\end{equation*}
	is an equivalence of \categories.
	The inverse is given by right Kan extension along the inclusion $ \incto{\Eucop}{\Mfldop} $.
\end{corollary}

\begin{proof}
	Since $ \Sh(\Euc;C) $ and $ \Sh(\Mfld;C) $ are the tensor products of presentable \categories
	\begin{equation*}
		\Sh(\Euc;C) \equivalent \Sh(\Euc) \tensor C \andeq \Sh(\Mfld;C) \equivalent \Sh(\Mfld) \tensor C \comma
	\end{equation*}
	it suffices to treat the case where $ C = \Space $ is the \category of spaces.
	In this case, since the \topos $ \Sh(\Mfld) $ is hypercomplete (\Cref{app.prop:hypercompleteness}), the claim follows from the fact that $ \incto{\Euc}{\Mfld} $ is a basis for the topology on $ \Mfld $; see \cite[Appendix A]{arXiv:2001.00319} or \cite[Corollary 3.12.13]{exodromy}.
\end{proof}

\resumesections % Goes back from appendix numbering to normal section numbering
\newpage

%-------------------------------------------------------------------%
%  Characteristic classes                                           %
%-------------------------------------------------------------------%

\part{Characteristic classes}

\numberwithin{equation}{part}
%!TEX root = ../diffcoh.tex

%-------------------------------------------------------------------%
%-------------------------------------------------------------------%
%  Part II Introduction                                             %
%-------------------------------------------------------------------%
%-------------------------------------------------------------------%

\label{part:charclasses}\label{char_class_part}

The objective of this portion of the notes is to construct, study, and use refinements of standard characteristic classes to differential cohomology. 

Historically, differential characteristic classes were studied by Cheeger and Simons \cite{MR827262}. This view is covered in \cref{DifferentialCharacteristicClasses}.

The modern approach uses the machinery of sheaves on manifolds developed in \Cref{part:basics} of these notes. 
Given a Lie group $G$, we consider three different, but related, sheaves $\Mfldop \to \Spc $ on $\Mfld$:
\begin{enumerate}
	\item The constant sheaf at the classifying space $\BG$ of $ G $ (\Cref{ntn:classifyingspaceBG}).
	We simply denote this sheaf by $ \BG $.

	\item The sheaf $ \BunG $ sending a manifold $ M $ to the groupoid of principal $ G $-bundles on $ M $ (\Cref{ex:BunG}).

	\item The sheaf $ \BunGnabla $ sending a manifold $ M $ to the groupoid of principal $ G $-bundles on $ M $ with connection (\Cref{ex:BunGnabla})..
\end{enumerate}
Characteristic classes live in the de Rham cohomology of these sheaves. 

\begin{definition}
	Let $S$ be a sheaf on manifolds. 
	The \emph{de Rham cohomology} of $S$ is $\Omegabullet(S)$.
\end{definition}

\noindent For example, the de Rham cohomology $\Omegabullet(\BG)$ of the constant sheaf $\BG$ is where ordinary characteristic classes live.

\begin{remark}
	Given a manifold $M$, one can recover the differential cohomology $\Hcech^k(M)$ by taking the $k$-th de Rham cohomology. %%%%%finish!
\end{remark}

The de Rham cohomology of $\BunGnabla $ is studied in \cref{WorkofFreedHopkins}. The de Rham cohomology $\Omegabullet(\BunGnabla )$ classifies characteristic classes for $G$-bundles with connections. In \cref{WorkofFreedHopkins}, we give a proof of the main theorem of \cite{FreedHopkins}. The theorem is as follows,

\begin{theorem}[(Freed--Hopkins)]\label{FreedHopkinsThm}
	The Chern--Weil homomorphism induces an isomorphism
	\[
		(\Sym^\bullet\gdual)^G \isomorphism \Omegabullet (\BunGnabla ) \period
	\]
\end{theorem}

\noindent Thus the Chern--Weil construction, reviewed in \cref{ChernWeilTheory}, produces all  characteristic classes for bundles with connection. The set up for the proof of Theorem \ref{FreedHopkinsThm} uses tools similar to the Cartan model for equivariant de Rham cohomology, which we review in \cref{EquivariantdeRhamCohomology}. 
Related work on Borel equivariant differential cohomology was done by Redden \cite{MR3649665}.

The de Rham cohomology of $\BbulletG$ is a bit more complicated. The tools we use to compute $\Omegabullet \BbulletG$ originate in Bott's paper \cite{BottsPaper}. In \cref{BottsMethod}, we review the techniques used in \cite{BottsPaper} including continuous cohomology and the van Est theorem. 
The takeaway of \cref{BottsMethod} is the following theorem of Bott:

\begin{theorem}[(Bott)]
	The continuous cohomology $\Hcont^{p-q}(G;\Sym^q(\gdual))$ is isomorphic to the de Rham cohomology group $\H^{p}(\BbulletG;\Omega^q)$:
	\[
		\H^{p}(\BbulletG;\Omega^q) \isomorphic \Hcont^{p-q}(G;\Sym^q(\gdual)) \period
	\]
\end{theorem}

\noindent We will really only use Bott's theorem in degrees $p-q\leq 0$. 

In \cref{LiftsofChernClasses}, the results of \cite{BottsPaper} are applied to provide lifts of Chern classes to differential cohomology. 
In particular, we will see there exists multiple lifts of each Chern class $c_i$ to $\H^{2n}(\BbulletGL{n}{\CC}; \ZZ_\CC(n))$. The collection of lifts is determined by the following result, credited by Hopkins to Bott:

\begin{theorem}
	There is a pullback square
	\begin{equation*}
		\begin{tikzcd}
			\H^{2n}(\BbulletGL{k}{\CC}; \ZZ_\CC(n))\arrow[rr]\arrow[d] & & \H^{2n}(\BU_k;\ZZ)\arrow[d] \\
			\H^n(\BU_k\times \BU_k;\CC)\arrow[rr, "\textup{diagonal}^*"'] & & \H^{2n}(\BU_k;\CC) \period
		\end{tikzcd}
	\end{equation*}
\end{theorem}

\noindent A real analogue of this theorem provides lifts of the Pontryagin classes.

\begin{remark}\label{rmk-OnOffDiagonal}
Note that differential cohomology $\H^i(-;\bb{Z}(j))$ is bigraded. 
The differential lifts of characteristic classes discussed in \cref{WorkofFreedHopkins} live in bidegree where $i=j$. 
We refer to these classes as ``on-diagonal." 
The classes defined in \cref{LiftsofChernClasses} live in bidegree where $i=2j$, 
and we call these ``off-diagonal" classes. 
Notationally, for a class $c$, we use $\hat{c}$ to denote an on-diagonal differential lift 
and $\ctilde$ for an off-diagonal lift.
\end{remark}

As an application of this construction, in \cref{VirasoroAlgebra} we explain how a differential lift of the first Pontryagin class $\ptilde_1\in\H^4(\BSL(\RR);\ZZ(2))$ can be used to produce the Virasoro group. The Virasoro group is a certain central extension of $\Diffplus(\Circ)$ by $\Uup_1$,
\[
	\Uup_1\to\Vir\to\Diffplus(\Circ) \period
\]
The construction of $\Vir$ uses the fiber integration for differential cohomology covered in \cref{FiberIntegration} and pullback along the classifying map of a certain bundle. 
This process is outlined in \cref{LiftsofChernClasses} and covered in depth in \cref{VirasoroAlgebra}. 
Note that there are multiple lifts of $p_1$ to differential cohomology. 
We obtain criterion for which lift $\ptilde_1$ could correspond to the Virasoro algebra central extension, but we do not pin down which lift works.

As far as we know, the material in \cref{LiftsofChernClasses} and \cref{VirasoroAlgebra} does not appear elsewhere in the literature, aside from the underpinning in \cite{BottsPaper}. 
The new ideas here are due to Dan Freed, Mike Hopkins, and Constantin Teleman.

\numberwithin{equation}{subsection}
\newpage
%!TEX root = ../diffcoh.tex

%-------------------------------------------------------------------%
%-------------------------------------------------------------------%
%  Chern–Weil theory                                                %
%-------------------------------------------------------------------%
%-------------------------------------------------------------------%

\section{Chern--Weil theory}\label{ChernWeilTheory}
\textit{by Greg Parker}

As this talk is a review of standard material, 
many technical results are stated without proof. 
For more detailed review, including proofs, the reader should consult \cite[Appendix C]{MilnorStasheff} for a review of connections and Chern--Weil theory for vector bundles, 
\cite[Chapter II]{KobayashiNomizu} or \cite[Chapter 2]{JohnRoe} for the theory of connections on principal bundles, and 
\cite[Chapter XII]{KobayashiNomizuII} for Chern--Weil theory for principal bundles.

%-------------------------------------------------------------------%
%-------------------------------------------------------------------%
%  Motivation                                                       %
%-------------------------------------------------------------------%
%-------------------------------------------------------------------%

\subsection{Motivation}

To begin, let's recall
\begin{theorem}[(Gauss--Bonnet)]
	Let $(\Sigma,g)$ be a compact, oriented, Riemannian $ 2 $-manifold without boundary. Let $\chi(\Sigma)$ be its Euler characteristic. 
	Then
	\[
		\int_\Sigma \kappa \d A=2\pi\chi(\Sigma) \period
	\]
\end{theorem}
Here $\kappa$ is the Gaussian curvature defined as follows. 
If $R_{ij}\d x_i\d x_j$ is the Riemann curvature tensor, locally
\[R=\begin{pmatrix}
0 & R_{12}\\
-R_{21} & 0
\end{pmatrix} \d x_1\wedge \d x_2\]
and $\kappa= R_{12}$. So we can rewrite the above as
\[
	\langle[\sqrt{\det(R)}],[\Sigma]\rangle=\int_\Sigma\sqrt{\det(R)}=2\pi\chi(\Sigma)=\langle 2\pi e(\Tan\Sigma),[\Sigma]\rangle \comma
\]
where $e(\Tan\Sigma)$ is the Euler class of $\Sigma$ and the brackets on the right-hand side  denote the pairing 
\[
	\H^2(\Sigma;\RR)\otimes \H_2(\Sigma;\RR)\to\RR\period
\]

Thus we observe $\sqrt{\det(R)}$, a polynomial in the curvature, captures information about the topology of $\Sigma$ and its tangent bundle $\Tan\Sigma$. Chern--Weil theory (which was actually the original formulation/theory of characteristic classes) generalizes the above to higher dimension and arbitrary bundles.

%-------------------------------------------------------------------%
%-------------------------------------------------------------------%
%  Connections and curvature                                        %
%-------------------------------------------------------------------%
%-------------------------------------------------------------------%

\subsection{Connections and curvature}

In order to formulate things correctly, we will need to recall some facts about connections and curvature, 
both for vector bundles and for principal bundles. 

\begin{convention}
	Throughout this talk, 
	let $M$ be a closed $n$-manifold, $\pi\colon E\to M$ a rank $k$ real or complex vector bundle with structure group $G = \Or_k $ or $ G = \Uup_k $. 
	Denote the real (or Hermitian) inner product by $\langle -, -\rangle$. 
	Let $\g$ denote the Lie algebra of $G$.

	Also, let $K$ be a Lie group and $p\colon Q\to X$ be a principal $K$-bundle. 
	Let $\kfrak$ denote the Lie algebra of $K$.
\end{convention}

%-------------------------------------------------------------------%
%  For vector bundles                                               %
%-------------------------------------------------------------------%

\subsubsection{For vector bundles}\label{subsubsec:connections_for_vbs}

We would like to differentiate sections $\psi\colon M\to E$. 
The problem is $\psi_{x(t)}$ for a path $x(t)\subset M$ all live in different vector spaces: $E_{x(t)}$, respectively, so we must find a way to ``connect'' them. 

View $\psi_{x(t)}$ as a path in the total space. The derivative (intuitively) is the vertical component of $\frac{\partial\psi}{\d t}$. Think of $f\colon\RR\to\RR$, then $\frac{\d f}{\d t}$ is the $y$-coordinate of the graph inside $\RR^2$. 
To define this precisely we need to choose a splitting
\[
	\Tan E\simeq VE\oplus HE
\]
into the ``vertical'' and ``horizontal'' subbundles.
The vertical piece $VE=\ker \d\pi$ is canonical, and the horizontal piece $HE$ is not. Such a splitting is called a \emph{connection}. 
Once we choose a connection, 
we get an isomorphism $\d\pi\colon HE\to TM$. 
So given $e\in E_{x(t)}$ we can lift $\dot{x}$ (a vector field along $x(t)$) to one $X^e_H\subset HE$. 
Then the flow is a path in $E$ projecting to $x(t)$, which is the \emph{parallel transport}, denoted $\varphi_te\in E_{x(t\times I)}$. Then
\[
	\varphi_{-t}\psi(t)\in E_{x(0)}
\]
for all $t$, so we can differentiate. 
The covariant derivative (with respect to our chosen connection) in the $\dot{x}(0)$ direction at $x(0)$ is 
$\frac{\d}{\d t}\vert_{t=0}\varphi_{-t}\psi_{x(t)}$. 
Thus we get an operator
\[
	\d_{A}\text{ or }\nabla^A\colon\Gamma(M;E)\to\Gamma(M;\Tstar M\otimes E)
\]
associated to a connection $A$, called the \emph{covariant derivative}. 
Here, $\nabla^A$ eats a vector field $X\in\Gamma(M;\TanM)$ and gives the derivative in that direction at each point. It satisfies
\begin{itemize}
	\item \textit{$\Cinf$-linearity in the direction of the derivative:} $\nabla_{fX}^A\psi=f\nabla_X^A\psi$, and 

	\item \textit{Leibniz rule:} $\nabla^Af\psi=\d f\otimes\psi+f\nabla\psi$.
\end{itemize}

The existence of connections is preserved under various bundle constructions.  

\begin{proposition}\label{properties of connections}
	Given $\nabla^A$ on $E$, $\nabla^B$ on $F$ we get connections:
	\begin{itemize}
		\item $\nabla^{A^*}$ on the dual bundle $E^*$.

		\item $\nabla^{AB}$ on the tensor product $E\otimes F $ by the formula 
		\[
			\nabla^{AB}(\varphi\otimes\psi)=\nabla^A\varphi\otimes\psi+\varphi\otimes\nabla^B\psi \period 
		\]

		\item if $F\colon M\to N$ and $E\to N$ then $F^*(\nabla^A)$ is a connection on  $f^*E$ by
			\[
				(F^*\nabla^A)_X\psi(m)\colonequals \nabla^A_{F_*X}\psi(f(m))\in E_{F(m)}=F^*E_m \period 
			\]
	\end{itemize}
\end{proposition}

\begin{proposition}\label{difference_of_connections}
	Two connections differ by a $1$-form valued in $\End(E)$. 
	In particular, the set of connections form an affine, and hence contractible, space.
\end{proposition}

\begin{remark}
	Thus one might expect invariants defined using them (if discrete) to not depend on the choice of connection.
\end{remark}

\begin{proof}
	Let $A$ and $A'$ be two connections on the bundle $\pi\colon E\to M$. 
	For $f\in\Cinf(X)$ and $\psi$ a section of $\pi$, we have
	\begin{align*}
		(\nabla^A-\nabla^{A'})(f\psi)&=\d f\otimes\psi+f\nabla^A\psi-\d f\otimes\psi-f\nabla^{A'}\psi \\
		&=f(\nabla^A-\nabla^{A'})\psi
	\end{align*}
	is $\Cinf$-linear with values in $\Gamma(\Tstar M\otimes E)$ so $\nabla^A-\nabla^{A'}\in\Omega^1(\End(E))$.
\end{proof}

\begin{example}
	On the trivial rank $k$-bundle $\underline{\RR}_M$ on $M$, the exterior derivative 
	\[
		\d\colon\Gamma(M;\underline{\RR}_M)\to\Omega^1(M)
	\]
	 is a connection.
\end{example}

\begin{example}
	In a local trivialization (by \Cref{difference_of_connections}) we can always write $\nabla=\d+A$, where $A\in\Omega^1(\End(E))$. 
	That is $A=A_1\d x_1+\cdots A_n\d x_n$ for $A_i$ matrices, and
	\[
		\nabla_i\psi=\frac{\del \psi}{\del x_i}+A_i\psi \period 
	\]
\end{example}

\begin{example} 
	On $\End(E)=E^*\otimes E$, the induced $\nabla$ from \Cref{properties of connections} is 
	\[
		\nabla B=\d B+[A,B]
	\]
	in a trivialization.
\end{example}

Define a connection as \emph{compatible with $\ang{-,-}$} if
\[
	\d\ang{\psi,\varphi} = \ang{\nabla\psi,\varphi} + \ang{\psi,\nabla\varphi} \period
\]
Note that for compatible $\nabla$, $A$ will be in $\Omega^1(\ofrak(E))$ or $\ufrak(E)$. 

\begin{remark}
	A fancy way of saying this is $\ang{-,-}\in E^*\otimes E^*$ has $\nabla =0$.
\end{remark}

\begin{lemma}
	Every bundle $E$ has a connection compatible with $\ang{-,-}$.
\end{lemma}

\begin{proof}
	Locally, connections of the form $ \d+A$ are compatible with $\ang{-,-}$ if $A$ is in $\Omega^1(\ofrak(E))$ or $\Omega^1(\ufrak(E))$. 
	This gives existence locally. Using a partition of unity, one obtains the desired connection globally.
\end{proof}

%-------------------------------------------------------------------%
%  For principal bundles                                            %
%-------------------------------------------------------------------%

\subsubsection{For principal bundles}\label{subsec-principal}

For a principal $K$-bundle $p\colon Q\to X$, 
the space of vertical tangent vectors $\ker(\d p)$ of $Q$ gives a short exact sequence
\begin{align}\label{ses:principal}
0\to\ker(p_*)\to \Tan Q\to p^*\Tan X\to 0
\end{align}
of vector bundles over $P$. 
As in the vector bundle situation, a connection will be a way of considering \emph{horizontal} tangent vectors.

\begin{definition}
	A \emph{principal connection} on $p\colon Q\to X$ is a splitting of the exact sequence \eqref{ses:principal}. 
\end{definition}

The kernel $\ker(\d p)$ can be identified with the trivial bundle with fiber the tangent space of the fiber $K$ of $p$. 
That is, we have an isomorphism
\begin{equation*}
	\ker(\d p)\isomorphic Q\times \kfrak \period
\end{equation*} 
A splitting of \cref{ses:principal} is equivalent to a section of the map $\ker (\d p)\to \Tan Q$. 
Using the identification $\ker(\d p)\isomorphic Q\times \kfrak$, 
a section $\Tan Q\to \ker(\d p)$ is equivalent to a section of $\Tstar Q\otimes( \kfrak\times Q)$; i.e., a one form with coefficients in $\kfrak$. 

\begin{definition}
	Let $p\colon Q\to X$ a principal $K$-bundle with principal connection. 
	The \emph{connection $ 1 $-form} is the principal connection is the one form  $\omega\in \Omega^1(Q;\kfrak)$ 
	corresponding the splitting of \cref{ses:principal}.
\end{definition}

Note that $\kfrak$ acts on $\kfrak$ in two ways: 
by right $m_r$, 
and by conjugation $\Ad_\kfrak$. 

\begin{lemma}\label{connection_1-form_equivariant}
	Let $p\colon Q\to X$ a principal $K$-bundle with principal connection $ 1 $-form $\omega$. 
	Then $\omega$ is $K$-equivariant,
	\[
		\Ad_\kfrak(m_r(\omega))=\omega
	\]
	and for $\xi\in\kfrak$ with associated vector field $X_\xi$, we have $\omega(X_\xi)=\xi$. 
\end{lemma}

\begin{remark}
	A connection on a principal bundle gives rise to a vector bundle connection on any associated vector bundle. 
	Likewise, a $K$-compatible connection on a vector bundle $E$ gives rise to a connection on the $K$-frame bundle, 
	and these operations are inverses. 
	The horizontal distribution on $\Tan Q$ complementing $\ker(p_*)$ in \cref{ses:principal} 
	is obtained from a vector bundle connection as directions of the infinitesmal parallel transport at a point. 
	In the opposite direction, the parallel transport of frames on $Q$ naturally gives a parallel transport of section of the vector bundle. 
	Alternatively, in local coordinates the connection form is just $\d+\omega$ for $\omega$ the ad-equivariant connection form on $Q$.
\end{remark}

%-------------------------------------------------------------------%
%-------------------------------------------------------------------%
%  Curvature                                                        %
%-------------------------------------------------------------------%
%-------------------------------------------------------------------%

\subsection{Curvature}

%-------------------------------------------------------------------%
%  For vector bundles                                               %
%-------------------------------------------------------------------%

\subsubsection{For vector bundles} 

Given two vector fields $X,Y\in\Gamma(M;\TanM)$, 
the maps $\nabla_X$ and $\nabla_Y$ need not commute; i.e.,
\[\nabla_X\nabla_Y-\nabla_Y\nabla_X\neq 0\]
Geometrically, since these were defined by flowing along horizontal lifts, $\tilde{X},\tilde{Y}$, this is a question about non-commuting flows; i.e., $[\tilde{X},\tilde{Y}]$. 
In particular, if the horizontal bundle $HE$ is integrable, then $[\tilde{X},\tilde{Y}]=0$ so the flows (and hence $\nabla_X,\nabla_Y$) commute. Thus the \emph{curvature} 
\[
	F_A(X,Y)(\psi)\colonequals[\nabla_X,\nabla_Y](\psi)-\nabla_{[X,Y]}(\psi)
\]
is a measure of the integrability of $HE\subseteq E$. 
Here $A$ is such that locally we have $\nabla=d+A$. 

We get a local description of the curvature by 
	\[
		F_A=\d A+A\wedge A\period 
	\]
	In other words, 
	\[
		F_A=F_A^{ij}\d x^i\wedge \d x^j
	\]
	with
	\[
		F^{ij}_A=\del_iA_j-\del_jA_i+A_iA_j-A_jA_i \period
	\]

\begin{claim}
The curvature $F_A$ defines a 2-form with values in the endomorphism bundle, \[F_A\in \Omega^2(M;\End(E))\period\]
In particular, the curvature is $\Cinf$-linear, $F_A(f\psi)=fF_A\psi$. 
\end{claim}

\begin{proof}
This follows from the Leibniz rule for connections. 
\end{proof}

For $F_A\in\Omega^2(\End(E))$. 
We have 
\[
	\d_{A}\colon\Omega^2(\End(E))\to\Omega^3(\End(E))
	\]
by $\alpha\otimes B\mapsto \d\alpha\otimes B+\alpha\otimes\nabla B$. 

\begin{theorem}[(Bianchi identity)]
The exterior derivative of the curvature vanishes, $\d_{A}F_A=0$.
\end{theorem}

%-------------------------------------------------------------------%
%  For principal bundles                                            %
%-------------------------------------------------------------------%

\subsubsection{For principal bundles}\label{sssec:CW_principal}

The wedge product of $\omega\in \Omega^1(Q;\kfrak\otimes\kfrak)$ with itself is an element of $\Omega^2(Q;\kfrak)$. 
The Lie bracket on $\g$ induces a map 
\[[-]\colon \Omega^2(Q;\kfrak\otimes\kfrak)\to\Omega^2(Q;\kfrak)\period\]

\begin{definition}
Let $p\colon Q\to X$ a principal $K$-bundle with principal connection $ 1 $-form $\omega$. 
The \emph{curvature} of $\omega$ is 
\[\Omega=\d\omega+[\omega\wedge \omega]\]
in $\Omega^2(Q;\kfrak)$.
\end{definition}

Consider $\kfrak$ as a $K$-module with the adjoint action. 
Let $\kfrak_Q\to X$ denote the adjoint bundle ${\kfrak_Q=Q\times_K\kfrak}$. 

\begin{lemma}
Let $p\colon Q\to X$ a principal $K$-bundle with connection. 
Let $\Omega$ be its curvature. 
Then $\Omega$ descends to a $2$-form $\tilde{\Omega}\in\Omega^2(X;\g_Q)$. 
\end{lemma}

\begin{example}
Take $K=\GL_n$ so that $Q$ has an associated rank $n$ vector bundle $V\to X$. 
The adjoint bundle can be identified with the endomorphism bundle $\End(V)$. 
Under this identification, 
a principal connection on $Q$ corresponds to a connection on the vector bundle $V\to X$, and 
the curvature $\tilde{\Omega}$ from a principal connection on $Q$ corresponds to the curvature of $V\to X$. 
\end{example}

\begin{theorem}[(Bianchi identities)]
	We have $\d\Omega+[\omega\wedge\Omega]=0$ and $\d\Omega=0$. 
\end{theorem}

%-------------------------------------------------------------------%
%-------------------------------------------------------------------%
%  Invariant polynomials                                            %
%-------------------------------------------------------------------%
%-------------------------------------------------------------------%

\subsection{Invariant polynomials}

\subsubsection{For vector bundles}

In Gauss--Bonnet we used $\sqrt{\det}$ to turn the $R\in\Omega^2(\so(\Tan\Sigma))$ into an $\RR$-valued form to integrate. In general, since $F_A$ isn't basis-invariant we want a map $P\colon\g\to\RR$ (for $G=\SO_k$ or $\SU_k$) invariant under $\Ad$. If $P$ is a polynomial, we say it is an \emph{invariant polynomial}. 
The space of $\Ad$-invariant polynomials on $\g$ is $\Sym^\bullet(\gdual)^{G}$. 

\begin{example}
	Both $\Tr$ and $\Det$ are $\Ad$-invariant. 
\end{example}

\noindent Thus given $P,A$ we get an $\RR$-valued form $P(F_A)\in\Omega^*(M;\RR)$.

\begin{proposition}
	The form $P(F_A)$ is closed, $\d P(F_A)=0$. 
	Hence we get a homomorphism
	\[
		\Sym^\bullet(\gdual)^G\to \HdR^*(M;\RR) \period  \index[terminology]{Chern--Weil homomorphism} 
	\]
\end{proposition}
 The map above is called the \emph{Chern--Weil homomorphism}, or sometimes just the Weil homomorphism.

\begin{proof}
Write $P(\xi)=\sum_I P_I\xi_{i_1},\dots,\xi_{i_N}$. 
Since $P$ is $\Ad$-invariant, for $g_t=\exp(t\eta)$, we have
		\[
			P(\xi)=P(\Ad_{g_t}\xi)
		\]
		so
		\begin{align*}
			0&=\frac{\d}{\d t}P(\xi)\\
			&=\frac{\d}{\d t}\sum_I P_I(\Ad_{g_t}\xi)_{i_1}\cdots(\Ad_{g_t}\xi)_{i_N}\\
			&=\sum_{I,k}P_I\xi_{i_1}\cdots\xi_{i_{k-1}}[\eta,\xi]_{i_k}\cdots\xi_{i_N} \period
		\end{align*}
Writing $F_A=\sum F_A^i$, we have
		\begin{align*}
			\d P(F_A)&=\d\paren{\sum_IP_IF^{i_1}\wedge\cdots\wedge F^{i_N}} \\
			&=\sum_{I,k}P_IF^{i_1}\wedge\cdots \wedge dF^{i_k}\wedge\cdots \wedge F^{i_N}+\sum_{I,k}P_IF^{i_1}\wedge\cdots\wedge [A,F_A]_{i_k} \wedge\cdots\\
			&=\sum_{I,k}P_IF^{i_1}\wedge\cdots\wedge (\d_{A}F_A)_{i_k}\wedge\cdots\wedge F_{i_N}\\
			&=0 \period \qedhere
		\end{align*}
\end{proof}

\begin{proposition}[(invariance)]
	The class $[P(F_A)]$ satisfies the following properties.
	\begin{enumerate}[{\upshape (1)}]
		\item $[P(F_A)]$ is independent of $A$.

		\item $[P(F_A)]$ is independent of $\ang{-,-}$.

		\item If $E\simeq E'$ then $[P(F_A)]=[P(F_{A'})]$. 
		The \emph{characteristic class} of $E$ is $[P(F_A)]\in \H^*$.
	\end{enumerate}
\end{proposition}

\begin{proof}[Proof Sketch]
	For (1), take $A,A'$ and set $\nabla_A-\nabla_{A'}=B$. 
	Define
	\begin{equation*}
		A_t \colon E\times I\to M\times I
	\end{equation*}
	by $\nabla_A+tB$. 
	Then $P(F_{A_t})\in\Omega^\bullet(M\times I;\RR)$, and $i_0\colon M\to M\times\{0\}$ has $i_0^*P(F_{A_t})=P(F_A)$ for some $i_1,A'$. 
	But $i_0,i_1$ are homotopic.

	The proof of (2) is similar. 

	For (3), use the pullback connection plus the independence of $A$.
\end{proof}

%-------------------------------------------------------------------%
%  For principal bundles                                            %
%-------------------------------------------------------------------%

\subsubsection{For principal bundles}

We have an analogous story for principal bundles, 
using the corresponding notions of curvature and Bianchi identities. 

\begin{proposition}\label{Chern--Weil map for Principal}
	Let $Q\to X$ be a principal $K$-bundle with curvature $\Omega$. 
	The assignment $P\mapsto P(\Omega)$ determines a map
	\[
		\Sym(\kfrak^\vee)^{K} \to\OmegadR^\bullet(X)
	\]
	that descends to a map on cohomology. 
\end{proposition}

%-------------------------------------------------------------------%
%-------------------------------------------------------------------%
%  Examples                                                         %
%-------------------------------------------------------------------%
%-------------------------------------------------------------------%

\subsection{Examples}

Now the fun part: choose different $P$ and see what we get.

%-------------------------------------------------------------------%
%  Chern classes                                                    %
%-------------------------------------------------------------------%

\subsubsection{Chern classes}
\label{chern_const}

Consider the polynomial $P=\Det(\lambda\id{}-\frac{1}{2\pi i}X)\colon\ufrak_k\to\RR$. 
Then expanding out, we get 
\[P=\lambda^k-c_1(X)\lambda^{k-1}+c_2(X)\lambda^{k-2}+\cdots\]
for $c_k$ polynomials in $X$. Define the characteristic class $c_k$ in $\H^{2k}$ obtained from $P$ to be the $k$-th Chern class. Explicitly
\begin{align*}
c_k(F_A)&=\frac{\Tr(F_A^{\wedge k})}{(2\pi i)^k}\\
&=1-\frac{1}{2\pi i}\Tr(F_A)+\frac{\Tr(F_A\wedge F_A)-\Tr(F_A)^2}{8\pi^2}-\cdots \period
\end{align*}
\begin{remark}
It's immediate that $c_1=0$ for an $\SU_n$-bundle since $\su_n$ is traceless. In fact, one can show $c_k$ are a basis for $\Ad$-invariant polynomials so this is a complete list.
\end{remark}

%-------------------------------------------------------------------%
%  Pontryagin classes                                               %
%-------------------------------------------------------------------%

\subsubsection{Pontryagin classes}
\label{pontrjagin}
Consider the polynomial $P$ from 
\[
	\Det(\lambda\id{}-\frac{1}{2\pi}X)\colon\ofrak_k\to\RR \period
\]
Expanding out, we get
\[
	P=\lambda^k-g_1(X)\lambda^{k-1}+\cdots \period
\]
Since $\ofrak_k$ is skew-symmetric $g_{\odd}=0$ and $g_{2k}=p_k(E)$ is the $k$-th Pontryagin class. 
For example, we have
\[
	p_1=-\frac{\Tr(F_A\wedge F_A)}{8\pi^2}
\]
and
\[
	p_2=\frac{\Tr(F_A\wedge F_A)^2-2\Tr(F_A\wedge\cdots\wedge F_A)}{128\pi^4} \period
\]

%-------------------------------------------------------------------%
%  Euler class                                                      %
%-------------------------------------------------------------------%

\subsubsection{Euler class}
\label{euler_const}

If $k$ is even, there is the Pfaffian $\Pf\colon\ofrak(2k)\to\RR$ with $\Pf(X)^2=\det(X)$. Then the \emph{Euler class} is 
\[e(E)=\Pf(F_A)\period\]

%-------------------------------------------------------------------%
%  Other classes                                                    %
%-------------------------------------------------------------------%

\subsubsection{Other classes}

If $g(X)=a_0+a_1X+a_2X^2+\cdots$ is a power series, then $\det(g(X))$ is invariant. For example,
\begin{itemize}
\item we get the total Chern class from
\[g=1+\frac{z}{2\pi i}\]
\item we get the L-genus from
\[g=\frac{z}{\tanh(z)}\]
\item we get the Todd genus from
\[g=\frac{z^2}{1-\exp(-z^2)}\period\]
\end{itemize}

%-------------------------------------------------------------------%
%-------------------------------------------------------------------%
%  Axioms                                                           %
%-------------------------------------------------------------------%
%-------------------------------------------------------------------%

\subsection{Axioms}

There are a set of axioms that Chern classes satisfy. 
Moreover, these axioms uniquely determine the Chern classes. 
See, for example, \cite[\S 4]{MilnorStasheff} for a discussion of this perspective. 
The axioms are
\begin{enumerate}
\item $c_0(E)=1$, $c_i(E)=0$ for $i>\rank(E)$
\subitem $c_i=\Tr(\wedge^iF_A)$ and $\wedge^i=0$ for $i\geq\rank(E)+1$.
\item Naturality with pullbacks
\item Whitney sum, $c(E\oplus F)=c(E)\cup c(F)$
\item Normalization $c_1(\cal{O}(1))=-1$ on $\CCP^1$.
\end{enumerate}
One can check that the Chern classes, as we have defined them above, 
satisfy these axioms, see \cite[Appendix C]{MilnorStasheff}. 
Thus, the Chern--Weil definition gives the same Chern classes as other definitions.
\begin{remark}
Although, a priori, $c_k$ has real coefficients $\HdR^{2k}(M;\RR)$, the normalization shows it is actually in the image of the map
\[\H^{2k}(M;\ZZ)\to \H^{2k}(M;\RR)\period\]
\end{remark}

%-------------------------------------------------------------------%
%-------------------------------------------------------------------%
%  An application                                                   %
%-------------------------------------------------------------------%
%-------------------------------------------------------------------%

\subsection{An application}

Here's an application of Chern--Weil theory to something harder to see with other definitions of characteristic classes.

\begin{lemma}
	Let $E\to M$ be a complex vector bundle that admits a reduction of structure group to locally constant transition functions (i.e., $E$ is a local system with group $\CC^n)$, then the Chern class $c_k(E)\in \H^{2k}(M;\ZZ)$ is torsion.
\end{lemma}

\begin{proof}
	The vector bundle $E$ admits a flat connection 
	\[A\mapsto gAg^{-1}+g^{-1}\d g\]
	so we can take $ \d+A $  with $A=0$, and this is preserved by changing trivializations.
\end{proof}

\newpage
%!TEX root = ../diffcoh.tex

%-------------------------------------------------------------------%
%-------------------------------------------------------------------%
%  Equivariant de Rham cohomology                                   %
%-------------------------------------------------------------------%
%-------------------------------------------------------------------%

\section{Equivariant de Rham cohomology}\label{EquivariantdeRhamCohomology}
\textit{by Greg Parker}

%-------------------------------------------------------------------%
%-------------------------------------------------------------------%
%  Motivation                                                       %
%-------------------------------------------------------------------%
%-------------------------------------------------------------------%

\subsection{Motivation}
Let $G$ be a Lie group and $M$ be a smooth manifold with a $G$ action. We want a cohomology theory that takes into account the $G$-action. If the action is free, then we can take
\[\HG^\bullet(M)\colonequals \H^\bullet(M/G)\period\]
If the action is not free, we take the homotopy quotient $\EG\times_G M$ and set the \emph{equivariant cohomology of $M$} to be 
\[\HG^\bullet(M)\colonequals \H^\bullet(\EG\times_G M)\period\]
Here $\EG\to \BG$ is the universal bundle, so that $\EG$ is a contractible space with a free $G$-action. 
\begin{quest}
How should one define equivariant cohomology using differential forms?
\end{quest}

\noindent To answer this question, we will roughly follow \cite[Chapter 1-4]{MR1689252}. 
The reader is encouraged to read \cite{MR1689252} for more details and applications. 

As motivation, again consider a free action. That is, take $P\to X$ to be a principal $G$-bundle. 
We want to distinguish forms in $\Omegabullet(P)$ that pullback from $X=P/G$. 
Let $\g $ be the Lie algebra of $G$. 
%The linearized action gives $X_\xi\in\Gamma(P,TP)$ for $\xi\in\g$, which together span the vertical subbundle $\ker(\pi_*:TP\to TM)$. 

For $\alpha\in\Omegabullet(P)$, we can locally write
\[\alpha=\sum_I\alpha_I\dup x_{i_1}\wedge\cdots\wedge\dup x_{i_N}\period\]
The form $\alpha$ is pulled back from $M$ if, for all $i$ the following conditions hold:
\begin{enumerate}
	\item The form $\dup x_i$ is vertical: $i_\xi\alpha=0$ for all $\xi\in\g$. 

	\item $\alpha$ does not depend on vertical coordinates: $i_\xi\dup \alpha=0$ for all $\xi\in\g$. 
\end{enumerate}
Forms satisfying these two conditions are called \emph{basic}. 
Let $\Omegabullet(P)_{\basic}$ denote the subcomplex of basic forms. 
Then, we have 
\[
	\H^\bullet(\Omega(P)_{\basic})\isomorphic\HdR^\bullet(X)\isomorphic\HdR^\bullet(P/G)=\HG^\bullet(P)\period
\]

%-------------------------------------------------------------------%
%-------------------------------------------------------------------%
%  G*-Algebras                                                      %
%-------------------------------------------------------------------%
%-------------------------------------------------------------------%

\subsection{\texorpdfstring{$G^*$}{G*}-Algebras}

Given an element $\xi\in \g$, there are multiple maps on $\Omegabullet(M)$:

\begin{itemize}

\item a degree $-1$ map by contraction, $\xi\mapsto i_\xi$ and 

\item a degree $0$ map by Lie derivative, $\xi\mapsto L_\xi$.

\end{itemize}

We can package these actions of $\g $, together with the differential $d$, on $\Omegabullet(M)$ as a representation of a certain Lie superalgebra $\gtilde$. 
Take 
\[\gtilde\colonequals \g_{-1}\oplus\g_0\oplus\RR\]
where, for each element $\xi\in\g$, we have corresponding elements of $\g_{-1}$ and $\g_0$ that we denote by their action on $\Omegabullet(M)$. That is, by $i_\xi$ and $L_\xi$, respectively. 
The generator of $\RR$ is denoted $d$. 
The bracket of the Lie superalgebra $\gtilde$ is defined by
\begin{align*}
	[i_\xi,i_\eta] &= 0\\
	[L_\xi,i_\eta] &= i_{[\xi,\eta]}\\
	[L_\xi,L_\eta] &= L_{[\xi,\eta]}\\
	[d,i_\xi] &= L_\xi\\
	[d,L_\xi] &= 0\\
	[d,d] &= 2d^2 = 0
\end{align*}
for all $\xi,\eta\in\g$.

The following is \cite[Definition 2.3.1]{MR1689252}. 
\begin{definition}
	A \emph{$G^*$-algebra} is a graded algebra $A$ with an action $G\to\Aut(A)$ of $G$ and an action $\gtilde\to\End(A)$ of $\gtilde$, so that the following hold:
	\begin{enumerate}
		\item $ \left. \frac{\d}{\d t} \right|_{t=0} \exp(t_\xi)=L_\xi$. 
		\item $gL_\xi g^{-1}=L_{\Ad_g\xi}$ and $gi_\xi g^{-1}=i_{\Ad_g\xi}$.

		\item $gd=dg$.
	\end{enumerate}
\end{definition}

Note that the tensor product of two $G^*$-algebras is again a $G^*$-algebra. %prove?

\begin{example}
	The complex $\Omegabullet(M)$ is a $G^*$-algebra with multiplication by the wedge product.
\end{example}

\noindent Considering a $G^*$-algebra $A$ with its differential from the action of $d\in\gtilde$, 
we define ${\H^\bullet(A)\colonequals \H_\bullet(A,d)}$. 

\begin{definition}
Let $A$ be a $G^*$-algebra. 
A \emph{basic form} in $A$ is an element $\alpha\in A$ so that 
\[i_\xi\alpha=L_\xi=0\]
for all $\xi\in\g $. 
\end{definition}

We will need to add an assumption on our $G^*$-algebra, referred to as \emph{Condition C} in \cite[\S 2.3.4]{MR1689252}. 
Condition C will ensure the existence of a certain $G$-invariant subspace that acts like the vertical subbundle (\Cref{subsec-principal}) in the locally free case, see \cite[Definition 2.3.3]{MR1689252}. 

\begin{definition}
	Let $\xi_1,\dots,\xi_k$ be a basis for $\g $. 
	A $G^*$-algebra $A$ is \emph{satisfies Condition C} if there exists elements $\theta^1,\dots,\theta^k\in A$ of degree 1 so that for all $i,j=1,\dots, k$,
	\[\iota_{\xi_i}\theta^{j}=\delta_{ij}\]
	and the subspace spanned by $\{\theta^i\}$ is invariant under $G$.
\end{definition}

In particular, if the action of $G$ on $M$ is free, then $\Omegabullet(M)$ is a $G^*$-algebra satisfying Condition C.

We say a $G^*$-algebra $A$ is \emph{acyclic} if the chain complex $(A,d)$ is.

\begin{definition}\label{def-Conditionc}
Let $M$ be a manifold with $G$ action and 
let $E$ be a $G^*$-algebra that is acyclic and satisfies Condition C. 
Define the equivariant de Rham cohomology by
\[
	\HGdR^\bullet(M)\colonequals \H^\bullet((\Omega(M)\otimes E)_{\basic}) \period
\]
\end{definition}

The following is \cite[Theorem 2.5.1]{MR1689252}. 
In particular, by \cite[Prop. 2.5.4]{MR1689252}, such $G^*$-algebras $E$ as in \Cref{def-Conditionc} exist in the context we care about. 

\begin{theorem}[(equivariant de Rham)]
	There is an isomorphism 
	\[
		\HGdR^*(M)\isomorphic\HG^*(M) \period
	\]
\end{theorem}

We discuss the idea of the proof here. For a full proof, see \cite[\S 2.5]{MR1689252}. 
\begin{proof}[Proof Idea]
Approximate $\EG$ with a sequence of finite-dimensional manifolds $E_k$ and take \[E = \lim_k \Omega(E_k)\period\]
By the free case,
\[
	\H^*(M\times E_k/G)=\H^*(\Omega(M\times E_k)_{\basic})
\] for $* \ll k$. 
To finish, one shows that
\[
	\Omega(M\times E_k)_{\basic}=\Omega(M)\otimes\Omega(E_k)_{\basic}
\]
in the limit.
\end{proof}

\begin{remark}
By \cite[\S 4.4]{MR1689252}, the definition of $\HGdR^*$ is independent of $E$ satisfying the assumptions (acyclic and Condition C).
\end{remark}

%-------------------------------------------------------------------%
%-------------------------------------------------------------------%
%  Cartan model                                                     %
%-------------------------------------------------------------------%
%-------------------------------------------------------------------%

\subsection{Cartan model}

Now we can look for a specific $E$ that gives a nice algebraic structure, so it might be more computable.

For a vector space $V$, the \emph{Koszul algebra} is $(\exterior^\bullet(V)\otimes \Sym^\bullet(V),\d)$ where
\begin{equation*}
	\d(\alpha\otimes 1)=1\otimes \alpha \andeq \d(1\otimes\alpha)=0
\end{equation*}
extended as a derivation.  \index[terminology]{Koszul complex}The \emph{Weil Algebra} is the Koszul algebra of $\gdual$, 
\[W=\exterior^\bullet(\gdual)\otimes\Sym^\bullet(\gdual)  \index[terminology]{Weil Algebra}\]
as a $G^*$ algebra: For a basis (as an algebra) $\theta^i, z^j$ we have
\begin{align*}
	i_a\theta_b &= \delta_{ab}\\
	L_a\theta_b &= -[\theta_a,\theta_b]=-c_{ab}^k\theta_k\\
	L_az_b &= -c_{ab}^k z_k \\
	i_az_b &= -c_{ab}^k\theta_k \period
\end{align*}

The following is \cite[Theorem 3.2.1]{MR1689252}.
\begin{proposition}
	The Weil algebra $W$ is acyclic and satisfies Condition C.
\end{proposition}

\begin{proof}[Proof of Acyclicity]
	Define a chain homotopy $ Q $ from $\id{}$ to $0$ by setting
	\begin{equation*}
		Q(\alpha\otimes 1) \colonequals 0 \andeq Q(1\otimes \alpha) \colonequals \alpha\otimes 1 \period \qedhere
	\end{equation*}
\end{proof}

\noindent In particular, we can use $W$ as a model for $E$.

The $G^*$-algebra $W$ has a rather nice subalgebra of basic forms. 
By \cite[Theorem 3.2.2]{MR1689252}, the basic cohomology ring of the Weil algebra $W$ is $\Sym^\bullet(\gdual)^G$. 
Thus 
\[\H_*((W\otimes \Omega^*(M))_{\basic},\dup \vert_{\basic})\]
 calculates $\HG^*(M)$. 
One can use this description of the equivariant de Rham cohomology of the Weil algebra to obtain a description, 
called the \emph{Cartan model}, 
of the equivariant de Rham cohomology of any $G^*$-algebra. 

\begin{theorem}[(Cartan model)]  \index[terminology]{Cartan model}\label{thm-Cartan}
	For a $G^*$-algebra $A$, there is an isomorphism (the Mathai--Quillen isomorphism) 
	\begin{align*}
		\varphi\colon (W\otimes A)_{\basic} &\isomorphism (\Sym^\bullet(\gdual)\otimes A)^G \\ 
		\shortintertext{sending}
		\dup \vert_{\basic} &\mapsto\dup_G=1\otimes \dup_A-\mu^a\otimes i_a \period
	\end{align*}
\end{theorem}

\noindent In particular, $\HG^*(M)$ can be computed from $(\Sym^\bullet(\gdual)\otimes\Omega^*(M))^G,\dup_G)$.

\newpage
%!TEX root = ../diffcoh.tex

%-------------------------------------------------------------------%
%-------------------------------------------------------------------%
%  Classifying spaces for G-bundles                                 %
%-------------------------------------------------------------------%
%-------------------------------------------------------------------%

\section{Classifying spaces for \texorpdfstring{$G$}{G}-bundles}\label{sec:classifying_spaces_for_G-bundles}
\textit{by Peter Haine}

Let $ G $ be a Lie group.
The purpose of this chapter is to study two differential cohomology variants of the classifying space $ \BG $ of $ G $.
We first introduced these objects in \Cref{ex:BunG,ex:BunGnabla}: they are the sheaves $ \BunG $ and $ \BunGnabla $ that send a manifold $ M $ to the groupoid of principal $ G $-bundles on $ M $ and principal $ G $-bundles on $ M $ equipped with a connection, respectively.
Our main goal is to show that $ \BunG $ and $ \BunGnabla $ are differential \textit{refinements} of the space $ \BG $ in the sense that there are natural equivalences
\begin{equation*}	
	\Gammalowersharp(\BunGnabla) \equivalence \Gammalowersharp(\BunG) \equivalence \BG
\end{equation*}
in the \category $ \Spc $ (\Cref{cor:formula_for_GammalowersharpBunG,cor:formula_for_Gammalowersharp_BunGnabla}).
We also show that $ \BunG $ and $ \BunGnabla $ refine the classifying space of the underlying \textit{discrete} group $ \Gdisc $ associated to $ G $ in the sense that there are natural equivalences
\begin{equation*}	
	\Gammalowerstar(\BunGnabla) \equivalence \Gammalowerstar(\BunG) \equivalence \BGdisc
\end{equation*}
in the \category $ \Spc $ (\Cref{lem:BGdisc}).

In order to prove that $ \BunG $ and $ \BunGnabla $ refine the classifying spaces $ \BG $ and $ \BGdisc $, we use explicit \textit{presentations} of $ \BunG $ and $ \BunGnabla $.
The sheaf $ \BunG $ can be realized as the geometric realization in $ \Sh(\Mfld) $ of the bar construction
\begin{equation*}
	\begin{tikzcd}[sep=1.5em]
	    \cdots \arrow[r, shift left=0.75ex] \arrow[r, shift right=0.75ex] \arrow[r, shift right=2.25ex] \arrow[r, shift left=2.25ex] & G \cross G \arrow[l] \arrow[l, shift left=1.5ex] \arrow[l, shift right=1.5ex] \arrow[r] \arrow[r, shift left=1.5ex] \arrow[r, shift right=1.5ex] & G \arrow[l, shift left=0.75ex] \arrow[l, shift right=0.75ex] \arrow[r, shift left=0.75ex] \arrow[r, shift right=0.75ex] & \ast  \arrow[l] 
	\end{tikzcd}
\end{equation*}
on the Lie group $ G $. 
That is, $ \BunG $ is the quotient $ \ast \modmod G $ in $ \Sh(\Mfld) $ of the point by $ G $ (\Cref{prop:Cech_nerve_trivG}).
Similarly, the sheaf $ \BunGnabla $ can be realized as the geometric realization in $ \Sh(\Mfld) $ of \textit{action groupoid}
\begin{equation*}
	\begin{tikzcd}[sep=1.5em]
	    \cdots \arrow[r, shift left=0.75ex] \arrow[r, shift right=0.75ex] \arrow[r, shift right=2.25ex] \arrow[r, shift left=2.25ex] & G \cross G \cross \Omega^1(-;\g) \arrow[l] \arrow[l, shift left=1.5ex] \arrow[l, shift right=1.5ex] \arrow[r] \arrow[r, shift left=1.5ex] \arrow[r, shift right=1.5ex] & G \cross \Omega^1(-;\g) \arrow[l, shift left=0.75ex] \arrow[l, shift right=0.75ex] \arrow[r, shift left=0.75ex] \arrow[r, shift right=0.75ex] & \Omega^1(-;\g) \period \arrow[l] 
	\end{tikzcd}
\end{equation*}
of the \textit{adjoint action} of $ G $ on the sheaf $ \Omega^1(-;\g) $ (see \Cref{def:adjoint_action,rem:adjoint_action_explicitly}). 
That is, $ \BunGnabla $ is the quotient $ \Omega^1(-;\g) \modmod G $ (\Cref{prop:Cech_nerve_trivGnabla}).
We then deduce that $ \BunG $ and $ \BunGnabla $ refine $ \BG $ and $ \BGdisc $ by using these presentations to compute the spaces $ \Gammalowersharp(\BunG) $ and $ \Gammalowersharp(\BunGnabla) $ as well as $ \Gammalowerstar(\BunG) $ and $ \Gammalowerstar(\BunGnabla) $.

Many texts take this colimit as the \textit{definition} of $ \BunG $, and then try to argue that this colimit indeed classifies principal $ G $-bundles.
This is difficult for a few reasons, the main one being that colimits in \categories of sheaves are not very explicit.
Instead, we work in the other direction and use general categorical techniques to show that $ \BunG $ admits such a presentation.
After these categorical reductions, the result boils down to the straightforward claim that the group of automorphisms of the trivial $ G $-bundle $ M \cross G \to M $ is naturally isomorphic to $ \Cinf(M,G) $; see \Cref{obs:gauge_group_trivial}.

We also want to highlight two advantages of the sheaf $ \BunG $ over the classifying space $ \BG $:
\begin{enumerate}
	\item If $ G $ is not discrete, then $ \BG $ generally has homotopy in all degrees; on the other hand, $ \BunG $ is a sheaf of \textit{groupoids}.

	\item Homotopy classes of maps $ \fromto{M}{\BG} $ are in bijection with \textit{isomorphism classes} of principal $ G $-bundles on $ M $; on the other hand, the space $ \Map_{\Sh(\Mfld)}(M,\BunG) $ is by definition the groupoid of principal $ G $-bundles on $ M $.
\end{enumerate}
Thus, in many ways the sheaf $ \BunG $ is more simple to work with than the classifying space $ \BG $.
The price we pay for this simplicity is that we have to work with a sheaf, rather than a single space.

Looking forward, there are also a number of applications of the work from this chapter:
\begin{enumerate}
	\item \Cref{DifferentialCharacteristicClasses} uses the sheaves $ \BunG $ and $ \BunGnabla $ to lift characteristic classes to classes living in differential cohomology.

	\item \Cref{WorkofFreedHopkins} presents work of Freed--Hopkins that uses Chern--Weil theory to compute the de Rham complex of the sheaf $ \BunGnabla $.

	\item \Cref{LiftsofChernClasses} computes differential cohomology groups of $ \BbulletGL{m}{\CC} $ and $ \BbulletGL{m}{\RR} $. 
\end{enumerate}

\Cref{subsec:groupoids_effective_epis} sets up the general categorical framework that we need to provide presentations for $ \BunG $ and $ \BunGnabla $.
In particular, we discuss monoid and group objects in \categories, actions of monoid objects, and effective epimorphisms.
The reader comfortable with these ideas can safely skip this section.
In \cref{sec:structure_of_BunG}, we use these tools to provide a presentation for $ \BunG $.
In \cref{sec:structure_of_BunGnabla}, we provide a presentation for $ \BunGnabla $.

%-------------------------------------------------------------------%
%-------------------------------------------------------------------%
%  Groupoid objects and effective epimorphisms                      %
%-------------------------------------------------------------------%
%-------------------------------------------------------------------%

\subsection{Groupoid objects and effective epimorphisms}\label{subsec:groupoids_effective_epis}

The purpose of this section is to develop a dictionary between sheaves of groupoids on $ \Mfld $ and simplicial diagrams in $ \Sh(\Mfld;\Set) $.
The point of developing this dictionary is that it will give us our desired presentations of the sheaves of groupoids $ \BunG $ and $ \BunGnabla $.

%-------------------------------------------------------------------%
%  Motivation: monoid objects & group objects in ordinary categories%
%-------------------------------------------------------------------%

\subsubsection{Motivation: monoid objects \& group objects in ordinary categories}\label{subsec:monoids_and_groups}

In order to introduce monoid, groupoid, and group objects in \acategory, we start with some motivation.
Recall that in an ordinary category $ \X $ with finite products, a \textit{monoid object} in $ \X $ is the data of an object $ M \in \X $ along with unit and multiplication maps
\begin{equation*}
	u \colon \fromto{\ast}{M} \andeq m \colon \fromto{M \cross M}{M}
\end{equation*}
such that the unitality and associativity diagrams
\begin{equation*}
	\begin{tikzcd}[sep=3em]
		M \isomorphic M \cross \ast \arrow[r, "\id{M} \cross u"] \arrow[dr, equals] & M \cross M \arrow[d, "m" description] & \ast \cross M \isomorphic M \arrow[l, "u \cross \id{M}"'] \arrow[dl, equals] \\
		& M & 
	\end{tikzcd}
\end{equation*}
and 
\begin{equation*}
	\begin{tikzcd}[sep=3em]
		M \cross M \cross M \arrow[r, "m \cross \id{M}"] \arrow[d, "\id{M} \cross m"'] & M \cross M \arrow[d, "m"] \\
		M \cross M \arrow[r, "m"'] & M 
	\end{tikzcd}
\end{equation*}
commute.
To explain the definition of a monoid in \acategory, observe that this data can be neatly packaged as a truncated simplicial diagram 
\begin{equation*}
    \begin{tikzcd}[sep=9em]
        M \cross M \cross M \arrow[r, shift left=1.5ex, "m \cross \idnosub"{description, near end}] \arrow[r, shift right=1.5ex, "\idnosub \cross m"{description, near end}] \arrow[r, shift right=4.5ex, "\pr_{1,2}"'] \arrow[r, shift left=4.5ex, "\pr_{2,3}"] &%
        M \cross M \arrow[l, "\idnosub \cross u \cross \idnosub"{description, near end}] \arrow[l, shift left=3ex, "u \cross \idnosub \cross \idnosub"{description, near end}] \arrow[l, shift right=3ex, "\idnosub \cross \idnosub \cross u"{description, near end}] \arrow[r, "m"{description, near end}] \arrow[r, shift left=3ex, "\pr_1"] \arrow[r, shift right=3ex, "\pr_2"'] &%
        M \arrow[l, shift left=1.5ex, "\idnosub \cross u"{description, near end}] \arrow[l, shift right=1.5ex, "u \cross \idnosub"{description, near end}] \arrow[r, shift left=1.5ex] \arrow[r, shift right=1.5ex] &%
        \ast \period \arrow[l, "u" description] 
      \end{tikzcd}
\end{equation*} 
The simplicial identities encode the associativity and unitality of composition, as well as some information that is redundant in ordinary category theory.

Readers familiar with the bar construction may recognize this truncated simplicial diagram as the first few stages of the bar construction. 
To make a tighter connection between simplicial objects and monoids, let us recall the precise description of
monoids via their bar constructions.

\begin{notation}
	Write $ \Catalg_1 $ for the $ (1,1) $-category with objects $ 1 $-categories and morphisms functors.%
	\footnote{We've used the somewhat strange notation $ \Catalg_1 $ here to highlight that we are working with the ordinary category of categories, where we do not even keep track of natural isomorphisms.
	Thus this is the `algebraic theory' of categories where we do not identify equivalent categories.}
\end{notation}

\begin{notation}
	We write
	\begin{equation*}
		\Nerve \colon \incto{\Catalg_{1}}{\Fun(\Deltaop,\Set)}
	\end{equation*}
	for the fully faithful \emph{nerve} functor from $ 1 $-categories to simplicial sets.
	Write $ \Mon $ for the $ 1 $-category of monoids.
	We regard $ \Mon $ as a full subcategory of $ \Catalg_1 $ via the functor sending a monoid $ M $ to the one-object category with endomorphisms the monoid $ M $.
	The composite fully faithful functor
	\begin{equation*}
		\begin{tikzcd}
			\Mon \arrow[r, hooked] & \Catalg_1 \arrow[r, hooked, "\Nerve"] & \Fun(\Deltaop,\Set)
		\end{tikzcd}
	\end{equation*}
	is called the \emph{bar construction}.
	We denote this composite by
	\begin{equation*}
		\Bar \colon \incto{\Mon}{\Fun(\Deltaop,\Set)} \period
	\end{equation*}
\end{notation}

\begin{lemma}[{\cite[Theorem 1.1.52]{MR4259746}}]\label{lem:Nerve_essential_image}
	The essential image of the nerve
	\begin{equation*}
		\Nerve \colon \incto{\Catalg_{1}}{\Fun(\Deltaop,\Set)}
	\end{equation*}
	consists of those simplicial sets $ X \colon \Deltaop \to \Set $ satisfying the \emph{Segal condition}:
	\begin{enumerate}[label=\stlabel{lem:Nerve_essential_image}, ref=\arabic*]
	    \item\label{lem:Nerve_essential_image.1} For each $ n > 0 $ and $ t \in [n] $, the square 
	    \begin{equation*}
	        \begin{tikzcd}[sep=2em]
	            X([n]) \arrow[r] \arrow[d] & X(\{t < \cdots < n\}) \arrow[d] \\
	            X(\{0 < \cdots < t\}) \arrow[r] & X(\{t\})
	        \end{tikzcd}
	    \end{equation*}
	    is a pullback square.
	\end{enumerate} 

	Consequently, the essential image of the bar construction 
	\begin{equation*}
		\Bar \colon \incto{\Mon}{\Fun(\Deltaop,\Set)} 
	\end{equation*}
	consists of those simplicial sets $ X \colon \Deltaop \to \Set $ satisfying \enumref{lem:Nerve_essential_image}{1} and:
	\begin{enumerate}[label=\stlabel{lem:Nerve_essential_image}, ref=\arabic*]
		\setcounter{enumi}{1}

	    \item\label{lem:Nerve_essential_image.2} We have $ X_0 \equivalent \ast $.
	\end{enumerate}
\end{lemma}

\begin{remark}
	By induction, the Segal condition \enumref{lem:Nerve_essential_image}{1} is equivalent to the condition that for each $ n > 0 $, the natural map
	\begin{equation*}	
		X([n]) \to X(\{0 < 1\}) \crosslimits_{X(\{1\})} X(\{1 < 2\}) \crosslimits_{X(\{2\})} \cdots \crosslimits_{X(\{n-1\})} X(\{n-1 < n\})
	\end{equation*}
	is an isomorphism.
	Said differently, let $ \Spine^n \subset \Delta^n $ be the simplicial subset given by the union of edges between successive vertices
	\begin{equation*}
		\Spine^n \colonequals \Union_{i=0}^{n-1} \Delta^{\{i<i+1\}} \period
	\end{equation*}
	\Cref{lem:Nerve_essential_image} says that a simplicial set $ X $ is in the essential image of the nerve if and only if every map $ \fromto{\Spine^n}{X} $ admits a \textit{unique} extension to an $ n $-simplex $ \fromto{\Delta^n}{X} $.
\end{remark}

\begin{remark}
	Note that the spine $ \Spine^2 $ of $ \Delta^2 $ is the horn $ \Lambda_1^2 $.
	Thus, the Segal condition \enumref{lem:Nerve_essential_image}{1} for $ n = 2 $ is equivalent to the statement that each horn $ \fromto{\Lambda_1^2}{X} $ admits a \textit{unique} extension to a $ 2 $-simplex $ \fromto{\Delta^2}{X} $.
	More generally, the Segal condition \enumref{lem:Nerve_essential_image}{1} is equivalent to the condition that for each pair $ n > 0 $ and $ 0 < i < n $, every inner horn $ \fromto{\Lambda_i^n}{X} $ admits a \textit{unique} extension to an $ n $-simplex $ \fromto{\Delta^n}{X} $.
\end{remark}

\begin{observation}[(multiplication and unit via simplicial maps)]
	If $ X \colon \Deltaop \to \Set $ is a simplicial set satisfying \enumref{lem:Nerve_essential_image}{1} and \enumref{lem:Nerve_essential_image}{2}, the face map 
	\begin{equation*}
		d_1 \colon X_1 \times X_1 \equivalent X_2 \to X_1
	\end{equation*}
	provides a multiplication on $ X_1 $ with unit given by the degeneracy map
	\begin{equation*}
		s_0 \colon \ast \equivalent X_0 \to X_1 \period
	\end{equation*}
	The simplicial identities involving $ X_0 $, $ X_1 $, $ X_2 $, and $ X_3 $ encode the associativity and unitality of the multiplication.
	The simplicial identities involving the $ (n+1) $-simplices for $ n \geq 3 $ encode higher-order associativity conditions for multiplying $ n $ elements.
	While these higher coherences are automatic for monoids in the category of sets, when working in higher category theory it is important to encode these coherences.
\end{observation}	

Since groups form a full subcategory of the category of monoids, the bar construction also identifies the category of groups with a full subcategory of the category of simplicial sets.
For this it is better to use an alternative characterization of the existence of inverses: a monoid $ M $ is a group if and only if the \textit{shear maps}
\begin{align*}
    M \times M &\to M \times M \andeq M \times M \to M \times M\\ 
    (x,y) &\mapsto (x,xy) \phantom{\andeq} (x,y) \mapsto (xy,y)
\end{align*}
are bijections.
Translating this into simplicial sets one sees that the category of groups is equivalent to the full subcategory of $ \Fun(\Deltaop,\Set) $ spanned by the simplicial sets $ X $ satisfying \enumref{lem:Nerve_essential_image}{1}, \enumref{lem:Nerve_essential_image}{2}, and:
\begin{enumerate}
   \setcounter{enumi}{2}

    \item\label{cor:Bar_fully_faithful.3} The induced squares
    \begin{equation*}
        \begin{tikzcd}[sep=2em]
            X(\{0 < 1 < 2\}) \arrow[r, "d_1"] \arrow[d, "d_2"'] & X(\{0 < 2\}) \arrow[d] \\
            X(\{0 < 1\}) \arrow[r] & X(\{0\})
        \end{tikzcd} 
        \andeq
         \begin{tikzcd}[sep=2em]
            X(\{0 < 1 < 2\}) \arrow[r, "d_1"] \arrow[d, "d_0"'] & X(\{0 < 2\}) \arrow[d] \\
            X(\{1 < 2\}) \arrow[r] & X(\{0\})
        \end{tikzcd} 
    \end{equation*}
    are pullback squares.
\end{enumerate}
We emphasize that condition \eqref{cor:Bar_fully_faithful.3} is not implied by the Segal condition \enumref{lem:Nerve_essential_image}{1}.

\begin{notation}
	Write $ \Gpdalg_1 \subset \Catalg_1 $ for the full subcategory spanned by the groupoids.
\end{notation}

\begin{corollary}\label{cor:Nerve_essential_image_Gpd}
	The essential image of the nerve functor restricted to groupoids
	\begin{equation*}
		\Nerve \colon \incto{\Gpdalg_{1}}{\Fun(\Deltaop,\Set)}
	\end{equation*}
	consists of those simplicial sets $ X \colon \Deltaop \to \Set $ satisfying:
	\begin{enumerate}[label=\stlabel{cor:Nerve_essential_image_Gpd}, ref=\arabic*]
	    \item\label{cor:Nerve_essential_image_Gpd.1} For each object $ S \in \Deltaop $ and partition $ S = T \cup T' $ such that $ T \cap T' = \{t\} $ consists of a single element, the induced square
        \begin{equation*}
            \begin{tikzcd}[sep=2em]
                X(S) \arrow[r] \arrow[d] & X(T') \arrow[d] \\
                X(T) \arrow[r] & X(\{t\})
            \end{tikzcd}
        \end{equation*}
        is a pullback square.
	\end{enumerate} 

	Consequently, the essential image of the bar construction restricted to groups
	\begin{equation*}
		\Bar \colon \incto{\Grp}{\Fun(\Deltaop,\Set)} 
	\end{equation*}
	consists of those simplicial sets $ X \colon \Deltaop \to \Set $ satisfying \enumref{cor:Nerve_essential_image_Gpd}{1} and:
	\begin{enumerate}[label=\stlabel{cor:Nerve_essential_image_Gpd}, ref=\arabic*]
		\setcounter{enumi}{1}

	    \item\label{cor:Nerve_essential_image_Gpd.2} We have $ X_0 \equivalent \ast $.
	\end{enumerate}
\end{corollary}

%-------------------------------------------------------------------%
%  Monoid & group objects in ∞-categories                           %
%-------------------------------------------------------------------%

\subsubsection{Monoid \& group objects in \texorpdfstring{$\infty$}{∞}-categories}

The conclusions of \Cref{lem:Nerve_essential_image,cor:Nerve_essential_image_Gpd} give the correct definition of a monoid object in an arbitrary \category.

\begin{definition}[(monoid object)]\label{def:monoid_object}
	Let $ \X $ be \acategory with finite products.
	An \textit{associative monoid} or \textit{$ \Eone $-monoid} in $ \X $ is a simplicial object $ M \colon \fromto{\Deltaop}{\X} $ such that 
	\begin{enumerate}[label=\stlabel{def:monoid_object}, ref=\arabic*]
	    \item\label{def:monoid_object.1} The object $ M_{0} $ is a terminal object of $ \X $.

	    \item\label{def:monoid_object.2} \textit{Segal condition:} For each $ n > 0 $ and $ t \in [n] $, the square 
	    \begin{equation*}
	        \begin{tikzcd}[sep=2em]
	            M([n]) \arrow[r] \arrow[d] & M(\{t < \cdots < n\}) \arrow[d] \\
	            M(\{0 < \cdots < t\}) \arrow[r] & M(\{t\})
	        \end{tikzcd}
	    \end{equation*}
	    is a pullback square in $ \X $.
	\end{enumerate}
	In this case, we call $ M_1 \in \X $ the \textit{underlying object} of $ M $.
    We often identify a monoid object by its underlying object.
    We write
    \begin{equation*}
    	\Mon(\X) \subset \Fun(\Deltaop,\X)
    \end{equation*}
    for the full subcategory spanned by the monoid objects.
\end{definition}

\begin{definition}[{(groupoid object \cites[Definitions \HTTthmlink{6.1.2.7} \& \HTTthmlink{7.2.2.1}]{HTT})}]\label{def:group}
    Let $ \X $ be \acategory with finite limits.
    A \textit{groupoid object} in $ \X $ is a simplicial object $ G \colon \Deltaop \to \X $ such that for each object $ S \in \Deltaop $ and partition $ S = T \cup T' $ such that $ T \cap T' = \{t\} $ consists of a single element, the induced square
    \begin{equation*}
        \begin{tikzcd}[sep=2em]
            G(S) \arrow[r] \arrow[d] & G(T') \arrow[d] \\
            G(T) \arrow[r] & G(\{t\})
        \end{tikzcd}
    \end{equation*}
    is a pullback square in $ \X $.
    We write
    \begin{equation*}
    	\Gpd(\X) \subset \Fun(\Deltaop,\X)
    \end{equation*}
    for the full subcategory spanned by the groupoid objects.
    
    A groupoid object $ G $ is a \emph{group object} if $ G_0 \equivalent \ast $.
    We write
    \begin{equation*}
    	\Grp(\X) \colonequals \Mon(\X) \intersect \Gpd(\X) 
    \end{equation*}
    for the full subcategory of $ \Fun(\Deltaop,\X) $ spanned by the group objects.
\end{definition}

Another way of phrasing \Cref{def:monoid_object,def:group} is that in the approach we take, we take the \textit{definition} of a group object to be its delooping.
The reader should consult \cite[\S\S\HAsubseclink{4.1.2}, \HAsubseclink{5.2.6}, \& \HAsubseclink{6.1.2}]{HA} for a more complete treatment.

%-------------------------------------------------------------------%
%  Čech nerves & effective epimorphisms                             %
%-------------------------------------------------------------------%

\subsubsection{Čech nerves \& effective epimorphisms}

In this subsection, we go in a slightly different direction and introduce a class of maps called \textit{effective epimorphisms}.
These maps generalize surjections between sets.
They are important to us because specifying a groupoid object in $ \Sh(\Mfld) $ is equivalent to specifying an effective epimorphism (\Cref{thm:effective_epis_and_groupoids}).
In particular, we take this alternative perspective to give presentations of $ \BunG $ and $ \BunGnabla $.

We begin with some notation.

\begin{notation}
	Write $ \Deltaplus $ for the \emph{augmented simplex category}.
	That is $ \Deltaplus $ is the category of (possibly empty) linearly ordered finite sets.
	We write $ [-1] \in \Deltaplus $ for the empty linearly ordered set.
	The usual simplex category $ \DDelta $ is the full subcategory of $ \Deltaplus $ spanned by the \textit{nonempty} linearly ordered sets.

	Given an integer $ n \geq -1 $, write
	\begin{equation*}
		\DDelta_{+,\leq n} \subset \Deltaplus 
	\end{equation*}
	for the full subcategory containing $ [i] $ for $ -1 \leq i \leq n $ and closed under isomorphism.
	% We write
	% \begin{equation*}
	% 	\DDelta_{\leq n} \colonequals \DDelta_{+,\leq n} \intersect \DDelta \period
	% \end{equation*} 
\end{notation}

\begin{observation}
	Note that $ \DDelta_{+,\leq 0} $ is equivalent to the category with two objects $ [-1] $ and $ [0] $ and a single non-identity arrow $ \fromto{[-1]}{[0]} $.
\end{observation}

\begin{definition}[(Čech nerve)]\label{def:Cech_nerve}
    Let $ \X $ be \acategory with pullbacks, and let $ e \colon W \to X $ be a morphism in $ \X $.
    The \textit{Čech nerve} $ \Cech_{+}(e) $ of $ e $ is the augmented simplicial object in $ \X $ given by the right Kan extension of the functor 
    \begin{equation*}
    	\fromto{\DDelta_{+,\leq 0}^{\op}}{\X} 
    \end{equation*}
    that picks out the morphism $ e \colon \fromto{W}{X} $ along the inclusion $ \DDelta_{+,\leq 0}^{\op} \subset \Deltaplusop $.
    Concretely, $ \Cech_{+}(e)_n $ is the augmented simplical object
    \begin{equation*}
        \begin{tikzcd}[sep=1.5em]
            \cdots \arrow[r, shift left=0.75ex] \arrow[r, shift right=0.75ex] \arrow[r, shift right=2.25ex] \arrow[r, shift left=2.25ex] & \displaystyle W \crosslimits_X W \crosslimits_X W \arrow[l] \arrow[l, shift left=1.5ex] \arrow[l, shift right=1.5ex] \arrow[r] \arrow[r, shift left=1.5ex] \arrow[r, shift right=1.5ex] & \displaystyle W \crosslimits_X W \arrow[l, shift left=0.75ex] \arrow[l, shift right=0.75ex] \arrow[r, shift left=0.75ex] \arrow[r, shift right=0.75ex] & W \arrow[l] \arrow[r, "e"] & X 
          \end{tikzcd}
    \end{equation*} 
   	where $ \Cech_{+}(e)_n $ is the $ (n+1) $-fold fiber product of $ W $ over $ X $, each degeneracy map is a diagonal morphism, and each face map is a projection.
    
    We write $ \Cech(e) \colon \fromto{\Deltaop}{\X} $ for the restriction of $ \Cech_{+}(e) $ to $ \Deltaop \subset \Deltaplusop $.
    We also refer to $ \Cech(e) $ as the Čech nerve of $ e $.
\end{definition}

\begin{observation}
	Let $ \X $ be \acategory with pullbacks, and let $ e \colon W \to X $ be a morphism in $ \X $.
    The Čech nerve $ \Cech(e) $ is always a groupoid object of $ \X $.
\end{observation}

\begin{definition}[(effective epimorphism)]\label{def:effective_epi}
    Let $ \X $ be \acategory with pullbacks.
    A morphism $ e \colon W \to X $ in $ \X $ is an \textit{effective epimorphism} if the induced map $ |\Cech(e)| \to X $ is an equivalence.
    We often denote an effective epimorphism by the two-headed arrow `$ \twoheadrightarrow $'.
\end{definition}

\begin{example}
	A morphism $ e $ in the $ 1 $-category of sets is an effective epimorphism if and only if $ e $ is a surjection.
\end{example}

\begin{example}[\HTT{Corollary}{7.2.1.15}]\label{ex:effective_epi_Spc}
	A morphism $ e \colon W \to X $  in the \category $ \Spc $ is an effective epimorphism if and only if $ \uppi_0(e) \colon \fromto{\uppi_0(W)}{\uppi_0(X)} $ is a surjection.
\end{example}

The important fact that we need is that in $ \Sh(\Mfld) $ a map that is \textit{locally} a $ \uppi_0 $-surjection is an effective epimorphism:

\begin{lemma}\label{lem:ShMfld_criterion_for_effective_epi}
	Let $ f \colon \fromto{E}{E'} $ be a morphism in $ \Sh(\Mfld) $.
	If for each $ n \geq 0 $, the induced morphism
	\begin{equation*}
		\fromto{\uppi_0 E(\RR^n)}{\uppi_0 E'(\RR^n)}
	\end{equation*}
	on connected components is a surjection, then $ f $ is an effective epimorphism.
\end{lemma}

\begin{proof}
	Combine the equivalence
	\begin{equation*}
		\Sh(\Mfld) \equivalent \Sh(\Euc)
	\end{equation*}
	of \Cref{lem:hypersheavesonCart} with \Cref{ex:effective_epi_Spc} and \cite[\HTTthm{Remark}{6.5.1.15} \& \HTTthm{Proposition}{7.2.1.14}]{HTT}.
\end{proof}

\begin{warning}[(categorical epimorphisms vs. effective epimorphisms)]
	Recall that a morphism $ e \colon \fromto{W}{X} $ in a $ 1 $-category $ \X $ is called an \emph{epimorphism} if the square
	\begin{equation}\label{sq:epimorphism_def}
		\begin{tikzcd}
			W \arrow[r, "e"] \arrow[d, "e"'] & X \arrow[d, equals] & \\
			X \arrow[r, equals] & X
		\end{tikzcd}
	\end{equation}
	is a pushout square.
	In many $ 1 $-categories, the notions of an `epimorphism' and an `effective epimorphism' coincide.
	For example, a map of sets $ e \colon \fromto{W}{X} $ is an effective epimorphism if and only if $ e $ is an epimorphism if and only if $ e $ is a surjection.

	This is no longer the case in the setting of \categories. 
	For example, given a morphism $ e \colon \fromto{W}{X} $ in $ \Spc $, the square \eqref{sq:epimorphism_def} is a pushout if and only if $ e $ is \emph{acyclic} \cite[Theorem 2.1]{MR3987558}.
	That is, the reduced integral homology groups of all of the fibers of $ e $ vanish.
	Since $ \uppi_0 \colon \fromto{\Spc}{\Set} $ preserves colimits and epimorphisms in $ \Set $ are surjections, every acyclic map in $ \Spc $ is an effective epimorphism.
	However, acyclicity is a much stronger condition.
\end{warning}

The following results explain why the data of a groupoid object in the \category $ \Sh(\Mfld) $ is \textit{equivalent} to the data of an effective epimorphism.

\begin{notation}
	Let $ \X $ be \acategory with pullbacks.
	We write
	\begin{equation*}
		\Eff(\X) \subset \Fun([1],\X)
	\end{equation*}
	for the full subcategory of the arrow \category of $ \X $ spanned by those arrows $ \fromto{W}{X} $ that are effective epimorphisms.
	Write $ \Eff_{\ast}(\X) \subset \Eff(\X) $ for the full subcategory spanned by those effective epimorphisms $ \surjto{\ast}{X} $ with source the terminal object of $ \X $.
\end{notation}

\begin{theorem}[{\cite[\HTTpage{587}]{HTT}}]\label{thm:effective_epis_and_groupoids}
	Let $ \X $ be \atopos.
	The formation of the Čech nerve defines an equivalence of \categories
	\begin{equation*}
		\Cech \colon \equivto{\Eff(\X)}{\Gpd(\X)}
	\end{equation*}
	between effective epimorphisms in $ \X $ and groupoid objects in $ \X $.
	The inverse is given by sending a groupoid object $ G $ in $ \X $ to the induced effective epimorphism $ \surjto{G_0}{|G|} $.

	Moreover, the Čech nerve restricts to an equivalence
	\begin{equation*}
		\Cech \colon \equivto{\Eff_{\ast}(\X)}{\Grp(\X)}
	\end{equation*}
	between effective epimorphisms with source the terminal object of $ \X $ and group objects of $ \X $.
\end{theorem}

\begin{remark}
	Given a morphism $ x \colon \fromto{\ast}{X} $, note that the $ 1 $-simplices of the Čech nerve of $ x \colon \fromto{\ast}{X} $ are given by the loop object $ \Omega_x X $, i.e., the pullback
	\begin{equation*}
		\begin{tikzcd}[sep=2.5em]
			\Omega_x X \arrow[r] \arrow[d] \arrow[dr, phantom, "\lrcorner"{description, very near start}] & \ast \arrow[d, "x"] \\
			\ast \arrow[r, "x"'] & X \period
 		\end{tikzcd}
	\end{equation*}
	In particular, under the equivalence of \Cref{thm:effective_epis_and_groupoids}, an effective epimorphism $ x \colon \surjto{\ast}{X} $ corresponds to the loop object $ \Omega_x X $ equipped with a natural group structure.
\end{remark}

\begin{remark}
	The second part of \Cref{thm:effective_epis_and_groupoids} can be formulated in a slightly different manner.
	Write $ \X_{\ast} $ for the \category of pointed objects of $ \X $ and $ \X_{\ast}^{\cn} \subset \X_{\ast} $ for the full subcategory spanned by the \textit{connected} objects.
	Then \HA{Theorem}{5.2.6.15} asserts that the formation of the Čech nerve of the basepoint $ \fromto{\ast}{X} $ lifts to an equivalence of \categories
	\begin{equation*}
		\Omega \colon \equivto{\X_{\ast}^{\cn}}{\Grp(\X)} \period
	\end{equation*}
\end{remark}

%-------------------------------------------------------------------%
%  Actions of monoid objects                                        %
%-------------------------------------------------------------------%

\subsubsection{Actions of monoid objects}

In this subsection, we explain how to use \Cref{def:monoid_object} to give the correct definition of actions of monoid objects in \categories.
We begin with some motivation from ordinary category theory. 
Let $ \X $ be an ordinary category with finite products, and let $ M $ be a monoid object in $ \X $ with multiplication $ m \colon \fromto{M \cross M}{M} $ and unit $ u \colon \fromto{\ast}{M} $.
Recall that a \textit{(left) action} of $ M $ on an object $ X \in \X $ is a map $ a \colon \fromto{M \cross X}{X} $ such that the diagrams 
\begin{equation*}
	\begin{tikzcd}[sep=3em]
		M \isomorphic \ast \cross M \arrow[r, " u \cross \id{M}"] \arrow[dr, equals] & M \cross X \arrow[d, "a"] \\
		& X 
	\end{tikzcd}
	\andeq
	\begin{tikzcd}[sep=3em]
		M \cross M \cross X \arrow[r, "a \cross \id{M}"] \arrow[d, "\id{M} \cross a"'] & M \cross X \arrow[d, "a"] \\
		M \cross X \arrow[r, "a"'] & X 
	\end{tikzcd}
\end{equation*}
commute.
Like with monoid objects themselves, this data can be neatly packaged as a truncated simplicial diagram 
\begin{equation*}
    \begin{tikzcd}[sep=9em]
        M \cross M \cross X \arrow[r, "m \cross \idnosub \cross \idnosub "{description, near end}] \arrow[r, shift left=3ex, "\idnosub \cross a"] \arrow[r, shift right=3ex, "\pr_{2,3}"'] &%
        M \cross X \arrow[l, shift left=1.5ex, "\idnosub \cross u \cross \idnosub"{description, near end}] \arrow[l, shift right=1.5ex, "u \cross \idnosub \cross \idnosub"{description, near end}] \arrow[r, shift left=1.5ex, "a"] \arrow[r, shift right=1.5ex, "\pr_2"'] &%
        X \period \arrow[l, "u \cross \id{}" description] 
      \end{tikzcd}
\end{equation*} 
With this in mind, the following is the homotopy coherent definition of an action of a monoid object:

\begin{definition}[(action of a monoid object)]\label{def:monoid_action}
	Let $ \X $ be \acategory with finite products, and let $ M \colon \fromto{\Deltaop}{\X} $ be a monoid object in $ \X $.
	A \emph{(left) action} of $ M $ on an object of $ \X $ consists of:
	\begin{enumerate}[label=\stlabel{def:monoid_action}, ref=\arabic*]
		\item A simplicial object $ A \colon \fromto{\Deltaop}{\X} $.

		\item A map of simplicial objects $ p \colon \fromto{A}{M} $ such that for each $ [n] \in \DDelta $, the maps
		\begin{equation*}
			f([n]) \colon \fromto{A([n])}{M([n])} \andeq \fromto{A([n])}{A(\{n\})}
		\end{equation*}
		exhibit the $ n $-simplices $ A([n]) $ as the product $ M([n]) \cross A(\{n\}) $.
	\end{enumerate}
	In this case, we say that $ M $ \textit{acts} on the $ 0 $-simplices $ A_0 $.
	Given an object $ X \in \X $, an \textit{action of $ M $ on $ X $} is an action $ A \colon \fromto{\Deltaop}{\X} $ of $ M $ equipped with an identification $ A_0 \equivalent X $.

	We write
	\begin{equation*}
		\LMod_M(\X) \subset \Fun(\Deltaop,\X)_{/M}
	\end{equation*}
	for the full subcategory spanned by the (left) $ M $-actions.
	(See \HA{Proposition}{4.2.2.9}.)
\end{definition}

\begin{example}
	The terminal object $ \ast \in \X $ admits a unique $  M $-action: this is just the simplicial object $ M \colon \fromto{\Deltaop}{\X} $. 
\end{example}

\begin{definition}[(quotient by an action)]\label{def:action_quotient}
	Let $ \X $ be \acategory with finite products and geometric realizations, let $ M \colon \fromto{\Deltaop}{\X} $ be a monoid object in $ \X $, let $ X \in \X $, and let $ A \colon \fromto{\Deltaop}{\X} $ be an action of $ M $ on $ X $.
	The \emph{quotient} $ X \modmod M $ of $ X $ by the action of $ M $ is the geometric realization
	\begin{equation*}
		X \modmod M \colonequals \real{A} \period
	\end{equation*}
\end{definition}

\begin{example}[(classifying spaces of topological groups)]\label{ntn:classifyingspaceBG}
	Let $ G $ be a topological group.
	Since the underlying homotopy type functor $ \Piinf \colon \fromto{\Top}{\Spc} $ preserves finite products, the underlying homotopy type $ \Piinf(G) $ is naturally a group object of $ \Spc $.
	The \textit{classifying space} $ \BG $ is the quotient $ \ast \modmod \Piinf(G) $ of $\ast$ by the group object $ \Piinf(G) \in \Grp(\Spc) $.
	That is, $ \BG $ is the geometric realization of the simplicial space
	\begin{equation*}
		\begin{tikzcd}[sep=1.5em]
		    \cdots \arrow[r, shift left=0.75ex] \arrow[r, shift right=0.75ex] \arrow[r, shift right=2.25ex] \arrow[r, shift left=2.25ex] & \Piinf(G) \cross \Piinf(G) \arrow[l] \arrow[l, shift left=1.5ex] \arrow[l, shift right=1.5ex] \arrow[r] \arrow[r, shift left=1.5ex] \arrow[r, shift right=1.5ex] & \Piinf(G) \arrow[l, shift left=0.75ex] \arrow[l, shift right=0.75ex] \arrow[r, shift left=0.75ex] \arrow[r, shift right=0.75ex] & \ast \period \arrow[l] 
		\end{tikzcd}
	\end{equation*}
\end{example}

\begin{notation}[(classifying object)]\label{ntn:classifying_object}
	Let $ \X $ be \acategory with finite products and geometric realizations, and let $ M \colon \fromto{\Deltaop}{\X} $ be a monoid object in $ \X $.
	We write
	\begin{equation*}
		\Bup_{\X} M \colonequals \ast \modmod M = \real{M} 
	\end{equation*}
	and call $ \Bup_{\X} M $ the \emph{classifying object} of $ M $.

	Note that the unit defines a map $ \fromto{\ast}{M} $ from the constant simplicial object at the terminal object of $ \X $ to $ M $.
	Passing to geometric realizations, we see that $ \Bup_{\X} M $ admits a natural point $ \fromto{\ast}{\Bup_{\X} M} $.
	Whenever we regard $ \Bup_{\X} M $ as a pointed object, we use this natural point.
\end{notation}

\begin{remark}[(classifying spaces of Lie groups)]
	Let $ G $ be a Lie group.
	In most of this text, $ \BG $ denotes the classifying space of a Lie group in the classical sense: in \Cref{ntn:classifying_object}, $ \BG $ is the space $ \Bup_{\Spc}\Piinf(G) $.
	We are also interested in the classifying object of $ G $ regarded as an object of $ \Sh(\Mfld) $: \Cref{prop:Cech_nerve_trivG} shows that the latter classifying object coincides with the sheaf $ \BunG $.
	We've included the subscript $ \X $ in our notation for classifying objects in order to distinguish these two objects.
\end{remark}

\begin{observation}
	Let $ \X $ be \acategory with finite products and geometric realizations, and let $ M \colon \fromto{\Deltaop}{\X} $ be a monoid object in $ \X $.
	By definition, $ M $ is the terminal in $ \LMod_M(\X) $, hence the quotient by $ M $ (i.e., geometric realization) naturally refines to a functor
	\begin{equation*}
		(-) \modmod M \colon \fromto{\LMod_M(\X)}{\X_{/\Bup_{\X} M}} \period
	\end{equation*}
\end{observation}

\begin{proposition}[\SAG{Proposition}{E.6.4.4}]
	Let $ \X $ be \acategory with finite limits and geometric realizations, and let $ G $ be a group object of $ \X $.
	Assume that geometric realizations preserve finite products in $ \X $.
	Then the functor
	\begin{equation*}
		(-) \modmod G \colon \fromto{\LMod_{G}(\X)}{\X_{/\Bup_{\X} M}}
	\end{equation*}
	admits a right adjoint
	\begin{align*}
		\Fup_G \colon \X_{/\Bup_{\X} M} &\to \LMod_{G}(\X) \\
	\intertext{defined by}
		[f \colon \fromto{Y}{\Bup_{\X} G}] &\mapsto \Cech(\ast \to \Bup_{\X} G) \crosslimits_{\Bup_{\X} G} Y \period
	\end{align*}
	In particular, the composite
	\begin{equation*}
		\begin{tikzcd}[sep=2.5em]
			\X_{/\Bup_{\X} M} \arrow[r, "\Fup_G"] & \LMod_{G}(\X) \arrow[r, "\forget"] & \X
		\end{tikzcd}
	\end{equation*}
	is given by
	\begin{equation*}
		[f \colon \fromto{Y}{\Bup_{\X} G}] \mapsto \fib(f) \period
	\end{equation*}
\end{proposition}

\begin{theorem}[(actions via classifying objects)]\label{lem:actions_via_classifying_spaces}
	Let $ \X $ be \atopos and let $ G \colon \fromto{\Deltaop}{\X} $ be a group object in $ \X $.
	Then: 
	\begin{enumerate}[label=\stlabel{lem:actions_via_classifying_spaces}, ref=\arabic*]
		\item The functor
		\begin{equation*}
			(-) \modmod G \colon \fromto{\LMod_{G}(\X)}{\X_{/\Bup_{\X} M}}
		\end{equation*}
		is an equivalence of \categories with inverse $ \Fup_G $.

		\item Let $ X \in \X $, and let $ A \colon \fromto{\Deltaop}{\X} $ be an action of $ G $ on $ X $.
		There is a natural equivalence
		\begin{equation*}
			X \equivalence \fib(X \modmod G \to \Bup_{\X} G) \period
		\end{equation*}
	\end{enumerate}
\end{theorem}

\begin{proof}
 	See the proof of \SAG{Theorem}{E.6.5.1}.%
 	\footnote{Lurie's result \SAG{Theorem}{E.6.5.1} is in the context of profinite homotopy theory (rather than \topos theory).
 	However, Lurie's proof only uses a few categorical properties that are true in both of these settings.
 	Specifically, the conclusion of \Cref{lem:actions_via_classifying_spaces} holds in any \category $ \X $ with finite limits and geometric realizations in which geometric realizations are universal and geometric realizations of groupoid objects are effective.}
\end{proof}

%-------------------------------------------------------------------%
%-------------------------------------------------------------------%
%  The structure of BunG                                            %
%-------------------------------------------------------------------%
%-------------------------------------------------------------------%

\subsection{The structure of \texorpdfstring{$ \BunG $}{BunG}}\label{sec:structure_of_BunG}

In this section, we use the abstract material introduced in \cref{subsec:groupoids_effective_epis} to show that $ \BunG $ is the quotient $ \ast \modmod G $ in $ \Sh(\Mfld) $ (\Cref{prop:Cech_nerve_trivG}).
This presentation then lets us show that the space $ \Gammalowersharp(\BunG) $ recovers the classifying space $ \BG $ (\Cref{cor:formula_for_GammalowersharpBunG}).

The first step is to define an effective epimorphism $ \surjto{\ast}{\BunG} $.

\begin{construction}[($ \trivG $)]\label{cons:trivG}
	Let $ G $ be a Lie group.
	Define a global section
	\begin{equation*}
		\trivG \colon \fromto{\ast}{\BunG}
	\end{equation*}
	of sheaves of groupoids on $ \Mfld $ as follows. 
	For each manifold $ M $, the map
	\begin{equation*}
		\trivG(M) \colon \fromto{\ast}{\BunG(M)}
	\end{equation*}
	picks out the trivial $ G $-bundle $ \pr_{M} \colon M \cross G \to M $.
\end{construction}

\begin{lemma}\label{lem:trivG_effective_epi}
	Let $ G $ be a Lie group.
	Then the global section $ \trivG \colon \fromto{\ast}{\BunG} $ is an effective epimorphism in $ \Sh(\Mfld;\Spc) $.
\end{lemma}

\begin{proof}
	By \Cref{lem:ShMfld_criterion_for_effective_epi}, it suffices to check that for each $ n \geq 0 $, the groupoid $ \BunG(\RR^n) $ is connected.
	This is a consequence of the fact that every principal $ G $-bundle on a contractible manifold is trivializable.
\end{proof}

In light of \Cref{thm:effective_epis_and_groupoids,lem:trivG_effective_epi}, the Čech nerve of \smash{$ \trivG $} provides a presentation of \smash{$ \BunG $}.
Our next goal is to identify this Čech nerve with the simplical object of $ \Sh(\Mfld) $ defined by the group Lie $ G $.

\begin{notation}
	Let $ G $ be a Lie group, that is, a group object in the category of manifolds.
	In order to distinguish $ G $ thought of as a manifold from $ G $ thought of as a group object, we write 
	\begin{equation*}
		\Bar(G) \in \Fun(\Deltaop,\Sh(\Mfld))
	\end{equation*}
	for the group object of $ \Sh(\Mfld) $ defined by $ G $.
\end{notation}

Computing the pullback $ \ast \cross_{\BunG} \ast $ amounts to computing the automorphisms of trivial bundles. 

\begin{observation}[(automorphism groups of trivial bundles)]\label{obs:gauge_group_trivial}
	Let $ M $ be a manifold and $ G $ a Lie group.
	Given a principal $ G $-bundle $ \fromto{E}{M} $, write $ \Aut(E) $ for the group of automorphisms of $ E $ as a principal $ G $-bundle over $ M $.

	The map
	\begin{align*}
		\Cinf(M,G) & \to \Aut(M \cross G) \\ 
		f &\mapsto [(m,g) \mapsto (m,g \cdot f(m))]
	\end{align*}
	is an isomorphism of groups.
	The inverse is given by sending an isomorphism of principal $ G $-bundles $ \phi \colon \isomto{M \cross G}{M \cross G} $ to the composite
	\begin{equation*}
		\begin{tikzcd}[row sep=0.5em]
			M \isomorphic M \cross \{e\} \arrow[r, hooked] & M \cross G \arrow[r, "\phi", "\sim"'{yshift=0.15em}] & M \cross G \arrow[r, "\pr_G"] & G \\ 
			m \arrow[rrr, |->] & & & \phi(m,e) \period
		\end{tikzcd}
	\end{equation*}
\end{observation}

\begin{proposition}\label{prop:Cech_nerve_trivG}
	Let $ G $ be a Lie group.
	Then:
	\begin{enumerate}[label=\stlabel{prop:Cech_nerve_trivG}, ref=\arabic*]
		\item\label{prop:Cech_nerve_trivG.1} There is a natural equivalence of simplicial objects 
		\begin{equation*}
			\isomto{\Bar(G)}{\Cech(\trivG)} 
		\end{equation*}
		in $ \Sh(\Mfld) $.

		\item\label{prop:Cech_nerve_trivG.2} There is a natural equivalence $ \equivto{\ast \modmod G}{\BunG} $ in $ \Sh(\Mfld) $.
	\end{enumerate}
\end{proposition}

\begin{proof}
	For \enumref{prop:Cech_nerve_trivG}{1}, we begin by computing the $ n $-fold pullbacks $ \ast \cross_{\BunG} \cdots \cross_{\BunG} \ast $ along the effective epimorphism $ \trivG \colon \surjto{\ast}{\BunG}$.
	By definition, given a manifold $ M $, an object of the groupoid
	\begin{equation*}
		\Big\lparen \ast \crosslimits_{\BunG} \cdots \crosslimits_{\BunG} \ast \Big\rparen(M)
	\end{equation*}
	consists of a tuple
	\begin{equation*}
		(M \cross G, \ldots, M \cross G, \phi_1 \colon M \cross G \isomorphism M \cross G,\ldots,\phi_{n-1} \colon M \cross G \isomorphism M \cross G) 
	\end{equation*}
	of $ n $ trivial $ G $-bundles along with $ (n-1) $ isomorphisms of trivial $ G $-bundles
	\begin{equation*}
		\phi_1,\cdots,\phi_{n-1} \colon \isomto{M \cross G}{M \cross G} \period
	\end{equation*}
	Moreover, the only morphisms in this groupoid are the identities; that is, $ \ast \cross_{\BunG} \cdots \cross_{\BunG} \ast $ is a sheaf of sets.
	\Cref{obs:gauge_group_trivial} provides natural isomorphisms
	\begin{align*}
		\ast \crosslimits_{\BunG} \cdots \crosslimits_{\BunG} \ast &\isomorphic \Cinf(-,G) \cross \cdots \cross \Cinf(-,G) \\ 
		&= \yo(G) \cross \cdots \cross \yo(G) \period
	\end{align*}
	Thus for each $ n \geq 0 $ we have provided natural isomorphisms
	\begin{equation*}
		\Cech(\trivG)_n \isomorphic \Bar(G)_n \period
	\end{equation*}
	It is immediate from the definitions that these isomorphisms are compatible with the simplicial identities, proving \enumref{prop:Cech_nerve_trivG}{1}. 

	To conclude, note that \enumref{prop:Cech_nerve_trivG}{1} and \Cref{lem:trivG_effective_epi} imply \enumref{prop:Cech_nerve_trivG}{2}.
\end{proof}

Using this presentation, we compute the homotopification of $ \BunG $:

\begin{corollary}\label{cor:formula_for_GammalowersharpBunG}
	Let $ G $ be a Lie group.
	Then:
	\begin{enumerate}[label=\stlabel{cor:formula_for_GammalowersharpBunG}, ref=\arabic*]
		\item\label{cor:formula_for_GammalowersharpBunG.1} There is a natural equivalence of spaces $ \Gammalowersharp(\BunG) \equivalent \BG $.
		
		\item\label{cor:formula_for_GammalowersharpBunG.2} There is a natural equivalence $ \Lhi(\BunG) \equivalent \Gammaupperstar(\BG) $ of sheaves on $ \Mfld $.
	\end{enumerate}
\end{corollary}

\begin{proof}
	First we show \enumref{cor:formula_for_GammalowersharpBunG}{1}.
	Since $ \BunG \equivalent |\Bar(G)| $ and $ \Gammalowersharp $ preserves colimits and on $ \Mfld $ agrees with the underlying homotopy type functor $ \Piinf \colon \fromto{\Mfld}{\Spc} $, we have natural equivalences
	\begin{align*}
		\Gammalowersharp(\BunG) &\equivalent \Gammalowersharp|\Bar(G)| \\ 
		&\equivalent |\Gammalowersharp\Bar(G)| \\ 
		&\equivalent |\Bar(\Piinf(G))| \\ 
		&\equivalent \BG \period
	\end{align*}

	To conclude note that \enumref{cor:formula_for_GammalowersharpBunG}{2} follows from \enumref{cor:formula_for_GammalowersharpBunG}{1} and the definition $ \Lhi = \Gammaupperstar\Gammalowersharp $.
\end{proof}

\begin{remark}
	Due to \Cref{prop:Cech_nerve_trivG,cor:formula_for_GammalowersharpBunG}, Freed and Hopkins denoted the sheaf $ \BunG $ by $ \Bup_{\bullet}G $.
\end{remark}

\begin{warning}
	As \Cref{cor:formula_for_GammalowersharpBunG} demonstrates, the sheaf $ \BunG $ is \textit{not} generally \RRinvariant.
	This might seem surprising: one often quotes the classical result that `principal $ G $-bundles are homotopy-invariant'.
	However, what this classical result says is that the \textit{set of isomorphism classses} of principal $ G $-bundles is homotopy-invariant.
	On the other hand, the \textit{groupoid} of principal $ G $-bundles is not homotopy-invariant!
\end{warning}

%-------------------------------------------------------------------%
%-------------------------------------------------------------------%
%  The structure of Bun_G^∇                                         %
%-------------------------------------------------------------------%
%-------------------------------------------------------------------%

\subsection{The structure of \texorpdfstring{$ \BunGnabla $}{BunG∇}}\label{sec:structure_of_BunGnabla}

The purpose of this section is to give a presentation of $ \BunGnabla $ as the quotient of $ \Omega^1(-;\g) $ by the `adjoint action' of $ G $ (see \Cref{def:adjoint_action,rem:adjoint_action_explicitly}).
To do this, we use the equivalence
\begin{equation*}
	(-) \modmod G \colon \fromto{\LMod_{G}(\Sh(\Mfld))}{\Sh(\Mfld)_{/\BunG}}
\end{equation*}
provided by \Cref{lem:actions_via_classifying_spaces,prop:Cech_nerve_trivG} and compute the fiber of the forgetful map as 
\begin{equation*}
	\fib(\fromto{\BunGnabla}{\BunG}) \equivalent \Omega^1(-;\g) 
\end{equation*}
(\Cref{lem:fib_of_BunGnabla_to_BunG}).

In order to compute this fiber, we first need to define a map
\begin{equation*}
	\trivGnabla \colon \surjto{\Omega^1(-;\g)}{\BunGnabla} \period
\end{equation*} 
As the notation suggests, $ \trivGnabla $ sends a form $ \omega \in \Omega^1(M;\g) $ to the trivial $ G $-bundle $ M \cross G $ with a connection involving $ \omega $.
In order to give a formula for this connection, we start in \cref{subsec:Maurer-Cartan} by explaining some background material on why every connection on a trivial $ G $-bundle takes a particular form (\Cref{lem:connection_trivial_bundle}).
In \cref{subsec:presentation_of_BunGnabla} we prove that
\begin{equation*}
	\BunGnabla \equivalent \Omega^1(-;\g) \modmod G
\end{equation*}
(\Cref{prop:Cech_nerve_trivGnabla}).
In \cref{subsec:homotopification_of_BunGnabla}, we use this presentation to show that $ \Gammalowersharp(\BunGnabla) \equivalent \BG $ (\Cref{cor:formula_for_Gammalowersharp_BunGnabla}).
Finally, in \cref{subsec:G-bundles_flat_connection} we show that the global sections of $ \BunG $ and $ \BunGnabla $ recover the classifying space of $ G $ equipped with the \textit{discrete} topology (\Cref{lem:BGdisc}).

%-------------------------------------------------------------------%
%  Maurer–Cartan forms & connections on trivial bundles             %
%-------------------------------------------------------------------%

\subsubsection{Maurer--Cartan forms \& connections on trivial bundles}\label{subsec:Maurer-Cartan}

Every Lie group admits a canonical $ 1 $-form valued in its Lie algebra called the \textit{Maurer--Cartan form}.
In this subsection, we introduce Maurer--Cartan forms and use them to explain why all connections on a trivial $ G $-bundle have a very particular form (\Cref{lem:connection_trivial_bundle}).

To define the Maurer--Cartan form, let us fix some notation.

\begin{notation}\label{ntn:Lie_group_stuff}
	Let $ G $ be a Lie group.
	\begin{enumerate}[label=\stlabel{ntn:Lie_group_stuff}, ref=\arabic*]
		\item We write $ e \in G $ for the identity element, and $ \g \colonequals \Tan_{e}G $ for the Lie algebra of $ G $.

		\item For each $ g \in G $, we write
		\begin{align*}
			\Lup_g \colon G \to G &\andeq  \Rup_g \colon G \to G \\
			\phantom{\Lup_g \colon} h \mapsto gh &\phantom{\andeq\Rup_g \colon} h \mapsto hg 
		\end{align*}
		for the maps given by left and right multiplication by $ g $, respectively.
	
		\item We write $ \Ad \colon \fromto{G}{\Aut(\g)} $ for the \textit{adjoint action} of $ G $ on $ \g $.
		That is, $ \Ad $ is the derivative of the conjugation action $ \fromto{G}{\Aut(G)} $ at the identity element $ e \in G $.
	\end{enumerate}
\end{notation}

Geometrically, the Maurer--Cartan form is defined by using the fact that the tangent bundle of a Lie group naturally splits:

\begin{observation}\label{obs:TG_splits}
	Let $ G $ be a Lie group.
	There is a natural splitting $ \equivto{\Tan G}{G \cross \g} $ defined by
	\begin{equation*}
		(g,v) \mapsto (g, (\Lup_{g^{-1}})_{\ast}(v)) \period
	\end{equation*}
	Here note that $ (\Lup_{g^{-1}})_{\ast}(v) $ is an element of $ \Tan_{g^{-1} g}G = \g $.
\end{observation}

\begin{definition}[(Maurer--Cartan form)]\label{def:Maurer-Cartan_form}
	Let $ G $ be a Lie group.
	The \textit{Maurer--Cartan form}\index[terminology]{Maurer--Cartan form} of $ G $ is the the $ \g $-valued $ 1 $-form $ \MC \in \Omega^1(G;\g) $\index[notation]{MC@$\MC$} defined by the composite
	\begin{equation*}
		\begin{tikzcd}
			\Tan G \arrow[r, "\sim"{yshift=-0.25em}] & G \cross \g \arrow[r, "\pr_{\g}"] & \g 
		\end{tikzcd}
	\end{equation*}
	of the splitting of \Cref{obs:TG_splits} with the projection.
	Explicitly, the Maurer--Cartan form is defined by
	\begin{equation*}
		\MC_{g}(v) \colonequals (\Lup_{g^{-1}})_{\ast}(v) \period
	\end{equation*}
\end{definition}

\begin{remark}
	The Maurer--Cartan form of is often written as $ g^{-1}\d g $.
\end{remark}

It follows immediately from the definitions that the Maurer--Cartan form is the unique left-invariant $ \g $-valued $ 1 $-form on $ G $ that is the identity on $ \Tan_e G = \g $:

\begin{proposition}[(characterization of the Maurer--Cartan form)]
	Let $ G $ be a Lie group.
	The Maurer--Cartan form of $ G $ is the unique $ \g $-valued $ 1 $-form $ \MC \in \Omega^1(G;\g) $ satisfying the following properties:
	\begin{enumerate}[label=\stlabel{def:Maurer-Cartan_form}, ref=\arabic*]
		\item The map $ \MC_e \colon \fromto{\Tan_e G}{\g} $ is the identity map.

		\item For each $ g \in G $, we have
		\begin{equation*}
			\MC_g = \Ad_{g}(\Rup_{g}^{\ast} \MC_{e}) \period
		\end{equation*}
		Here $ \Rup_{g}^{\ast} $ denotes the pullback of forms under right translation by $ g $.
	\end{enumerate}
\end{proposition}

\begin{observation}\label{obs:MC_unique_connection}
	Rephrasing \Cref{def:Maurer-Cartan_form}, the Maurer--Cartan form $ \MC \in \Omega^1(G;\g) $ is the \textit{unique} connection on the trivial $ G $-bundle $ \fromto{G}{\ast} $.
\end{observation}

Using the Maurer--Cartan form, we can also see that every connection on a trivial $ G $-bundle on a manifold has a very particular form.
Namely, they are all obtained by pulling back the Maurer--Cartan form from $ G $ and adding a form that lives `horizontally' over $ M $.

\begin{observation}
	Let $ M $ be a manifold and $ G $ a Lie group.
	Since the map $ \pr_{M} \colon \fromto{M \cross G}{M} $ admits a section, the pullback map $ \prupperstar_M \colon \Omega^1(M;\g) \to \Omega^1(M \cross G;\g) $ is injective.
	Hence the map
	\begin{align*}
		\Omega^1(M;\g) &\to \Omega^1(M \cross G;\g) \\ 
		\omega &\mapsto \prupperstar_M(\omega) + \prupperstar_G(\MC)
	\end{align*}
	is also injective.
\end{observation}

\begin{lemma}[(connections on trivial bundles)]\label{lem:connection_trivial_bundle}
	Let $ M $ be a manifold and $ G $ a Lie group.
	Write $ i_e \colon \incto{M}{M \cross G} $ for the identity section $ \goesto{m}{(m,e)} $.
	Then:
	\begin{enumerate}[label=\stlabel{lem:connection_trivial_bundle}, ref=\arabic*]
		\item Given a $ 1 $-form $ \omega \in \Omega^1(M;\g) $, the $ \g $-valued $ 1 $-form $ \prupperstar_M(\omega) + \prupperstar_G(\MC) $ is a connection $ 1 $-form on the trivial $ G $-bundle $ M \cross G $.

		\item If $ \theta \in \Omega^1(M \cross G;\g) $ is a connection $ 1 $-form on the trivial $ G $-bundle $ M \cross G $, then
		\begin{equation*}
			\theta = \prupperstar_M\iupperstar_e(\theta) + \prupperstar_G(\MC) \period
		\end{equation*}

		\item The image of the injection
		\begin{align*}
			\Omega^1(M;\g) &\inclusion \Omega^1(M \cross G;\g) \\ 
			\omega &\mapsto \prupperstar_M(\omega) + \prupperstar_G(\MC)
		\end{align*}
		is the subset of connection $ 1 $-forms on the trivial $ G $-bundle. 
	\end{enumerate}
\end{lemma}

%-------------------------------------------------------------------%
%  The presentation of BunG∇                                        %
%-------------------------------------------------------------------%

\subsubsection{The presentation of \texorpdfstring{$ \BunGnabla $}{BunG∇}}\label{subsec:presentation_of_BunGnabla}

We now appeal to \Cref{lem:actions_via_classifying_spaces} to show $ \BunGnabla $ admits a presentation as $ \Omega^1(-;\g) \modmod G $.
To do this, we start by giving an explicit description of the fiber of the forgetful map $ \fromto{\BunGnabla}{\BunG} $.

\begin{notation}[($ \BunGnablatriv $)]\label{ntn:BunGnablatriv}
	Let $ G $ be a Lie group. 
	Write $ \BunGnablatriv $ for the pullback
	\begin{equation*}
		\begin{tikzcd}[sep=2.5em]
			\BunGnablatriv \arrow[r, ->>] \arrow[d] \arrow[dr, phantom, "\lrcorner"{description, very near start}] & \BunGnabla \arrow[d, "\forget"] \\
			\ast \arrow[r, "\trivG"', ->>] & \BunG \period
 		\end{tikzcd}
	\end{equation*}
\end{notation}

\begin{observation}[(explicit description of $ \BunGnablatriv $)]\label{obs:explicit_description_of_BunGnablatriv}
	Let $ G $ be a Lie group. 
	By the explicit description of pullbacks of groupoids, for each manifold $ M $, objects of \smash{$ \BunGnablatriv(M) $} consist of triples
	\begin{equation*}
		(P,\theta \in \Omega^1(P;\g),\phi \colon \isomto{M \cross G}{P}) \comma
	\end{equation*} 
	where $ \fromto{P}{M} $ is a $ G $-bundle, $ \theta $ is a connection on $ P $, and $ \phi $ is a trivialization of $ P $.
	A morphism $ \fromto{(P_1,\theta_1,\phi_1)}{(P_1,\theta_1,\phi_1)} $ is an isomorphism of $ G $-bundles
	\begin{equation*}
		f \colon \isomto{P_1}{P_2}
	\end{equation*}
	such that $ \fupperstar(\omega_2) = \omega_1 $ and the triangle of isomorphisms
	\begin{equation*}
		\begin{tikzcd}[row sep=3em, column sep=1em]
			& M \cross G \arrow[dl, "\phi_1"'] \arrow[dr, "\phi_2"] & \\
			P_1 \arrow[rr, "f"'] & & P_2 
		\end{tikzcd}
	\end{equation*}
	commutes.
\end{observation}

\begin{lemma}\label{lem:BunGnablatriv_is_0-truncated}
	Let $ G $ be a Lie group.
	Then the sheaf $ \BunGnablatriv $ is $ 0 $-truncated (i.e., equivalent to a sheaf of sets).
\end{lemma}

\begin{proof}
	Let $ M $ be a manifold and $ f \colon \fromto{(P_1,\theta_1,\phi_1)}{(P_1,\theta_1,\phi_1)} $ a morphism in $ \BunGnablatriv(M) $.
	Note that the condition $ \phi_2 = f\phi_1 $ uniquely determines $ f $: we necessarily have $ f = \phi_2 \phi_1^{-1} $.
	Hence the condition $ \fupperstar(\omega_2) = \omega_1 $ is equivalent to the condition $ \phiupperstar_1(\omega_1) = \phiupperstar_2(\omega_2) $.
	Thus we have 
	\begin{equation*}
		\Map_{\BunGnablatriv(M)}((P_1,\theta_1,\phi_1),(P_1,\theta_1,\phi_1)) = 
		\begin{cases}
			\{\phi_2 \phi_1^{-1}\} & \phiupperstar_1(\omega_1) = \phiupperstar_2(\omega_2) \\ 
			\emptyset & \text{otherwise} \period 
		\end{cases}
		\qedhere
	\end{equation*}
\end{proof}

Using this description of $ \BunGnablatriv $, we now provide an equivalence $ \equivto{\Omega^1(-;\g)}{\BunGnablatriv} $.

\begin{convention}
	Let $ V $ be an $ \RR $-vector space and $ n \geq 0 $ an integer. 
	Throughout this chapter, we regard $ \Omega^n(-;V) $ as a sheaf of \textit{sets} (hence a sheaf of spaces) on $ \Mfld $.
\end{convention}

\begin{construction}[($ \trivGnabla $)]\label{cons:trivGnabla}
	Let $ G $ be a Lie group.
	Define a morphism
	\begin{equation*}
		\trivGnabla \colon \fromto{\Omega^1(-;\g)}{\BunGnabla}
	\end{equation*}
	of sheaves of groupoids on $ \Mfld $ as follows. 
	For each manifold $ M $, the map $ \trivGnabla(M) $ is given by sending a $ 1 $-form $ \omega \in \Omega^1(M;\g) $ to the trivial $ G $-bundle $ \pr_{M} \colon M \cross G \to M $ equipped with the connection
	\begin{align*}
		\prupperstar_M(\omega) + \prupperstar_G(\MC) \period
	\end{align*}
\end{construction}

\begin{lemma}\label{lem:fib_of_BunGnabla_to_BunG}
	Let $ G $ be a Lie group.
	Then the natural commutative square
	\begin{equation*}
		\begin{tikzcd}[sep=2.5em]
			\Omega^1(-;\g) \arrow[r, "\trivGnabla"] \arrow[d] & \BunGnabla \arrow[d, "\forget"] \\
			\ast \arrow[r, "\trivG"', ->>] & \BunG
 		\end{tikzcd}
	\end{equation*}
	is a pullback square in $ \Sh(\Mfld) $.
	In particular, $ \trivGnabla $ is an effective epimorphism.
\end{lemma}

\begin{proof}
	Note that for each manifold $ M $, the map of groupoids
	\begin{align*}
		t_M \colon \fromto{\Omega^1(M;\g)}{\BunGnablatriv(M)}
	\end{align*}
	induced by the universal property of the pullback is given by the assignment
	\begin{align*}
		\omega \mapsto (M \cross G, \prupperstar_M(\omega) + \prupperstar_G(\MC), \idnosub \colon\isomto{M \cross G}{M \cross G}) \period
	\end{align*}
	The map $ t_M $ is fully faithful because the groupoids $ \Omega^1(M;\g) $ and $ \BunGnablatriv(M) $ are both equivalent to sets (\Cref{lem:BunGnablatriv_is_0-truncated}).
	\Cref{lem:connection_trivial_bundle} implies that $ t $ is essentially surjective; hence $ t $ is an equivalence, as desired.
\end{proof}

We now explain why \Cref{lem:fib_of_BunGnabla_to_BunG} gives rise to a presentation of $ \BunGnabla $ as a quotient $ \Omega^1(-;\g) \modmod G $.

\begin{observation}\label{obs:adjoint_action}
	Let $ G $ be a Lie group.
	By \Cref{prop:Cech_nerve_trivG}, $ \BunG \equivalent \ast \modmod G $ in $ \Sh(\Mfld) $.
	Thus \Cref{lem:actions_via_classifying_spaces} provides an equivalence of \categories
	\begin{equation*}
		(-) \modmod G \colon \equivto{\LMod_G(\Sh(\Mfld))}{\Sh(\Mfld)_{/\BunG}} \period
	\end{equation*}
	The inverse sends an object $ f \colon \fromto{E}{\BunG} $ to the fiber $ \fib(f) $ equipped with a $ G $-action.
	By \Cref{lem:fib_of_BunGnabla_to_BunG}, applying this inverse equivalence to the forgetful map
	\begin{equation*}
		\fromto{\BunGnabla}{\BunG}
	\end{equation*}
	defines a $ G $-action on $ \Omega^1(-;\g) $.
\end{observation}

\begin{definition}[(adjoint action)]\label{def:adjoint_action}
	Let $ G $ be a Lie group.
	We refer to the $ G $-action on $ \Omega^1(-;\g) $ described in \Cref{obs:adjoint_action} as the \emph{adjoint action}.
\end{definition}

\begin{remark}[(the adjoint action, explicitly)]\label{rem:adjoint_action_explicitly}
	Unwinding the definitions shows that the adjoint action admits the following explicit description.
	Given a manifold $ M $, map $ \phi \colon \fromto{M}{G} $, and $ 1 $-form $ \omega \in \Omega^1(M;\g) $, write $ \Ad_{\phi} \omega \in \Omega^1(M;\g) $ for the $ 1 $-form defined by
	\begin{equation*}
		m \mapsto \Ad_{\phi(m)} \omega_m \period
	\end{equation*}
	Then the adjoint action of $ G $ on $ \Omega^1(-;\g) $ is given by
	\begin{align*}
		\Cinf(M,G) \cross \Omega^1(M;\g) &\to \Omega^1(M;\g) \\ 
		(\phi,\omega) &\mapsto \Ad_{\phi} \omega \period
	\end{align*}
\end{remark}

\begin{corollary}\label{prop:Cech_nerve_trivGnabla}
	Let $ G $ be a Lie group.
	Then there is a natural equivalence
	\begin{equation*}
		\equivto{\Omega^1(-;\g) \modmod G}{\BunGnabla} \period
	\end{equation*}
	from the quotient of $ \Omega^1(-;\g) $ by the adjoint action to $ \BunGnabla $.
\end{corollary}

\begin{proof}
	Immediate from \Cref{lem:actions_via_classifying_spaces} and the definition of the adjoint action.
\end{proof}

\begin{remark}[($ \Bup_{\nabla}G $)]
	Due to \Cref{prop:Cech_nerve_trivGnabla}, Freed and Hopkins denoted the sheaf $ \BunGnabla $ by $ \Bup_{\nabla}G $ \cite[Example 5.11]{FreedHopkins}.
\end{remark}

\begin{remark}[($ \EnablaG $)]
	The sheaf \smash{$ \BunGnablatriv \isomorphic \Omega^1(-;\g) $} has the following alternative description. 
	Write $ \EnablaG $ for the sheaf on $ \Mfld $ that assigns a manifold $ M $ the groupoid of triples
	\begin{equation}\label{eq:EnablaG_data}
		(P, s \colon M \to P, \theta \in \Omega^1(P;\g)) \comma
	\end{equation}
	where $ P \to M $ is a principal $ G $-bundle on $ M $, $ s $ is a global section of $ P $, and $ \theta $ is a connection on $ P $.
	The morphisms in $ \EnablaG(M) $ are isomorphisms of principal bundles preserving the specified sections and $ 1 $-forms.
	Since the data of a section of a principal $ G $-bundle is equivalent to the a trivialization, there are isomorphisms
	\begin{equation*}
		\EnablaG \isomorphic \BunGnablatriv \isomorphic \Omega^1(-;\g) \period 
	\end{equation*}
	Freed and Hopkins use this alternative description in \cite{FreedHopkins}.

	% Consider the map
	% \begin{equation*}
	% 	\alpha \colon \fromto{\EnablaG(M)}{\Omega^1(M;\g)}
	% \end{equation*}
	% that sends a triple \eqref{eq:EnablaG_data} to the $ \g $-valued $ 1 $-form $ \supperstar(\theta) $.
	% Also consider the map
	% \begin{equation*}
	% 	\beta \colon \fromto{\Omega^1(M;\g)}{\EnablaG(M)}
	% \end{equation*}
	% sending a $ 1 $-form $ \omega \in \Omega^1(M;\g) $ to the trivial $ G $-bundle $ \pr_{M} \colon M \cross G \to M $ equipped with the identity section and connection $ 1 $-form
	% \begin{align*}
	% 	\prupperstar_M(\omega) + \prupperstar_G(\MC) \period
	% \end{align*}
	% By definition, the composite $ \alpha\beta $ is the identity on $ \Omega^1(M;\g) $. 
	% \Cref{lem:connection_trivial_bundle} implies that the composite $ \beta\alpha $ is also equivalent to the identity.
	% That is, $ \alpha $ and $ \beta $ define an equivalence of sheaves 
	% \begin{equation*}
	% 	\Omega^1(-;\g) \equivalent \EnablaG \period
	% \end{equation*}
	% The map $ \trivGnabla \colon \fromto{\Omega^1(-;\g)}{\BunGnabla} $ is the composite of the equivalence $ \beta \colon \equivto{\Omega^1(-;\g)}{\EnablaG} $ with the map $ \fromto{\EnablaG}{\BunGnabla} $ forgetting the additional data of a trivializing section.
\end{remark}

%-------------------------------------------------------------------%
%  The homotopification of BunG∇                                    %
%-------------------------------------------------------------------%

\subsubsection{The homotopification of \texorpdfstring{$ \BunGnabla $}{BunG∇}}\label{subsec:homotopification_of_BunGnabla}

We now use the presentation $ \BunGnabla \equivalent \Omega^1(-;\g) \modmod G $ to compute the homotopification of $ \BunGnabla $.
We begin by showing that the homotopification of $ \Omega^1(-;\g) $ is trivial.

\begin{notation}\label{ntn:Cinf_ring}
	We write $ \Cinf $ for the sheaf of $ \RR $-algebras on $ \Mfld $ given by the assignment
	\begin{equation*}
		M \mapsto \Cinf(M,\RR) \comma
	\end{equation*}
	with pointwise addition and multiplication.
	We regard $ \Cinf $ as an object of the \category $ \Sh(\Mfld;\Vect(\RR)) $.
\end{notation}

\begin{observation}
	As a sheaf of sets, $ \Cinf $ is just the sheaf represented by $ \RR $.
	We have introduced the notation \Cref{ntn:Cinf_ring} to distinguish when we want to think of the sheaf of rings represented by $ \RR $ or the sheaf of sets represented by $ \RR $.
\end{observation}

\begin{lemma}\label{lem:homotopification_kills_Cinf-modules}
	Let $ E \in \Sh(\Mfld;\Vect(\RR)) $.
	If $ E $ admits the structure of a $ \Cinf $-module, then the composite
	\begin{equation*}
		\begin{tikzcd}[sep=3em]
			\Sh(\Mfld;\Vect(\RR)) \arrow[r, "\forget"] & \Sh(\Mfld;\Spc) \arrow[r, "\Gammalowersharp"] & \Spc
		\end{tikzcd}
	\end{equation*}
	carries $ E $ to the terminal object.
\end{lemma}

\begin{proof}
	To prove the claim, it suffices to show that the identity map $ \fromto{\Gammalowersharp(E)}{\Gammalowersharp(E)} $ is nullhomotopic.
	Write $ i_0,i_1 \colon \fromto{\ast}{\Cinf} $ for the global sections specified by $ 0,1 \in \RR $, respectively.
	Since $ E $ admits the structure of a $ \Cinf $-module, there exist a point $ 0 \colon \fromto{\ast}{E} $ and a multiplication map 
	\begin{equation*}
		m \colon \fromto{\Cinf \cross E}{E}
	\end{equation*}
	such that the diagram
	\begin{equation}\label{diag:Cinf-module_structure}
		\begin{tikzcd}[sep=3em]
			E \isomorphic \ast \cross E \arrow[r, "i_0 \cross \id{E}"] \arrow[dr, "0_{E}"'] & \Cinf \cross E \arrow[d, "m" description] & \ast \cross E \isomorphic E \arrow[l, "i_1 \cross \id{E}"'] \arrow[dl, equals] \\
			& E & 
		\end{tikzcd}
	\end{equation}
	commutes.
	Since the manifold $ \RR $ is contractible and $ \Gammalowersharp \colon \fromto{\Mfld}{\Spc} $ sends a manifold to its underlying homotopy type, we have $ \Gammalowersharp(\Cinf) \equivalent \ast $.
	Hence $ \Gammalowersharp(i_0) $ and $ \Gammalowersharp(i_1) $ are equivalences and
	\begin{equation*}
		\Gammalowersharp(i_0) \equivalent \Gammalowersharp(i_1) \period
	\end{equation*}
	Since $ \Gammalowersharp $ preserves finite products (\Cref{cor:formula_for_Gammalowersharp_preserves_finite_products}), the commutativity of the diagram \eqref{diag:Cinf-module_structure} shows that there are equivalences
	\begin{align*}
		\id{\Gammalowersharp(E)} &\equivalent \Gammalowersharp(m) \of (\Gammalowersharp(i_1) \cross \Gammalowersharp(\id{E})) \\
		&\equivalent \Gammalowersharp(m) \of (\Gammalowersharp(i_0) \cross \Gammalowersharp(\id{E})) \\ 
		&\equivalent \Gammalowersharp(0_{E}) \period
	\end{align*}
	Hence the map $ \id{\Gammalowersharp(E)} $ is nullhomotopic, as desired.
\end{proof}

\begin{remark}[(on the proof of \Cref{lem:homotopification_kills_Cinf-modules})]
	One can strengthen \Cref{lem:homotopification_kills_Cinf-modules} to the following claim: if $ E \in \Sh(\Mfld;\D(\RR)) $ admits the structure of a $ \Cinf $-module, then $ \Gammalowersharp(E) \equivalent 0 $.
	The point is that from the colimit formula for $ \Gammalowersharp $ (\Cref{cor:formula_for_Gammalowersharp}), one can show that 
	\begin{equation*}
		\Gammalowersharp \colon \fromto{\Sh(\Mfld;\D(\RR))}{\D(\RR)}
	\end{equation*}
	admits a canonical \textit{lax} monoidal structure with respect to the tensor products on both sides.
	In particular, $ \Gammalowersharp $ preserves algebras and modules over algebras.
	Hence $ \Gammalowersharp(E) $ is a module over the algebra $ \Gammalowersharp(\Cinf) $ in $ \D(\RR) $.
	The fact that the topological space $ \RR $ is contractible implies that $ \Gammalowersharp(\Cinf) $ is zero; since every module over the zero algebra is also zero, $ \Gammalowersharp(E) \equivalent 0 $.

	This sketch is the real idea behind the proof of \Cref{lem:homotopification_kills_Cinf-modules}.
	However, to spell out this argument in full detail requires a number of technical digressions that we do not need to use elsewhere. 
	Thus we decided to give a more direct proof of the specific result we need in order to compute $ \Gammalowersharp(\BunGnabla) $.
\end{remark}

\begin{example}\label{ex:homotopification_of_Omega^n}
	Let $ V $ be an $ \RR $-vector space and $ n \geq 0 $ an integer.
	The sheaf $ \Omega^n(-;V) $ is a $ \Cinf $-module with multiplication defined by
	\begin{align*}
		\Cinf(M,\RR) \cross \Omega^n(M;V) &\to \Omega^n(M;V) \\ 
		(f,\omega) &\mapsto f\omega \period 
	\end{align*}
	\Cref{lem:homotopification_kills_Cinf-modules} shows that $ \Gammalowersharp(\Omega^n(-;V)) \equivalent \ast $.
\end{example}

\begin{warning}
	The subsheaf $ \Omegacl^n \subset \Omega^n $ of closed $ n $-forms is not a $ \Cinf $-module.
	The reason is that multiplication by a function is not compatible with the de Rham differential: given a function $ f \in \Cinf(M,\RR) $ and form $ \omega \in \Omega^n(M) $, we have
	\begin{equation*}
		\d(f\omega) = \d f \wedge \omega - f\d\omega \period
	\end{equation*}
	Moreover, even if $ \omega $ is closed, $ \d f \wedge \omega $ need not be zero.
\end{warning}

We now compute the homotopification of $ \BunGnabla $.

\begin{corollary}\label{cor:formula_for_Gammalowersharp_BunGnabla}
	Let $ G $ be a Lie group.
	Then:
	\begin{enumerate}[label=\stlabel{cor:formula_for_Gammalowersharp_BunGnabla}, ref=\arabic*]
		\item\label{cor:formula_for_Gammalowersharp_BunGnabla.1} The forgetful morphism $ \fromto{\BunGnabla}{\BunG} $ induces an equivalence
		\begin{equation*}
			\equivto{\Gammalowersharp(\BunGnabla)}{\Gammalowersharp(\BunG)} \period
		\end{equation*}

		\item\label{cor:formula_for_Gammalowersharp_BunGnabla.2} There is a natural equivalence $ \Gammalowersharp(\BunGnabla) \equivalent \BG $.
	\end{enumerate}
\end{corollary}

\begin{proof}
	First we prove \enumref{cor:formula_for_Gammalowersharp_BunGnabla}{1}.
	By definition, the map $ \fromto{\BunGnabla}{\BunG} $ given by forgetting connection $ 1 $-forms is induced on geometric realizations by a map
	\begin{equation*}
		G^{\cross (\bullet-1)} \cross \Omega^1(-;\g) \to G^{\cross (\bullet-1)}
	\end{equation*}
	from the simplicial object defining the adjoint action of $ G $ on $ \Omega^1(-;\g) $ to the bar construction of $ G $.
	Moreover, on each term this map is the projection.
	Hence using the fact that the left adjoint $ \Gammalowersharp $ preserves finite products (\Cref{cor:formula_for_Gammalowersharp_preserves_finite_products}), we compute
	\begin{align*}
		\Gammalowersharp(\BunGnabla) &\equivalent \Gammalowersharp\real{ G^{\cross(\bullet -1)} \cross \Omega^1(-;\g)} \\ 
		&\equivalent \real{\Gammalowersharp(G^{\cross(\bullet -1)}) \cross \Gammalowersharp(\Omega^1(-;\g))} \\ 
		&\equivalent \real{\Gammalowersharp(G^{\cross(\bullet -1)}) \cross \ast} && (\textup{\Cref{ex:homotopification_of_Omega^n}}) \\
		&\equivalent \Gammalowersharp\real{\Bar(G)} \\ 
		&\equivalent \Gammalowersharp(\BunG) && (\textup{\Cref{prop:Cech_nerve_trivG}}) \period
	\end{align*}
	Item \enumref{cor:formula_for_Gammalowersharp_BunGnabla}{2} now follows from \Cref{cor:formula_for_GammalowersharpBunG}. 
\end{proof}

Later on, we make use of the following suspension spectra variant of \Cref{cor:formula_for_Gammalowersharp_BunGnabla}:

\begin{corollary}\label{cor:formula_for_Gammalowersharp_Suspension_spectrum_BunG}
	Let $ G $ be a Lie group.
	There are natural equivalences of spectra
	\begin{equation*}
		\Gammalowersharp(\Sigma_+^\infty\BunGnabla ) \equivalent \Gammalowersharp(\Sigma_+^\infty\BunG)  \equivalent \Sigma_+^\infty \BG \period
	\end{equation*}
\end{corollary}

\begin{proof}
	Combine the compatibility of $ \Gammalowersharp $ and $ \Sigma_+^\infty $ (\Cref{cor:suspension_spectra_ShMfld}) with \Cref{cor:formula_for_Gammalowersharp_BunGnabla}. 
\end{proof}

%-------------------------------------------------------------------%
%  G-bundles with flat connection                                   %
%-------------------------------------------------------------------%

\subsubsection{\texorpdfstring{$ G $}{G}-bundles with flat connection}\label{subsec:G-bundles_flat_connection}

In this final subsection, we prove that the global sections of $ \BunG $ and $ \BunGnabla $ recover the classifying space of the group $ G $ equipped with the \textit{discrete} topology.

\begin{notation}
	Let $ G $ be a topological group.
	We write $ \Gdisc \in \Grp(\Set) $ for the discrete group obtained by forgetting the topology on $ G $.
	We write $ \BGdisc \in \Spc $ for the classifying space of $ \Gdisc $.
\end{notation}

\begin{observation}\label{obs:GammalowerstarG_is_Gdisc}
	Let $ G $ be a Lie group.
	Note that the global sections
	\begin{equation*}
		\Gammalowerstar(G) = \Cinf(\ast,G)
	\end{equation*}
	recover the discrete group $ \Gdisc $.
\end{observation}

\begin{proposition}\label{lem:BGdisc}
	Let $ G $ be a Lie group.
	Then: 
	\begin{enumerate}[label=\stlabel{lem:BGdisc}, ref=\arabic*]
		\item\label{lem:BGdisc.1} The forgetful map $ \fromto{\BunGnabla}{\BunG} $ induces an equivalence 
		\begin{equation*}
			\equivto{\Gammalowerstar(\BunGnabla)}{\Gammalowerstar(\BunG)}
		\end{equation*}
		on global sections.

		\item\label{lem:BGdisc.2} There is a natural equivalence $ \Gammalowerstar(\BunG) \equivalent \BGdisc $.
	\end{enumerate}
\end{proposition}

\begin{proof}
	To prove \enumref{lem:BGdisc}{1}, note that since the Maurer--Cartan form is the \textit{unique} connection form on the trivial $ G $-bundle on the point (\Cref{obs:MC_unique_connection}), the map
	\begin{equation*}
		\fromto{\BunGnabla(\ast)}{\BunG(*)} 
	\end{equation*} 
	forgetting the connection form is an equivalence.
	For \enumref{lem:BGdisc}{2}, using the fact that $ \Gammalowerstar $ preserves limits and colimits, we compute
	\begin{align*}
		\Gammalowerstar(\BunG) &\equivalent \Gammalowerstar\real{\Bar(G)} &&\text{(\Cref{prop:Cech_nerve_trivG})} \\ 
		&\equivalent \real{\Gammalowerstar\Bar(G)} \\
		&\equivalent \real{\Bar(\Gammalowerstar(G))} \\ 
		&\equivalent \real{\Bar(\Gdisc)} &&\text{(\Cref{obs:GammalowerstarG_is_Gdisc})} \\
		&= \BGdisc \period && \qedhere 
	\end{align*} 
\end{proof}

\newpage
%!TEX root = ../diffcoh.tex

%-------------------------------------------------------------------%
%-------------------------------------------------------------------%
%  On-diagonal differential characteristic classes                  %
%-------------------------------------------------------------------%
%-------------------------------------------------------------------%

\section{On-diagonal differential characteristic classes}\label{DifferentialCharacteristicClasses}
\textit{by Arun Debray}

In \cref{ChernWeilTheory}, we constructed Chern, Pontryagin, and Euler classes of vector bundles in the
de Rham cohomology of manifolds $M$. The catalyst for this chapter is the observation that these classes are always in the image of the map $\H^*(M;\ZZ)\to\HdR^*(M)$. That is, we have the diagram
\begin{equation}
% https://q.uiver.app/?q=WzAsMyxbMCwwLCJcXEheayhNO1xcWlopIl0sWzEsMSwiXFxIZFJeayhNKSJdLFswLDIsIlxcT21lZ2FjbF5rKE0pIl0sWzAsMSwiY15cXFpaKFYpXFxtYXBzdG8gYyhWKSIsMCx7ImxhYmVsX3Bvc2l0aW9uIjo5MH1dLFsyLDEsIlAoRl9cXG5hYmxhKVxcbWFwc3RvIGMoVikiLDIseyJsYWJlbF9wb3NpdGlvbiI6OTB9XV0=
\begin{tikzcd}
	{\H^k(M;\ZZ)} \\
	& {\HdR^k(M)}, \\
	{\Omegacl^k(M)}
	\arrow["{c^\ZZ(V)\mapsto c(V)}"{pos=0.5}, from=1-1, to=2-2]
	\arrow["{P(\curvature{A})\mapsto c(V)}"'{pos=0.5}, from=3-1, to=2-2]
\end{tikzcd}
\end{equation}
which looks suspiciously like two sides of the differential cohomology hexagon. We therefore ask whether it is
possible to fill in the middle: can one choose a class $\chat\in \Hhat^*(M; \ZZ)$ whose image under the curvature
map is the Chern--Weil form, and whose image under the characteristic class map is the lift of the characteristic
class to $\ZZ$-valued cohomology?

The answer is yes, and in fact this was one of Cheeger--Simons' original applications of their theory of
differential characters \cite[\S 2]{MR827262}. In this section, we will follow the proof of
Bunke--Nikolaus--Völkl \cite[\S 5.2]{MR3462099}, who work universally on the classifying stack $\BunGnabla $ from
\Cref{ex:BunGnabla}. After that, we review our examples, constructing differential lifts of Chern,
Pontryagin, and Euler classes, and discuss how the Whitney sum formula behaves in the differential context.
Finally, we use the differential refinement of Chern--Weil theory to give a clean general description of secondary
invariants. These invariants in particular include Chern--Simons invariants, which we will use again and again in
\cref{part:applications}.

%-------------------------------------------------------------------%
%-------------------------------------------------------------------%
%  Lifting the Chern–Weil map to differential cohomology             %
%-------------------------------------------------------------------%
%-------------------------------------------------------------------%

\subsection{Lifting the Chern--Weil map to differential cohomology}

Begin with a Lie group $G$ and an invariant polynomial $P\in \Sym^n(\gdual)^G$. From $P$, the Chern--Weil machine
constructs a closed form $P(\Omega)\in \Omegacl^\bullet(\BunGnabla )$.\footnote{We end up with a form on the
universal object $\BunGnabla $ because Chern--Weil forms are natural in the connection. For more information,
see \cite[(7.21)]{FreedHopkins} and the surrounding text.}

We next need to choose an integer lift $c^\ZZ$ of $c$. There is both an existence and a uniqueness question: an
arbitrary cohomology class need not be in the lattice
\begin{equation*}
	\im(\H^k(\BG;\ZZ)\to \H^k(\BG;\RR)) \comma
\end{equation*}
and if there is torsion in $\H^k(\BG;\ZZ)$, the lift is not unique.\footnote{For example, there is torsion in
$\H^*(\BO_n; \ZZ)$ and $\H^*(\BSO_n;\ZZ)$ \cite{Bro82}.}

\begin{theorem}[{(Cheeger--Simons \cite[Theorem 2.2]{MR827262}, Bunke--Nikolaus--Völkl \cite[\S 5.2]{MR3462099})}]
\label{differential_CW_lift}
Given this data, there is a unique natural class $\chat\in\Hhat^k(\BunGnabla ; \ZZ)$ whose image under the
characteristic class map is $c^\ZZ$ and whose image under the curvature map is $P(\Omega)$.
\end{theorem}

\noindent Naturality is with respect to $G$, keeping track of the data $c^\ZZ$.

\begin{proof}
	The invariant polynomial $P$ gives us a map of sheaves of sets on $\Mfld$:
	\begin{equation}
	\label{first_sheaf_set}
		\begin{tikzcd}
			\Omega^1(-;\g) \arrow[rrr, "{\omega\mapsto \d\omega + [\omega, \omega]}", from=1-1, to=1-4] & & & \Omega^2(-;\g) \arrow[r, "P"] & \Cyc^{2p}(\Omega) \comma	
		\end{tikzcd}
	\end{equation}
	where $\Cyc$ is the sheaf of differential cycles from \cref{differential_cycles}.
	\index[terminology]{differential cycles}
	If $\yo(G)$ denotes the sheaf of groups associated to $G$ by the Yoneda embedding, then the maps
	in~\eqref{first_sheaf_set} are $\yo(G)$-equivariant, where $\Cyc^{2p}(\Omega)$ is given the trivial
	$\yo(G)$-action. Take the groupoid quotient
	\begin{equation}
		\begin{tikzcd}
			(\Omega^1(-;\g)) \modmod \yo(G) & (\Omega^2(-;\g)) \modmod \yo(G) & \Cyc^{2p}(\Omega)\modmod \yo(G) \comma
			\arrow[from=1-1, to=1-2]
			\arrow[from=1-2, to=1-3]
		\end{tikzcd}
	\end{equation}
	then take the nerve and sheafify, giving $\BunGnabla $ and $\BbulletG$ as we discussed in
	\Cref{prop:Cech_nerve_trivG,prop:Cech_nerve_trivGnabla}:
	\begin{equation}\label{BnablaBbullet}
		\BunGnabla \longrightarrow {\BbulletG}\times i(\Cyc^{2p}(\Omega))\comma
	\end{equation}
	where $i\colon\Set\to\sSet$ builds the constant simplicial set out of a set. There is an equivalence of simplicial
	sheaves
	\begin{equation}
	\label{what_is_iz2p}
		i(\Cyc^{2p}(\Omega)) \isomorphism \Omega^\infty \EM(\Cyc^{2p}(\Omega)[0])\comma
	\end{equation}
	where $\EM\colon \D(\ZZ) \to\Sp$ is the Eilenberg--MacLane functor and $[0]$ means we regard the sheaf
	$\Cyc^{2p}(\Omega)$ of abelian groups as a sheaf of complexes concentrated in degree zero.

	Take~\eqref{BnablaBbullet}, compose with the projection onto $i(\Cyc^{2p}(\Omega))$, and apply~\eqref{what_is_iz2p}
	to obtain a map
	\begin{equation*}
		\varphi_P\colon \BunGnabla \to \Omega^\infty \EM(\Cyc^{2p}(\Omega)[0]) \period
	\end{equation*}
	Let 
	\begin{equation*}
		\psi_P\colon \Sigma^\infty_+ \BunGnabla \to \EM(\Cyc^{2p}(\Omega)[0])
	\end{equation*}
	be the adjoint of $\varphi_P$ under the $(\Sigma^\infty,\Omega^\infty)$ adjunction.

	Now apply the homotopification functor $\Lhi\colon\Sh(\Mfld; \Spc)\to\Sh_\RR(\Mfld; \Spc)$ from
	\cref{homotopification_definition} to $\psi_P$. We claim this produces a map
	\begin{equation}
		\Gamma^*(\Sigma_+^\infty\BG)\overset{\chi_P}{\longrightarrow} \Gamma^*(\HRR[2p])\period
	\end{equation}
	To see this, use the identifications $\Lhi\simeq\Gammaupperstar\Gammalowersharp$
	(\cref{homotopification_definition}) and $\Gammalowersharp(E)\simeq \real{E(\Deltaalgdot)}$
	(\cref{cor:formula_for_Gammalowersharp}). The identification
	\begin{equation*}
		\Gammalowersharp(\EM(\Cyc^{2p}(\Omega)[0]))\simeq \HRR[2p]
	\end{equation*}
	is a dressed-up version of the de Rham theorem, and the equivalence
	\begin{equation*}
		\Gammalowersharp(\Sigma_+^\infty\BunGnabla ) \simeq \Sigma_+^\infty\BG
	\end{equation*}
	is \Cref{cor:formula_for_Gammalowersharp_Suspension_spectrum_BunG}.

	Next we have to identify $\chi_P$. On cohomology, the Chern--Weil construction uses $P$ to naturally assign a
	degree-$2p$ real cohomology class to a principal $G$-bundle; this soups up to a map
	$\xi_P\colon\Sigma_+^\infty\BG\to \HRR[2p]$. Looking back at the definition of $\varphi_P$, we see that
	$\Gammalowersharp(\psi_P) = \xi_P$; therefore $\chi_P = \Gammaupperstar(\xi_P)$.

	Here is where $c^\ZZ$ comes in. It is data of a lift
	\begin{equation}
	\begin{gathered}
	\begin{tikzcd}
		& \HZZ[2p] \\
		\Sigma_+^\infty\BG & \HRR[2p],
		\arrow[from=1-2, to=2-2]
		\arrow["\xi_P"', from=2-1, to=2-2]
		\arrow["c^\ZZ", from=2-1, to=1-2]
	\end{tikzcd}
	\end{gathered}
	\end{equation}
	giving us a diagram
	\begin{equation*}
		\begin{tikzcd}[column sep=3em]
			\Sigma_+^\infty\BunGnabla \arrow[rr, bend left, "\psi_P"] \arrow[d] & & \EM(\Cyc^{2p}(\Omega)[0]) \arrow[d] \\
			\Gamma^*(\Sigma_+^\infty\BG) \arrow[rr, bend right, "\Gammaupperstar(\xi_P)"'] \arrow[r, "\Gammaupperstar(c^\ZZ)"] & \Gammaupperstar(\HZZ[2p]) \arrow[r] & \Gammaupperstar(\HRR[2p]) \period
		\end{tikzcd}
	\end{equation*}
	The vertical arrows are both of the form $F\to\Lhi(F)$, and are unit maps for the adjunction $(\Lhi,
	\text{inclusion})$ from \cref{nul:Gammaadjunctions}.

	The map from the upper left to the lower right factors through the pullback
	\begin{equation*}
		\begin{tikzcd}[column sep=3em]
			\Sigma_+^\infty\BunGnabla \arrow[rr, bend left, "\psi_P"] \arrow[r, "\chat"'] \arrow[d] & \HZZhat(2p) \arrow[r] \arrow[d] \arrow[dr, phantom, "\lrcorner"{very near start}] & \EM(\Cyc^{2p}(\Omega)[0]) \arrow[d] \\
			\Gamma^*(\Sigma_+^\infty\BG) \arrow[rr, bend right, "\Gammaupperstar(\xi_P)"'] \arrow[r, "\Gammaupperstar(c^\ZZ)"] & \Gammaupperstar(\HZZ[2p]) \arrow[r] & \Gammaupperstar(\HRR[2p]) \period
		\end{tikzcd}
	\end{equation*}
	and $\chat$ is the desired differential refinement.
\end{proof}

\begin{remark}
	Cheeger--Simons' original proof did not use this language: they did not have $\BunGnabla $ available. Instead, they
	use \emph{$n$-classifying spaces} $\beta_\nabla^{(n)} G$. These are spaces such that all connections on principal
	$G$-bundles $P\to M$ pull back from $\beta_\nabla^{(n)} G$, provided $\dim M < n$, and the pullback need not be
	unique. Narasimhan--Ramanan \cite{NR61, NR63} proved $n$-classifying spaces exist for all $n$ and $G$, provided
	$\uppi_0(G)$ is finite.
\end{remark}

\begin{example}[(differential Chern classes)]
	\label{differential_Chern}
	\index[terminology]{Chern class!on-diagonal differential refinement}
	\index[terminology]{differential Chern class}
	Borel \cite[\S 29]{Bor53} shows that
	\begin{equation*}
		\H^*(\BU_n;\ZZ)\cong \ZZ[c_1,\dotsc, c_n] \comma
	\end{equation*}
	so integer lifts are unique, and using Grothendieck's axioms, one can show that the images of these Chern classes in de Rham
	cohomology are equal to the Chern classes we constructed in \S\ref{chern_const}. Therefore we obtain
	\emph{on-diagonal differential Chern classes} $\chat_k(P, \conn)\in\Hhat^{2k}(M; \ZZ)$ associated to principal
	$\Uup_n$-bundles $P\to M$ with connection $\conn$. 
	See \Cref{rmk-OnOffDiagonal} for a note on terminology.

	Several authors construct differential Chern classes by other methods, including
	Brylinski--McLaughlin \cite{BM96},
	Berthomieu \cite{Ber10},
	Bunke \cite{Bun10,Bunke:notes}, and
	Ho \cite{Ho15}. Schreiber \cite{Urs} constructs $\chat_1$.
	\index[notation]{ccheckk@$\chat_k$}
\end{example}

\begin{example}[(differential Pontryagin classes)]
\label{differential_Pontryagin}
\index[terminology]{Pontryagin class!on-diagonal differential refinement}
\index[terminology]{differential Pontryagin class}
	Brown \cite[Theorem 1.6]{Bro82} shows there is torsion in $\H^*(\BO_n; \ZZ)$, so choosing $p_k^\ZZ$ is not
	automatic. Let $c\colon\BO_n\to\BU_n$ be the complexification map, and for a principal
	$\Or_n$-bundle $P\to M$ define
	\begin{equation}
		p_k(P) \colonequals (-1)^k c_{2k}(c(P))\in \H^{2k}(M; \ZZ)\period
	\end{equation}
	The images of these classes in de Rham cohomology are equal to the Pontryagin classes we defined in
	\S\ref{pontrjagin},
	so \cref{differential_CW_lift} produces for us
	\emph{on-diagonal differential Pontryagin classes} $\phat_k(P,  \conn)\in\Hhat^{4k}(M;\ZZ)$ associated to
	principal $\Or_n$-bundles with connection $\conn$.

	Brylinski--McLaughlin \cite{BM96} and Grady--Sati \cite[Proposition 3.6]{GS21} construct $\phat_k$ in a
	different way.
	\index[notation]{pcheck@$\phat_k$}
\end{example}

\begin{example}[(differential Euler classes)]
	\index[terminology]{Euler class!on-diagonal differential refinement}
	\index[terminology]{differential Euler class}
	Brown \cite[Theorem 1.6]{Bro82} shows that there is also torsion in $\H^*(\BSO_n; \ZZ)$, so we must
	choose a lift $e^\ZZ$ of the Euler class we constructed in \cref{euler_const}. There, we defined $e$ only for $n$
	even; for odd $n$, we set $e\colonequals 0$.

	Let $V\to\BSO_n$ denote the tautological bundle. Since $V$ is oriented, it has a $\ZZ$-cohomology
	Thom class $\tau(E)\in \widetilde{\H}{}^n(V, V\setminus 0; \ZZ)$. We let $e^\ZZ$ be the pullback of $\tau(E)$ by
	the zero section of $V$. The image of this class is $e$, so the class defined by the Pfaffian when $n$ is even, and
	$0$ when $n$ is odd. For all $n$, however, $e^\ZZ\ne 0$; it is $2$-torsion when $n$ is odd.

	Therefore we obtain a \emph{on-diagonal differential Euler class} $\ehat(P, \conn)\in\Hhat^n(M; \ZZ)$
	associated to a principal $\SO_n$-bundle with connection $\conn$, and it can be nonzero for all $n$, not just
	even $n$.

	Brylinski--McLaughlin \cite{BM96} and Bunke \cite[Example 3.85]{Bunke:notes} construct $\ehat$ in a different way.
\end{example}

\begin{remark}[(from principal bundles to vector bundles: an important nuance)]
\label{metric_differential_cc}
	We would like to use the characteristic classes we just constructed to define differential lifts of characteristic
	classes of vector bundles with connection. The way this usually works for characteristic classes is that a vector
	bundle has an associated principal $G$-bundle, and we consider characteristic classes for $G$. For example, a
	rank-$n$ complex vector bundle has a principal $\GL_n( \CC)$-bundle of frames. The maximal compact of $\GL_n( \CC)$
	is $\Uup_n$, so inclusion $\Uup_n\to\GL_n( \CC)$ induces a homotopy equivalence of classifying spaces, which
	means characteristic classes of principal $\Uup_n$-bundles give you characteristic classes of principal $\GL_n(
	\CC)$-bundles give you characteristic classes of complex vector bundles. Both of these steps are necessary: the
	Chern--Weil map is only guaranteed to be an isomorphism for compact groups, and without additional structure such
	as a metric, the structure group of a vector bundle is noncompact.

	In differential cohomology, this becomes a stumbling block: homotopy equivalences do not always induce isomorphisms
	on differential cohomology, so what we learn about principal $\Uup_n$-bundles does not necessarily help us with
	complex vector bundles. Therefore a priori, the differential characteristic classes we defined above only make
	sense for vector bundles with a metric and a compatible connection as in \S\ref{subsubsec:connections_for_vbs},
	because these correspond to connections on principal $\Uup_n$-bundles, rather than principal $\GL_n(
	\CC)$-bundles.
	\index[terminology]{connection!compatibility with the metric}

	In addition to complex vector bundles and $\Uup_n$ versus $\GL_n( \CC)$, which is about differential Chern
	classes, there are two more cases to worry about.
	\begin{enumerate}
		\item Real vector bundles and $\Or_n$ and $\GL_n( \RR)$, and differential Pontryagin classes.
		
		\item Oriented real vector bundles, $\SO_n$, and $\GL_n( \RR)_0$ (i.e.\ the connected component of $\GL_n(
		\RR)$ containing the identity), for the differential Euler class.
	\end{enumerate}
	First, Chern classes. For $\GL_n( \CC)$ the Chern--Weil map is not an isomorphism, but it is surjective \cites[\S
	4]{MR827262}[\S 11.8.1]{Pro07}, so differential Chern classes can be defined in the absence of a metric.

	Next, Pontryagin classes. The construction in \cref{differential_Pontryagin} implies differential Pontryagin
	classes of $V,  \conn$ are equal to differential Chern classes of $V\otimes\CC$ with connection induced from
	$\conn$, so differential Pontryagin classes can be defined in the absence of a metric.

	But Euler classes are different! If $A\in\GL_n( \RR)$ and $X\in\so_n$, then
	\begin{equation}
		\mathrm{pf}(AXA^{-1}) = \det(A)\mathrm{pf}(X)\comma
	\end{equation}
	so the Pfaffian is not $\GL_n( \RR)_0$-invariant. Therefore the differential Euler class requires an oriented
	vector bundle, a Euclidean metric, and a compatible connection.
	\index[terminology]{differential Euler class!metric-compatibility requirement}
\end{remark}

We will use these classes in a few different ways in \cref{applications_part}, including obstructing conformal
immersions in \cref{conformal_immersions} and constructing non-topological invertible field theories in
\cref{invertible_field_theories}. Cheeger--Simons \cite{Cheeger-Simons} discuss some additional applications,
including characteristic classes associated to foliations and a geometric refinement of the Atiyah--Singer index
theorem. There are also differential refinements of the Todd genus, $\Ahat$-genus \cite[Definition
3.9]{GS21}, and so forth.
\index[terminology]{Todd genus!differential refinement}
\index[terminology]{Ahat genus@$\Ahat$-genus!differential refinement}

%-------------------------------------------------------------------%
%-------------------------------------------------------------------%
%  The Whitney sum formula for differential characteristic classes  %
%-------------------------------------------------------------------%
%-------------------------------------------------------------------%

\subsection{The Whitney sum formula for on-diagonal differential characteristic classes}

The Whitney sum formula expresses the Chern, Pontryagin, and Euler classes of a direct sum $E\oplus F$ of vector
bundles in terms of the respective characteristic classes of $E$ and of $F$.
\index[terminology]{Whitney sum formula!in ordinary cohomology}
Let $c \colonequals 1 + c_1 + c_2 + \cdots$ denote the total Chern class\index[terminology]{total Chern class} and $p \colonequals 1 +
p_1 + p_2 + \cdots$ denote the total Pontryagin class.\footnote{Though these appear to be infinite sums, they are
finite when evaluated on any vector bundle, because $c_k(E)$ and $p_k(E)$ vanish when $k >
\mathrm{rank}(E)$.}\index[terminology]{total Pontryagin class} For complex vector bundles $E,F\to X$, we have
\begin{subequations}\label{ordinary_Whitney_sum}
\begin{equation}
	\begin{aligned}
		c(E\oplus F) &= c(E)c(F)\\
		c_k(E\oplus F) &= \sum_{i+j=k} c_i(E) c_j(F) \period
	\end{aligned}
\end{equation}
For oriented real vector bundles $E,F\to X$,
\begin{equation}
	e(E\oplus F) = e(E)e(F) \period
\end{equation}
Both of these equations take place in the ring $\H^*(X;\ZZ)$. However, for Pontryagin classes, the corresponding
formula only holds modulo $2$-torsion. That is, in the ring $\H^*(X;\ZZ[1/2])$,
\begin{equation}
\begin{aligned}
	p(E\oplus F) &= p(E)p(F)\\
	p_k(E\oplus F) &= \sum_{i+j=k} p_i(E) p_j(F) \period
\end{aligned}
\end{equation}
\end{subequations}
The formula for the Pontryagin classes of a direct sum with $\ZZ$ coefficients is known by work of
Thomas \cite{Tho62} and Brown \cite[Theorem 1.6]{Bro82}, but it is a little more complicated.

On to differential cohomology. Given vector bundles with connection $(E, A^E)$ and $(F, A^F)$ over a
space $X$, the direct sum $E\oplus F$ has an induced connection $A^E\oplus A^F$. One can prove the Whitney
sum formulas~\eqref{ordinary_Whitney_sum} by studying the effect of the maps 
\[\mathrm B(\GL_{n_1}(\CC)\times\GL_{n_2}(
\CC))\to \mathrm{BGL}_{n_1 + n_2}(\CC)\]
(resp.\ $\BSO_{n_i}$, $\mathrm{BGL}_{n_i}(\RR)$) on cohomology. Naturality of
\cref{differential_CW_lift} then implies
\begin{subequations}
\label{differential_Whitney}
\begin{align}
	\chat(E\oplus F,  \conn^E\oplus \conn^F) &= \chat(E,  \conn^E)\chat(F,  \conn^F)\\
	\label{differential_Whitney_Euler}
	\ehat(E\oplus F,  \conn^E\oplus \conn^F) &= \ehat(E,  \conn^E)\ehat(F,  \conn^F)\\
	\phat(E\oplus F,  \conn^E\oplus \conn^F) &= \phat(E,  \conn^E)\phat(F,  \conn^F),
	\label{differential_Whitney_Pontryagin}
\end{align}
\end{subequations}
where $E$ and $F$ are complex or oriented where needed. For~\eqref{differential_Whitney_Euler} we must assume $E$
and $F$ come with Euclidean metrics which $\conn^E$ and $\conn^F$ are compatible with, because of
\cref{metric_differential_cc}; and as usual~\eqref{differential_Whitney_Pontryagin} takes
place in $\Hhat^*(X;\ZZ[1/2])$.\index[terminology]{Whitney sum formula!in differential cohomology}

The formulas~\eqref{differential_Whitney} are less useful than they might seem: in some places you might want to
use it, the connection you care about on $E\oplus F$ is not a direct sum connection. This happens, for example, in
the proof of \cref{differential_conformal_immersion} in \cref{applications_part}. Fortunately, the differential
Whitney sum formula is true in more generality.
\begin{definition}
\index[terminology]{compatible connection}
\label{compatible_connection}
Choose connections $\conn^E$ on $E$, $\conn^F$ on $F$, and $\overline{\conn}$ on $E\oplus F$. The projections
$E\oplus F\rightrightarrows E, F$ induce connections $\overline{\conn}{}^E$, resp.\ $\overline{\conn}{}^F$ on $E$,
resp.\ $F$ from $\overline A$. 
Let $\curvature{\overline{\conn}}\in\Omega_X^2(\End(E\oplus F))$ be the curvature of
$\overline{\conn}$. We say $\conn$ is \textit{compatible} with $\conn^E\oplus\conn^F$ if
\begin{enumerate}
	\item $\overline{\conn}{}^E =  \conn^E$ and $\overline{\conn}{}^F =  \conn^F$, and
	\item given vector fields $v,w$ on $X$, $\curvature{\overline{\conn}}(v, w)\in\Gamma(\End(E\oplus F))$ is block
	diagonal.
\end{enumerate}
\end{definition}
There are two notions of compatibility floating around: compatibility with a metric, and compatibility with the
direct-sum connection. They are different.

\begin{theorem}[{(Cheeger--Simons \cite[Theorem 4.7]{Cheeger-Simons})}]
	\hfill
	\begin{enumerate}
		\item If $\overline{\conn}$ is compatible with $\conn^E\oplus\conn^F$, then
		\begin{equation*}
			\phat(E\oplus F, \overline{\conn}) = \phat(E\oplus F,  \conn^E\oplus \conn^F) \period
		\end{equation*}
	
		\item If $E$ and $F$ are oriented and Euclidean, then
		\begin{equation*}
			\ehat(E\oplus F,\overline{\conn}) = \ehat(E\oplus F,  \conn^E\oplus \conn^F) \period 
		\end{equation*} 

		\item If $E$ and $F$ are complex, then
		\begin{equation*}
			\chat(E\oplus F,\overline{\conn}) = \chat(E\oplus F,  \conn^E\oplus \conn^F) \period
		\end{equation*} 
	\end{enumerate}

	Therefore analogues of~\eqref{differential_Whitney} hold with $\overline A$ in place of $ \conn^E\oplus \conn^F$.
\end{theorem}

\noindent The proof uses a variation formula for the Chern--Simons form similar to \cref{variation_Chern_Simons}.
\index[terminology]{Chern--Simons form!variation formula}

Recall that the Whitney sum formula can be used to show that the Euler class obstructs the existence of
a section of an oriented vector bundle. In the same way, the differential Euler class obstructs flat sections.

\begin{lem}
	Let $V\to M$ be an oriented Euclidean vector bundle with compatible connection $\conn$ admitting a flat section.
	Then $\ehat(V, \conn) = 0$.
\end{lem}

\begin{proof}
	The flat section splits $V = V'\oplus\underline{\RR}$ such that $\conn$ is compatible with the direct sum
	connection, where $\underline{\RR}$ carries the standard connection $\d$.  Because $\ehat(\underline{\RR}, \d) =
	0$, the Whitney sum formula finishes the proof for us.
\end{proof}

%-------------------------------------------------------------------%
%-------------------------------------------------------------------%
%  Secondary invariants and Chern–Simons forms                      %
%-------------------------------------------------------------------%
%-------------------------------------------------------------------%

\subsection{Secondary invariants and Chern--Simons forms}
\label{secondary_invariants}
% first, the motivating idea
Degree-$n$ characteristic classes provide invariants of closed, oriented $n$-manifolds by integration, and these
invariants provide useful topological information: integrating the Euler class produces the Euler characteristic,
and integrating products of Pontryagin classes produces oriented bordism invariants. In this section we discuss
the analogous invariants defined by integrating on-diagonal differential characteristic classes; since the
differential cohomology of a point is not concentrated in degree zero, we do not have to stick to $n$-manifolds.

Let $G$ be a compact Lie group and $c^\ZZ\in\H^n(\BG;\ZZ)$. \Cref{differential_CW_lift} gives us an on-diagonal differential lift
$\chat\in\Hhat^n(\BunGnabla ;\ZZ)$ of $c^\ZZ$. Let $M$ be a closed, oriented $(n-1)$-manifold, and let $P \to M$ be a
principal $G$-bundle with connection $\conn$. In \cref{FiberIntegration}, we constructed an integration map on
differential cohomology. Integration has degree $-(n-1)$, so if $\alpha_c(P,  \conn)$ denotes the integral of
$\chat(P,  A)$, then $\alpha_c(P,  \conn)$ is an element of $\RR/\ZZ$:
\begin{equation}
\begin{aligned}
	\int_M\colon \Hhat^n(M; \ZZ) &\longrightarrow \Hhat^1(\pt; \ZZ)\cong\RR/\ZZ\\
	\chat(P,  \conn) &\longmapsto \alpha_c(P,  \conn) \period
\end{aligned}
\end{equation}
The quantity $\alpha_c(P, \conn)$, as an $\RR/\ZZ$-valued invariant of principal bundles with connection, is
called the \emph{secondary invariant} associated to $c$. In this context, the $\ZZ$-valued purely topological
invariant $\int_M c^\ZZ(P)$ on $n$-manifolds is called the \emph{primary invariant}.
\index[terminology]{secondary invariant}
\index[terminology]{primary invariant}

In examples, secondary invariants tend to be very geometric, despite our general abstract definition.

\begin{example}[(holonomy of a connection on a principal $\Uup_1$-bundle)]
	Let $P\to M$ be a principal $\Uup_{1}$-bundle with connection $\conn$ and consider the differential
	first Chern class $\chat_1(P, \conn)$, built from the curvature form of $\conn$. Given an embedded, oriented loop
	$i\colon \Circ\hookrightarrow M$, we can pull back $\chat_1(P, \conn)$ to $\Circ$ and integrate, defining an element
	of $\RR/\ZZ$. Cheeger--Simons \cite[Example 1.5]{MR827262} show that this $\RR/\ZZ$-valued quantity is the log of
	the holonomy of $P$ around $\Circ$. That is, holonomy is the secondary invariant associated to the first Chern class
	or the curvature for principal $\Uup_1$-bundles.
	\index[terminology]{holonomy!as a secondary invariant}
\end{example}

\begin{example}[(Chern--Simons invariants)]
	\label{secondary_chern_simons}
	Chern--Simons invariants are important examples of secondary invariants: they will appear several times in several
	different ways in \cref{applications_part}. In some settings, any secondary invariant constructed via Chern--Weil
	theory is called a Chern--Simons invariant, but by far the most commonly considered example is in dimension $3$.

	Choose a compact Lie group $G$ and an element $\lambda\in\H^4(\BG;\ZZ)$, which we call the
	\emph{level}. Given a closed $3$-manifold $Y$, a principal $G$-bundle $P\to Y$, and a connection $\conn$ on $P$,
	the \emph{Chern--Simons invariant} $\CS_\lambda(P,\conn)\in\RR/\ZZ$ \cite{cs} is defined to be value of
	the secondary invariant associated to $\lambda$ on $(P,
	\conn)$.\index[terminology]{level}\index[terminology]{Chern--Simons invariant}

	The standard construction of $\CS_\lambda(P, \conn)$, which is the construction Chern--Simons gave, is more
	geometric. We will discuss this in \cref{config_spaces}. The approach here, using differential cohomology, is due
	to Cheeger--Simons \cite{MR827262}.
\end{example}

Chern \cite{Che44} defines a differential form in a sphere bundle related to the secondary invariant built from the
Euler class.

% then a digression about eta-invariants
\begin{remark}[(secondary invariants and differential generalized cohomology)]
We can try to run the same story with a generalized cohomology theory $E$. To do so, we need a differential
refinement $\Ehat$ of $E$, an integration map for $\Ehat$-cohomology (possibly on manifolds with some additional
structure) and an on-diagonal differential characteristic class $\chat\in \Ehat^*(\BunGnabla )$. Together these data are a lot to
ask for, but everything goes through in $\KU$-theory, for example.

Definitions of differential refinements of $\KU$ and $\KO$ were first sketched by Freed \cite[Examples 1.12 and
1.13]{Fre00} and Freed--Hopkins \cite{FH00}. Hopkins--Singer \cite[\S 4.4]{HopkinsSinger} first constructed
differential $\KU$-theory, and Grady--Sati \cite{GS21} first systematically study differential
$\KO$-theory. There are differential lifts of the Atiyah--Bott--Shapiro integration maps in $\KU$- and $\KO$-theory
on closed spin\textsuperscript{$c$}, resp.\ spin manifolds.

We can therefore study secondary invariants for $\KU$- and $\KO$-theories. The final piece of data we need is a
differential characteristic class, and we choose $1\in\KU^0(X)$ or $\KO^0(X)$. The primary invariant associated
with this data on a spin or spin\textsuperscript{$c$} manifold admits a geometric interpretation as the index of
the spinor Dirac operator \cite{AS68}. The secondary invariant has a related description \cite{Lot94}, as the
\emph{$\eta$-invariant} of the Dirac operator, defined and studied by Atiyah--Patodi--Singer \cite{APS1, APS2,
APS3}. See Freed~\cite[\S 7.4]{Fre21} for more information.

There are several additional models for differential $\KU$-theory constructed by
Klonoff \cite{Klo08},
Bunke--Schick \cite[\S 2]{BS09},
Simons--Sullivan \cite{MR2732065},
Bunke--Nikolaus--Völkl \cite[\S 6]{MR3462099},
Schlegel \cite[\S 4.2]{Sch13},
Tradler--Wilson--Zenalian \cite{TWZ13, TWZ16},
Hekmati--Murray--Schlegel--Vozzo \cite{HMSV15},
Park \cite{Par17},
Gorokhovski--Lott \cite{GL18},
Schlarmann \cite{Sch19},
Park--Parzygnat--Redden--Stoffel \cite{PPRS21},
Cushman \cite{Cus21},
and Gomi--Yamashita \cite{GY21}; Cushman and Gomi-Yamashita also construct models for differential $\KO$-theory.

%Several of these teams of researchers also show how to lift the Atiyah--Bott--Shapiro
%integration map on \Ktheory, defined on closed spin\textsuperscript{$c$} manifolds, to differential \Ktheory.
See Bunke--Schick \cite{BunksSchick} for a survey.
\end{remark}

\newpage
%!TEX root = ../diffcoh.tex

%-------------------------------------------------------------------%
%-------------------------------------------------------------------%
%  Chern–Weil forms after Freed–Hopkins                             %
%-------------------------------------------------------------------%
%-------------------------------------------------------------------%

\section{Chern--Weil forms after Freed--Hopkins}\label{WorkofFreedHopkins}

\textit{by Dexter Chua}

Let $ G $ be a Lie group.
The main theorem of the Freed--Hopkins paper \emph{Chern--Weil forms and abstract homotopy theory} \cite{FreedHopkins} is that Chern--Weil forms are the only natural way to get a differential form from a principal $G$-bundle.
That is, Freed and Hopkins computed the de Rham complex of the sheaf $ \BunGnabla $ in terms of Chern--Weil forms.

Theorems along these lines are of interest historically. It is an important ingredient in the heat kernel proof of the Atiyah--Singer index theorem. 
Essentially, the idea of the proof is to use the heat equation to show that there is \emph{some} formula for the index of a vector bundle in terms of the derivatives of the metric, and then by invariant theory, this must be given by the Chern--Weil forms we know and love. One then computes this for sufficiently many examples to figure out exactly which characteristic class it is, as Hirzebruch originally did for his signature formula.

%-------------------------------------------------------------------%
%-------------------------------------------------------------------%
%  The statement                                                    %
%-------------------------------------------------------------------%
%-------------------------------------------------------------------%

%reference 7.2e for B_\nablaG definition

\subsection{The statement}\label{section:statement}

To state the theorem, we work in the category $\Shv(\Mfld;\Spc)$. 

\begin{recollection}\label{rec:Omega_BunG_and_BunGnabla}
  There are a number of important sheave on $ \Mfld $ that we have already encountered:
  \begin{enumerate}[label=\stlabel{rec:Omega_BunG_and_BunGnabla}, ref=\arabic*]
    \item Any $M \in \Mfld$ defines a representable (discrete) sheaf, which we denote by $M$ again (\Cref{ex:representable_sheaf}).

    \item For $p \geq 0$, we have a discrete sheaf
    \begin{equation*}
      \Omega^p \in \Shv(\Mfld;\Spc)\period
    \end{equation*}
    This is in fact a sheaf of vector spaces (\Cref{ex:Omegai_is_a_sheaf}).
    Moreover, there are linear natural transformations ${\d\colon\Omega^p \to \Omega^{p + 1}}$. 
    Thus, we get a sheaf of chain complexes $\Omega^{\bullet}$ whose value on a manifold $ M $ is the de Rham complex of $ M $.

    \item Fix $G$ a Lie group. 
    We write $ \BunGnabla  \colon \fromto{\Mfldop}{\Spc_{\leq 1}} $\index[notation]{Groupoid of Principal G-bundles@$\BunGnabla $} for the sheaf sending a manifold $ M $ to be the groupoid of principal $G$-bundles on $M$ with connection and isomorphisms (\Cref{ex:BunGnabla}).
  \end{enumerate}
\end{recollection}

\begin{notation}
    We also write
    \begin{equation*}
      \Omegabullet \colon \fromto{\PSh(\Mfld;\Spc)^{\op}}{\D(\RR)}
    \end{equation*}
    for the right Kan extension of $ \Omegabullet \colon \fromto{\Mfldop}{\D(\RR)} $ along the Yoneda embedding.
    Thus, for any presheaf $\Fcal \in \PSh(\Mfld;\Spc) $, we can think of
    \begin{equation*}
      \Omegabullet (\Fcal) \colonequals \lim_{M \to \Fcal} \Omega^{\bullet}(M)
    \end{equation*}
    as the de Rham complex of $\Fcal$.\index[terminology]{de Rham complex!of a sheaf on manifolds}
\end{notation}

The main theorem is:

\begin{theorem}[{\cite[Theorem 7.20]{FreedHopkins}}]
  Let $ G $ be a Lie group.
  The Chern--Weil homomorphism induces an isomorphism:
  \begin{equation*}
    (\Sym^\bullet\gdual)^G \isomorphism \Omegabullet (\BunGnabla ) \period
  \end{equation*}
  Here $ (\Sym^\bullet\gdual)^G $ is regarded as a complex with zero differential and $ (\Sym^i \gdual)^G $ is in degree $ 2i $.
  In particular, the de Rham differential on $ \BunGnabla $ is zero.
\end{theorem}

\noindent This implies that the Chern--Weil construction is the only natural way of obtaining differential forms from a principal $G$-bundle.

To prove the theorem, we use the universal principal $G$-bundle 
\begin{equation*}
  \trivGnabla \colon \surjto{\Omega^1(-;\g)}{\BunGnabla}
\end{equation*}
introduced in \Cref{cons:trivGnabla}, as well as the presentation of $ \BunGnabla $ as the quotient
\begin{equation*}
    \BunGnabla \equivalent \Omega^1(-;\g) \modmod G
\end{equation*}
of $ \Omega^1(-;\g) $ by the adjoint action of $ G $ (\Cref{prop:Cech_nerve_trivGnabla}).

Our proof naturally breaks into two steps. 
First, we compute $\Omegabullet(\Omega^1(-;\g))$, and then we need to know how to compute $\Omegabullet(\Fcal\modmod G)$ from $\Omegabullet (\Fcal)$ for any discrete sheaf $\Fcal$.

We first do the second part.

\begin{lemma}\label{lemma:omega-mod-g}
  Let $\Fcal \in \Shv(\Mfld;\Spc)$ be a discrete sheaf with a $G$-action $\alpha \colon G \times \Fcal \to \Fcal$. 
  Then $\Omegabullet(\Fcal \modmod G)$ is the subcomplex of $\Omegabullet(\Fcal)$ consisting of the $\omega$ satisfying the following properites:
  \begin{enumerate}[label=\stlabel{lemma:omega-mod-g}, ref=\arabic*]
    \item\label{lemma:omega-mod-g.1} For all $g \in G$, we have $\alpha^* \omega|_{\{g\} \times \Fcal} = \omega$.

    \item\label{lemma:omega-mod-g.2} For all $\xi \in\g$, we have $\iota_\xi \omega = 0$.
  \end{enumerate}
\end{lemma}

\noindent The first condition says $\omega$ should be $G$-invariant, and the second condition says $\omega$ is suitably ``horizontal''.

\begin{remark}\label{remark:iota}
  Let us explain what we mean by $\iota_\xi \omega$. 
  In general, for $M$ a manifold and $X$ a vector field on $M$, for each manifold $ N $, contraction with $ X $ on $ M $ defines a map
  \begin{equation*}
    \iota_X\colon \Omega^p(M \times N) \to \Omega^{p - 1}(M \times N) \period
  \end{equation*}
  Then by left Kan extension, this induces a map
  \begin{equation*}
    \iota_X\colon \Omega^p(M \times \Fcal) \to \Omega^{p - 1}(M \times \Fcal)
  \end{equation*}
  for all $\Fcal \in \Shv(\Mfld;\Spc)$. 

  Now if $\Fcal$ has a $G$-action and $\xi \in\g$, then $\xi$ induces an invariant vector field on $G$, which we also call $\xi$. We then define $\iota_\xi\colon \Omega^p(\Fcal) \to \Omega^{p - 1}(\Fcal)$ by the following composition
    \begin{equation*}
      \begin{tikzcd}
        \Omega^p(\Fcal) \ar[r, "\alpha^*"] & \Omega^p(G \times \Fcal) \ar[r, "\iota_\xi"] & \Omega^{p - 1}(G \times \Fcal) \ar[r] & \Omega^{p - 1}(\{e\} \times \Fcal) = \Omega^{p - 1}(\Fcal) \comma
      \end{tikzcd}
    \end{equation*}
  where the last map is induced by the inclusion.

  This gives us a very explicit method to compute the natural transformation $\iota_\xi \omega$ for $\omega \in \Omega^p(\Fcal)$ and $\xi \in\g$. Given a test manifold $M$ and $\phi \in \Fcal(M)$, which we think of as a natural transformation $\phi \colon M \to \Fcal$, we form the composite
    \begin{equation*}
      \begin{tikzcd}
        G \times M \ar[r, "1 \times \phi"] & G \times \Fcal \ar[r, "\alpha"] & \Fcal \ar[r, "\omega"] & \Omega^p
      \end{tikzcd}
    \end{equation*}
  This defines a differential form $\eta \in \Omega^p(G \times M)$. Then we have
  \begin{equation*}
    (\iota_\xi \omega)_M(\phi) = \iota_{\xi} \eta|_{\{e\} \times M}\period
  \end{equation*}
\end{remark}

\begin{proof}
  We have
  \begin{equation*}
    \Omega^p (\Fcal \modmod G) = \Omega^p(|(\Fcal\modmod G)_\bullet|) = \Tot(\Omega^p((\Fcal\modmod G)_\bullet))\period
  \end{equation*}
  Since $(\Fcal\modmod G)_\bullet$ is a simplicial discrete sheaf, its totalization can be computed by
    \begin{equation*}
      \Omega^p(\Fcal\modmod G) = \ker\, \left(
      \begin{tikzcd}[sep=3.5em]
        \Omega^p(\Fcal) \arrow[r, "\pr^* - \alpha^*"] & \Omega^p(G \times \Fcal)
      \end{tikzcd}
      \right) \comma
    \end{equation*}
  where $\pr \colon G \times \Fcal \to \Fcal$ is the projection.

  To prove the lemma, we have to show that $\pr^* \omega = \alpha^* \omega$ if and only if the conditions in the lemma are satisfied. This follows from the more general claim below with $\eta = \alpha^* \omega - \pr^* \omega$.

  \begin{claim}\label{claim:when_forms_are_zero}
    Let $M$ be a manifold and $\Fcal$ a sheaf. 
    Then $\eta \in \Omega^p(M \times \Fcal)$ is zero if and only if:
    \begin{enumerate}[label=\stlabel{claim:when_forms_are_zero}, ref=\arabic*]
      \item\label{claim:when_forms_are_zero.1} For all $x \in M$, we have $\eta|_{\{x\} \times \Fcal} = 0$. 
      
      \item\label{claim:when_forms_are_zero.2} For any vector field $X$ on $M$, we have $\iota_X \eta = 0$.
    \end{enumerate}
  \end{claim}

  The conditions \enumref{lemma:omega-mod-g}{1} and \enumref{claim:when_forms_are_zero}{1} match up exactly. 
  Unwrapping the definition of $\iota_\xi$ and noting that $\iota_X \pr^* \omega = 0$ always, the only difference between \enumref{lemma:omega-mod-g}{2} and \enumref{claim:when_forms_are_zero}{2} is that in \enumref{lemma:omega-mod-g}{2}, we only test on invariant vector fields on $G$, instead of all vector fields, and we only check the result is zero after restricting to a fiber $\{e\} \times \Fcal$. 
  The former is not an issue because the condition $\Cinf(G)$-linear and the invariant vector fields span as a $\Cinf(G)$-module. 
  The latter also doesn't matter because we have assumed that $\alpha^* \omega$ is invariant. 

  To prove the claim, if $\Fcal$ were a manifold, this is automatic, since the first condition says $\eta$ vanishes on vectors in the $N$ direction while the second says it vanishes on vectors in the $M$ direction.

  If $\Fcal$ were an arbitrary sheaf, we know $\eta$ is zero when pulled back along any map 
  \begin{equation*}(1 \times \phi) \colon M \times N \to M \times \Fcal\end{equation*}
   where $N$ is a manifold, by naturality of the conditions. But since $M \times \Fcal$ is a colimit of such maps, $\eta$ must already be zero on $M \times \Fcal$.
\end{proof}

Now it remains to describe $ \Omegabullet(\Omega^1(-;\g))$. More generally, for any vector space $V$, we can calculate $\Omegabullet(\Omega^1(-;V))$. We first state the result in the special case where $V = \RR $.

\begin{theorem}
  For each $ p \geq 0 $ there is an equivalence 
  \begin{equation*}
    \Omega^p(\Omega^1) \cong \RR \period
  \end{equation*}
  For $p = 2q$, it sends $\omega$ to $(\d\omega)^q$. 
  For $p = 2q + 1$, it sends $\omega$ to $\omega \wedge (\d \omega)^q$.
\end{theorem}

The general case is no harder to prove, and the result is described in terms of the Koszul complex.
\begin{definition}
  Let $V$ be a vector space. The \emph{Koszul complex} $\Kos^\bullet V$ is a differential graded algebra whose underlying algebra is
  \begin{equation*}
    \Kos^\bullet V = \exterior^\bullet(V) \tensor \Sym^\bullet(V) \period  \index[notation]{Kos@$\Kos^\bullet V$} \index[terminology]{Koszul complex}
  \end{equation*}
  For $v \in V$, we write $v$ for the corresponding element in $\exterior^1 V$, and $\tilde{v}$ for the corresponding element in $\Sym^1 V$. We set $|v| = 1$ and $|\tilde{v}| = 2$. 
  The differential is then defined by
  \begin{equation*}
    \d(v) = \tilde{v} \andeq \d(\tilde{v}) = 0\period
  \end{equation*}
\end{definition}

\begin{theorem}\label{thm:main}
  For any vector space $V$, we have an isomorphism of differential graded algebras
  \begin{equation*}
    \eta\colon \Kos^\bullet(\Vdual) \isomorphism \Omegabullet (\Omega^1(-;V))\period
  \end{equation*}
  In particular,
  \begin{equation*}
    \Omegabullet (\Omega^1(-;\g)) = \Kos^\bullet(\gdual) \period
  \end{equation*}

  Explicitly, for $\ell \in \Vdual = \exterior^1(\Vdual) $, the element $\eta(\ell) \in \Omega^1(\Omega^1(-;V))$ is defined by
  \begin{equation*}
    \eta(\ell)(\alpha \tensor v) = \langle v, \ell\rangle\, \alpha
  \end{equation*}
  for $\alpha \in \Omega^1$ and $v \in V$. This is then extended to a map of differential graded algebras.

  In other words, the theorem says every natural transformation
  \begin{equation*}
    \omega_M\colon \Omega^1(M; V) \to \Omega^p(M)
  \end{equation*}
  is (uniquely) a linear combination of transformations of the form
  \begin{equation*}
    \sum \alpha_i \tensor v_i \mapsto \sum_{I, J} M_{I, J}(v_{i_1}, \ldots, v_{i_k}, v_{j_1}, \ldots, v_{j_\ell})\, \alpha_{i_1} \wedge \cdots \wedge \alpha_{i_k} \wedge \d \alpha_{j_1} \wedge \cdots \wedge \d \alpha_{j_\ell}
  \end{equation*}
  where $M_{I, J}$ is anti-symmetric in the first $k$ variables and symmetric in the last $\ell$.
\end{theorem}

Using this, we conclude:

\begin{theorem}
  The Chern--Weil homomorphism gives an isomorphism
  \begin{equation*}
    (\Sym^\bullet\gdual)^G \isomorphism \Omegabullet(\BunGnabla ) \comma  \index[terminology]{Chern--Weil homomorphism}
  \end{equation*}
  and the differential on $\Omegabullet(\BunGnabla )$ is zero.
\end{theorem}

\noindent Note that this $\Sym^\bullet\gdual$ is different from that appearing in the Koszul complex.

\begin{proof}
  We apply the criteria in \cref{lemma:omega-mod-g}. The first condition is the $G$-invariance condition, and translates to the $(-)^G$ part of the statement. So we have to check that the forms satisfying the second condition are isomorphic to $\Sym^\bullet\gdual$.

  To do so, we have to compute the action of $\iota_\xi$ on $\Omega^1(-;\g)$ following the recipe in \cref{remark:iota}. Fix $\omega \in \Omega^p(\Omega^1(-;\g))$ and $\xi \in\g$.

  Let $\phi \colon M \to \EnablaG$ be a trivial principal $G$-bundle with connection $A \in \Omega^1(M;\g)$.
  The induced principal $G$-bundle on $G \times M$ under the action then has connection $\theta + \Ad_{g^{-1}} A$. So by definition,
  \begin{equation*}
    (\iota_\xi \omega)_M(A) = \left.\iota_\xi \left(\omega(\theta + \Ad_{g^{-1}} A)\right)\right|_{\{e\} \times M}\period
  \end{equation*}

  To compute the action on $\Kos^\bullet\gdual$, it suffices to compute it on $\exterior^1\gdual$ and $\Sym^1\gdual$.

  \begin{enumerate}
    \item If $\lambda \in\gdual = \exterior^1\gdual$, then $\lambda(A) = \langle A, \lambda\rangle$, and
      \begin{equation*}
        \iota_\xi \langle \theta + \Ad_{g^{-1}} A, \lambda\rangle = \langle \iota_\xi \theta + \iota_\xi \Ad_{g^{-1}} A, \lambda\rangle.
      \end{equation*}
      We know $\iota_\xi \theta = \xi$, and $\iota_\xi \Ad_{g^{-1}} A = 0$ since $\Ad_{g^{-1}} A$ vanishes on all vectors in the $G$ direction. So we know
      \begin{equation*}
        \iota_\xi \lambda = \langle \xi, \lambda\rangle \in \exterior^0\gdual.
      \end{equation*}

    \item Next, $\tilde{\lambda}(A) = \langle \d A, \lambda \rangle$. We compute
      \begin{align*}
        \iota_\xi \langle \d (\theta + \Ad_{g^{-1}}A), \lambda \rangle|_{\{e\} \times M} &= \left.\iota_\xi \left\langle -\frac{1}{2}[\theta, \theta] + \Ad_{\d g^{-1}} \wedge A + \Ad_{g^{-1}} \d A, \lambda\right\rangle\right|_{\{e\} \times M} \\
        &= \langle -\Ad_\xi A, \lambda\rangle \\ 
        &= \langle A, -\Ad_\xi^* \lambda\rangle \period
      \end{align*}
      So
      \begin{equation*}
        \iota_\xi \tilde{\lambda} = -\Ad^*_\xi \lambda \in \exterior^1\gdual \period
      \end{equation*}
  \end{enumerate}

  First observe that in $\exterior^\bullet\gdual$, the only elements killed by $\iota_\xi$ are those in $\exterior^0\gdual \cong \RR $. To take care of the $\Sym$ part, set
  \begin{equation*}
    \Omega_\lambda = \tilde{\lambda} + \frac{1}{2}[\lambda, \lambda]\period
  \end{equation*}
  Since $\tilde{\lambda}(A) = \langle \d A, \lambda\rangle$, we see that $\Omega_\lambda(A) = \langle \Omega_A, \lambda\rangle$, where $\Omega_A$ is the curvature, and one calculates $\iota_\xi \Omega_\lambda = 0$. By a change of basis, we can identify
  \begin{equation*}
    \Kos^\bullet(\gdual) \cong \exterior^\bullet(\gdual) \tensor \Sym^\bullet \langle \Omega_\lambda\colon \lambda \in\gdual\rangle \comma
  \end{equation*}
  and $\iota_\xi$ vanishes on the second factor entirely. 
  So we are done.
\end{proof}

More generally, the same proof shows that:

\begin{theorem}
  If $M$ is a smooth manifold, the de Rham complex of $M \times (\Omega^1(-;V))$ is the total complex of $\Omegabullet(M; \Kos^\bullet(\Vdual))$. 
  
 In particular, if $M$ has a $G$-action, then $(M \times (\Omega^1(-;\g)) )\modmod G$ is exactly the Cartan model for equivariant de Rham cohomology.  \index[terminology]{Cartan Model}
\end{theorem}

\noindent See \Cref{thm-Cartan} for more on the Cartan model. 

This would follow immediately if we had a result that said that 
\begin{equation*}
  \Omegabullet(M \times \Fcal) \cong \Omegabullet(M) \tensorhat \Omegabullet (\Fcal) \period
\end{equation*}
Here $ \tensorhat $ denotes the \textit{completed} tensor product.
Since $\Omegabullet(\Omega^1(-;\g))$ is finite dimensional, in the case of interest the completed tensor product is the usual tensor product.

%-------------------------------------------------------------------%
%-------------------------------------------------------------------%
%  The proof                                                        %
%-------------------------------------------------------------------%
%-------------------------------------------------------------------%

\subsection{The proof}\label{section:proof}

We now prove of \Cref{thm:main}. The $p = 0$ case is trivial, so assume $p > 0$.

Recall that we have to show that any natural transformation
\begin{equation*}
  \omega_M\colon \Omega^1(M; V) \to \Omega^p(M)
\end{equation*}
is (uniquely) a linear combination of transformations of the form
\begin{equation*}
  \sum \alpha_i \tensor v_i \mapsto \sum_{I, J} M_{I, J}(v_{i_1}, \ldots, v_{i_k}, v_{j_1}, \ldots, v_{j_\ell})\, \alpha_{i_1} \wedge \cdots \wedge \alpha_{i_k} \wedge \d \alpha_{j_1} \wedge \cdots \wedge \d \alpha_{j_\ell} \period
\end{equation*}
The uniqueness part is easy to see since we can extract $M_{I, J}$ by evaluating $\omega_M(\alpha)$ for $M$ of dimension large enough. So we have to show every $\omega_M$ is of this form.

The idea of the proof is to first use naturality to show that for $x \in M$, the form $\omega_M(\alpha)_x$ depends only on the $N$-jet of $\alpha$ at $x$ for some large but finite number $N$ (of course, \emph{a posteriori}, $N = 1$ suffices). Once we know this, the problem is reduced to one of finite dimensional linear algebra and invariant theory.

\begin{lemma}
  For $\omega \in \Omega^p(\Omega^1(-;V))$ and $\alpha \in \Omega^1(M; V)$, the value of $\omega_M(\alpha)$ at $x \in M$ depends only on the $N$-jet of $\alpha$ at $p$ for some $N$. In fact, $N = p$ suffices.
\end{lemma}
We elect to introduce the constant $N$, despite it being equal to $p$, because the precise value does not matter.

\begin{proof}
  Suppose $\alpha$ and $\alpha'$ have identical $p$-jets at $x$. Then there are functions $f_0, f_1, \ldots, f_p$ vanishing at $p$ and $\beta \in \Omega^1(M; V)$ such that
  \begin{equation*}
    \alpha' = \alpha + f_0 f_1 \cdots f_p \beta\period
  \end{equation*}

  The first step is to replace the $f_i$ with more easily understood coordinate functions. Consider the maps
    \begin{equation*}
      \begin{tikzcd}[column sep=6em]
        M \ar[r, "{1_M \times (f_0, \ldots, f_p)}"] & M \times \RR ^{p + 1} \ar[r, "\pr_1"] & M.
      \end{tikzcd}
    \end{equation*}
  Let $\tilde{\alpha}, \tilde{\beta}$ be the pullbacks of the corresponding forms under $\pr_1$, and $t_0, \ldots, t_p$ the standard coordinates on $\RR ^{p + 1}$. Then $\alpha, f_0 f_1\cdots f_p \beta$ are the pullbacks of $\tilde{\alpha}, t_0 t_1\cdots t_p \tilde{\beta}$ under the first map.

  So it suffices to show that $\omega_{M \times \RR ^{p + 1}}(\tilde{\alpha})$ and $\omega_{M \times \RR ^{p + 1}}(\tilde{\alpha} + t_0 t_1 \cdots t_p \tilde{\beta})$ agree as $p$-forms at $(x, 0)$.

  The point now is that by multilinearity of a $p$-form, it suffices to evaluate these $p$-forms on $p$-tuples of standard basis vectors (after choosing a chart for $M$), and there is at least one $i$ for which the $\partial_{t_i}$ is not in the list. So by naturality we can perform this evaluation in the submanifold defined by $t_i = 0$, in which these two $p$-forms agree. 
\end{proof}

By naturality, we may assume $M = W$ is a vector space and $x$ is the origin. The value of $\omega_W(\alpha)$ at the origin is given by a map
\begin{equation*}
  \tilde{\omega}_W\colon \Jet^N(W; \Wdual \tensor V) \to \exterior^p(\Wdual) \comma
\end{equation*}
where $\Jet^N(W; \Wdual \tensor V)$ is the space of $N$-jets of elements of $\Omega^1(W; V)$. This is a finite dimensional vector space, given explicitly by
\begin{equation*}
  \Jet^N(W; \Wdual \tensor V) = \bigoplus_{j = 0}^N \Sym^j(\Wdual) \tensor \Wdual \tensor V \period
\end{equation*}
Under this decomposition, the $j$-th piece captures the $j$-th derivatives of $\alpha$. 
Throughout the proof, we view $\Sym^j(\Wdual)$ as a \emph{quotient} of $(\Wdual)^{\tensor j}$, hence every function on $\Sym^j(\Wdual)$ is in particular a function on $(\Wdual)^{\tensor j}$.

At this point, everything else follows from the fact that $\tilde{\omega}_W$ is functorial in $W$, and in particular $\GL(W)$-invariant.

\begin{lemma}
  The map $\tilde{\omega}_W\colon \Jet^N(W; \Wdual \tensor V) \to \exterior^p(\Wdual)$ is a polynomial function.
\end{lemma}

This lemma is true in much greater generality --- it holds for any set-theoretic natural transformation between ``polynomial functors'' $\Vect \to \Vect$. Here a set-theoretic natural transformation is a natural transformations of the underlying set-valued functors. This is a polynomial version of the fact that a natural transformation between additive functors is necessarily additive, because being additive is a \emph{property} and not a structure.

\begin{proof}
	Write
	\begin{equation*}
		F(W) \colonequals \bigoplus_{j = 0}^N \Sym^j \Wdual \tensor \Wdual \tensor V \andeq G(W) \colonequals \exterior^p W \period
	\end{equation*}
	We think of these as functors $\Vect \to \Vect$ (with $V$ fixed). 
  The point is that for $f \in \Hom_\Vect(W, W')$, the functions $F(f), G(f)$ are polynomial in $f$. 
  This together with naturality will force $\tilde{\omega}_W$ to be polynomial as well.

	To show that $\tilde{\omega}_W$ is polynomial, we have to show that if $v_1, \ldots, v_n \in F(W)$, then $\tilde{\omega}_W(\sum \lambda_i v_i)$ is a polynomial function in $\lambda_1, \ldots, \lambda_n$. Without loss of generality, we may assume each $v_i$ lives in the $(j_i - 1)$-th summand (so that the summand has $j_i$ tensor powers of $\Wdual$).

	Fix a number  $j$ such that $j_i \mid j$ for all $i$. We first show that $\tilde{\omega}_W(\sum \lambda_i^j v_i)$ is a polynomial function in the $\lambda_i$'s.

	Let $f \colon W^{\oplus n} \to W^{\oplus n}$ be the map that multiplies by $\lambda_i^{j / j_i}$ on the $i$-th factor, and let $\Sigma \colon W^{\oplus n} \to W$ be the sum map. Consider the commutative diagram
	\begin{equation*}
		\begin{tikzcd}[sep=2.5em]
			F(W^{\oplus n}) \ar[r, "F(f)"] \ar[d, "\tilde{\omega}_{W^{\oplus n}}"'] & F(W^{\oplus n}) \ar[r, "F(\Sigma)"]\ar[d, "\tilde{\omega}_{W^{\oplus n}}"] & F(W)\ar[d, "\tilde{\omega}_{W}"] \\
			G(W^{\oplus n}) \ar[r, "G(f)"'] & G(W^{\oplus n}) \ar[r, "G(\Sigma)"'] & G(W)
		\end{tikzcd}
	\end{equation*}
	Let $\tilde{v}_i \in F(W^{\oplus n})$ be the image of $v_i$ under the inclusion of the $i$-th summand. 
  Then $x = \sum \tilde{v}_i$ gets sent along the top row to $\sum \lambda_i^j v_i$. 
  On the other hand, $\tilde{\omega}_{W^{\oplus n}}(x)$ is some element in $G(W^{\oplus n})$, and whatever it might be, the image along the bottom row gives a polynomial function in the $\lambda_i^{j/j_i}$, hence in the $\lambda_i$. 
  So we are done.

	We now know that for any finite set $v_1, \ldots, v_n$, we can write
	\begin{equation*}
		\tilde{\omega}_W(\lambda_1^j v_1 + \cdots + \lambda_n^j v_n) = \sum_{r_1, \ldots, r_m} a_R \lambda_1^{r_1} \cdots \lambda_n^{r_n} \period
	\end{equation*}
	We claim each $r_i$ is a multiple of $j$ (if the corresponding $a_R$ is non-zero). 
	Indeed, if we set
	\begin{equation*}
		\lambda_i \colonequals (\mu_i^j - \nu_i^j)^{1/j} \comma
	\end{equation*}
	then the result must be a polynomial in the $\mu_i$ and $\nu_i$ as well, since it is of the form
	\begin{equation*}
		\tilde{\omega}_W(\sum \mu_i^j v_i - \nu_i^j v_i) \period
	\end{equation*}
	But
	\begin{equation*}
		\sum a_R (\mu_1^j - \nu_1^j)^{r_1/ j} \cdots (\mu_n^j - \nu_n^j)^{r_n/ j}
	\end{equation*}
	is polynomial in $\mu_i, \nu_i$ if and only if $j \mid r_i$.

	Now by taking $j$-th roots, we know $\tilde{\omega}_W(\sum \lambda_i v_i)$ is polynomial in the $\lambda_i$ when $\lambda_i \geq 0$. That is, it is polynomial when restricted to the cone spanned by the $v_i$'s. But since the $v_i$'s are arbitrary, this implies it is polynomial everywhere.
\end{proof}

\begin{lemma}
  Any non-zero $\GL(W)$-invariant linear map
  \begin{equation*}
    (\Wdual)^{\tensor M} \to \exterior^p(\Wdual)
  \end{equation*}
  has $M = p$ and is a multiple of the anti-symmetrization map. 
  In particular, any such map is anti-symmetric.
\end{lemma}

\begin{proof}
  For convenience of notation, replace $\Wdual$ with $W$. Since the map is in particular invariant under $\RR ^\times \subseteq \GL(W)$, we must have $M = p$. By Schur's lemma, the second part of the lemma is equivalent to claiming that if we decompose $W^{\tensor p}$ as a direct sum of irreducible $\GL(W)$ representations, then $\exterior^p W$ appears exactly once. In fact, we know the complete decomposition of $W^{\tensor p}$ by Schur--Weyl duality.

  Let $\{V_\lambda\}$ be the set of irreducible representations of $S_p$. Then as an $S_p \times \GL(W)$-representation, we have
  \begin{equation*}
    W^{\tensor p} = \bigoplus_\lambda V_\lambda \tensor W_\lambda \comma
  \end{equation*}
  where $W_\lambda = \Hom_{S_p} (V_\lambda, W^{\tensor p})$ is either zero or irreducible, and are distinct for different $\lambda$. Under this decomposition, $\exterior^p W$ corresponds to the sign representation of $S_p$.
\end{proof}

So we know $\tilde{\omega}_W$ is a polynomial in $\bigoplus_j \Sym^j(\Wdual) \tensor \Wdual \tensor V$, and is anti-symmetric in the $\Wdual$. So the only terms that can contribute are when $j = 0$ or $j = 1$. In the $j = 1$ case, it has to factor through $\exterior^2 \Wdual \tensor V$. So $\tilde{\omega}_W$ is polynomial in $(\Wdual \tensor V) \oplus (\exterior^2 \Wdual \tensor V)$. This exactly says $\omega_W(\alpha)$ is given by wedging together $\alpha$ and $\d \alpha$ (and pairing with elements of $\Vdual$).

\newpage
%!TEX root = ../diffcoh.tex

%-------------------------------------------------------------------%
%-------------------------------------------------------------------%
%  Bott's method                                                    %
%-------------------------------------------------------------------%
%-------------------------------------------------------------------%

\section{Bott's method}\label{BottsMethod}
\textit{by Araminta Amabel}

For $G$ a Lie group, recall the sheaf of groupoids $\BbulletG$ from \Cref{ex:BunG}. 
The goal of this section is to prove Bott's theorem \cite[Theorem 1]{BottPaper}:

\begin{theorem}
	There is an isomorphism
	\[
		\H^{p}(\BbulletG;\Omega^q)=\Hcont^{p-q}(G;\Sym^q(\gdual)) \comma
	\]
	where the right-hand side is the continuous cohomology group.
\end{theorem}

%-------------------------------------------------------------------%
%  Motivation and set up                                            %
%-------------------------------------------------------------------%

\subsection{Motivation and set up}

Let $G$ be a Lie group. 
Recall the Chern--Weil homomorphism 
\[
	\phi \colon \Sym(\gdual)^G\to \H^*(\BG;\RR) \period  \index[terminology]{Chern--Weil homomorphism}
\]
Here, $\gdual$ denotes the linear dual of $\g$. 
We view $\gdual$ as a $G$-module under the adjoint action. 
If $G$ is compact, then this map $\varphi$ is an isomorphism.

Given any principal $G$-bundle on $X$ with connection, we get an induced map
\[
	\Sym(\gdual)^G\to \Omega^*(X) \period
\]
 Taking $X=\BG$ with principal $G$-bundle $\EG\to \BG$, recovers the universal case, $\varphi$. Note that this construction depends on a choice of connection, but this dependence no longer matters once we descend to cohomology. Bott's method will allow us to construct a similar map with no mention of a connection.
 %fix

%-------------------------------------------------------------------%
%-------------------------------------------------------------------%
%  Continuous cohomology                                            %
%-------------------------------------------------------------------%
%-------------------------------------------------------------------%

\subsection{Continuous cohomology}

The following definition can be found in \cite[\S 2]{MR494071}. 

\begin{definition} \index[terminology]{Continuous Cohomology} \index[notation]{Continuous cohomology@$\Hcont^\bullet(G;W)$}
	Let $G$ be a topological group. 
	Let $W$ be a $G$-space. 
	Then the \emph{continuous cohomology} of $G$ with coefficients in $W$ is the cohomology $\Hcont^p(G;W)$ of the cochain complex
	\[\Mapcont(G^{\times p},W)\]
	of continuous maps, with differential 
	\[\del\colon \Mapcont(G^{\times p},W)\to \Mapcont(G^{\times p+1},W)\]
	sending a map $f\colon G^{\times p}\to W$ to the map $(\del f)\colon G^{\times p+1}\to W$ by
	\begin{align*}
		(\del f)(g_1,\dots,g_{p+1}) \colonequals f(g_2,\dots,g_{p+1}) &+ \left(\sum_{i=1}^{p} (-1)^i f(g_1,\dots, g_ig_{i+1},\dots,g_{p+1})\right) \\ 
		&+(-1)^{p+1}f(g_1,\dots,g_p)\cdot g_{p+1} \period
	\end{align*}
\end{definition}

Note that on the third term in $(\del f)$, we are using the action of $G$ on $W$.

\begin{example}
Let $G$ be a topological group and $W$ a $G$-module. 
The zeroeth continuous cohomology of $G$ with values in $W$ is the fixed points, 
\[\Hcont^{0}(G;W)\simeq W^G\period
\]
\end{example}

The following theorem of van Est can be found in \cite{MR0059285}. 

\begin{theorem}[(van Est)] 
Let $G$ be a connected Lie group and $K\subset G$ a maximal compact subgroup. Then there is an equivalence
\[ \Hcont^\bullet(G;A)\simeq \HLie^\bullet(\g,\kfrak;A)\]
for any $G$-space $A$.
\end{theorem} \index[terminology]{van Est Theorem}

See \cite[\S 5]{MR494071} for a discussion of this result, and \cite{MR147577} for generalizations. 

\begin{cor}
	Let $G$ be a compact, connected Lie group. For $i>0$,
	\[
		\Hcont^i(G;A)=0 \period
	\]%this is corollary of Van Est theorem 
\end{cor}%i think can prove directly too?

%-------------------------------------------------------------------%
%-------------------------------------------------------------------%
% Relating continuous cohomology to ordinary cohomology             %
%-------------------------------------------------------------------%
%-------------------------------------------------------------------%

\subsection{Relating continuous cohomology to ordinary cohomology}

We would like to produce a map
\[ 
	\H^\bullet(\BG;\RR)\to \Hcont^\bullet(G;\RR)
\]
when $G$ is a connected Lie group. We will produce this map as the edge map of a spectral sequence.

For $K$ a Lie group, let $ \kfrak \colonequals \Lie(K)$.

\begin{lem}\label{ss}
	Let $G$ be a connected Lie group with maximal compact subgroup $K$. 
	There is a spectral sequence whose $E_1$ term is 
	\[
		E_1^{p,q}=\left(\exterior^p((\g/\kfrak)^\vee)\otimes \Sym^q(\gdual)\right)^\kfrak
	\]
	converging to 
	\[
		E_\infty^{p,q}=\Sym^{q-p}(\kdual)\period
	\]
\end{lem}

\begin{proof}
	Note that $\g$ splits as
	\[
		\g\simeq \g/\kfrak\oplus\kfrak\period
	\]
	Thus we can rewrite the $E_1$ page as
	\[
		E_1^{p,q}=\left(\exterior^p((\g/\kfrak)^\vee)\otimes \bigoplus_{a+b=q}\Sym^a((\g/\kfrak))^\vee\otimes \Sym^b(\kfrak)^\vee)\right)^\kfrak\period
	\]

	Note that the terms $\exterior^p((\g/\kfrak)^\vee)$ and $\Sym^p((\g/\kfrak)^\vee)$ are Koszul dual. 
	During the course of the spectral sequence, these Koszul dual terms cancel each other. 
	The $E_\infty$ page is thus
	\[
		E_\infty^{p,q}=\Sym^{q-p}(\kdual)\period \qedhere
	\]
\end{proof}

We can compute the $E_2$ term of this spectral sequence directly. The $E_1$ page comes from the relative Chevalley--Eilenberg complex,
\[E_1^{p,q}=\HLie^p(\g,\kfrak;\Sym^q(\gdual))\period
\]

The $d_1$ differential is the Chevalley--Eilenberg differential. Thus the $E_2$ page is just relative Lie algebra cohomology,
\[E_2^{p,q}=\HLie^p(\g,\kfrak;\Sym^q(\gdual)) \period\]
By the van Est theorem, this relative Lie algebra cohomology can be recognized in terms of continuous cohomology,
\[
	\HLie^p(\g,\kfrak;\Sym^q(\gdual))\simeq \Hcont^p(G;\Sym^q(\gdual)) \period
\]

\begin{cor}
	Let $G$ be a connected Lie group with maximal compact subgroup $K$.
	There is a map $\H^\bullet(\BG;\RR)\to \Hcont^\bullet(G;\RR)$.
\end{cor}

\begin{proof}
	One of the edge maps of the spectral sequence from \Cref{ss} goes from the $E_\infty$ term to the $E_2^{p,0}$ column. 
	Since $K$ is compact, the $E_\infty$ term can be identified with $\H^\bullet(\BK;\RR)$ be the Chern--Weil homomorphism. The $E_2^{p,0}$ column is  
	\[
		\Hcont^p(G;\Sym^0(\gdual))\simeq \Hcont^p(G;\RR) \period \qedhere
	\]
\end{proof}

\newpage
%!TEX root = ../diffcoh.tex

%-------------------------------------------------------------------%
%-------------------------------------------------------------------%
%  Lifts of Chern classes                                           %
%-------------------------------------------------------------------%
%-------------------------------------------------------------------%

\section{Lifts of Chern classes}\label{LiftsofChernClasses}
\textit{Talk by Mike Hopkins}\\
\textit{Notes by Araminta Amabel}
 %maybe should put in reminder that have lifts with connection all the way to Z(2n). use virasoro as motivation for getting them half way up? or add in explanation why they should only go half way up
%define B_\bullett G if haven't yet. say that can take \H^*(\BbulletG;\ZZ(n))$ tthings

%-------------------------------------------------------------------%
%-------------------------------------------------------------------%
%  Introduction                                                     %
%-------------------------------------------------------------------%
%-------------------------------------------------------------------%

\subsection{Introduction}

Let $\ZZ(n)$ be the Deligne complex  \index[notation]{Z(n)@$\ZZ(n)$}
\begin{equation*}
	\ZZ\to\Omega^0\to\cdots\to \Omega^{n-1}
\end{equation*}
We'll also let $\ZZ(\infty)$ denote the untruncated complex,
\begin{equation*}
	\ZZ\to\Omega^0\to\cdots
\end{equation*}
Similarly, we define $\RR(n)$ where $n=1,\dots,\infty$ to be the complex  \index[notation]{R(n)@$\RR(n)$}
\begin{equation*}
	\RR\to\Omega^0\to\cdots\to \Omega^{n-1}
\end{equation*}
and $\ZZ_\CC(n)$ to be the complex  \index[notation]{ZC(n)@$\ZZ_\CC(n)$}
\begin{equation*}
	\ZZ\to\Omega^0_\CC\to\Omega^1_\CC\to\cdots\to\Omega^{n-1}_\CC \period
\end{equation*}
One can also think of $\ZZ(n)$ as the homotopy pullback
\begin{equation*}
	\begin{tikzcd}
		\ZZ(n)\arrow[r]\arrow[d] & \ZZ\arrow[d]\\
		\Sigma^{-n}\Omegacl^n\arrow[r] & \RR \period
	\end{tikzcd}
\end{equation*}
One take away is that there are a lot more characteristic classes in differential cohomology than you would expect.

%-------------------------------------------------------------------%
%  Virasoro group motivation                                        %
%-------------------------------------------------------------------%

\subsubsection{Virasoro group motivation}

The Virasoro group is a certain central extension of $\Diffplus(\Circ)$ by $\Uup_1$,  \index[terminology]{Virasoro!group}
\begin{equation*}
	\Uup_1\to\Vir\to\Diffplus(\Circ) \period
\end{equation*}
Let $\Gamma=\Diffplus(\Circ)$ be the group of orientation preserving diffeomorphisms of $\Circ$. Central extensions of $\Gamma$ are classified by elements of $\H^3(\BGamma;\ZZ(1))$; i.e., by homotopy classes of maps $\BGamma\to \K(\ZZ(1), 3)$. We have a fiber sequence
\begin{equation*}
	\K(\ZZ(1),2)\to \EGamma\to \BGamma\to \K(\ZZ(1), 3) \period
\end{equation*}
Consider the fibration $\EGamma\times_\Gamma \Circ\to \BGamma$ with fiber $\Circ$. Integration along the fibers gives a map
\begin{equation*}
	\H^4(\EGamma\times_\Gamma \Circ;\ZZ(2))\to \H^3(\BGamma;\ZZ(1)) \period 
\end{equation*}
There is a map $\EGamma\times_\Gamma \Circ\to \BSL(\RR)$. Thus given a class $\ptilde_1\in \H^4(\BSL(\RR);\ZZ(2))$, we can pull it back to get a class in $\H^4(\EGamma\times_\Gamma \Circ;\ZZ(2))$. 
Integrating along the fiber produces a class in $\H^3(\BGamma;\ZZ(1))$. 
Thus, classes in $\H^4(\BSL(\RR);\ZZ(2))$ produce central extensions of $\Diffplus(\Circ)$. 

%-------------------------------------------------------------------%
%  Hopes                                                            %
%-------------------------------------------------------------------%

\subsubsection{Hopes}

Let $G$ be a Lie group. 
Recall the sheaf of groupoids $\BbulletG$ from \Cref{ex:BunG}.
\begin{enumerate}
	\item If $V\to X$ is a real vector bundle, we want lifted Pontryagin classes $\ptilde_n(V)\in \H^{4n}(X;\ZZ(2n))$. 

	\subitem To obtain such lifts, it suffices to construct $\ptilde_n\in \H^{4n}(\BbulletGL{n}{\RR};\ZZ(2n))$ such that $\ptilde_n$ maps to $p_n$ under the map
	\begin{equation*}
		\H^{4n}(\BbulletGL{n}{\CC};\ZZ(2n))\to \H^{4n}(\BGL_n(\CC);\ZZ)
	\end{equation*}

	\item If $W\to X$ is a complex vector bundle, we want (off-diagonal) Chern classes
	\begin{equation*}
		\ctilde_m(W)\in \H^{2m}(X;\ZZ_\CC(m)) \period
	\end{equation*} 

	\subitem To obtain such lifts, it suffices to construct $\ctilde_n\in \H^{2n}(\BbulletGL{n}{\CC};\ZZ_\CC(n))$ such that $\ctilde_n$ maps to $c_n$ under the map
	\begin{equation*}
		\H^{2n}(\BbulletGL{n}{\CC};\ZZ_\CC(n))\to \H^{2n}(\BGL_n(\CC);\ZZ) \period
	\end{equation*}

	\item \textit{Cartan formula:} Given a short exact sequence of vector bundles
	\begin{equation*}
		0\to V\to W\to U\to 0
	\end{equation*}
	an expression of the differential characteristic classes of $W$ in terms of the differential characteristic classes for $U$ and $V$. Every short exact sequence of vector bundles is split, but this splitting might not be smooth. 
	Thus it's possible that such a formula exists for split short exact sequences, $V\oplus U$. %right? what lol

	\item \textit{Projective bundle formula:} More generally, higher characteristic classes being determined by those for line bundles.
\end{enumerate}

%-------------------------------------------------------------------%
%  Statement of results                                             %
%-------------------------------------------------------------------%

\subsubsection{Statement of results}

The following are things Hopkins has worked out and attributes to ideas found in papers of Bott, \cites{BottPaper,BottNotes}

\begin{theorem}
	There is a pullback square
	\begin{equation*}
		\begin{tikzcd}
			\H^{2n}(\BbulletGL{m}{\CC}; \ZZ_\CC(n))\arrow[rr]\arrow[d] & & \H^{2n}(\BU_m;\ZZ)\arrow[d] \\
			\H^n(\BU_m\times \BU_m;\CC)\arrow[rr, "\textup{diagonal}^*"'] & & \H^{2n}(\BU_m;\CC) \period
		\end{tikzcd}
	\end{equation*}
\end{theorem}

\noindent This is \Cref{cor-BGLC} below.

So if we wanted to lift the first Chern class $c_1$, we could take 
\begin{equation*}
	\frac{1}{2}(c_1\otimes 1+ 1\otimes c_1)\in \H^2(\BU_1\times \BU_1;\CC) \period
\end{equation*}
But, could also add to this any terms that are in the kernel of the diagonal map. So there are many possible off-diagonal lifts of $c_1$ to something with $\ZZ_\CC(1)$ coefficients. 

Using the $e^{2\pi i}$ induced isomorphism $\K(\ZZ_\CC(1);2)\similarrightarrow \BGL_1(\CC)$ produces the lift of $c_1$ corresponding to $\frac{1}{2}(c_1\otimes 1+1\otimes c_1)$. %explain this

\begin{remark}
This also works for products of copies of $\GL_n(\CC)$. For example, let 
\begin{equation*}
	G=\GL_n(\CC)\times\cdots\times \GL_n(\CC) \period
\end{equation*}
Then we have a pullback square
\begin{equation*}
	\begin{tikzcd}
		\H^{2n}(\BbulletG; \ZZ_\CC(n))\arrow[rr]\arrow[d] & & \H^{2n}(\BG;\ZZ)\arrow[d] \\
		\H^n(\BG\times \BG;\CC)\arrow[rr, "\textup{diagonal}^*"'] & & \H^{2n}(\BG;\CC) \period
	\end{tikzcd}
\end{equation*}
\end{remark}

\noindent Let $P_{a|b}\subset \GL_{a+b}(\CC)$ be the subset of matrices of the form 
\begin{equation*}
	\begin{pmatrix}
		A & B \\
		0 & C
	\end{pmatrix} \comma
\end{equation*}
where $A$ is an $(a\times a)$-matrix and $B$ is a $(b\times b)$-matrix. %this is some paragbolic?
Note that there is a map
\begin{equation*}
	\GL_a(\CC)\times \GL_b(\CC)\to P_{a|b}
\end{equation*}
sending $(A,B)$ to the block matrix with $A$ and $B$ on the diagonal.

\begin{conjecture}
	The induced map
	\begin{equation*}
		\H^{2n}(P_{a|b};\ZZ_\CC(n))\to \H^{2n}( \GL_a(\CC)\times \GL_b(\CC);\ZZ_\CC(n))
	\end{equation*}
	is an isomorphism.
\end{conjecture} 

\begin{proof}[Proof Outline]
	Completing \Cref{exercise-Pab} below, one should find that 
	$\H^{2n}(P_{a\vert b};\ZZ_\CC(n))$ fits into a pullback diagram
	\begin{equation*}
		\begin{tikzcd}
			\H^{2n}(\Bun_{P_{a\vert b}}; \ZZ_\CC(n))\arrow[rr]\arrow[d] & & \H^{2n}(\Bup(P_{a\vert b}\cap U_{a+b});\ZZ)\arrow[d]\\
			\H^n(\Bup(P_{a\vert b}\cap U_{a+b}) \times \Bup(P_{a\vert b}\cap U_{a+b});\CC)\arrow[rr, "\textup{diagonal}^*"'] & & \H^{2n}(\Bup(P_{a\vert b}\cap U_{a+b});\CC)
		\end{tikzcd}
	\end{equation*}
	and $\H^{2n}(\GL_a(\CC)\times \GL_b(\CC));\ZZ_\CC(n))$ fits into a pullback diagram
	\begin{equation*}
		\begin{tikzcd}
			\H^{2n}(\Bun_{\GL_a(\CC)\times \GL_b(\CC)}; \ZZ_\CC(n))\arrow[rr]\arrow[d] & & \H^{2n}(\BU_{a+b};\ZZ)\arrow[d]\\
			\H^n(\BU_{a+b}\times \BU_{a+b};\CC)\arrow[rr, "\textup{diagonal}^*"'] & & \H^{2n}(\BU_{a+b};\CC)
		\end{tikzcd}
	\end{equation*}
	Since every short exact sequence of vector bundles splits, the inclusion $\BU_{a+b}\hookrightarrow \BP_{a\vert b}$ is a homotopy equivalence. Thus so is the inclusion $\BU_{a+b}\hookrightarrow \Bup(P_{a\vert b}\cap U_{a+b})$. 
	Hence the lower left corners of the above two pullback diagrams are isomorphic.
\end{proof}

\noindent Thus if we have a Cartan-like formula for split short exact sequences, we can get a Cartan-like formula for any short exact sequence.

The following is an example of \Cref{cor-BGZ} below.

\begin{theorem}\label{thm-BGZ}
	There is a pullback square
	\begin{equation*}
		\begin{tikzcd}
			\H^{2n}(\BbulletGL{m}{\RR};\ZZ(n))\arrow[r]\arrow[d] & \H^{2n}(\BO_m;\ZZ)\arrow[d]\\
			\H^{2n}(\BGL_m(\CC);\RR)\arrow[r] & \H^{n}(\BO_m;\RR) \period
		\end{tikzcd}
	\end{equation*}
\end{theorem}

\begin{example}
	Take $n=1$ and choose $m$ large. 
	The first Pontryagin class $p_1$ lives in $\H^4(\BO_m;\ZZ)$. 
	By \Cref{thm-BGZ}, off-diagonal differential lifts of $p_1$ are given by a choice of class in 
	\begin{equation*}
		\H^4(\BGL_m(\CC);\RR)\simeq \RR\oplus \RR
	\end{equation*}
	that agrees with the image of $p_1$ in
	\begin{equation*}
		\H^4(\BO_m;\RR)\simeq\RR \period
	\end{equation*}
	Pictorially, there is a pullback diagram
	\begin{equation*}
		\begin{tikzcd}
			\H^4(\BbulletGL{m}{\RR};\ZZ(2)) \arrow[r, "f"] \arrow[d] & \H^4(\BO_m;\ZZ)\arrow[d] \\
			\RR\oplus \RR\arrow[r] & \RR \period
		\end{tikzcd}
	\end{equation*}
	Since this is pullback diagram, the kernel of $f$ is the same as the kernel of the bottom horizontal map. That is, $\ker(f)=\RR$. Thus there is a 1-parameter family of differential lifts of $p_1$. 

	One way to choose such a lift $\ptilde_1$ is to ask for $\ptilde_1$ to be primitive; i.e.,
	\begin{equation*}
		\ptilde_1(V\oplus U)=\ptilde_1(V)+\ptilde_1(U)
	\end{equation*}
	Up to a scalar $\lambda$, there is only one choice of primitive element of $\H^4(\BGL_m;\RR)$ that agrees with $p_1$ in $\H^4(\BO_m;\RR)$. That class is 
	\begin{equation*}
		\frac{1}{2}(\lambda c_1^2-2c_2) \period
	\end{equation*}
\end{example}

%-------------------------------------------------------------------%
%-------------------------------------------------------------------%
%  Computations                                                     %
%-------------------------------------------------------------------%
%-------------------------------------------------------------------%

\subsection{Computations}

Suppose that $G$ is a finite-dimensional Lie group. 
We are interested in computing
\begin{equation*}
	\H^{2n}(\BbulletG;\ZZ(n)) \period
\end{equation*}
We start with $\H^{2n}(\BbulletG;\RR(n))$.

\begin{proposition}
	For all $k$ one has $\H^k(\BbulletG;\RR(\infty))=0$.
\end{proposition}

\begin{proof}
	By definition, $\RR(\infty)$ is the complex
	\begin{equation*}
		\RR \to \Omega^0 \to \Omega^1 \to \cdots
	\end{equation*}
	which is acyclic by the Poincaré Lemma.
\end{proof}

\begin{corollary}\label{zero_below_2n}
	For $k<2n$ one has $ \H^k(\BbulletG;\RR(n))=0 $.
\end{corollary}

\begin{proof}
	We will show that for $k<2n$ the map
	\begin{equation*}
		\H^k(\BbulletG;\RR(n+1))\to \H^k(\BbulletG;\RR(n))
	\end{equation*}
	is surjective. For this we have the long exact sequence associated to the short exact sequence
	\begin{equation*}
		0\to\Sigma^{-(n+1)}\Omega^n\to\RR(n+1)\to\RR(n)\to 0 \period
	\end{equation*}
	It gives us an exact sequence
	\begin{equation*}
		\H^k(\BbulletG;\RR(n+1))\to \H^k(\BbulletG;\RR(n))\to \H^{k-n}(\BbulletG;\Omega^n) \period
	\end{equation*}
	By Bott's theorem \cite[Theorem 1]{BottPaper}, we have %ref thm in these notes
	\begin{equation*}
		\H^{k-n}(\BbulletG;\Omega^n)=\Hcont^{k-2n}(G;\Sym^n(\gdual)) \comma
	\end{equation*}
	where the right-hand side is the continuous cohomology group.
	Since $ k-2n < 0 $, this group is zero.
\end{proof}

\begin{corollary}\label{R_to_Omega}
	The map
	\begin{equation*}
		\H^{2n}(\BbulletG;\RR(n))\to \H^n(\BbulletG;\Omega^n)
	\end{equation*}
	is an isomorphism.
\end{corollary}

\begin{proof}
	This map is part of the long exact sequence
	\begin{equation*}
		\cdots\to \H^{2n}(\BbulletG;\RR(n+1))\to \H^{2n}(\BbulletG;\RR(n))\to \H^{n}(\BbulletG;\Omega^n)\to \H^{2n+1}(\BbulletG;\RR(n+1))\to\cdots
	\end{equation*}
	and the two end terms are zero by \Cref{zero_below_2n}.
\end{proof}

\begin{corollary}\label{cor-CWoff}
	We have an isomorphism
	\begin{equation*}
		\H^{2n}(\BbulletG;\RR(n))\isomorphic\Sym^n(\gdual)^G \period
	\end{equation*}
\end{corollary}

\begin{proof}
	By \Cref{R_to_Omega}, we have an isomorphism
	\begin{equation*}
		\H^{2n}(\BbulletG;\RR(n))\similarrightarrow \H^n(\BbulletG;\Omega^n) \period
	\end{equation*}
	Bott's theorem gives an isomorphism
	\begin{equation*}
		\H^{n}(\BbulletG;\Omega^n)\isomorphic \Hcont^{n-n}(G;\Sym^n(\gdual)) \period
	\end{equation*}
	One has
	\begin{equation*}
		\Hcont^{0}(G;\Sym^n(\gdual))\isomorphic (\Sym^n(\gdual))^G \period \qedhere
	\end{equation*}
\end{proof}

\begin{corollary}\label{cor-BGZ}
\label{2n_to_n}
	For every $n$ there is a pullback square
	\begin{equation*}
		\begin{tikzcd}
			\H^{2n}(\BbulletG;\ZZ(n))\arrow[r]\arrow[d] & \H^{2n}(\BG;\ZZ)\arrow[d]\\
			\Sym^n(\gdual)^G\arrow[r] & \H^{2n}(\BG;\RR) \period
		\end{tikzcd}
	\end{equation*}
\end{corollary}

\begin{proof}
	For this consider the pullback square
	\begin{equation*}
		\begin{tikzcd}
			\ZZ(n)\arrow[r]\arrow[d] & \ZZ\arrow[d]\\
			\RR(n)\arrow[r] & \RR
		\end{tikzcd}
	\end{equation*}
	The associated Mayer--Vietoris sequence shows that the kernel of the map from the upper left corner of
	\begin{equation*}
		\begin{tikzcd}
			\H^{2n}(\BbulletG;\ZZ(n))\arrow[r]\arrow[d] & \H^{2n}(\BG;\ZZ)\arrow[d]\\
			\H^{2n}(\BbulletG;\RR(n))\arrow[r] & \H^{2n}(\BG;\RR)
		\end{tikzcd}
	\end{equation*}
	to the pullback is $\H^{2n-1}(\BbulletG;\RR)$, which is zero by Chern--Weil. %cite earlier too
\end{proof}

Tensoring with $\CC$ gives:

\begin{corollary}
	For every $n$ there is a pullback square
	\begin{equation*}
		\begin{tikzcd}
			\H^{2n}(\BbulletG;\ZZ_\CC(n))\arrow[r]\arrow[d] & \H^{2n}(\BG;\ZZ)\arrow[d]\\
			\Sym^n_\CC(\gdual\otimes\CC)^{G_\CC}\arrow[r] & \H^{2n}(\BG;\CC)
		\end{tikzcd}
	\end{equation*}
	where $G_\CC$ is the complexification of the Lie group $G$.
\end{corollary}
%this needs a complex version of Bott's Theorem
%\begin{equation*}\H^p(\BbulletG;\Omega^q_\CC)\simeq \Hcont^{p-q}(\BG;\Sym_\CC^n(\gdual\otimes\CC)^{G_\CC}\end{equation*}
%I think that just needs
%\begin{equation*}\CC\otimes\Sym^n(\gdual)\simeq\Sym_\CC^n(\gdual\otimes\CC)^{G_\CC}\end{equation*}

\begin{remark}
	When $G$ is connected, the map
	\begin{equation*}
		\Sym^n(\gdual)^G\to\Sym^n(\gdual)^\g
	\end{equation*}
	is an isomorphism. 
	Otherwise, there is a residual action of $\uppi_0G$ and one has an isomorphism
	\begin{equation*}
		\Sym^n(\gdual)^G\to\left(\Sym^n(\gdual)^\g\right)^{\uppi_0G} \period
	\end{equation*}
\end{remark}

We now turn to evaluating these groups.

\begin{example}
	Let's take $G = \GL_n(\CC)$. 
	Then since $ \GL_n(\CC) $ is connected, we have
	\begin{equation*}
		\Sym^n(\gdual)^G=\Sym^n(\gdual)^\g
	\end{equation*}
	which depends only on $\g$. 
	Since $\g$ is complex, we have
	\begin{equation*}
		\g\otimes\CC\isomorphic \g\oplus\g
	\end{equation*}
	and so 
	\begin{equation*}
		\CC\otimes\Sym^n(\gdual)^G = \Sym^n_\CC(\gdual\oplus\gdual)^{\g\oplus\g}
	\end{equation*}
	Now $\g$ is also the complexification of the Lie algebra $\ufrak_n$ of the unitary group $\Uup_n$. Thus the above is isomorphic to
	\begin{equation*}
		\CC\otimes(\Sym^n(\ufrak_n\oplus\ufrak_n))^{\Uup_n\times \Uup_n}
	\end{equation*}
	which, by Chern--Weil, is 
	\begin{equation*}
		\H^{2n}(\BU_n\times \BU_n;\CC) \period
	\end{equation*}
\end{example}

\begin{corollary}\label{cor-BGLC}
	There is a pullback diagram
	\begin{equation*}
		\begin{tikzcd}
			\H^{2n}(\BbulletGL{m}{\CC};\ZZ_\CC(n))\arrow[r]\arrow[d] & \H^{2n}(\BU_m;\ZZ)\arrow[d] \\
			\H^{2n}(\BU_m\times \BU_m;\CC)\arrow[r] & \H^{2n}(\BU_m;\CC)
		\end{tikzcd}
	\end{equation*}
\end{corollary}

\begin{example}
	Let's now take the case $G = \GL_n(\RR)$. The main thing now is to compute
	\begin{equation*}
		\Sym^n(\gdual)^G=\left(\Sym^n(\gdual)^\g\right)^{\ZZ/2}
	\end{equation*}
	Using Weyl's unitary trick again, we can complexify and recognize
	\begin{equation*}
		\g_\CC\isomorphic (\ufrak_n)\otimes\CC
	\end{equation*}
	and we find by Chern--Weil that
	\begin{equation*}
		\left(\Sym^n(\gdual)^\g\right)_\CC\isomorphic \H^{2n}(\BU_m;\CC) \period 
	\end{equation*}
	The action of $\Gal(\CC/\RR)$ is complex conjugation on both $\CC$ and on $\Uup_m$ so
	\begin{equation*}
		\H^{2n}(\BU_m;\CC)^{\Gal(\CC/\RR)} \isomorphic \begin{cases}
		\H^{2n}(\BU_m;i\RR), & n\text{ odd} \\
		\H^{2n}(\BU_m;\RR), & n \text{ even} \period
		\end{cases}
	\end{equation*}
	In this case, the action of $\uppi_0\GL_m$ is trivial.
\end{example}

\begin{remark}
	Maybe the easiest way to be convinced of the action of complex conjugation and of $\uppi_0\GL_m$ is to remember the formula for the Chern classes in terms of $\Sym^\bullet(\gdual)$. 
	For $x\in\gl_n(\CC)$, the total Chern class
	\begin{equation*}
		1+c_1t+\cdots+ c_nt^n
	\end{equation*}
	is given by the homogeneous terms in the characteristic polynomial
	\begin{equation*}
		\det\frac{t}{2\pi i}\begin{pmatrix}
		e_{1,1} & \cdots & e_{1,n}\\
		\vdots & \ddots & \vdots\\
		e_{n,1} & \cdots & e_{n,n}
		\end{pmatrix}-1
	\end{equation*}
	where $e_{i,j}\in\gl_n\CC^*$ is the function associating to a matrix its $(i,j)$ entry. 
	If we apply this to a matrix with real entries, we see that the $k$-th chern class lies in $\frac{1}{(2\pi i)^k}\RR$ and that it is invariant under conjugation by any matrix.
\end{remark}

\begin{exercise}\label{exercise-Pab}
	Let $P_{a,b}\subset \GL_{a+b}(\CC)$ be the subgroup which sends vectors whose last $b$ coordinates are zero to vectors whose last $b$ coordinates are zero, as above. One may compute $
		\Sym^\bullet(\pfrak^\vee_{a,b})^{P_{a,b}}$
	by first computing $
		\Sym^\bullet(\pfrak^\vee_{a,b})^{G_{a,b}}$ 
	and appealing to the unitary trick. 
This is the relevant computation for working out a Cartan formula for an exact sequence which does not necessarily split.
\end{exercise}

\newpage
%!TEX root = ../diffcoh.tex

%-------------------------------------------------------------------%
%-------------------------------------------------------------------%
%  The Virasoro algebra                                             %
%-------------------------------------------------------------------%
%-------------------------------------------------------------------%

\section{The Virasoro algebra}\label{VirasoroAlgebra}
\textit{by Arun Debray}

The contents of this section can be summarized as follows:
\begin{itemize}
	\item The Virasoro group is a particular central extension of $\Diffplus(\Circ)$ by $\TT$.
	\item A theorem of Segal \cite[Corollary 7.5]{Seg81} proves that
	\begin{equation}
		\Cent_{\TT}(\Diffplus(\Circ)) \isomorphism \Cent_{\TT}(\PSL_2(\RR))\times\Cent_\RR(\wittRR) \comma
	\end{equation}
	where $\wittRR = \mathrm{Lie}(\Diffplus(\Circ))$ is the Witt algebra. The map is: restrict the central extension to
	$\PSL_2(\RR)\subset\Diffplus(\Circ)$ for the first component, and differentiate for the second component.
	\item It is possible to construct this central extension using an off-diagonal differential lift of $p_1$ and
	a transgression map.
%%	\item The (isomorphism classes of) central extensions created from differential lifts of $p_1$ are expected to
%	be trivial when restricted to $\PSL_2(\RR)$. Proving this might be a good warm-up to studying $\Diffplus(\Circ)$.
%	\item The Virasoro extension is also trivial when restricted to $\PSL_2(\RR)$. Therefore, to identify which
%	off-diagonal differential lift $\ptilde_1$ of $p_1$ induces the Virasoro central extension, it suffices to look at the
%	induced central extension of Lie algebras of $\wittRR$. There is an $\RR$ worth of differential lifts of $p_1$, and
%	$\Cent_\RR(\wittRR)\cong\RR$.
\end{itemize}

%-------------------------------------------------------------------%
%-------------------------------------------------------------------%
%  Review of central extensions                                     %
%-------------------------------------------------------------------%
%-------------------------------------------------------------------%

\subsection{Review of central extensions}

\begin{definition} \index[terminology]{Central Extension}
	Let $G$ be a group and $A$ be an abelian group. A \emph{central extension} of $G$ by $A$ is a short exact sequence
	of groups
	\begin{equation}\label{centralext}
		\begin{tikzcd}[sep=1.5em]
			1\arrow[r] & A\arrow[r] & \widetilde G\arrow[r] & G\arrow[r] & 1,
		\end{tikzcd}
	\end{equation}
	such that $A\subset Z(\widetilde G)$. An equivalence of central extensions is a map of short exact sequences which
	is the identity on $G$ and on $A$. These form an abelian group we denote $\Cent_A(G)$.  \index[notation]{Cent@$\Cent$}
\end{definition}

When $G$ and $A$ have additional structure, we will ask that central extensions respect that structure: for
example, when both are Lie groups (possibly infinite-dimensional), we want~\eqref{centralext} to be a short exact
sequence of Lie groups. 

For discrete $G$ and $A$, central extensions are classified by $\H^2(G;A)$. Explicitly, given a cocycle $b\colon
G\times G\to A$, we build the central extension by setting $\widetilde G = G\times A$ as sets, with the twisted
multiplication
\begin{equation}
	(g_1, a_1)\cdot_b (g_2, a_2) \colonequals (g_1g_2, a_1 + a_2 + b(g_1, g_2)) \period
\end{equation}
Associativity follows from the cocycle condition; if two cocycles are related by a coboundary, their induced
central extensions are equivalent.

Generalizing this to Lie groups is not straightforward --- you can't just use smooth cochains unless $A$ is a
topological vector space. We are interested in central extensions by $\TT$, so we'll have to be craftier. The fix is
due to Segal \cite{Seg70}, and was later rediscovered by Brylinski \cite{Bry00}, following Blanc \cite{Bla85}. We
rephrase it in language familiar to this seminar.

Let $A$ be an abelian Lie group. Throughout today's talk, $\underline A$ denotes the simplicial sheaf on $\msf{Man}$
whose value on a test manifold $M$ is the space of \emph{smooth } maps $M\to A$.\footnote{By contrast, the
simplicial sheaf just denoted ``$A$'' treats $A$ as having the discrete topology. This is a little bit
counterintuitive but is standard notation.}
\begin{theorem}[{(Segal \cite{Seg70}, Brylinski \cite{Bry00})}]
Let $G$ and $A$ be abelian Lie groups. Then, equivalences classes of central extensions in which $\widetilde G\to
G$ is a principal $A$-bundle are classified by $\H^2(\BbulletG; \underline A)$.
\end{theorem}
The idea of the characterization is that $\BbulletG$ admits a simplicial resolution
% \[
% 	\BbulletG \simeq \left(\xymatrix{
% 		{*} & G\ar@<-0.4ex>[l]\ar@<0.4ex>[l] & G\times G\ar@<-0.6ex>[l]\ar[l]\ar@<0.6ex>[l] & G\times G\times
% 		G\ar@<-0.8ex>[l]\ar@<-0.4ex>[l]\ar@<0.4ex>[l]\ar@<0.8ex>[l] &\cdots}\right)
% \]
\begin{equation*}
	\BbulletG \simeq \left(
	\begin{tikzcd}[sep=1.5em]
	    \cdots \arrow[r, shift left=0.75ex] \arrow[r, shift right=0.75ex] \arrow[r, shift right=2.25ex] \arrow[r, shift left=2.25ex] & G \cross G \arrow[l] \arrow[l, shift left=1.5ex] \arrow[l, shift right=1.5ex] \arrow[r] \arrow[r, shift left=1.5ex] \arrow[r, shift right=1.5ex] & G \arrow[l, shift left=0.75ex] \arrow[l, shift right=0.75ex] \arrow[r, shift left=0.75ex] \arrow[r, shift right=0.75ex] & * \arrow[l]
	\end{tikzcd}\right)
\end{equation*}
which is the content of the bar construction (see \Cref{prop:Cech_nerve_trivG}), and we want to compute $\uppi_0$ of the simplicial set of maps
\begin{equation}
	\begin{tikzcd}[sep=2.5em]
	    \cdots \arrow[r, shift left=0.75ex] \arrow[r, shift right=0.75ex] \arrow[r, shift right=2.25ex] \arrow[r, shift left=2.25ex] & G \cross G \arrow[d, blue] \arrow[l] \arrow[l, shift left=1.5ex] \arrow[l, shift right=1.5ex] \arrow[r] \arrow[r, shift left=1.5ex] \arrow[r, shift right=1.5ex] & G \arrow[d] \arrow[l, shift left=0.75ex] \arrow[l, shift right=0.75ex] \arrow[r, shift left=0.75ex] \arrow[r, shift right=0.75ex] & * \arrow[d] \arrow[l] \\
	    \cdots \arrow[r, shift left=0.75ex] \arrow[r, shift right=0.75ex] \arrow[r, shift right=2.25ex] \arrow[r, shift left=2.25ex] & A \arrow[l] \arrow[l, shift left=1.5ex] \arrow[l, shift right=1.5ex] \arrow[r] \arrow[r, shift left=1.5ex] \arrow[r, shift right=1.5ex] & * \arrow[l, shift left=0.75ex] \arrow[l, shift right=0.75ex] \arrow[r, shift left=0.75ex] \arrow[r, shift right=0.75ex] & * \arrow[l]
	\end{tikzcd}
% \xymatrix{
% 	{*}\ar[d] &G\ar@<-0.4ex>[l]\ar@<0.4ex>[l] & G\times G\ar@[blue][d]\ar@<-0.6ex>[l]\ar[l]\ar@<0.6ex>[l] &\dots
% 	\ar@<-0.8ex>[l]\ar@<-0.4ex>[l]\ar@<0.4ex>[l]\ar@<0.8ex>[l]\\
% 	{*} & {*}\ar[l] & A\ar[l] & \dots\ar[l]
% }
\end{equation}
The blue map corresponds to the $2$-cocycle for the extension in ordinary group cohomology.

\begin{remark}
\label{cext_lie_alg}
	Differentiating a central extension of Lie groups produces a \emph{central extension of Lie algebras}
	\[
	0\to	\afrak \to \gtilde \to \g \to 0 \comma 
	\]
	which is what you would expect ($\afrak$ is an abelian Lie algebra contained in the center of
	$\gtilde$).

	Central extensions of Lie algebras are classified by second \emph{Lie algebra cohomology} $\HLie^2(\g;\afrak)$.  \index[notation]{HLie@$\HLie^\bullet$} \index[terminology]{Lie algebra!cohomology}
	Cocycles are alternating bilinear maps $\omega\colon\exterior^2\g\to\afrak$ satisfying a version of the Jacobi
	identity:
	\begin{equation}
	\label{LA_cext_Jacobi}
		\omega(X, [Y, Z]) + \omega(Y, [Z, X]) + \omega(Z, [X, Y]) = 0\period
	\end{equation}
	From such an $\omega$, we build a central extension which, as a vector space, is $\g\oplus\afrak$, but with
	Lie bracket
	\begin{equation}
	\label{liealgext}
		[(X_1, A_1), (X_2, A_2)] \colonequals [X_1, X_2] + \omega(X_1, X_2)\period
	\end{equation}
	A $1$-cochain is a map $\lambda\colon\g\to\afrak$, and its differential is $d\lambda(X,
	Y)\colonequals \lambda([X, Y])$.

	So we have a map
	\begin{equation*}
		\H^2(\BbulletG, \underline A)\to \HLie^2(\g; \afrak) \period
	\end{equation*}
	The van Est theorem says this is an
	equivalence in certain nice situations (not ours, unfortunately).
\end{remark}

%-------------------------------------------------------------------%
%-------------------------------------------------------------------%
%  The Virasoro algebra and the Virasoro group                      %
%-------------------------------------------------------------------%
%-------------------------------------------------------------------%

\subsection{The Virasoro algebra and the Virasoro group}

Let $\Gamma\colonequals \Diffplus(\Circ)$,  \index[notation]{Diff Pluss@$\Diffplus(\Circ)$}
the group of orientation-preserving diffeomorphisms of the circle. This is an
infinite-dimensional \emph{Fréchet Lie group}, meaning it is locally modeled on a Fréchet space and has a group
structure in which multiplication and inversion are smooth.

\begin{definition} \index[terminology]{Witt algebra} \index[notation]{Witt@$\wittRR$}
\label{Witt_algebra}
	The \emph{Witt algebra} $\wittRR$ is the infinite-dimensional real Lie algebra of polynomial vector fields on $\Circ$.
	Explicitly, it is generated by $\xi_n\colonequals -x^{n+1}\frac{\partial}{\partial x}$ for $n\in\ZZ$, with bracket
	\begin{equation}
		[\xi_m, \xi_n] \colonequals (m-n)\xi_{m+n} \period
	\end{equation}
\end{definition}

\noindent Skating over issues of regularity, the Witt algebra is the Lie algebra of $\Gamma$.\footnote{If we were to treat
regularity more carefully, we would allow some infinite linear combinations of the $\xi_n$, corresponding to the
Fourier series of a smooth vector field.}

The Virasoro algebra $\virRR$ is a central extension of $\wittRR$ by $\RR$.  \index[notation]{Vir@$\virRR$} \index[terminology]{Virasoro!algebra}
There is also a Virasoro group
$\widetilde\Gamma$, a central extension of $\Gamma$; the Virasoro algebra is its Lie algebra, and is easier to
define (since Lie algebra $\HLie^2$ just works to produce central extensions, whereas we had to modify group
cohomology). Specifically, consider the $2$-cocycle $c\colon \exterior^2\wittRR\to\RR$ given by
\begin{equation}
	c(\xi_m, \xi_n) \colonequals \frac{1}{12}(m^3-m)\delta_{m+n, 0}c \comma
\end{equation}
where $c$ is a chosen basis for $\RR$. The $1/12$ is not there for any deep reason, just as a normalization
constant. Anyways, as in~\eqref{liealgext} this defines for us an extension
\begin{equation*}
	1\to\RR\to\virRR\to\wittRR\to 1 \comma
\end{equation*}
called the \emph{Virasoro algebra}. The element $c$ inside $\virRR$ is called the \emph{central charge}.  \index[terminology]{Central Charge}

The Virasoro group $\widetilde\Gamma$ is the extension of $\Gamma$ by $\TT$ which is, as a space,
$\TT\times\Gamma$, with multiplication
\begin{equation}
	(z_1, f)\cdot (z_2, g) \colonequals (z_1 + z_2 + B(f, g), f\circ g) \comma
\end{equation}
where $B\colon \Gamma\times\Gamma\to\TT$ is the \emph{Bott cocycle}  \index[terminology]{Bott cocycle}
\begin{equation}
\label{bott_cocyc}
	B(f, g)\colonequals \oint_{\Circ} \log(f\circ g)' \mrm{d}(\log g)' \period
\end{equation}

\begin{remark}
	The identification $\Circ\cong\RRP^1$ embeds $\PGL_2^+(\RR) = \PSL_2(\RR)\subset\Gamma$ as the real fractional linear
	transformations; hence also
	\begin{equation*}
		\psl_2(\RR) = \slfrak_2(\RR)\subset\wittRR \comma
	\end{equation*}
	as the Lie algebra generated by
	$\xi_{-1}$, $\xi_0$, and $\xi_1$. Restricted to $\PSL_2(\RR)$, the Virasoro central extension is trivializable,
	which will be useful later.
\end{remark}

\begin{remark}
	Some authors' definitions will differ. 
	For example, defining the Witt and Virasoro algebras as complex Lie algebras, or
	defining the Virasoro group as the universal cover of ours. 
	%In particular, I used a few different sources when preparing this talk, and though I think they were consistent, it's possible that I missed something and they used different conventions.
	%%I (Minta) checked these and think they are all ok!
\end{remark}

\begin{remark}[(applications)]  \index[terminology]{Conformal Field Theory}
\label{virasoro_applications}
The Virasoro group and algebra appear in two-dimensional conformal field theory (CFT). Usually, in quantum field
theory, one specifies a (Riemannian or Lorentzian) metric on spacetime, and the information in the theory depends
on the metric. A \emph{conformal field theory} is a quantum field theory in which all information only depends on
the conformal class of the metric. Two-dimensional CFTs in particular connect to many areas of mathematics and
physics.
\begin{itemize}
	\item The mathematical formalization of 2d CFT, using vertex algebras, has connections to representation
	theory, and, famously, to monstrous moonshine.

	\item One way to think of string theory is as a 2d CFT on the worldsheet, one of whose fields is a map into
	(10- or 26-dimensional) spacetime.

	\item In condensed-matter physics, Wess--Zumino--Witten models (particular 2d CFTs) are used in modeling the
	quantum Hall effect. See also \Cref{ex-WZW}.

	\item Maybe closest to the hearts of the attendees of this seminar: the Stolz--Teichner conjecture suggests that
	cocycles for TMF on a space $X$ are given by families of 2d supersymmetric quantum field theories parametrized
	by $X$. Superconformal field theories are particularly nice examples of these, and have been used to shine
	light on this conjecture.\footnote{The appearances of SCFTs, rather than just CFTs, in superstring theory and
	in the Stolz--Teichner conjecture aren't as related to the Virasoro group and algebra; they have a larger
	symmetry algebra, though it's closely related.}
\end{itemize}
So how does the Virasoro appear in CFT? Let's suppose we're on a Riemann surface $\Sigma$ in a local holomorphic
coordinate $z$. If you write out commutators for the Lie algebra $\mathfrak c$ of infinitesimal conformal
transformations, you might notice they look like those for the Witt algebra --- in fact, if you complexify it, you
obtain precisely $\wittCC \oplus \wittCC $. So this acts on the system as a symmetry; you can think of it as two
different Witt group symmetries.

The fact that we obtain a central extension is standard lore from quantum mechanics. The state space in a quantum
system is a complex Hilbert space, but if $\lambda\in\bb{C}^\times$, the states $\vert\psi\rangle$ and $\lambda\vert\psi\rangle$ are
thought of as the same, in that measurements cannot distinguish them. Nonetheless, the formalism of quantum
mechanics uses the Hilbert space structure.

The takeaway, though, is that a symmetry of the system, as in acting on the states and all that, only has to be a
projective representation on the state space! So to describe an honest Lie group or Lie algebra acting on the state
space, we need to take a central extension of the symmetry group or Lie algebra. This leads us to the
(complexified) Virasoro algebra and Virasoro group. Thus, the symmetry algebra of conformal field theory is (at
least) a product of two copies of the Virasoro algebra, and the space of states is a representation of the Virasoro
algebra.
\end{remark}

%-------------------------------------------------------------------%
%-------------------------------------------------------------------%
%  Constructing the central extension with differential cohomology  %
%-------------------------------------------------------------------%
%-------------------------------------------------------------------%

\subsection{Constructing the central extension with differential cohomology}
\label{diffcoh_virasoro}

The key fact bridging differential cohomology and central extensions is:

\begin{lem}
\label{cext_lemma}
	There is an equivalence of simplicial sheaves $\ZZ(1)\simeq \Sigma^{-1}\underline{\TT}$.
\end{lem}

\begin{proof}
	By definition, $\ZZ(1)$ is the sheaf $0\to \ZZ\to\Omega^0\to 0$, and $\Omega^0 = \underline{\RR}$. The chain map
	\begin{equation*}
		\begin{tikzcd}
			0\arrow[r] & \ZZ\arrow[d]\arrow[r] & \underline{\RR}\arrow[d, "\mod{\ZZ}"] \arrow[r] & 0\\
			0\arrow[r] & 0\arrow[r] & \underline{\TT}\arrow[r] & 0
		\end{tikzcd}
	\end{equation*}
	is a quasi-isomorphism.
\end{proof}

\begin{cor}
\label{cext_corollary}
	For any Lie group $G$, possibly infinite-dimensional, we have an isomorphism
	\begin{equation*}
		\H^2(\BbulletG;\underline{\TT})\cong \H^3(\BbulletG;\ZZ(1)) \period
	\end{equation*}
	In particular, the group $ \H^3(\BbulletG;\ZZ(1)) $ classifies central extensions of $G$ by $\TT$ which are principal $\TT$-bundles
	over $G$.
\end{cor}

Thus, we would like to construct the Virasoro central extension via a differential cohomology class in
$\H^3(\BbulletGamma;\ZZ(1))$. 
This builds on the hard work of the previous few talks.
Let $\GL_n^+(\RR)$ denote the group of orientation-preserving, invertible $n\times n$ matrices. The restriction map
\begin{equation*}
	\H^*(\BGL_n(\RR);\RR)\to\H^*(\BGL_n^+(\RR);\RR)
\end{equation*}
is not an isomorphism, but is close enough to one that we still
have Pontryagin classes for
$\GL_n^+(\RR)$-bundles, or for oriented vector bundles.

In \Cref{2n_to_n}, Hopkins described how $\H^4(\BbulletGLp{n}{\RR};\ZZ(2))$ fits into a pullback square
\begin{equation}\label{diffp1square}
	\begin{tikzcd}
		\H^4(\BbulletGLp{n}{\RR};\ZZ(2)) \arrow[r, blue] \arrow[d] & \H^4(\BGL_n^+(\RR);\ZZ) \arrow[d]\\
		\H^4(\BbulletGLp{n}{\RR};\RR(2)) \arrow[r, red] & \H^4(\BGL_n^+(\RR);\RR) \period
	\end{tikzcd}
\end{equation}

\begin{definition}
	An \emph{off-diagonal differential lift of $p_1$} is a class
	\begin{equation*}
		\ptilde_1\in \H^4(\BbulletGLp{n}{\RR};\ZZ(2))
	\end{equation*}
	whose image under the blue map is the usual $p_1\in \H^4(\BGL_n^+(\RR);\ZZ)$.
\end{definition}

By \Cref{cor-CWoff}, we have an isomorphism
\begin{equation*}
	\H^4(\BbulletG; \RR(2))\cong\Sym^2(\gdual)^G \period
\end{equation*} 
For $\GL_n^+(\RR)$, this is
an $\RR^2$, spanned by the invariant polynomials $\Tr(A)^2$ and $\Tr(A)^2$, which we call $c_1^2$ and $c_2$,
respectively. 
The group $\H^4(\BG;\RR)$ can be dispatched with ordinary Chern--Weil theory: we repeat the same story, but
retracting $G$ onto its maximal compact. 
Here, we get
\begin{equation*}
	\H^4(\BSO(n);\RR)\cong\RR \comma
\end{equation*}
spanned by $\Tr(A^2)$, as $\Tr(A)^2 = 0$. 
Accordingly, the red map in~\eqref{diffp1square} is a rank-$1$ map $\RR^2\to\RR$. Since~\eqref{diffp1square} is
a pullback square, there is an $\RR$ worth of differential lifts of $p_1$: explicitly, $\lambda\in\RR$ gives you the
lift of $p_1$ which maps in the lower left to $(1/2)(\lambda c_1^2 - 2c_2)$. However, if you want
\begin{equation*}
	\ptilde_1(E_1\oplus E_2) = \ptilde_1(E_1) + \ptilde_1(E_2) \comma
\end{equation*}
you force $\lambda = 1$, which is a quick calculation with
the Whitney formula. 
(All this was in \Cref{sec:Delignecup}.)\footnote{As we have mentioned a few times already, the Whitney sum formula
is not true for Pontryagin classes in $\ZZ$-valued cohomology. So what's going on here? Brown~\cite[Theorem
1.6]{Bro82} identified the obstructions to the Whitney sum formula holding for the Pontryagin classes of a pair of
vector bundles, and for $p_1$ the obstructions vanish for oriented vector bundles. This is one reason to use
$\GL_n^+(\RR)$ instead of $\GL_n(\RR)$.}

In \cref{FiberIntegration}, we also discussed the fiber integration map for an $\Hhat $-oriented fiber bundle
\begin{equation*}
	F\to E\to B \comma  
\end{equation*}
which has the form
\begin{equation*}
	\H^k(E;\ZZ(\ell))\to \H^{k-\dim(F)}(B;\ZZ(\ell-\dim(F))) \period
\end{equation*}
Combining all this, consider the universal
oriented sphere bundle $\EbulletGamma\times_\Gamma \Circ\to \BbulletGamma$, which is an $\Hhat $-oriented fiber
bundle with fiber $\Circ$. 
Therefore, given a differential lift of $p_1$, we can apply it to the vertical tangent
bundle $V\to \EbulletGamma\times_\Gamma \Circ$, and get a class
\begin{equation*}
	\ptilde_1(V)\in \H^4(\EbulletGamma\times_\Gamma \Circ;\ZZ(2)) \period
\end{equation*}
Then we can push it forward to a class in $\H^3(\BbulletGamma;\ZZ(1))$, which determines an
isomorphism class of central extensions of $\Gamma$ as above. The goal is to determine the choice of $\lambda$ such
that this central extension gives the Virasoro group. I'll suggest some ways forward.

The first thing we need is a way to get a handle on the group of extensions of $\Gamma$. Recall that
$\PSL_2(\RR)\subset\Gamma$ as the real fractional linear transformations of $\RRP^1 = \Circ$; hence a central extension
of $\Gamma$ restricts to a central extension of $\PSL_2(\RR)$.

\begin{theorem}[{(Segal \cite[Corollary 7.5]{Seg81})}]
	A central extension of $\Gamma$ by $\TT$ is determined by the pair of (1) its restriction to $\PSL_2(\RR)$ and (2)
	the induced Lie algebra central extension of $\wittRR$ by $\RR$. 
	Said differently, there is an isomorphism of abelian
	groups
	\begin{equation*}
		\Cent_{\TT}(\Gamma) \isomorphism \Cent_{\TT}(\PSL_2(\RR))\times\Cent_{\RR}(\wittRR) \period
	\end{equation*}
\end{theorem}

We can identify both of these groups. First, $\uppi_1\PSL_2(\RR)\cong\ZZ$, and the universal cover
$\SLtilde_2(\RR)\to\PSL_2(\RR)$\footnote{The notation is because it's also the universal cover of $\SL_2(\RR)$,
which is the connected double cover of $\PSL_2(\RR)$.} is the \emph{universal central extension} of $\PSL_2(\RR)$:
for any abelian group $A$, central extensions of $\PSL_2(\RR)$ by $A$ are in bijection with maps $\varphi\colon \ZZ\to
A$, given by
\begin{equation}
	\begin{tikzcd}
		0\arrow[r] & \ZZ\arrow[r]\arrow[d, "\varphi"'] & \SLtilde_2(\RR)\arrow[d]\arrow[r] & \PSL_2(\RR)\arrow[r] \arrow[d, equals] & 0\\
		0\arrow[r] & A\arrow[r] & (\SLtilde_2(\RR))_\varphi\arrow[r] & \PSL_2(\RR)\arrow[r] & 0 \period
	\end{tikzcd}
\end{equation}
So
\begin{equation*}
	\Cent_{\TT}(\PSL_2(\RR))\cong\Hom(\ZZ, \TT) = \TT \period
\end{equation*} 
The computation that $\H^2(\wittRR;\RR)\cong\RR$ is standard, e.g. \cite[\S 6.2.1]{Obl17}.

Thus the map from off-diagonal differential lifts of $p_1$ to central extensions of $\Gamma$ is a map $\RR\to\RR\times\TT$. 

One can then ask the following question, 
which was posed to us by Dan Freed and Mike Hopkins.
\begin{question}
	\label{VirasoroQuestion}
	Does there exist an off-diagonal differential lift $\ptilde_1$ of the first Pontryagin class that hits the Virasoro algebra central extension in $\RR\times\TT$?
\end{question}

Note that the Virasoro central extension is in the first factor of $\RR\times\TT$.  
Indeed, it induces the Virasoro algebra central
extension, and hence is nontrivial on the first factor and trivial when restricted to $\PSL_2(\RR)$.

The answer to \Cref{VirasoroQuestion} is: yes!
\begin{theorem}[\cite{DLW21}]
	Let $\ptilde_1\in \H^4(\Bun_{\GL_n^+(\RR)};\ZZ(2))$ be the unique off-diagonal differential lift of $p_1$ which
	satisfies the Whitney sum formula. The central extension of $\Gamma$ defined by transgressing $\ptilde_1$
	according to the procedure above is the Virasoro group with central charge $-12$.
\end{theorem}
``Central charge $-12$'' means that this group extension is classified by $-12$ times the Bott
cocycle~\eqref{bott_cocyc}.

\newpage

%-------------------------------------------------------------------%
%  Applications                                                     %
%-------------------------------------------------------------------%

\part{Applications}

\numberwithin{equation}{part}
%!TEX root = ../diffcoh.tex

%-------------------------------------------------------------------%
%-------------------------------------------------------------------%
%  Part III Introduction                                            %
%-------------------------------------------------------------------%
%-------------------------------------------------------------------%

\label{part:applications}\label{applications_part}

In \cref{applications_part}, we survey some applications of differential cohomology to questions in geometry and
physics. Some of these applications belong to the pattern that what ordinary cohomology tells us about topological
objects, differential cohomology tells us about their geometric analogues: this includes both the use of
differential cohomology to obstruct conformal immersions as well as the classification of invertible field
theories, both of which we say more about below. For other applications, the analogy with ordinary cohomology is
subtler; some use the differential characteristic classes we built in \cref{char_class_part}, such as the study of
loop groups and the Virasoro group.

\subsection*{Chern--Simons invariants}
Chern--Simons invariants, which we define and study in \cref{config_spaces}, are the key to many of these
applications. Let $G$ be a compact Lie group; choose a class $\lambda\in \H^4(\BG;\ZZ)$ and let $\langle-,
-\rangle$ be the degree-$2$ $G$-invariant symmetric polynomial on $\g$ associated to the image of $\lambda$
in de Rham cohomology. The Chern--Simons invariant associated to $\lambda$ is defined for a $3$-manifold $Y$, a
principal $G$-bundle $\pi\colon P\to Y$, and a connection $\conn$ on $P$ with curvature $\curvature{\conn}$. If we
assume that $\pi$ has a section, so that we can descend $\curvature{\conn}$ to a form on $Y$, the Chern--Simons
invariant is
\begin{equation}
\label{CS_3_intro}
	\CS_\lambda(P, \conn) = \int_Y \langle \conn\wedge\curvature{\conn}\rangle - \frac 16\langle \conn\wedge [\conn,
	\conn]\rangle \in \RR/\ZZ \period
\end{equation}
We first met this invariant in a different guise in \cref{secondary_chern_simons}. In \cref{differential_CW_lift}
we showed $\lambda$ and $\langle-, -\rangle$ determine a differential refinement
$\hat\lambda\in\Hhat^4(\BunGnabla ; \ZZ)$, and said but did not prove that the Chern--Simons invariant is the
secondary invariant associated to $\hat\lambda$. We will prove the latter fact in \cref{config_spaces}.

Chern--Simons invariants and their generalizations play a central role in most of the applications of differential
cohomology which we survey: they bridge the geometry of connections with the algebraic topology of (differential)
characteristic classes, and therefore have something to say about both worlds.

% config spaces
For example, in \cref{config_spaces} we follow Evans-Lee--Saveliev \cite{deletedsquare} and use Chern--Simons
invariants as a tool to determine when two homotopy-equivalent lens spaces are not
diffeomorphic; to do so, we also spend time developing a little of the theory of Chern--Simons invariants. The
classification of lens spaces up to diffeomorphism or homotopy equivalence is classical \cites{Rei35}[\S
5]{Whi41}{Bro60}, which makes it a good testing ground to determine how powerful manifold invariants are. For
example, Longoni--Salvatore \cite{LS05} proved the surprising result that the homotopy type of the two-point
configuration space of a lens space can distinguish homotopy-equivalent lens spaces.  Evans-Lee--Saveliev build on
Longoni--Salvatore's work, providing more comprehensive tools for understanding when the homotopy type of the
two-point configuration space of $L(p, q)$ is a stronger invariant than the homotopy type of $L(p, q)$. They extend
Chern--Simons invariants to two-point configuration spaces and use them to give a numerical criterion
(\cref{f_h_restr}) for a map of configuration spaces to be a homotopy equivalence. They combine this criterion with
a few other tools, including Massey products, to provide many pairs of homotopy-equivalent lens spaces whose
two-point configuration spaces are not homotopy equivalent.

% conformal immersions
In \cref{conformal_immersions}, we use on-diagonal differential characteristic classes to obstruct conformal
immersions, following Chern--Simons \cite{cs}. Recall that characteristic classes in ordinary cohomology can
obstruct immersions into $\RR^n$ as follows: if $M$ is a smooth $m$-manifold that immerses into $\RR^n$ with normal
bundle $\nu$, then $TM\oplus\nu \cong T\RR^n|_M\cong\underline\RR^n$, and $\nu$ is rank $n-m$, so all of its
characteristic classes in degree greater than $n-m$ vanish. This places constraints on the characteristic classes
of $M$. For example, let $w_i$ denote the $i$-th Stiefel--Whitney
class;\index[terminology]{Stiefel--Whitney class} if $\CCP^2$ immersed in $\RR^5$, then the normal bundle $\nu$
would be
one-dimensional, so
\begin{equation}
	w_2(T\CCP^2\oplus\nu) = w_2(T\CCP^2) + \underbracket{w_1(T\CCP^2)w_1(\nu) + w_2(\nu)}_{=0} =
	w_2(\underline\RR^5) = 0 \comma
\end{equation}
but $w_2(T\CCP^2)\ne 0$, which prevents such an immersion. One can run the same argument using Cheeger--Simons'
differential characteristic classes, which we discussed in \cref{DifferentialCharacteristicClasses}: since these
characteristic classes are defined for vector bundles with connection, they can obstruct isometric embeddings of a
Riemannian manifold $M$ by placing constraints on $TM$ with its Levi-Civita
connection.\index[terminology]{Levi-Civita connection} Chern--Simons \cite{cs} improve on this argument in two
ways, giving it considerably more power: they prove that the on-diagonal differential Pontryagin classes of the
Levi-Civita connection only depend on the conformal class of the metric (\cref{diff_p_conformally_invariant}), so
can be used to obstruct conformal immersions. They then use the Chern--Simons form to obtain additional
obstructions: in some cases, the Chern--Simons form is closed, and conformal immersions restrict what its de Rham
class can be. The proofs of these obstructions make use of the close relationship between differential
characteristic classes and Chern--Simons forms. Chern--Simons' obstructions are strong enough to prove that
$\RRP^3$ with the round metric cannot conformally immerse in $\RR^4$ (\cref{nope_RP3}).

Our third application of Chern--Simons invariants is to physics: there is a classical field theory whose Lagrangian
is the Chern--Simons invariant~\eqref{CS_3_intro}. We discuss this theory in \cref{classical_CS}, focusing on how
various pieces of the theory can be described using differential cohomology. Schwarz \cite{Sch77} and
Witten \cite{Wit89} quantized this theory, producing a topological field theory called Chern--Simons
theory\index[terminology]{Chern--Simons theory} which has been a major object of study in both mathematics and
physics. See \cref{quantum_CS} for references and more information on the quantum theory.

\subsection*{Quantum physics}
Speaking of physics, several of the applications of differential cohomology that we survey are in physics or are
closely related to it. In these applications, differential cohomology tends to appear because quantization imposes
integrality conditions on objects in field theories; in many cases these can be lifted to integrality data,
allowing differential cohomology to enter the picture.

\Cref{field_theory} is dedicated to this idea, working with the example of electromagnetism. We first discuss
classical Maxwell theory, describing how information in this theory can be expressed with differential forms. Then
we walk through Dirac's argument \cite{Dir31} that the presence of magnetic monopoles forces electric and magnetic
charges to be quantized, i.e.\ valued in a discrete subgroup of $\RR$. As a consequence, the fields in the quantum
theory are cocycles for differential cohomology, and the action can be rewritten using the
differential-cohomological cup product and integration. For electromagnetism, the appearance of differential
cohomology is relatively explicit and simple, making it a good example, but the concept of quantization of abelian
gauge fields leading to differential cohomology appears in numerous other places in quantum physics, and can
involve fancier objects such as differential \Ktheory.

The next chapter, \cref{invertible_field_theories}, is about a different application of differential (generalized)
cohomology to physics: the classification of invertible field theories. This is one of the applications which is a
geometric analogue of a use of ordinary (generalized) cohomology for something topological. Following Atiyah and
Segal, a topological field theory (TFT) is a symmetric monoidal functor
\[Z\colon\Bord_n\to\mathsf C\comma\]
where $\Bord_n$ is a bordism (higher) category and $\mathsf C$ is some symmetric monoidal (higher) category, often
$\mathsf{Vect}_{\CC}$. The simplest nontrivial TFTs are the invertible TFTs, which are the TFTs whose values on all
objects and morphisms in $\Bord_n$ are invertible in $\mathsf C$, meaning that objects are invertible under the
tensor product, and morphisms are invertible under composition.
%Invertible TFTs were classified by
%Freed--Hopkins--Teleman \cite{FHT10} using the work of Galatius--Madsen--Tillmann--Weiss \cite{GMTW09} and
%Nguyen \cite{Ngy17} characterizing the stable homotopy types of bordism categories.
We are interested in \emph{reflection-positive} invertible TFTs; this extra requirement is a physically motivated
version of unitarity. The classification of reflection-positive invertible TFTs is due to
Freed--Hopkins \cite{FH21}, who show that, up to isomorphism, reflection-positive invertible TFTs are classified by
the torsion subgroup of $[\mathrm{MTH}, \Sigma^n \mathrm I_\ZZ]$ (see \cref{top_IFT} for definitions of these
spectra). In typical examples, the partition functions of these theories are bordism invariants defined by
integrating characteristic classes in (generalized) cohomology. Freed--Hopkins (\textit{ibid.}) go further and
conjecture that the entirety of $[\mathrm{MTH}, \Sigma^n \mathrm I_\ZZ]$ classifies invertible field theories that
need not be topological, which would be defined on some yet-to-be-constructed geometric bordism category. Again,
partition functions can often be described by integrating characteristic classes, but this time in differential
(generalized) cohomology, and typically in one dimension lower, so as to obtain a secondary invariant. We discuss
this conjecture and several examples: classical Chern--Simons theory as mentioned above, the classical
Wess--Zumino--Witten model, and an example using differential $\mathrm{KO}$-theory.

\subsection*{Representations of loop groups}

In \cref{loop_groups}, we turn to the representation theory of loop groups. These are infinite-dimensional Lie
groups, but unusually nice ones: as long as you are careful about what you mean by a representation, their
representation theory closely resembles that of compact Lie groups! The representations we care about are
projective representations, so genuine representations of a central extension by the circle group
\begin{equation}
\label{Tcent_}
	1\to \TT\to \LGtilde\to \LG\to 1 \comma
\end{equation}
satisfying a ``positive energy'' condition: restricting the representation to $\TT$, its weight subspaces for
negative weights are trivial. The reader may wonder how invariant this definition is, and is right to be concerned:
it is a significant theorem of Pressley--Segal \cite[Theorem 13.4.2]{loop} that when $G$ is simply connected and
compact, every positive energy representation of $\LG$ admits an intertwining projective $\Diffplus(\TT)$-action,
meaning that the notion of positive energy is preserved under reparametrizations of $\TT$. One of the major goals
of \cref{loop_groups} is to discuss the key ideas in this theorem and its proof: we introduce and motivate the
positive energy condition, we discuss the nice properties of positive energy representations, and we sketch the
proof of Pressley--Segal's theorem. Along the way, we discuss some connections with physics. In \cref{PS_diffcoh},
we discuss two different connections to differential cohomology: first, the central extensions of the sort we
consider are principal $\TT$-bundles over $\LG$, hence determine classes in $\H^2(\LG;\ZZ)$. It turns out that
every element of this cohomology group comes from a central extension, and moreover, as principal $\TT$-bundles
they carry canonical connections, allowing for a lift to $\Hhat^2(\LG;\ZZ)$. This class is related to the ``level''
that one starts with via transgression maps
\begin{equation*}
	\Hhat^4(\BunGnabla ;\ZZ)\to\Hhat^3(G;\ZZ)\to\Hhat^2(\LG;\ZZ) \period
\end{equation*}
Central
extensions that are principal $\TT$-bundles correspond to off-diagonal classes in the differential cohomology group $\H^3(\Bun_{\LG};
\ZZ(1))$, as in \cref{VirasoroAlgebra}, and we say a little about this perspective too.

Our final chapter, \cref{segal_sugawara}, takes the above story and makes it explicit, albeit at the level of Lie
algebras. The Lie algebra of a central extension $\LGtilde$, denoted $\Lgtilde{}$, is an example of a
\emph{Kac--Moody algebra}, and is a central extension of the loop algebra of the Lie algebra of $\g$. The
Pressley--Segal theorem cooks up an intertwining projective $\Diffplus(\TT)$-action on the representations of
$\LGtilde$, so at the level of Lie algebras we might expect a compatible Virasoro algebra action on the
representations of Kac--Moody algebras. This is true, and Segal--Sugawara show we can do better, explicitly
identifying how the central $\CC$ in the Virasoro algebra acts in terms of the level of the central
extension~\eqref{Tcent_}. Both this chapter and the previous chapter on loop groups are closely related to
two-dimensional conformal field theory: the data of the category of positive-energy representations of $\LGtilde$
can be used to build a two-dimensional conformal field theory called the Wess--Zumino--Witten model. This CFT is
further related to Chern--Simons theory, a 3d TFT. All of this data --- the central extension of $\LG$, the
specific Wess--Zumino--Witten model, the specific Chern--Simons theory --- is indexed by groups such as
$\H^2(\LG;\ZZ)$, $\H^3(G;\ZZ)$, and $\H^4(\BG;\ZZ)$, which when $G$ is simple and simply connected are all
canonically isomorphic to $\ZZ$. These groups are related to each other by transgression maps, and this corresponds
to the relationship between, e.g.\ loop groups and the WZW model, or the WZW model and Chern--Simons theory. These
cohomology classes have differential refinements, as do the transgression maps relating them.

These are not the only applications of differential cohomology to topology, geometry, or physics, but we hope they
illustrate the diversity of things that can be done with differential cohomology, and that they make for an
interesting and enjoyable read.

\numberwithin{equation}{subsection}
\newpage
%!TEX root = ../diffcoh.tex

%-------------------------------------------------------------------%
%-------------------------------------------------------------------%
%  Chern–Simons invariants                                          %
%-------------------------------------------------------------------%
%-------------------------------------------------------------------%

\section{Chern--Simons invariants}
\textit{by YiYu (Adela) Zhang}
\label{config_spaces}

%-------------------------------------------------------------------%
%-------------------------------------------------------------------%
%  Motivation/Review                                                %
%-------------------------------------------------------------------%
%-------------------------------------------------------------------%

Our first application is to the theory of Chern--Simons forms and invariants, tools in geometry which are closely
tied to differential cohomology. We first mentioned these in \Cref{secondary_chern_simons}, where we said that
Chern--Simons invariants are defined to be the secondary invariants associated to the on-diagonal differential
characteristic classes constructed in Chern--Weil theory. But they also have a much more geometric description,
given by integrating a specific form built from the connection and curvature forms. These two descriptions are part
of the reason Chern--Simons invariants are so useful: one can use homotopy-theoretic methods in differential
cohomology to learn facts about geometry, and vice versa. This will be a common theme throughout this part of the
book, and Chern--Simons forms will appear several times.

We begin in \cref{ssec:lens_CS}, defining and discussing Chern--Simons forms associated to a principal bundle
$\pi\colon P\to M$ with connection, and relating them to the differential lifts of Chern--Weil characteristic
classes from \cref{DifferentialCharacteristicClasses}. In \cref{cs_invariants}, we focus on the case when $\pi$ is
a principal $\SU_2$-bundle over a $3$-manifold, where we can descend the Chern--Simons invariant from an integral
on $P$ to an integral on $M$.  Finally, in \cref{config_ssec}, we show an application of Chern--Simons invariants,
as a tool to determine when two-point configuration spaces of lens spaces are homotopy equivalent.

\subsection{Chern--Simons forms}
\label{ssec:lens_CS}

Let $G$ be a compact Lie group and $\pi\colon P\rightarrow M$ a principal $G$-bundle. 
Fix a degree-$k$ invariant polynomial $f \in\Sym^k(\gdual)^G$. Given a connection $\conn$ on $P$ with
curvature $\curvature{\conn}$, we will write $f(\curvature{\conn})\in \HdR^{2k}(M)$ for the associated Chern--Weil form.
%We will write $f (A)=f (F_A)$ for a connection $A$.

\begin{recall}
\label{recall_connection}
  A \textit{connection} $\conn$ on the principal $G$-bundle $\pi\colon P\to M$ is a $\g$-valued $1$-form on $P$
  which is $G$-equivariant in the sense that $(R_g)^*\conn =\Ad_{g^{-1}}\conn$, and it
  is ``the identity'' on tangent vectors along the fiber, i.e.\ $\conn(X_{\xi})=\xi$ for $\xi\in\g$ and $X_{\xi}$ its
  fundamental vector field.

  The \textit{curvature} of $\conn$, which we usually denote $\curvature{\conn}$, is the form $\d\omega + [\omega,
  \omega]\in\Omega^2(M;\g)$.
\end{recall}
\index[terminology]{connection!on a principal $G$-bundle}
Analogous to connections on vector bundles, a $G$-connection corresponds to a splitting
\begin{equation*}
    \Tan P\cong H\oplus V \comma
\end{equation*}
where
$V$ is the vertical tangent bundle (the kernel of $\pi_*\colon \Tan P\to \Tan M$), and $H$ is the horizontal tangent
bundle. A priori, there is only a short exact sequence
\index[terminology]{vertical tangent bundle}
\index[terminology]{horizontal tangent bundle}
\begin{equation}
	\begin{tikzcd}[ampersand replacement=\&]
        0 \& V
        \&  \Tan P
        \& H
        \& 0 \semicolon
        \arrow[from=1-1, to=1-2]
        \arrow[from=1-2, to=1-3]
        \arrow[from=1-3, to=1-4]
        \arrow[from=1-4, to=1-5]
\end{tikzcd}
\end{equation}
a connection is a $G$-equivariant splitting. Because the fibers of a principal $G$-bundle are $G$-torsors,
there is an isomorphism $V\cong\underline{\g}$, and the $G$-action is the fiberwise adjoint action, leading to the
definition of connection given in \ref{recall_connection}.

Recall from \cref{sssec:CW_principal} that the \textit{adjoint bundle} to a principal $G$-bundle $P\to M$, denoted
$\g_P$, is the associated vector bundle to the adjoint representation $G\to\Aut(\g)$. The affine space of
connections on $P$ can be identified with $\cA_P=\Omega^1(M;\g_P)$, i.e., $ 1 $-forms on $M$ with values in the
adjoint bundle. Given two connections $\conn_0$, $\conn_1\in\cA_P$, the straight-line path
$\conn_t:I\rightarrow \cA_P$ determines a connection $\overline\conn$ on the $G$-bundle $P\times [0,1]$ over
$M\times[0,1]$. Let $\curvature{\overline\conn}$ be the curvature of $\overline\conn$.
\begin{definition}
The \textit{Chern--Simons form} associated to $\conn_0,\conn_1\in\cA_P$ and $f $ is given by
\index[terminology]{Chern--Simons form}
\index[notation]{CSf@$\CS_f(\conn)$}
\[\CS_{f }(\conn_1, \conn_0)=\int_{[0,1]}f(\curvature{\overline\conn})\in\Omega^{2k-1}(M)\period\]
\end{definition}
Let $\curvature{\conn_i}$ denote the curvature of $\conn_i$; then, by Stokes' theorem,
\begin{equation}
\label{stokes_CS}
\d\CS_{f }(\conn_1, \conn_0)= f(\curvature{\conn_1})-f(\curvature{\conn_0}) \period
\end{equation}
That is, the de Rham class $[f(\curvature{\conn_1})]$ is independent of the choice of connection, a fact that we first saw in
\cref{ChernWeilTheory}.
\begin{remark}
The path from $\conn_0$ to $\conn_1$ matters --- if we choose a different path, the Chern--Simons form will
differ by an exact term. This is beyond the scope of this chapter.
\end{remark}
Suppose instead we take the $G$-bundle $\pi^*P\rightarrow P$, which has a tautological section and hence a
tautological (flat) connection $\conn_0$. Then we can define a Chern--Simons form on $P$ (not on $M$!) for a single
connection $\conn$:
\begin{equation}
\label{total_space_CS}
\CS_{f }(\conn)=\CS_{f }(\pi^*\conn, \conn_0)\in \Omega^{2k-1}(P) \period
\end{equation}
Since $\conn_0$ is flat,~\eqref{stokes_CS} implies
\begin{equation}
\label{chern_simons_differential}
	\d\CS_{f }(\conn )= f(\pi^*\curvature{\conn})=\pi^*f(\curvature{\conn}) \period
\end{equation}
At this point, we want you to recall the differential cohomology hexagon from \cref{thm:SimonsSullivanunique}.
%
%Chern and Simons (\cite{cs}) showed that when $[f (A)]$ is an integral class, then there is a $u\in C^{2k-1}(M,\RR/\ZZ)$ such that $\pi^*u$ is the reduction of $\CS_{f }(A)\mod \ZZ$.
%
%-------------------------------------------------------------------%
%-------------------------------------------------------------------%
%  Relation to Cheeger–Simons differential characters               %
%-------------------------------------------------------------------%
%-------------------------------------------------------------------%
%
%\subsubsection{Relation to Cheeger--Simons differential characters}
%
%Now we will briefly explain the relation between Chern--Simons forms and Cheeger--Simons differential cohomology, following \cite[Chapter 2]{b}.
%
%
%Recall the definition of the Cheeger--Simons differential cohomology
%$$\Hhat^{k}(M;\ZZ)=\big\{ \chi\colon \Zsm_{k-1}\rightarrow \RR/\ZZ\big| \exists \alpha\in\curvature{\conn}^{k}(M)_{\ZZ},
%\chi(\partial c)=\int_c\alpha\mod\ZZ\big\}$$
%and the differential cohomology diagram. (c.f. Peter's introductory talk.)
\index[terminology]{differential cohomology hexagon}
\begin{equation}
\begin{gathered}
    \begin{tikzcd}[column sep={13ex,between origins}, row sep={11ex,between origins}]
      0 \arrow[dr] & & & & 0 \\
      & \H^{*-1}(M;\RR/\ZZ) \arrow[rr, "-\Bock"] \arrow[dr] & & \H^*(M;\ZZ) \arrow[dr] \arrow[ur] & \\
      \HdR^{*-1}(M) \arrow[ur] \arrow[dr] & & \Hhat^*(M;\ZZ) \arrow[ur, "\ch" description] \arrow[dr, "\curv" description] & & \HdR^*(M) \\
      & \displaystyle\frac{\Omega^{*-1}(M)}{\Omegacl^{*-1}(M)_{\ZZ}} \arrow[rr, "\d"'] \arrow[ur,
	  "\iota" description] & & \Omegacl^*(M)_{\ZZ} \arrow[ur] \arrow[dr] & \\
      0 \arrow[ur] & & & & 0 
    \end{tikzcd}
\end{gathered}
\end{equation}
The squares and triangles are commutative, and the diagonals are short exact sequences.
%
%The curvature map\index[terminology]{curvature map}
%\begin{equation*}
%  \curv \colon \Hhat^{k}(M;\ZZ)\rightarrow \curvature{\conn}^{k}(M)
%\end{equation*}
%sends $\chi$ to $\alpha$. 
%The characteristic class map\index[terminology]{characteristic class map}
%\begin{equation*}
%  \ch\colon\Hhat^{k}(M;\ZZ)\rightarrow \H^{k}(M;\ZZ)
%\end{equation*}
%is obtained by lifting $\chi$ to $\tilde{\chi} \colon \Zsm_{k-1}\rightarrow\RR$ and sending $\chi$ to the integral class defined by
%\begin{equation*}
%  c\mapsto -\tilde{\chi}(\partial c)+\int_c \alpha \,, \quad \text{ for } c\in \Csm_{k}(M;\ZZ) \period
%\end{equation*}
%We specifically need $\iota$, which is the descent of the map $\curvature{\conn}^{k-1}(M)\rightarrow\Hhat^{k}(M,\ZZ)$ defined
%by
%\begin{equation*}
%  \iota(\omega)(z) \colonequals \exp(2\pi i\int_z\omega)
%\end{equation*}
%for $z\in \Zsm_{k-1}$. Since the Chern--Simons form is a closed $(2k-1)$-form on $P$, we can ask about its image
%under $\iota$.
\begin{proposition}
  \label{iota_chern_simons}
  Suppose $c^\ZZ\in\H^{2k}(\BG;\ZZ)$ is an integral lift of the Chern--Weil characteristic class of $f$ and
  $\chat\in\Hhat^{2k}(\BunGnabla ; \ZZ)$ is the differential refinement of $c^\ZZ$ and $f$ guaranteed by
  \cref{differential_CW_lift}. Then for any principal $G$-bundle $\pi\colon P\to M$ with connection $\conn$,
  \begin{equation}
  	\iota(\CS_f(\conn)) = \pi^*\chat(P, \conn)\in \Hhat^{2k-1}(P;\ZZ).
  \end{equation}
\end{proposition}

\begin{proof}
  As usual, we can prove this for all principal bundles with connection at once by working universally on $\pi\colon
  (\BunGnablatriv, \conn)\to \BunGnabla $; here $\BunGnablatriv$ denotes the classifying stack of trivialized principal $G$-bundles
  with connection (see \Cref{ntn:BunGnablatriv}) and $\pi$ is the universal principal $G$-bundle with connection $A$ in the setting of stacks on
  $\Mfld$. By construction, if $\curvature{\conn}$ denotes the curvature of $\conn$,
  \begin{equation*}
    \curv(\chat) = f(\curvature{\conn})\in\Omega^{2k}(\BunGnabla ) \comma
  \end{equation*}
  so by~\eqref{chern_simons_differential},
  \begin{equation*}
  	\d\CS_f(\conn) = \pi^*\curv(\chat)\in\Hhat^{2k}(\BunGnablatriv;\ZZ).
  \end{equation*}
  The hexagon does all the hard work for us: the homotopification of $\BunGnablatriv$ is contractible (\Cref{ex:homotopification_of_Omega^n}), so
  \begin{equation*}
    \H^{2k-1}(\BunGnablatriv; \RR/\ZZ) = 0
  \end{equation*}
  and thus the curvature map is injective.  Since $\d = {\curv}\circ\iota$, we
  can conclude.
\end{proof}

Now suppose that $\pi\colon P\rightarrow M$ admits a section $\sigma\colon M\rightarrow P$. 
Then we further deduce that
\begin{equation*}
	\chat(P, \conn) = \sigma^*\pi^*\chat(P, \conn) = \iota(\sigma^*\CS_f(\conn)) \comma
%  \CWhat_{A}(f ,u)=\sigma^*(\pi^*\CWhat_{A}(f ,u))=\iota(\sigma^*\CS_{A}(f )) \period
\end{equation*}
meaning that
\begin{equation}
\label{CS_secondary_CW}
	\int_M \chat(P, \conn) = \int_M \sigma^*(\CS_f(\conn))\in\RR/\ZZ \period
\end{equation}
That is, as promised in \cref{secondary_chern_simons}, this Chern--Simons invariant is the secondary invariant
associated to $\chat$.
\index[terminology]{secondary invariant}
%
%Then (the generalization of) the \textit{Chern--Simons functional} is the evaluation at the fundamental class
%\begin{equation}
%  \CWhat_{A}(f ,u)([M])=\exp(2\pi i\int_M\sigma^*\CS_{A}(f )) \comma
%\end{equation}
%or its more familiar form 
%\begin{equation*}
%  \log(\CWhat_{A}(f ,u)([M]))=\int_M\sigma^*\CS_{A}(f ) \mod \ZZ \period
%\end{equation*}
This is conceptually nice, but how do we obtain computable topological invariants from this formula?

%-------------------------------------------------------------------%
%-------------------------------------------------------------------%
%  Chern–Simons invariants for 3-manifolds                          %
%-------------------------------------------------------------------%
%-------------------------------------------------------------------%

\subsection{Chern--Simons invariants for 3-manifolds}
\label{cs_invariants}
\index[terminology]{Chern--Simons invariant!of a $3$-manifold}

As an example, we examine the case where $P$ is a principal $\SU_2$-bundle over a path-connected $3$-manifold $M$,
\begin{equation*}
  f (\conn)=\frac{1}{8\pi^2}\tr(\curvature{\conn}\wedge \curvature{\conn}) \comma
\end{equation*}
and $c^\ZZ\in\H^4(\BSU_2;\ZZ)$ is the
second Chern class.\index[terminology]{Chern class}
%is the second Chern class in $\H^4(M;\RR)$. This is the classical Chern--Simons theory.
We mostly follow the exposition in \cite{KK}.

The quaternionic projective space $\HHP^\infty$ is a $\BSU_2$, so $\BSU_2$ is $3$-connected; hence every principal
$\SU_2$-bundle over a $3$-manifold is trivializable. Fix a trivialization; then there is a trivial (flat)
connection $\conn_0$, which allows us to identify $\cA_P$ with $\curvature{\conn}^1(M; \mathfrak{su}_2)$. Recall
that $\SU_2$ acts on $\cA_P$ by
\[g\cdot \conn =g\conn g^{-1}-dg\ g^{-1}\period\]
This action preserves flatness: if $\curvature{\conn}$ is the curvature of $\conn$, then the curvature of $g\cdot\conn$ is
$g\curvature{\conn} g^{-1}$. The \textit{gauge group} of $P$ is the group of bundle automorphisms of $P$ which cover the
identity on $M$.\index[terminology]{gauge group!of a principal bundle} In this case, the gauge group is
$\mathcal{G}\cong\Mapsm(M,\SU_2)$ and it acts on $P\cong M\times \SU_2$ by left multiplication, so the
$\mathcal{G}$-action preserves flat connections.

On the other hand, each flat connection $A$  gives rise to a \textit{holonomy representation} $\uppi_1(M)\rightarrow
G$: parallel transport along a loop $\gamma$ at $m_0$ gives an automorphism of the fiber $\SU_2$ at $m_0$, which
depends only on the homotopy class $[\gamma]\in\uppi_1(M,m_0)$.\index[terminology]{holonomy
representation} With a bit of work, one can recover the well-known fact that 
\begin{equation}
	\{\text{flat connections on }P\}/\mathcal{G}\hookrightarrow R(M) \colonequals\Hom(\uppi_1(M),\SU_2)/\text{conjugation} \period
\end{equation}
Since $P$ is trivial, this injection becomes a bijection. In fact, this can be upgraded to a homeomorphism, with
the right-hand side the character variety of $M$.\index[terminology]{character variety}

Now look at the $3$-form
\begin{equation*}
  \CS_{f }(\conn)=\CS_{f }(\conn, \conn_0)=\int_{[0,1]}\frac{1}{8\pi^2}\tr(\curvature{\conn}\wedge \curvature{\conn}) \comma
\end{equation*}
where as usual $\curvature{\conn}$ is the curvature of $\conn$.
Integrating over $M$ gives us the \textit{Chern--Simons functional} on $\cA_P$:\index[terminology]{Chern--Simons functional}
\index[notation]{cs@$\tilde{\cs}$}
\begin{equation}
\label{CS_fun}
  \tilde{\cs}(\conn)=\int_{M\times[0,1]}\frac{1}{8\pi^2}\tr(\curvature{\conn}\wedge \curvature{\conn}) = \frac{1}{8\pi^2}\int_M 
  \tr\paren{\conn\wedge \d\conn +\frac{2}{3}\conn\wedge [\conn\wedge \conn]} \period
\end{equation}
This map is smooth and functorial in $P\rightarrow M$, and up to $\ZZ$ factors, it is independent of the
trivialization of $P$. Therefore $\tilde{\cs}$ descends to a functional
\begin{equation}
	\cs \colon R(M)\cong\cA_P/\mathcal{G}\rightarrow\RR/\ZZ \period
\end{equation}
\index[notation]{cs@$\cs$}
The reason is that if $\sigma\in\mathcal{G}$, there is a straight-line path in $\cA_P$ from $\conn$ to
$\sigma\cdot\conn$, which we can interpret as a connection $\overline\conn$ on $[0,1]\times P\to[0,1]\times M$
with curvature $\curvature{\overline\conn}$. When we quotient by $\mathcal G$, we obtain a loop in $\cA_P/\mathcal{G}$, or
a connection on $P\times \Circ\to M\times \Circ$. Then
\begin{equation*}
  \cs(\sigma\cdot \conn)-\cs(\conn)=\int_{M\times \Circ}\frac{1}{8\pi^2}\tr(\curvature{\overline\conn}\wedge \curvature{\overline\conn}) =
  \int_{M\times\Circ} c_2(P\times\Circ) \comma
\end{equation*}
which is an integer because $c_2$ is an integer-valued characteristic class.

The function $\cs\colon R(M)\to\RR/\ZZ$ is a homotopy invariant of $M$. In practice, it is relatively computable,
as we will see for lens spaces.

%-------------------------------------------------------------------%
%  Chern–Simons invariants of Lens spaces                           %
%-------------------------------------------------------------------%

\subsubsection{Chern--Simons invariants of lens spaces}
Let $p$ and $q$ be coprime positive integers and $\zeta$ be a primitive $p$-th root of unity. Then
$\ZZ/p$ acts on $\CC^2$ by
\begin{equation}
	(z_1, z_2)\mapsto (\zeta z_1, \zeta^q z_2) \period
\end{equation}
Restricting to the unit $\Sph{3}\subset\CC^2$, this is a free action, and the quotient is called a \textit{lens space}
and denoted $L(p, q)$ \cite[\S 20]{Tie08}.
\index[terminology]{lens space}
\index[notation]{$L(p, q)$}

Lens spaces form a nice collection of examples of $3$-manifolds, and given an invariant of $3$-manifolds, one can
test how powerful it is by checking how well it distinguishes inequivalent lens spaces. For example, $L(5, 1)$ and
$L(5, 2)$ have the same homology and fundamental group, but are not homotopy equivalent \cite{Ale19}; and there are
homotopy-equivalent lens spaces which are not homeomorphic \cites{Rei35}[\S 5]{Whi41}{Bro60}. The full
classifications of lens spaces up to homotopy equivalence and homeomorphism are known, due to work of
Whitehead \cite[\S 5]{Whi41}, resp.\ Reidemeister \cite{Rei35} and Brody \cite{Bro60}.

Let's test the power of Chern--Simons invariants on lens spaces.

\begin{theorem}[{\cite[Theorem 5.1]{KK}}]\label{theorem:kk}
  The image of $\cs\colon R(L(p, q))\to\RR/\ZZ$ is the set\index[terminology]{lens space}
  \[
    \left\{-\frac{n^2r}{p} \,\bigg|\, n=0,1,\ldots,\left\lfloor \frac{p}{2}\right\rfloor \right\} \comma
  \]
  where $r$ is an integer satisfying $qr\equiv -1\bmod p$.
\end{theorem}

You can think of $\image(\cs)$ as the set of Chern--Simons invariants of a $3$-manifold.

\begin{remark}  
  Two lens spaces $L(p,q)$ and $L(p',q')$ have the same set of Chern--Simons invariants if and only if $p=p'$ and
  $q'q^{-1}\equiv a^2 \mod p$ for some $a\in\ZZ$, i.e., there is an orientation preserving homotopy equivalence
  between the two \cite[\S 5]{Whi41}. Hence Chern--Simons invariants detect the homotopy type of lens spaces.
\end{remark}

\begin{proof}[Proof sketch of \cref{theorem:kk}]
  The lens space $L(p,q)$ can be obtained by gluing the boundary of two solid tori $X$, $K$ together via an element
  \begin{equation*}
    \begin{pmatrix}
      p&q \\
      r&s 
  \end{pmatrix}\in \SL_2(\ZZ)
  \end{equation*}
  Let $x=\Circ\times\{1\}$ represent a generator of $\uppi_1(X)$ and $y$ a meridian of $\partial X$. Let $\mu, \lambda$ be the corresponding generators of $\partial K$, so $\mu=px+qy $ and $ \lambda=rx+sy$. 

  Now we utilize some general results about $3$-manifolds with a single torus boundary in \cite{KK}. 
  Suppose we have a path $f _t$ in $\Hom(\uppi_1(X),\SU_2)$ with 
  \begin{equation*}
    f _t(\mu) =
    \begin{pmatrix}
      e^{2\pi i\alpha(t)}& \\
      &e^{-2\pi i\alpha(t)} 
    \end{pmatrix} 
    \andeq
    f _t(\lambda) =
    \begin{pmatrix}
      e^{2\pi i\beta(t)}& \\
      &e^{-2\pi i\beta(t)} 
    \end{pmatrix},
  \end{equation*}
  where $\alpha, \beta \colon [0,1] \rightarrow\RR$. 
  The corresponding path of flat connections takes the form
  \begin{equation*}
  A_t=\begin{pmatrix}
     i\alpha(t)& \\
      &-i\alpha(t) 
  \end{pmatrix}\, \d x+ \begin{pmatrix}
      i\beta(t)& \\
      &-i\beta(t) 
  \end{pmatrix}\, \d y
  \end{equation*}
  near the torus boundary. 
  If $f _0 $ and $ f _1$ send $\mu$ to $1$,
  then \cite[Theorem 4.2]{KK}
  \begin{equation}
  	\cs(f _1)-\cs(f _0)=-2\int_0^1\beta\alpha'\, dt \bmod \ZZ.
  \end{equation}
  On the other hand, a holonomy representation on $X$ extends to one on the Dehn filling $M$ (in our case, the lens
  space itself) if and only if it sends $\mu$ to 1 (\textit{ibid.}, proof of Theorem 4.2).\index[terminology]{Dehn
  filling}

  Back to the sketch. We take $\gamma_t$ to be a path sending $x$ to $e^{2\pi i\theta}$ with $\theta\in[0, 1/2]$.
  (Every representation of $\uppi_1(X)$ is conjugate to a representation in the image of the path.)
  Then $\gamma_{t_1}$ extends to a representation $f _t$ of $\uppi_1(L(p,q))=\ZZ/p$ if and only if $pt_1\in\ZZ$, so
  we can obtain $\lfloor p/2 \rfloor+1$ conjugacy classes of representations of $\ZZ/p$, which correspond to $t_1 =
  n/p$ for $0\le n\le \lfloor p/2\rfloor$.

  On the other hand, $\alpha(t)=pt$ and $\beta(t)=rt$, so
  \begin{equation*}
    \cs(f _{t_1})=-2\int_0^{t_1}\beta\alpha'dt=-rpt_1^2\period
  \end{equation*}
  Plug in $t_1 = n/p$ and conclude.
\end{proof}

%-------------------------------------------------------------------%
%  Application: confiuguration spaces of Lens spaces                %
%-------------------------------------------------------------------%

\subsection{Application: configuration spaces of lens spaces}
\label{config_ssec}
To strengthen our Chern--Simons invariants, let's use them to study a related invariant of lens spaces:
the homotopy type of $F_2(L(p, q))$, the space of two-point subsets of
$L(p, q)$. Longoni--Salvatore \cite{LS05} showed that this distinguishes $L(7, 1)$ and $L(7, 2)$, which are
homotopy equivalent; the fact that the homotopy type of $F_2(X)$ knows more than the homotopy type of $X$ was a
surprising result. Differential cohomology enters the story with work of Evans-Lee--Saveliev \cite{deletedsquare}
using Chern--Simons invariants to provide a more comprehensive way to test whether the two-point configuration
spaces of two homotopy-equivalent lens spaces are homotopy equivalent.
%In this subsection, we follow Evans-Lee--Saveliev \cite{deletedsquare}, extending Chern--Simons invariants to
%two-point configuration spaces of lens spaces. Recall that $L(7, 1)$ and $L(7, 2)$ are homotopy equivalent but not
%homeomorphic, and Longoni--Salvatore \cite{LS05} showed that the homotopy type of their two-point configuration
%spaces can tell them apart. Now the question is: does the homotopy type of two-point configuration space
%distinguish lens spaces up to homeomorphism?

Choose a lens space $L=L(p,q)$ and a CW structure on it with a single $i$-cell $e_i$ for $0\le i\le 3$. Let $X =
L\times L$. The \textit{two-point configuration space}\index[terminology]{configuration space} of $L$ is
\begin{equation}
  X_0 \colonequals \Conf_2(L) \cong X \smallsetminus \Delta \comma
\end{equation}
where $\Delta\subset X$ is the diagonal, i.e.\ the subspace of elements $(x, x)$ with $x\in L$. Taking the
product CW structure on $X$, $X_0\subset L$ is a subcomplex, and the inclusion $X_0\hookrightarrow L\times
L$ induces an isomorphism of fundamental groups.

Using this CW structure, one can compute that
\begin{equation*}
  \H_3(X)\cong\ZZ\oplus\ZZ\oplus\ZZ/p \semicolon
\end{equation*}
the classes $[e_0\times e_3]=[e_0\times L]$ and $[e_3\times e_0]=[L\times e_0]$ generate the two $\ZZ$ summands and $[e_1\times
e_2 + e_2\times e_1]$ generates the $\ZZ/p$ summand.

\begin{lemma}
  There is a closed, oriented $3$-manifold $S$ with a map $f\colon M\to X$ such that $f_*[M] = [e_1\times e_2
  + e_2\times e_1]$.
\end{lemma}

\begin{proof}
  This is a special case of the \textit{Steenrod realization problem} asking when a given degree-$n$ homology class
  can be represented as a map from a closed, oriented $n$-manifold. This can be reformulated as a question about
  oriented bordism $\Omega_n^{\SO}(X)$, a generalized homology
  theory,\index[notation]{OmeganSO@$\Omega_n^{\SO}$}\index[terminology]{bordism} and the natural
  transformation
  \begin{align*}
    \Omega_n^{\SO} &\to \H_n \\ 
    \shortintertext{sending}
    (M, f\colon M\to X) &\mapsto f_*[M] \period
  \end{align*}
  In this
  form, the question was answered negatively in general by Thom \cite[Théorème III.9]{ThomThesis}, but when $X$ is a
  manifold, $\Omega_3^{\SO}(X)\to\H_3(X)$ is surjective (\textit{ibid.}, Théorème III.3).
\end{proof}

\noindent Evans-Lee--Saveliev \cite[\S 3]{deletedsquare} give an explicit example of such a representative manifold $S$.

With this choice of generators of $\H_3(X)$, the inclusion $\H_3(X_0)=\ZZ\oplus\ZZ/p \hookrightarrow \H_3(X)$
sends a generator of the free summand to $(1,1,0)$ and a generator of the torsion summand to $[S]=(0,0,1)$.

Given a representation
\begin{equation*}
  \alpha \colon \uppi_1(X)=\ZZ/p\times\ZZ/p\rightarrow \SU_2
\end{equation*}
and a closed, oriented
$3$-manifold $M$ with a map $f\colon M\to X$, we get a representation $f^*\alpha$ of $\uppi_1(M)$. Hence we can define an extension of the Chern--Simons invariants 
\begin{subequations}
\begin{equation}
  \cs_X \colon R(X_0)\rightarrow \Hom(\H_3(X_0), \RR/\ZZ)
\end{equation}
by
\begin{equation}
  \cs_X(\alpha) = \cs_M(f^*\alpha) = \frac{1}{8\pi^2}\int_M \tr\paren{\conn\wedge \d\conn +\frac{2}{3}\conn
  \wedge [\conn \wedge \conn]} \period
\end{equation}
\end{subequations}
A priori this depends on our choice of $(M, f)$, but it is actually independent of this choice, and is also
functorial in $X$. Thus we obtain a homotopy invariant for each pair of conjugacy class of representation and third
homology class.

Now we compute. Fix an $\SU_2$-representation $\alpha$, which is conjugate to one sending the generators of
$\uppi_1(X)$ to $e^{2\pi i k/p}$ and $e^{2\pi i \ell/p}$; we will call this representation $\alpha(k, \ell)$. Under
the two maps $L\rightrightarrows X$ realizing our two nontorsion generators of $\H_3(X)$, $\alpha(k, \ell)$ pulls
back to the representations sending a generator of $\uppi_1(L)$ to $e^{2\pi i k/p}$ and $e^{2\pi i\ell/p}$. By
\cref{theorem:kk}, the Chern--Simons invariants of these representations are $-k^2r/p$ and $-\ell^2r/p$, where $r$
can be any integer such that $qr\equiv -1\bmod p$.

Evaluating the Chern--Simons invariant for $S\to X$ is harder. Evans-Lee--Saveliev show that the choice of $S$
they constructed is Seifert fibered\index[terminology]{Seifert fiber space} over $S^2$ (\textit{ibid.}, Lemma 4.4),
allowing them to use a theorem of Auckly \cite[\S 2]{auckly} computing the Chern--Simons invariants of such
$3$-manifolds. The upshot is that the Chern--Simons invariant of $f^*\alpha$ on $S$ is $2k\ell/p$. Pulling back
along $X_0\hookrightarrow X$, our nontorsion generator of $\H_3(X_0)$ has Chern--Simons invariant $r(k^2 +
\ell^2)/p$, and our torsion generator has invariant $2k\ell/p$.

Let
\begin{equation*}
  X_0=\Conf_2(L(p,q)) \andeq X'_0=\Conf_2(L(p,q')) 
\end{equation*}
and suppose that $f \colon X_0\rightarrow X'_0$ is a homotopy equivalence.
Then the induced isomorphism $\ZZ/p\times \ZZ/p\to
\ZZ/p\times\ZZ/p$ on fundamental groups corresponds to a matrix
\begin{equation*}
	f_1=\begin{pmatrix}
   a& c\\
    b& d 
\end{pmatrix}\in\GL_2(\ZZ/p) \period
\end{equation*}
The induced isomorphism on $\H_3=\ZZ\oplus\ZZ/p$ has the form
\begin{equation*}
    h_3=\begin{pmatrix}
   \epsilon& 0\\
    a & b
\end{pmatrix} \comma
\end{equation*}
where $\epsilon=\pm 1$ and $b\in (\ZZ/p)^{\times}$. Using naturality of Chern--Simons
invariants, we can deduce the following numerical constraints: 

\begin{proposition}[{\cite[Proposition 5.2]{deletedsquare}}]
  \label{f_h_restr}
  If $f$ is a homotopy equivalence, then $\epsilon q'\equiv qa^2 \bmod p$ and
  \begin{equation*}
    f_1 = \begin{pmatrix}
       a &0 \\
        0& \pm a 
    \end{pmatrix}, 
    \begin{pmatrix}
      0 & a\\
         \pm a &0 
    \end{pmatrix};\ 
    h_3=\begin{pmatrix}
       \epsilon& 0\\
        0 & \pm a^2 
    \end{pmatrix} \period
  \end{equation*}
\end{proposition}
Composing with the swap map $(x, y)\mapsto (y, x)$ if necessary, we can and do make $f$ diagonal, rather than
antidiagonal.

To learn more information about lens spaces, we have to combine \cref{f_h_restr} with other invariants. These
invariants are further away from differential cohomology, so we will be terser and point the reader towards
references with more information. Specifically, we will combine the Chern--Simons invariants results from above
with information about Massey products\index[terminology]{Massey product} in the cohomology of the universal cover
$\tilde X_0$ of $X_0$.
% info abour H*(\tilde X_0)
\begin{proposition}[{\cite[Lemma 6.1]{deletedsquare}}]
\label{deleted_square_coh}
  There is an isomorphism
  \begin{equation*}
    \H^*(\tilde X_0)\cong \ZZ[a_1,\dotsc,a_{p-1}, b]/(a_i^2, b^2) \comma
  \end{equation*}
  where $\lvert a_i\rvert = 2$ and $\lvert b\rvert = 3$.
\end{proposition}

\begin{proof}[Proof sketch]
  The universal cover of $X$ is $\Sph{3}\times \Sph{3}$; therefore the universal cover of $X_0$ is a subspace of $\Sph{3}\times
  \Sph{3}$, specifically the complement of the orbit of the diagonal of $\Sph{3}\times \Sph{3}$ under the $\uppi_1(X)$-action.
  Therefore there is a map $\pi\colon \tilde X_0\hookrightarrow \Sph{3}\times \Sph{3}\to \Sph{3}$ given by inclusion followed by
  projection onto the first factor; it is a surjective submersion, and the fiber is a $(p-1)$-punctured $\Sph{3}$. 
  Set up the Serre spectral sequence;\index[terminology]{Serre spectral sequence} there are only a few differentials not
  zeroed out by degree considerations, and they vanish because $\pi$ has a section. 
  Thus the spectral sequence
  collapses. 
  There are no nontrivial extension questions, so the cohomology ring of $\tilde X_0$ is the tensor
  product of
  \begin{equation*}
    \H^*(\Sph{3};\ZZ/2)\cong\ZZ/2[b]/(b^2) \andeq \H^*(\Sph{3}\setminus \{x_1,\dotsc,x_{p-1})\}\cong\ZZ/2[a_1,\dotsc,a_{p-1}]/(a_i^2) \period \qedhere
  \end{equation*}
\end{proof}

Let $a_0 = -a_1 - \dotsb -a_{p-1}$. Miller \cite[\S 2.1]{Mil11} calculates the $\uppi_1(X_0)$-action on $\H^2(\tilde
X_0)$. Specifically, for $k,\ell\in\ZZ/p$, let $\tau_{k,\ell}$ denote the element corresponding to $(k,\ell)$ under
the identification
\begin{equation*}
    \uppi_1(X_0)\cong\ZZ/p\times\ZZ/p
\end{equation*}
above. Then,
\begin{equation*}
	\tau_{k,\ell}\cdot a_i = a_{i+k-\ell} \period
\end{equation*}
This puts an additional constraint on a homotopy equivalence $f\colon X_0\to X_0'$: $f$ must intertwine the action
map $\uppi_1(X_0)\to\Aut(\H^2(\tilde X_0))$. With $\alpha,\epsilon$ as above, this implies $f =
\alpha\cdot\mathrm{id}$ and that the following diagram commutes \cite[Proposition 6.3]{deletedsquare}:
\begin{equation*}
% https://q.uiver.app/?q=WzAsNCxbMCwwLCJcXEheMihcXHRpbGRlIFhfMCkiXSxbMSwwLCJcXEheMihcXHRpbGRlIFhfMCcpIl0sWzAsMSwiXFxIXjIoXFx0aWxkZSBYXzApIl0sWzEsMSwiXFxIXjIoXFx0aWxkZSBYXzAnKSJdLFswLDEsIlxcdGlsZGUgZl4qIl0sWzIsMywiXFx0aWxkZSBmXioiXSxbMSwzLCJcXHRhdV97ayxcXGVsbH0iXSxbMCwyLCJcXHRhdV97XFxhbHBoYSBrLCBcXGFscGhhXFxlbGx9IiwyXV0=
\begin{tikzcd}
	\H^2(\tilde X_0') & \H^2(\tilde X_0) \\
	\H^2(\tilde X_0') & \H^2(\tilde X_0) \period
	\arrow["{\tilde f^*}", from=1-1, to=1-2]
	\arrow["{\tilde f^*}"', from=2-1, to=2-2]
	\arrow["{\tau_{k,\ell}}", from=1-2, to=2-2]
	\arrow["{\tau_{\alpha k, \alpha\ell}}"', from=1-1, to=2-1]
\end{tikzcd}
\end{equation*}
This provides an additional constraint on $f$.

% what is a Massey product?
Next we need information about Massey products in $\H^*(\tilde X_0; \ZZ/2)$. The Massey product is a secondary
cohomology operation; the corresponding primary operation is the cup product.\index[terminology]{Massey
product}\index[terminology]{secondary cohomology operation} As a quick review, a Massey product \cites[\S
2]{UM57}{Mas58} is defined for $x,y,z\in \H^*(X; A)$ when $A$ is a ring, $x\cupprod y = 0$, and $y\cupprod z = 0$: one
chooses cocycles $\overline x$, $\overline y$, and $\overline z$ representing $x$, $y$, and $z$ respectively, and
chooses cochains $A$ and $B$ such that
\begin{equation*}
  \delta A = \overline x\cupprod\overline y \andeq \delta B = \overline y\cupprod\overline z \period
\end{equation*}
The Massey product $\ang{x, y, z}$ is defined to be the set of cohomology classes
$[A\cupprod\overline z - \overline x\cupprod B]$ for all possible choices of $A$ and $B$. Massey products are
functorial, which follows directly from their definition.

Assume $p$ is odd and $0 < q < p/2$. It follows from \cref{deleted_square_coh} that there are identifications of
abelian groups
\begin{equation}\label{zeta_p}
	\FF_2(\zeta_p) \colonequals \FF_2[t]/(1 + t + \cdots + t^{p-1}) \isomorphism \H^m(\tilde X_0;\ZZ/2),\ m = 2,5 \semicolon
\end{equation}
for $m = 2$, this map sends $t^k\mapsto a_k\bmod 2$, and for $m = 5$, $t^k\mapsto a_kb\bmod 2$.\footnote{We chose
the notation $\FF_2(\zeta_p)$ because this is the cyclotomic field associated to a primitive $p$-th root
of unity $\zeta_p$ over $\FF_2$.\index[terminology]{cyclotomic field}} If $x,y,z\in
\H^2(\tilde X_0;\ZZ/2)$ satisfy $xy = yz = 0$ (so that their Massey product is defined), then
\begin{equation*}
  \ang{x,y,z} \subset \H^5(\tilde X_0;\ZZ/2) \comma
\end{equation*}
so we may describe these Massey products as (possibly multivalued) maps
\begin{equation}\label{zeta_Massey}
	\ang{-, -, -}\colon \FF_2(\zeta_p)\times\FF_2(\zeta_p)\times\FF_2(\zeta_p)\to
	\FF_2(\zeta_p) \period
\end{equation}
Miller \cite[Theorem 3.33]{Mil11} calculates these Massey products. 
For example, 
\begin{equation*}
  t^n\cdot \ang{t^k,t^\ell, t^j} = \ang{t^{k+n}, t^{\ell+n}, t^{j+n}} \andeq \ang{t^k, t^\ell, t^j} = \ang{t^j, t^\ell, t^k} \period
\end{equation*}
These two relations allow us to inductively reduce to the case when at least one of $j$, $k$, or $\ell$ is $0$; the description of the
Massey products in that case is a little more complicated, and can be found in \cite[Theorem 7.1]{deletedsquare}.

This leads us to our last obstruction. The two different maps 
\begin{equation*}
  \tilde f^*\colon \H^m(\tilde X_0';\ZZ_2)\to\H^m(\tilde X_0;\ZZ/2)
\end{equation*}
for $m = 2,5$, become the same map $\tilde f^*\colon
\FF_2(\zeta_p)\to\FF_2(\zeta_p)$ under the identification~\eqref{zeta_p}. Therefore we obtain the constraint that
this $\tilde f^*$ must intertwine the Massey product map~\eqref{zeta_Massey}.

Our three constraints (coming from Chern--Simons invariants, cohomology of $\tilde X_0$, and Mass\-ey products) each
boil down to numerical constraints on $p$ and $q$, and these are amenable to computer calculation. This is how
Evans-Lee--Saveliev showed that these constraints can detect some homotopy-equivalent but not homeomorphic lens
spaces that Longoni--Salvatore's techniques miss. These pairs include $L(11, 2)$ and $L(11, 3)$; $L(13, 2)$ and
$L(13, 5)$; and $L(17, 3)$ and $L(17, 5)$.

\newpage
%!TEX root = ../diffcoh.tex

%-------------------------------------------------------------------%
%-------------------------------------------------------------------%
%  Conformal immersions                                             %
%-------------------------------------------------------------------%
%-------------------------------------------------------------------%

\section{Conformal immersions}
\textit{by Charlie Reid}
\label{conformal_immersions}

Let $M$ be a smooth, $m$-dimensional manifold and suppose $M$ immerses in $\RR^n$ with normal bundle $\NormM$. Then
there is a short exact sequence
% https://q.uiver.app/?q=WzAsNSxbMCwwLCIwIl0sWzEsMCwiVE0iXSxbMiwwLCJUXFxSUl5uIl0sWzMsMCwiTk0iXSxbNCwwLCIwIl0sWzAsMV0sWzEsMl0sWzIsM10sWzMsNF1d
\begin{equation}
\begin{tikzcd}
	0 & \TanM & {\Tan\RR^n|_M} & \NormM & 0,
	\arrow[from=1-1, to=1-2]
	\arrow[from=1-2, to=1-3]
	\arrow[from=1-3, to=1-4]
	\arrow[from=1-4, to=1-5]
\end{tikzcd}
\end{equation}
so the Pontryagin classes of $\TanM$ and $\NormM$ satisfy\footnote{Because the Whitney sum formula for Pontryagin classes
only holds up to $2$-torsion, this formula should be thought of as taking place in cohomology with $\ZZ[1/2]$ or
$\RR$ coefficients.}
\begin{equation}
	p(\TanM)p(\NormM) = p(\Tan\RR^n|_M) = 1 \period
\end{equation}
The total Pontryagin class is the sum of $1$ and a nilpotent element ($p_1(M) + p_2(M) + \cdots$), hence is
invertible. This means $p(\NormM)$ is uniquely determined if it exists: there is a formula for $p_k(\NormM)$ in terms of
$p(\TanM)$. If $M$ immerses in $\RR^n$, then $\NormM$ is rank $n-m$, so $p_k(\NormM) = 0$ for $k > n-m$, and because of the
formula, this is actually a constraint on the Pontryagin classes of $\TanM$. Thus Pontryagin classes can be used to
prove nonimmersion results for smooth manifolds by showing this constraint is not met.

In \cref{DifferentialCharacteristicClasses}, we saw that given a connection on the tangent bundle, Pontryagin
classes lift to differential cohomology. It therefore seems worthwhile to imitate the above argument and use
on-diagonal differential Pontryagin classes given by the Levi-Civita connection to obstruct isometric immersions of
Riemannian manifolds. Chern and Simons \cite{cs} did this, though with a few key differences.
\begin{enumerate}
	\item Chern and Simons were able to show (\textit{ibid.}, Theorem 4.5) that if $g$ and $g'$ are two conformally
	equivalent metrics on a manifold $M$, with Levi-Civita connections $\conn$, resp.\ $\conn'$, then $\phat(M,
	\conn) = \phat(M, \conn')$. Therefore the differential Pontryagin classes of $M$ are conformal invariants,
	and can be used to study conformal immersions.
	\item There is an additional integrality result which has no analogue in the purely topological case
	(\textit{ibid.}, Theorem 5.14): when a Pontryagin class' Chern--Weil form vanishes, the corresponding
	Chern--Simons form is closed, and one-half of its de Rham class is contained within the lattice
	$\image(\H^*(-;\ZZ)\to\H^*(-;\RR))$. After some more work, this leads to another necessary
	condition for the existence of a conformal immersion.
\end{enumerate}
As an example, $\RRP^3$ smoothly immerses in $\RR^4$ \cite{Boy03}, and given the round metric,
$\RRP^3$ locally conformally immerses in $\RR^4$. But Chern--Simons show (\textit{ibid.}, \S 6) that there is
no conformal immersion $\RRP^3\hookrightarrow\RR^4$.

In \cref{ssec:conformal_invariance}, we prove that the on-diagonal differential Pontryagin classes of the
Levi-Civita connection are conformal invariants of the Riemannian metric. Then, in \cref{ssec:obstruct}, we use
on-diagonal differential Pontryagin classes to obstruct conformal immersions. Finally, in \cref{ssec:div2}, we
produce the integrality obstruction using the Chern--Simons form and use it to show $\RRP^3$ with the round metric
cannot conformally immerse in $\RR^4$.

\begin{remark}
	The story we just told is a little anachronistic: Chern--Simons' work came before Cheeger--Simons' paper on
	differential characters, and was not stated in this language. But Chern and Simons were aware that their ideas
	could be rephrased as calculations in the ring of differential characters, as they write in the introduction to
	their paper. In any case, the paper \cite{cs} is best known for an entirely different reason: for introducing the
	Chern--Simons form of a connection!
\end{remark}

%-------------------------------------------------------------------%
%-------------------------------------------------------------------%
%  Conformal invariance of differential Pontryagin classes          %
%-------------------------------------------------------------------%
%-------------------------------------------------------------------%

\subsection{Conformal invariance of differential Pontryagin classes}
\label{ssec:conformal_invariance}
Let $G$ be a compact Lie group. Recall that given a degree-$k$ invariant polynomial $f$ on the Lie algebra $\g$ and a
characteristic class $c^\ZZ\in\H^{2k}(\BG;\ZZ)$, we obtain a differential characteristic class $\chat\in
\Hhat^{2k}(\BunGnabla ;\ZZ)$ (as proven in \cref{differential_CW_lift}) and a Chern--Simons form
$\CS_f(\conn)\in\Omega^{2k-1}(P)$ given a principal $G$-bundle $\pi\colon P\to M$ and a connection $\conn$ on $P$
(as defined in~\eqref{total_space_CS}). We are specifically interested in the Pontryagin polynomials $P_k$ from
\cref{pontrjagin}, which we lifted to on-diagonal differential Pontryagin classes $\phat_k$ in
\cref{differential_Pontryagin}.

Our aim in this section is to prove:
\begin{theorem}[{\cite[Theorem 4.5]{cs}}]\label{diff_p_conformally_invariant}
	Let $M$ be a manifold and $g_0,g_1$ be conformally equivalent Riemannian metrics on $M$. If $\conn_0$ and
	$\conn_1$ denote the Levi-Civita connections for $g_0$ and $g_1$, then for all $k$, we have
	\begin{equation*}
		\phat_k(M, \conn_0) = \phat_k(M, \conn_1)
	\end{equation*}
	and $\CS_{P_k}(\conn_0) - \CS_{P_k}(\conn_1)$ is exact.
\end{theorem}

The first ingredient in the proof is a variation formula.
\index[terminology]{variation formula!for Chern--Simons forms}
\index[terminology]{Chern--Simons form}

\begin{lemma}[{(variation formula \cite[Proposition 3.8]{cs})}]
	\label{variation_Chern_Simons}
	Suppose $\conn_t$ is a smooth path of connections on a principal $G$-bundle $P\to M$ and $\curvature{\conn_t}$ is the
	curvature of $\conn_t$. Then
	\begin{equation}
		\label{var_form_CS}
		\left.\frac{\mathrm d}{\mathrm dt} \CS_f(\conn_t)\right|_{t=0} = k\cdot  f(\conn'\wedge
		\curvature{\conn_0}^{k-1}) + \omega,
	\end{equation}
	where $\omega$ is exact and $\conn' = \left.\frac{\mathrm d}{\mathrm dt}(\conn_t)\right|_{t=0}$.
\end{lemma}

\begin{proof}
	It suffices to work universally in $\BunGnablatriv$, the stack of trivial principal $G$-bundles with connection: the
	forgetful map
	\begin{equation*}
		\BunGnablatriv\to\BunGnabla
	\end{equation*}
	is the universal principal $G$-bundle with connection in the world of stacks
	on $\Mfld$ (see \cref{subsec:presentation_of_BunGnabla}). The de Rham complex of $\BunGnablatriv$ is acyclic \cite[Theorem 7.19]{FreedHopkins}, so it suffices to
	apply the de Rham differential to~\eqref{var_form_CS} and then show both sides are equal.\footnote{One can avoid
	the use of the abstract object $\BunGnablatriv$ by using Narasimhan--Ramanan's $n$-classifying spaces \cite{NR61,
	NR63}.}

	For the left-hand side, we know
	\begin{align*}
		\d\left(\left.\frac{\mathrm d}{\mathrm dt} \CS_f(\conn_t)\right|_{t=0}\right) &= \left.\frac{\d}{\d t}\left(
		\d(\CS_f(\conn_t))\right)\right|_{t=0}\\
		&= \left.\frac{\d}{\d t}\left(f((\curvature{\conn_t})^k)\right)\right|_{t=0}\\
		&= k\cdot f(\curvature{\conn_0}'\wedge \curvature{\conn_0}^{k-1}),
	\end{align*}
	where $\curvature{\conn_0}' = \left.\frac{\mathrm d}{\mathrm dt}(\curvature{\conn_t})\right|_{t=0}$.

	For the right-hand side,
	\begin{align*}
		\d(k\cdot f(\conn' \wedge \curvature{\conn_0}^{k-1})) &= k\cdot f(\d\conn' \wedge \curvature{\conn_0}^{k-1}) -
		k(k-1) f(\conn' \wedge \d\curvature{\conn}\wedge\curvature{\conn}^{k-2})\\
		&= k f(\d\conn'\wedge\curvature{\conn_0}^{k-1}) - k(k-1) f(\conn'\wedge [\curvature{\conn_0}, \conn_0]\wedge
		\curvature{\conn}^{k-2})\\
		&= k f(\d\conn'\wedge \curvature{\conn_0}^{k-1}) + k f([\conn', \conn_0]\wedge \curvature{\conn_0}^{k-1}).
	\end{align*}
	This uses two important facts from Chern--Weil theory: that $\d\curvature{\conn_0} = [\curvature{\conn_0},
	\conn_0]$ together with the value of the invariant polynomial for a commutator \cite[(2.9)]{cs}. Now
	\begin{align*}
		\d\conn' &= \left.\frac{\d}{\d t}\left(\d\conn_t\right)\right|_{t=0}\\
		&= \left.\frac{\d }{\d t}\left(\curvature{\conn_t} - \frac 12[\conn_t, \conn_t]\right)\right|_{t=0}\\
		&= \curvature{\conn_0}' - [\conn', \conn_0],
	\end{align*}
	so $\d(k\cdot f(\conn'\wedge\curvature{\conn_t}^{k-1}))= k\cdot (\curvature{\conn_0}'\wedge
	\curvature{\conn_t}^{k-1})$ and we are done.
	\end{proof}
	\begin{proof}[Proof of \cref{diff_p_conformally_invariant}]
	Now for $f$ we take $P_k$, the invariant polynomial that we used in \cref{pontrjagin} to define the $k$-th
	Pontryagin class.  This is the pullback of the $2k$-th Chern polynomial under the complexification map
	$\mathfrak o(n)\to\mathfrak u(n)$; we tend not to use the pullback of the $(2k+1)^{\mathrm{st}}$ Chern polynomial
	as much because it is $2$-torsion and its Chern--Simons form is exact \cite[Proposition 4.3]{cs}.

	It suffices to show that $\delta\colonequals \CS_{P_k}(\conn_0) - \CS_{P_k}(\conn_1)$ is exact; this implies
	it is a closed form with integral periods, so the image $\overline\delta$ of $\delta$ in
	$\Omega^{4k-1}(M)/\Omegacl^{4k-1}(M)_\ZZ$ vanishes. This is the lower-left corner of the differential cohomology
	hexagon, and as we saw in \cref{iota_chern_simons}, applying
	\begin{align*}
		\iota\colon\Omega^{4k-1}(M)/\Omegacl^{4k-1}(M)_\ZZ &\to\Hhat^{4k}(M; \ZZ) \\
	\shortintertext{sends}
		\overline\delta &\mapsto\phat_k(P,\conn_0) - \phat_k(P, \conn_1) \comma
	\end{align*}
	so showing $\overline\delta = 0$ is good enough.

	Now to show $\delta$ is exact. It is always possible to connect $g_0$ and $g_1$ by a path $g_t$,
	$t\in(-\varepsilon, 1+\varepsilon)$ of conformally equivalent metrics. Moreover, this path may be chosen to satisfy
	\[g_t = e^{2th}g_0\]
	for some real-valued smooth function $h$. Choose such a path and let $\conn_t$ be the Levi-Civita connection of
	$g_t$.  Differentiating in $t$ commutes with the de Rham differential, so is suffices to show that $\frac{\d}{\d
	t}\CS_{P_k}(\conn_t)$ is exact; without loss of generality, we prove this for $t = 0$.
	\Cref{variation_Chern_Simons} means we only have to show
	\begin{equation*}
		P_k(\conn_0'\wedge \curvature{\conn_0}^{2k-1}) = 0.
	\end{equation*}
	For a little while we work locally on the bundle $\pi\colon B(M)\to M$ of frames: the fiber at $x\in M$ is the
	$\GL_n(\RR)$-torsor of orthonormal bases $(e_1, \dotsc, e_n)$ of $T_xM$. There are canonical one-forms
	$\omega_i\in\Omega^1(B(M))$ defined at a point $(x, (e_1, \dotsc, e_n))$ so that
	\begin{equation*}
		\d\pi = \sum_{i=1}^n \omega_i \cdot e_i.
	\end{equation*}
	Let $E_i$ be the horizontal vector field dual to $\omega_i$; here ``horizontal'' is with respect to the connection
	$\conn_0$. Then on frames orthogonal to $g_0$, there is a decomposition \cite[Lemma 4.4]{cs}
	\begin{equation*}
		\conn_{ij}' = \underbracket{\delta_{ij}\d(h\circ\pi)}_{\alpha} + \underbracket{E_i(h\circ\pi)\omega_j -
		E_j(h\circ\pi)\omega_i}_{\beta}.
	\end{equation*}
	We will address each piece separately. First, one directly checks that for $\varphi =
	(\varphi_{ij})\in\Omega^k(F(M))$,
	\begin{equation}
	\label{cycle_indices}
		P_k(\varphi\wedge \curvature{\conn}^{k-1}) = \sum_{i_1,\dotsc,i_k = 1}^n
		\varphi_{i_1i_2}\wedge (\curvature{\conn})_{i_2i_3}\wedge\dotsb\wedge (\curvature{\conn})_{i_ni_1}.
	\end{equation}
	Plugging in $\varphi = \alpha$, we obtain
	\begin{equation*}
		P_k(\alpha\wedge \curvature{\conn}^{2k-1}) = \d(f\circ\pi)\wedge P_{2k-1}(\curvature{\conn}^{2k-1}) = 0,
	\end{equation*}
	because $A$ is compatible with the metric. Now plugging $\beta$ into~\eqref{cycle_indices},
	\begin{equation}
	\label{beta_qty}
		P_k(\beta\wedge\curvature{\conn}^{2k-1}) = \sum_{i_1,\dotsc,i_{2k} = 1}^n
		(E_{i_1}(f\circ\pi)\omega_{i_2} - E_{i_2}(f\circ\pi)\omega_{i_1})\wedge(\curvature{\conn})_{i_2i_3}\wedge\dotsb\wedge
		(\curvature\conn)_{i_{2k}i_1}.
	\end{equation}
	The Jacobi identity implies $\sum \omega_i\wedge(\curvature{\conn})_{ij} = 0$, so~\eqref{beta_qty} vanishes as
	well. Lastly, we need to descend from $B(M)$ to $M$, and $0$ descends to $0$.
\end{proof}

\subsection{Obstructing conformal immersions with differential Pontryagin classes}
\label{ssec:obstruct}

%-------------------------------------------------------------------%
%  Differential Chern Classes                                       %
%-------------------------------------------------------------------%

%\subsubsection{Differential Chern Classes}
\index[terminology]{Chern class!differential refinement}
\index[terminology]{differential Chern class}
Recall the on-diagonal differential lifts of Chern classes we constructed in
\cref{DifferentialCharacteristicClasses}, specifically \cref{differential_Chern}, defined as follows: define the
\emph{Chern polynomials} $C_k\in I^k(\Uup_n)$ by 
\index[terminology]{Chern polynomial}
\begin{equation}
	\det(\lambda I - \frac{1}{2\pi i}A) = \sum_{k=0} ^n C_k(A)\lambda^{n-k};
\end{equation}
apply Chern--Weil theory to $C_k$, producing a characteristic class $c_k$. The integer
cohomology of $\Uup_n$ is torsion-free \cite[\S 29]{Bor53} and its image in de Rham cohomology contains $c_k$,
so there is a unique lift to $\chat_k$ to degree-$2k$ differential cohomology.

We will also need the \emph{inverse Chern polynomials} $C_k^\perp$, which are defined to satisfy
\index[terminology]{inverse Chern polynomials}
\begin{equation}
	(1+C_1+\cdots+C_n)(1+C_1^\perp + C_2^\perp + \cdots) = 1 \period
\end{equation}
For example, $C_1^\perp = -C_1$, $C_2^\perp = -C_2 - C_1C_1^\perp$, $C_3^\perp = -C_3 - C_2C_1^\perp -
C_1C_2^\perp$, and so on. Chern--Weil theory associates de Rham characteristic classes $c_k^\perp\in\HdR^{2k}$ to
these, and like ordinary Chern classes, these classes lift uniquely to differential cohomology classes
$\chat_k^\perp\in \Hhat^{2k}(\BunUnabla{n};\ZZ)$. They satisfy analogous formulas to the inverse Chern
polynomials: for example
\index[notation]{chatperp@$\chat_k^\perp$}
%
%We don't care about the $D$'s. Since $\BGL_n(\CC)$ has no torsion in it's cohomology, these polynomials are enough
%to define differential characteristic classes $\chat_k(V)$ of a complex vector bundle with connection of rank $n$,
%which refine the ordinary Chern classes $c_k(V)$. We also define the total Chern class
%$$\chat(V) = 1 + \chat_1(V) + \cdots + \chat_n(V)\in \Hhat^{\even}(M)$$
%We also define inverse Chern polynomials by . We can solve for them and write explicit inductive formulas:
%\begin{align*}
%	C_1^\perp &= -C_1 \comma \\ 
%	C_2^\perp &= -C_2 -C_1C_1^\perp \comma \\ 
%	C_3^\perp &= -C_3 -C_2C_1^\perp - C_1C_2^\perp \comma \\ 
%	&\vdots
%\end{align*}
%We then define differential characteristic classes $\chat_k^\perp$ using these polynomials. 
%I haven't mentioned it yet, but the Weil map and it's differential refinement are actually ring homomorphisms, so for example
\begin{equation}
\label{inverse_Chern_relation}
	\chat_2^\perp = -\chat_2 - \chat_1 \chat_1^\perp \period
\end{equation} 
%\todo[inline]{This uses that the Chern--Weil map is a ring homomorphism. Did we say this anywhere?}

%-------------------------------------------------------------------%
%  Differential Pontryagin Classes                                  %
%-------------------------------------------------------------------%

%\subsubsection{Differential Pontryagin Classes}
\index[terminology]{Pontryagin class!differential refinement}
\index[terminology]{differential Pontryagin class}

In \cref{differential_Pontryagin}, we defined on-diagonal differential Pontryagin classes $\phat_k$ in much the
same way as we defined differential Chern classes. Using the \emph{inverse Pontryagin polynomials} $P_k^\perp$,
defined to satisfy \index[terminology]{inverse Pontryagin polynomials}
\begin{equation}
	(1+P_1+\cdots+P_n)(1+P_1^\perp + P_2^\perp + \cdots) = 1 \comma
\end{equation}
we define on-diagonal inverse differential Pontryagin classes $\phat_k^\perp\in \Hhat^{4k}(\BunOnabla{n};\ZZ)$.  Because there is torsion in $\H^*(\BO_n;\ZZ)$, a priori the lift to differential cohomology
requires a choice, but there is a canonical way to do this: complexify to pass to on-diagonal inverse differential
Chern classes. This means that analogues of~\eqref{inverse_Chern_relation} and its higher-rank generalizations hold
for on-diagonal inverse Pontryagin classes. For example, $\phat_2^\perp = -\phat_2 - \phat_1\phat_1^\perp$.
\index[notation]{phatperp@$\phat_k^\perp$}

\begin{theorem}
\label{differential_conformal_immersion}
Let $M$ be a Riemannian manifold and $\phi \colon M^n\to \RR^{n+k}$ be a conformal immersion of $M$ into Euclidean
space. Then the image of $\phat^\perp_i(M, \conn^{\mathrm{LC}})$ in $\Hhat^{4i}(M;\ZZ[1/2])$ vanishes for all
$i>k/2$.
\end{theorem}
\index[terminology]{conformal immersion}
\begin{proof}
Since the classes $\phat_k$ are conformally invariant (\cref{diff_p_conformally_invariant}), so too are the classes
$\phat_k^\perp$. Therefore, without loss of generality, we can assume $\phi$ is isometric. Let $\NormM$ denote the
orthogonal normal bundle: there is an orthogonal direct sum
	$$\TanM \oplus \NormM = \Tan\RR^{n+k} = \underline{\RR}^{n+k}.$$
	The Levi-Civita connection $\conn_{\Tan\RR^{n+k}}^{\mathrm{LC}}$ on $\RR^{n+k}$ compresses to the Levi-Civita
	connection $\conn_{\TanM}^{\mathrm{LC}}$ on $M$, and to a connection $\conn_{\NormM}$ on $\NormM$. Since
	$\conn_{\Tan\RR^{n+k}}^{\mathrm{LC}}$ is flat, it is compatible with
	$\conn_{\TanM}^{\mathrm{LC}}\oplus\conn_{\NormM}$ (\cref{compatible_connection}). Hence
	\begin{equation*}
		\phat(\TanM, \conn_{\TanM}^{\mathrm{LC}})*\phat(\NormM, \conn_{\NormM})= \phat(\Tan\RR^{n+k},
		\conn_{\Tan\RR^{n+k}}^{\mathrm{LC}}) = 1\comma
	\end{equation*}
	implying
	\begin{equation*}
		\phat^\perp(\TanM, \conn_{\TanM}^{\mathrm{LC}}) = \phat(\NormM, \conn_{\NormM})\period
	\end{equation*}
	Since $\NormM$ has rank $k$, $\phat_i(\NormM, \conn_{\TanM}^{\mathrm{LC}})$ vanishes for $i > k/2$.
\end{proof}
\begin{remark}
As always, we use $\ZZ[1/2]$ coefficients because the Whitney sum for Pontryagin classes is more complicated over
the integers. See Thomas \cite{Tho62} and Brown \cite[Theorem 1.6]{Bro82}. The extra factors ultimately come from
Chern classes, so they too admit differential refinements, and a $\ZZ$-valued differential Whitney sum formula
exists. Using this, it is possible to upgrade \cref{differential_conformal_immersion} to take place in
$\Hhat^*(M;\ZZ)$.
\index[terminology]{Whitney sum formula}
\end{remark}

\subsection{Dividing by \texorpdfstring{2}{$2$}}
\label{ssec:div2}
We foreshadowed that Chern--Simons theory will allow us to prove that $\RRP^3$ with the round metric does not
conformally immerse in $\RR^4$, but to actually prove this we need another obstruction. This one is an evenness
result: we will use the Chern--Simons form to define a de Rham cohomology class of on the frame bundle of $\RRP^3$,
and prove that a conformal immersion would imply this class is in the image of the map induced by the inclusion
$2\ZZ\to\RR$. A direct calculation shows this is not the case, and we conclude.
\begin{lemma}[{\cite[Proposition 3.15]{cs}}]
\label{universal_pullback_coboundary}
If $\pi\colon P\to M$ is a principal $G$-bundle with connection $\conn$,
there is a cochain $u\in C^{2k-1}(M; \RR/\ZZ)$ such that $\delta(u) = f(\curvature\conn)\bmod\ZZ$ and in
$C^*(P; \RR/\ZZ)$, $\CS_f(\conn)\bmod\ZZ \pi^*(u)$ is a coboundary.
\end{lemma}
\begin{proof}
Since $[f(\curvature\conn)]$ is in the image of the map from integer cohomology to de Rham cohomology,
$f(\curvature\conn)\bmod\ZZ$ is a coboundary, so choose $u\in C^{2k-1}(M; \RR/\ZZ)$ with $\delta u =
f(\curvature\conn)\bmod\ZZ$. Then
\begin{align*}
	\delta(\pi^*(u)) &= \pi^*(\delta u) = \pi^*(f(\curvature\conn))\bmod\ZZ\\
		&= \delta(\CS_f(\conn))\bmod\ZZ = \delta(\CS_f(\conn)\bmod\ZZ) \period
\end{align*}
That is, $\delta(\pi^*(u) - \CS_f(\conn)\bmod\ZZ)$ vanishes.
\end{proof}
Let $\pi\colon P\to M$ be a principal $G$-bundle with connection $\conn$. In the previous chapter,
specifically~\eqref{chern_simons_differential}, we showed that
$\d\CS_f(\conn) = \pi^*f(\curvature\conn)$. Therefore if $f(\curvature\conn) = 0$, $\CS_f(\conn)$ is
closed and defines a class $[\CS_f(\conn)]\in \H^{2k-1}(P;\RR)$.

\begin{corollary}[{\cite[Theorem 3.16]{cs}}]
	\label{when_CS_is_pulledback}
	Assume that $f(\curvature\conn) = 0$. 
	Then there is a class
	\begin{equation*}
		\overline u\in \H^{2k-1}(M;\RR/\ZZ)
	\end{equation*}
	such that in $\H^{2k-1}(P;\RR/\ZZ)$, we have
	\begin{equation*}
		[\CS_f(\conn)] \bmod \ZZ = \pi^*(\overline u) \period
	\end{equation*}
\end{corollary}

\begin{proof}
By hypothesis of \cref{universal_pullback_coboundary}, $\delta(u) = f(\curvature\conn) = 0$, so we can choose
$\overline u$ to be the class of $u$ in cohomology.
\end{proof}
\begin{example}
\label{cpx_stiefel}
Let $\St_n(\CC^{n+k})$ denote the \textit{Stiefel manifold} of isometric immersions
$\CC^n\hookrightarrow\CC^{n+k}$. Sending an immersion to its image defines a map $\pi$ to the \textit{Grassmannian
manifold} $\Gr_n(\CC^{n+k})$ parametrizing codimension-$k$ subspaces of $\CC^{n+k}$, and this map is a principal
$\Uup_n$-bundle.
\index[terminology]{Stiefel manifold}
\index[terminology]{Grassmannian}
This bundle has a natural connection. It is equivalent to describe the connection on the associated rank-$n$
complex vector bundle $\pi'\colon S\to\Gr_n(\CC^{n+k})$, which is the tautological bundle. If $\rho\colon
(-\varepsilon, \varepsilon)\to S$ is a smooth curve, $\rho(t)$ is an element of the vector space
$\pi(\rho(t))\in\Gr_n(\CC^{n+k})$; we specify the connection by declaring the covariant derivative of $\rho(t)$
along $\pi\circ\rho$ to be the orthogonal projection of $\rho'(t)$ into the subspace $\pi(\rho(t))$. Call this
connection $\conn^{\mathrm{can}}$.

There is a canonically defined rank-$k$ complex vector bundle $Q\to\Gr_n(\CC^{n+k})$, whose fiber at an
$n$-dimensional subspace $V\subset\CC^{n+k}$ is $V^\perp\subset\CC^{n+k}$. Thus $S\oplus Q = \underline\CC^{n+k}$,
so in a similar manner as in the proof of \cref{differential_conformal_immersion}, $[C_i^\perp(\conn^{\mathrm{can}})] = 0$, i.e.\ $C_i^\perp(\conn^{\mathrm{can}})$ is exact. The Grassmannian is a compact,
irreducible Riemannian symmetric space, so since $C_i^\perp(\conn^{\mathrm{can}})$ is an invariant, exact
differential form, it must vanish. Therefore \Cref{when_CS_is_pulledback} tells us
$[\CS_{C_i^\perp}(\conn^{\mathrm{can}})]\bmod\ZZ$ pulls back from $\overline u\in
\H^{2k-1}(\Gr_n(\CC^{n+k});\RR/\ZZ)$. Because the cohomology of complex Grassmannians is concentrated in even
degrees, $\overline u = 0$, meaning $[\CS_{C_i^\perp}(\conn^{\mathrm{can}})]$ is in the image of the map
\begin{equation*}
	\H^*(\St_n(\CC^{n+k});\ZZ)\to\H^*(\St_n(\CC^{n+k});\RR) \period
\end{equation*}

\end{example}
By passing to real vector bundles, we will gain an additional factor of $2$.
We will say a real-valued cohomology
class is \textit{contained in the even integer lattice} if it is in the image of the composite
\begin{equation}
	\H^*(-;\ZZ)\overset{\cdot 2}{\longrightarrow} \H^*(-; \ZZ)\longrightarrow\H^*(-;\RR) \period
\end{equation}

\begin{lemma}[{\cite[Lemma 5.12]{cs}}]
	\label{double_cpxif}
	Let $c\colon\St_n(\RR^{n+k})\to\St_n(\CC^{n+k})$ be the complexification map. 
	For $\ell > 0 $, the image of
	\begin{equation*}
		c^*\colon \H^\ell(\St_n(\CC^{n+k});\ZZ)\to\H^\ell(\St_n(\RR^{n+k});\ZZ)
	\end{equation*}
	is contained in the even integer lattice.
\end{lemma}

\begin{proof}
First suppose $k = 0$, for which $\St_n(\CC^n)\cong\Uup_n$ and $\St_n(\RR^n)\cong\Or_n$; $c$ is the usual
complexification map. It suffices to show that the mod $2$ reductions of all positive-degree classes in the image
of $c^*$ vanish.

At this point we need a tool called the \textit{inverse transgression map}.\index[terminology]{inverse
transgression map} We will say more about this map in \cref{transgression_detail} at the end of this chapter; for
this proof, we need only that inverse transgression is a map $\tau\colon\H^\ell(\BG;\ZZ)\to\H^{\ell-1}(G;\ZZ)$
satisfying two key properties:
\begin{enumerate}
	\item $\tau$ is natural in $G$, and
	\item for $A = \ZZ$ or $\ZZ/2$ and $x\in \H^*(\BG; A)$, $\tau(x^2) = 0$.
\end{enumerate}
Let $\Bup c\colon\BO_n\to\BU(n)$ be the map induced from complexification on classifying spaces. We
know $(\Bup c)^*(c_i) \bmod 2 = w_i^2$ \cite[Theorem 1.5]{Bro82}, so
\begin{equation}
	c^*(\tau(c_i))\bmod 2 = \tau((\Bup c)^*(c_i)\bmod 2) = 0 \period
\end{equation}
This suffices because $\{\tau(c_i)\}$ generates $\H^*(\Uup_n;\ZZ)$ \cite[Théorèmes 8.2 et 8.3]{Bor54}.

For more general $k$, recall that
\begin{equation*}
	\St_n(\RR^{n+k}) \cong \Or_{n+k}/\Or_k \andeq \St_n(\CC^{n+k}) \cong \Uup_{n+k}/\Uup_k \period
\end{equation*}
Let $\pi$ denote the quotient $\Or_{n+k}\to\St_n(\RR^{n+k})$ as well as its complex analogue.
Then $\pi$ commutes with complexification, so it suffices to show that
\begin{equation*}
	\pi^*\colon\H^*(\St_n(\RR^{n+k});\ZZ/2)\to\H^*(\Or_{n+k};\ZZ/2)
\end{equation*}
is injective, and this is due to Borel \cite[\S10]{Bor53}.
\end{proof}
	
This extra factor of two provides an additional obstruction to the existence of a conformal immersion, and this is
what we will use to show $\RRP^3$ cannot conformally immerse in $\RR^4$. 

\begin{theorem}[{\cite[Theorem 5.14]{cs}}]
	\label{extra_factor_of_two}
	Let $M$ be an $n$-dimensional Riemannian manifold,
	\begin{equation*}
		B(M)\to M
	\end{equation*}
	be the principal $\Or_n$-bundle of 
	frames, and $\conn$ be the Levi-Civita connection on $B(M)$. Suppose $M$ conformally immerses in $\RR^{n+k}$;
	then, for $i\ge \lfloor k/2\rfloor$, $\CS_{P_i^\perp}(\conn)$ is contained in the even integer lattice.
\end{theorem}

\begin{proof}
	Let $\varphi\colon M\to\RR^{n+k}$ be a conformal immersion. By \cref{diff_p_conformally_invariant}, we can assume
	$\varphi$ is an isometric immersion. We then have a Gauss map $\Phi\colon M\to\Gr_n(\RR^{n+k})$ sending $x\mapsto
	T_xM\subset T_x\RR^{n+k} = \RR^{n+k}$, as well as its analogue on total spaces $\Phi\colon B(M)\to\St_n(\RR^{n+k})$
	defined analogously.

	For $i > \lfloor k/2\rfloor$, we know by \cref{cpx_stiefel,double_cpxif} that
	\begin{equation*}
		[\CS_{P_i^\perp}(\conn^{\mathrm{can}})]\in \H^{2i-1}(\St_n(\RR^{n+k});\RR) 
	\end{equation*}
	is contained in the even integer
	lattice. This property is natural in principal bundles with a connection, and $\conn =
	\Phi^*(\conn^{\mathrm{can}})$, so this is also true for $\CS_{P_i^\perp}(\conn)$.
\end{proof}

We use this to define an $\RR/\ZZ$-valued invariant which obstructs conformal immersions of an orientable
Riemannian $3$-manifold $Y$ into $\RR^4$. The frame bundle $B(Y)\to Y$ admits a section $\chi$; define
\begin{equation}
	\Phi(Y) \colonequals \int_Y \frac 12 \chi^*\CS_{P_1}(\conn)\in\RR/\ZZ \comma
\end{equation}
where $\conn$ is the Levi-Civita connection. A priori this depends on the section, but one can calculate (e.g.\
\cite[\S 6]{cs}) that if $\chi$ and $\chi'$ are two sections, the difference of their pullbacks of the
Chern--Simons invariant consists of torsion and an integer number of copies of an integral cohomology class; the
torsion disappears when we integrate, and the integer-valued cohomology class does not affect the answer mod $\ZZ$.
\Cref{extra_factor_of_two} (and the fact that $P_1^\perp = -P_1$) implies that if $Y$ conformally immerses in
$\RR^4$, then $\Phi(Y) = 0$.

And now the moment we've all been waiting for.

\begin{theorem}[{\cite[\S 6, Example 1]{cs}}]
	\label{nope_RP3}
	The manifold $\RRP^3$ with the round metric does not conformally immerse into $\RR^4$.
\end{theorem}

\begin{proof}
	We will calculate $\CS_{P_1}(\conn)$ for $\conn$ the Levi-Civita connection on $\RRP^3$. The identification
	$\RRP^3 = \SO(3)$ gives us an orthonormal basis $\{v_1, v_2, v_3\}$ of $\mathfrak{so}(3)$, the space of
	left-invariant vector fields; in the Levi-Civita connection, $\nabla_{v_1}v_2 = v_3$, $\nabla_{v_2}v_3 = v_1$, and
	$\nabla_{v_1}v_3 = -v_2$. If $\pi\colon B_{\Or}(\RRP^3)\to\RRP^3$ denotes the bundle of orthonormal frames, the
	above basis gives us a section $\chi$ of $\pi$. We have a formula for $\chi^*\CS_{P_1}(\conn)$~\eqref{CS_fun};
	expanding in coordinates and using the covariant derivatives of the $v_i$s, and we obtain
	\begin{equation}
		\chi^*\left(\frac 12\CS_{P_1}(\conn)\right) = -\frac{1}{2\pi^2} \,\mathrm{vol}\comma
	\end{equation}
	where $\mathrm{vol}$ is the volume form on $\RRP^3$. As a Riemannian manifold, $\RRP^3$ with the round metric is
	the quotient of $\Sph{3}$ with the round metric under the antipodal map, so the volume of $\RRP^3$ is one-half that of
	$\Sph{3}$, i.e.\ $\mathrm{Vol}(\RRP^3) = \pi^2$. Thus $\Phi(\RRP^3) = 1/2$.
\end{proof}
%-------------------------------------------------------------------%
%  ℝP³                                                              %
%-------------------------------------------------------------------%
%
%\subsubsection{\texorpdfstring{$\RRP^3$}{ℝP³}}
%TODO: in this section, we will go into the additional integrality result (probably just sketching the proof) and
%then show $\RRP^3$ does not conformally immerse in$\RR^4$.
%
%In the Chern--Simons paper, they prove that $\RRP^3$ cannot conformally immerse in $\RR^4$. One would hope to show this by showing $\phat_1^\perp(\Tan\RRP^3)$ is non vanishing. Note that $\phat_1^\perp = -\phat_1$. Also note that we are on a $3$ manifold so the characteristic class and curvature of $\phat_1$ both vanish, but there is a little more information in the differential character. One can compute this it by choosing a section of the orthogonal frame bundle of $\RRP^3$, and integrating the pullback of the Chern--Simons form by this section. Unfortunately this integral gives $1$, and we are supposed to consider it as an element of $\RR/\ZZ$. In the paper they go further and show that it is actually well defined mod $2\ZZ$. At least thats what I understood.
%
There are numerous examples in the literature of calculations of this sort to obtain conformal nonimmersion
results: see \cite{HL74, APS2, Mil75, Don77, Tsu81, Bac82, Tsu84, Ouy94, MM01, MZ10, PT10, Li15} for some examples.
\begin{remark}[(transgression and inverse transgression)]
\label{transgression_detail}
Here we go into a little more detail about the transgression and inverse transgression maps, the latter of which
appeared in the proof of \cref{double_cpxif}. We follow \cites[\S 9]{Bor55}{cs}.
\begin{definition}
Let $F\overset i\to E\overset{\pi}{\to} B$ be a fiber bundle, $x\in \H^k(F; A)$, and $y\in \H^{k+1}(B; A)$. We say
that \textit{$x$ transgresses to $y$}\index[terminology]{transgression} when there is a cochain $c\in Z^k(F; A)$
such that $[i^*(c)] = x$ and $\delta c = \pi^*b$ for some cocycle $b$ in the cohomology class of $y$.
\end{definition}
Given $x$, $y$ may not exist, and may not be unique if it exists. Transgression is natural under pullback of fiber
bundles, so when studying transgression in principal $G$-bundles, it makes sense to work universally in $G\to
\EG\to\BG$.

Transgression has something to say about the Serre spectral sequence\index[terminology]{Serre
spectral sequence} for the fiber bundle $F\to E\to B$. We can identify $x$ and $y$ with their images on the
$E_2$-page, in $E_2^{0,k}$ and $E_2^{k+1,0}$ respectively. Transgression as defined above is equivalent to asking
that
\begin{enumerate}
	\item no differential $d_r$ for $r < k+1$ kills $x$ or $y$, so that their images in the $E_{k+1}$-page are
	nonzero; and
	\item $d_{k+1}(x) = y$.
\end{enumerate}
The Serre spectral sequence is first-quadrant, so $d_{k+1}$ is the last differential that could kill $x$ or $y$.
In the bundle $G\to \EG\to \BG$, all positive-degree elements must be killed by differentials, because $\EG$ is
contractible; this is another indication that transgression is important here.\footnote{Similarly, when $A$ is an
abelian group, there is a fibration $\mathrm K(A, n)\to E\to \mathrm K(A, n+1)$, where $E$ is contractible, and a
theorem of Borel \cite[Theorem 13.1]{Bor53} on transgression is a crucial part of Serre's calculation \cite{Ser53}
of the cohomology of Eilenberg--MacLane spaces.\index[terminology]{Eilenberg--MacLane space}} When $G$ is a
connected Lie group, transgression is often as nice as it can be: $\H^*(G;A)$ is an exterior algebra on odd-degree
generators $x_1,\dotsc,x_n$, $\H^*(\BG;A)$ is a polynomial algebra on even-degree generators $y_1,\dotsc,y_n$, and
$x_i$ transgresses to $y_i$. Here $A$ may be $\QQ$, $\ZZ/p$, or $\ZZ$ depending on $G$; for example, when $G =
\Uup_n$, we can use $\ZZ$ coefficients. In these settings we can begin to see how to define the inverse
transgression map: ignoring gradings, the only differences between the rings $\H^*(\BG;A)$ and $\H^*(G;A)$ are the
relations $x_i^2 = 0$, so we can think of transgression as a map $\H^*(G;A)\to \H^{*+1}(\BG; A)$ whose image is
everything not containing terms of the form $y_i^m$ for $m > 1$. Thus we can define an inverse transgression
map\index[terminology]{inverse transgression map} $\tau$ by sending $y_i\mapsto x_i$ and $y_i^2 = 0$.

Chern--Simons \cite[\S 5]{cs} define $\tau$ differently, and more directly: given $y\in \H^{k+1}(\BG;A)$, let $b$ be
a cocycle representative for $y$ which vanishes when pulled back to any point of $\BG$; since $\EG$ is
contractible, $\pi^*(b) = \delta c$ for some $c\in Z^k(\EG; A)$. Then $\tau(y)$ is defined to be the cohomology
class of the restriction of $c$ to a fiber; one has to check this is well-defined, but it is. When $\H^*(G; A)$ is
an exterior algebra on odd-degree generators, this definition recovers the definition from the previous paragraph,
but this definition is more general. 
It is natural in $G$, and $\tau(y^2) = 0$ follows because if we choose $b,c$
as above, then
\begin{equation*}
	\delta(b\cupprod c) = \pi^*(b\cupprod b) \comma
\end{equation*}
and restricted to a fiber, $b\cupprod c$ vanishes.

From here it is natural to wonder whether the inverse transgression map admits a differential refinement
$\hat\tau\colon\Hhat^4(\BunGnabla ;\ZZ)\to\Hhat^3(G;\ZZ)$. This is true, and there are constructions of this map due
to Carey--Johnson--Murray--Stevenson--Wang \cite[\S 3]{CJMSW05} and Schreiber \cite[1.4.1.2]{Urs}.

Chern--Simons (\textit{ibid.}, \S 3) also discuss transgression in the context of the Chern--Simons form and when
the fiber bundle is a principal $G$-bundle $P\to M$ with connection $\conn$. Fixing an invariant polynomial $f$,
they use the Maurer--Cartan form\index[terminology]{Maurer--Cartan form} on $G$ to define a class in $\HdR^*(G)$
which transgresses to $[f(\curvature\conn)]\in\HdR^*(M)$.
\end{remark}

\newpage
%!TEX root = ../diffcoh.tex

%-------------------------------------------------------------------%
%-------------------------------------------------------------------%
%  Charge quantization                                              %
%-------------------------------------------------------------------%
%-------------------------------------------------------------------%

\section{Charge quantization}
\textit{Talk by Dan Freed}\\
\textit{Notes by Arun Debray}
%These notes were typeset by Arun Debray, following a lecture by Dan Freed. Any mistakes or typos are probably due to me (i.e. Arun).
\label{field_theory}

There are a few different applications of differential cohomology to quantum physics; today, we'll focus on charge
quantization, using Maxwell theory as an example. First, in \S\ref{classical_Maxwell}, we introduce classical
Maxwell theory, formulated in the language of differential forms. Then, in \S\ref{quantum_Maxwell}, we pass to the
quantum theory. This imposes integrality conditions on differential forms, leading to the appearance of
differential cohomology. This lecture is based on~\cite[Part 3]{Fre02a}.

The history of the use of differential cohomology to implement charge quantization is closely tied to the
development of the theory of differential cohomology itself. Alvarez~\cite{Alv85} was the first to use
differential cohomology in this context, though he does not use the words ``differential
cohomology.''\footnote{Alvarez also uses differential cohomology to characterize quantized topological terms. This
is a related but different application of differential cohomology to physics, and is more closely related to the
discussion of invertible field theories in the next chapter. See Deligne--Freed~\cite[Chapter 6]{DF99} for a
mathematical exposition of topological terms and their relationship to differential cohomology, along with more
recent work of
Davighi, Gripaios, and Randal-Williams~\cite{DG18a, DGR20} and Córdova--Freed--Lam--Seiberg~\cite{CFLS20a, CFLS20b}
from a more physics-based perspective.}
Gawędzki~\cite{Gaw88} then explicitly brings in differential cohomology in the form of Deligne cohomology.

The original motivation to consider generalized differential cohomology came from charge quantization in string
theory: work of Minasian--Moore~\cite{MM97}, Sen~\cite{Sen98}, and Witten~\cite{Wit98} argued that D-brane charges
and Ramond--Ramond field strengths are valued in \Ktheory,\footnote{See~\cite{FW99, MW00, DMW02} for some
related work.} leading to a search for a $\Kup$-theoretic analogue of differential
cohomology.\index[terminology]{D-brane charge}\index[terminology]{Ramond--Ramond field} Freed--Hopkins~\cite{FH00}
first provided a definition of differential \Ktheory for this purpose, and Freed~\cite{Fre00} considers more
general differential generalized cohomology theories. Hopkins--Singer~\cite{HopkinsSinger}, who comprehensively
studied differential generalized cohomology theories, write that they originally began their project to investigate
string-theoretic phenomena.\footnote{Similarly, twisted differential cohomology was first motivated by the
appearance of examples of twisted differential \Ktheory in string theory~\cites[\S 5.3]{Wit98}{BV00}{Fre00},
and has since become an object of study in its own right~\cites[\S
4.1.2]{Urs}{GS18}{BN19}{GS19c}{GS19a}{GS19b}{FSS20d}.}

%-------------------------------------------------------------------%
%-------------------------------------------------------------------%
%  Classical Maxwell theory                                         %
%-------------------------------------------------------------------%
%-------------------------------------------------------------------%

\subsection{Classical Maxwell theory}
\label{classical_Maxwell}
\index[terminology]{Maxwell theory!classical}

Let $(N, g_N)$ be a Riemannian $3$-manifold without boundary and $M=\RR\times N$. Let $t$ be the $\RR$
coordinate, so we give $M$ the Lorentz metric
\begin{equation}
	g_M= \d t^2 - g_N \period
\end{equation}
Choose differential forms $E\in\Omega^1(N)$ and $B\in\Omega^2(N)$, respectively the electric and magnetic fields;
also choose the \emph{charge density} $\rho_E\in\Omegac^3(N)$, and the \emph{current}
\index[terminology]{charge density}
\index[terminology]{current}
$J_E\in\Omegac^2(N)$.% 
\footnote{Here $\Omegac^k(X)$ denotes the space of compactly supported $k$-forms on $X$.} 
If $\star_N$ denotes the Hodge star on $N$, then \emph{Maxwell's equations}, as you might see them on a t-shirt, are
\index[terminology]{Maxwell's equations}
\begin{align*}
	\d B &= 0 \\
	 \frac{\partial B}{\partial t} + \d E &= 0\\
	\d  {\star_N} E &= \rho_E \\
	 {\star_N} \frac{\partial E}{\partial t} - \d {\star_N} B &= J_E \period
\end{align*}
Writing $F= B - \d t\wedge E\in\Omega^2(M)$ and $j_E= \rho_E + \d t\wedge J_E\in\Omega^3(M)$, we
obtain a more concise form of Maxwell's equations:
\begin{equation}
	\d F = 0 \andeq \d {\star_M} F = j_E \period
\end{equation}
Now we include topology. We just saw that $j_E$ is exact, so it cannot define an interesting de Rham cohomology
class, but $F$ is closed, so may be interesting. Define the \emph{charge} at time $t$ to be the de Rham class
\begin{equation}
	Q_E = [j_E|_{\{t\}\times N}]\in \Hc^3(N;\RR) \period
\end{equation}
This is in the kernel of the map $\Hc^3(N;\RR)\to \H^3(N;\RR)$; hence, on a compact manifold, $Q_E = 0$.

Let $W$ be the worldline of a charged particle with electric charge $q_E\in\RR$. Then $j_E = q_E\cdot\delta_W$,
where $\delta_W$ is the ``current sitting at $W$.'' We have two ways of making sense of this.
\begin{itemize}
	\item First, we could take $\delta_W$ to be a current in the de Rham sense, akin to a differential form but
	built with distributions instead of smooth functions. Amusingly, this is a current in both the Maxwell and de
	Rham senses. This is a typical example of a current in electromagnetism.
	\item Alternatively, we could take $\delta_W$ to be an honest $3$-form Poincaré dual to $W$. In this case we
	can choose $\delta_W$ to be supported in an arbitrary neighborhood of $W$.
\end{itemize}
One more ingredient in Maxwell theory, though not strictly necessary, is an action principle. This follows the
Lagrangian formulation of physics: we aim to find a variational problem whose solutions are the Maxwell equations.
We add an assumption from classical physics: that $[F] = 0$ in $\HdR^2(M)$; this means there are no magnetic
monopoles.
\index[terminology]{action principle}

This assumption also implies $F = \d A$ for some $1$-form $A$ called the \emph{electromagnetic potential}.
\index[terminology]{electromagnetic potential} This is not unique, but its class in $\Omega^1(M)/\Omegacl^1(M)$
(i.e.\ up to closed $1$-forms) is unique.  Then, the \emph{classical action} \index[terminology]{classical action}
of Maxwell theory is
\index[terminology]{electromagnetic potential}
\index[terminology]{Lagrangian!for Maxwell theory}
\begin{equation}
\label{classical_action}
	S = \int_M -\frac{1}{2} \d A\wedge {\star} \d A + A\wedge j_E \period
\end{equation}
Since $M$ is noncompact, this could be infinite, but we're just interested in its first variation anyways, which is
well-behaved.
\begin{exercise}
Show that the Euler--Lagrange equation for \eqref{classical_action} is $\d{\star} F = j_E$. (We already assumed $\d
F = 0$, the other half of Maxwell's equations.)
\index[terminology]{Euler--Lagrange equation}
\end{exercise}
One caveat: defining the action requires $A$ to be in $\Omega^1(M)$, not
$\Omega^1(M)/\Omegacl^1(M)$. This ends up not a problem; adding a closed form to $A$ does not change
the Euler--Lagrange equation.

%-------------------------------------------------------------------%
%-------------------------------------------------------------------%
%  Quantum Maxwell theory                                           %
%-------------------------------------------------------------------%
%-------------------------------------------------------------------%

\subsection{Quantum Maxwell theory}
\label{quantum_Maxwell}
\index[terminology]{Maxwell theory!quantum}

In the quantum theory, we allow magnetic monopoles. Dirac \cite{Dir31} argues that this forces electric and
magnetic charges to be \emph{quantized}, i.e.\ taking values in a discrete subgroup of $\RR$. This is how
differential cohomology enters the picture.

So assume $N = \RR^3$ with the usual Euclidean metric, and introduce a magnetic monopole of charge $q_B\in\RR$ at
the origin. Then we have a \emph{magnetic current} \index[terminology]{magnetic current} $j_B = q_B\cdot\delta_0$.
The condition that $\d F = 0$ is modified to
\begin{equation}
\label{dFmon}
	\d F = q_B\cdot\delta_0 \period
\end{equation}

The input to the path integral is the exponentiated action $\exp(iS/\hbar)$ (where $S$ is
as in~\eqref{classical_action}. However, this is not quite consistent with~\eqref{dFmon} --- there is a problem at
the origin. On $\RR\times(\RR^3\smallsetminus 0)$, we can write $F = \d A$, and therefore realize $F$ as the curvature of
a connection $A$ on a principal $\RR/q_B\ZZ$-bundle $P$. The characteristic class of $P$ is
\begin{equation}
	[P]\in \H^2(\RR\times (\RR^3\smallsetminus 0); q_B\ZZ)\cong \H^2(S^2;q_B\ZZ) = q_B\ZZ \comma
\end{equation}
and $[P]$ is a generator of this abelian group.

The space of fields in the quantum theory is the groupoid of principal $\RR/q_B\ZZ$-bundles with connection. Now we
can revisit the action~\eqref{classical_action} --- it doesn't have to make sense as is (e.g.\ $A$ isn't exactly a
$1$-form), but we do want $\exp(iS/\hbar)$ to make sense.

Let's work on a general $4$-manifold $X$. To avoid causality issues, let's make $X$ a Riemannian manifold, rather
than a Lorentz one. Assume $j_E$ is Poincaré dual to some loop $\gamma\subset X$. If there is a $q_E$ charge moving
along this loop, then
\begin{equation}
	\int_M A\wedge j_E = \oint_\gamma q_EA = q_E\Hol_\gamma(A) \period
\end{equation}
Now $\Hol_\gamma(A)\in\RR/q_B\ZZ$, so the quantity
\begin{equation}
	\exp\left(\frac{i}{\hbar} q_E\Hol_\gamma(A)\right)
\end{equation}
is well-defined if and only if
\begin{equation}
	\frac{1}{\hbar} q_Eq_B\in 2\pi\ZZ \period
\end{equation}
This is Dirac's quantization condition. Thus integrality enters a story told with differential forms; this is
already suggestive of differential cohomology!
\index[terminology]{Dirac quantization}

To say it more explicitly, the space of quantum fields is the stack $\Bunnabla_{\RR/q_B\ZZ}(X)$; the set of
isomorphism classes of objects is $\Hhat^2(X; q_B\ZZ)$. The curvature map lands in those $2$-forms with periods in
$q_B\ZZ$, giving us a short exact sequence we've seen before:
\index[terminology]{curvature map}
\begin{equation*}
	\begin{tikzcd}
		0 \arrow[r] & \H^1(X;\RR/q_B\ZZ) \arrow[r] & \Hhat^2(X;q_B\ZZ) \arrow[r, "\curv"] & \Omegacl^2(X)_{q_B\ZZ} \arrow[r] & 0 \period
	\end{tikzcd}
\end{equation*}
The classical fields $\Omega^1(X)/\Omegacl^1(X)$ sit as a subspace in $\Hhat^2(X; q_B\ZZ)$; the cokernel
is $\H^2(X;q_B\ZZ)$ modulo torsion, indicating the new information in the quantum theory.

Another interesting upshot is that since the kernel of the curvature map corresponds to the flat connections, i.e.\
those on which $F$ is boring, the electric flux really lives in $\Hhat^2(X; q_B\ZZ)$. This is new. The flat
connections are new, too --- even if you don't usually get to observe them, they manifest in the physics, e.g.\
through the Aharonov--Bohm effect. And all of this is still ``semiclassical,'' i.e.\ about the input to the path
integral, before we try to evaluate said path integral.
\index{Aharonov--Bohm effect}

\begin{remark}
	One important clarification: $F$ is not a differential cohomology class; it's the curvature of an actual bundle
	with connection, not an equivalence class. So really we need a cochain model: bundles and connections glue, but
	equivalence classes don't. 
	Cheeger--Simons characters aren't built in this way, so for physics applications one must
	do something different.
\end{remark}

Now we revisit the electric charge, a closed $3$-form.
%\footnote{You might be wondering where the compact support
%condition went. To work in Euclidean signature, rather than Minkowski signature, we must \emph{Wick-rotate} the
%theory, a nontrivial procedure which ultimately removes the requirement for compact
%supports.\index[terminology]{Wick rotation}}
Because
$(i/\hbar)j_Ej_B\in 2\pi\ZZ$, we'd like to impose that $[j_E]\in \HdR^3(X)$ is also in the image of the map
$\H^3(X;q_E\ZZ)\to \H^3(X;\RR)$, i.e.\ that we're in the homotopy pullback, which is $\Hhat^3(X;q_E\ZZ)$. Again,
though, we want a local object in the end, not just its isomorphism class.

We can also rewrite one term in the exponentiated action in terms of differential cohomology, as
\begin{equation}
	\exp\left(\frac i\hbar \int_X \Fhat\cdot\jhat\right) \period
\end{equation}
Here $\Fhat$ and $\jhat$ are the differential cohomology refinements of $F$ and $j_E$, respectively.  The
product $\cdot$ is the cup product from \cref{sec:Delignecup}, which is a map
\begin{equation}\label{differential_action}
	\Hhat^2(X; q_B\ZZ)\otimes\Hhat^3(X;q_E\ZZ)\longrightarrow \Hhat^5(X; q_Eq_B\ZZ) \period
\end{equation}
Since $X$ is a $4$-manifold, the integration map has degree $-4$, so is of the form
\begin{equation}
	\int_X\colon \Hhat^5(X; q_Eq_B\ZZ)\longrightarrow \Hhat^1(\pt; q_Eq_B\ZZ)\cong \RR/q_Eq_B\ZZ \period
\end{equation}

\begin{exercise}
	Show that if $\Fhat$ is topologically trivial, meaning that it comes from a connection on a trivial vector
	bundle, or equivalently that its image under the characteristic class map vanishes, then
	$\Fhat\cdot\jhat$ is also topologically trivial.
\end{exercise}

\begin{remark}
	There are many variations of this story in field theory and string theory, generally for abelian gauge fields.
	For example, $F$ might have some other degree, or even be inhomogeneous. Dirac charge quantization still
	applies, and will refine $F$ to an appropriate differential cohomology group.\index[terminology]{string theory}

	More recently, people realized that this story sometimes yields generalized differential cohomology theories.
	Understanding which cohomology theory one obtains is a bit of an art --- physics tells you some constraints,
	but not an algorithm. For example, this happens in superstring theory: the Ramond--Ramond
	field\index[terminology]{Ramond--Ramond field} is realized in differential \Ktheory \cite{FH00,
	MW00}\index[terminology]{differential K-theory@differential \Ktheory}, and the $B$-field in a differential
	refinement of (a truncation of) $\ko$ \cite{DFM11a, DFM11b}. These and other refinements of Dirac
	quantization to generalized differential cohomology are also studied in \cite{BM06a, BM06b, DFM07, Fre08,
	Sat10, SZ10, Sat11, SSS12, KM13, KV14, DMDR14, FSS15, FR16, GS19, Sat19, FRRB20}. The choice of generalized
	cohomology theory is not always an exact science: for example, there are different proposals for the $C$-field
	in M-theory. Witten \cite[\S 2.3]{Wit97} argues that the $C$-field should be quantized in $w_1$-twisted
	degree-$4$ ordinary differential cohomology, which passes consistency checks for various possible
	anomalies \cites[\S 4]{Wit97}[\S 4]{Wit16}{FH21a}; there is also the ambitious ``hypothesis H'' of
	Fiorenza--Sati--Schreiber \cite{Sat18, FSS20a, FSS19} proposing that the $C$-field in M-theory is quantized
	using a differential refinement of $\image(J)$-twisted stable cohomotopy instead. Work of Fiorenza, Sati,
	Schreiber, and their collaborators \cite{SS19, FSS20a, GS20, SS20b, SS20a, BSS21, FSS19, FSS20b, SS21, SS20c,
	FSS20c} and Roberts \cite{Rob20} recovers as consequences of hypothesis H several things physicists predicted
	to be true about M-theory.\index[terminology]{hypothesis
	H}\index[terminology]{C-field@$C$-field}\index[terminology]{M-theory}
\end{remark}

If we consider Maxwell theory with both electric and a magnetic currents, the theory has an
``anomaly,''\index[terminology]{anomaly} meaning that some quantity that we'd like to obtain as a complex number is
actually an element of a complex line that's not trivialized (and in some cases cannot be trivialized canonically
for all manifolds of a given dimension). Differential cohomology also provides a perspective on the anomaly. The
expression $\Fhat\cdot\jhat_E$ in \eqref{differential_action} is valid if there's electric current but
not magnetic current; if $\jhat_B\neq 0$, then $F$ isn't closed, hence isn't the curvature of a line bundle.
But $\jhat_B$ is also quantized, hence represents a differential cohomology class, and we can ask for
$\Fhat$ to trivialize $\jhat_B$. Now the action is
\begin{equation}
	\exp\left(\frac i\hbar \int_X \Fhat\cdot \jhat_E\jhat_B\right) \period
\end{equation}
Since $\Fhat\cdot\jhat_E\jhat_B\in\Hhat^6$, integrating brings us to $\Hhat^2(\pt;
q_Eq_B\ZZ)$, yielding the complex line which signals the anomaly. More on this anomaly can be found in
Freed--Moore--Segal \cite{FMS07a, FMS07b}.
% redundant in view of invertible field theories section
%\begin{example}
%In the last few minutes, we'll discuss a different example of differential cohomology in physics. Suppose $M$ is an
%oriented Riemannian $3$-manifold and $P\to M$ is a principal $SU_2$-bundle with connection $\Theta$. The second
%Chern class admits a differential refinement $\ccech_2(\Theta)\in \Hhat^4(M)$, and
%\begin{equation}
%	\int_M \ccech_2(\Theta) \in \Hhat^1(\pt)\cong\RR/\ZZ.
%\end{equation}
%Hence this is the sort of thing you can add to an action. It's an example of a \emph{Chern--Simons term} in 3d QFT.
%\index[terminology]{Chern--Simons term}
%The de Rham class underlying $\ccech_2$ is sometimes called the level of the theory, and the fact that it must
%refine to differential cohomology is saying the level is quantized. In general, quantization of coupling constants
%%provides another instance of differential cohomology in physics.
%
%Chern--Simons terms are usually described without differential cohomology, using the Chern--Simons form associated to
%$\Theta$, but writing that term on $M$, rather than on $P$, requires a choice of a section of $P$, and we don't
%always have that.
%\end{example}

\newpage
%!TEX root = ../diffcoh.tex

%-------------------------------------------------------------------%
%-------------------------------------------------------------------%
%  Invertible field theories                                        %
%-------------------------------------------------------------------%
%-------------------------------------------------------------------%

\section{Invertible field theories}
\textit{by Arun Debray}
\label{invertible_field_theories}

Freed--Hopkins \cite[\S 5.4]{FH21} conjecture a different application of generalized differential cohomology to
field theory, describing reflection-positive invertible field theories which are not necessarily topological. In
this chapter we go over this conjecture. This story is similar to an established theorem, Freed--Hopkins'
classification of reflection-positive invertible \emph{topological} field theories \cite{FH21}, so we begin in
\cref{top_IFT} by going over that classification; then in \cref{non_top_field_theory} we generalize to the
nontopological setting.
\subsection{Topological invertible field theories}
\label{top_IFT}
\begin{definition}
Let $\rho(n)\colon H_n\to\Or_n$ be a Lie group homomorphism. An \textit{$H_n$-structure} on a smooth manifold
$M$ is a principal $H_n$-bundle $P\to M$ together with an isomorphism of principal $\Or_n$-bundles \[\theta\colon
P\times_{H_n}\Or_n\isomorphism \mathcal B_{\Or}(M),\] where $\mathcal B_{\Or}(M)$ is the frame
bundle of $M$.
\index[terminology]{Hn-structure@$H_n$-structure}
\index[terminology]{frame bundle}
\end{definition}
%
%
%Let $\rho\colon H\to\Or$ be a homomorphism of topological groups. We may then speak of ``manifolds with
%$H$-structure,'' meaning manifolds $M$ with a lift of the stable tangent bundle map $TM\colon M\to B\Or$ across
%$\rho$:
%\begin{equation}
%% https://q.uiver.app/?q=WzAsMyxbMCwxLCJNIl0sWzEsMSwiQlxcbWF0aHJtIE8iXSxbMSwwLCJCIl0sWzIsMSwiXFx4aSJdLFswLDEsIlRNIiwyXSxbMCwyLCIiLDEseyJzdHlsZSI6eyJib2R5Ijp7Im5hbWUiOiJkYXNoZWQifX19XV0=
%\begin{tikzcd}
%	& BH \\
%	M & {B\mathrm O}
%	\arrow["\rho", from=1-2, to=2-2]
%	\arrow["TM"', from=2-1, to=2-2]
%	\arrow[dashed, from=2-1, to=1-2]
%\end{tikzcd}
%\end{equation}
For example, an $\SO_n$-structure is equivalent data to an orientation, a $\Spin_n$-structure is equivalent to a
spin structure, and so forth.

An $H_n$-structure on a manifold $M$ induces an $H_n$-structure on $\partial M$, and we may therefore consider
bordism groups $\Omega_n^{H}$ of $H_n$-manifolds, as Lashof \cite{Las63} did, and their categorified analogues:
bordism $(\infty, n)$-categories $\Bord_n^H$ of $n$-manifolds with $H_n$-structure, such as the bordism categories
constructed by Lurie \cite{Lur09}, Schommer-Pries \cite{SP17}, and Calaque--Scheimbauer \cite{CS19}.
\index[terminology]{bordism groups}
\index[notation]{OmegaH@$\Omega_n^H$}
\index[terminology]{bordism category}
\index[notation]{BordnH@$\Bord_n^H$}

Recall that a topological field theory (TFT) is a symmetric monoidal functor
\index[terminology]{topological field theory}
\index[terminology]{TFT|see {topological field theory}}
\begin{equation}
	Z\colon\Bord_n^H\to\fC \comma
\end{equation}
where $\fC$ is some symmetric monoidal $(\infty, n)$-category. The $\infty$-category of TFTs is symmetric monoidal
under ``pointwise tensor product:''\index[terminology]{tensor product!of topological field theories}
\[(Z_1\otimes Z_2)(M) \colonequals Z_1(M)\otimes Z_2(M)\period\]

\begin{definition}[{(Freed--Moore \cite{FM06})}]
	A TFT 
	\begin{equation*}
		Z\colon\Bord_n^H\to\fC
	\end{equation*}
	is \textit{invertible} if there is some other TFT $Z^{-1}$ such that $Z\otimes
	Z^{-1}$ is isomorphic to the trivial theory (i.e.\ the constant functor valued in $\mathbf{1}_\fC$).
	\index[terminology]{topological field theory!invertible}
	\index[terminology]{invertible topological field theory|see {topological field theory!invertible}}
\end{definition}

\noindent Equivalently, $Z$ carries objects of $M$ to $\otimes$-invertible objects in $\fC$ and $k$-morphisms to
composition-invertible $k$-morphisms in $\fC$ for all $k$. In many cases it suffices to check invertibility on a
subset of objects, such as certain spheres \cite{432876} or tori \cite{SP18}.

\begin{example}[(Euler theories)]
Let $\lambda\in\CC^\times$. The \emph{Euler theory}
\begin{equation*}
	Z_\lambda\colon\Bord_{n, n-1}^{\Or}\to\Vect_\CC
\end{equation*}
is an invertible TFT which to every object assigns the vector space $\CC$, and to every bordism $X\colon M_1\to M_2$ assigns
multiplication by $\lambda^{\chi(X, M_1)}$. These compose properly because the Euler characteristic satisfies a
gluing formula.\index[terminology]{Euler theory}
\end{example}
Freed--Hopkins--Teleman \cite{FHT10} classified invertible TFTs using work of
Galatius--Madsen--Tillmann--Weiss \cite{GMTW09} and Nguyen \cite{Ngu17}. Freed--Hopkins \cite{FH21} went further:
they studied \textit{reflection-positive invertible TFTs}, which have additional structure. This structure is
related to the notion of unitarity in quantum field theory, so invertible TFTs appearing in the study of unitary
QFTs should have reflection-positive structures.
\index[terminology]{unitarity!in quantum field theory}
\index[terminology]{topological field theory!reflection-positive invertible}

Let $\MTH$ denote the Thom spectrum of $-\mrm{B}\rho\colon \mrm{B}H\to \BO$.\footnote{There is an important
subtlety here: we started with $\rho(n)\colon H_n\to\Or_n$, not the stabilized version $\rho\colon H\to\Or$.
Freed--Hopkins \cite[Theorem 2.19]{FH21} show that the additional data associated to reflection positivity allows
one to define $\rho$ and $H$ such that $\rho(n)\colon H_n\to\Or_n$ is the pullback of $\rho\colon H\to\Or$ along
the inclusion $\Or_n\hookrightarrow\Or$.}
\index[terminology]{Thom spectrum}%
\index[notation]{MTH@$\MTH$}%
Thom's collapse map identifies the homotopy groups
of $\MTH$ with the bordism groups of manifolds with $H_n$-structure \cites[Théorème
IV.8]{ThomThesis}{Pon55}[Theorem C]{Las63}.\footnote{The use of $-\rho$ ensures that we obtain an $H$-structure on
the stable tangent bundle.  Homotopy theorists more traditionally study the Thom spectrum of $\rho$, denoted
$\mathrm{MH}$, which corresponds to bordism of manifolds with an $H$-structure on the stable \emph{normal} bundle.
Often $\MTH\simeq\mathrm{MH}$, as is the case for $\mathrm{MTO}$, $\mathrm{MTSO}$, $\mathrm{MTSpin}$,
$\mathrm{MTSpin}^c$, $\mathrm{MTString}$, and $\mathrm{MTU}$, but not always:
$\mathrm{MTPin}^+\not\simeq\mathrm{MPin}^+$.} Let $\IZ$ denote the \textit{Anderson dual of the sphere
spectrum} \cite{And69, Yos75}, which satisfies the universal property that there is a short exact sequence
\index[terminology]{Anderson dual of the sphere spectrum}
\index[notation]{IZ@$\IZ$}
\begin{equation}
	% https://q.uiver.app/?q=WzAsNSxbMCwwLCIwIl0sWzEsMCwiXFxidWxsZXQiXSxbMiwwLCJcXGJ1bGxldCJdLFszLDAsIlxcYnVsbGV0Il0sWzQsMCwiMCJdLFswLDFdLFsxLDJdLFsyLDNdLFszLDRdXQ==
\begin{tikzcd}
	0 & \Ext(\uppi_{n-1}(X), \ZZ) & {[X, \Sigma^n \IZ]} & \Hom(\uppi_n(X), \ZZ) & 0 \comma
	\arrow[from=1-1, to=1-2]
	\arrow[from=1-2, to=1-3]
	\arrow[from=1-3, to=1-4]
	\arrow[from=1-4, to=1-5]
\end{tikzcd}
\end{equation}
which noncanonically splits.

\begin{theorem}[{(Freed--Hopkins \cite{FH21})}]
	There is an isomorphism of abelian groups from $\uppi_0$ of the space reflection-positive, invertible,
	$n$-dimensional, topological field theories to the torsion subgroup of $[\MTH, \Sigma^{n+1}\IZ]$.
\end{theorem}

\begin{remark}
	Any classification of TFTs $Z\colon\Bord_n^H\to\fC$ depends on what we take $\fC$ to be. For this theorem,
	Freed--Hopkins make an ansatz about the choice of $\fC$. Example $\fC$ meeting this ansatz are known in category
	number $2$ and below: see \cite[Theorem 1.52]{Vienna} and \cite[Proposition 4.21]{DG18}.
\end{remark}

If $B$ admits a CW structure with finitely many cells in each dimension, so that the homotopy groups of $\MTH$ are
finitely generated, then
\begin{equation*}
	\Tors([\MTH, \Sigma^{n+1}\IZ])\cong \Tors(\Hom(\uppi_n(\MTH),\CC^\times)) \period
\end{equation*}
Thus we have identified $\Tors([\MTH, \Sigma^{n+1}\IZ])$ with the group of torsion
$\CC^\times$-valued bordism invariants for $n$-dimensional $H$-manifolds.  Given such a bordism invariant
$\varphi$, it is possible to choose a reflection-positive invertible TFT $Z$ in the component of
$\uppi_0(\mathrm{ITFT}s)$ corresponding to $\varphi$ such that the partition function of $Z$ is equal to $\varphi$.

\begin{example}[{(classical Dijkgraaf--Witten theory \cite{DW90, FQ93})}]
	\label{classical_DW}
	\index[terminology]{Dijkgraaf--Witten theory!classical}
	Let $G$ be a group and $\lambda\in \H^n(\BG;\QQ/\ZZ)$. Then $\lambda$ defines a bordism invariant of oriented
	$n$-manifolds $M$ with a principal $G$-bundle $P$ by integrating, then exponentiating:
	\begin{equation}
		(M, P)\longmapsto \exp\paren{2\pi i \int_M \lambda(P)}\in\CC^\times\comma
	\end{equation}
	where $\lambda(P)$ denotes the pullback of $\lambda$ along the map $M\to \BG$ defined by $P$. Stokes' theorem
	\index[terminology]{Stokes' theorem}
	implies this is a bordism invariant, and it is torsion; therefore~\eqref{classical_DW} is the partition
	function of a unique (up to isomorphism) reflection-positive invertible TFT. This TFT is called
	\textit{classical Dijkgraaf--Witten theory}. The state space assigned to any codimension-$1$ manifold is
	noncanonically isomorphic to $\CC$; see Freed--Quinn \cite[\S 1]{FQ93} for a fuller description and
	Yonekura \cite[\S 4]{Yon19} for another construction.
\end{example}

\begin{example}[(Arf theory)]
	\label{arf_TFT}
	\index[terminology]{Arf theory}
	\index[terminology]{Arf invariant}
	We have $\Omega_2^{\Spin}\cong\ZZ/2$, and the \textit{Arf invariant} is a complete invariant
	\begin{equation*}
		\Arf\colon \Omega_2^{\Spin}\to\{\pm 1\}
	\end{equation*}
	\cite[Proposition (4.1)]{Ati71}. 
	Using Freed--Hopkins' classification, there is a
	reflection-positive invertible TFT $Z_A\colon\Bord_2^{\Spin}\to\fC$, called the \textit{Arf theory}, whose
	partition function is the Arf invariant, and $Z_A$ is unique up to isomorphism. Gunningham \cite[Example
	2.19]{Gun16} showed that we can take $\fC$ to be $\categ{sAlg}_\CC$, the Morita bicategory of complex
	superalgebras.
	\index[notation]{sAlgC@$\categ{sAlg}_\CC$}
	\index[terminology]{Morita bicategory}

	As in \cref{classical_DW}, we can recast this example as integration, this time in generalized cohomology.
	Atiyah--Bott--Shapiro \cite{ABS64} showed that spin manifolds admit pushforward maps for $\KO$-theory. On a
	spin surface, the partition function of the Arf theory (i.e.\ the Arf invariant) is the pushforward
	\index[terminology]{Atiyah--Bott--Shapiro map}
	\begin{align*}
			\exp2 \pi i\int_\Sigma^{\KO}\colon \KO^0(\Sigma) &\longrightarrow \KO^{-2}(\pt)\cong\{\pm 1\}\\
		1 &\longmapsto Z_A(\Sigma)
	\end{align*}
	That is, the $\KO$-theoretic pushforward lands in $\ZZ/2$, and exponentiation brings us to $\{\pm
	1\}\subset\CC^\times$.

	Something similar also works in positive codimension! Let $C$ be a closed spin $1$-manifold.
	\begin{align}
	\label{curve_KO_int}
		\int_C^{\KO}\colon \KO^0(C) &\longrightarrow \KO^{-1}(\pt)\cong\ZZ/2\\
		1 &\longmapsto Z_A(C) \period
	\end{align}
	This $\ZZ/2$ is different --- we interpret it as the group of isomorphism classes of complex super lines
	$\{\CC, \Pi\CC\}$ under tensor product. That is, an invertible field theory valued in $\categ{sAlg}_\CC$
	assigns to a codimension-$1$ manifold a $\otimes$-invertible complex super vector space; up to isomorphism this
	is either the even line or the odd line, and~\eqref{curve_KO_int} tells us which one the Arf theory assigns to
	$C$. For example, the bounding spin circle is assigned an even line, and the nonbounding spin circle is
	assigned an odd line.\index[terminology]{super line}\index[terminology]{super vector
	space}\index[terminology]{bounding spin circle}\index[terminology]{nonbounding spin circle}
\end{example}

When we turn to non-topological invertible field theories, these integrals will use differential (generalized)
cohomology.
\subsection{Non-topological invertible field theories}
\label{non_top_field_theory}
Using reflection-positive invertible TFTs, we saw the torsion subgroup of $[\MTH, \Sigma^{n+1}\IZ]$.
Freed--Hopkins \cite[\S 5.4]{FH21} go further and conjecture that the entire group classifies reflection-positive
invertible field theories that are not necessarily topological. At present, it is not clear how to define these
field theories. But Freed--Hopkins predict what the partition functions of these theories should be, which is a
differential-cohomological lift of the topological story, where we had bordism invariants. We follow
Freed \cite[Lecture 9]{Fre19} and Freed--Hopkins \cite[\S 5.4]{FH21} in this section.
\begin{definition}
A \textit{differential $H_n$-structure} on a smooth manifold $M$ is
\index[terminology]{differential Hn-structure@differential $H_n$-structure}
\begin{enumerate}
	\item a Riemannian metric on $M$,
	\item an $H$-structure in the sense above, i.e.\ a principal $H_n$-bundle $P\to M$ with an isomorphism
	$\theta\colon P\times_{H_n}\Or_n\isomorphism\mathcal B_{\Or}(M)$, and
	\item a connection $\conn$ on $P$ whose induced connection under $\theta$ is the Levi-Civita connection for
	the metric.
\end{enumerate}
\end{definition}
A differential $H_n$-structure on $M$ induces a differential $H_n$-structure on a collar neighborhood of $\partial
M$, so analogously to $\Bord_n^H$, there should be a ``geometric bordism category'' $\Bord_n^{H_n, \nabla}$. Then
one should be able to define field theories as symmetric monoidal functors from $\Bord_n^{H_n, \nabla}$ to
something like a category of topological vector spaces, and define invertibility as above.
\index[terminology]{bordism category!geometric}
\index[terminology]{geometric bordism category|see {bordism category!geometric}}
Following ideas of Atiyah, Kontsevich, and Segal \cite{Seg11}, various geometric versions of bordism categories
have been constructed or sketched by
Cheung \cite{Che07},
Ayala \cite{Aya09},
Hohnhold--Stolz--Teichner \cite[\S 6.2]{HST10},
Hohnhold--Kreck--Stolz--Teichner \cite[\S 5.2]{HKST11},
Stolz--Teichner \cite{ST11},
Tachikawa \cite[\S 1]{Tac13},
Schommer-Pries--Stapleton \cite[\S 7]{SS14},
Kandel \cite{Kan16},
Grady--Sati \cite[\S 5.2]{GS17},
Ulrickson \cite[\S 2.1.2]{Ulr17},
Müller--Szabo \cite[\S 2.1]{MS18},
Grady--Pavlov \cite[\S 4.2]{GP20},
Ludewig--Stoffel \cite[\S 3]{LS20}, and
Kontsevich--Segal \cite{KS21}; Müller--Szabo use their model to study examples of invertible,
non-topological field theories.
\begin{conjecture}[{(Freed--Hopkins \cite[Conjecture 8.37]{FH21})}]
\label{nontopinv}
	There is an isomorphism of abelian groups from $\uppi_0$ of the space reflection-positive, invertible,
	$n$-dimensional field theories to $[\MTH, \Sigma^{n+1}\IZ]$.
\end{conjecture}
Key to this conjecture is formulating a good definition of invertible, non-topological field theory. In the rest of
this section, we assume the conjecture is true, which in particular means finding a definition.

% conjectural description of partition function: a differential H_n-structure means we can integrate classes in
% \check IZ(BH). maybe should phrase parallelly
This conjecture includes a prediction for the value of the partition function of an invertible field theory given
by $\varphi\in\Map(\MTH, \Sigma^{n+1}\IZ)$. An $H$-manifold $M$ gives a point in $\MTH$, i.e.\ a map
$M\colon \Sigma^n\mathbb S\to\MTH$. Composing with $\varphi$ and desuspending, we have a map $\mathbb
S\to\Sigma\IZ$; its homotopy class is an element of $\IZ^1(\pt) = \uppi_{-1}\IZ = 0$, so this construction is not
very interesting. But conjecturally, a differential refinement of this procedure takes a manifold $M$ with a
differential $H_n$-structure and obtains an element $\varphi(M)\in \IZhat^1(\pt)\cong\RR/\ZZ$; then the partition
function of the corresponding invertible field theory is predicted to be $\exp(2\pi i\varphi(M))$. See
Hopkins--Singer \cite[\S 5.1]{HopkinsSinger} for a construction which adopts this perspective; they in particular
construct the differential refinement $\IZhat$ of $\IZ$, by using that $\HZZ\to\IZ$ is a rational equivalence.
Yamashita--Yonekura \cite{YY21, Yam21} take another approach, directly constructing a differential refinement of
$\mathrm{Map}(\MTH, \Sigma^2\IZ)$ and using it to access the partition functions of these conjectured field
theories.
\index[terminology]{Anderson dual of the sphere spectrum!differential refinement}

Often there is a simpler description. Assume $\varphi$ can be identified with the element of the group $\Hom(\Omega_{n+1}^H,
\ZZ)$ given by integrating a (generalized) cohomology class $c$. Then the partition function of the theory
associated to $\varphi$ is the secondary invariant associated to $c$, as defined in \cref{secondary_invariants}.
\index[terminology]{secondary invariant}
\begin{example}[(classical Chern--Simons theory)]
\label{classical_CS}
\index[terminology]{Chern--Simons theory!classical}
The Chern--Simons invariants we discussed above in \cref{cs_invariants} fit together into an invertible,
non-topological field theory which is a differential analogue of \cref{classical_DW}. Fix a compact Lie group and a
\textit{level} $\lambda\in\H^4(\BG;\ZZ)$. Assume $\lambda$ is not torsion. Since $G$ is compact, the Chern--Weil map
is an isomorphism, so as in \cref{DifferentialCharacteristicClasses}, $\lambda$ refines to a class
$\hat \lambda\in \Hhat^4(\BunGnabla ;\ZZ)$.

The level $\lambda$ defines an element of $\Hom(\Omega_4^{\SO}(\BG);\ZZ)$: send an oriented $4$-manifold $X$ with
principal $G$-bundle $P\to M$ to the integer $\int_M \lambda(P)$, where $\lambda(P)$ denotes the pullback of
$\lambda$ along the homotopy class of maps $M\to \BG$ defined by $P$. Again, Stokes' theorem is why this is a
bordism invariant. According to \cref{nontopinv}, this bordism invariant determines (up to isomorphism) an
invertible field theory for $3$-manifolds with a differential $\SO_3\times G$-structure. This field theory is
classical Chern--Simons theory \cite{Fre95, Fre02, Gom01}
\begin{equation}
	\alpha_{(G, \lambda)}\colon\Bord_3^{\SO\times G, \nabla}\longrightarrow \Line_\CC \period
\end{equation}
Let $Y$ be a closed $3$-manifold with a differential $\SO\times G$-structure, which means an orientation, a
Riemannian metric, a principal $G$-bundle $P\to Y$, and a connection $\conn$ for $P$. The data of $(P,\conn)$
gives us a map $Y\to\BunGnabla $, allowing us to pull $\hat \lambda$ back to $Y$, and the orientation allows us to
integrate differential cohomology classes, as in \cref{FiberIntegration}. The
partition function of $\alpha_{(G, \lambda)}$ is $\exp\paren{2\pi i \int_Y \hat \lambda(P, \conn)}$, which is
exactly the exponentiated Chern--Simons invariant of $(P, \conn)$, as we established in~\eqref{CS_secondary_CW}:
\begin{align*}
	\exp2\pi i\int_Y\colon \Hhat^4(Y) &\longrightarrow \Hhat^1(\pt)\to \CC^\times\\
	\hat \lambda(P, \conn) &\longmapsto \exp\paren{2\pi i\CS_\lambda(P, \conn)} \comma
\end{align*}
That is, $\Hhat^1(\pt)\cong\RR/\ZZ$, and exponentiating gets us to $\CC^\times$.

On a closed, oriented surface $\Sigma$ with a Riemannian metric, principal $G$-bundle $P\to\Sigma$, and connection
$\conn$, $\alpha_{(G, \lambda)}$ again assigns the pushforward of $\hat \lambda(P, \conn)$, but this time the
pushforward map has signature
\begin{equation}
\int_\Sigma\colon \Hhat^4(\Sigma)\longrightarrow \Hhat^2(\pt)\cong\Line_\CC \comma
\end{equation}
which sends $\hat \lambda(P, \conn)$ to the Chern--Simons line constructed in, e.g., \cite[\S 4]{Fre95}. This
story continues in extended TFT, assigning higher-categorical objects to lower-dimensional manifolds, such as
in \cite{Gom01a}.

See also Fiorenza--Sati--Schreiber \cite{FSS15a}, Yamashita--Yonekura \cite[Example 4.56]{YY21}, and
Yamashita \cite[\S 3.4.2]{Yam21} for additional constructions of classical Chern--Simons theory as an invertible
field theory, and Freed--Neitzke \cite{FN20} for an application to special functions.
\end{example}

\begin{remark}[(quantizing Chern--Simons theory)]
\label{quantum_CS}
One of the interesting things you can do with the classical Chern--Simons theory is to quantize it. This amounts to
summing $\alpha_{(G, \lambda)}$ over the space of all principal $G$-bundles with connection on a given closed,
oriented $3$-manifold. This procedure, known as taking the path integral, is still only heuristically
defined,\footnote{When $G$ is finite, Freed--Quinn \cite{FQ93} define a path integral of \emph{topological} field
theories whose fields include a principal $G$-bundle. Applied to classical Dijkgraaf--Witten theory from
\cref{classical_DW}, the resulting TFT, called \textit{(quantum) Dijkgraaf--Witten theory}, is a commonly studied
model organism in topological field theory.} but enough is known about it in the physics literature that we can ask
mathematical questions about the quantized theory. In physics, this quantum Chern--Simons theory was first studied
by Schwarz \cite{Sch77} and Witten \cite{Wit89}.%
\index[terminology]{Chern--Simons theory}%
\index[terminology]{Dijkgraaf--Witten theory}%
\index[terminology]{path integral}

Something strange happens in this quantization procedure, though: Witten (\textit{ibid.}) gives a physical argument
that quantum Chern--Simons theory is in fact a topological field theory! Therefore it should be possible to
formalize it mathematically as a symmetric monoidal functor 
\begin{equation}
	Z_{G,k}\colon\Bord_3^{\SO}\longrightarrow \fC \comma
\end{equation}
where $\fC$ is some symmetric monoidal $(\infty, 3)$-category. It is not known how to do this in
general,\footnote{There are a few different perspectives on what $Z_{G,k}(\mathrm{pt}_+)$ should be. For $G$
finite, the answer is known by work of Freed--Hopkins--Lurie--Teleman \cite[\S 4.2]{FHLT10} and Wray \cite[\S
9]{Wra10}; for $G$ a torus, the answer is due to Freed--Hopkins--Lurie--Teleman (\textit{ibid}.). For general $G$,
two different approaches are provided by Freed--Teleman (see \cite{432876}) and Henriques \cite{Hen17b, Hen17a}.
See also \cite{FT20}.}
but it is known how to extend it to a theory of $1$-, $2$-, and $3$-manifolds, valued in the $2$-category of
$\mathbb C$-linear categories, by work of Reshetikhin--Turaev \cite{RT90, RT91}, Walker \cite{Wal91},
Bakalov--Kirillov \cite{BK01}, Kerler--Lyubashenko \cite{KL01}, and
Bartlett--Douglas--Schommer-Pries--Vicary \cite{BDSV15}.\footnote{These constructions require some additional
structure on our manifolds, such as a choice of trivialization of the first Pontryagin class. As theories of merely
oriented manifolds, Chern--Simons theories are \emph{anomalous}. See \cites[\S 9.3]{FHLT10}[]{432876} for more
information.} Much more can be said about this TFT and its connections to various parts of
geometry, topology, representation theory, and physics; see Freed \cite{Fre09} for a
general survey on Chern--Simons theory and the references therein for more information.
\end{remark}

\begin{example}[(classical Wess--Zumino--Witten theory)]\label{ex-WZW}
	This example is related to the previous example, but with a slightly different flavor. Let $G$ be a compact Lie
	group and $\hat h\in \Hhat^3(G;\ZZ)$. If $h \colonequals \cc(\hat{h})$ (the image of $\hat{h}$ under the characteristic class map of \Cref{cons:charclassmap}),
	then $h$ defines a bordism invariant of oriented $3$-manifolds $M$ with a map $\psi\colon M\to G$:
	\begin{align*}
		\Omega_3^{\SO}(G) &\longrightarrow \ZZ\\
		(M,\psi) &\longmapsto\int_M \psi^*(h) \period
	\end{align*}
	\Cref{nontopinv} therefore says there is a two-dimensional invertible field theory $\beta_{G, h}$ whose
	partition function is the secondary invariant associated to $\hat h$. This theory is called \textit{classical
	Wess--Zumino--Witten (WZW) theory}; it was originally studied by Witten \cite{Wit83}, following
	Wess--Zumino \cite{WZ71}. See Freed \cite[Appendix A]{Fre95} for a discussion of the classical theory
	specifically.\footnote{There are considerably more general objects studied in quantum physics under the name
	``Wess--Zumino--Witten theory'' or ``Wess--Zumino--Witten term.'' See \cites[\S 6]{DF99}{Fre08}[\S
	5.6]{Urs}{FSS15b}{LOT20}{Yon20} for some examples taking an algebro-topological viewpoint.}%
	\index[terminology]{secondary invariant}%
	\index[terminology]{Wess--Zumino--Witten model}%
	\index[terminology]{Wess--Zumino--Witten model!classical}

	As part of a trend you may have noticed by now, the original description of the classical WZW partition function
	$\int_M \psi^*(\hat h)$ was not phrased in this way; the connection with differential cohomology is due to
	Gawędzki \cite{Gaw88}. For a moment assume that $G$ is connected, simple, and simply connected, so that
	$\H^3(G;\ZZ)\cong\ZZ$. 
	Let $\MC\in\Omega^1(G; \g)$ denote the Maurer--Cartan form (see \Cref{def:Maurer-Cartan_form}).,
	As mentioned in \cref{transgression_detail}, the transgression map
	\begin{equation*}
		\tau^{-1}\colon \H^3(G;\ZZ)\to \H^4(\BG;\ZZ)
	\end{equation*}
	is an isomorphism; since $G$ is compact, the Chern--Weil machine associates to $\tau^{-1}(h)$ (or rather, its image in
	$\RR$-valued cohomology) a degree-two invariant polynomial $f$. In this case, the Wess--Zumino--Witten action is%
	\index[terminology]{transgression}
	\begin{equation}
		\beta_{G,h}(M, \psi) = \int_M -\frac 16\psi^*(f(\MC\wedge [\MC, \MC])).
	\end{equation}
	The differential refinement of $\tau\colon \H^4(\BG;\ZZ)\to\H^3(G;\ZZ)$ constructed by
	Carey--Johnson--Mur\-ray--Stevenson--Wang \cite[\S 3]{CJMSW05} and Schreiber \cite[1.4.1.2]{Urs} can be thought of as
	starting with a classical Chern--Simons theory and obtaining a classical Wess--Zumino--Witten theory in one
	dimension lower.
\end{example}

\begin{remark}[(quantizing the Wess--Zumino--Witten model)]
	\label{quantum_WZW}
	Just as in \cref{quantum_CS}, it is possible to quantize the classical WZW model, at least at a physical level of
	rigor: one sums over the space of maps to $G$. The result is called the quantum Wess--Zumino--Witten model, or just
	the Wess--Zumino--Witten or WZW model.\index[terminology]{Wess--Zumino--Witten model!quantum} This theory is a
	conformal field theory,~\index[terminology]{conformal field theory} meaning its value on a manifold depends only
	on the conformal class of the Riemannian metric. Some of what we do in the next two chapters, involving the
	representation theory of loop groups, is related to the WZW model.

	Given a level $h\in\Hhat^4(\BunGnabla ;\ZZ)$, there is a (quantum) Chern--Simons theory and a quantum WZW model
	(obtained by transgressing $h$ to $\Hhat^3(G;\ZZ)$), and the two are related: the WZW model is a boundary theory
	for the Chern--Simons theory. There are different ways of formulating this precisely: one uses \textit{relative
	field theory}\index[terminology]{relative field theory} \cite{FT12}. In this formalism, the bulk theory $\alpha$
	is a symmetric monoidal functor out of a bordism category, and its boundary theory $Z$ is a natural transformation
	from (a truncation of) $\alpha$ to the trivial field theory. Among other things, this implies that the partition
	function of $Z$ on an $(n-1)$-manifold $M$ is not a number, but an element of the state space $\alpha(M)$; when
	$\alpha$ is Chern--Simons theory and $Z$ is the WZW model, this fact was first noticed by Witten \cite{Wit89}.
	See Gwilliam--Rabinovich--Williams \cite{GRW20} for another approach to this bulk-boundary correspondence, in the
	language of factorization algebras.
\end{remark}

\begin{example}[(exponentiated $\eta$-invariants)]
We now give a differential analogue of \cref{arf_TFT}: in that example, we used the Atiyah--Bott--Shapiro
pushforward \cite{ABS64} in $\KO$-theory to produce a torsion bordism invariant, hence an invertible topological
field theory. Here we will use the same pushforward to produce a nontorsion bordism invariant, hence an invertible,
non-topological field theory. This theory is discussed by Freed \cite[Example 9.24]{Fre19}.
\index[terminology]{Atiyah--Bott--Shapiro map}

The bordism invariant in question is the \emph{$\Ahat$-genus} $\Ahat\colon
\Omega_4^{\Spin}\to\ZZ$,\footnote{$\Ahat$ is pronounced ``$A$-hat'' or ``$A$-roof.'' This gives rise to the
following joke: A man walks into a bar with a dog and says to the bartender, ``This is a talking dog. I'll bet you
a drink he can answer a question.''

The bartender says, ``Sure. Ok dog, what's your favorite spin bordism invariant?''

``Arf!''

``\dots''

``Clifford, how about a different one?''

``A-roof!''

(they get thrown out)

The dog looks at the man and says, ``Ok fine, next time I'll say `index of the Dirac operator.' ''}
\index[terminology]{A-genus@$\Ahat$-genus}
\index[notation]{Ahat@$\Ahat$}
which, like the Arf invariant, is a pushforward in $\KO$-theory: for a closed spin $4$-manifold $X$, we have
\begin{align*}
	\int_X^{\KO}\colon \KO^0(X) &\longrightarrow \KO^{-4}(\pt)\cong\ZZ\\
	1 &\longmapsto \Ahat(X) \period
\end{align*}
This is nonvanishing on the K3 surface,\index[terminology]{K3 surface} hence nontorsion. By Freed--Hopkins'
conjecture, this bordism invariant corresponds to some invertible, non-topological field theory on $3$-dimensional
differential spin manifolds (i.e.\ $3$-manifolds with a spin structure and a Riemannian metric):
\begin{equation*}
	\alpha'\colon\Bord_3^{\Spin, \nabla}\longrightarrow \sLine_\CC \period
\end{equation*}
And analogously to the Arf theory, we can describe the value of $\alpha'$ on closed $2$- and $3$-manifolds with
differential spin structure using the pushforward in differential $\KO$-theory. Grady--Sati \cite[\S 4.3]{GS21}
construct this pushforward for a closed spin manifold; using this, the partition function of $\alpha'$ on a closed
spin Riemannian $3$-manifold $Y$ is
\index[terminology]{differential KO-theory@differential $\KO$-theory}
\begin{align*}
	\exp2\pi i\int_Y^{\KOhat}\colon \KOhat^0(Y) &\longrightarrow \KOhat^{-3}(\pt)\to
	\CC^\times\\
	1 &\longmapsto \alpha'(Y) \comma
\end{align*}
where as usual $\KOhat^{-3}(\pt)\cong\RR/\ZZ$, and we exponentiate to obtain the partition function in
$\CC^\times$. The isomorphism type of the state space assigned to a closed spin Riemannian $2$-manifold $\Sigma$ is
in a similar way the image of $1$ under the pushforward
$\KOhat^0(\Sigma)\to\KOhat^{-2}(\pt)\cong\ZZ/2$, corresponding to the two isomorphism classes of
complex super lines, $\CC$ and $\Pi\CC$.

Like in \cref{classical_CS}, the partition function of $Y$ also has a more geometric description. A differential
spin structure is the data needed to define the Dirac operator on the spinor bundle of $Y$, and index-theoretic
methods allow one to extract an \textit{exponentiated $\eta$-invariant} from this Dirac operator, as constructed by
Atiyah--Patodi--Singer \cite{APS1, APS2, APS3}. The Dai--Freed theorem \cite{DF94} proves this exponentiated
$\eta$-invariant satisfies a gluing law which can be interpreted as implying that $\alpha'$ is symmetric monoidal.
\index[terminology]{Dai--Freed theorem}
\index[terminology]{Dirac operator}
\index[terminology]{eta-invariant@$\eta$-invariant}
\end{example}
For more examples of invertible, non-topological field theories and their relationship to differential cohomology,
see Monnier \cites[\S 4]{Mon15}[\S 5]{Mon17}{Mon18}, Monnier--Moore \cite{MM19},
Córdova--Freed--Lam--Seiberg \cite[\S\S 6.2 \& 7]{CFLS20a},
Yamashita--Yonekura \cite[\S\S 4.2, 6]{YY21}, and Yamashita \cite[\S 3.4]{Yam21}.

\newpage
%!TEX root = ../diffcoh.tex

\section{Loop groups and intertwining of positive-energy representations}
\textit{by Sanath Devalapurkar}
\label{loop_groups}

We will give an introduction to the representation theory of loop groups of
compact Lie groups: we will discuss what positive energy representations
are, why they exist, how to construct them (via a Schur--Weyl style
construction and a Borel--Weil style construction), and how to show that they
don't depend on choices. Motivation will come from both mathematics and
quantum mechanics.

The theory of positive-energy representations of loop groups is modeled
on the representation theory of compact Lie groups. Some parts of the talk will make more sense if you are familiar
with the compact Lie group story, but this is not a requirement: in this section, we try to emphasize the ``big
picture'' over details, and we hope that this choice makes it readable for you. Likewise, we will not assume any
familiarity with loop groups or infinite-dimensional topology, nor will we dig into those details.

In \cref{ssec_loop_overview}, we state the main theorem (\cref{main-thm}) and discuss some motivation for caring
about representations of loop groups. In \cref{rep_loop}, we begin thinking about projective representations of
loop groups and the corresponding central extensions. In \cref{PS_proof_sketch}, we provide an extended proof
sketch of \cref{main-thm}, and discuss some connections to physics. Finally, in \cref{PS_diffcoh}, we discuss how
this relates to differential cohomology. There are two ways to lift the construction of central extensions of loop
groups to differential cohomology; one follows the Chern--Weil story we've used several times already in this part,
and the other more closely resembles the story we told about off-diagonal Deligne cohomology and the Virasoro
algebra in \cref{VirasoroAlgebra}.

%-------------------------------------------------------------------%
%-------------------------------------------------------------------%
%  Introduction													 %
%-------------------------------------------------------------------%
%-------------------------------------------------------------------%

\subsection{Overview}
\label{ssec_loop_overview}

The objective of this chapter is to explain the following theorem of Pressley--Segal \cite[Theorem 13.4.2]{loop}:
\begin{theorem}\label{main-thm}
	Let $G$ be a simply connected compact Lie group. Then any positive energy representation $E$ of the loop group
	$\LG$ admits a projective intertwining action of $\Diffplus(\Circ)$.
\end{theorem}
If this means nothing to you, that's okay: the goal of this talk is to explain all the components of this theorem
(\cref{rep_loop}) and sketch a proof (\cref{PS_proof_sketch}). Then, in \cref{PS_diffcoh}, we discuss how the
representation theory of loop groups is related to differential cohomology.

Here's a rough sketch of what \Cref{main-thm} is about. The
representation theory of a semisimple compact Lie group $G$ is very well-behaved: the Peter--Weyl
theorem \cite{PW27} \index[terminology]{Peter--Weyl theorem} allows one to provide any finite-dimensional
$G$-representation with a $G$-invariant Hermitian inner product, and this inner product decomposes the
representation into a direct sum of irreducibles. Moreover, the irreducibles are in bijection with dominant
weights,\index[terminology]{dominant weight} where by the Borel--Weil theorem (see \cite{Ser54}), the
representation associated to a dominant weight is given  as the global sections of a line bundle associated to a
homogeneous space of $G$ (a particular flag variety). \index[terminology]{Borel--Weil
theorem}\index[terminology]{homogeneous space}\index[terminology]{flag variety}

Most representations of loop groups will not satisfy analogues of this property,
so we'd like to hone down on the ones which do. These are the ``positive energy
representations''; these essentially satisfy properties necessary to be able to
write down highest/lowest weight vectors. \Cref{main-thm} then states
that positive energy representations are preserved under reparametrizations of
the circle (which give automorphisms of the loop group $\LG$). One can therefore
think of \Cref{main-thm} as a consistency result.

Before proceeding, I'd like to give some motivation for caring about the
representation theory of loop groups.
\begin{enumerate}
	\item One motivation comes from the connection between representation theory and homotopy theory. The
	Atiyah--Segal completion theorem \cites[Theorem 7.2]{Ati61}[\S 4.8]{AH61}[Theorem
	2.1]{AS69}\index[terminology]{Atiyah--Segal completion theorem} relates representations of a compact Lie group
	$G$ to $G$-equivariant \Ktheory, and likewise the representation theory of the loop group $\LG$ is related
	to (twisted) $G$-equivariant elliptic cohomology. This has been explored in \cite{Bry90, Dev96, Liu96, And00,
	And03, Gro07, Lur09a, Gan14, Lau16, Kit19, Rez20, BET21}.\index[terminology]{elliptic cohomology}
	\item Another motivation comes from the hope that geometry on the free loop space $\LM$ of a manifold $M$ is
	supposed to correspond to correspond to ``higher-dimensional geometry'' over $M$.\index[terminology]{loop
	space!in terms of higher-dimensional geometry} For instance, if $M$ has a Riemannian metric, one can
	think of the scalar curvature of $\LM$ at a loop as the integral of the Ricci curvature of $g$ over the loop.
	Similarly, spin structures on $M$ are closely related to orientations on $\LM$ \cites{Wit85}[\S
	3]{Ati85}{Wit88}[\S 2]{McL92}[Theorem 9]{ST05}[Corollary E, \S 1.2]{Wal16}, and string structures on $M$ are
	closely related to spin structures on $\LM$ \cites{Kil87}[Theorem 6.9]{NW13}.\footnote{There are a number of
	other works providing additional proofs of this fact or pointing out subtleties in the definitions,
	including \cites{PW88}{CP89}[\S 3]{McL92}{KY98}{ST05}{KM13a}{Wal15}{Cap16}{Wal16a}{Kri20}.}
	\index[terminology]{string structure}
\end{enumerate}
In light of this hope, it is rather pacifying to
have a strong analogy between representation theory of compact Lie groups and of loop groups. In fact, all of
these motivations are related by a story that still seems to be mysterious at the moment.

There's also motivation from physics for studying the representation theory of
loop groups. The wavefunction of a free particle on the circle $\Circ$ must be an
$\Lup^2$-function on $\Circ$ (because the probability of finding the particle
somewhere on the circle is $1$). There is an action of the loop group $\LU_1$
on $\Lup^2(\Circ;\CC)$ given by pointwise multiplication (a pair $\gamma \colon \Circ\to
\Uup_1$ and $f\in \Lup^2(\Circ;\CC)$ is sent to the $\Lup^2$-function $f_\gamma(z) =
\gamma(z) f(z)$). In particular, $\LU_1$ gives a lot of automorphisms of the
Hilbert space $\Lup^2(\Circ;\CC)$; this is relevant to quantum mechanics, where
observables are (Hermitian) operators on the Hilbert space of states. Having a
particularly (mathematically) natural source of symmetries is useful. In
\cite{segal-survey}, Segal in fact says: ``In fact it is not much of an
exaggeration to say that the mathematics of two-dimensional quantum field theory
is almost the same thing as the representation theory of loop groups''.

%-------------------------------------------------------------------%
%-------------------------------------------------------------------%
%  Representations of loop groups								   %
%-------------------------------------------------------------------%
%-------------------------------------------------------------------%

\subsection{Representations of loop groups}
\label{rep_loop}

\begin{definition}
	Let $G$ be a compact connected Lie group. The loop group $ \LG \colonequals \Cinf(\Circ, G)$ is the group of smooth unbased loops in $G$.
\end{definition}
If $G$ is positive-dimensional, $\LG$ is not finite-dimensional. A fair amount of the theory of finite-dimensional
manifolds generalizes to infinite-dimensional spaces locally modeled by nice classes of topological vector spaces,
and in this sense $\LG$ is an infinite-dimensional Lie group, in fact quite a nice one. Reading this chapter does
not require any additional familiarity with infinite-dimensional topology, but if you're interested, you can learn
more in \cites{Ham82}{Mil84}[\S 3.1]{loop}

There will be a lot of circles floating around, and so we will distinguish these
by subscripts. Some of these will be denoted by $\TT$, for ``torus''.

\begin{remark}[(classification of compact Lie groups)]
We quickly review the classification of compact Lie groups. This may clarify the generality in which some of the
results in this section hold.
\begin{itemize}
	\item Let $G$ be a compact Lie group and $G_0\subset G$ denote the connected component containing the identity.
	Then there is a short exact sequence
	\begin{equation*}
		1\to G_0\to G\to \uppi_0(G)\to 1 \period
	\end{equation*}

	\item Let $G$ be a compact, connected Lie group. Then there is a short exact sequence
	\begin{equation*}
		1\to F \to\tilde G\to G\to 1 \comma
	\end{equation*}
	where $F$ is finite and $\tilde G$ is a product of a torus $\TT^n$ and a simply connected group.

	\item Let $G$ be a compact, connected, simply connected Lie group. Then $G$ is a product of \emph{simple}
	simply connected Lie groups.\index[terminology]{simple Lie group}

	\item Let $G$ be a compact, simply connected, simple Lie group. Then $G$ is isomorphic to one of
	$\SU_n$, $\mathrm{Spin}_n$, $\mathrm{Sp}_n$, $\mathrm G_2$, $\mathrm F_4$, $\mathrm E_6$,
	$\mathrm E_7$, or $\mathrm
	E_8$.\index[notation]{SUn@$\SU_n$}\index[notation]{Spinn@$\mathrm{Spin}_n$}
	\index[notation]{Spn@$\mathrm{Sp}_n$}\index[notation]{G2@$\mathrm G_2$} \index[notation]{F4@$\mathrm F_4$}
	\index[notation]{E6@$\mathrm E_6$}\index[notation]{E7@$\mathrm E_7$}\index[notation]{E8@$\mathrm E_8$}
\end{itemize}
Most of the results in this section require $G$ to be connected and simply connected; a few will also require $G$
to be simple. In particular, when $G$ is simple, $\H^4(\BG;\ZZ)\cong\ZZ$.\footnote{This isomorphism can be made
canonical by specifying that under the Chern--Weil map, the Killing form $B\colon\g\times\g\to\RR$ defines a
positive element of $\HdR^4(\BG)\cong\RR$.\index[terminology]{Killing form}}
\end{remark}

\begin{remark}
	The loop group $\LG$ is an infinite-dimensional Lie group, and it has an
	action of $\Circ$ by rotation. We will denote this ``rotation'' circle by
	$\TTrot$.  This action will turn out to be very useful shortly.\index[notation]{Trot@$\TTrot$}
\end{remark}

The action of $\TTrot$ allows one to consider the semidirect product $\LG \rtimes \TTrot$. The following proposition is then an exercise in manipulating
symbols:
\begin{proposition}
	An action of $\LG \rtimes \TTrot$ on a vector space $V$ is the same data as
	an action $R$ of $\TTrot$ on $V$ and an action $U$ of $\LG$ on $V$
	satisfying
	\begin{equation*}
		R_\theta U_\gamma R_\theta^{-1} = U_{R_\theta \gamma} \period
	\end{equation*}
\end{proposition}
Most interesting representations $U$ of $\LG$ on a vector space $V$ are not, strictly speaking, representations:
instead of $U_\gamma U_{\gamma'} = U_{\gamma\gamma'}$, they satisfy the weaker condition that
\begin{equation}
	U_\gamma U_{\gamma'} = c(\gamma, \gamma') U_{\gamma \gamma'} \comma
\end{equation}
where $c(\gamma, \gamma')\in \CC^\times$. This is precisely:
\begin{definition}
	A \emph{projective representation} of $\LG$ on a Hilbert space $V$ is a
	continuous homomorphism $\LG \to \PU(V)$.\index[terminology]{projective representation}
\end{definition}

\begin{remark}
	Why Hilbert spaces? From a mathematical perspective, this is because Hilbert
	spaces are well-behaved infinite-dimensional vector spaces. From a physical
	perspective, this is because Hilbert spaces are spaces of states. In fact,
	this also explains why most interesting representations are projective: the
	state of a quantum system is not a vector in the Hilbert space, but rather a
	vector in the projectivization of the Hilbert space. This corresponds to the
	statement that shifting the wavefunction by a phase does not affect physical
	observations.
\end{remark}

Assume $V$ is an infinite-dimensional, separable Hilbert space. Then $\PU(V)$ is a $\Kup(\ZZ, 2)$, so
projective representations determine cohomology classes in $\H^2(\LG;\ZZ)$.
\begin{lem}
\label{csc_loop}
When $G$ is compact and simply connected, $\H^2(\LG;\ZZ)\cong\H^3(G;\ZZ)$.
\end{lem}
\begin{proof}
Since $G$ is simply connected, $\uppi_1(G) = 0$, and $\uppi_2$ vanishes for any Lie group. Therefore the Hurewicz
theorem identifies $\uppi_3(G)$ and $\H_3(G;\ZZ)$. Let $\Omega G$ denote the based loop space of
$G$\index[notation]{OmegaG@$\Omega G$}, i.e.\ the subspace of $\LG$ consisting of loops beginning and ending at the
identity. Essentially by definition, there is an isomorphism $\uppi_k(G)\to\uppi_{k-1}(\Omega G)$ for $k > 1$, so we
learn $\uppi_1(\Omega G) = 0$ and $\uppi_2(\Omega G)\cong\uppi_3(G)$.

To get to $\LG$, we use that as topological spaces, $\LG\cong G\times\Omega G$ \cite[\S 4.4]{loop}. Thus
$\uppi_1(\LG) = 0$ and $\uppi_2(\LG) \cong \uppi_3(G)$, and the Hurewicz and universal coefficient theorems allow us to
conclude.\index[terminology]{Hurewicz theorem}\index[terminology]{universal coefficient theorem}
\end{proof}
Another way to construct this isomorphism is as follows: there is an evaluation map
\begin{align*}
	\ev\colon \Circ \times \LG &\to G \\ 
	(x, \ell) &\mapsto \ell(x) \semicolon
\end{align*}
then the isomorphism in~\cref{csc_loop} is: pull back by $\ev$, then integrate in the $\Circ$ direction.

It turns out that when $G$ is compact and simply connected, every class in $\H^2(\LG;\ZZ)$ arises from a projective
representation as above \cite[Theorem 4.4.1]{loop}. There is a central extension\footnote{This central extension is
also a fiber bundle, and by Kuiper's theorem \cite{Kui65}, the total space $\Uup(V)$ is contractible (see
also \cites[Lemme 3]{DD63}[Proposition A2.1]{AS04}).\index[terminology]{Kuiper's theorem} This fiber bundle is
homotopy equivalent to two other interesting fiber bundles: the universal principal $\Uup_1$-bundle
$\Uup_1\to\mathrm{EU}_1\to\mathrm{BU}_1$, and the loop space-path space bundle $\Omega \Kup(\ZZ, 2)\to P
\Kup(\ZZ, 2)\to \Kup(\ZZ, 2)$.}
\begin{equation}
\label{PUcentral_ext}
	1\to \TTce\to \Uup(V)\to \PU(V)\to 1 \comma
\end{equation}
and so any projective representation $\rho$ of $\LG$ determines a central
extension by pulling \eqref{PUcentral_ext} back:
\begin{equation}
	1\to \TTce\to \LGtilde_\rho \to \LG \to 1 \period
\end{equation}
Conversely, any central extension of $\LG$ gives rise to a projective
representation of $\LG$. In particular:

%\begin{remark}\label{central}
%	Isomorphism classes of central extensions of $\LG$ are in bijection with
%	elements of $\H^2(\LG;\ZZ)$. If $G$ is simple and simply connected, then
%	$\H^2(\LG; \ZZ) \cong \H^3(G; \ZZ) \cong \ZZ$.
%\end{remark}

\begin{definition}
	Let $G$ be a simple and simply connected compact Lie group. The
	\emph{universal central extension} $\LGtilde$ of $\LG$ is the central
	extension corresponding to the generator of $\H^2(\LG; \ZZ) \cong \ZZ$.\index[terminology]{universal central
	extension}
\end{definition}
We first met universal central extensions in a different context, in \cref{diffcoh_virasoro}.

The following result is key.
\begin{theorem}[{\cite[Theorem 4.4.1]{loop}}]
Let $G$ be simply connected. Then there is a unique action of $\Diffplus(\TTrot)$ on $\LGtilde$ which covers
	the action on $\LG$. Moreover, $\LGtilde$ deserves to be called
	``universal'', because there is a unique map of extensions from $\LGtilde$ to
	any other central extension of $\LG$.
\end{theorem}

\begin{remark}
	As a consequence, the action of $\TTrot$ on $\LG$ lifts canonically to
	$\LGtilde$. Every projective unitary representation of $\LG$ with an
	intertwining action of $\TTrot$ is equivalently a unitary representation
	of $\LGtilde\rtimes \TTrot$. For the remainder of this talk, we will assume $G$ is simply connected and
	abusively say write ``representation of $\LG$'' to mean a representation of $\LGtilde \rtimes \TTrot$.
\end{remark}

\begin{notation}
	It is a little inconvenient to constantly keep writing $\LGtilde\rtimes
	\TTrot$, so we will henceforth denote it by $\LGtildeplus$. The subgroup
	$\TTrot$ of $\LGtildeplus$ is also known as the ``energy circle'' (for
	reasons to be explained below).
	\index[terminology]{energy circle}
	\index[notation]{LGplus@$\LGtildeplus$}
\end{notation}

One of the nice properties of tori is that their representations take on a
particularly simple form, thanks to the magic of Fourier series. The action of
$\Circ$ on a finite-dimensional vector space is the same data as a $\ZZ$-grading.
The case of topological vector spaces is slightly more subtle: if $\Circ$ acts on
a topological vector space $V$, then one can consider the closed ``weight''
subspace $V_n$ of $V$ where the action of $\Circ$ is by the
character\footnote{Some conventions are different: the action might be by
$z\mapsto z^n$. We're following \cite{loop}.} $z\mapsto z^{-n}$. Then the direct
sum $\bigoplus_{n\in \ZZ} V_n$ is a dense subspace of $V$; it is known as the
subspace of \emph{finite energy} vectors in $V$. This is simply the usual weight
decomposition adapted to the topological setting.

\begin{definition}
	The action of $\Circ$ on a topological vector space $V$ is said to satisfy the
	\emph{positive energy condition} if the weight subspace $V_n = 0$ for
	$n<0$. Equivalently, the action of $\Circ$ is represented by $e^{-iA\theta}$,
	where $A$ is an operator with positive spectrum.
\end{definition}

\begin{remark}
	The motivation for this definition comes from quantum mechanics: the
	wavefunction of a free particle on a circle is $e^{inx}$ (up to
	normalization), and requiring that the energy (which is essentially the
	weight $n$) to be positive is mandated by physics.
\end{remark}

\begin{definition}
\label{pos_en}
	A representation of $\LG$ (which, recall, means a representation of
	$\LGtildeplus$) is said to satisfy the \emph{positive energy condition} if it
	satisfies the positive energy condition when viewed as a representation of
	the energy/central circle $\TTrot$.
	\index[terminology]{positive energy condition}
\end{definition}

\begin{remark}
	It doesn't make sense for a representation of $\LG$ to be positive energy if
	you take ``representation of $\LG$'' to mean a literal representation of
	$\LG$; one needs to interpret that phrase as meaning a representation of
	$\LGtildeplus$.
\end{remark}

We can now see the utility of \Cref{main-thm}: the positive energy
condition involves the canonical parametrization of the circle, and to ensure that our definition would agree with
that of an alien civilization's, we should ensure that the pullback $\fupperstar V$ of any positive energy
representation $V$ of $\LG$ along an orientation-preserving diffeomorphism $f\in \Diffplus(\TTrot)$ is another
positive energy representation. That is precisely the content of \Cref{main-thm}.

At the beginning of this chapter, we said that positive energy representations of loop groups satisfy analogues of
many properties of representations of compact Lie groups. To make that statement precise, we need to introduce some
definitions that impose sanity conditions on the representations we want to study.

\begin{definition}
	Let $V$ be a representation of a topological group $G$ (possibly
	infinite-di\-men\-sion\-al). Then $V$ is said to be:
	\begin{itemize}
	\item \emph{irreducible} if it has no closed $G$-invariant subspace;\index[terminology]{irreducible
	representation}
	\item \emph{smooth} if the following condition is satisfied: let
		$V_\mathrm{sm}$ denote the subspace of vectors $v\in V$ such that
		the orbit map $G\to V$ sending $g$ to $gv$ is continuous; then
		$V_\mathrm{sm}$ is dense in $V$.\index[terminology]{smooth representation}
	\end{itemize}
	Two $G$-representations $V$ and $W$ are \emph{essentially equivalent} if
	there is a continuous $G$-equivariant map $V\to W$ which is injective and
	has dense image.\index[terminology]{essential equivalence}
\end{definition}

\begin{warning}
	Essential equivalence is \emph{not} an equivalence relation!
\end{warning}

The representation theory of compact Lie groups is really nice: every finite-dimensional complex representation of
a compact Lie group $G$ is semisimple (i.e.\ it is a direct sum of irreducible
representations),\index[terminology]{semisimple!representation} and unitary, and extends to a representation of the
complexification $G_\CC$ of $G$.\footnote{A complexification of a real Lie group $G$ is a complex Lie group,
generally noncompact, whose Lie algebra is isomorphic to $\g\otimes\CC$. When $G$ is compact, $G_\CC$ is unique up
to isomorphism.}\index[terminology]{complexification} These properties have analogues for positive energy
representations of loop groups.

\begin{theorem}[{\cite[Theorem 9.3.1]{loop}}]
\label{like_cpt_Lie}
	Let $V$ be a smooth positive energy representation of $\LG$. Then up to
	essential equivalence:
	\begin{itemize}
	\item $V$ is completely reducible into a discrete direct sum of
		irreducible representations,
	\item $V$ is unitary,
	\item $V$ extends to a holomorphic projective representation of
		$\mathrm L(G_\CC)$, and
	\item $V$ admits a projective intertwining action of $\Diffplus(\Circ)$,
		where this $\Circ$ is the energy/rotation circle. (This is \Cref{main-thm}.)
	\end{itemize}
\end{theorem}

The proof of this result takes up the bulk of the second part of Pressley--Segal.

\begin{remark}\label{pos-energy}
	The group $G$ includes into $\LG$ as the subgroup of constant loops. Let $G$
	be simple and simply connected. If $T$ is a maximal torus of $G$, then one
	has
	\begin{equation*}
		\TTrot \times T\times \TTce \subseteq \LGtildeplus \period
	\end{equation*}
	Consequently, if
	$V$ is a representation of $\LGtildeplus$, then $V$ can be decomposed (up to
	essential equivalence) as a $\TTrot \times T\times
	\TTce$-representation:
	\begin{equation}
		V = \bigoplus_{(n,\lambda, h) \in \TTrot^\vee \times T^\vee\times
	\TTce^\vee} V_{(n,\lambda,h)}\period
	\end{equation}
	Here, $n$ is the energy of $V$; $\lambda$ is a weight of $V$ (regarded
	as a representation of $T$); and $h$ is a character of $\TTce$. The notation $(-)^\vee \colonequals
	\Hom(-, \CC^\times)$ denotes the character dual:\index[terminology]{character dual} because
	$\TTrot\times T\times\TTce$ is a compact abelian group, its unitary representations are direct sums of
	one-dimensional representations. Therefore as a $\TTrot\times T\times\TTce$-representation, $V$ splits as a
	direct sum of one-dimensional representations, which are indexed by the character dual
	\begin{equation*}
		(\TTrot \times T \times \TTce)^\vee = \TTrot^\vee \times T^\vee\times \TTce^\vee \period
	\end{equation*}
	
	If $V$ is
	irreducible, then $\TTce$ must act by scalars by Schur's lemma, and so
	only one value of $h$ can occur; this is called the \emph{level} of $V$. It\index[terminology]{level}
	turns out that if $V$ is a smooth positive energy representation, then each
	weight space $V_{n,\lambda,h}$ is finite-dimensional. In fact, a
	representation of $\LG$ of level $h$ is the same as a representation of
	$\LGtilde_h\rtimes \TTrot$, where $\LGtilde_h$ is the central extension of
	$\LG$ corresponding to $h\in \ZZ \cong \H^2(\LG;\ZZ)$.
\end{remark}

\begin{remark}
	By \Cref{pos-energy}, an irreducible positive energy representation
	$V$ of $\LG$ is uniquely determined by the level $h$ and its lowest energy
	subspace $V_0$: the representation $V$ is generated as a
	$\LGtildeplus$-representation by $V_0$.
\end{remark}

\begin{remark}
	Since $G$ is simply connected, there are transgression isomorphisms
	\begin{equation*}
		\H^4(\BG;\ZZ)\to\H^3(G;\ZZ)\to\H^2(\LG;\ZZ) \comma
	\end{equation*}
	meaning we can understand the level as (up to homotopy)
	a map $\BG\to \Kup(\ZZ,4)$. This $\Kup(\ZZ,4)$ is closely tied to the
	twisting $\Kup(\ZZ,4)\to \BGL_1(\tmf)$ of $\tmf$ constructed in \cite[Theorem 1.1]{ABG10}: see \cite{And00,
	Gro07, BET21}.
\end{remark}

As a side note, we observe the following:

\begin{proposition}
	Let $V$ be a smooth positive energy representation of $\LG$. Then $V$ is
	irreducible as a representation of $\LGtilde$.
\end{proposition}

\begin{proof}
	Assume $V$ is not irreducible as a $\LGtilde$-representation. Projection onto
	a proper $\LGtilde$-invariant summand defines a bounded self-adjoint operator
	$T:V\to V$ which commutes with $\LGtilde$, but (by hypothesis) not with the action
	of $\TTrot$. Choose $R\in\TTrot$; then define for each $n\in \ZZ$ the bounded operator
	\begin{equation}
		T_n = \int_{\TTrot} z^n R_z T R_z^{-1}\, \d z\period
	\end{equation}
	$T_n$ commutes with the action of $\LGtilde$, and $T_n$ sends the weight space
	$V_m$ to $V_{m+n}$. Because $T$ does not commute with $\TTrot$, the
	operator $T_n$ must be nontrivial for at least one $n<0$. Suppose that $m$
	is the lowest energy of $V$ (i.e., the smallest $m$ such that the weight
	space $V_m \neq 0$).\footnote{Because $V$ is positive energy, $m\geq 0$ --- but
	that doesn't matter for now.} Then $T_n(V_m) = 0$ if $n<0$. Since $V$ is
	irreducible as a representation of $\LGtildeplus$, it is generated as a
	representation by $V_m$.  But then $T_n(V) = 0$ for all $n<0$. The adjoint
	to $T_n$ is $T_{-n}$, and so $T_n(V) = 0$ for all $n\neq 0$.

	This implies that $T$ commutes with the action of $\TTrot$, which is a
	contradiction: the $T_n$ are the Fourier coefficients of the loop
	$\Circ\to \End(V)$ sending $z$ to $R_z T R_z^{-1}$, so we find that this loop
	must be constant.  Consequently, $T$ must commute with the action of
	$\TTrot$, as desired. 
\end{proof}

%-------------------------------------------------------------------%
%-------------------------------------------------------------------%
%  A proof sketch of Theorem \ref{main-thm}						    %
%-------------------------------------------------------------------%
%-------------------------------------------------------------------%

\subsection{A proof sketch of \texorpdfstring{\Cref{main-thm}}{Theorem \ref*{main-thm}}}
\label{PS_proof_sketch}

The goal of this section is to go through the proof of \Cref{main-thm}.
As with all proofs in representation theory, we may first reduce to
the irreducible case, thanks to the first part of \cref{like_cpt_Lie}.
%Before arguing the general case, I want to make an
%observation.
\begin{observation}\label{observe}
	Recall that Schur--Weyl duality sets up a one-to-one correspondence between
	representations of $\SU_n$ and representations of the symmetric groups, by
	studying the decomposition of the tensor power $V^{\otimes d}$ of the
	standard representation $V$ under the action of $\Sigma_d$.\index[terminology]{Schur--Weyl duality}
\end{observation}
One may hope that some analogue of \Cref{observe} is true for
representations of loop groups: suppose we could construct a giant
representation of $\LSU_n$ whose $h$-fold tensor product contains all the
irreducible positive energy representations of level $h$, such that this big
representation admits an intertwining action of $\Diffplus(\Circ)$. Then (with a
little bit of work), we would obtain an intertwining action of $\Diffplus(\Circ)$ on
all irreducible positive representations of $\LSU_n$, which would prove \Cref{main-thm} in this particular case. We would like to then reduce from the
case of a general $G$ to the case of $\SU_n$. The Peter--Weyl theorem says that
a simply connected $G$ is a closed subgroup of $\SU_n$ for some $n$, suggesting
that a technique like this might work.\index[terminology]{Peter--Weyl theorem}

Pressley--Segal's approach is similar, but not the same.
\begin{itemize}
	\item Their base case consists not just of $\LSU_n$, but the loop groups of all simply connected,
	simply laced compact Lie groups.\footnote{Recall that $G$ is simply laced if all its nonzero
	roots have the same length; in other words, if the Dynkin diagram\index[terminology]{Dynkin diagram} of $G$
	does not have multiple edges (so the Dynkin diagram is of ADE type). The simple, simply connected, simply laced
	Lie groups are $\SU_n$ for all $n$, $\mathrm{Spin}_n$ for $n$ even, $\mathrm E_6$, $\mathrm E_7$, and
	$\mathrm E_8$.}\index[terminology]{simply laced} In \cite[Lemma 13.4.4]{loop}, they extend from simply laced
	groups to all simply connected Lie groups; the reason they cannot just use an embedding $j\colon
	G\hookrightarrow \SU_n$ is that, given a representation $V$ of $\LGtilde$, Pressley--Segal need not
	just the embedding $j$, but also the condition that there is an irreducible representation $V'$ of the bigger
	group with $V$ a summand in $j^*V'$.
	\item Now assume $G$ is simply connected and simply laced. Instead of constructing a huge tensor product,
	Pressley--Segal reduce to the case of level $1$ representations in a different way. Let $m_n\colon\LG\to\LG$ be
	the map precomposing a loop $\Circ\to G$ with the $n$-th-power map $\Circ\to\Circ$.
	Then \cite[Proposition 9.3.9]{loop} every irreducible representation $V$ of $\LGtilde$ is contained in $m_h^*F$
	for some level $1$ representation $F$. This allows Pressley--Segal to carry the $\Diffplus(\Circ)$-action from
	$F$ to $V$.
	\item Finally, when $G$ is simply laced and $F$ is level $1$, Pressley--Segal construct the
	$\Diffplus(\Circ)$-action directly using the ``blip construction'' \cite[\S 13.2, \S 13.3]{loop}.
\end{itemize}
\begin{remark}
Pressley--Segal write that ``one hopes that a more satisfactory proof of \cref{main-thm} can be
found,'' \cite[p.\ 271]{loop}, so perhaps there's a proof out there that more closely resembles the
Schur--Weyl-style argument.
\end{remark}

%This isn't exactly the
%approach that's taken in \cite{loop}, but something very close to it is in fact
%what is done. I haven't had time to flesh out the details of the approach
%outlined above (nor have I found a reference which takes this
%approach\footnote{Probably because it doesn't work for some obvious reason.} ---
%except for maybe \cite{loop-rep}, but I've gotten very confused over whether
%certain results in that paper are stated for a general compact Lie group or only
%for $\SU_n$), but I think it might be interesting.

Now we will see how the story goes for $\LSU_n$.
\begin{construction}
\label{Fock_constr}
	Let $G = \SU_n$. Define $H \colonequals \Lup^2(\Circ, V)$, where $V$ is the standard
	representation.\index[terminology]{standard representation} Let $\Har^2(\Circ, V)\subseteq H$ denote the \emph{Hardy
	space}\index[terminology]{Hardy space} of
	$\Lup^2$-functions on $\Circ$ with only nonnegative Fourier coefficients, and let
	$P$ denote orthogonal projection of $H$ onto $\Har^2(\Circ, V)$. Then $H = PH \oplus
	P^\perp H$. The \emph{Fock space}\index[terminology]{Fock space} $\Fock_P$ is the Hilbert
	space completion of the alternating algebra:
	\begin{equation}
	\label{Fock_P}
	\Fock_P = \exteriorhat (PH \oplus \overline{P^\perp H}) \cong
	\Directsumhat_{i,j \geq 0} \exterior^i(PH) \oplus \exterior^j(\overline{P^\perp H})\period
	\end{equation}
	Here $\overline V$ denotes the complex conjugate vector space to $V$, and $\exteriorhat$ and
	$\Directsumhat $ denote Hilbert space completions. The Fock space
	turns out to be the ``giant representation'' we were after: it's the
	fundamental representation of $\LSU_n$.
\end{construction}

\begin{remark}[(the Fock space in physics)]
The process of building a Fock space out of a Hilbert space $H$, as in \eqref{Fock_P}, has a quantum-mechanical
interpretation. Suppose that $H$ is the space of states describing the mechanics of a particle: for example,
$\Lup^2(\Circ, \CC)$ corresponds to a particle moving on a circle. The corresponding Fock space is the space of
states for systems with any number of particles. In \cref{Fock_constr}, we used the alternating algebra, which
means that the particles are fermions: the relation $f\wedge f = 0$ is the Pauli exclusion principle, imposing that
two fermions cannot be in the same state. For a bosonic many-body system, one would use the (Hilbert space
completion of the) symmetric algebra.%
\index[terminology]{fermion}%
\index[terminology]{boson}%
\index[terminology]{Pauli exclusion principle}
The process of building a Fock space from a single-particle Hilbert space is called second
quantization.\index[terminology]{second quantization}

In our setting, $\Lup^2(\Circ, V)$ corresponds to a system with a fermion moving on a circle, together with some
kind of $G$-symmetry. The subspace $\exterior^i(PH) \oplus \exterior^j(\overline{P^\perp H})$ consists of $i$ fermionic
particles and $j$ fermionic antiparticles. This explains why we take the conjugate space to $P^\perp H$: it is so
that the antiparticles have positive energy.\index[terminology]{antiparticle}
\end{remark}
A loop on $G$ acts on $H$ by pointwise multiplication, and $f\in \Diffplus(\Circ)$
acts on $H$ by sending $\xi \colon \Circ\to V$ to $\xi(f^{-1}(z)) \cdot
|(f^{-1})'(z)|^{1/2}$. (The square root factor is a normalization factor to
ensure unitarity of the action.) In fact, this gives an action of $\LG\rtimes
\Diffplus(\Circ)$ on $H$, and one can ask when this descends to a projective
representation of $\LG \rtimes \Diffplus(\Circ)$ on the Fock space $\Fock_P$. Segal
wrote down a \emph{quantization condition} for when a unitary operator on $H$ descends to a projective
transformation of $\Fock_P$: namely, $u$ descends to $\Fock_P$ if and only if the commutator $[u,P]$ is
Hilbert--Schmidt.\footnote{Recall that a bounded operator $A$ on a Hilbert space is \textit{Hilbert--Schmidt} if
$\Tr(A^{*} A)$ is finite.} One checks that the action of $\LG \rtimes \Diffplus(\Circ)$ on $H$ satisfies Segal's
quantization criterion, and so descends to a projective representation of $\LG \rtimes \Diffplus(\Circ)$ on the
Fock space $\Fock_P$.

Almost by definition, the action of $\Circ = \TTrot$ on $\Fock_P$ is of positive
energy, and so $\Fock_P$ is a representation of positive energy. 
It turns out that:

\begin{theorem}[{\cites[Section 10.6]{loop}[Chapter I.5]{loop-rep}}]\label{sun-summand}
	The irreducible summands of $\Fock_P^{\otimes h}$ give all the irreducible
	positive energy representations of $\LSU_n$ of level $h$.
\end{theorem}

%\begin{remark}
%	The space denoted $\cal{H}$ in that section of \cite{loop} is the Hilbert space
%	completion of the Fock space $\Fock_P$.
%\end{remark}
%
We will expand on this construction of the irreducible level $h$ representations of $\LSU_n$ in
\cref{segal_sugawara}, when we discuss the Segal--Sugawara construction.

%This approach isn't exactly the one taken in \cite{loop}, but something
%essentially like it is.
The first reduction comes from:
\begin{lemma}[{\cite[Lemma 13.4.3]{loop}}]
\label{summand}
	Let $V$ and $W$ be positive energy representations of $\LGtilde$. Suppose
	that $V$ is irreducible, and that $V\oplus W$ admits an intertwining action
	of $\Diffplus(\Circ)$. Then $V$ admits an intertwining action of $\Diffplus(\Circ)$.
\end{lemma}

We will prove this shortly; first, we will indicate how to use this to prove the
general case.

\begin{remark}\label{reduction}
	It suffices to prove by \cref{summand} that for every irreducible
	positive energy representation $V$ of $\LG$, there is some $G'$ and an
	embedding $i\colon \LG\to \LG'$ where \Cref{main-thm} is true for $G'$, and
	an irreducible representation $V'$ of $\LG'$ such that $V$ is a summand of
	$i^{*} V'$.
\end{remark}

To use this reduction, we first need to establish that \Cref{main-thm} is true for a
class of Lie groups $G$. In fact:

\begin{theorem}
	\Cref{main-thm} is true if $G$ is simple, simply connected, and
	simply laced.
\end{theorem}

The proof of this result is quite similar to that of \cref{sun-summand}:
one constructs the analogue of the Fock space for $\LG$ (which, like in the
$\SU_n$ case, has an intertwining action of $\Diffplus(\Circ)$), and then shows
that every irreducible positive energy representation is a summand of some twist
of this representation of $\LG$. See \cite[\S 13.4]{loop} for more details.

\begin{construction}
	Let $\Omega G$ denote the \emph{based} loop space of $G$, regarded as the
	homogeneous quotient $\LG/G \simeq \LG_\CC/\Lup^+ G_\CC$.
	Since $ G $ is simple any simply connected,
	\begin{equation*}
		\H^2(\Omega G;\ZZ) \cong \H^3(G; \ZZ) \cong \ZZ \comma
	\end{equation*}
	so every integer gives rise to a complex line bundle on
	$\Omega G$. The holomorphic sections $\Gamma$ of the line bundle
	corresponding to the generator is called the \emph{basic representation} of
	$\LG$.\footnote{Of course, the abelian group $\ZZ$ has two generators. Here we have a canonical one: as
	discussed above, we have a canonical generator for $\H^4(\BG;\ZZ)$, hence $\H^3(G;\ZZ)$ via transgression, and
	therefore also for $\H^2(\Omega G;\ZZ)$.}
	\index[terminology]{based loop space}
	\index[notation]{OmegaG@$\Omega G$}
	\index[terminology]{basic representation of $\LG$}
\end{construction}

\begin{example}
	If $G = \SU_n$, then $\Gamma$ is the Fock space
	described above.
\end{example}

Then:

\begin{proposition}[{\cite[Proposition 9.3.9]{loop}}]\label{basic}
	Let $G$ be a simple, simply connected, and simply laced Lie group. Then any
	irreducible positive energy representation of level $h$ of $\LG$ is a summand
	in $i_h^{*} \Gamma$, where $i_h:\LG\to \LG$ is the map induced by the degree
	$h$ map $\Circ\to \Circ$.
\end{proposition}

The level $1$ representation $\Gamma$ admits an intertwining
action of $\Diffplus(\Circ)$ via the ``blip construction.''%
\index[terminology]{blip construction} We will not go into the details here; see \cite[\S 13.3]{loop}. Assuming
this, combining \cref{basic} with \cref{summand} shows that \cref{main-thm} is true for $\LG$ when $G$ is simply
laced (and simple and simply connected).

According to \cref{reduction}, it now suffices to show:
\begin{proposition}\label{new-reduction}
	For every irreducible positive energy representation $V$ of $\LG$, there is a
	simply laced $G'$ and an embedding $i:\LG\to \LG'$, as well as an irreducible
	representation $V'$ of $\LG'$ such that $V$ is a summand of $i^{*} V'$. 
\end{proposition}
This is proved in \cite[Lemma 13.4.4]{loop} in the following manner.

One first classifies all the irreducible representations of $\LG$. Using the loop
group analogue of Schur--Weyl duality worked well when $G = \SU_n$, but that
won't do in the general case. Instead, one utilizes a loop group analogue of
Borel--Weil (see \cite[Section 4.2]{segal-survey}). Recall how this works for finite-dimensional, compact Lie
groups: fix a maximal torus $T$ of $G$, and then, for every antidominant weight $\lambda$ of $T$ (i.e., $\langle
h_\alpha, \lambda\rangle \leq 0$ for every positive root $\alpha$), there is an associated line bundle
$\cal{L}_\lambda$ on $G/T \cong G_\CC/B^+$. The space of holomorphic sections of $\cal{L}_\lambda$ is an
irreducible representation of $G$ of lowest weight $\lambda$, and all irreducible representations of $G$ arise this
way.\index[terminology]{maximal torus}

In the loop group case, one again begins by fixing a maximal torus $T$ of $G$
(one should think of $\TTrot\times T \times \TTce$ as a maximal torus of
$\LG$). Consider the homogeneous space $\LG/T$. There is a fiber sequence
\begin{equation}
	G/T\to \LG/T \to \Omega G \comma
\end{equation}
and the set of isomorphism classes of complex line bundles on $\LG/T$ is
\begin{equation}
	\H^2(\LG/T;\ZZ) \cong \H^2(\Omega G;\ZZ) \oplus \H^2(G/T; \ZZ) = \ZZ \oplus
\widehat{T} \comma
\end{equation}
where $\widehat{T}$ is the character group of $T$. You can prove this using the Serre spectral sequence, which as
usual is easier because $G$ is simple and simply connected. Anyways, we learn that line bundles on $\LG/T$ are
indexed by $(h,\lambda)\in \ZZ\oplus \widehat{T}$.

\begin{theorem}[{(Borel--Weil for loop groups \cite[Theorem 9.3.5]{loop})}]
	One has:
	\begin{itemize}
	\item The space $\Gamma(\cal{L}_{h,\lambda})$ of holomorphic sections is
		zero or irreducible of positive energy of level $h$; moreover, every
		projective irreducible representation of $\LG$ arises this way.
	\item The space $\Gamma(\cal{L}_{h,\lambda})$ is nonzero if and only if
		$(h,\lambda)$ is antidominant,\footnote{Recall that if $G$ is the
		simply laced group $\SU_n$, then the weight lattice is
		$\bigoplus_{1\leq i\leq n+1} \ZZ\chi_i/\ZZ\sum_i \chi_i$, and the
		roots are $\chi_i - \chi_j$ with $i\neq j$. The positive roots,
		corresponding to the usual Borel subgroup of upper-triangular matrices, are
		$\chi_i - \chi_j$ for $i<j$.  Therefore, $(h, \lambda = \lambda_1,
		\cdots, \lambda_n)$ is antidominant if $\lambda$ is antidominant,
		i.e., $\lambda_1 \leq \cdots \leq \lambda_n$, and if $\lambda_n -
		\lambda_1 \leq h$.} i.e.,
		$$0\geq \lambda(h_\alpha) \geq -\frac{h}{2}\langle h_\alpha,
		h_\alpha\rangle$$
		for each positive coroot $h_\alpha$ of $G$. (In particular,
		$\lambda$ is antidominant as a weight of $T\subseteq G$.)
		\index[terminology]{antidominant}%
		\index[terminology]{Borel subgroup}%
		\index[terminology]{coroot}
	\end{itemize}
\end{theorem}
The upshot is that irreducible representations correspond to antidominant
weights. To prove \Cref{new-reduction}, it suffices to show that
all antidominant weights of $\LG$ are restrictions of antidominant weights of
$\LG'$ for some simply laced $G'$. The argument now proceeds case-by-case, as $G$
ranges over all simple simply connected simply laced compact Lie groups. The
proof is not very enlightening, so we will not go into more detail here.

%-------------------------------------------------------------------%
%-------------------------------------------------------------------%
%  Appendix: Random thoughts										%
%-------------------------------------------------------------------%
%-------------------------------------------------------------------%

\begin{remark}[(relationship with Wess--Zumino--Witten theory)]
Segal \cite{segal-book} studies the theory of positive energy representations of $\LG$ from a different
perspective, that of conformal field theory.\index[terminology]{conformal field theory} Specifically, the category
of level $h$ positive energy representations of $\LG$ has the structure of a \textit{modular tensor
category}\index[terminology]{modular tensor category} Given a modular tensor category $\mathsf C$, one can build
\begin{enumerate}
	\item a $3$-dimensional topological field theory $Z_{\mathsf C}$ \cite{RT90, RT91, Wal91, BK01, KL01, BDSV15},
	and
	\item a $2$-dimensional conformal field theory \cite{MS89}.
\end{enumerate}
These two theories are related: the 2d CFT is a boundary theory for the 3d TFT \cite{Wit89, FT12}. When $\fC$ is
the category of level $h$ representations of $\LG$, the TFT is Chern--Simons theory (see \cref{quantum_CS}) and the
CFT is the Wess--Zumino--Witten model (see \cref{quantum_WZW}).\footnote{One might wonder if every modular tensor
category arises in this way, as a category of positive-energy representations of a loop group. This is the
Moore--Seiberg conjecture,\index[terminology]{Moore--Seiberg conjecture} and is open at the time of writing. See,
e.g., \cite{HRW08}.}

You do not need \cref{main-thm} to construct the modular tensor category structure on $\mathsf{Rep}_k(\LG)$, and
the TFT and CFT provide a very large amount of data associated to that structure. It may be possible to coax
\cref{main-thm} out of that extra structure. For example, Segal \cite[\S 12]{segal-book} discusses this for abelian
Lie groups.
\end{remark}

\subsection{OK, but what does this have to do with differential cohomology?}
\label{PS_diffcoh}
There is differential cohomology hiding in the background of the story of central extensions of loop groups. There
are two ways in which it appears: one which is related to the story of on-diagonal differential characteristic
classes built from Chern--Weil theory, and another which relates central extensions to off-diagonal Deligne
cohomology similarly to the discussion of the Virasoro group in \cref{VirasoroAlgebra}. This, together with the
appearance of $\Diffplus(\Circ)$ in the representation theory of loop groups, suggests that loop groups and the
Virasoro group should interact somehow, as we will see in the next chapter.

% first: Chern--Weil, curvature, etc
\subsubsection{The on-diagonal story}
Suppose $G$ is simple and simply connected, so that $\H^4(\BG;\ZZ)$, $\H^3(G;\ZZ)$, and $\H^2(\LG;\ZZ)$ are all
isomorphic to $\ZZ$, and the transgression maps
\begin{equation*}
	\H^4(\BG;\ZZ)\to\H^3(G;\ZZ)\to\H^2(\LG;\ZZ)
\end{equation*}
are isomorphisms.%
\index[terminology]{transgression}
The level $h$ canonically refines to $\hat h\in\Hhat^4(\BunGnabla ;\ZZ)$
(\cref{differential_CW_lift}), and the transgression map refines to a map
\begin{equation*}
	\Hhat^4(\BunGnabla ;\ZZ)\to\Hhat^3(G;\ZZ)
\end{equation*}
\cite[\S 3]{CJMSW05}, as we discussed in \cref{transgression_detail}. Does
the story continue to a differential refinement $\Hhat^3(G;\ZZ)\to\Hhat^2(\LG;\ZZ)$? That is, a projective
representation $\LG\to\PU(V)$ determines a central extension $\LGtilde$ of $\LG$, which is a principal
$\TT$-bundle over $\LG$. Does this $\TT$-bundle come with a canonical connection?

Of course, this is a loaded question, and we'll see that the answer is yes. But first, a (relatively) down-to-Earth
plausibility argument. Given a central extension
% TODO\colon \shortexact?
\begin{subequations}
\begin{equation}
	1\to \TTce\to \LGtilde \to \LG \to 1 \comma
\end{equation}
we can differentiate it to obtain a central extension of Lie algebras
\begin{equation}
\label{Loop_alg_cext}
	0\to \RR\to \Lgtildeno \to \Lg \to 0 \period
\end{equation}
\end{subequations}
Recall from \cref{cext_lie_alg} that the central extension \eqref{Loop_alg_cext} can be described by a cocycle for
the Lie algebra cohomology group\index[terminology]{Lie algebra!cohomology} $\HLie^2(\Lg; \RR)$. Cocycles
are alternating maps $\omega\colon \Lg\times\Lg\to\RR$ satisfying the cocycle condition \eqref{LA_cext_Jacobi}.
Choose a cocycle $\omega$; then, $\Lgtildeno$ is the vector space $\Lg\oplus\RR$ with the Lie bracket
\begin{equation}
	[(\xi, a), (\eta, b)] \colonequals ([\xi, \eta],\omega(\xi, \eta)) \period
\end{equation}
For example, an element of $\H^4(\BG;\RR)$ corresponds via the Chern--Weil machine to an invariant symmetric
bilinear form $\langle-, -\rangle\colon \g\times\g\to\RR$, and it defines a degree-$2$ Lie algebra
cocycle for $\Lg$ by \cite[\S 4.2]{loop}
\begin{equation}
	\omega(\xi, \eta) \colonequals \frac{1}{2\pi}\int_{\Circ} \langle \xi(\theta), \eta'(\theta)\rangle\,\d\theta \period
\end{equation}
Suppose that $\omega$ comes from a central extension of $\LG$ which is a principal $\TT$-bundle $\pi\colon
\LGtilde\to \LG$. 
Then $\Tan\LGtilde$ fits into a short exact sequence
\begin{equation}
\label{loop_ext_tangent_ses}
	0\to \Tan\TT\to \Tan\LGtilde\to \pi^*\Tan\LG\to 0 \period
\end{equation}
At the identity of $\LGtilde$ this is \eqref{Loop_alg_cext}, and left translation carries this identification to
every tangent space. The data of $\omega$ includes a splitting of \eqref{Loop_alg_cext}, and left translation turns
this into a splitting of \eqref{loop_ext_tangent_ses}. A connection on $\pi\colon\LGtilde\to\LG$ is a
$\TT$-invariant splitting, and since $\TT$ acts trivially on its Lie algebra, we have just built a connection with
curvature $\omega$. Thus the class of \eqref{Loop_alg_cext} in $\H^2(\LG;\ZZ)$ refines to a class in
$\Hhat^2(\LG;\ZZ)$. Pressley--Segal \cite[Theorem 4.4.1]{loop} show that this is a necessary and sufficient
condition on $\omega$ for any compact, simply connected Lie group $G$, and that $\omega$ determines the
extension.\footnote{When $G$ is not simply connected, the theorem is not quite as nice: see \cite[Theorem
4.6.9]{loop} and \cite{Wal17}.}

\begin{remark}
	It may be possible to do this ``all at once'' by finding a canonical connection $\conn$ on the principal
	$\TT$-bundle $\pi\colon \Uup(V)\to \PU(V)$ where $V$ is an infinite-dimensional separable Hilbert space;
	this would lift the tautological class
	\begin{equation*}
		c_1(\Uup(V))\in\H^2(\PU(V);\ZZ) = \H^2(\Kup(\ZZ, 2); \ZZ)
	\end{equation*}
	to
	\begin{equation*}
		\chat_1(\Uup(V), \conn) \in \Hhat^2(\PU(V);\ZZ) \period
	\end{equation*}
	Then a projective representation would pull back $\chat_1(\Uup(V), \conn)$ (and $\conn$) to $\LG$.
\end{remark}

To summarize a little differently, given $\hat h\in\Hhat^4(\BunGnabla ;\ZZ)$, we can obtain a Chern--Weil form
$\langle-, -\rangle$, hence a cocycle $\omega\in\HLie^2(\Lg;\RR)$. Because $\curv(\hat h)$
satisfies an integrality condition, so does $\omega$, which turns out to be the same condition needed to define a
central extension $\LGtilde\to\LG$ with a connection. That is, we built a map
$\Hhat^4(\BunGnabla ;\ZZ)\to\Hhat^2(\LG;\ZZ)$. We would like to describe it more directly.

The first step is the transgression map $\Hhat^4(\BunGnabla ;\ZZ)\to\Hhat^3(\BG;\ZZ)$ constructed by \cite[\S
3]{CJMSW05}. To get from $3$ to $2$, Gawędzki \cite[\S 3]{Gaw88} constructs for any closed manifold $M$ a
transgression map\index[terminology]{transgression}
\begin{equation}
	\Hhat^3(M;\ZZ)\to\Hhat^2(\mathrm LM;\ZZ)
\end{equation}
from the perspective that differential cohomology is isomorphic to the
hypercohomology\index[terminology]{hypercohomology} of the Deligne complex\index[terminology]{Deligne complex}%
\footnote{Gawędzki actually works with a different complex,
namely $0\to \TT\to i\Omega^1\to\dotsb\to i\Omega^{n-1}\to 0$, where the map $\TT\to i\Omega^1$ is $\d\circ{\log}$.
This is equivalent to $\Sigma\ZZ(n)$ \cite[Remark 3.6]{BML94}, and the proof is a straightforward generalization of
\cref{cext_lemma}.} 
\begin{equation*}
	0\to\ZZ\to\Omega^0\to\dotsb\to\Omega^{n-1}\to 0 \period
\end{equation*}
Another option is to construct the transgression as follows: first pull back by the evaluation
map $\Circ\times \mathrm LM\to M$, then integrate over the $\Circ$ factor using the map we constructed in
\cref{FiberIntegration}.

% second: Deligne-Beĭlinson
\subsubsection{The off-diagonal story}
In \cref{VirasoroAlgebra}, we saw in \cref{cext_corollary} that central extensions of a Lie group $\Gamma$
(possibly infinite-dimensional) which are principal $\TT$-bundles are classified by $\H^3(\BbulletGamma;\ZZ(1))$.
The central extensions of loop groups we constructed in this chapter are principal $\TT$-bundles. Therefore there
is in principle a way to start with a class $h\in \H^4(\BG;\ZZ)$ and obtain a class $\phi(h)\in\H^3(\Bun_{\LG};\ZZ(1))$, and that is what we are going to do next.

Recall that truncating defines a map of complexes of sheaves of abelian groups $\ZZ(n)\to\ZZ$, inducing for us a
map
\begin{equation}\label{4_to_2}
	\H^4(\BbulletG;\ZZ(2)) \to \H^4(\BbulletG;\ZZ)\isomorphism \H^4(\BG;\ZZ) \period
\end{equation}
%(TODO: make sure the equality in the second map is right. We don't say why this is anywhere, yet.)
\begin{lem}
For $G$ a compact Lie group, \eqref{4_to_2} is an isomorphism.
\end{lem}
\begin{proof}
Recall from \cref{2n_to_n} that \eqref{4_to_2} is part of the pullback square
\begin{equation}
\begin{gathered}
	\begin{tikzcd}[column sep=3.5em]
		\H^{4}(\BbulletG;\ZZ(2))\arrow["\eqref{4_to_2}", r]\arrow[d] & \H^{4}(\BG;\ZZ)\arrow[d]\\
		\Sym^2(\gdual)^G\arrow[r] & \H^{4}(\BG;\RR),
	\end{tikzcd}
\end{gathered}
\end{equation}
where the bottom map is the Chern--Weil map. Since $G$ is compact, the Chern--Weil map is an isomorphism,
so \eqref{4_to_2} is as well.
\end{proof}
Therefore our level $h\in \H^4(\BG;\ZZ)$ is equivalent data to an off-diagonal characteristic class $\tilde h\in
\H^4(\BbulletG;\ZZ(2))$. The next step is the construction of yet another transgression map, this time due to
Brylinski--McLaughlin \cite[\S 5, on p.\ 618]{BML94}:
\begin{equation}
	\H^4(\BbulletG;\ZZ(2))\longrightarrow \H^3(\Bun_{\LG};\ZZ(1)) \period
\end{equation}
Their construction models elements of these two differential cohomology groups simplicially: they identify
$\H^4(\BbulletG;\ZZ(2))$ as the abelian group of equivalence classes of gerbes with a connective structure over a
simplicial manifold model for $\BbulletG$, and $\H^3(\Bun_{\LG};\ZZ(1))$ as equivalence classes of line
bundles over a simplicial model for $\Bun_{\LG}$ (\textit{ibid.}, Theorem 5.7).

We have obtained some class in $\H^3(\Bun_{\LG};\ZZ(1))$ from a level $h\in \H^4(\BG;\ZZ)$, hence some
central extension. That this coincides with the central extension obtained from $h$ by the other methods in this
chapter is due to Brylinski--McLaughlin (\textit{ibid.}, \S 5). See also Brylinski \cite[\S 6.5]{Bry08} for related
discussion and Waldorf \cite[\S 3.1]{Wal10} for another construction of this transgression map.

\newpage
%!TEX root = ../diffcoh.tex

\section{The Segal--Sugawara construction}
\textit{by Peter Haine}
\label{segal_sugawara}

Let $ G $ be a simply connected, simple, compact Lie group with Lie algebra $ \g $.
In \Cref{loop_groups}, we looked at central extensions
\begin{equation*}
	\begin{tikzcd}[sep=1.5em]
		1 \arrow[r] & \Circ \arrow[r] & \LGtilde \arrow[r] & \LG \arrow[r] & 1
	\end{tikzcd}
\end{equation*}
of the \textit{loop group} $ \LG \colonequals \Cinf(\Circ,G) $.
The group $ \Diffplus(\Circ) $ of orientation-preserving diffeomorphisms of the circle acts on $ \LG $ by precomposition.
So we might expect an action of the Virasoro group $ \Vir $ on $ \LGtilde $.
We saw that even though there is not an action of $ \Vir $ on $ \LGtilde $, roughly, the Virasoro group acts on any \textit{positive energy} representation of $ \LGtilde $.
However, the Virasoro action on positive energy representations of $ \LGtilde $ is very inexplicit, and we can
only guarantee the existence of the Virasoro action up to ``essential equivalence,'' which is not actually an
equivalence relation.\index[terminology]{essential equivalence}
In particular, the Pressley--Segal Theorem \cite[Theorem 13.4.3]{loop} (\cref{main-thm}) does not explicitly
explain how the central circle $ \Circ \subset \LGtilde $ acts.

The goal of this chapter is to explain the Lie algebra version of the Pressley--Segal Theorem, which gives an
explicit representation of the Virasoro algebra on any positive energy representation of the \textit{Kac--Moody
algebra} $ \Lgtildeno $ associated to a simple Lie algebra $ \g $ (over the complex numbers).  We'll be able to do
this by writing down explicit universal formulas for ``elements'' of the universal enveloping algebra $
\Univ(\Lgtildeno) $ that satisfy the Virasoro relations.  The catch is that these universal formulas involve
infinite sums, so they do not actually make sense as elements of $ \Univ(\Lgtildeno) $, but they do make sense
whenever we act on a representation where only finitely many of the terms don't act by zero; this is what the
positive energy condition guarantees.

Like in the previous chapter, we are not assuming you're familiar with all of these words. In
\cref{subsec:reminders}, we review some important definitions from \Cref{VirasoroAlgebra}. In \cref{sec:KacMoody},
we define the loop algebra of a Lie algebra, which up to regularity issues is the Lie-algebraic analogue of the
loop group of a Lie group. We also introduce \textit{Kac--Moody algebras}, the analogues of the central extensions
of loop groups we constructed in \cref{rep_loop}. In \cref{sec:SegalSugawara}, we introduce the Segal--Sugawara
construction, first at a high level, then digging into the details.

%-------------------------------------------------------------------%
%-------------------------------------------------------------------%
%  Introduction                                                     %
%-------------------------------------------------------------------%
%-------------------------------------------------------------------%

%-------------------------------------------------------------------%
%  Reminders on Virasoro \& Witt algebras                           %
%-------------------------------------------------------------------%

\subsection{Reminders on Virasoro \& Witt algebras}\label{subsec:reminders}

\begin{definition}
	The (complex) \textit{Witt algebra} is the complex Lie algebra $ \wittCC $ of polynomial vector fields on $
	\Circ $.  Explicitly, $ \wittCC $ has generators $ L_m \colonequals ie^{im\theta} \frac{\dup}{\dup\theta} $ for
	$ m \in \ZZ $ with Lie bracket
	\begin{equation*}
		[L_m,L_n] \colonequals (m-n) L_{m+n}
	\end{equation*}
	for all $ m,n \in \ZZ $.\index[terminology]{Witt algebra!complex}
\end{definition}
This is the complexification of the Witt algebra we discussed in \Cref{Witt_algebra}.

\begin{nul}
	Ignoring regularity issues, the Witt algebra is the complexification of the Lie algebra of the group $
	\Diffplus(\Circ) $ of orientation-preserving diffeomorphisms of the circle.\footnote{For the readers who care
	about regularity: the Lie algebra of $\Diffplus(\Circ)$ is the Lie algebra of all smooth vector fields on
	$\Circ$, and $\wittCC$ is a dense subset of the complexification. See \cites[\S 3.3]{loop}{MO:267249}.}
\end{nul}

\begin{nul}
	Recall from \Cref{cext_lie_alg} that central extensions of Lie algebras are classified by Lie algebra
	cohomology.\index[terminology]{Lie algebra cohomology} We have that
	\begin{equation*}
		\HLie^2(\wittCC;\CC) \isomorphic \CC \comma
	\end{equation*}
	so there is a $ 1 $-dimensional space of central extensions of the Witt algebra.
\end{nul}

\begin{definition}
	The (complex) \textit{Virasoro algebra} $ \virCC $ is the central extension\index[terminology]{Virasoro
	algebra!complex}
	\begin{equation}
	\label{cpxVira}
		\begin{tikzcd}[sep=1.5em]
			1 \arrow[r] & \CC \chg \arrow[r] & \virCC \arrow[r] & \wittCC \arrow[r] & 1
		\end{tikzcd}
	\end{equation}
	of $ \wittCC $ with generators $ L_m $ for $ m \in \ZZ $ and a central element $ \chg $, and nontrivial Lie bracket given by
	\begin{equation}
	\label{Witt_cocycle}
		[L_m,L_n] \colonequals (m-n) L_{m+n} + \delta_{m,-n} \frac{m^3 - m}{12} \chg
	\end{equation}
	for all $ m,n \in \ZZ $.

	We call the central element $ \chg \in \virCC $ the \textit{central charge}.\index[terminology]{central
	charge}\index[notation]{chg@$\chg$}
\end{definition}
Said a little differently,~\eqref{Witt_cocycle} spells out a cocycle for $\HLie^2(\wittCC;\CC)$, which determines
the central extension~\eqref{cpxVira}.

\begin{nul}
	Again, ignoring regularity issues, the Virasoro algebra is the complexification of the Lie algebra of the Virasoro group $ \Vir $.
\end{nul}

%-------------------------------------------------------------------%
%  Talk overview                                                    %
%-------------------------------------------------------------------%

%\subsubsection{Talk overview}\label{subsec:overview}

%\begin{goal}
%	\end{goal}

%-------------------------------------------------------------------%
%-------------------------------------------------------------------%
%  Loop algebras & Kac–Moody algebras                               %
%-------------------------------------------------------------------%
%-------------------------------------------------------------------%

\subsection{Loop algebras and Kac--Moody algebras}\label{sec:KacMoody}

The first thing we need to explain in order to state the Segal--Sugawara construction is what the Kac--Moody algebra $ \Lgtildeno $ is.
As the notation suggests, $ \Lgtildeno $ is the Lie algebra analog of the central extension $ \LGtilde $ of the loop group $ \LG $ (with suitable finiteness hypotheses).
Before talking about Kac--Moody algebras, we need to talk about loop algebras.

%-------------------------------------------------------------------%
%  Loop algebras                                                    %
%-------------------------------------------------------------------%

\subsubsection{Loop algebras}\label{subsec:loops}

\begin{recollection}
	Let $ \g $ be a Lie algebra over a ring $ R $, and let $ S $ be an $ R $-algebra.
	The basechange $ \g \tensor_R S $ of $ \g $ to $ S $ is the Lie algebra over $ S $ with underlying $ S $-module the basechange $ \g \tensor_R S $ of the underlying $ R $-module of $ \g $ to $ S $ with Lie bracket extended from pure tensors from the formula
	\begin{equation*}
		[X_1 \tensor s_1, X_2 \tensor s_2]_{\g \tensor_R S} \colonequals [X_1,X_2]_{\g} \tensor s_1 s_2 \period 
	\end{equation*}
\end{recollection}

\begin{definition}
	Let $ \g $ be a complex Lie algebra.
	The \textit{loop algebra} $ \Lg $ of $ \g $ is the Lie algebra\index[terminology]{loop algebra}
	\index[notation]{Lg@$\Lg$}
	\begin{equation*}
		\Lg \colonequals \g \tensor_{\CC} \CC[t^{\pm 1}] \comma
	\end{equation*}
	regarded as a Lie algebra over $ \CC $ (rather than $ \CC[t^{\pm 1}] $).
\end{definition}

\begin{notation}
	Let $ \g $ be a complex Lie algebra, $ X \in \g $, and $ m $ an integer. 
	We write 
	\begin{equation*}
		X\ang{m} \colonequals X \tensor t^m \in \Lg \period \index[notation]{Xangm@$X\ang m$}
	\end{equation*}
\end{notation}

\begin{nul}
	If $ \{u_i\}_{i \in I} $ is a Lie algebra basis for $ \g $, then $ \{u_i\ang{m} \}_{(i,m) \in I \cross \ZZ} $ is a basis for $ \Lg $.
\end{nul}

\begin{remark}\label{remark:prodloopalg}
	The loop algebra functor $ \Lup \colon \LieAlg_{\CC} \to \LieAlg_{\CC} $ preserves finite products.
\end{remark}

\begin{recollection}
	A finite dimensional Lie algebra $ \g $ is \textit{simple} if $ \g $ is not abelian and the only ideals of $ \g
	$ are $ \g $ and $ 0 $.\index[terminology]{simple Lie algebra}
\end{recollection}

\begin{theorem}[{(Garland \cite[\S\S1 \& 2]{MR601519})}]
	If $ \g $ is a simple Lie algebra over $ \CC $, then
	\begin{equation*}
		\HLie^2(\Lg;\CC) \isomorphic \CC \period 
	\end{equation*}
\end{theorem}

In particular, if $ \g $ is simple there is a $ 1 $-dimensional space of central extensions of $ \Lg $.

%-------------------------------------------------------------------%
%  Recollection on bilinear forms & semisimplicity                  %
%-------------------------------------------------------------------%

\subsubsection{Recollection on bilinear forms \& semisimplicity}\label{subsec:recbilin}

\begin{notation}
	Let $ \g $ be a complex Lie algebra.
	We write $ \ad \colon \fromto{\g}{\End_{\CC}(\g)} $ for the adjoint representation, defined
	by\index[terminology]{adjoint representation!for Lie algebras}
	\begin{equation*}
		\ad(X) \colonequals [X,-] \period
	\end{equation*}
\end{notation}

\begin{example}
	A Lie algebra $ \g $ is abelian if and only if the adjoint representation of $ \g $ is trivial.
\end{example}

\begin{recollection}[(Killing form)]
	Let $ \g $ be a finite-dimensional Lie algebra.
	The \textit{Killing form} on $ \g $ is the bilinear form\index[terminology]{Killing form}
	\begin{align*}
		\Kil_{\g} \colon \g \cross \g &\to \CC \\
		(X,Y) &\mapsto \tr(\ad(X) \of \ad(Y)) \period
	\end{align*}
	The Killing form is symmetric and \textit{invariant} in the sense that
	\begin{equation*}
		\Kil_{\g}([X,Y],Z) = \Kil_{\g}(X,[Y,Z])
	\end{equation*}
	for all $ X,Y,Z \in \g $.
\end{recollection}

\begin{example}
	If $ \g $ is a simple Lie algebra, then every invariant symmetric bilinear form on $ \g $ is a $ \CC $-multiple of the Killing form $ \Kil_{\g} $.
	See \cite{MSE:298287} for a nice exposition of this fact. It is also related to Chern--Weil theory, which tells
	us that the space of invariant symmetric bilinear forms is isomorphic to $\H^4(\BG;\RR)$, and when $G$ is a
	compact, simple, simply connected Lie group, $\H^4(\BG;\RR)\cong\RR$. This is because $\H^4(\BG;\ZZ)\cong\ZZ$,
	which we have discussed and used in previous chapters.
\end{example}

\begin{example}
	Let $ \afrak $ be a finite-dimensional abelian Lie algebra over $ \CC $.
	Since the adjoint representation of $ \afrak $ is trivial, the Killing form of $ \afrak $ is identically zero.
	Also note that every bilinear form on the underlying vector space of $ \afrak $ is an invariant bilinear form on $ \afrak $. 
\end{example}

\begin{proposition}[{\cite[Chapter II, Theorems 2 \& 4]{MR1808366}}]\label{prop:semisimplicity}
	Let $ \g $ be a finite dimensional complex Lie algebra.
	The following conditions are equivalent:
	\begin{enumerate}[label=\stlabel{prop:semisimplicity}, ref=\arabic*]	
		\item\label{prop:semisimplicity.0} The center of $ \g $ is trivial.

		\item\label{prop:semisimplicity.1} The only abelian ideal in $ \g $ is $ 0 $.

		\item\label{prop:semisimplicity.2} The Lie algebra $ \g $ is isomorphic to a product of simple Lie algebras.

		\item\label{prop:semisimplicity.3} \emph{Cartan--Killing criterion:} the Killing form of $ \g $ is
		nondegenerate.\index[terminology]{Cartan--Killing criterion}
	\end{enumerate}
\end{proposition}

\begin{definition}
	Let $ \g $ be a finite-dimensional complex Lie algebra.
	If the equivalent conditions~\enumref{prop:semisimplicity}{0}--\enumref{prop:semisimplicity}{3} are satisfied,
	we say that $ \g $ is \textit{semisimple}.\index[terminology]{semisimple!Lie algebra}
\end{definition}

%-------------------------------------------------------------------%
%  Kac–Moody algebras                                               %
%-------------------------------------------------------------------%

\subsubsection{Kac--Moody algebras}\label{subsec:KacMoody}

Now we define the Lie algebra analogue of the central extensions $ \LGtilde $ of the loop group $ \LG $ that we
studied in \cref{loop_groups}. Those central extensions were parametrized by an element of $\H^4(\BG;\ZZ)$, and
these similarly require the additional data of an invariant symmetric bilinear form on $ \g $, i.e.\ an element of
$\H^4(\BG;\RR)$. The Killing form provides a canonical choice. The forms not in the image of
$\H^4(\BG;\ZZ)\to\H^4(\BG; \RR)$ correspond to loop algebra central extensions which do not lift to loop groups.

\begin{definition}[\cite{Kac68, Moo68}]
	Let $ \g $ be a Lie algebra over $ \CC $ with invariant symmetric bilinear form $ B \colon \fromto{\g \cross \g}{\CC} $.
	The \textit{Kac--Moody algebra} of $ \g $ with respect to the form $ B $ is the central extension
	\index[terminology]{Kac--Moody algebra}
	\begin{equation*}
		\begin{tikzcd}[sep=1.5em]
			1 \arrow[r] & \CC c \arrow[r] & \Lgtilde{B} \arrow[r] & \Lg \arrow[r] & 1
		\end{tikzcd}
	\end{equation*}
	with central element $ c $ and with Lie bracket extended from the relation
	\begin{align*}
		[X\ang{m},Y\ang{n}]_{\Lgtilde{B}} &\colonequals [X\ang{m},Y\ang{n}]_{\Lg} + \delta_{m,-n} m B(X,Y) c \\ 
		&= [X,Y]_{\g}\ang{m+n} + \delta_{m,-n} m B(X,Y) c
	\end{align*} 
	for all $ X,Y \in \g $.
\end{definition}

\begin{nul}
	If $ \{u_i\}_{i \in I} $ is a Lie algebra basis for $ \g $, then $ \{u_i\ang{m} \}_{(i,m) \in I \cross \ZZ} \union \{c\} $ is a basis for $ \Lgtilde{B} $.
\end{nul}

\begin{remark}
	The Kac--Moody algebra $ \Lgtildeno $ is usually denoted by $ \hat{\g} $ and is also known as the
	\textit{affine Lie algebra} of $ \g $.\index[terminology]{affine Lie algebra}
\end{remark}

\begin{remark}\label{remark:prodKacMoodyalg}
	Let $ \g_1 $ and $ \g_2 $ be complex Lie algebras equipped with invariant symmetric bilinear forms
	\begin{equation*}
		B_1 \colon \fromto{\g_1 \cross \g_1}{\CC} \andeq B_2 \colon \fromto{\g_2 \cross \g_2}{\CC} \period
	\end{equation*}
	Write $ B $ for the bilinear form on the product Lie algebra $ \g_1 \cross \g_2 $ defined by
	\begin{equation*}
		B((x_1,x_2),(y_1,y_2)) \colonequals B_1(x_1,y_1) + B(x_2,y_2) \period
	\end{equation*}
	Then we have a canonical isomorphism
	\begin{equation*}
		\widetilde{\Lup}_{B}(\g_1 \cross \g_2) \isomorphic \Lgtilde{B_1}_1 \cross \Lgtilde{B_2}_2 \period
	\end{equation*}
\end{remark}

%-------------------------------------------------------------------%
%-------------------------------------------------------------------%
%  The Segal–Sugawara construction                                  %
%-------------------------------------------------------------------%
%-------------------------------------------------------------------%

\subsection{The Segal--Sugawara construction}\label{sec:SegalSugawara}

We now have enough of of the background on Lie algebras to give a vague statement of the Segal--Sugawara construction.

\begin{definition}\label{def:posenergy}
	Let $ \g $ be a Lie algebra over $ \CC $ and $ B $ an invariant symmetric bilinear form on $ \g $.
	A representation $ \rho \colon \fromto{\Lgtilde{B}}{\End_{\CC}(V)} $ has \textit{positive energy} if for all $
	v \in V $ and $ X \in \g $ there exists an integer $ m > 0 $ such that\index[terminology]{positive energy!for
	Kac--Moody algebra representations}
	\begin{equation*}
		\rho(X\ang{m})v = 0 \period
	\end{equation*}
\end{definition}

\begin{remark}
	In the theory of Kac--Moody algebras, positive energy representations are more often called \textit{admissible}.
	We have chosen the term ``positive energy'' to align with the loop group
	terminology (see \Cref{pos_en}).\index[terminology]{admissible representation}
\end{remark}

\begin{theorem}[(Segal--Sugawara construction, vague formulation)]\label{thm:SegalSugawaravague}
	Let $ \g $ be an abelian or simple Lie algebra over $ \CC $ and let $ B \colon \fromto{\g \cross \g}{\CC} $ be
	a nondegenerate invariant symmetric bilinear form on $ \g $.\index[terminology]{Segal--Sugawara
	construction!vague formulation}
	Write $ \Cas_B(\g) \in \Univ(\g) $ for the Casimir element of $ \g $ with respect to the bilinear form $ B
	$.\index[terminology]{Casimir element}\index[notation]{CasBg@$\Cas_B$}
	Let
	\begin{equation*}
		\rho \colon \fromto{\Lgtilde{B}}{\End_{\CC}(V)}
	\end{equation*} 
	be a positive energy representation of $ \Lgtilde{B} $ such that
	\begin{enumerate}[label=\stlabel{thm:SegalSugawaravague}, ref=\arabic*]
		\item the central element $ c \in \Lgtilde{B} $ acts by multiplication by a complex number $ \ell $,

		\item and the complex number $ -\ell $ is not equal to
		\begin{equation*}
			\lambda_B(G) \colonequals \frac{\tr(\ad(\Cas_{B}(\g)))}{2\dim(\g)} \period
		\end{equation*}
	\end{enumerate}
	Then there is an explicit action of the Virasoro algebra on $ V $ where the central charge $ \chg \in \virCC $
	acts by multiplication by\index[terminology]{Virasoro algebra}
	\begin{equation*}
		\frac{\ell \dim(\g)}{\ell + \lambda_B(\g)} \period
	\end{equation*}

	As special cases:
	\begin{enumerate}[label=\stlabel{thm:SegalSugawaravague}, ref=\arabic*]
		\setcounter{enumi}{2}
		\item If $ \g $ is abelian, then $ \lambda_B(\g) = 0 $ for any nondegenerate invariant symmetric bilinear form $ B $, and the central charge $ \chg \in \virCC $ acts by multiplication by $ \dim(\g) $.

		\item If $ \g $ is simple and $ B $ is the normalization of the Killing form such that the long roots of $
		\g $ have square length $ 2 $, then $ \lambda_B(\g) $ is a positive integer known as the \emph{dual Coxeter
		number} of $ \g $.\index[terminology]{dual Coxeter number}
	\end{enumerate}
\end{theorem}

\begin{nul}
	The complex number $ \ell $ in \Cref{thm:SegalSugawaravague} is known as the \textit{level} of the positive energy representation $ \rho $.
\end{nul}

\begin{goal}
	The goal for the rest of the talk is to explain this construction, the Casimir element $ \Cas_B(\g) $, and give a better description of the normalized trace $ \lambda_B(\g) $ as an eigenvalue of $ \ad(\Cas_B(\g)) $.
\end{goal}

%-------------------------------------------------------------------%
%  Motivating case: the Heisenberg algebra                          %
%-------------------------------------------------------------------%

\subsubsection{Motivating case: the Heisenberg algebra}\label{subsec:Heisenberg}

As motivation for the Segal--Sugawara construction, we start with the most simple case, where $ \g $ is the $ 1 $-dimensional abelian Lie algebra.
Since the constant $ \lambda_B(\g) $ will be zero in this case, we can do this without yet introducing the Casimir element.

\begin{definition}
	The \textit{Heisenberg algebra}\index[terminology]{Heisenberg!algebra}\index[terminology]{Lie algebra!Heisenberg} is the Kac--Moody algebra
	\begin{equation*}
		\heis \colonequals \widetilde{\Lup}\CC \index[notation]{Heis@$\heis$}
	\end{equation*}
	of the $ 1 $-dimensional abelian Lie algebra $ \CC $ with respect to the bilinear form $ \fromto{\CC \cross \CC}{\CC} $ given by multiplication.
\end{definition}

\begin{nul}
	Write $ u \in \CC $ for the element $ 1 $, which we regard as a basis for $ \CC $ as a $ 1 $-dimensional abelian Lie algebra.
	Then the Heisenberg algebra has generators $ \{c\} \union \{u\ang{m}\}_{m \in \ZZ} $, where $ c $ is central and the nontrivial bracket relation is given by
	\begin{equation*}
		[u\ang{m},u\ang{n}] \colonequals \delta_{m,-n}m c \period
	\end{equation*}
\end{nul}

\begin{definition}
	Let $ \mu, \hbar \in \CC $.
	Write $ u \in \CC $ for the element $ 1 $, which we regard as a basis for $ \CC $ as a $ 1 $-dimensional abelian Lie algebra.
	The \textit{Fock representation} $ \Fock(\mu,\hbar) $ is the representation of the Heisenberg algebra on the
	polynomial ring\index[terminology]{Fock representation}
	\begin{equation*}
		\Fock(\mu,\hbar) \colonequals \CC[x_1,x_2,\ldots]
	\end{equation*}
	in infinitely many variables, where
	\begin{align*}
		c &\mapsto \hbar \id{} \\ 
		u\ang{n} &\mapsto \begin{cases}
			\frac{\partial}{\partial x_n} \comma & n > 0 \\
			-\hbar x_{-n} \comma & n < 0 \\ 
			\mu \id{} \comma & n = 0 \period
		\end{cases} \\
	\end{align*}
\end{definition}

The following fact about the irreducibility of Fock representations is easy:

\begin{lemma}[{\cite[Lemma 2.1]{MR1021978}}]
	Let $ \mu,\hbar \in \CC $.
	If $ \hbar \neq 0 $, then the $ \heis $-representation $ \Fock(\mu,\hbar) $ is irreducible.
\end{lemma}

\begin{nul}
	If $ \hbar = 0 $, then the constants $ \CC \subset \Fock(\mu,0) $ are invariant.
\end{nul}

\begin{properties}\label{properties:Fock}
	The following are some important properties of the Fock representations of the Heisenberg algebra.
	\begin{enumerate}[label=\stlabel{properties:Fock}, ref=\arabic*]
		\item\label{properties:Fock.1} The elements $ u(0) $ and $ c $ of $ \heis $ act by multiplication.

		\item\label{properties:Fock.2} For every polynomial $ p \in \Fock(\mu,\hbar) $, there exists an integer $ n \gg 0 $ such that $ u\ang{n} p = 0 $: let $ n $ be any positive such that the variable $ x_n $ does not appear in $ p $.
		That is, the Fock representation $ \Fock(\mu,\hbar) $ is ``positive energy'' in the sense of \Cref{def:posenergy}.

		\item\label{properties:Fock.3} For each integer $ n > 0 $, the element $ u\ang{n} \in \heis $ acts locally nilpotently on $ \Fock(\mu,\hbar) $.
	\end{enumerate}
\end{properties}

Now we can give the Segal--Sugawara construction for the Fock representations of the Heisenberg algebra.

\begin{construction}[(Virasoro action of Fock representations)]
	For each integer $ m \in \ZZ $, define an infinite sum of elements of $ \Univ(\heis) $ by
	\begin{equation*}
		L_m^S \colonequals \frac{1}{2} \sum_{j \in \ZZ} \normord{u\ang{-j} u\ang{j+m}} \period
	\end{equation*}
	Here, $ \normord{u\ang{-j} u\ang{j+m}} $ denotes the \textit{normal ordering} on $ u\ang{-j} u\ang{j+m} $,
	defined by\index[terminology]{normal ordering}\index[notation]{uangj@$\normord{u\ang{-j}u\ang{j+m}}$}
	\begin{equation*}
		\normord{u\ang{-j} u\ang{j+m}} \colonequals \begin{cases}
			u\ang{-j} u\ang{j+m} \comma & -j \leq j + m \\ 
			u\ang{j+m} u\ang{-j} \comma & -j \geq j + m \period
		\end{cases}
	\end{equation*}
	Explicitly,
	\begin{equation*}
		L_m^S = \begin{cases}
			\displaystyle \frac{1}{2} u\ang{n}^2 +  \sum_{j > 0 } u\ang{n - j} u\ang{n + j} \comma & m = 2n \\ 
			\displaystyle \sum_{j > 0} u\ang{n + 1 - j} u\ang{n + j} \comma & m = 2n + 1 \period
		\end{cases}
	\end{equation*}
	
	The operators $ L_m^S $ are not well-defined elements of $ \Univ(\heis) $, but since the Fock representations of $ \heis $ are positive energy \enumref{properties:Fock}{2}, the operators $ L_m^S $ make sense as operators on $ \Fock(\mu,\hbar) $.
\end{construction}

\begin{theorem}[{(Segal--Sugawara for $ \Fock(\mu,1) $ \cite[Proposition 2.3]{MR1021978})}]
	Under the representation of $ \heis $ on the Fock space $ \Fock(\mu,1) $, the operators $ L_{m}^S $ on $ \Fock(\mu,1) $ satisfy the commutation relation
	\begin{equation*}
		[L_m^S,L_n^S] = (m - n)L_{m+n}^S + \delta_{m,-n} \frac{m^3 - m}{12} \period
	\end{equation*}
	Hence the assignment
	\begin{align*}
		\virCC &\to \End_{\CC}(\Fock(\mu,1)) \\ 
		L_m &\mapsto L_{m}^S \\ 
		\chg &\mapsto \id{}
	\end{align*}
	is a $ \virCC $-representation with central charge $ 1 $.
\end{theorem}

\begin{remark}
	To derive the Segal--Sugawara action on $ \Fock(\mu,\hbar) $ for $ \hbar \neq 0 $, let $ L_{m} $ act by $ \frac{1}{\hbar} L_{m}^S $.
\end{remark}

\begin{remark}
	Gordon's notes \cite{Gordon:InfiniteLie} give a nice exposition of the Segal--Sugawara construction for Fock representations and the representation theory of the Virasoro algebra.
\end{remark}

%-------------------------------------------------------------------%
%  The Casimir element                                              %
%-------------------------------------------------------------------%

\subsubsection{The Casimir element}\label{subsec:Casimir}

In the general case, the idea is to try to mimic the formulas that we wrote down defining the operators on the Fock representations that satisfy the Virasoro relations.
First, we need to explain the ``Casimir element'' and normalized trace $ \lambda_B(\g) $ appearing in \Cref{thm:SegalSugawaravague}.

\begin{definition}\label{rec:Casimir}
	Let $ \g $ be a finite-dimensional Lie algebra over $ \CC $ and let $ B $ be a nondegenerate invariant symmetric bilinear form on $ \g $.
	The \textit{Casimir element} $ \Cas_{B}(\g) $ of $ \g $ with respect to the form $ B $ is the element of the
	universal enveloping algebra $ \Univ(\g) $ given by the image of $ \id{\g} $ under the
	composite\index[terminology]{Casimir element}\index[notation]{Casb@$\Cas_B$}
	\begin{equation*}
		\begin{tikzcd}[sep=1.5em]
			\End_{\CC}(\g) \isomorphic \g \tensor_{\CC} \gdual \arrow[r, "\sim"{yshift=-0.25em}] & \g \tensor_{\CC} \g \arrow[r, hook] & \Tup_{\CC}(\g) \arrow[r, ->>] & \Univ(\g) \period
		\end{tikzcd}
	\end{equation*}
	Here the isomorphism $ \isomto{\g \tensor_{\CC} \gdual}{\g \tensor_{\CC} \g} $ is the identity on the first factor and the isomorphism $ \isomto{\gdual}{\g} $ induced by the form $ B $ on the second factor, and $ \Tup_{\CC}(\g) $ is the tensor algebra of $ \g $ over $ \CC $.
\end{definition}

The following are some key properties that we need to know about the Casimir element:

\begin{properties}\label{properties:Casimir}
	% Let $ \g $ be a finite-dimensional Lie algebra over $ \CC $ and let $ B $ be a nondegenerate invariant symmetric bilinear form on $ \g $.
	\hfill
	\begin{enumerate}[label=\stlabel{properties:Casimir}, ref=\arabic*]
		\item\label{properties:Casimir.1} The Casimir element $ \Cas_{B}(\g) $ is a central element of $ \Univ(\g) $.

		\item\label{properties:Casimir.2} If $ \{u_1,\ldots,u_d\} $ and $ \{u^1,\ldots,u^d\} $ are bases of $ \g $ that
		are dual with respect to the bilinear form $ B $ in the sense that $B(u_i,u^j) = \delta_{i,j}$, then
		\[\Cas_{B}(\g) = \sum_{i=1}^d u_i u^i\period\]
		\item\label{properties:Casimir.3} Assume that $ \g $ is simple. 
		Then the Casimir element of the Killing form of $ \g $ acts by the identity in the adjoint representation. 
		Hence for any nondegenerate invariant symmetric bilinear form $ B $ on $ \g $, the Casimir element $ \Cas_{B}(\g) $ acts by scalar multiplication in the adjoint representation of $ \g $.
		If $ B $ is the normalization of the Killing form on $ \g $ such that long roots have square length $ 2 $, then in the adjoint representation $ \Cas_{B}(\g) $ acts by multiplication by an even positive integer.

		\item\label{properties:Casimir.4} If $ \g $ is abelian, then since the adjoint representation of $ \g $ is trivial, for any nondegenerate invariant symmetric bilinear form $ B $ on $ \g $ we have $ \ad(\Cas_{B}(\g)) = 0 $.
		In particular, in the adjoint representation $ \Cas_{B}(\g) $ acts by scalar multiplication. 
	\end{enumerate}
\end{properties}

Even though there are no Lie algebras that are both abelian and simple, it is important for us that both types of Lie algebras have the property that the Casimir element associated to any nondegenerate invariant symmetric bilinear form acts by scalar multiplication in the adjoint representation.
In particular, if $ \g $ is abelian or simple, then $ \ad(\Cas_{B}(\g)) $ only has exactly one eigenvalue.

\begin{definition}
	Let $ \g $ be a finite dimensional abelian or simple Lie algebra over $ \CC $ and let $ B \colon \fromto{\g \cross \g}{\CC} $ be a nondegenerate invariant symmetric bilinear form on $ \g $.
	Define a complex number $ \lambda_B(\g) $ by 
	\begin{equation*}
		\lambda_B(\g) \colonequals \frac{1}{2} \Big(\text{eigenvalue of } \ad(\Cas_B(\g)) \Big) \period
	\end{equation*}
\end{definition}

\begin{nul}
	If $ \dim(\g) > 0 $, then 
	\begin{equation*}
		\lambda_B(\g) = \frac{\tr(\ad(\Cas_{B}(\g)))}{2\dim(\g)} \comma
	\end{equation*}
	which aligns with the vague formulation of the Segal--Sugawara construction (\Cref{thm:SegalSugawaravague}).
\end{nul}

\begin{example}
	If $ \g $ is simple and $ B $ is the normalization of the Killing form on $ \g $ such that long roots have square length $ 2 $, then $ \lambda_B(\g) $ is a positive integer \enumref{properties:Casimir}{3} known as the \textit{dual Coxeter number} of $ \g $.
\end{example}

\begin{example}
	If $ \afrak $ is an abelian Lie algebra, then for any nondegenerate invariant symmetric bilinear form $ B $ on $ \afrak $, we have $ \lambda_{B}(\afrak) = 0 $.
\end{example}

%-------------------------------------------------------------------%
%  The general case                                                 %
%-------------------------------------------------------------------%

\subsubsection{The general case}\label{subsec:SSgeneral}

Now let us try using ``the same'' formula to write down a Virasoro action on positive energy representations of $ \Lgtildeno $ as we did for the Heisenberg algebra.
The first modification is that we need to sum over a basis of $ \g $.

\begin{construction}
	Let $ \g $ be a finite-dimensional Lie algebra over $ \CC $ and let $ B $ be a nondegenerate invariant symmetric bilinear form on $ \g $.
	Given a positive energy representation $ \rho \colon \fromto{\Lgtilde{B}}{\End_{\CC}(V)} $, for each integer $ m \in \ZZ $ define
	\begin{equation*}
		T_m^{\rho} \colonequals \frac{1}{2} \sum_{i = 1}^d \sum_{j \in \ZZ} \normord{\rho(u_i\ang{-j}) \rho(u^i\ang{j+m})} \, \in \End_{\CC}(V) \period
	\end{equation*}
	Note that even though the formula defining $ T_m^{\rho} $ involves an infinite sum, since $ \rho $ is a positive energy representation, for each $ v \in V $, all but finitely many terms in the sum defining $ T_m^{\rho} $ annihilate $ v $.
	Hence $ T_m^{\rho} $ is well-defined as an element of $ \End_{\CC}(V) $. 
\end{construction}

We used the letter ``$ T $'' instead of ``$ L $'' because the commutation relation is not quite right:

\begin{lemma}[{\cite[Theorem 10.1]{MR1021978}}]
	Let $ \g $ be a finite dimensional abelian or simple Lie algebra over $ \CC $ and let $ B $ be a nondegenerate invariant symmetric bilinear form on $ \g $.
	For every positive energy representation $ \rho \colon \fromto{\Lgtilde{B}}{\End_{\CC}(V)} $, we have the following commutation relation in $ \End_{\CC}(V) $:
	\begin{align*}
		[T_m^{\rho},T_n^{\rho}] &= (\rho(c) + \lambda_B(\g))(m - n) T_{m+n}^{\rho} \\ 
		&\phantom{=} \qquad + \delta_{m,-n} \dim(\g) \frac{m^3 - m}{12} \rho(c) (\rho(c) + \lambda_B(\g))  \period 
	\end{align*}
\end{lemma}

\begin{idea}
	The naive guess that the operators $ T_m^{\rho} $ satisfy the Virasoro relations is not correct. 
	However, if we could invert $ \rho(c) + \lambda_B(\g) $, then the operators
	\begin{equation*}
		\frac{1}{\rho(c) + \lambda_B(\g)} T_m^{\rho} 
	\end{equation*}
	would satisfy the Virasoro relations.
	We can do this provided that the central element $ c \in \Lgtilde{B} $ acts by a scalar $ \ell $ on $ V $, and $ \ell \neq -\lambda_B(\g) $.
\end{idea}

\begin{theorem}[{(Segal--Sugawara construction \cite[Corollary 10.1]{MR1021978})}]\label{thm:SegalSugawara}
	Let $ \g $ be a finite dimensional abelian or simple Lie algebra over $ \CC $ and let $ B \colon \fromto{\g
	\cross \g}{\CC} $ be a nondegenerate invariant symmetric bilinear form.\index[terminology]{Segal--Sugawara
	construction}
	Let
	\begin{equation*}
		\rho \colon \fromto{\Lgtilde{B}}{\End_{\CC}(V)}
	\end{equation*} 
	be a positive energy representation of $ \Lgtilde{B} $ such that
	\begin{enumerate}[label=\stlabel{thm:SegalSugawara}, ref=\arabic*]
		\item the central element $ c \in \Lgtilde{B} $ acts by multiplication by a complex number $ \ell $,

		\item and $ \ell \neq - \lambda_B(\g) $.  
	\end{enumerate}
	Choose bases $ \{u_1,\ldots,u_d\} $ and $ \{u^1,\ldots,u^d\} $ of $ \g $ that are dual with respect to the bilinear form $ B $.

	Then the assignment
	\begin{equation*}
		L_m \mapsto L_m^{\rho} \colonequals \frac{1}{2(\ell + \lambda_B(\g))} \sum_{i=1}^d \sum_{j \in \ZZ} \normord{\rho(u_i\ang{-j}) \rho(u^i\ang{j+m})}
	\end{equation*}
	extends to a $ \virCC $-representation on $ V $ with central charge
	\begin{equation*}
		\frac{\ell \dim(\g)}{\ell + \lambda_B(\g)} \period
	\end{equation*}
	That is, in $ \End_{\CC}(V) $, the operators $ L_m^{\rho} $ satisfy the commutation relation
	\begin{equation*}
		[L_m^{\rho},L_n^{\rho}] = (m - n) L_{m+n}^{\rho} + \delta_{m,-n} \frac{m^3 - m}{12} \frac{\ell \dim(\g)}{\ell + \lambda_B(\g)} \period 
	\end{equation*}
\end{theorem}

\begin{remark}
	For $ a,b \in \ZZ $, the sum $ \sum_{i=1}^d u_i(a)u^i(b) $ is independent of the choice of basis $ \{u_1,\ldots,u_d\} $ of $ \g $.
	In particular, the operators $ L_m^{\rho} $ are independent of the choice of basis.
\end{remark}

\begin{remark}
	If $ \ell = -\lambda_{B}(\g) $, then the formulas we wrote down for the Segal--Sugawara operators $ L_m^{\rho}
	$ do not make sense, and there is a fundamental difficulty in dealing with the ``critical level'' $ \ell =
	-\lambda_{B}(\g) $.\index[terminology]{critical level}
	At the critical level, the theory seems to resemble the positive characteristic situation rather than the classical one; see \cite{MO:25592} for some discussion of this point. 
\end{remark}

\begin{remark}
	In light of \Cref{remark:prodKacMoodyalg}, the Segal--Sugawara construction can be extended to the case where $ \g $ is \textit{reductive}, i.e., $ \g $ decomposes as a product
	\begin{equation*}
		\g \isomorphic \afrak \cross \g_1 \cross \cdots \cross \g_r \comma
	\end{equation*}
	where $ \afrak $ is an abelian Lie algebra and $ \g_1,\ldots,\g_r $ are simple Lie algebras.
	In this case, the central charge of the resulting $ \virCC $-representation is 
	\begin{equation*}
		\dim(\afrak) + \sum_{i=1}^r \frac{\ell_i \dim(\g_i)}{\ell_i + \lambda_{B_i}(\g_i)} \period
	\end{equation*}
	Here the central element of $ \widetilde{\Lup}\afrak $ acts by multiplication by a nonzero complex number and the central element of each $ \Lgtildeno_i $ acts by multiplication by $ \ell_i \in \CC \smallsetminus \{ - \lambda_{B_i}(\g_i) \} $. 
	This is rather useful as all of the classical Lie algebras are reductive \cite[Theorem 5.49]{MR2440737}; see \cite[Remark 10.3]{MR1021978} for details.
\end{remark}

\begin{remark}
	The Segal--Sugawara construction is usually stated with the assumptions that $ \g $ is simple and $ B $ is the normalization of the Killing form such that the long roots of $ \g $ have square length $ 2 $ (so that $ \lambda_{B}(\g) $ is the \textit{dual Coxeter number}, often denoted by $ h^{\vee} $).
	This is somewhat unfortunate; because the Killing form of an abelian Lie algebra is trivial, to include the
	abelian case (and the reductive extension) the ``usual'' statement needs to be modified to include arbitrary nondegenerate invariant symmetric bilinear forms as in \Cref{thm:SegalSugawaravague}.
\end{remark}

\begin{remark}
	One of the motivations for the formula for the Segal--Sugawara operators $ L_m^{\rho} $ comes from the theory of vertex algebras.
	See \cite[\S3]{MR2079371}, in particular \cite[Proposition 3.3.1]{MR2079371}, for more details on the relation
	to vertex algebras.\index[terminology]{vertex algebra}
\end{remark}

%-------------------------------------------------------------------%
%-------------------------------------------------------------------%
%-------------------------------------------------------------------%
%  References                                                       %
%-------------------------------------------------------------------%
%-------------------------------------------------------------------%
%-------------------------------------------------------------------%

\newpage

\DeclareFieldFormat{labelalphawidth}{#1} % Removes square brackets
\DeclareFieldFormat{shorthandwidth}{#1}
\printbibliography[heading=references]

%-------------------------------------------------------------------%
%-------------------------------------------------------------------%
%-------------------------------------------------------------------%
%  Indices                                                          %
%-------------------------------------------------------------------%
%-------------------------------------------------------------------%
%-------------------------------------------------------------------%

\newpage
\printindex[notation]

\newpage
\printindex[terminology]

\end{document}